\documentclass[11pt,a4paper]{article}
\usepackage[T1]{fontenc}
\usepackage[utf8]{inputenc}
\usepackage[english]{babel}
\usepackage{lmodern}
\usepackage{amsmath,amssymb,amsthm}
\usepackage[margin=25mm,headheight=14pt]{geometry}
\usepackage{booktabs,longtable,array}
\usepackage{enumitem,needspace,listings}
\usepackage{xcolor}
\usepackage[colorlinks=true,linkcolor=blue!45!black,urlcolor=blue!45!black,citecolor=blue!45!black]{hyperref}
\hypersetup{pdftitle={Nonsolvable critical groups with the spectrum of S10},pdfauthor={D. V. Lytkina and V. D. Mazurov},pdfsubject={Nonsolvable finite groups critical with respect to the spectrum of S10},pdfkeywords={finite group, spectrum, symmetric group, solvable radical, almost simple group}}
\newtheorem{theorem}{Theorem}
\newtheorem{lemma}[theorem]{Lemma}
\newtheorem{proposition}[theorem]{Proposition}
\newtheorem{corollary}[theorem]{Corollary}
\theoremstyle{definition}
\newtheorem{definition}[theorem]{Definition}
\newtheorem*{remark}{Remark}
\newtheorem*{classical}{Classical result used below}
\DeclareMathOperator{\Soc}{Soc}
\DeclareMathOperator{\Aut}{Aut}
\DeclareMathOperator{\Out}{Out}
\DeclareMathOperator{\sgn}{sgn}
\DeclareMathOperator{\lcm}{lcm}
\newcommand{\Om}{\Omega}
\newcommand{\F}{\mathbb F}

\setlist[enumerate,1]{label=\arabic*),leftmargin=*,itemsep=3pt}
\setlist[enumerate,2]{label=\alph*),leftmargin=*}
\title{Nonsolvable critical groups with the spectrum of \(S_{10}\)}
\author{D.\,V.~Lytkina\thanks{The work of D.V. Lytkina was supported by the Russian Science Foundation (grant no. 25-41-10008).}\\[3pt]
\small Novosibirsk State University\\
\small 1 Pirogova St., Novosibirsk, 630090, Russia\\
\small \texttt{daria.lytkin@gmail.com}
\and
V.\,D.~Mazurov\thanks{The work of V.D. Mazurov was carried out within the framework of the state assignment of the Sobolev Institute of Mathematics, project FWNF-2026-0017.}\\[3pt]
\small Sobolev Institute of Mathematics\\
\small 4 Koptyug Ave., Novosibirsk, 630090, Russia\\
\small \texttt{vic.mazurov@gmail.com}}
\date{}
\begin{document}
\maketitle
\begin{abstract}
We classify the nonsolvable finite groups that are critical with respect to the spectrum of the symmetric group \(S_{10}\). Up to isomorphism, there are exactly two such groups: \(S_{10}\) and \((A_9\times C_3)\rtimes C_2\), where the involution acts on \(A_9\) by conjugation by a transposition and on \(C_3\) by inversion. The quotient of any nonsolvable group with this spectrum by its solvable radical is almost simple. All possible simple socles are treated without assuming that the radical is nilpotent. The final cases, \(L_3(4)\) and \(A_6\), are excluded using the canonical subgroup \(K=[R(G)\cap G^\infty,G^\infty]\), actions on chief factors, and the required element orders \(9,20,30\). In the \(A_6\) case, the final contradiction is the existence of an element of order \(16\). Reproducible certificates accompany the character-table and matrix calculations.
\end{abstract}
\tableofcontents
\clearpage
\section{Introduction}
The \emph{spectrum} of a finite group \(G\) is the set \(\omega(G)\) of its element orders. We consider the spectrum of the symmetric group of degree \(10\):
\[
 \Om=\omega(S_{10})=\{1,2,3,4,5,6,7,8,9,10,12,14,15,20,21,30\}.
\]
A group with spectrum \(\Om\) is called \(\Om\)-critical if none of its proper sections has this spectrum. Sections and the remaining notation are defined below.
\begin{theorem}\label{thm:final-unconditional-classification}
Up to isomorphism, there are exactly two nonsolvable finite \(\omega(S_{10})\)-critical groups:
\[
 S_{10}\qquad\text{and}\qquad H=(A_9\times C_3)\rtimes C_2.
\]
In the definition of \(H\), a generator of \(C_2\) acts on \(A_9\) by conjugation by a transposition in \(S_9\), and on \(C_3\) by inversion.
\end{theorem}
The proof proceeds as follows. First, we calculate the spectra of the two groups in the statement and prove that they are critical. Next, the exclusions \(25,35\notin\Om\) show that the quotient of any nonsolvable group with spectrum \(\Om\) by its solvable radical is almost simple. We then treat every possible simple socle. The socles \(A_9,A_{10}\) lead to the groups in the statement; the other socles are excluded without assuming criticality.

To analyze the solvable radical, we introduce the perfect residual \(L=G^\infty\), its radical \(P=R(G)\cap L\), and the characteristic subgroup \(K=[P,L]=[K,L]\). The quotient \(L/K\) is a perfect central cover of the simple socle. We study the joint actions of the radical and the perfect residual on chief factors inside \(K\). Choosing an uppermost chief factor on which \(K\) acts nontrivially gives a nilpotent normal subgroup in the induced action group. Theorems on minimal polynomials and critical subgroups then control its structure. The final contradictions use both the required orders in \(\Om\) and the absence of the orders \(16,18,24,28,35,42,45,63\).

All lemmas, reductions, and proofs needed for this argument are included here. Classical theorems used without proof are stated explicitly. The data on simple groups, their central covers, and modules are given to the extent needed, with references to published sources. The appendices contain verification calculations. The proof of the main theorem is completed in Section~\ref{sec:classification-proof}.

\section{Notation and preliminary facts}

All groups considered in this text are finite. Groups are called isospectral if they have the same spectrum; a group is subspectral to another if its spectrum is contained in that of the other group. We denote by \(\omega(G)\) the set of element orders of \(G\), by \(R(G)\) its largest solvable normal subgroup, and by \(\Soc(G)\) the subgroup generated by all minimal normal subgroups of \(G\). We write conjugation as \(a^x=x^{-1}ax\) and commutators as \([a,b]=a^{-1}b^{-1}ab\). For subgroups \(A,B\), the notation \([A,B]\) denotes the subgroup generated by all commutators \([a,b]\), \(a\in A\), \(b\in B\).

\begin{definition}
A section of a group \(G\) is a quotient group \(K/L\), where \(L\trianglelefteq K\leqslant G\). The section is proper if \((K,L)\ne(G,1)\). For a set of positive integers \(\Om\), a group \(G\) is called \(\Om\)-critical if \(\omega(G)=\Om\), but the spectrum of every proper section differs from \(\Om\).
\end{definition}

\begin{lemma}\label{lem:sections}
For every section \(K/L\) of a finite group \(G\),
\[
 \omega(K/L)\subseteq\omega(K)\subseteq\omega(G).
\]
Consequently, to prove that a group \(G\) with spectrum \(\Om\) is \(\Om\)-critical, it suffices to exclude:
\begin{enumerate}
\item all proper subgroups \(K<G\) with \(\omega(K)=\Om\);
\item all proper quotients \(G/L\), \(1\ne L\trianglelefteq G\), with \(\omega(G/L)=\Om\).
\end{enumerate}
\end{lemma}
\begin{proof}
The second inclusion is obvious. Suppose that \(xL\) has order \(d\) and \(x\) has order \(e\). Then \(d\mid e\), and \(x^{e/d}\) has order \(d\). Hence \(d\in\omega(K)\). If the section is proper, then either \(K<G\), or \(K=G\) and \(L\ne1\), proving the final assertion.
\end{proof}

A normal subgroup is called minimal normal if it is nontrivial and contains no proper nontrivial normal subgroup of the ambient group. A chief factor of \(G\) is a quotient \(N/M\), where \(M,N\trianglelefteq G\), \(M<N\), and there is no normal subgroup of \(G\) strictly between them. A solvable chief factor is elementary abelian and is regarded as a simple module for the group induced by conjugation on that factor. We write \(F(A)\) for the Fitting subgroup of a solvable group \(A\), that is, its largest nilpotent normal subgroup; \(O_p(A)\) is its largest normal \(p\)-subgroup. The subgroup \(\Phi(A)\) is the intersection of the maximal subgroups of \(A\). For a \(p\)-group it equals \(A'A^p\). If \(U\le A\) and \(U\Phi(A)=A\), then \(U=A\): otherwise \(U\) would be contained in a maximal subgroup, which also contains \(\Phi(A)\). For an abelian \(p\)-group, \(\Omega_1(A)\) denotes the subgroup of elements whose orders divide \(p\). The exponent \(\exp A\) is the least common multiple of the element orders of \(A\). A group \(E\) is called an extraspecial \(p\)-group if \(Z(E)=E'=\Phi(E)\) has order \(p\).

The derived series is defined by \(G^{(0)}=G\), \(G^{(i+1)}=[G^{(i)},G^{(i)}]\). Its stable term is denoted by \(G^\infty\) and called the perfect residual. A group \(A\) is perfect if \(A=A'\). A perfect central cover of a simple group \(S\) is a perfect group \(D\) with \(D/Z(D)\cong S\). The abbreviated notation \(m.S\), \(S.A\) for covers and outer extensions is specified whenever it is used. We write \(\pi(A)\) for the set of prime divisors of \(|A|\), and \(\mu(A)\) for the elements of \(\omega(A)\) maximal under divisibility.

By a simple type we mean an isomorphism class of simple modules over the specified field, identified with a chosen representative when convenient.

In the analysis of the radical, we use the global notation
\[
 L=G^\infty,\qquad P=R(G)\cap L,\qquad K=[P,L].
\]
In individual local lemmas, the letters \(P,Q,H,T\) may denote Sylow subgroups, kernels, or groups acting on modules; their meanings are specified in the corresponding statements. The notation \(H=(A_9\times C_3)\rtimes C_2\) refers to the group in the main theorem. In local arguments, the same letter may denote a Hall subgroup; each such use is introduced explicitly.

\subsubsection{Classical results}
The following classical results are used without reproducing their general proofs. The versions needed are stated here or immediately before their first specific application.

\begin{classical}[Hall's and Schur--Zassenhaus theorems]
For every set of primes \(\pi\), a finite solvable group has Hall \(\pi\)-subgroups, and all such subgroups are conjugate. A subgroup \(A\le G\) is called a Hall subgroup if \((|A|,[G:A])=1\). If \(A\trianglelefteq G\) is a Hall subgroup, then it has a complement \(B\) in \(G\): \(G=AB\), \(A\cap B=1\). If at least one of \(A,G/A\) is solvable, all complements are conjugate by elements of \(A\).
\end{classical}

\begin{classical}[Maschke's and Clifford's theorems and Schur's lemma]
If the characteristic of a field \(k\) does not divide the order of a finite group \(A\), every \(kA\)-module is completely reducible. If \(N\trianglelefteq G\), the restriction of a simple \(kG\)-module to \(N\) is completely reducible, and \(G\) permutes its isotypic components transitively. The endomorphism ring of a simple module is a division ring; over an algebraically closed field, every endomorphism of a finite-dimensional simple module is scalar.
\end{classical}

\begin{classical}[Gasch\"utz's complement criterion]
Let \(A\trianglelefteq G\) be an abelian \(p\)-subgroup, and let \(U\) be a Sylow \(p\)-subgroup of \(G\). If \(A\) has a complement in \(U\), then it has a complement in \(G\). The converse also holds. We use the formulation of Theorem~1.1 in \cite{GaschutzComplement}.
\end{classical}

\begin{classical}[Universal central extension]
A finite perfect group \(S\) has a universal central extension \(\widehat S\to S\) with kernel \(M(S)\), called the Schur multiplier. Every perfect central extension \(D\to S\) is a quotient of it; in particular, the kernel of \(D\to S\) is a quotient of \(M(S)\).
\end{classical}

\begin{classical}[Jordan's theorem]\label{fact:jordan}
If \(P\leqslant S_n\) is a primitive permutation group containing a cycle of prime length \(p\), where \(p\leqslant n-3\), then \(A_n\leqslant P\).
\end{classical}

Here a cycle means a single cyclic permutation of length \(p\), fixing the other \(n-p\) points. A transitive group is called primitive if it preserves no partition of the point set into equal blocks of size strictly between \(1\) and \(n\). This formulation appears, for example, as Theorem~1.1 in Jones's paper \cite{Jones}.

\begin{classical}[Simplicity and automorphisms of alternating groups]
For \(n\geqslant5\), the group \(A_n\) is nonabelian simple. For \(n\ge7\), we have \(\Aut(A_n)\cong S_n\). For the two groups used in the main argument,
\[
 \Aut(A_5)\cong S_5,\qquad
 |\Aut(A_5)|=120,\qquad |\Aut(A_6)|=1440.
\]
In particular, both of the latter orders are coprime to \(7\).
\end{classical}

For automorphism data one may use the ATLAS \cite{Atlas}; for \(A_6\), we have \(\Aut(A_6)\cong P\Gamma L_2(9)\), and \(|\operatorname{Out}(A_6)|=4\). We also use Cauchy's and Sylow's theorems in their usual formulations: if a prime \(p\) divides the order of a finite group, the group has an element of order \(p\); Sylow \(p\)-subgroups have order equal to the full \(p\)-part of the group order, and their images under an epimorphism are Sylow \(p\)-subgroups of the image. Every \(p\)-subgroup is contained in a Sylow \(p\)-subgroup, and all Sylow \(p\)-subgroups are conjugate.

\begin{lemma}\label{lem:centralizer-alt}
For \(n\geqslant4\), we have \(C_{S_n}(A_n)=1\). If \(n\geqslant5\) and \(N\trianglelefteq S_n\), then \(N=1,A_n\) or \(S_n\).
\end{lemma}
\begin{proof}
If \(s\) centralizes \(A_n\), it centralizes every \(3\)-cycle and therefore preserves its support. Hence \(s\) preserves every three-element subset of the point set. For each point \(i\), the intersection of all three-element subsets containing \(i\) is \(\{i\}\). Thus \(s\) fixes every point.

Now \(N\cap A_n\trianglelefteq A_n\), so simplicity of \(A_n\) implies that this intersection is \(1\) or \(A_n\). In the first case, \([N,A_n]\leqslant N\cap A_n=1\), whence \(N\leqslant C_{S_n}(A_n)=1\). In the second case, \(A_n\leqslant N\), and \([S_n:A_n]=2\) gives \(N=A_n\) or \(S_n\).
\end{proof}

\section{The spectra and criticality of the two groups}
\subsection{Spectra and the two critical groups}

For a permutation whose disjoint cycles have lengths \(\lambda_1,\ldots,\lambda_k\), including fixed points, the order is \(\lcm(\lambda_1,\ldots,\lambda_k)\), and the sign is \((-1)^{n-k}\). Thus the spectra of symmetric and alternating groups are calculated here by enumerating partitions of \(n\), rather than the permutations themselves.

\begin{proposition}\label{prop:spectra}
We have
\begin{align*}
 \Om:=\omega(S_{10})
 &=\{1,2,3,4,5,6,7,8,9,10,12,14,15,20,21,30\},\\
 \omega(S_9)
 &=\{1,2,3,4,5,6,7,8,9,10,12,14,15,20\},\\
 \omega(A_9)
 &=\{1,2,3,4,5,6,7,9,10,12,15\},\\
 \omega(A_{10})
 &=\{1,2,3,4,5,6,7,8,9,10,12,15,21\}.
\end{align*}
In particular, \(25,35\notin\Om\), \(\Om\setminus\omega(S_9)=\{21,30\}\), and \(\Om\setminus\omega(A_{10})=\{14,20,30\}\).
\end{proposition}
\begin{proof}
Let \(d=\prod p^{e_p}\). The minimum number of points on which a symmetric group contains an element of order \(d\) is
\[
 \ell(d)=\sum_{e_p>0}p^{e_p},\qquad \ell(1)=0.
\]
Indeed, disjoint cycles of lengths \(p^{e_p}\) give the upper bound. Conversely, for any collection of cycles, each maximal prime power \(p^{e_p}\) must divide some cycle length. If several coprime powers occur in one cycle, its length is at least their product, and that product is at least their sum. Summing over the cycles yields the lower bound \(\ell(d)\).

Now \(\ell(d)\leqslant10\) gives precisely the stated spectrum of \(S_{10}\), and \(\ell(d)\leqslant9\) gives the spectrum of \(S_9\). For odd \(d\), all cycles in the minimal construction have odd length, so the permutation is even. If \(d\) is even, the minimal construction contains exactly one cycle of even length. To obtain an even permutation, one must add another even cycle; the minimum additional length is \(2\). The same lower bound is necessary because an even permutation of even order has at least two cycles of even length. Consequently, for even \(d\), existence in \(A_n\) is equivalent to \(\ell(d)+2\leqslant n\). This gives both stated spectra of alternating groups. An independent verification by enumerating all partitions is given in the appendix.
\end{proof}

The elements of \(\Om\) maximal under divisibility, and their cycle realizations in \(S_{10}\), are:
\[
\begin{array}{c|rrrrrrr}
d&8&9&12&14&20&21&30\\ \hline
\text{cycle lengths}>1&8&9&4+3&7+2&5+4&7+3&5+3+2 .
\end{array}
\]

\begin{remark}
The group \(A_{10}\) contains an element of order \(8\): \((1\,2\,3\,4\,5\,6\,7\,8)(9\,10)\in A_{10}\) has order \(8\). This causes no difficulty in proving criticality of \(S_{10}\): \(A_{10}\) has no elements of orders \(14,20,30\).
\end{remark}

\begin{proposition}\label{prop:S10-critical}
The group \(S_{10}\) is \(\Om\)-critical.
\end{proposition}
\begin{proof}
Suppose that \(K\leqslant S_{10}\) and \(\omega(K)=\Om\). An element of order \(9\) on ten points must be a \(9\)-cycle with one fixed point. Hence either \(K\) is transitive or its orbits have sizes \(9\) and \(1\). The second case is impossible: an element of order \(21\) has cycle structure \(7+3\), and therefore has no fixed points.

We prove that \(K\) is primitive. In a nontrivial block system on ten points, the number of blocks and the block size are \(2,5\), in some order. An element of order \(7\) must act trivially on the blocks because there are fewer than seven of them, and then trivially inside each block for the same reason. This contradicts its order.

On ten points, an element of order \(7\) is a single \(7\)-cycle fixing three points. Jordan's theorem gives \(A_{10}\leqslant K\). An element of order \(14\) has cycle structure \(7+2+1\) and is odd. Hence \(K=S_{10}\).

It remains to check proper quotients. By Lemma~\ref{lem:centralizer-alt}, the nontrivial normal subgroups of \(S_{10}\) are \(A_{10}\) and \(S_{10}\); the corresponding quotients are \(C_2\) and \(1\). Their spectra differ from \(\Om\). Lemma~\ref{lem:sections} now applies.
\end{proof}

Define the second group in a way that removes any ambiguity in the semidirect-product notation:
\begin{equation}\label{eq:H}
 H=\{(s,t)\in S_9\times S_3:\sgn(s)=\sgn(t)\}.
\end{equation}
In its natural embedding in \(S_{12}\), the two factors act on disjoint sets of \(9\) and \(3\) points. The subgroup \(A_9\times A_3\cong A_9\times C_3\) is normal of index \(2\), and \((1\,2)(11\,12)\) defines a complement of order \(2\). It acts by an outer automorphism on \(A_9\) and by inversion on \(C_3\). Thus this is precisely the group \(H=(A_9\times C_3):2\) in the main theorem.

\begin{proposition}\label{prop:H-spectrum}
The group \(H\) is nonsolvable, is not isomorphic to \(S_{10}\), and satisfies \(\omega(H)=\Om\).
\end{proposition}
\begin{proof}
It contains the nonabelian simple subgroup \(A_9\), and is therefore nonsolvable. Its order is
\[
 |H|=\frac{9!\,3!}{2}=3\cdot9!=1\,088\,640
      =2^7\,3^5\,5\,7,
\]
which differs from \(|S_{10}|=3\,628\,800\).

If \(s\in A_9\), then \(t\in A_3\), and the pair has order either \(|s|\) or \(\lcm(|s|,3)\). By Proposition~\ref{prop:spectra}, these pairs give the set
\[
 \omega(A_9)\ \cup\
 \{\lcm(d,3):d\in\omega(A_9)\}
 =
 \{1,2,3,4,5,6,7,9,10,12,15,21,30\}.
\]
If \(s\) is odd, then \(t\) is a transposition, and \(|s|\) is even, so \(|(s,t)|=|s|\). The orders of odd permutations in \(S_9\) are
\[
 \{2,4,6,8,10,12,14,20\}.
\]
The latter set is obtained from the same partitions of \(9\), subject to the oddness condition. The union of the two displayed sets is exactly \(\Om\).
\end{proof}

\begin{proposition}\label{prop:H-critical}
The group \(H\) is \(\Om\)-critical.
\end{proposition}
\begin{proof}
Let \(C=\{1\}\times A_3\cong C_3\) be the kernel of the projection \(\pi:H\to S_9\). Suppose that \(K\leqslant H\) and \(\omega(K)=\Om\). Put \(P=\pi(K)\).

In the model~\eqref{eq:H}, the order of a pair is the least common multiple of the component orders, and the second component has order \(1,2\) or \(3\). Therefore an element of order \(9\) in \(H\) has first component of order \(9\); similarly, for orders \(7\) and \(20\), the first component has the same order. Hence \(9,7,20\in\omega(P)\).

The presence of a \(9\)-cycle makes \(P\) transitive on nine points. The only possible nontrivial block system would consist of three blocks of three points. An element of order \(7\) would act trivially both on the blocks and within them. Thus \(P\) is primitive. A permutation of order \(20\) on nine points has structure \(5+4\); its fourth power is a \(5\)-cycle. Jordan's theorem with \(p=5\leqslant9-3\) gives \(A_9\leqslant P\). The permutation of structure \(5+4\) itself is odd, so \(P=S_9\).

Now \(K\cap C\) is \(1\) or \(C\). If \(K\cap C=1\), then \(K\cong S_9\), contrary to \(21\in\omega(K)\) and \(21\notin\omega(S_9)\). If \(K\cap C=C\), then \(|K|=|C|\,|S_9|=|H|\), so \(K=H\). Proper subgroups have been excluded.

Consider a nontrivial normal subgroup \(N\trianglelefteq H\). Write \(A=A_9\times\{1\}\trianglelefteq H\). By simplicity of \(A_9\), we have \(N\cap A=A\) or \(1\). In the first case, \(A\leqslant N\), and \(H/N\) is a quotient of \(H/A\cong S_3\); its spectrum does not contain \(5\).

In the second case, \([N,A]\leqslant N\cap A=1\). By Lemma~\ref{lem:centralizer-alt} and the model~\eqref{eq:H}, the centralizer of \(A\) in \(H\) is \(C\). Hence \(1\ne N\leqslant C\), so \(N=C\), and \(H/N\cong S_9\) has no element of order \(21\). Thus proper quotients are also excluded. Lemma~\ref{lem:sections} completes the proof.
\end{proof}

\section{The quotient by the radical is almost simple}
\subsection{The action of \(C_5\times C_5\) and the solvable radical}

\begin{lemma}\label{lem:fixed}
Suppose that \(E\cong C_5\times C_5\) acts linearly on a nonzero finite-dimensional space \(V\) over a field of characteristic \(7\). Then there exists \(1\ne e\in E\) with \(C_V(e)\ne0\).
\end{lemma}
\begin{proof}
Extend the field to an algebraic closure \(k\). For every \(e\in E\), the minimal polynomial of the induced operator divides \(X^5-1\), which has no repeated roots in characteristic \(7\): its derivative is \(5X^4\). Consequently, all operators from \(E\) are diagonalizable over \(k\).

Since they commute, there is a common eigenvector \(0\ne v\in V\otimes k\). Indeed, the eigenspaces of one operator are invariant under all the others, so an eigenspace can be chosen successively inside the one already chosen. We have \(v^e=\chi(e)v\), where \(\chi:E\to k^*\) is a homomorphism. Its image lies in the group of fifth roots of unity, of order \(5\). Thus \(|\ker\chi|\geqslant5\). Choose \(1\ne e\in\ker\chi\). Then \(e-1\) has nonzero kernel over \(k\).

Matrix rank is unchanged by a field extension (it is determined by nonzero minors). Hence the kernel of \(e-1\) is nonzero over the original field as well. Thus \(C_V(e)\ne0\).
\end{proof}

\begin{lemma}\label{lem:five}
If \(P\) is a finite \(5\)-group, \(|P|\geqslant25\), and \(P\) has no element of order \(25\), then \(P\) contains \(C_5\times C_5\).
\end{lemma}
\begin{proof}
Every nontrivial element of \(P\) has order \(5\): an element of larger order would have a power of order \(25\). A nontrivial finite \(p\)-group has nontrivial center: this follows from the class equation because noncentral conjugacy classes have sizes divisible by \(p\). Choose \(z\in Z(P)\) of order \(5\) and \(x\in P\setminus\langle z\rangle\). Then \(z\) and \(x\) commute, both have order \(5\), and their cyclic subgroups intersect trivially. Therefore \(\langle z,x\rangle\cong C_5\times C_5\).
\end{proof}

\begin{lemma}\label{lem:chief}
Let \(R\trianglelefteq G\) be solvable, and suppose that \(7\mid |R|\). Then there exist subgroups \(K<N\leqslant R\), normal in \(G\), such that \(N/K\) is a nontrivial elementary abelian \(7\)-group.
\end{lemma}
\begin{proof}
Choose a maximal chain of subgroups normal in \(G\) between \(1\) and \(R\). Each factor is a solvable minimal normal subgroup of the corresponding quotient of \(G\). Let \(M\) be such a factor. The subgroup \(M'\) is characteristic in \(M\), and hence normal in the corresponding quotient of \(G\). Minimality and solvability of \(M\) imply \(M'=1\): the equality \(M'=M\ne1\) is incompatible with solvability.

Thus \(M\) is abelian. Its Sylow subgroups are characteristic, so minimality implies that \(M\) is a \(p\)-group for a single prime \(p\). The subgroup of elements whose orders divide \(p\) is also characteristic and nontrivial; again by minimality it is \(M\). Therefore \(M\) is elementary abelian. Since the product of the factor orders equals \(|R|\), at least one factor is a \(7\)-group.
\end{proof}

\begin{lemma}\label{lem:radical}
Let \(R\trianglelefteq G\) be solvable, and suppose that
\[
 25,35\notin\omega(G),\qquad 25\mid |G/R|.
\]
Then \(7\nmid |R|\).
\end{lemma}
\begin{proof}
Suppose that \(7\mid |R|\). By Lemma~\ref{lem:chief}, choose \(K<N\leqslant R\), normal in \(G\), such that \(V=N/K\) is a nonzero elementary abelian \(7\)-group. It is normal in \(\overline G=G/K\).

Take a Sylow \(5\)-subgroup \(P\) of \(\overline G\). Its image under the epimorphism \(\overline G\to G/R\) is a Sylow \(5\)-subgroup of \(G/R\), so \(|P|\geqslant25\). By Lemma~\ref{lem:sections}, \(25,35\notin\omega(\overline G)\). Hence Lemma~\ref{lem:five} yields \(E\leqslant P\), \(E\cong C_5\times C_5\).

The conjugation action of \(E\) on \(V\) is linear over \(\F_7\). By Lemma~\ref{lem:fixed}, there exist \(1\ne e\in E\) and \(1\ne v\in V\) with \(v^e=v\). The orders of \(e\) and \(v\) are \(5\) and \(7\). The product of commuting elements of coprime orders has order equal to the product of those orders. Therefore \(|ev|=35\), contrary to \(35\notin\omega(\overline G)\).
\end{proof}

\begin{remark}
Lemma~\ref{lem:radical} assumes neither that \(G\) splits over \(R\), nor that \(\overline G\) splits over \(V\). The subgroup \(E\) is chosen in \(\overline G\) itself; it is not lifted from a quotient as a prescribed complement. The hypothesis \(25\notin\omega(G)\) is essential to the method used to choose \(E\).
\end{remark}

\subsection{Permuting the simple direct factors}

\begin{lemma}\label{lem:faithful}
Suppose that \(T=\Soc(X)\) is a direct product of nonabelian simple groups. Then conjugation gives
\[
 C_X(T)=1,\qquad X\hookrightarrow\Aut(T)
\]
\end{lemma}
\begin{proof}
Since \(T\) is a direct product of nonabelian simple groups, \(Z(T)=1\). The centralizer \(D=C_X(T)\) is normal in \(X\). If \(D\ne1\), choose a minimal nontrivial normal subgroup \(M\) of \(X\) contained in \(D\). By the definition of the socle, \(M\leqslant T\). Then \(M\leqslant T\cap C_X(T)=Z(T)=1\), a contradiction. Thus the kernel of the conjugation action is trivial.
\end{proof}

\begin{lemma}\label{lem:aut}
Let \(T=T_1\times\cdots\times T_m\), with every \(T_i\) isomorphic to \(A_5\) or \(A_6\). If \(\alpha\in\Aut(T)\) has order \(7\), then the permutation it induces on the simple direct factors has a cycle of length \(7\).
\end{lemma}
\begin{proof}
First, we explain why automorphisms permute the factors \(T_i\). They are precisely the minimal normal subgroups of \(T\). Indeed, if \(M\trianglelefteq T\) is minimal and the projection of \(M\) onto \(T_i\) is nontrivial, then \([M,T_i]\leqslant M\cap T_i\). This subgroup is nontrivial: otherwise the nontrivial projection of \(M\) would centralize \(T_i\), whereas that projection is normal in the simple group \(T_i\), hence equals \(T_i\), and \(Z(T_i)=1\). Consequently, \(M\cap T_i\ne1\), whence \(T_i\leqslant M\), and minimality gives \(M=T_i\). Automorphisms preserve minimal normal subgroups and thus permute the \(T_i\).

The kernel of the action of \(\Aut(T)\) on the factors equals \(\prod_i\Aut(T_i)\): if each factor is preserved, the restrictions to the individual factors independently define an automorphism of the whole direct product. The order of this kernel is not divisible by \(7\), since
\[
 |\Aut(A_5)|=120,\qquad |\Aut(A_6)|=1440.
\]
Thus \(\alpha\) does not lie in the kernel and induces a nontrivial permutation of order \(7\). Its cycles have lengths \(1\) or \(7\); at least one has length \(7\).
\end{proof}

\begin{lemma}\label{lem:diagonal}
Let \(T=T_1\times\cdots\times T_m\trianglelefteq X\), and let \(x\in X\) have order \(7\) and cyclically permute seven factors \(T_1,\ldots,T_7\). If \(T_1\) contains an element of order \(5\), then \(35\in\omega(X)\).
\end{lemma}
\begin{proof}
Let \(a\in T_1\), \(|a|=5\). With a suitable numbering, \(a^{x^j}\in T_{j+1}\) for \(j=0,\ldots,6\). Put
\[
 d=a\,a^x\,a^{x^2}\cdots a^{x^6}.
\]
All seven factors lie in different direct factors and commute. Their projections have order \(5\), so \(|d|=5\). Since \(x^7=1\),
\[
 d^x=a^x\,a^{x^2}\cdots a^{x^6}a^{x^7}
     =a^x\,a^{x^2}\cdots a^{x^6}a=d.
\]
Thus \(d\) and \(x\) commute and \(|dx|=5\cdot7=35\).
\end{proof}

\subsection{Excluding composite socles}

\begin{theorem}\label{thm:main}
Let \(\omega(G)=\Om=\omega(S_{10})\), \(R=R(G)\), and \(X=G/R\). Then the following isomorphism is impossible:
\[
 \Soc(X)\cong A_5^a\times A_6^b,
 \qquad a,b\geqslant0,\quad a+b\geqslant2.
\]
 This assertion does not require \(G\) to be \(\Om\)-critical.
\end{theorem}
\begin{proof}
Suppose otherwise and write \(T=\Soc(X)=T_1\times\cdots\times T_m\), where \(m=a+b\geqslant2\). Every \(T_i\) contains an element of order \(5\); in particular, \(25\mid |T|\) and \(25\mid |X|\). By Proposition~\ref{prop:spectra}, \(25,35\notin\omega(G)\). Lemma~\ref{lem:radical} gives \(7\nmid |R|\).

Since \(7\in\omega(G)\), the prime \(7\) divides \(|G|\), and hence also \(|X|=|G|/|R|\). By Cauchy's theorem, there exists \(x\in X\) of order \(7\). By Lemma~\ref{lem:faithful}, conjugation by \(x\) induces an automorphism of \(T\) of the same order \(7\). Lemma~\ref{lem:aut} shows that it cyclically permutes seven simple factors. Lemma~\ref{lem:diagonal} gives \(35\in\omega(X)\). However, Lemma~\ref{lem:sections} gives
\[
 \omega(X)\subseteq\omega(G)=\Om,\qquad35\notin\Om.
\]
This contradiction completes the proof.
\end{proof}

\begin{corollary}\label{cor:cases}
If \(\omega(G)=\omega(S_{10})\), then all of the following possibilities for \(\Soc(G/R(G))\) are excluded:
\[
 A_5^m\ (m\geqslant2),\qquad A_5\times A_6,\qquad
 A_6\times A_6,
\]
as are all mixed products of three or more copies of \(A_5\) and \(A_6\).
\end{corollary}
\begin{proof}
Every case listed is covered directly by Theorem~\ref{thm:main}.
\end{proof}

\begin{remark}
For \(2\leqslant m<7\), the final step of the proof is even shorter: there is no permutation of order \(7\) on \(m\) factors, so Lemma~\ref{lem:aut} already gives a contradiction. For \(m\geqslant7\), one uses the diagonal element from Lemma~\ref{lem:diagonal}. Thus no separate bound on \(m\) is needed.
\end{remark}

\begin{proposition}\label{prop:products}
For direct products of copies of \(A_5\) and \(A_6\), the condition \(\omega(T)\subseteq\Om\) holds:
\begin{enumerate}
\item for \(A_5^m\) with any \(m\geqslant1\);
\item for \(A_6\), \(A_5\times A_6\), and \(A_6\times A_6\).
\end{enumerate}
If there are at least three factors, one of which is \(A_6\), the condition does not hold.
\end{proposition}
\begin{proof}
The cycle structures give
\[
 \omega(A_5)=\{1,2,3,5\},\qquad
 \omega(A_6)=\{1,2,3,4,5\}.
\]
The order of an element of a direct product is the least common multiple of the component orders. For \(A_5^m\), it always divides \(30\), and every divisor of \(30\) lies in \(\Om\). For two factors chosen from \(A_5,A_6\), all possible least common multiples lie in \(\{1,2,3,4,5,6,10,12,15,20\}\subseteq\Om\). If there are at least three factors and one is isomorphic to \(A_6\), take an element of order \(4\) in that factor and elements of orders \(3\) and \(5\) in two other factors. Their product has order \(\lcm(4,3,5)=60\notin\Om\).
\end{proof}

Thus the spectrum condition on the direct product itself excludes neither \(A_5^m\) nor the pairs listed above. Theorem~\ref{thm:main} uses the further requirement that the whole group contain an element of order \(7\).

\subsubsection{The socle of the quotient by the radical}

\begin{lemma}\label{lem:radical-free-socle-structure}
Let \(G\) be a nonsolvable finite group, \(R=R(G)\), and \(X=G/R\). Then \(X\ne1\), \(R(X)=1\), and its socle has the form
\[
 T=\operatorname{Soc}(X)=T_1\times\cdots\times T_m,
 \qquad m\ge1,
\]
where all \(T_i\) are nonabelian simple groups. If \(\omega(G)=\Om\), then \(\omega(T)\subseteq\Om\).
\end{lemma}
\begin{proof}
If \(X=1\), then \(G=R\) is solvable, contrary to the choice of \(G\). Suppose that \(1\ne D\trianglelefteq X\) is solvable, and let \(\widetilde D\) be its full preimage in \(G\). Then
\[
 R<\widetilde D\trianglelefteq G,\qquad
 \widetilde D/R\cong D.
\]
An extension of a solvable group by a solvable group is solvable. Indeed, if \(D^{(a)}=1\) and \(R^{(b)}=1\), then \(\widetilde D^{(a)}\leqslant R\), whence \(\widetilde D^{(a+b)}=1\). Thus \(\widetilde D\) would be a solvable normal subgroup of \(G\) strictly containing \(R(G)\), which is impossible. Therefore \(R(X)=1\).

Let \(M\) be a minimal nontrivial normal subgroup of \(X\). Such a subgroup exists because \(X\) is finite and nontrivial. It has no nontrivial proper characteristic subgroup: any such subgroup would be normal in \(X\), contradicting minimality of \(M\).

Choose a minimal nontrivial normal subgroup \(L\) of the group \(M\) itself, and list its distinct images under all automorphisms of \(M\):
\[
 L_1=L,L_2,\ldots,L_k.
\]
All the \(L_i\) are isomorphic and minimal normal in \(M\). For \(i\ne j\), their intersection is normal in \(M\), and minimality gives \(L_i\cap L_j=1\). Moreover,
\[
 [L_i,L_j]\leqslant L_i\cap L_j=1,
\]
so distinct \(L_i\) commute. Their product is characteristic in \(M\), since automorphisms of \(M\) permute all the \(L_i\). Characteristic simplicity of \(M\) therefore gives \(M=L_1\cdots L_k\).

If \(L\) were abelian, all the \(L_i\) would be abelian; since they commute pairwise, \(M\) would be abelian. But then \(M\) would be a nontrivial solvable normal subgroup of \(X\), contrary to \(R(X)=1\). Thus all the \(L_i\) are nonabelian. Their centers are characteristic in \(L_i\) and normal in \(M\). Minimality and nonabelianness of \(L_i\) imply \(Z(L_i)=1\).

For every \(i\), the product of the other \(L_j\) centralizes \(L_i\), so
\[
 L_i\cap\prod_{j\ne i}L_j\leqslant Z(L_i)=1.
\]
Consequently,
\[
 M=L_1\times\cdots\times L_k.
\]
Every normal subgroup \(U\) of one factor \(L_i\) is normal in all of \(M\): it is normalized by \(L_i\), and the other direct factors centralize it. Minimality of \(L_i\) as a normal subgroup of \(M\) gives \(U=1\) or \(U=L_i\). Thus every \(L_i\) is nonabelian simple.

List all minimal normal subgroups of \(X\) as \(M_1,\ldots,M_s\). By definition, \(T=M_1\cdots M_s\). Distinct \(M_i\) intersect trivially and commute, since \([M_i,M_j]\leqslant M_i\cap M_j\). By what we have proved, each \(M_i\) is centerless. Hence
\[
 M_i\cap\prod_{j\ne i}M_j\leqslant Z(M_i)=1,
\]
and the product \(M_1\cdots M_s\) is direct. Decomposing each \(M_i\) as above, we obtain
\[
 T=T_1\times\cdots\times T_m,\qquad m\geqslant1,
\]
where all \(T_i\) are nonabelian simple. Here \(m\) counts the simple direct factors, not necessarily the minimal normal subgroups of \(X\): a single \(M_i\) may contain several simple direct factors.

By Lemma~\ref{lem:sections},
\[
 \omega(T)\subseteq\omega(X)\subseteq\omega(G)=\Om.
\]
Indeed, \(T\) is a subgroup of the quotient \(X\), and the spectrum of a quotient is contained in that of the original finite group. The equality \(\omega(X)=\Om\) is neither required nor asserted here.
\end{proof}

\subsection{The simple direct factors of the socle}

We use the following classification theorem of Zavarnitsine \cite[Table 1, rows \(S_5,S_7\)]{Zav}.

\begin{classical}[Simple groups with prime divisors at most \(7\)]
Every finite nonabelian simple group \(S\) such that \(\pi(S)\subseteq\{2,3,5,7\}\) is isomorphic to one of the following 18 groups:
\begin{align*}
&A_5,\ A_6,\ \operatorname{PSp}_4(3),\
 L_2(7),\ L_2(8),\ U_3(3),\\
&A_7,\ A_8,\ A_9,\ A_{10},\
 L_2(49),\ L_3(4),\
 U_3(5),\ U_4(3),\\
&J_2,\ \operatorname{PSp}_4(7),\
 \operatorname{Sp}_6(2),\ \mathrm P\Omega_8^+(2).
\end{align*}
Here \(L_n(q)=\operatorname{PSL}_n(q)\), \(U_n(q)=\operatorname{PSU}_n(q)\). The list is up to isomorphism; \(\operatorname{PSp}_4(3)\cong U_4(2)\).
\end{classical}

To exclude a composite socle, the full spectra of the listed groups are not needed. It suffices to divide the list into three classes:
\[
\begin{array}{c|l}
\text{class}&\text{groups}\\ \hline
\mathcal A&A_5,\ A_6,\ \operatorname{PSp}_4(3)\\
\mathcal B&L_2(7),\ L_2(8),\ U_3(3)\\
\mathcal C&A_7,A_8,A_9,A_{10},L_2(49),L_3(4),\\
&U_3(5),U_4(3),J_2,\operatorname{PSp}_4(7),
  \operatorname{Sp}_6(2),\mathrm P\Omega_8^+(2).
\end{array}
\]
Groups in \(\mathcal A\) have elements of order \(5\), those in \(\mathcal B\) have elements of order \(7\), and those in \(\mathcal C\) have elements of both orders. This follows from the factorizations of the orders in the cited table and Cauchy's theorem.

We also need elements of order \(4\) in \(L_2(7)\) and \(U_3(3)\), and of order \(9\) in \(L_2(8)\) and \(\operatorname{PSp}_4(3)\). We justify these facts.

The group \(L_2(q)\) has a cyclic subgroup of order \((q+1)/(2,q-1)\). Indeed, multiplication by the elements of norm \(1\) in \(\mathbb F_{q^2}^*\) gives a cyclic subgroup of order \(q+1\) in \(\operatorname{SL}_2(q)\); on passing to the projective group, its order is divided by \((2,q-1)\). For \(q=7\) and \(q=8\), this gives \(4\) and \(9\), respectively.

In \(U_3(3)\), choose an element \(i\) of order \(4\) in \(\mathbb F_9\) and take the matrix \(\operatorname{diag}(i,i^{-1},1)\). It preserves the standard Hermitian form, since \(i\,i^3=i^4=1\), and has determinant \(1\). The center of \(\operatorname{SU}_3(3)\) has order \((3,4)=1\), so the order \(4\) is preserved in \(U_3(3)\).

For \(\operatorname{PSp}_4(3)\), explicit matrices suffice. Over \(\mathbb F_3\), put
\[
 J=\begin{pmatrix}
 1&1&0&0\\0&1&1&0\\0&0&1&1\\0&0&0&1
 \end{pmatrix},\qquad
 B=\begin{pmatrix}
 0&0&0&1\\0&0&2&1\\0&1&0&0\\2&2&0&0
 \end{pmatrix}.
\]
Direct multiplication gives \(J^{\mathsf T}BJ=B\). The matrix \(B\) is skew-symmetric and has determinant \(1\), so it defines a nondegenerate alternating form. If \(N=J-I\), then \(N^4=0\), \(N^3\ne0\); in characteristic \(3\),
\[
 J^3=I+N^3\ne I,\qquad J^9=I.
\]
Thus \(J\) has order \(9\) in the symplectic group of this form. Factoring out the center of order \(2\) does not change the order \(9\).

\begin{theorem}\label{thm:first-premise}
Let \(\omega(G)=\Om\), \(R=R(G)\). If
\[
 T=\Soc(G/R)=T_1\times\cdots\times T_m,\qquad m\geqslant2,
\]
where all \(T_i\) are nonabelian simple groups, then each \(T_i\) is isomorphic to \(A_5\) or \(A_6\). This assertion does not require criticality of \(G\).
\end{theorem}
\begin{proof}
The inclusions \(\omega(T_i)\subseteq\omega(T)\subseteq\Om\) and Cauchy's theorem imply \(\pi(T_i)\subseteq\{2,3,5,7\}\). Thus every \(T_i\) belongs to the displayed list of 18 simple groups.

Elements in distinct direct factors commute. Consequently, if one factor contains an element of order \(r\), and another contains one of order \(s\), their product contains an element of order \(\lcm(r,s)\). We use
\[
 28,\ 35,\ 45,\ 63\notin\Om.
\]

If one factor belongs to \(\mathcal C\), it has elements of orders \(5\) and \(7\). Every other factor belongs to \(\mathcal A\), \(\mathcal B\), or \(\mathcal C\), so it contains an element of order \(5\) or \(7\). Choosing the other order in the first factor gives \(35\), a contradiction. Hence \(\mathcal C\) is excluded entirely.

If one factor belongs to \(\mathcal A\) and another to \(\mathcal B\), we again obtain \(5\cdot7=35\). Thus all factors belong to the same one of \(\mathcal A\), \(\mathcal B\).

Suppose that all belong to \(\mathcal B\). Every group in this class contains elements of order \(7\) and also of order \(4\) or \(9\), as shown above. Any two factors therefore give either \(7\cdot4=28\) or \(7\cdot9=63\). This is impossible.

Thus all factors belong to \(\mathcal A\). If one is \(\operatorname{PSp}_4(3)\), take an element of order \(9\) in it and an element of order \(5\) in another factor. This gives \(45\), the final contradiction. Only \(A_5\) and \(A_6\) remain.
\end{proof}

\begin{corollary}\label{cor:simple-socle}
For every nonsolvable finite group \(G\) with spectrum \(\Om\), the group \(X=G/R(G)\) has a nonabelian simple socle \(S\), and
\[
 S\leqslant X\leqslant\Aut(S).
\]
In particular, no such group has a composite socle.
\end{corollary}
\begin{proof}
By Lemma~\ref{lem:radical-free-socle-structure}, we have \(R(X)=1\), \(X\ne1\), and the socle of \(X\) is a nontrivial direct product of nonabelian simple groups. If there are at least two direct factors, Theorem~\ref{thm:first-premise} shows that it is a product of copies of \(A_5\) and \(A_6\). Theorem~\ref{thm:main} excludes this case. Thus the socle is simple. The embedding in \(\Aut(S)\) follows from Lemma~\ref{lem:faithful}.
\end{proof}

\begin{lemma}\label{lem:standalone-simple-candidates}
If \(G\) is nonsolvable and \(\omega(G)=\Om\), then the simple socle of \(G/R(G)\) belongs to the list
\[
\begin{gathered}
 A_5,A_6,A_7,A_8,A_9,A_{10},L_2(7),L_2(8),L_3(4),\\
 U_3(3),U_3(5),U_4(2),U_4(3),J_2,S_6(2),O_8^+(2).
\end{gathered}
\]
\end{lemma}
\begin{proof}
By Corollary~\ref{cor:simple-socle}, the socle \(S\) is simple, and \(\omega(S)\subseteq\Om\). Cauchy's theorem gives \(\pi(S)\subseteq\{2,3,5,7\}\), so the preceding list of eighteen simple groups applies. In that list, \(\operatorname{PSp}_4(3)\cong U_4(2)\). We exclude two further groups because they contain elements of order \(25\).

In \(L_2(49)\), the cyclic nonsplit torus has order \((49+1)/2=25\). Its construction using norm-\(1\) elements of a quadratic field extension was given above in the discussion of \(L_2(q)\).

For \(\operatorname{PSp}_4(7)\), put \(q=7\) and regard \(V=\mathbb F_{q^4}\) as a four-dimensional space over \(\mathbb F_q\). Choose \(0\ne c\in\mathbb F_{q^4}\) with \(c^{q^2}=-c\). The form
\[
 B(x,y)=\operatorname{Tr}_{\mathbb F_{q^4}/\mathbb F_q}
             (cxy^{q^2})
\]
is nondegenerate and skew-symmetric. Indeed, the trace is invariant under the \(q^2\)-power map, so \(B(y,x)=-B(x,y)\); nondegeneracy follows from that of the trace pairing and invertibility of \(y\mapsto y^{q^2}\). Multiplication by elements \(\alpha\) of norm \(1\) from \(\mathbb F_{q^4}\) to \(\mathbb F_{q^2}\) preserves \(B\), since \(\alpha^{q^2+1}=1\). These multiplications form a cyclic subgroup of order \(q^2+1=50\) in \(\operatorname{Sp}_4(7)\). After factoring out the scalars \(\{I,-I\}\), we obtain a cyclic subgroup of order \(25\) in \(\operatorname{PSp}_4(7)\). The spectrum \(\Om\) forbids order \(25\), excluding both groups and giving the stated list.
\end{proof}

\section{\texorpdfstring{The socles \(A_9,A_{10}\)}{The socles A9,A10}}
\subsection{Almost simple quotients and socles of large degree}

We use the published classification of Zhurtov--Lytkina--Mazurov \cite{ZhLM}, stating explicitly the consequences of its tables needed below.

\subsubsection{Trivial radical}

\begin{proposition}\label{prop:radical-zero}
If \(\omega(G)=\Om\), \(G\) is nonsolvable, and \(R(G)=1\), then \(G\cong S_{10}\).
\end{proposition}
\begin{proof}
By Corollary~\ref{cor:simple-socle}, \(G\) is almost simple. Lemma 1 and Table 2 in \cite{ZhLM} give the complete list of almost simple groups subspectral to \(S_{10}\). Among these groups, order \(21\) occurs only in
\[
 L_3(4).3,\quad L_3(4).6,\quad A_{10},\quad S_{10}.
\]
For the first two, the corresponding sets of element orders maximal under divisibility are \(\{4,6,15,21\}\) and \(\{8,12,15,21\}\); \(A_{10}\) has no element of order \(20\). Thus the simultaneous occurrence of \(20,21\in\omega(G)\) implies \(G\cong S_{10}\).
\end{proof}

\subsubsection{A criterion for the absence of a radical above a simple quotient}

\begin{lemma}\label{lem:no-covers}
Let \(R\trianglelefteq P\) be solvable, \(P/R\cong S\), where \(S\) is nonabelian simple, and \(\omega(P)\subseteq\Om\). Suppose that there is no group \(E\) subspectral to \(S_{10}\) with a nontrivial normal elementary abelian subgroup \(V\) such that \(E/V\cong S\) and the action of \(E\) on \(V\) is irreducible. Then \(R=1\).
\end{lemma}
\begin{proof}
If \(R\ne1\), choose \(K\) maximal among the proper subgroups of \(R\) that are normal in \(P\). Then \(V=R/K\) is minimal normal in \(E=P/K\). Since \(V\) is solvable, the proof of Lemma~\ref{lem:chief} shows that \(V\) is elementary abelian. Its minimal normality is equivalent to irreducibility of the conjugation action. Moreover,
\[
 E/V\cong P/R\cong S,\qquad
 \omega(E)\subseteq\omega(P)\subseteq\Om.
\]
This gives a cover \(E\) forbidden by the hypothesis. Hence \(R=1\). No splitting assumption on \(P\) is made in this proof.
\end{proof}

\begin{theorem}\label{thm:A10}
If \(\omega(G)=\Om\) and \(\Soc(G/R(G))\cong A_{10}\), then \(G\cong S_{10}\).
\end{theorem}
\begin{proof}
Put \(R=R(G)\), \(X=G/R\). By Lemma~\ref{lem:faithful},
\[
 A_{10}\leqslant X\leqslant\Aut(A_{10})\cong S_{10},
\]
so \(X=A_{10}\) or \(S_{10}\). We use the standard theorem on automorphisms of alternating groups: \(\Aut(A_n)\cong S_n\) for \(n\geqslant7\).

Let \(P\) be the full preimage of \(A_{10}\) in \(G\). Then \(P/R\cong A_{10}\). By Theorem 1 and the \(A_{10}\) row of Table 1 in \cite{ZhLM}, no nontrivial irreducible elementary abelian extension of \(A_{10}\) is subspectral to \(S_{10}\). Lemma~\ref{lem:no-covers} gives \(R=1\). Now \(G=A_{10}\) or \(S_{10}\). The first case is excluded by order \(20\), so \(G\cong S_{10}\).
\end{proof}

\subsubsection{A single minimal normal kernel}

\begin{proposition}\label{prop:minimal-alt}
Let \(\omega(G)=\Om\), and suppose that \(R=R(G)\ne1\) is minimal normal in \(G\). Then \(G/R\) is not isomorphic to any simple alternating group.
\end{proposition}
\begin{proof}
If \(G/R\cong A_n\), its prime divisors imply \(5\leqslant n\leqslant10\). As in Lemma~\ref{lem:chief}, \(R\) is elementary abelian and the action of \(G\) on \(R\) is irreducible. Thus Theorem 1 in \cite{ZhLM} applies. Its Table 1 gives the following obstructions:
\[
\begin{array}{c|l}
G/R&\text{order absent from an admissible cover}\\ \hline
A_5&7\text{ if the characteristic is not }7;\
      8\text{ for }7^4.A_5\\
A_6&7\\
A_7&9\\
A_8&9\\
A_9&20\\
A_{10}&\text{no nontrivial admissible cover}.
\end{array}
\]
In each row, either the required cover does not exist, or its spectrum omits an element of \(\Om\). This is a contradiction.
\end{proof}

\subsubsection{A single top \(3\)-kernel above \(S_9\)}

\begin{lemma}\label{lem:split-spectrum}
Let \(V\trianglelefteq E\) be elementary abelian of characteristic \(p\), and put \(Q=E/V\). For the action of \(Q\) on \(V\) induced by conjugation,
\[
 \omega(V\rtimes Q)\subseteq\omega(E).
\]
No splitting of \(E\) is assumed.
\end{lemma}
\begin{proof}
Let \(q\in Q\) have order \(m\). In the split extension, every element of the coset \(Vq\) has order \(m\) or \(pm\). Order \(m\) already lies in \(\omega(E)\) by Lemma~\ref{lem:sections}. If \(Vq\) contains an element of order \(pm\), the operator \(N_q=1+q+\cdots+q^{m-1}\) on \(V\) is nonzero: in additive notation, \((v,q)^m=N_qv\).

Choose any preimage \(x\in E\) of \(q\). Since \(x^m\in V\), the order of \(x\) is \(m\) or \(pm\). The second case proves the assertion. In the first, \(x^m=1\), and choosing \(v\) with \(N_qv\ne0\) gives an element \(vx\) of order \(pm\) in \(E\). This proves the inclusion for every coset.
\end{proof}

\begin{proposition}\label{prop:S9-three}
Let \(G\) be \(\Om\)-critical, \(R=R(G)\), and \(G/R\cong S_9\). Suppose that there exists \(K\trianglelefteq G\), \(K<R\), such that \(R/K\) is a minimal normal \(3\)-subgroup of \(G/K\). Then \(K=1\) and \(G\cong H\).
\end{proposition}
\begin{proof}
Put \(E=G/K\), \(V=R/K\). The group \(V\) is elementary abelian over \(\mathbb F_3\). Its conjugation action factors through \(E/V\cong S_9\), since \(V\) is abelian. Minimal normality means that \(V\) is irreducible as an \(S_9\)-module. Lemma~\ref{lem:split-spectrum} gives
\[
 \omega(V\rtimes S_9)\subseteq\omega(E)\subseteq\Om.
\]

The restriction of \(V\) to \(A_9\trianglelefteq S_9\) is completely reducible. To see this, choose a minimal nonzero \(A_9\)-submodule \(U\). Its conjugates under \(S_9\) are also simple as \(A_9\)-modules; their sum is \(S_9\)-invariant, and hence equals \(V\) by irreducibility. A sum of simple submodules is completely reducible. Therefore each simple \(A_9\)-constituent \(U\) is also a quotient of \(V\) as an \(A_9\)-module, and \(U\rtimes A_9\) is a section of \(V\rtimes S_9\).

By Theorem 1 and Table 1 in \cite{ZhLM}, the only admissible irreducible \(3\)-kernel above \(A_9\) is the trivial one-dimensional module, giving \(C_3\times A_9\). Thus \(A_9\) acts trivially on all simple summands of \(V\), and hence on \(V\). The action of \(S_9\) factors through \(C_2\). An irreducible \(C_2\)-module over \(\mathbb F_3\) is one-dimensional: the involution is diagonalizable with eigenvalues \(1,-1\). Consequently, \(V\cong C_3\).

Let \(P\) be the full preimage of \(A_9\) in \(E\). The same row of the cover classification gives \(P\cong C_3\times A_9\). The subgroup \(A=P'\cong A_9\) is characteristic in \(P\), and hence normal in \(E\).

If \(V\leqslant Z(E)\), choose a preimage \(x\) of an element of order \(8\) in \(S_9\). Its order is \(8\) or \(24\). The latter is forbidden. In the former case, \(x\) commutes with a nontrivial \(v\in V\), and \(|xv|=24\), again forbidden. Thus odd permutations invert \(V\).

Choose a preimage \(y\) of a transposition in \(S_9\). We have \(y^2\in V\). The element \(y\) centralizes its square and also inverts every element of \(V\). Hence \(y^2=(y^2)^{-1}\); since \(|V|=3\), this gives \(y^2=1\). Conjugation by \(y\) on \(A\) corresponds to conjugation by a transposition on \(A_9\): the isomorphism \(A\to P/V\) is induced by projection. Therefore
\[
 E=(A\times V)\rtimes\langle y\rangle\cong
 (A_9\times C_3):2=H
\]
with the action defined in \eqref{eq:H}. In particular, \(\omega(E)=\Om\). Criticality of \(G\) excludes a proper isospectral quotient \(G/K\), so \(K=1\) and \(G\cong H\).
\end{proof}

\subsubsection{Excluding a top \(2\)-kernel above \(S_9\)}

For this case, the three simple \(8\)-dimensional \(\mathbb F_2 A_9\)-modules must be distinguished. The completeness of this list follows from the decomposition matrix \cite{A9mod2}. Matrices for the standard generators are taken from the ATLAS \cite{A9Matrices}. The required matrices and the results of Gaussian elimination are given below, allowing the numerical part of the proof to be checked directly.

\begin{lemma}\label{lem:A9-rank}
Let \(W\) be any simple \(8\)-dimensional \(\mathbb F_2 A_9\)-module, and let
\[
 t=(1\,2)(3\,4)(5\,6)(7\,8)\in A_9.
\]
Then
\[
 \operatorname{rank}_W(t-I)=4.
\]
\end{lemma}
\begin{proof}
There are exactly three such modules, denoted \(8a,8b,8c\) in the ATLAS. The module \(8a\) is isomorphic to the quotient of the permutation module \(\mathbb F_2^9\) by the one-dimensional space of constant vectors. Since \(9\) is odd, the permutation module is the direct sum of this space and the \(8\)-dimensional module. On \(\mathbb F_2^9\), the operator \(t-I\) has rank \(4\): each of the four transpositions contributes one independent image vector. The operator is zero on constant vectors, so it also has rank \(4\) on \(8a\).

We check the other two modules. The standard generators in the version of the ATLAS used here are
\[
 a=(1\,2\,3),\qquad b=(3\,4\,5\,6\,7\,8\,9).
\]
We use a right action: the left factor acts first, followed by the right factor. For both modules, the matrix of \(a\) is
\[
 A=\operatorname{diag}(D,D,D,D),\qquad
 D=\begin{pmatrix}0&1\\1&1\end{pmatrix}.
\]
The rows of the matrices of \(b\) on \(8b\) and \(8c\), respectively, are:
\[
\begin{array}{c|cc}
\text{row}&B_b&B_c\\ \hline
1&\mathtt{10000000}&\mathtt{10000000}\\
2&\mathtt{00100000}&\mathtt{00100000}\\
3&\mathtt{00000100}&\mathtt{00001000}\\
4&\mathtt{10100001}&\mathtt{11000010}\\
5&\mathtt{01001100}&\mathtt{00000001}\\
6&\mathtt{00001010}&\mathtt{00101100}\\
7&\mathtt{01010100}&\mathtt{00011011}\\
8&\mathtt{10010001}&\mathtt{00110100}.
\end{array}
\]
Each string of eight zeros and ones represents a row of a matrix over \(\mathbb F_2\).

Direct multiplication of permutations gives
\[
 w=bab\,a^2b^2ab^2a^2b^2a
   =(1\,2)(3\,4)(5\,6)(7\,8)=t.
\]
Thus the matrix of \(t\) is obtained by substituting \(a=A,\ b=B_b\) or \(B_c\) into this word. For \(w(A,B_b)-I\) and \(w(A,B_c)-I\), Gaussian elimination over \(\mathbb F_2\) gives the following nonzero rows of reduced row echelon form:
\[
\begin{array}{c|cc}
 &8b&8c\\ \hline
1&\mathtt{10000100}&\mathtt{10000001}\\
2&\mathtt{01000011}&\mathtt{01000111}\\
3&\mathtt{00010001}&\mathtt{00110011}\\
4&\mathtt{00001010}&\mathtt{00001000}.
\end{array}
\]
The other four rows are zero. In each column the four rows are linearly independent, so the rank is \(4\). This checks all three modules.
\end{proof}

\begin{lemma}\label{lem:S9-order16}
Let \(0\ne V\) be an \(\mathbb F_2 S_9\)-module whose restriction to \(A_9\) is a direct sum of simple \(8\)-dimensional modules. Then
\[
 16\in\omega(V\rtimes S_9).
\]
\end{lemma}
\begin{proof}
Take \(x=(1\,2\,3\,4\,5\,6\,7\,8)\in S_9\). We have
\[
 x^4=(1\,5)(2\,6)(3\,7)(4\,8).
\]
This involution is conjugate in \(A_9\) to \(t\) from Lemma~\ref{lem:A9-rank}. Indeed, all permutations of type \(2^4 1\) are conjugate in \(S_9\); if a conjugating permutation is odd, one can multiply it by an odd transposition in the centralizer of \(t\) to make it even.

Let \(M\) be the operator induced by \(x\) on \(V\), and put \(d=\dim V\). By Lemma~\ref{lem:A9-rank}, on each simple summand the rank of \(x^4-I\) equals half its dimension. Consequently,
\[
 \operatorname{rank}(M^4-I)=d/2.
\]
Put \(N=M-I\). In characteristic \(2\),
\[
 N^8=M^8-I=0,\qquad N^4=M^4-I.
\]
Let the sizes of the Jordan blocks of the nilpotent operator \(N\) be \(\ell_1,\ldots,\ell_s\). All \(\ell_i\leqslant8\), and
\[
 \sum_i\ell_i=d,\qquad
 \operatorname{rank}N^4
   =\sum_i\max(\ell_i-4,0).
\]
For \(1\leqslant\ell\leqslant8\),
\[
 \max(\ell-4,0)\leqslant\ell/2,
\]
with equality only for \(\ell=8\). Since the total rank equals \(d/2\), all blocks have size \(8\). In particular, \(N^7\ne0\). In characteristic \(2\),
\[
 I+M+M^2+\cdots+M^7=(M-I)^7=N^7\ne0.
\]
Choose \(v\in V\) with \(N^7v\ne0\). In the semidirect product,
\[
 (v,x)^8=(N^7v,1)\ne1.
\]
The right-hand side has order \(2\), while the image of \((v,x)\) in \(S_9\) has order \(8\). Hence \((v,x)\) has order \(16\).
\end{proof}

\begin{theorem}\label{thm:S9-complete}
Let \(\omega(G)=\Om\), \(R=R(G)\), and \(G/R\cong S_9\). Then \(G\) has a quotient isomorphic to \(H=(A_9\times C_3):2\). In particular, if \(G\) is \(\Om\)-critical, then \(G\cong H\).
\end{theorem}
\begin{proof}
If \(R=1\), then \(G\cong S_9\), contrary to \(21\in\omega(G)\). Thus \(R\ne1\). Choose \(K\) maximal among the proper subgroups of \(R\) normal in \(G\). Then \(V=R/K\) is a nontrivial minimal normal subgroup of \(E=G/K\). It is solvable and hence elementary abelian of characteristic \(p\in\{2,3,5,7\}\). The action of \(E\) on \(V\) factors through \(E/V\cong S_9\), and \(V\) is an irreducible \(\mathbb F_pS_9\)-module.

As in the proof of Proposition~\ref{prop:S9-three}, the restriction of \(V\) to \(A_9\) is completely reducible. Every simple summand \(W\) gives a section \(W\rtimes A_9\) of \(V\rtimes S_9\). Lemmas~\ref{lem:split-spectrum} and~\ref{lem:sections} give
\[
 \omega(W\rtimes A_9)
 \subseteq\omega(V\rtimes S_9)
 \subseteq\omega(E)\subseteq\Om.
\]
Theorem 1 and Table 1 in \cite{ZhLM} give exactly two possibilities for these simple summands:
\begin{enumerate}
\item \(p=2\) and \(\dim W=8\);
\item \(p=3\), with \(W\) the one-dimensional trivial module.
\end{enumerate}
There are no admissible summands for \(p=5\) or \(p=7\); the trivial module of characteristic \(2\) is also inadmissible because \(C_2\times A_9\) contains an element of order \(18\).

For \(p=2\), Lemma~\ref{lem:S9-order16} applies. It gives \(16\in\omega(V\rtimes S_9)\), and Lemma~\ref{lem:split-spectrum} transfers this order to \(\omega(E)\subseteq\Om\). This contradicts \(16\notin\Om\).

Thus \(p=3\). The proof of Proposition~\ref{prop:S9-three}, up to its final use of criticality, establishes \(E=G/K\cong H\): first \(V\cong C_3\), then the preimage of \(A_9\) is \(C_3\times A_9\), and the outer involution must invert \(C_3\) to avoid order \(24\). Thus the required quotient \(G/K\) exists without assuming criticality. If \(G\) is critical, then \(\omega(G/K)=\omega(H)=\Om\) forces \(K=1\).
\end{proof}

\subsection{The complete analysis of the quotient \(A_9\)}

Here the condition \(\omega(G)\subseteq\Om\) is weaker than equality of spectra. This allows induction on the radical without requiring every intermediate quotient to retain the orders \(20\) and \(21\).

\begin{lemma}\label{lem:A9-good}
Let \(W\ne0\) be a simple \(\mathbb F_2A_9\)-module with \(\omega(W\rtimes A_9)\subseteq\Om\). Then \(W\) is isomorphic to one of \(8b,8c\), and every element of \(A_9\) of order \(5\) or \(9\) acts on \(W\) without nonzero fixed vectors.
\end{lemma}
\begin{proof}
Table 1 in \cite{ZhLM} gives dimension \(8\) for \(W\), leaving the three modules in Lemma~\ref{lem:A9-rank}. We exclude \(8a=\mathbb F_2^9/\langle(1,\ldots,1)\rangle\). For \(z=(1\,2\,3\,4)(5\,6\,7)(8\,9)\), we have \(|z|=12\), and
\[
 (1+z+\cdots+z^{11})[e_1]=[e_1+e_2+e_3+e_4]\ne0.
\]
Thus \(8a\rtimes A_9\) contains an element of order \(24\notin\Om\).

For the other two modules, use the matrices \(A,B_b,B_c\) from Lemma~\ref{lem:A9-rank}. The word
\[
 f=(a^2bab^2)^3=(3\,4\,5\,7\,8)
\]
has order \(5\). Direct Gaussian elimination gives
\[
 \operatorname{rref}\bigl((A^2B_bAB_b^2)^3-I\bigr)=I_8,
\]
and the same equality with \(B_c\) in place of \(B_b\). In other words, both operators \(f-I\) are invertible. All \(5\)-cycles are conjugate in \(A_9\): an odd conjugating permutation can be adjusted by a transposition of two fixed points. Hence the assertion holds for every element of order \(5\).

If an element of order \(9\) fixed a nonzero vector, that vector of order \(2\) would commute with it, giving order \(18\notin\Om\).
\end{proof}

Below we use the Schur--Zassenhaus theorem (a normal Hall subgroup has a complement; if that subgroup is solvable, all complements are conjugate), Maschke's theorem, and the usual correspondence between irreducible complex representations of a group of order coprime to the characteristic and absolutely irreducible representations in that characteristic. This correspondence is compatible with group automorphisms; see \cite{CurtisReiner}.

\begin{lemma}\label{lem:fixedfree-characters}
Let \(q\) be a prime, \(P\) a finite \(q\)-group, and \(a\) an automorphism of \(P\) of order coprime to \(q\), with \(C_P(a)=1\). Then the only irreducible ordinary character of \(P\) fixed by \(a\) is the principal character. The same holds for isomorphism types of absolutely simple \(P\)-modules over an algebraically closed field of characteristic different from \(q\).
\end{lemma}
\begin{proof}
We show that every \(a\)-invariant conjugacy class of \(P\) contains a fixed element. If the class \(g^P\) is invariant, then \(C_{P\rtimes\langle a\rangle}(g)\) maps onto \(\langle a\rangle\). Its kernel \(C_P(g)\) is a normal Hall \(q\)-subgroup. A complement to this kernel is also a complement to \(P\) in \(P\rtimes\langle a\rangle\). By Schur--Zassenhaus, it is conjugate to \(\langle a\rangle\) by an element of \(P\). The corresponding conjugate of \(g\) is centralized by \(a\), and hence equals \(1\). Thus exactly one conjugacy class is fixed.

The ordinary character table is invertible and intertwines the permutation action on conjugacy classes with the action on irreducible characters. These two permutation matrices have equal traces, so their numbers of fixed classes and fixed characters coincide. The principal character is fixed, and there are no other fixed irreducible characters. The module assertion follows from the correspondence above in characteristic not dividing \(|P|\).
\end{proof}

\begin{lemma}\label{lem:isotypic-norm}
Let \(P\trianglelefteq X\) be a \(2\)-group, let \(V\) be a module over a field of odd characteristic \(p\), and let \(a\in X\) have order \(m\). If, after extending the field to an algebraic closure, \(a\) cyclically permutes \(m\) distinct isotypic components of the restriction of \(V\) to \(P\), then
\[
 pm\in\omega(V\rtimes X).
\]
\end{lemma}
\begin{proof}
The restriction is completely reducible by Maschke's theorem. Choose a nonzero vector \(v\) in one of these components. The vectors \(v,av,\ldots,a^{m-1}v\) lie in distinct direct summands, so their sum is nonzero. Thus the operator \(1+a+\cdots+a^{m-1}\) is nonzero after extending the field, and hence over the original field as well. For a suitable \(u\in V\), we obtain \((u,a)^m\ne1\); this element has order \(p\). Since the image of \((u,a)\) has order \(m\), the element itself has order \(pm\).
\end{proof}

\begin{lemma}\label{lem:A9-three-radical}
Let \(Q\trianglelefteq E\) be a \(3\)-group, \(E/Q\cong A_9\), and \(\omega(E)\subseteq\Om\). Then \(E\cong A_9\times Q\).
\end{lemma}
\begin{proof}
Put \(U=Q/\Phi(Q)\), where \(\Phi(Q)\) is the Frattini subgroup; the case \(Q=1\) is obvious. The action on \(U\) factors through \(A_9\). Lemma~\ref{lem:split-spectrum}, applied to \(E/\Phi(Q)\), gives \(\omega(U\rtimes A_9)\subseteq\Om\). Every simple composition factor \(W\) of \(U\) gives a section \(W\rtimes A_9\) of this group. By Table 1 in \cite{ZhLM}, only the one-dimensional trivial module is admissible in characteristic \(3\). Hence the image of \(A_9\) on \(U\) is conjugate to a subgroup of upper unitriangular matrices over \(\mathbb F_3\). This image is a \(3\)-group, whereas \(A_9\) is perfect, so the image is trivial.

The kernel of \(\Aut(Q)\to\Aut(Q/\Phi(Q))\) is a \(3\)-group. We recall a short proof. Fix a basis of \(Q/\Phi(Q)\) consisting of \(d\) elements. Every ordered set of lifts of that basis generates \(Q\); otherwise it would lie in a maximal subgroup and would not generate the quotient by \(\Phi(Q)\). There are \(|\Phi(Q)|^d\) such ordered sets. The kernel acts freely on them, since an automorphism fixing a generating set is the identity. Thus its order is a power of \(3\).

Let \(L\) be the perfect residual of \(E\). Then \(L=L'\), and the image of \(L\) in \(E/Q=A_9\) equals \(A_9\). The image of \(E\) in \(\Aut(Q)\) lies in this \(3\)-kernel, so the image of the perfect group \(L\) is trivial. Hence \([L,Q]=1\) and \(E=LQ\). The subgroup \(Z=L\cap Q\) is central in \(L\). If \(Z\ne1\), choose \(Z_0<Z\) of index \(3\). It is normal in \(L\), and \(L/Z_0\) is a perfect central cover of \(A_9\) with kernel \(C_3\), subspectral to \(S_{10}\). By the same table, the only admissible cover with this kernel is \(C_3\times A_9\), which is not perfect. This is a contradiction. Thus \(Z=1\) and \(E=L\times Q\).
\end{proof}

\begin{theorem}\label{thm:A9-structure}
Let \(R=R(G)\), \(G/R\cong A_9\), and \(\omega(G)\subseteq\Om\). Then there is a normal \(2\)-subgroup \(P\) such that
\[
 R/P=Q\text{ is a \(3\)-group},\qquad
 G/P\cong A_9\times Q.
\]
Moreover, if \(a\in G\) has order \(5\) or \(9\) and its image in \(G/P\) lies in the indicated direct factor \(A_9\), then \(C_P(a)=1\).
\end{theorem}
\begin{proof}
We argue by induction on \(|R|\). For \(R=1\), the assertion is obvious. Suppose that \(R\ne1\), and let \(V\leqslant R\) be a minimal normal subgroup of \(G\). It is elementary abelian of characteristic \(p\in\{2,3,5,7\}\). Put \(G_0=G/V\). Its radical equals \(R/V\), so induction gives \(P_0\trianglelefteq G_0\) with
\[
 P_0\text{ a \(2\)-group},\qquad
 G_0/P_0\cong A_9\times Q_0,
\]
where \(Q_0\) is a \(3\)-group. The action on \(V\) makes it a simple \(\mathbb F_pG_0\)-module, and Lemma~\ref{lem:split-spectrum} gives
\[
 \omega(V\rtimes G_0)\subseteq\omega(G)\subseteq\Om.
 \tag{*}
\]
For every cyclic subgroup of order \(5\) or \(9\) in the factor \(A_9\), its preimage in \(G_0\) has a complement to \(P_0\) by Schur--Zassenhaus. Thus \(G_0\) contains elements of these orders whose images lie entirely in \(A_9\). By induction, they act fixed-point-freely on \(P_0\).

\emph{The cases \(p=5,7\).}
First suppose that the action of \(P_0\) on \(V\) is nontrivial. After extending the field, choose a nontrivial simple \(P_0\)-constituent type. An element \(a_5\in G_0\) of order \(5\) whose image lies entirely in \(A_9\) does not fix this type, by Lemma~\ref{lem:fixedfree-characters}. Its orbit has length \(5\), and Lemma~\ref{lem:isotypic-norm} gives \(5p\in\omega(V\rtimes G_0)\). This is the forbidden order \(25\) or \(35\). If \(P_0\) acts trivially, the action factors through \(A_9\times Q_0\), and \(V\rtimes(A_9\times Q_0)\) is a quotient of \(V\rtimes G_0\). Every simple composition factor \(W\) of the restriction of \(V\) to \(A_9\) gives a section \(W\rtimes A_9\) with admissible spectrum. But Table 1 in \cite{ZhLM} allows no such modules in characteristics \(5,7\). Thus both characteristics are excluded.

\emph{The case \(p=3\).}
We show that \(P_0\) acts trivially on \(V\). Suppose otherwise. Let \(\Sigma\) be the set of simple \(P_0\)-constituent types after extending the field; it contains a nontrivial type. Let \(X\) be the preimage of \(A_9\times1\) in \(G_0\), so that \(X/P_0\cong A_9\). The set \(\Sigma\) is \(X\)-invariant, and \(P_0\) acts trivially on it because inner automorphisms preserve simple-module types. Choose \(a_9\in X\) of order \(9\). If \(a_9^3\) fixed every type in \(\Sigma\), the kernel of the action of \(X\) on \(\Sigma\), containing \(P_0\) and \(a_9^3\), would have a nontrivial normal image in the simple group \(A_9\). Thus this kernel would equal \(X\). Then \(a_9\) would fix every type in \(\Sigma\), contrary to Lemma~\ref{lem:fixedfree-characters} and the existence of a nontrivial type. Hence \(a_9^3\) moves some type. Its orbit under \(a_9\) has length \(9\), and Lemma~\ref{lem:isotypic-norm}, together with (*), gives the forbidden order \(27\).

Thus \(P_0\) centralizes \(V\). Its preimage \(Y\) in \(G\) is a central extension of \(P_0\) by \(V\). Schur--Zassenhaus gives \(Y=V\times P\), where \(P\) is the unique Sylow \(2\)-subgroup of \(Y\). Uniqueness follows because complements are conjugate by elements of the central subgroup \(V\). Therefore \(P\trianglelefteq G\), \(P\cong P_0\), and \(R/P\) is a \(3\)-group. Lemma~\ref{lem:A9-three-radical} gives \(G/P\cong A_9\times(R/P)\). The centralizer condition transfers from \(P_0\) to the isomorphic subgroup \(P\). The indicated direct factor \(A_9\) maps into the corresponding direct factor of \(G_0/P_0\): its homomorphism into the solvable group \(Q_0\) is trivial, and its image in \(A_9\) is isomorphic to \(A_9\). Orders \(5,9\) are preserved on passing to \(G_0\), since they are already preserved in the subsequent image in the corresponding \(A_9\) factor.

\emph{The case \(p=2\).}
The preimage \(P\) of \(P_0\) in \(G\) is a normal \(2\)-subgroup, and \(G/P\cong A_9\times Q_0\). The normal \(2\)-subgroup \(P_0\) acts trivially on the simple module \(V\) of characteristic \(2\). Indeed, the action of a \(2\)-group on the vector set gives a nonzero space \(C_V(P_0)\), and normality makes this space \(G_0\)-invariant; irreducibility makes it all of \(V\). Thus \(V\) is a simple module for \(A_9\times Q_0\). By Clifford's theorem, its restriction to the normal factor \(A_9\) is completely reducible. Every simple summand \(W\) gives a section \(W\rtimes A_9\) of the group in (*). Lemma~\ref{lem:A9-good} shows that elements of orders \(5,9\) whose images lie entirely in \(A_9\) have no nonzero fixed vectors on \(V\).

Now let \(a\in G\) have one of these orders and image entirely in \(A_9\) in \(G/P\). Its image \(a_0\in G_0\) has the same order, and induction gives \(C_{P_0}(a_0)=1\). If \(y\in C_P(a)\), its image in \(P_0\) is \(1\), whence \(y\in V\). But \(C_V(a)=0\), so \(y=1\). The induction is complete.
\end{proof}

\begin{corollary}\label{cor:A9-no20}
If \(G/R(G)\cong A_9\) and \(\omega(G)\subseteq\Om\), then \(20\notin\omega(G)\). In particular, \(G/R(G)\cong A_9\) is impossible when \(\omega(G)=\Om\).
\end{corollary}
\begin{proof}
Take \(P,Q\) from Theorem~\ref{thm:A9-structure}, and suppose that \(|x|=20\). The radical \(R\) is a \(\{2,3\}\)-group, so the image of \(x\) in \(G/R=A_9\) still has order divisible by \(5\). There are no elements of order \(20\) in \(A_9\), so this image has order \(5\) or \(10\). Consequently, \(x^{10}\in R\) has order \(2\) and lies in \(P\), since \(R/P\) is a \(3\)-group. But \(a=x^4\) has order \(5\) and image entirely in \(A_9\leqslant G/P\): its projection into \(Q\) is trivial. Thus \(1\ne x^{10}\in C_P(a)\), contrary to Theorem~\ref{thm:A9-structure}.
\end{proof}

\Needspace{10\baselineskip}
\begin{corollary}\label{cor:A9-complete}
Let \(G\) be nonsolvable and \(\Om\)-critical, with \(\Soc(G/R(G))\cong A_9\). Then \(G\cong H\).
\end{corollary}
\begin{proof}
By Corollary~\ref{cor:simple-socle}, the quotient \(G/R(G)\) is almost simple. Since \(\Aut(A_9)=S_9\), it equals \(A_9\) or \(S_9\). Corollary~\ref{cor:A9-no20} excludes the first case; Theorem~\ref{thm:S9-complete} applies in the second.
\end{proof}

\section{\texorpdfstring{The socles \(S_6(2),O_8^+(2),U_3(5),U_4(3)\)}{The socles S6(2),O8+(2),U3(5),U4(3)}}
\subsection{Excluding the socles \(S_6(2)\) and \(O_8^+(2)\)}

Here \(S_6(2)=\operatorname{Sp}_6(2)\) is the simple symplectic group, not the symmetric group \(S_6\). The reduction below treats an arbitrary solvable radical.

\begin{lemma}\label{lem:Sp6-modules}
Let \(S=\operatorname{Sp}_6(2)\), \(p\in\{2,3,5,7\}\), and let \(W\ne0\) be a simple \(\mathbb F_pS\)-module. If \(\omega(W\rtimes S)\subseteq\Om\), then \(p=2\), \(\dim W=8\), and \(W\) is isomorphic to the unique simple \(8\)-dimensional \(\mathbb F_2S\)-module.
\end{lemma}
\begin{proof}
We prove the assertion directly from the ordinary character table and the four decomposition matrices \cite{Sp6Data}. After extending the field to an algebraic closure, it suffices to check all absolutely simple constituents, since the dimension of a fixed-point space is unchanged by field extension. For an element \(g\) of order \(m\) coprime to the characteristic, we use
\[
 c_\varphi(g)=\frac1m\sum_{j=0}^{m-1}\varphi(g^j).
\]
If \(c_\varphi(g)>0\), the semidirect product contains an element of order \(pm\).

\emph{Characteristic \(2\).}
For a representative of the unique class \(9A\), the square and all powers coprime to \(9\) again belong to \(9A\), while the cube belongs to \(3B\). Hence
\[
 c_\varphi(9A)=
 \frac{\varphi(1)+2\varphi(3B)+6\varphi(9A)}9.
\]
The decomposition matrix gives the complete table:
\[
\begin{array}{c|rrrrrrrr}
\dim W&1&6&8&14&48&64&112&512\\ \hline
c_\varphi(9A)&1&0&0&2&6&6&12&58.
\end{array}
\]
Thus every module except those of dimensions \(6,8\) gives the forbidden order \(18\). The \(6\)-dimensional module is unique and is the natural symplectic module. We exclude it without calculating characters on \(2\)-singular elements. Write the symplectic space as the orthogonal direct sum of a plane and a four-dimensional space. On the plane choose a symplectic operator of order \(3\). On the four-dimensional space with basis \(e_1,e_2,e_3,e_4\), define a nondegenerate alternating form by
\[
 (e_1,e_3)=(e_2,e_4)=1,
\]
with all other pairwise values zero, subject to symmetry in characteristic \(2\). The basis permutation \((1\,2\,3\,4)\) preserves this form and has order \(4\). The direct sum of the two operators is an element \(z\in S\) of order \(12\), with
\[
 (1+z+\cdots+z^{11})e_1=e_1+e_2+e_3+e_4\ne0.
\]
Hence the semidirect product contains order \(24\). Only the unique \(8\)-dimensional module remains; a realization over \(\mathbb F_2\) is given below.

\emph{Characteristic \(3\).}
For the classes \(8A,8B\), the fixed-space formula becomes
\[
 c_A=\frac{\varphi(1)+\varphi(2B)+2\varphi(4A)+4\varphi(8A)}8,
 \qquad
 c_B=\frac{\varphi(1)+\varphi(2B)+2\varphi(4D)+4\varphi(8B)}8.
\]
The full list of absolutely simple modules gives:
\[
\begin{array}{r|rr@{\qquad}r|rr}
\dim W&c_A&c_B&\dim W&c_A&c_B\\ \hline
1&1&1&91&13&13\\
7&1&1&98&10&10\\
14&2&2&189&23&23\\
21&3&3&189&25&23\\
27&5&3&189&25&27\\
34&4&6&196&22&24\\
35&7&5&405&47&49\\
49&5&5&&&
\end{array}
\]
Every one of these modules gives the forbidden order \(24\).

\emph{Characteristic \(5\).}
There is exactly one class \(7A\) in \(S\), so
\[
 c_\varphi(7A)=\frac{\varphi(1)+6\varphi(7A)}7.
\]
All \(27\) irreducible Brauer characters give positive values:
\[
\begin{array}{r|r@{\quad}r|r@{\quad}r|r}
\dim W&c_\varphi&\dim W&c_\varphi&\dim W&c_\varphi\\ \hline
1&1&70&10&168&24\\
7&1&83&11&210&30\\
15&3&105&15&210&30\\
21&3&105&15&280&40\\
21&3&105&15&280&40\\
27&3&120&18&315&45\\
35&5&133&19&371&53\\
35&5&141&21&405&57\\
56&8&168&24&420&60
\end{array}
\]
Consequently, every module gives the forbidden order \(35\).

We check \(p=7\) separately. There is exactly one class of elements of order \(5\) in \(S\). For a representative \(f\) and an irreducible Brauer character \(\varphi\), the fixed-space dimension is
\[
 c_\varphi=\frac{\varphi(1)+4\varphi(f)}5.
\]
\Needspace{22\baselineskip}
The ordinary character table of \(S\) and its decomposition matrix modulo \(7\) in \cite{Sp6Data} give the following data; repeated rows correspond to distinct characters:
\[
\begin{array}{rrr|rrr}
\varphi(1)&\varphi(f)&c_\varphi&
\varphi(1)&\varphi(f)&c_\varphi\\ \hline
1&1&1&168&-2&32\\
7&2&3&189&-1&37\\
15&0&3&189&-1&37\\
21&1&5&189&-1&37\\
21&1&5&201&1&41\\
26&1&6&210&0&42\\
35&0&7&210&0&42\\
35&0&7&280&0&56\\
56&1&12&280&0&56\\
70&0&14&311&1&63\\
84&-1&16&315&0&63\\
94&-1&18&336&1&68\\
105&0&21&378&-2&74\\
105&0&21&420&0&84\\
105&0&21&&&
\end{array}
\]
All \(29\) values of \(c_\varphi\) are positive. To reproduce the column \(\varphi(f)\), note that outside the principal \(7\)-block all characters have defect zero and agree on \(7\)-regular elements with the corresponding ordinary characters. In the principal block, using the ATLAS numbering of ordinary characters, the six Brauer characters equal
\[
 \chi_1,\quad \chi_3,\quad \chi_6-\chi_1,\quad
 \chi_{15}-\chi_6+\chi_1,\quad
 \chi_{22}-\chi_3,\quad
 \chi_{28}-\chi_{15}+\chi_6-\chi_1
\]
on \(7\)-regular elements. This follows directly from the cited decomposition matrix. Thus an element of order \(5\) has a nonzero fixed vector in every simple module of characteristic \(7\), giving order \(35\) in the semidirect product. Extending the field does not change this conclusion.

\emph{Clarification of the source.} The statement about a noncyclic Sylow \(5\)-subgroup of \(S\) in \cite[p.~71]{ZhLM} is not used here: \(|S|=2^9\cdot3^4\cdot5\cdot7\), so that Sylow subgroup is cyclic of order \(5\). The calculation above replaces that step.
\end{proof}

\begin{lemma}\label{lem:Sp6-fixedfree}
On the module in Lemma~\ref{lem:Sp6-modules}, every element of order \(5\) or \(9\) acts without nonzero fixed vectors.
\end{lemma}
\begin{proof}
There is exactly one class in \(S\) of each of the orders \(5\) and \(9\); see \cite{Sp6Data}. For standard generators \(a,b\), representatives can be chosen as
\[
 f=(ab^3)^3,\qquad u=ab.
\]
The rows of the matrices \(A,B\) defining these generators on the unique \(8\)-dimensional ATLAS module are
\[
\begin{array}{c|cc}
\text{row}&A&B\\ \hline
1&\mathtt{01000000}&\mathtt{00100000}\\
2&\mathtt{10000000}&\mathtt{00010000}\\
3&\mathtt{00001000}&\mathtt{00000100}\\
4&\mathtt{10011100}&\mathtt{00000010}\\
5&\mathtt{00100000}&\mathtt{00011010}\\
6&\mathtt{11101100}&\mathtt{10100000}\\
7&\mathtt{00000001}&\mathtt{01111110}\\
8&\mathtt{00000010}&\mathtt{10011111}.
\end{array}
\]
Direct Gaussian elimination over \(\mathbb F_2\) gives
\[
 \operatorname{rref}\bigl((AB^3)^3-I\bigr)=I_8,
 \qquad
 \operatorname{rref}(AB-I)=I_8.
\]
Thus both difference operators are invertible, so neither group element has a nonzero fixed vector. The appendix contains code checking these equalities together with the orders of both operators.
\end{proof}

\begin{theorem}\label{thm:Sp6-radical}
Let \(\omega(E)\subseteq\Om\), \(R=R(E)\), and \(E/R\cong\operatorname{Sp}_6(2)\). Then \(R\) is a \(2\)-group. Moreover, every element of \(E\) of order \(5\) or \(9\) acts fixed-point-freely on \(R\).
\end{theorem}
\begin{proof}
We argue by induction on \(|R|\), proving both assertions simultaneously. For \(R=1\), they are obvious. Let \(V\leqslant R\) be a minimal normal subgroup of \(E\). It is elementary abelian of characteristic \(p\in\{2,3,5,7\}\). Put \(E_0=E/V\), \(P_0=R/V\). By induction, \(P_0\) is a \(2\)-group and \(C_{P_0}(a)=1\) for all \(a\in E_0\) of order \(5\) or \(9\). The subgroup \(V\) is a simple \(\mathbb F_pE_0\)-module. Lemma~\ref{lem:split-spectrum} gives
\[
 \omega(V\rtimes E_0)\subseteq\omega(E)\subseteq\Om.
 \tag{**}
\]
Elements of orders \(5,9\) in \(E_0/P_0=S\) have lifts of the same orders in \(E_0\): apply Schur--Zassenhaus to the preimage of the corresponding cyclic subgroup.

If \(p=2\), the normal \(2\)-subgroup \(P_0\) acts trivially on the simple module \(V\). Indeed, counting the orbits of the \(2\)-group on the vector set gives \(C_V(P_0)\ne0\), and normality and irreducibility yield \(C_V(P_0)=V\). Thus \(V\) is a simple \(\mathbb F_2S\)-module, and \(V\rtimes S\) is a quotient of the group in (**). Lemmas~\ref{lem:Sp6-modules} and~\ref{lem:Sp6-fixedfree} apply. Hence \(R\) is a \(2\)-group. For \(a\in E\) of order \(5\) or \(9\), its image \(a_0\in E_0\) has the same order. If \(y\in C_R(a)\), its image lies in \(C_{P_0}(a_0)=1\). Thus \(y\in V\), and Lemma~\ref{lem:Sp6-fixedfree} gives \(y=1\).

Now suppose that \(p\) is odd. If \(P_0\) acts trivially on \(V\), then \(V\) is a simple \(\mathbb F_pS\)-module with \(\omega(V\rtimes S)\subseteq\Om\), contradicting Lemma~\ref{lem:Sp6-modules}. It therefore suffices to exclude a nontrivial action of \(P_0\).

After extending the field, consider the set \(\Sigma\) of simple \(P_0\)-constituent types of \(V\). It contains a nontrivial type. For \(p=5,7\), choose \(a_5\in E_0\) of order \(5\). By induction, \(C_{P_0}(a_5)=1\), and Lemma~\ref{lem:fixedfree-characters} shows that \(a_5\) fixes no nontrivial type. Thus \(\Sigma\) contains an orbit of length \(5\). Lemma~\ref{lem:isotypic-norm} and (**) give the forbidden order \(25\) or \(35\).

It remains to treat \(p=3\). Choose \(a_9\in E_0\) of order \(9\). The subgroup \(P_0\) acts trivially on \(\Sigma\), since inner automorphisms preserve simple-module types. If \(a_9^3\) fixed all of \(\Sigma\), the kernel of the action of \(E_0\) on \(\Sigma\) would contain \(P_0\) and \(a_9^3\). Its image in the simple group \(E_0/P_0=S\) would be nontrivial, and hence equal to \(S\). The kernel would then equal \(E_0\). But \(a_9\) cannot fix a nontrivial type in \(\Sigma\), since \(C_{P_0}(a_9)=1\) and Lemma~\ref{lem:fixedfree-characters} applies. This contradiction shows that \(a_9^3\) moves some type. Its orbit under \(a_9\) has length \(9\). Lemma~\ref{lem:isotypic-norm} and (**) give \(27\in\omega(E)\), which is impossible. All odd characteristics are excluded, completing the induction.
\end{proof}

\Needspace{12\baselineskip}
\begin{corollary}\label{cor:Sp6-O8-excluded}
If \(G\) is nonsolvable and \(\omega(G)=\Om\), then \(\Soc(G/R(G))\) is isomorphic to neither \(\operatorname{Sp}_6(2)\) nor \(O_8^+(2)\). Criticality of \(G\) is not required.
\end{corollary}
\begin{proof}
By Corollary~\ref{cor:simple-socle}, \(X=G/R(G)\) is almost simple. For the socle \(S=\operatorname{Sp}_6(2)\), there are no outer automorphisms, so \(X=S\). Theorem~\ref{thm:Sp6-radical} shows that \(R(G)\) is a \(2\)-group. But \(S\) has no element of order \(21\): \(\mu(S)=\{7,8,9,10,12,15\}\). An element of odd order \(21\) in \(G\) would map injectively into \(G/R(G)\), since its cyclic subgroup intersects the \(2\)-radical trivially. This is a contradiction.

Now let the socle be \(T=O_8^+(2)\). Among almost simple groups with this socle, only \(T\) itself is subspectral to \(S_{10}\), by Table 2 in \cite{ZhLM}. Thus \(X=T\). The group \(T\) has a subgroup \(S\cong\operatorname{Sp}_6(2)\) of index \(120\); see \cite{Sp6Data}. Let \(E\) be its full preimage in \(G\). Then \(R(E)=R(G)\): solvability and normality give \(R(G)\leqslant R(E)\), while the image of \(R(E)\) in the nonsolvable simple group \(S\) is solvable and normal, and hence trivial. Since \(\omega(E)\subseteq\Om\), Theorem~\ref{thm:Sp6-radical} again shows that \(R(G)\) is a \(2\)-group. In \(T\), we also have \(\mu(T)=\{7,8,9,10,12,15\}\), so the same argument excludes order \(21\) in \(G\).
\end{proof}

\subsection{Excluding the socle \(U_3(5)\)}

In this case, we can avoid enumerating modules of characteristic \(2\). The obstruction comes from the subgroup \(5^{1+2}:8\). The trivial \(2\)-module is admissible, so we must prove that the radical is central, rather than trivial.

\begin{lemma}\label{lem:extra5-norm8}
Let \(B=P\rtimes\langle a\rangle\), where \(P\) is extraspecial of order \(5^3\) and exponent \(5\), and \(|a|=8\). Suppose that \(a\) induces an automorphism of order \(4\) on \(Z=Z(P)\cong C_5\), while \(a^4\) acts on \(P/Z\) as multiplication by \(-1\). If \(V\) is a module of characteristic \(p\ne5\) on which \(Z\) acts nontrivially, then
\[
 1+a+a^2+\cdots+a^7\ne0\quad\text{on }V.
\]
In particular, \(8p\in\omega(V\rtimes B)\).
\end{lemma}
\begin{proof}
Extend the field to an algebraic closure. Since \(p\ne5\), the restriction to \(Z\) is a direct sum of weight spaces. Choose a nontrivial weight \(\lambda\) with \(V_\lambda\ne0\). The four nontrivial linear characters of \(C_5\) are cyclically permuted by an automorphism of order \(4\). Thus
\[
 V_\lambda,\quad aV_\lambda,\quad a^2V_\lambda,
 \quad a^3V_\lambda
\]
are four distinct direct summands, and \(a^4\) preserves \(V_\lambda\).

The action of \(P\) on \(V_\lambda\) is faithful. Indeed, its kernel \(N\trianglelefteq P\) intersects \(Z\) trivially because a nontrivial character of a group of order \(5\) is faithful. Every nontrivial normal subgroup of a finite \(5\)-group intersects its center nontrivially: conjugation on \(N\) gives \(|C_N(P)|\equiv|N|\pmod5\), whence \(|C_N(P)|\geqslant5\). Thus \(N=1\).

The operator \(a^4\) on \(V_\lambda\) is not scalar. Otherwise it would commute with the faithful image of \(P\), and conjugation by \(a^4\) would act trivially on \(P\), contrary to its action as \(-1\) on \(P/Z\). In particular, \(1+a^4\ne0\) on \(V_\lambda\), including in characteristic \(2\). Choose \(v\in V_\lambda\) with \(u=(1+a^4)v\ne0\). Then
\[
 (1+a+\cdots+a^7)v=u+au+a^2u+a^3u\ne0,
\]
since the four nonzero summands lie in distinct weight spaces. An operator nonzero after field extension is nonzero over the original field as well. The formula \((v,a)^8=((1+a+\cdots+a^7)v,1)\) now gives an element of order \(8p\).
\end{proof}

\begin{lemma}\label{lem:U35-modules}
Let \(S=U_3(5)=\operatorname{PSU}_3(5)\), \(p\in\{2,3,5,7\}\), and let \(W\ne0\) be a simple \(\mathbb F_pS\)-module. If \(\omega(W\rtimes S)\subseteq\Om\), then \(p=2\) and \(W\) is the one-dimensional trivial module.
\end{lemma}
\begin{proof}
We first construct the subgroup from Lemma~\ref{lem:extra5-norm8} in \(S\). Put \(F=\mathbb F_{25}\), \(\bar x=x^5\), and consider the Hermitian form with matrix
\[
 J=\begin{pmatrix}0&0&1\\0&1&0\\1&0&0\end{pmatrix}.
\]
The group \(\operatorname{SU}_3(5)\) contains
\[
 P=\left\{
 \begin{pmatrix}1&-\bar x&y\\0&1&x\\0&0&1\end{pmatrix}:
 x,y\in F,\ y+\bar y+x\bar x=0
 \right\}.
\]
For every \(x\), exactly five values of \(y\) are possible, so \(|P|=125\). Its center consists of the matrices with \(x=0\) and has order \(5\); the commutator is given by the nondegenerate alternating form \(\bar x'x-\bar x x'\) on the two-dimensional \(\mathbb F_5\)-space \(F\). Hence \(P\) is extraspecial; its exponent is \(5\), since all these matrices are unipotent of size \(3\) in characteristic \(5\).

Take \(t\in F^\times\) of order \(24\), and put
\[
 d=\operatorname{diag}(t,t^4,t^{-5}).
\]
The equalities \(\det d=1\), \(\bar d^{\,T}Jd=J\) are checked directly. Conjugation \(u\mapsto dud^{-1}\) multiplies the parameter \(x\) by \(t^9\) and the central parameter \(y\) by \(t^6\). This gives an operator of order \(8\) on \(P/Z(P)\), and one of order \(4\) on the center. Moreover, \((t^9)^4=t^{12}=-1\). The center of \(\operatorname{SU}_3(5)\) has order \(3\), is generated by \(d^8\), and intersects \(P\) trivially. Thus the image of \(d\) in \(S\) has order \(8\), and the image of \(P\langle d\rangle\) has the required structure \(P:8\).

If \(p=2,3,7\) and \(S\) acts nontrivially on \(W\), then the action is faithful by simplicity of \(S\). In particular, the center of \(P\) acts nontrivially. Lemma~\ref{lem:extra5-norm8} gives orders \(16,24,56\), respectively, none in \(\Om\). Trivial modules of characteristics \(3,7\) are also inadmissible: an element of order \(8\) commutes with a vector of order \(p\), giving \(24\) or \(56\). Only the trivial \(2\)-module remains.

We exclude \(p=5\). Here we use the following precise special case of Theorem 1 and Table 1 in \cite{ZhLM}: \emph{for every nonzero simple \(\mathbb F_5A_7\)-module \(U\), the group \(U\rtimes A_7\) is not subspectral to \(S_{10}\)}. The group \(S\) has a subgroup \(A_7\) of index \(50\); see \cite{U35Data}. Choose any simple composition factor \(U=W_i/W_{i-1}\) of the restriction of \(W\) to \(A_7\). Then
\[
 (W_i\rtimes A_7)/W_{i-1}\cong U\rtimes A_7
\]
is a section of \(W\rtimes S\), contradicting the cited result. Complete reducibility of the restriction to the nonnormal subgroup \(A_7\) is not assumed here.
\end{proof}

\begin{theorem}\label{thm:U35-radical}
Let \(\omega(E)\subseteq\Om\), \(R=R(E)\), and \(E/R\cong U_3(5)\). Then
\[
 R\cong C_2^d\quad(d\geqslant0),\qquad
 E\cong U_3(5)\times R.
\]
\end{theorem}
\begin{proof}
We prove both assertions by induction on \(|R|\). For \(R=1\), both assertions are immediate. Choose a minimal normal subgroup \(V\leqslant R\) of \(E\). It is elementary abelian of characteristic \(p\in\{2,3,5,7\}\). By induction,
\[
 E_0=E/V\cong S\times Q,\qquad
 Q=R/V\cong C_2^e,
 \qquad S=U_3(5).
\]
The subgroup \(V\) is a simple \(\mathbb F_pE_0\)-module, and Lemma~\ref{lem:split-spectrum} gives
\[
 \omega(V\rtimes E_0)\subseteq\omega(E)\subseteq\Om.
\]
If \(p\) is odd, all operators from \(Q\) are simultaneously diagonalizable over \(\mathbb F_p\), with eigenvalues \(\pm1\). Since \(Q\) is central in \(E_0\), its simultaneous weight spaces are \(E_0\)-invariant. Irreducibility of \(V\) means that there is just one weight; \(Q\) acts by scalars. Thus the restriction of \(V\) to \(S\) is irreducible: every \(S\)-invariant subspace is automatically \(Q\)-invariant. The group \(V\rtimes S\) is a subgroup of \(V\rtimes E_0\), contradicting Lemma~\ref{lem:U35-modules} for odd \(p\).

Thus \(p=2\). The normal \(2\)-subgroup \(Q\) acts trivially on a simple module in characteristic \(2\), as shown in the proof of Theorem~\ref{thm:Sp6-radical}. Lemma~\ref{lem:U35-modules} gives \(V\cong C_2\) with trivial action of all of \(E_0\). In particular, \(R\) is a \(2\)-group.

We show that \(E\) splits as a direct product. Consider \(U=R/\Phi(R)\). The conjugation action on \(U\) factors through \(S\). Applying Lemma~\ref{lem:split-spectrum} to \(E/\Phi(R)\) and then passing to the simple module sections, Lemma~\ref{lem:U35-modules} shows that all composition factors of \(U\) are trivial. Consequently, the image of \(S\) on \(U\) is a \(2\)-group of upper unitriangular matrices. Perfection of \(S\) makes this image trivial. The kernel of \(\Aut(R)\to\Aut(R/\Phi(R))\) is a \(2\)-group; the proof by counting generating sets is the same as in Lemma~\ref{lem:A9-three-radical}.

Let \(L\) be the perfect residual of \(E\). Then \(L=L'\), \(E=LR\), and the image of \(L\) in \(\Aut(R)\) is trivial, since it is perfect and lies in that \(2\)-kernel. Thus \([L,R]=1\), and \(Z=L\cap R\) is a central \(2\)-subgroup of the perfect group \(L\), with \(L/Z\cong S\).

Here we use the standard property of the universal central extension: \emph{if a perfect group \(L\) has a central subgroup \(Z\) with \(L/Z\cong S\) for a nonabelian simple group \(S\), then \(L\) is a quotient of the universal central extension of \(S\), and \(Z\) is a quotient of the Schur multiplier \(M(S)\)}. For \(S=U_3(5)\), we have \(M(S)\cong C_3\); see \cite{U35Data}. Thus the central \(2\)-subgroup \(Z\) is trivial. We obtain \(L\cong S\) and \(E=L\times R\).

If \(R\) contains an element of order \(4\), its product with an element of order \(7\) in \(L\) has order \(28\notin\Om\). Thus the \(2\)-group \(R\) has exponent at most \(2\), and is therefore elementary abelian. The induction is complete.
\end{proof}

\Needspace{12\baselineskip}
\begin{corollary}\label{cor:U35-excluded}
If \(G\) is nonsolvable and \(\omega(G)=\Om\), then \(\Soc(G/R(G))\not\cong U_3(5)\).
\end{corollary}
\begin{proof}
Suppose otherwise and put \(R=R(G)\). By Corollary~\ref{cor:simple-socle}, \(G/R\) is almost simple. Let \(E\) be the full preimage of its socle \(S=U_3(5)\). Then \(R(E)=R\), since every solvable normal subgroup of \(E\) has trivial image in the nonsolvable simple group \(S\). Theorem~\ref{thm:U35-radical} shows that \(R\) is a \(2\)-group. By the complete table of almost simple groups subspectral to \(S_{10}\) in \cite{ZhLM}, the only possibilities here are \(S\) and \(S:2\); the row labeled \(S:2\) also accounts for its simple socle \(S\). Neither group has elements of order \(21\): \(\mu(S:2)=\{7,8,12,20\}\). But an element of order \(21\) in \(G\) would retain its order modulo the normal \(2\)-subgroup \(R\). This is a contradiction.
\end{proof}

\subsection{Excluding the socle \(U_4(3)\)}

\begin{lemma}\label{lem:U43-no-modules}
Let \(S=U_4(3)=\operatorname{PSU}_4(3)\), \(p\in\{2,3,5,7\}\), and let \(W\ne0\) be a simple \(\mathbb F_pS\)-module. Then \(\omega(W\rtimes S)\not\subseteq\Om\).
\end{lemma}
\begin{proof}
The group \(S\) has a subgroup \(A_7\) of index \(1296\); see \cite{U43Data}. For \(p=5,7\), Table 1 in \cite{ZhLM} permits no nontrivial elementary abelian irreducible cover of \(A_7\) in these characteristics. A simple composition factor of the restriction of \(W\) to \(A_7\) gives a forbidden section \(U\rtimes A_7\), as in Lemma~\ref{lem:U35-modules}.

For \(p=3\), we use the following precise special case of the same table: \emph{if \(U\ne0\) is a simple \(\mathbb F_3A_7\)-module and \(U\rtimes A_7\) is subspectral to \(S_{10}\), then \(U\) is one-dimensional and trivial}. Thus all composition factors of the restriction of \(W\) to \(A_7\) are trivial. The image of \(A_7\) on \(W\) then lies in the group of upper unitriangular matrices over \(\mathbb F_3\). The perfect group \(A_7\) has no nontrivial solvable image, so it acts trivially. The kernel of the representation of \(S\) contains \(A_7\ne1\); simplicity of \(S\) makes the entire representation trivial. But \(C_3\times S\) contains an element of order \(24\), since \(8\in\omega(S)\).

It remains to treat \(p=2\). We use the ordinary ATLAS characters and the full MOC decomposition matrix in \cite{U43Data}. Choose an element \(g\) in class \(9A\). Its cube belongs to \(3A\); the six powers coprime to \(9\) are divided equally between \(9A,9B\), three in each. Thus, for an irreducible Brauer character \(\varphi\),
\[
 c_\varphi(g)=
 \frac{\varphi(1)+2\varphi(3A)+
       3\bigl(\varphi(9A)+\varphi(9B)\bigr)}9.
\]
All twelve absolutely simple modules of characteristic \(2\) give the following table:
\[
\begin{array}{r|r|r|r}
\dim W&\varphi(3A)&\varphi(9A)+\varphi(9B)&c_\varphi(g)\\ \hline
1&1&2&1\\
20&-7&-2&0\\
34&7&2&6\\
34&7&-4&4\\
70&-11&-1&5\\
70&-11&-1&5\\
70&-11&2&6\\
70&-11&2&6\\
120&12&0&16\\
640&-8&2&70\\
640&-8&2&70\\
896&32&-2&106
\end{array}
\]
Positive values give the forbidden order \(18\). Thus only the unique absolutely simple module \(W_{20}\) of dimension \(20\) remains. It is realizable over \(\mathbb F_2\), as also follows from the integral realization below.

Let \(M\) be the reduction modulo \(2\) of the \(21\)-dimensional integral ATLAS representation with ordinary character \(\chi_2\). The decomposition matrix shows that its composition factors are exactly \(1,W_{20}\), each once. In the standard generators \(a,b\), put
\[
 x=abab^4.
\]
Direct computation on \(M\) gives
\[
 |x|=8,\qquad
 \operatorname{rank}_M(1+x+\cdots+x^7)
 =\operatorname{rank}_M(x-I)^7=2.
\]
The matrices and the complete verification are given in Appendix~\ref{app:U43-certificate}. The two nonzero rows in reduced row echelon form are encoded by \(1150594,802764\) in the binary encoding used there.

If the operator \(N=1+x+\cdots+x^7\) were zero on \(W_{20}\), its rank on \(M\) would be at most \(1\). Indeed, for a composition series \(0<M_1<M\), there are two cases. If \(\dim M_1=1\), then \(N(M)\subseteq M_1\). If \(\dim M_1=20\), then \(N(M_1)=0\), and the image of \(N\) has dimension at most \(\dim(M/M_1)=1\). Both cases contradict rank \(2\). Hence \(N\ne0\) on \(W_{20}\), and \(W_{20}\rtimes S\) contains an element of order \(16\notin\Om\).

This excludes all absolutely simple constituents. Extending the original field cannot turn a zero operator into a nonzero one, so the conclusion also holds for simple modules over \(\mathbb F_2\) that are not assumed in advance to be absolutely simple.
\end{proof}

\begin{corollary}\label{cor:U43-excluded}
If \(G\) is nonsolvable and \(\omega(G)=\Om\), then \(\Soc(G/R(G))\not\cong U_4(3)\).
\end{corollary}
\begin{proof}
Put \(R=R(G)\), and suppose that the socle is isomorphic to \(S=U_4(3)\). By Corollary~\ref{cor:simple-socle}, \(G/R\) is almost simple. Let \(E\) be the full preimage of its socle in \(G\). Then \(E/R\cong S\), \(R\) is solvable, and \(\omega(E)\subseteq\Om\). If there existed a nontrivial irreducible elementary abelian extension \(F\) of \(S\), with kernel \(V\), subspectral to \(S_{10}\), Lemma~\ref{lem:split-spectrum} would give \(\omega(V\rtimes S)\subseteq\omega(F)\), contrary to Lemma~\ref{lem:U43-no-modules}. Therefore Lemma~\ref{lem:no-covers} gives \(R=1\). Proposition~\ref{prop:radical-zero} then yields \(G\cong S_{10}\), which is incompatible with socle \(U_4(3)\).
\end{proof}

\section{General reductions for radical actions} \subsection{Modules when the radical is nilpotent}

In this subsection we exclude the ten remaining socles under the additional assumption that the solvable radical is nilpotent. Nilpotence is not deduced from criticality; it is used as an explicit hypothesis.

\begin{lemma}\label{lem:nilpotent-chief} Let \(R=R(G)\) be nilpotent, let \(X=G/R\), and let \(T_p\) be a Sylow \(p\)-subgroup of \(R\). Then the following assertions hold. \begin{enumerate} \item Each group \(T_p\) is normal in \(G\), and distinct \(T_p\) commute. The action of \(G\) on \(V_p=T_p/\Phi(T_p)\) factors through \(X\), and \[
 \omega(V_p\rtimes X)\subseteq\omega(G).
\] The same inclusion holds for every chief factor \(U\) of \(G\) lying inside \(T_p\), with the action of \(X\) on \(U\). \item If \(a\in G\) has order coprime to \(p\), then the natural map \[
 C_{T_p}(a)\longrightarrow C_{V_p}(a)
\] is surjective. In particular, if \(T_p\ne1\) and \(C_{V_p}(a)\ne0\), then the order \(p|a|\) belongs to \(\omega(G)\). \item Factoring out \(R_{p'}=\prod_{q\ne p}T_q\) preserves the orders of all \(p\)-elements of \(G\). \end{enumerate} \end{lemma} \begin{proof} In a finite nilpotent group, the Sylow subgroups are characteristic and form a direct product. Inner automorphisms of \(T_p\) act trivially on \(T_p/\Phi(T_p)\), and the other Sylow subgroups of \(R\) centralize \(T_p\). The first inclusion now follows from Lemma~\ref{lem:split-spectrum} applied to \(G/(R_{p'}\Phi(T_p))\).

Let \(U=N/K\) be a chief factor inside \(T_p\). The normal \(p\)-subgroup induced by \(T_p\) acts trivially on the simple module \(U\) of characteristic \(p\): its fixed-point subspace is nonzero and \(G\)-invariant. The other \(T_q\) centralize \(U\). Thus \(R\) acts trivially on \(U\). By Lemma~\ref{lem:split-spectrum}, \(U\rtimes(G/N)\) is subspectral to \(G\), and \(U\rtimes X\) is a quotient of it.

For the second assertion, write \(T=T_p\), \(F=\Phi(T)\). If \(vF\) is fixed by \(a\), then in \(T\langle a\rangle\) the subgroups \(\langle a\rangle\) and \(\langle a\rangle^v\) are complements to the normal Hall subgroup \(F\) in \(F\langle a\rangle\). By the Schur--Zassenhaus theorem, there exists \(u\in F\) such that \(\langle a\rangle^v=\langle a\rangle^u\). Consequently, \(vu^{-1}\) normalizes \(\langle a\rangle\). Conjugation by an element of \(T\) preserves the image of \(a\) modulo \(T\). Since \(T\cap\langle a\rangle=1\), this normalizing element actually centralizes \(a\). It belongs to the required coset \(vF\). A nonzero fixed vector gives a nontrivial \(p\)-subgroup in \(C_T(a)\); by Cauchy's theorem, this subgroup contains an element of order \(p\) commuting with \(a\).

The final assertion follows because a cyclic \(p\)-subgroup intersects the \(p'\)-kernel \(R_{p'}\) trivially. \end{proof}

\Needspace{10\baselineskip} \begin{lemma}\label{lem:frob54-seven} If \(F=C_5:C_4\) and the complement acts faithfully on \(C_5\), then for every nonzero \(\mathbb F_7F\)-module \(V\), the group \(V\rtimes F\) contains an element of order \(35\) or \(28\). \end{lemma} \begin{proof} If \(C_5\) acts trivially, it commutes with a nonzero vector of order \(7\), giving \(35\). Otherwise, after extending the field, there is a nonzero weight space for a nontrivial character of \(C_5\). A generator \(b\) of the complement of order \(4\) cyclically permutes the four nontrivial characters. For a nonzero vector \(v\) in one of these spaces, the vectors \(v,bv,b^2v,b^3v\) belong to distinct direct summands. Hence \((1+b+b^2+b^3)v\ne0\), and \((v,b)\) has order \(28\). \end{proof}

\subsubsection{Required fixed-point tables}

The following lemma collects the module-theoretic facts used below. Dimensions refer to absolutely simple modules, with repetitions retained. The notation \(c_m\) denotes the dimension of the subspace fixed by the specified element of order \(m\).

\begin{lemma}\label{lem:nilpotent-tables} The following assertions hold. \begin{enumerate} \item In characteristic \(5\), every nonzero \(L_2(7)\)-module has a nonzero fixed vector for every element of order \(3\). \item In characteristic \(5\), every nonzero \(L_2(8)\)-module has a nonzero fixed vector for every element of order \(7\). \item Every simple \(\mathbb F_3J_2\)-module gives the forbidden order \(24\) in its semidirect product with \(J_2\). \item If a simple \(\mathbb F_3U_3(3)\)-module \(W\) satisfies \(\omega(W\rtimes U_3(3))\subseteq\Om\), then elements of order \(7\) act on \(W\) without nonzero fixed vectors. \item In characteristic \(3\), every nonzero \(L_3(4)\)-module has a nonzero fixed vector for every element of order \(7\). \end{enumerate} \end{lemma} \begin{proof} We use the ordinary character tables and decomposition matrices in \cite{NilpotentData}. For an element \(g\) of order \(m\) coprime to the characteristic, \[
 c_m=\frac1m\sum_{j=0}^{m-1}\varphi(g^j).
\] In the first two cases, the characteristic \(5\) does not divide the order of the simple group, so ordinary and absolutely simple modular characters correspond. The complete lists of dimensions and fixed-point dimensions are as follows: \[
\begin{array}{c|rrrrrr}
L_2(7):\ \dim W&1&3&3&6&7&8\\ \hline
c_3&1&1&1&2&3&2
\end{array}
\qquad
\begin{array}{c|rrrrrrrrr}
L_2(8):\ \dim W&1&7&7&7&7&8&9&9&9\\ \hline
c_7&1&1&1&1&1&2&1&1&1
\end{array}
\] In the first table, the character values on the class \(3A\) are \(1,0,0,0,1,-1\). In the second table, the characters of degree \(7\) vanish on all \(7\)-classes, the character of degree \(8\) has value \(1\), and for each character of degree \(9\) the sum of the values on the six nonidentity powers of an element of order \(7\) is \(-2\).

For \(J_2\) in characteristic \(3\), we obtain the following complete table: \[
\begin{array}{c|rrrrrrrrrrrrrr}
\dim W&1&13&13&21&21&36&57&57&63&90&133&189&189&225\\ \hline
c_8&1&1&1&3&3&6&7&7&11&12&17&23&23&25
\end{array}
\] Every entry is positive, giving \(24\).

For \(U_3(3)\) in characteristic \(3\), the complete table is \[
\begin{array}{c|rrrrrrrrr}
\dim W&1&3&3&6&6&7&15&15&27\\ \hline
c_8&1&0&0&0&0&1&2&2&5\\
c_7&1&0&0&0&0&1&3&3&3
\end{array}
\] Positive values of \(c_8\) are ruled out by the order \(24\). In all remaining columns, \(c_7=0\). Both calculations require only the classes of \(2\)-elements and the sums of values on the \(7\)-classes; the decomposition matrices uniquely recover the Brauer characters used.

Finally, for \(L_3(4)\) in characteristic \(3\), we have \[
\begin{array}{c|rrrrrrrrr}
\dim W&1&15&15&15&19&45&45&63&63\\ \hline
c_7&1&3&3&3&1&6&6&9&9.
\end{array}
\] The three characters of degree \(15\) are \(\chi_3-\chi_2,\chi_4-\chi_2,\chi_5-\chi_2\), the character of degree \(19\) is \(\chi_2-\chi_1\), and the others coincide with the ordinary characters of degrees \(1,45,45,63,63\) on \(3\)-regular elements. This gives the last row directly.

For arbitrary modules, nonzero fixed-point spaces on composition factors give a nonzero fixed-point space on the whole module: restriction to a cyclic group of order coprime to the characteristic is completely reducible. A nonzero operator and its rank are preserved under field extension. \end{proof}

\subsubsection{Elements of order \(3\) among the automorphisms of \(L_3(4)\)}

\begin{lemma}\label{lem:L34-outer3} Let \(S=L_3(4)\), \(S\leqslant X\leqslant\Aut(S)\), and let \(W\) be a simple \(X\)-module over an algebraically closed field of characteristic \(2\), with \(S\) acting nontrivially on \(W\). Then every element of order \(3\) in \(X\) has a nonzero fixed vector on \(W\). \end{lemma} \begin{proof} Put \(Y=\operatorname{PGL}_3(4)\). The standard description of automorphisms of linear groups gives \(Y\trianglelefteq\Aut(S)\), \(Y/S\cong C_3\), \(\Aut(S)/Y\cong C_2^2\). Therefore the element \(a\) of order \(3\) under consideration belongs to \(Y\).

We state the special case of Steinberg's theorem \cite{Steinberg} used here. \emph{All absolutely simple \(\operatorname{SL}_3(4)\)-modules in characteristic \(2\) have the form \(U\otimes V^{(2)}\), where \(U,V\in\{1,N,N^*,A\}\), \(N\) is the natural three-dimensional module, \(A=\mathfrak{sl}_3\) is the eight-dimensional adjoint module, and \((2)\) denotes twisting by the field automorphism \(z\mapsto z^2\). All sixteen modules are simple and pairwise nonisomorphic.} Checking the action of the center of order \(3\) gives exactly six modules for \(S\): \[
 1,\quad A,\quad A^{(2)},\quad
 N\otimes N^{(2)},\quad
 N^*\otimes(N^*)^{(2)},\quad A\otimes A^{(2)}.
\] Their dimensions \(1,8,8,9,9,64\) are independently confirmed by the complete decomposition matrix in \cite{NilpotentData}. All six representations extend to \(Y\): the actions of \(\operatorname{GL}_3(4)\) in the indicated tensor models are trivial on scalar matrices from \(\mathbb F_4^\times\).

We determine the eigenvalues of \(a\). Choose a preimage \(D\in\operatorname{GL}_3(4)\). Then \(D^3=cI\), where \(c\in\mathbb F_4^\times\). If \(c=1\), the eigenvalues of \(D\) belong to \(\{1,\omega,\omega^2\}\). Since \(a\ne1\), their multiplicities, up to permutation, are \((2,1,0)\) or \((1,1,1)\). If \(c\ne1\), the polynomial \(z^3-c\) is irreducible over \(\mathbb F_4\); the eigenvalues are \(\lambda,\lambda\omega,
\lambda\omega^2\), each occurring once. Here \(\omega^3=1\), \(\omega\ne1\).

The multiplicities of the eigenvalues \(1,\omega,\omega^2\) on the extended modules are, up to interchange of the last two entries, as follows: \[
\begin{array}{c|cc}
\text{module}&(2,1,0)\text{ on }N&
\text{three distinct eigenvalues on }N\\ \hline
A,\ A^{(2)}&(4,2,2)&(2,3,3)\\
N\otimes N^{(2)},\ (N\otimes N^{(2)})^*&(5,2,2)&(3,3,3)\\
A\otimes A^{(2)}&(24,20,20)&(22,21,21).
\end{array}
\] For the first row, we take the ratios of the eigenvalues on \(N\) in \(\operatorname{End}(N)\) and remove the one-dimensional space of scalars. The second and third rows follow by multiplying eigenvalues in tensor products. If \(c\ne1\), the common factor in the second row merely cyclically permutes the three equal multiplicities.

Every nontrivial simple \(S\)-module extends to \(Y\). All simple \(Y\)-modules lying above it are obtained by multiplying this extension by one of the three linear characters of \(Y/S\). This follows from Schur's lemma and the cyclicity of \(Y/S\): the action on the multiplicity space reduces to a simple module for a cyclic group and hence is one-dimensional. Multiplication by a linear character cyclically permutes the three multiplicities, which are positive in every row.

Finally, \(S\langle a\rangle\) is \(S\) or \(Y\) and is normal in \(X\). By Clifford's theorem, the restriction of a simple \(X\)-module is completely reducible; all its constituents are nontrivial on \(S\). On each constituent, \(a\) has a nonzero fixed vector. \end{proof}

\begin{theorem}\label{thm:nilpotent-exclusion} Let \(G\) be nonsolvable, let \(\omega(G)=\Om\), and suppose that \(R(G)\) is nilpotent. Then \(\Soc(G/R(G))\) is not isomorphic to any of the following ten groups: \[
 A_5,A_6,A_7,A_8,L_2(7),L_2(8),U_3(3),U_4(2),L_3(4),J_2.
\] \end{theorem} \begin{proof} Put \(R=R(G)\), \(X=G/R\), \(S=\Soc(X)\). By Corollary~\ref{cor:simple-socle}, the group \(X\) is almost simple. We use the notation \(T_p\) and Lemma~\ref{lem:nilpotent-chief}.

\medskip\noindent\textbf{1. \(S=A_5\).} The group \(X\) is \(A_5\) or \(S_5\), so \(T_7\ne1\). Nilpotence of \(R\) means that \(T_2\) commutes with \(T_7\). If \(T_2\) contains an element of order \(4\), the forbidden order \(28\) occurs. Thus \(T_2\) has exponent at most \(2\). If \(X=A_5\), the square of an element of order \(8\) would belong to \(T_2\) and have order \(4\), a contradiction. If \(X=S_5\), take a chief factor inside \(T_7\). Its semidirect product with \(X\) is subspectral to \(G\). But \(S_5\) contains \(5:4\), so Lemma~\ref{lem:frob54-seven} applies.

\medskip\noindent\textbf{2. \(S=A_6\).} The order of \(\Aut(A_6)\) is not divisible by \(7\), so \(T_7\ne1\). By Theorem 1 and Table 1 in \cite{ZhLM}, there is no nonzero simple \(\mathbb F_7A_6\)-module whose semidirect product with \(A_6\) is subspectral to \(S_{10}\). A chief factor inside \(T_7\) and any simple composition factor of its restriction to \(S\) contradict this assertion and Lemma~\ref{lem:nilpotent-chief}.

\medskip\noindent\textbf{3. \(S=A_7\) or \(A_8\).} The group \(X\leqslant S_7\) or \(S_8\) has no element of order \(9\). Factor \(G\) by \(R_{3'}\). The order \(9\) is preserved, the radical of the resulting group is \(T=T_3\), and the quotient by it is \(X\). For brevity, retain the name \(G\). Let \(E\) be the preimage of \(S\). Table 1 in \cite{ZhLM} allows only the trivial simple \(A_7\)- or \(A_8\)-module in characteristic \(3\). Hence all composition factors of \(T/\Phi(T)\) are trivial on \(S\). The image of the perfect group \(S\) on this space, being a \(3\)-group, is trivial.

As in Theorem~\ref{thm:U35-radical}, the kernel \(\Aut(T)\to\Aut(T/\Phi(T))\) is a \(3\)-group. The perfect residual \(L\) of \(E\) therefore centralizes \(T\), with \(E=LT\) and \(L/(L\cap T)\cong S\). The Schur multipliers are \(M(A_7)=C_6\), \(M(A_8)=C_2\); see \cite{NilpotentData}. The property of perfect central extensions stated in Theorem~\ref{thm:U35-radical} shows that \(L\) is \(A_7\), \(3.A_7\), or \(A_8\). None of these groups has an element of order \(9\); for \(3.A_7\) this follows from \(\mu(3.A_7)=\{12,15,21\}\) in \cite{ZhLM}. The group \(L\) contains an element of order \(5\) centralizing \(T\). Therefore \(T\) has exponent at most \(3\), since otherwise \(45\) would occur. Every \(3\)-element \(x\in E\) can be written as \(lt\), where \(t\in T\) and \(l\in L\) has \(3\)-power order: its image in \(S\) has order \(1\) or \(3\), and \(L\cap T\) is a \(3\)-group. Then \(|l|\leqslant3\), \(|t|\leqslant3\), \([l,t]=1\), whence \(|x|\leqslant3\). Since \([X:S]\leqslant2\), all \(3\)-elements of \(G\) belong to \(E\). The order \(9\) is impossible.

\medskip\noindent\textbf{4. \(S=L_2(7)\).} The order of \(X\leqslant\operatorname{PGL}_2(7)\) is not divisible by the prime \(5\), so \(T_5\ne1\). Choose \(x\in G\) of order \(9\). Its image in \(X\) has order \(1\) or \(3\); in the second case it belongs to \(S\), since \([X:S]\leqslant2\). By Lemma~\ref{lem:nilpotent-tables}, the action of this image on the nonzero space \(T_5/\Phi(T_5)\) has a nonzero fixed vector. If the image is trivial, this is immediate. Lemma~\ref{lem:nilpotent-chief} gives an element of order \(5\) in \(C_{T_5}(x)\), yielding \(45\).

\medskip\noindent\textbf{5. \(S=L_2(8)\).} Here too the order of \(X\leqslant S:3\) is not divisible by \(5\), so \(T_5\ne1\). Choose an element \(a\) of order \(7\) with nontrivial image in \(S\). Such a preimage exists: any preimage of an element of order \(7\) has order \(7,14\) or \(21\), and its \(7\)-part retains a nontrivial image. By Lemma~\ref{lem:nilpotent-tables}, the space \(T_5/\Phi(T_5)\) has a nonzero fixed vector for \(a\). Lifting this vector gives the forbidden order \(35\).

\medskip\noindent\textbf{6. \(S=U_4(2)\).} The order of \(X\leqslant\Aut(S)\) is not divisible by \(7\), so \(T_7\ne1\). The group \(S\) contains \(5:4\), the normalizer of its Sylow \(5\)-subgroup; see \cite{NilpotentData}. Apply Lemma~\ref{lem:frob54-seven} to a chief factor inside \(T_7\), whose action factors through \(X\).

\medskip\noindent\textbf{7. \(S=J_2\).} The group \(J_2\) has no element of order \(9\), and its index \([\Aut(J_2):J_2]=2\), so neither does \(X\). Thus \(T_3\ne1\). A chief factor inside \(T_3\) gives a simple \(X\)-module. By Lemma~\ref{lem:nilpotent-tables}, any simple composition factor of its restriction to \(S\) gives order \(24\), contradicting Lemma~\ref{lem:nilpotent-chief}.

\medskip\noindent\textbf{8. \(S=U_3(3)\).} The order of \(\Aut(S)\) is not divisible by \(5\), so \(T_5\ne1\). Hence \(T_7=1\): otherwise commuting elements of orders \(5\) and \(7\) in \(T_5,T_7\), respectively, would give \(35\). The group \(X\leqslant S:2\) has no element of order \(21\). Let \(|x|=21\). Since \(T_7=1\), the image of \(x\) in \(X\) retains the divisor \(7\) and has order \(7\). Therefore \(1\ne x^7\in T_3\), and \(a=x^3\) of order \(7\) centralizes \(x^7\). But by Lemmas~\ref{lem:nilpotent-chief} and~\ref{lem:nilpotent-tables}, every chief factor inside \(T_3\) has no nonzero fixed vectors for \(a\). Proceeding successively along a chief series inside \(T_3\), we obtain \(C_{T_3}(a)=1\). This is a contradiction.

\medskip\noindent\textbf{9. \(S=L_3(4)\).} We first prove \(T_5=1\). The group \(S\) contains a subgroup \(A_6\). By Table 1 in \cite{ZhLM}, the only permissible simple \(\mathbb F_5A_6\)-module is trivial. Consequently, every simple \(\mathbb F_5S\)-module with a permissible semidirect product has only trivial composition factors on \(A_6\). Perfection of \(A_6\) and simplicity of \(S\) force the entire action of \(S\) to be trivial. But then an element of order \(7\) in \(S\) gives \(35\). A chief factor inside \(T_5\), if one existed, would contradict this conclusion.

No almost simple group with socle \(S\) contains elements of orders \(9\) or \(20\). Indeed, \(\Aut(S)/\operatorname{PGL}_3(4)\) is a \(2\)-group. Odd-order semisimple elements of \(\operatorname{PGL}_3(4)\) lie in the tori \(C_3^2,C_{15},C_{21}\), which have no element of order \(9\). These three types correspond to decompositions of the three-dimensional space into irreducible summands of dimensions \(1+1+1\), \(2+1\), \(3\) for a semisimple preimage in \(\operatorname{GL}_3(4)\); the scalar subgroup of order \(3\) is factored out. Moreover, \(\Out(S)\cong C_2\times S_3\) has no element of order \(4\), and \(S\) contains no \(10\). If \(y\in\Aut(S)\) had order \(20\), its image in the outer automorphism group would have order \(1\) or \(2\); then \(y^2\in S\) would have order \(10\), which is forbidden in \(S\).

Consequently, \(T_3\ne1\). Also \(T_2\ne1\): if \(T_2=1\), then in view of \(T_5=1\) an element of order \(20\) would intersect \(R\) trivially and retain its order in \(X\). Choose \(x\in G\) of order \(9\). Its image \(\bar x\) has order \(1\) or \(3\). If the image is trivial, then \(x\in T_3\) centralizes \(T_2\), giving \(18\). Hence \(|\bar x|=3\). The absence of \(18\) and the lifting of fixed points give \[
 C_{T_2}(x)=1,\qquad
 C_{T_2/\Phi(T_2)}(\bar x)=0.
\] By Lemma~\ref{lem:L34-outer3}, all simple composition factors of the \(X\)-module \(T_2/\Phi(T_2)\) must be trivial on \(S\). Otherwise \(\bar x\) would have a nonzero fixed vector; restriction to a cyclic group of order \(3\) is completely reducible. Perfection of \(S\) then shows that \(S\) acts trivially on all of \(T_2/\Phi(T_2)\).

Let \(E\) be the preimage of \(S\) in \(G\), and let \(L\) be its perfect residual. The image of \(E\) in \(\Aut(T_2)\) belongs to the \(2\)-kernel of the map onto \(\Aut(T_2/\Phi(T_2))\). Thus the perfect group \(L\) centralizes \(T_2\) and maps onto \(S\). Choose in \(L\) an element \(a\) of order \(7\) with nontrivial image in \(S\). It centralizes a nonzero element of order \(2\) in \(T_2\). On the other hand, by Lemma~\ref{lem:nilpotent-tables}, the element \(a\) has a nonzero fixed vector on the nonzero space \(T_3/\Phi(T_3)\). Lifting it gives an element of order \(3\) in \(T_3\) centralizing \(a\). The Sylow subgroups of the radical commute, so the three elements of orders \(2,3,7\) commute pairwise and give \(42\notin\Om\).

All ten socles have been excluded. \end{proof}

\begin{corollary}\label{cor:nilpotent-classification} Let \(G\) be nonsolvable and \(\Om\)-critical, and suppose that \(R(G)\) is nilpotent. Then \[
 G\cong S_{10}\quad\text{or}\quad G\cong H.
\] \end{corollary} \begin{proof} By Corollary~\ref{cor:simple-socle}, the quotient by the radical is almost simple; the possible simple socles are listed in the proof of Theorem~\ref{thm:first-premise}. The socles \(S_6(2),O_8^+(2),U_3(5),U_4(3)\) are excluded by Corollaries~\ref{cor:Sp6-O8-excluded}, \ref{cor:U35-excluded},~\ref{cor:U43-excluded}. The ten remaining small socles are excluded by Theorem~\ref{thm:nilpotent-exclusion}. The cases \(A_{10}\), \(A_9\) give \(S_{10}\) and \(H\), respectively, by Theorem~\ref{thm:A10} and Corollary~\ref{cor:A9-complete}. Both groups are indeed critical; their radicals are \(1\) and \(C_3\), respectively, and hence are nilpotent. \end{proof}

\clearpage \subsection{A nonnilpotent radical: a precise reduction}

In this subsection the radical is not assumed nilpotent. We derive necessary conditions that will be used in the subsequent exclusions. Put \(R=R(G)\), \(F=F(R)\), where \(F(R)\) is the largest nilpotent normal subgroup of \(R\). It is characteristic in \(R\) and hence normal in \(G\). Choose a chief series of \(G\) inside the radical: \[
 1=R_0<R_1<\cdots<R_t=R.
\] Each \(V_i=R_i/R_{i-1}\) is an elementary abelian \(p_i\)-group and a simple \(\mathbb F_{p_i}G\)-module. Let \(C_G(V_i)\) denote the kernel of this action.

\begin{lemma}\label{lem:fitting-chief-kernel} For this series we have \[
 F(R)=R\cap\bigcap_{i=1}^t C_G(V_i).
\] In particular, the radical is nilpotent if and only if it acts trivially on all these chief factors. \end{lemma} \begin{proof} We first prove that \(F\) centralizes every \(V_i\). Work in \(G/R_{i-1}\). The image of \(F\) is nilpotent and normal, so each of its Sylow \(q\)-subgroups \(P_q\) is normal in the whole quotient. If \(q\ne p_i\), then \([P_q,V_i]\le P_q\cap V_i=1\). If \(q=p_i\), then \(C_{V_i}(P_q)\ne0\): when a \(p_i\)-group acts on the set of vectors, the number of fixed vectors is divisible by \(p_i\), and the zero vector is fixed. Since \(P_q\) is normal, this subspace is invariant under \(G\). Simplicity of \(V_i\) gives \(C_{V_i}(P_q)=V_i\). Thus \(F\) is contained in the right-hand side of the equality.

Denote the right-hand side by \(K\). It is normal in \(R\). For every \(i\), we have \[
 [K,K\cap R_i]\le K\cap R_{i-1}.
\] Consequently, the series \(K\cap R_i\), with repetitions removed, is a central series of \(K\). Thus \(K\) is nilpotent, whence \(K\le F(R)\). \end{proof}

\begin{lemma}\label{lem:mixed-chief-action} Let \(V=N/M\) be a chief factor inside \(R\), and let \(p=\operatorname{char}V\), \(Y=G/C_G(V)\). Then \(Y\) acts faithfully and irreducibly on \(V\), and \[
 \omega(V\rtimes Y)\subseteq\omega(G),\qquad O_p(Y)=1.
\] Here \(O_p(Y)\) is the largest normal \(p\)-subgroup. If \(X=G/R\) is almost simple with simple socle \(S\), exactly one of the following alternatives holds: \begin{enumerate} \item \(Y\) is solvable; \item \(C_G(V)\le R\), \(R(Y)=R/C_G(V)\), and \(Y/R(Y)\cong X\). \end{enumerate} If \(R\) is nonnilpotent, its image in \(Y\) is nontrivial for at least one \(V\). For such \(V\), the group \(T=R(Y)\) is nontrivial, and \(F(T)\) is a nontrivial \(p'\)-group. \end{lemma} \begin{proof} Faithfulness and irreducibility follow from the definitions. Apply Lemma~\ref{lem:split-spectrum} to \(G/M\) and its normal subgroup \(N/M\). We obtain \[
 \omega\bigl(V\rtimes(G/N)\bigr)\subseteq\omega(G/M).
\] Since \(N\le C_G(V)\), the group \(V\rtimes Y\) is a quotient of this split extension. This proves the spectral inclusion even when the original extension does not split.

Any normal \(p\)-subgroup of \(Y\) has a nonzero fixed-point subspace; normality and irreducibility imply that this subspace is all of \(V\). Faithfulness gives \(O_p(Y)=1\).

Put \(C=C_G(V)\). The normal subgroup \(CR/R\) of the almost simple group \(X\) is either trivial or contains \(S\). Indeed, a normal subgroup intersecting \(S\) trivially centralizes \(S\), and \(C_X(S)=1\). In the second case, \(G/CR\) is solvable because \(X/S\le\Out(S)\) is solvable; for all socles under consideration this also follows directly from their known outer automorphism groups. The group \(CR/C\), being an image of \(R\), is solvable, and therefore so is \(Y=G/C\). In the first case, \(C\le R\). Then \(R/C\) is precisely the radical of \(G/C\), since its quotient is isomorphic to the almost simple group \(X\), whose radical is trivial. This gives the second alternative.

The existence of a factor on which \(R\) acts nontrivially follows from Lemma~\ref{lem:fitting-chief-kernel}. Its image is a nontrivial solvable normal subgroup of \(Y\), so \(T\ne1\). A nontrivial finite solvable group has a nontrivial Fitting subgroup: take a minimal normal subgroup. Finally, \(O_p(T)\) is characteristic in \(T\) and normal in \(Y\), hence is trivial. The Sylow \(p\)-subgroup of the nilpotent group \(F(T)\) is also normal in \(Y\) and is therefore trivial. \end{proof}

The last lemma identifies the modules that must actually be studied: they include not only \(X\)-modules, but also faithful simple modules for groups \(Y\) with nontrivial solvable radical having no normal \(p\)-subgroup. A direct application of the module tables for \(S\) is not justified here.

The spectral inclusion in this lemma gives further explicit necessary conditions. If \(a\in Y\) has order \(m\), \((m,p)=1\), and \(pm\notin\Om\), then \(C_V(a)=0\); otherwise a nonzero fixed vector together with \(a\) would give order \(pm\) in \(V\rtimes Y\). The complete list of such orders \(m\in\Om\) is \[
\begin{array}{c|l}
 p& m\\ \hline
 2&9,21\\
 3&8,14,20\\
 5&7,8,9,12,14,21\\
 7&4,5,6,8,9,10,12,15,20,30.
\end{array}
\] The table applies to elements of the action image \(Y\), not automatically to elements of the same order in the whole group \(G\): their order may decrease on passage to the image.

\begin{lemma}\label{lem:mixed-five-seven} Under the assumptions of Lemma~\ref{lem:mixed-chief-action}, suppose that \(\omega(G)\subseteq\Om\), \(\{p,q\}=\{5,7\}\). If \(O_q(Y)\ne1\), then \(O_q(Y)\cong C_q\). If, in addition, \(Y\) is nonsolvable and \(p\mid |S|\), then \(O_q(Y)=1\). \end{lemma} \begin{proof} Put \(Q=O_q(Y)\). All nonidentity elements of \(Q\) have order \(q\), since \(25,49\notin\Om\). For each such element \(a\), we have \(C_V(a)=0\): otherwise \(V\rtimes Y\) would contain an element of order \(pq=35\). If \(Q\) is noncyclic, it contains \(C_q\times C_q\). Indeed, choose a central element \(z\) of order \(q\) and an element \(b\notin\langle z\rangle\); then \(\langle z,b\rangle\cong C_q^2\). Over an algebraic closure of the field of characteristic \(p\), this abelian group is diagonalizable. Each of its linear characters has nontrivial kernel, since the group of roots of unity of degree \(q\) is cyclic. Thus a nonzero vector is fixed by some element of order \(q\), which is impossible. Hence \(Q\cong C_q\).

Now suppose that \(Y\) is nonsolvable. Let \(E\) be the preimage of the socle \(S\) in \(Y\), and let \(L\) be the perfect residual of \(E\). The group \(L\) is perfect and maps onto \(S\). Its image in \(\Aut(Q)\cong C_{q-1}\) is trivial, so \([L,Q]=1\). If \(p\mid |S|\), then \(L\) contains an element of order \(p\). Together with \(Q\), it gives an element of order \(35\) already in \(Y\), a contradiction. \end{proof}

For example, for a factor of characteristic \(5\) with nonsolvable action image we obtain \(O_7(Y)=1\) if \(S\) is one of the groups \[
 A_5,A_6,A_7,A_8,U_4(2),L_3(4),J_2.
\] For characteristic \(7\), we obtain \(O_5(Y)=1\) if \(S\) is one of the groups \[
 A_7,A_8,L_2(7),L_2(8),U_3(3),U_4(2),L_3(4),J_2.
\] These are restrictions on the action image \(Y\), not assertions that the corresponding prime does not divide the order of the whole radical \(R\).

\begin{theorem}\label{thm:nonnilpotent-dichotomy} Let \(\omega(G)=\Om\), suppose that \(R=R(G)\) is nonnilpotent, and suppose that the socle \(S\) of \(G/R\) belongs to the list of ten groups. Let \(L\) denote the perfect residual of \(G\), that is, the last stable term of its derived series. Then at least one of the following possibilities holds: \begin{enumerate} \item there is a chief factor \(V\) inside \(R\) on which \(R\) acts nontrivially and for which \(G/C_G(V)\) is nonsolvable; both preceding lemmas apply to this factor; \item \([R,L]\le F(R)\). In this case \[
\begin{array}{c|c}
 S&\text{restriction to }R/F(R)\\ \hline
 A_5,A_6&7\nmid |R/F(R)|\\
 L_2(7),L_2(8),U_3(3)&5\nmid |R/F(R)|\\
 A_7,A_8,U_4(2),L_3(4),J_2&
 \pi(R/F(R))\subseteq\{2,3\}.
\end{array}
\] \end{enumerate} The possibilities are not asserted to be mutually exclusive. \end{theorem} \begin{proof} Suppose that the first possibility does not occur. On every chief factor \(V_i\) inside \(R\), either \(R\) acts trivially or the image of \(G\) on \(V_i\) is solvable. In the second case, \(L\) acts trivially, since the image of a perfect group in a solvable group is trivial. In either case \([R,L]\) centralizes \(V_i\). Lemma~\ref{lem:fitting-chief-kernel} gives \([R,L]\le F(R)\).

Put \(F=F(R)\), pass to \(G/F\), and denote images by bars. Then \([\overline R,\overline L]=1\). The image of \(L\) in \(G/R\) is \(S\): the group \(S\) is perfect and \((G/R)/S\) is solvable. Hence, for every prime \(q\mid |S|\), the group \(\overline L\) contains an element of order \(q\). If \(5\mid |S|\) and \(7\mid |\overline R|\), commuting elements of orders \(5\) and \(7\) give order \(35\) in the quotient \(G/F\). This is impossible. Similarly, \(7\mid |S|\) implies \(5\nmid |\overline R|\). The table is obtained by checking the prime divisors of the orders of the ten groups, taking \(\pi(G)=\{2,3,5,7\}\) into account. \end{proof}

\subsubsection{Properties of critical groups}

\begin{lemma}\label{lem:critical-radical-witnesses} Let \(G\) be \(\Om\)-critical, and let \(N\le R(G)\) be a minimal normal subgroup of \(G\) that is elementary abelian of characteristic \(p\). Put \[
 \mathcal M=\{8,9,12,14,20,21,30\},\qquad
 D(N)=\mathcal M\setminus\omega(G/N).
\] Then \(D(N)\ne\varnothing\), and every element of \(D(N)\) is divisible by \(p\). For distinct minimal normal \(p\)-subgroups \(N_1,N_2\), the sets \(D(N_1),D(N_2)\) are disjoint. In particular, the number of minimal normal subgroups inside \(R\) of each characteristic is bounded as follows: \[
\begin{array}{c|c|c}
 p&\{m\in\mathcal M:p\mid m\}&\text{number of subgroups}\\ \hline
 2&8,12,14,20,30&\le5\\
 3&9,12,21,30&\le4\\
 5&20,30&\le2\\
 7&14,21&\le2.
\end{array}
\] \end{lemma} \begin{proof} The set \(\mathcal M\) consists of the elements of \(\Om\) maximal under divisibility; all elements of \(\Om\) are their divisors. If \(D(N)\) were empty, the spectrum of the proper quotient \(G/N\) would equal \(\Om\), contrary to criticality.

If \(x\in G\) has order \(m\) coprime to \(p\), then \(\langle x\rangle\cap N=1\), and the image of \(x\) retains its order. Consequently, \(p\mid m\) for every \(m\in D(N)\).

Now suppose that \(m\in D(N_1)\cap D(N_2)\), and choose \(x\) of order \(m\). Its image has order strictly less than \(m\) in each quotient. Since \(N_j\) is elementary abelian, the intersection \(\langle x\rangle\cap N_j\) has order exactly \(p\). Thus \(1\ne x^{m/p}\in N_1\cap N_2\), which is impossible for distinct minimal normal subgroups. The resulting nonempty pairwise disjoint sets lie in the corresponding rows of the table, proving the bounds. \end{proof}

Thus criticality provides explicit spectral witnesses for the bottom layers. This lemma alone, however, does not bound the Fitting height and does not eliminate either branch of Theorem~\ref{thm:nonnilpotent-dichotomy}.

\subsubsection{Conditions on the action of the subgroup \(5{:}4\)}

\begin{lemma}\label{lem:central-twist-five-eight} There exists a solvable group \(E=7^4:(5{:}8)\) with nonnilpotent radical and \[
 \omega(E)=\{1,2,4,5,7,8,10\}\subset\Om.
\] The quotient of the complement \(5{:}8\) by its central involution is isomorphic to \(5{:}4\), but the action of the complement on \(7^4\) does not factor through this quotient. \end{lemma} \begin{proof} Let \(K=\mathbb F_{7^4}\) and \(V=(K,+)\). Choose \(\zeta\in K^*\) of order \(5\) and \(\lambda\in K^*\) of norm \(N_{K/\mathbb F_7}(\lambda)=-1\). Such \(\lambda\) exists: in the cyclic group \(K^*\), the norm is the power map with exponent \((7^4-1)/(7-1)\), whose image has order \(6\). Consider the \(\mathbb F_7\)-linear transformations \[
 b(v)=\zeta v,\qquad a(v)=\lambda v^7.
\] We have \(b^5=1\), \(a^4=-I\), \(|a|=8\), and also \(aba^{-1}=b^2\), since \(7\equiv2\pmod5\). Consequently, \[
 B=\langle b,a\rangle\cong C_5\rtimes C_8,
 \qquad |B|=40,
\] where the action of \(C_8\) on \(C_5\) has kernel of order \(2\). The element \(a^4\) is central in \(B\), and \(B/\langle a^4\rangle\cong5{:}4\).

We check that every nonidentity element of \(B\) acts on \(V\) without nonzero fixed vectors. The elements \(b^i\), \(i\ne0\), act by multiplication by \(\zeta^i\ne1\). The elements \(b^ia^4\) act by multiplication by \(-\zeta^i\ne1\). If \(j\not\equiv0\pmod4\), then \(b^ia^j\) is conjugate in \(B\) to \(a^j\), since \(2^j-1\not\equiv0\pmod5\). Some power of \(a^j\) equals \(a^4=-I\), so it too fixes no nonzero vector.

Put \(E=V\rtimes B\). For every \(h\ne1\), the transformation \(1-h\) is invertible on \(V\). Hence every element of the coset \(Vh\) is conjugate by an element of \(V\) to \(h\) and has the same order. The element orders in \(B\) are \(1,2,4,5,8,10\); passing to \(V\) adds order \(7\). This gives the stated spectrum. The group \(E\) is solvable. It is not nilpotent, since the Sylow \(7\)-subgroup \(V\) is not centralized by its \(7'\)-complement \(B\). Finally, \(a^4\) acts on \(V\) as \(-I\ne I\), so the action does not descend to \(5{:}4\). \end{proof}

This example is not a counterexample to the desired classification: it is solvable and has strictly smaller spectrum. It establishes a specific obstruction to the method under consideration. Lemma~\ref{lem:frob54-seven} rules out a nonzero \(\mathbb F_7[5{:}4]\)-module, but not a module for the central extension \(5{:}8\). One must first justify that the relevant central subgroup acts trivially; without this, there is no contradiction from the orders \(28\) and \(35\).

There are also nonsolvable groups subspectral to \(S_{10}\) with nonnilpotent radical, for example \(A_5\times S_4\). Its radical is \(S_4\), and its spectrum is \[
 \{1,2,3,4,5,6,10,12,15,20\}\subset\Om.
\] This follows from \(\omega(A_5)=\{1,2,3,5\}\), \(\omega(S_4)=\{1,2,3,4\}\), and the formula for the spectrum of a direct product. Thus even nonsolvability together with spectral inclusion does not force the radical to be nilpotent. Thus spectral inclusion alone cannot justify reducing to a nilpotent radical; the subsequent argument must use the full spectral equality.

\clearpage \subsection{Excluding the case without a joint action}

We now strengthen Theorem~\ref{thm:nonnilpotent-dichotomy}: for all ten remaining socles, a chief factor with jointly nontrivial actions must indeed exist. Moreover, it can be chosen inside the radical of the perfect residual. The joint actions are studied separately below.

Throughout this subsection, \(L\) is the perfect residual of \(G\), and \(R=R(G)\), \(P=R\cap L\), \(X=G/R\), \(S=\Soc(X)\). We assume that \(X\) is almost simple and that its outer part \(X/S\) is solvable. The latter holds for every socle under consideration.

\begin{lemma}\label{lem:perfect-residual-product} We have \(R(L)=P\), \(L/P\cong S\), and \[
 (LR)/P\cong S\times Q,\qquad Q=R/P.
\] If \(\omega(G)=\Om\), the spectrum of \(Q\) is contained in the corresponding row of the following table: \[
\begin{array}{c|l}
 S&\text{possible element orders of }Q\\ \hline
 A_5&1,2,3,4,5,6,10,15,30\\
 A_6&1,2,3,4,5,6,10\\
 A_7,A_8,L_2(7),L_3(4)&1,2,3\\
 L_2(8),U_4(2)&1,3\\
 U_3(3),J_2&1,2.
\end{array}
\] The assertion does not require \(R\) or \(P\) to be nilpotent. \end{lemma} \begin{proof} The image of the perfect residual under an epimorphism is the perfect residual of the image. For the almost simple group \(X\) with solvable outer part, this residual is \(S\). Hence \(L/P\cong S\). The subgroup \(R(L)\) is characteristic in \(L\), and \(L\) is characteristic in \(G\). Thus \(R(L)\) is normal in \(G\) and solvable, whence \(R(L)\le R\). The reverse inclusion \(P\le R(L)\) is immediate.

The normal subgroups \(L\) and \(R\) satisfy \([L,R]\le L\cap R=P\). After factoring out \(P\), their images commute and intersect trivially, giving the stated direct product. Consequently, for every \(m\in\omega(Q)\) and \(s\in\omega(S)\), the number \(\lcm(m,s)\) belongs to \(\Om\). In particular, \(m\in\Om\).

For completeness, the element orders maximal under divisibility in the simple groups used to check this condition are listed here: \[
\begin{array}{c|l@{\qquad}c|l}
 A_5&2,3,5&A_6&3,4,5\\
 A_7&4,5,6,7&A_8&4,6,7,15\\
 L_2(7)&3,4,7&L_2(8)&2,7,9\\
 U_3(3)&7,8,12&U_4(2)&5,9,12\\
 L_3(4)&3,4,5,7&J_2&6,7,8,10,12,15.
\end{array}
\] These are the spectra in Table~2 of \cite{ZhLM}; for alternating groups they also follow directly from cycle decompositions. Enumerating \(m\in\Om\) gives the first table. \end{proof}

To state precisely the case to be excluded, introduce the condition \[
\tag{A}\label{eq:no-mixed-in-perfect}
 \text{for every chief factor }V=N/M,
 \ N\le P,\quad
 [R,V]=1\ \text{or}\ [L,V]=1.
\] Here \([R,V]=1\) means that the conjugation action is trivial, that is, \([R,N]\le M\). Condition~\eqref{eq:no-mixed-in-perfect} concerns all chief factors inside \(P\), not only the minimal normal subgroups of \(G\).

\begin{lemma}\label{lem:no-mixed-nilpotent-perfect} Under condition~\eqref{eq:no-mixed-in-perfect}, the group \(P=R(L)\) is nilpotent. If \(K\le P\) is normal in \(G\), this condition is inherited by the quotient \(G/K\), with the images \(R,L,P\). \end{lemma} \begin{proof} Since \(P\le R\) and \(P\le L\), it centralizes every chief factor inside itself: in either alternative of condition (A), the relevant group contains \(P\). Thus any chief series of \(G\) inside \(P\) is a central series of \(P\). Hence \(P\) is nilpotent.

The radical and perfect residual of \(G/K\) are \(R/K\) and \(L/K\): the first group is the image of the radical on factoring out a solvable normal subgroup, and the second equality follows from the derived series. Chief factors inside \(P/K\) lift to chief factors of \(G\) inside \(P\). Their actions are preserved, so (A) is preserved. \end{proof}

\begin{lemma}\label{lem:coprime-fixed-general} Let a solvable automorphism group \(A\) act on a finite group \(D\), and let \(N\trianglelefteq D\) be invariant under \(A\). If \((|A|,|N|)=1\), then \[
 C_{D/N}(A)=C_D(A)N/N.
\] In particular, this holds for every invariant section of a \(p\)-group and every automorphism of order coprime to \(p\). \end{lemma} \begin{proof} Work in the semidirect product \(D\rtimes A\). If \(dN\) is fixed, the complements \(A\) and \(A^d\) lie in \(N\rtimes A\). By the Schur--Zassenhaus theorem, they are conjugate by an element \(n\in N\): here \(N\) is a normal Hall subgroup and the complement \(A\) is solvable. Thus \(dn^{-1}\) normalizes \(A\). An element of \(D\) normalizing \(A\) actually centralizes it: conjugation by this element does not change the image of \(A\) modulo \(D\), and \(A\cap D=1\). Therefore the coset \(dN\) has a fixed representative. The reverse inclusion is obvious. \end{proof}

\begin{lemma}\label{lem:perfect-automorphism-tools} The following two assertions hold. \begin{enumerate} \item If \(D=D'\), \(Z\le Z(D)\), and an automorphism \(\alpha\) induces the identity on \(D/Z\), then \(\alpha=1\). \item Suppose that a perfect group \(D\) acts on a \(p\)-group \(T\). If all composition factors of \(T/\Phi(T)\) are trivial on \(D\), then \([D,T]=1\). The same conclusion holds if \(D\) acts trivially on all factors of an invariant series inside \(T\) with elementary abelian factors. \end{enumerate} \end{lemma} \begin{proof} In the first case, \(d^{-1}d^\alpha\in Z\). Since these elements are central, the map \(d\mapsto d^{-1}d^\alpha\) is a homomorphism \(D\to Z\). Perfection of \(D\) makes it trivial.

In the second case, with respect to a suitable basis, the image of \(D\) on \(T/\Phi(T)\) lies in the group of upper unitriangular matrices over \(\mathbb F_p\) and is therefore a \(p\)-group. It is perfect, hence trivial. The kernel of the map \(\Aut(T)\to\Aut(T/\Phi(T))\) is also a \(p\)-group. Recall the proof: if \(d=\dim T/\Phi(T)\), this kernel acts freely on the set of ordered lifts of a fixed basis. Each such lift generates \(T\), and the number of lift systems is \(|\Phi(T)|^d\). The kernel order is therefore a power of \(p\). Using perfection again, we obtain \([D,T]=1\). If a series inside \(T\) is given, its intersections and images in \(T/\Phi(T)\) give a filtration with trivial factors, and the same argument applies. \end{proof}

\begin{lemma}\label{lem:coprime-radical-rigidity} Suppose that \eqref{eq:no-mixed-in-perfect} holds and that \(P\) is a \(p\)-group. If \(a\in R\) has prime order \(q\ne p\), then \([a,L]=1\). \end{lemma} \begin{proof} Put \(U=P/\Phi(P)\). The action of \(L\) on \(U\) factors through \(S=L/P\), and the action of \(a\) commutes with that of \(S\), since \([R,L]\le P\) and inner automorphisms of \(P\) act trivially on \(U\). By \(q\ne p\), we have a decomposition of \(S\)-modules \[
 U=U_0\oplus U_1,\qquad
 U_0=C_U(a),\quad U_1=[U,a].
\] There are no nonzero fixed vectors of \(a\) on \(U_1\).

We prove that \(S\) acts trivially on \(U_1\). Choose a composition series of the \(G\)-module \(U\). Its factors are chief factors of \(G\) inside \(P\). If \(L\) acts nontrivially on such a factor, then by (A) the group \(R\), and in particular \(a\), acts trivially on it. The remaining factors are trivial on \(L\). The operation \(W\mapsto[W,a]\) is exact on \(\langle a\rangle\)-modules, since \(q\ne p\). Applying it to this series gives a filtration of \(U_1\) with trivial \(S\)-factors: factors of the first type disappear, while those of the second type remain trivial on \(S\). The image of the perfect group \(S\) on \(U_1\) is therefore trivial.

In the group \(L/\Phi(P)\), the subgroup \(U_0\) is normal, since it is \(S\)-invariant. The group \[
 D=(L/\Phi(P))/U_0
\] is perfect and is a central extension of \(S\) with kernel \(U/U_0\cong U_1\). The automorphism \(a\) induces the identity on \(D/U_1\). By the first part of Lemma~\ref{lem:perfect-automorphism-tools}, it is the identity on \(D\) and hence on \(U_1\). But \(C_{U_1}(a)=0\), so \(U_1=0\).

Thus \(a\) acts trivially on \(P/\Phi(P)\). The kernel of the corresponding automorphism map for \(P\) is a \(p\)-group, whereas the order of the action of \(a\) divides \(q\). Hence \([P,a]=1\). Finally, \(a\) acts trivially on \(L/P\), so Lemma~\ref{lem:coprime-fixed-general} gives \(L=P C_L(a)\). Since \(P\le C_L(a)\), we obtain \([L,a]=1\). \end{proof}

\begin{lemma}\label{lem:normal-L34-no9} If \(\omega(B)\subseteq\Om\) and \(B\) contains a normal subgroup \(S\cong L_3(4)\), then \(9\notin\omega(B)\). \end{lemma} \begin{proof} Recall the necessary facts about automorphisms of \(S\), already used in Lemma~\ref{lem:L34-outer3} and Theorem~\ref{thm:nilpotent-exclusion}. The group \(\Aut(S)/\operatorname{PGL}_3(4)\) is a \(2\)-group, and odd-order elements of \(\operatorname{PGL}_3(4)\) lie in the tori \(C_3^2,C_{15},C_{21}\). Thus \(\Aut(S)\) has no element of order \(9\). Moreover, the subgroup of order \(3\) in \(\Out(S)\cong C_2\times S_3\) is unique and consists of diagonal outer automorphisms. Its two nonidentity elements have representatives \[
 d=\operatorname{diag}(\omega,1,1),\qquad d^{-1}
 \quad\text{in }\operatorname{PGL}_3(4),
 \qquad \omega^3=1,\ \omega\ne1.
\] Both representatives have order \(3\) and centralize the nonidentity involution \(t=I+E_{23}\) of \(S\). A representative of the trivial outer class can also be chosen to be the identity.

Now let \(x\in B\) have order \(9\), and let \(\alpha\) be the automorphism of \(S\) induced by it. Then \(\alpha^3=1\). Choose, in the same outer class, an automorphism \(d\) from the list above, so that \(d^3=1\) and \(d\) fixes an involution \(t\in S\). Since \(S\) is centerless, we identify its inner automorphisms with \(S\). There exists \(s\in S\) such that \(y=sx\) induces \(d\). Moreover, \[
 y^3=s\alpha(s)\alpha^2(s)x^3.
\] The element \(u=s\alpha(s)\alpha^2(s)\) lies in \(S\). Since \(d^3=\alpha^3=1\), it induces the trivial inner automorphism of \(S\), whence \(u=1\). Thus \(y^3=x^3\ne1\) has order \(3\), and \(|y|=9\). But \([y,t]=1\), so \(yt\) has order \(18\), a contradiction. No splitting of \(B\) over \(S\) or over its centralizer was assumed here. \end{proof}

\begin{theorem}\label{thm:mixed-chief-required} Let \(\omega(G)=\Om\), and suppose that the simple socle of \(G/R(G)\) is one of the groups \[
 A_5,A_6,A_7,A_8,L_2(7),L_2(8),U_3(3),U_4(2),L_3(4),J_2.
\] Then inside \(P=R(G)\cap L\) there is a chief factor \(V=N/M\) on which both \(R(G)\) and \(L\) act nontrivially. In particular, \(Y=G/C_G(V)\) is nonsolvable and has nontrivial radical, \(Y/R(Y)\cong G/R(G)\), and \[
 \omega(V\rtimes Y)\subseteq\Om,
 \qquad O_p(Y)=1,\quad p=\operatorname{char}V.
\] This assertion does not require criticality. \end{theorem} \begin{proof} Suppose otherwise, so that (A) holds. By Lemma~\ref{lem:no-mixed-nilpotent-perfect}, \(P\) is nilpotent. Its Sylow subgroups \(T_p\) are normal in \(G\) and commute pairwise. The radical of \(L\) is \(P\), and \(L/P\cong S\), so Lemma~\ref{lem:nilpotent-chief} applies to \(L\) with the inclusion \(\omega(L)\subseteq\Om\). It is this lemma, rather than equality of the spectra of \(L\) and \(G\), that justifies using the module tables.

Whenever \(R\) acts trivially on a chief factor \(V\) inside \(P\), we also have \[
 \omega(V\rtimes X)\subseteq\Om.
\] Indeed, Lemma~\ref{lem:split-spectrum} gives \(V\rtimes(G/N)\), and \(V\rtimes X\) is a quotient of it. This justifies applying the tables to the action of the whole almost simple group.

\medskip\noindent\textbf{1. \(S=A_5\).} The order of \(X\le S_5\) is not divisible by the prime \(7\). By Lemma~\ref{lem:perfect-residual-product}, neither is the order of \(Q=R/P\). Therefore \(T_7\ne1\). If \(L\) acts trivially on every chief factor inside \(T_7\), then Lemma~\ref{lem:perfect-automorphism-tools} gives \([L,T_7]=1\). An element of order \(5\) in \(L\) then gives the forbidden order \(35\). Thus choose a chief factor \(V\) inside \(T_7\) on which \(L\) acts nontrivially. By (A), \(R\) acts trivially on it, so it is a simple \(\mathbb F_7X\)-module.

If \(X=S_5\), its subgroup \(5{:}4\) and Lemma~\ref{lem:frob54-seven} give a contradiction. If \(X=A_5\), take \(x\in G\) of order \(8\). Its image in \(X\) has order \(1\) or \(2\). In either case, \(C_V(x)\ne0\). For a nontrivial image this can be explained without tables: an involution acting in odd characteristic without nonzero fixed vectors acts as \(-I\). But a nontrivial representation of the simple group \(A_5\) is faithful, and the image of its noncentral involution cannot be scalar. The fixed point lifts from \(V\) to \(T_7\) by Lemma~\ref{lem:coprime-fixed-general}, giving order \(7\cdot8=56\), a contradiction.

\medskip\noindent\textbf{2. \(S=A_6\).} Neither \(X\le\Aut(A_6)\) nor \(Q\) has order divisible by \(7\), so \(T_7\ne1\). But Table~1 in \cite{ZhLM} gives no permissible simple \(\mathbb F_7A_6\)-module. A chief factor of \(L\) inside \(T_7\), whose action factors through \(L/P=S\), gives a contradiction. This is the same module-theoretic fact used explicitly in case \(A_6\) of Theorem~\ref{thm:nilpotent-exclusion}.

\medskip\noindent\textbf{3. \(S=A_7\) or \(A_8\).} Factor \(G\) by the Hall \(3'\)-subgroup \(K\) characteristic in \(P\). The order \(9\) is preserved, since \(K\) has order coprime to \(3\). Retain the notation for images. Now \(P\) is a \(3\)-group. The permissible simple \(\mathbb F_3S\)-modules are trivial (Table~1 in \cite{ZhLM}). Hence the perfect group \(L\) centralizes \(P\) by Lemma~\ref{lem:perfect-automorphism-tools}.

The Schur multipliers of \(A_7,A_8\) are \(C_6,C_2\); thus the perfect group \(L\) with central \(3\)-kernel is \(A_7\), \(3.A_7\), or \(A_8\). None of these groups has an element of order \(9\); for \(3.A_7\) we use \(\mu(3.A_7)=\{12,15,21\}\), as in Theorem~\ref{thm:nilpotent-exclusion}. Every element of \(R\) induces an automorphism of \(L\) that acts trivially on \(L/P\). Since \(P\le Z(L)\), the first part of Lemma~\ref{lem:perfect-automorphism-tools} gives \([R,L]=1\).

The group \(L\) contains an element of order \(5\), so \(R\) contains no elements of order \(9\), since otherwise \(45\) would occur. The index \([G:LR]=[X:S]\le2\), so an element \(x\) of order \(9\) belongs to \(LR\). Write \(x=lr\). Its image in \(L/P=S\) has order \(1\) or \(3\); since \(P\) is a \(3\)-group, \(l\) is a \(3\)-element. Hence \(|l|\le3\), and since \([L,R]=1\), the element \(r=l^{-1}x\) is also a \(3\)-element. But \(|r|\le3\), whence \(|x|\le3\), a contradiction.

\medskip\noindent\textbf{4. \(S=L_2(7)\).} Neither \(X\le\operatorname{PGL}_2(7)\) nor \(Q\) has order divisible by \(5\), so \(T_5\ne1\). If \(L\) centralizes every chief factor inside \(T_5\), it centralizes \(T_5\), which together with an element of order \(7\) gives \(35\). Thus there is a chief factor \(V\) on which \(L\) acts nontrivially. By (A), \(R\) acts trivially on \(V\).

Choose \(x\in G\) of order \(9\). Its image in \(X\) has order \(1\) or \(3\), and in the second case lies in \(S\), since the outer group has order \(2\). By Lemma~\ref{lem:nilpotent-tables}, every nonzero \(S\)-module of characteristic \(5\) has a nonzero fixed vector for an element of order \(3\). Thus \(C_V(x)\ne0\). Lifting to the normal \(5\)-subgroup gives \(45\), a contradiction.

\medskip\noindent\textbf{5. \(S=L_2(8)\).} Neither \(X\le S{:}3\) nor \(Q\) has order divisible by \(5\), so \(T_5\ne1\). Choose an element \(a\in L\) of order \(7\) with nontrivial image in \(S\). Such an element is obtained by taking the \(7\)-part of a preimage; the absence of order \(49\) ensures that its order is \(7\). Every chief factor of \(L\) inside \(T_5\) is a \(\mathbb F_5S\)-module and has a nonzero fixed vector for \(a\) by Lemma~\ref{lem:nilpotent-tables}. Lifting it gives \(35\), a contradiction.

\medskip\noindent\textbf{6. \(S=U_4(2)\).} Neither \(X\le S{:}2\) nor \(Q\) has order divisible by \(7\), so \(T_7\ne1\). The group \(S\) contains \(5{:}4\). A chief factor of \(L\) inside \(T_7\) contradicts Lemma~\ref{lem:frob54-seven}.

\medskip\noindent\textbf{7. \(S=J_2\).} Lemma~\ref{lem:nilpotent-tables} excludes all simple \(S\)-modules in characteristic \(3\), since they give order \(24\). Hence \(T_3=1\). Nor is the order of \(Q\) divisible by \(3\), by Lemma~\ref{lem:perfect-residual-product}. Thus \(3\nmid |R|\), and the order \(9\) is preserved in \(X\le J_2{:}2\), where it does not occur.

\medskip\noindent\textbf{8. \(S=U_3(3)\).} Neither \(X\le S{:}2\) nor \(Q\) has order divisible by \(5\), so \(T_5\ne1\). Commutativity of the Sylow subgroups of \(P\) and the exclusion of \(35\) give \(T_7=1\). The table for \(Q\) implies that \(Q\) is an elementary abelian \(2\)-group. Thus \(7\nmid |R|\), and \(G/L\) is a \(2\)-group: its normal subgroup \(R/P\) and quotient \(G/LR\le\Out(S)\) are both \(2\)-groups.

Choose \(x\in G\) of order \(21\). There is no element of order \(21\) in \(X\), and the \(7\)-part is preserved, so \(x^7\in R\) has order \(3\). Since \(Q\) is a \(2\)-group, \(x^7\in T_3\). The element \(a=x^3\) of order \(7\) belongs to \(L\) and has nontrivial image in \(S\). On every chief factor of \(L\) inside \(T_3\), it acts without nonzero fixed vectors, by Lemma~\ref{lem:nilpotent-tables}. Proceeding along this series gives \(C_{T_3}(a)=1\). But \(1\ne x^7\in C_{T_3}(a)\), a contradiction.

\medskip\noindent\textbf{9. \(S=L_3(4)\).} We begin with \(T_5=T_7=1\). For \(T_5\), apply the argument of case \(L_3(4)\) of Theorem~\ref{thm:nilpotent-exclusion} to the group \(L\): the restriction of a permissible module to \(A_6<S\) in characteristic \(5\) has only trivial factors. Perfection of \(A_6\) makes its action trivial, and simplicity of \(S\) makes the whole action of \(S\) trivial. An element of order \(7\) then gives \(35\). For \(T_7\), the same subgroup \(A_6\) has no permissible simple module in characteristic \(7\). Consequently, \(P=T_2\times T_3\).

We show that \([L,T_2]=1\). Choose \(x\) of order \(9\). If the action of \(L\) on a chief factor \(V\) of \(G\) inside \(T_2\) is nontrivial, then by (A) the action of \(R\) is trivial. Thus \(V\) is a simple \(\mathbb F_2X\)-module nontrivial on \(S\). Since \(\Aut(S)\) has no element of order \(9\), the image of \(x\) in \(X\) has order \(1\) or \(3\). In the first case it fixes all of \(V\), and in the second it has a nonzero fixed vector by Lemma~\ref{lem:L34-outer3}. Lifting it to \(T_2\) gives \(18\), a contradiction. Thus all these factors are trivial on \(L\), and Lemma~\ref{lem:perfect-automorphism-tools} gives \([L,T_2]=1\).

The group \(T_3\) is nontrivial. Indeed, if \(T_3=1\), then \(P=T_2\), and the quotient \(G/P\) retains the order \(9\) and contains a normal subgroup \(L/P\cong S\). This contradicts Lemma~\ref{lem:normal-L34-no9}. If also \(T_2\ne1\), choose \(a\in L\) of order \(7\) with nontrivial image in \(S\). It centralizes \(T_2\) and, by Lemma~\ref{lem:nilpotent-tables}, has a nonzero fixed vector on every nonzero module of characteristic \(3\). Hence \(C_{T_3}(a)\ne1\). Commuting elements of orders \(2,3,7\) give \(42\), a contradiction. Thus \[
 P=T_3\ne1.
\]

By the table for \(Q\), all of \(R\) is a \(\{2,3\}\)-group. It contains an involution: otherwise \(R\) would be a \(3\)-group, and an element of order \(20\) would retain its order in \(X\). But \(\Aut(L_3(4))\) has no element of order \(20\), as proved in detail in Theorem~\ref{thm:nilpotent-exclusion}: the outer group has no element of order \(4\), while \(S\) has none of order \(10\). Choose an involution \(b\in R\). Lemma~\ref{lem:coprime-radical-rigidity}, with \(p=3,q=2\), now gives \([b,L]=1\). The group \(L\) contains an element \(a\) of order \(7\) and a nonidentity element \(c\in P\) of order \(3\) commuting with \(a\), by the same table of modules in characteristic \(3\). Since \(a,c\in L\), both are centralized by \(b\). Again \(42\) occurs.

A contradiction has been obtained for all ten socles. Thus (A) is impossible, and the required chief factor exists inside \(P\). Its action image \(Y\) is nonsolvable, since the image of the perfect group \(L\) is nontrivial. The remaining assertions follow from Lemma~\ref{lem:mixed-chief-action}. \end{proof}

\begin{corollary}\label{cor:perfect-radical-seven} If \(\omega(G)=\Om\) and the socle of \(G/R(G)\) is isomorphic to one of the groups \[
 A_6,A_7,A_8,L_2(8),U_3(3),U_4(2),J_2,
\] then the radical \(R(L)=R(G)\cap L\) of the perfect residual \(L\) is nonnilpotent. \end{corollary} \begin{proof} Suppose that \(P=R(L)\) is nilpotent. In cases 2, 3, 5, 6, 7, 8 of the proof of Theorem~\ref{thm:mixed-chief-required}, condition (A) was used only to obtain this nilpotence. All subsequent arguments in these cases apply without (A): the table for \(Q\) is unconditional, and the chief factors of \(L\) inside \(P\) are \(S\)-modules by Lemma~\ref{lem:nilpotent-chief}. Thus cases 2, 5, 6 exclude \(A_6,L_2(8),U_4(2)\), respectively, by a forbidden order or the absence of a permissible module; case 3 excludes \(A_7,A_8\) through the absence of order \(9\); case 7 excludes \(J_2\) by the orders \(24\) and \(9\); and case 8 excludes \(U_3(3)\) through the impossibility of order \(21\). This gives a contradiction in all seven cases. \end{proof}

Thus the possibility \emph{without} a joint action has been completely excluded. The condition \([R,L]\le F(R)\) from Theorem~\ref{thm:nonnilpotent-dichotomy} has not itself been ruled out: it may coexist with a joint action. The conclusion is that such a factor must exist inside \(R(L)\).

It remains to exclude faithful simple \(\mathbb F_pY\)-modules with \(R(Y)\ne1\) on which both the radical and the perfect residual act nontrivially, while preserving the inclusion \(\omega(V\rtimes Y)\subseteq\Om\). The previous restrictions \(O_p(Y)=1\) and Lemma~\ref{lem:mixed-five-seven} remain valid. The rest of the proof studies these joint actions on chief factors inside the canonical commutator kernel \(K\).

\clearpage \subsection{Module-theoretic restrictions on joint actions}

The choice of chief factor in Theorem~\ref{thm:mixed-chief-required} can be refined. This is important: an arbitrary action group may have a complicated radical, whereas for a \emph{highest} joint factor the radical of the perfect residual is already nilpotent. Below we exclude six combinations of the socle and the characteristic of such a factor. These exclusions are then strengthened to exclude two socles completely for an arbitrary radical.

\begin{lemma}\label{lem:highest-mixed} Under the assumptions of Theorem~\ref{thm:mixed-chief-required}, choose a chief series of \(G\) inside \(P=R\cap L\) and its highest factor \(V=N/M\) on which both groups \(R,L\) act nontrivially. Put \[
 p=\operatorname{char}V,\quad C=C_G(V),\quad
 Y=G/C,\quad E=Y^{(\infty)},\quad Q=R(E).
\] Then \[
 E\cong L/(L\cap C),\qquad E/Q\cong S,
 \qquad Q\text{ is nilpotent},\quad p\nmid |Q|.
\] Moreover, \(V\) is a faithful simple \(\mathbb F_pY\)-module, \(\omega(V\rtimes Y)\subseteq\Om\), \(O_p(Y)=1\), and \(R(Y)\ne1\). \end{lemma} \begin{proof} On every factor above \(N\) in the chosen series, at least one of the groups \(R,L\) acts trivially; hence so does their subgroup \(P\). Consequently, the part of the series above \(N\) gives a central series of \(P/N\), and \(P/N\) is nilpotent.

By Lemma~\ref{lem:mixed-chief-action}, the group \(C\) is contained in \(R\), the radical of \(Y\) is \(R/C\), and \(Y/R(Y)\cong X\). The perfect residual maps onto the perfect residual under an epimorphism, so \(E\) is the image of \(L\). Its solvable radical is the image of \(P\): this image is solvable, and the quotient by it is \(L/P\cong S\). Since \(N\le C\), the group \(Q\) is a quotient of the nilpotent group \(P/N\). It is characteristic in \(E\), and \(E\) is characteristic in \(Y\). Its Sylow \(p\)-subgroup is therefore normal in \(Y\) and lies in \(O_p(Y)=1\). Hence \(p\nmid |Q|\). The remaining assertions were proved in Lemma~\ref{lem:mixed-chief-action}. \end{proof}

\subsubsection{Minimal polynomials and central kernels}

We now allow any joint chief factor, not necessarily a highest one. Retain the notation \(Y=G/C_G(V)\), \(E=Y^{(\infty)}\), \(Q=R(E)\), \(E/Q\cong S\), with \(V\) a faithful simple \(\mathbb F_pY\)-module satisfying \(\omega(V\rtimes Y)\subseteq\Om\). The group \(Q\) is not assumed nilpotent. The group \(F(Q)\) is nilpotent of order coprime to \(p\): its Sylow \(p\)-subgroup is normal in \(Y\) and hence lies in \(O_p(Y)=1\).

We use the following minimal-polynomial theorem in the precise form needed below.

\begin{lemma}[Minimal polynomial theorem] \label{lem:HH-prime} Let \(s,q\) be distinct primes, and let \(H=T\rtimes\langle a\rangle\), where \(T\) is a normal \(q\)-group, \(|a|=s\), \([T,a]\ne1\). Suppose that \(H\) acts faithfully on a vector space \(W\) of positive characteristic \(\ell\ne q\). If the minimal polynomial of \(a\) on \(W\) differs from \(x^s-1\), then \[
 C_T(a)\ne1,\qquad T\text{ is nonabelian},\qquad
 q=2,\quad s\text{ is a Fermat prime}.
\] \end{lemma} \begin{proof}[Source of the result] This is Lemma~3 of \cite{StaroletovPrimeGraph}, where it is stated with reference to Lemma~3.6 of the work cited there. We use exactly this form of a Hall--Higman type theorem; its proof is not reproduced here. In particular, the statement \emph{does not require} \(\ell\ne s\). \end{proof}

If the minimal polynomial is \(x^s-1\), the operator \(1+a+\cdots+a^{s-1}\) is nonzero: otherwise the minimal polynomial would divide this polynomial of degree \(s-1\). In the semidirect product, this gives an element of order \(\ell s\). When \(\ell=s\), we obtain \(s^2\), since \(1+x+\cdots+x^{s-1}=(x-1)^{s-1}\) in characteristic \(s\). We shall apply the result both in distinct characteristics and in equal characteristics.

\begin{lemma}\label{lem:cyclic-kernel-prime} In the situation just described, let \(\{p,q\}=\{5,7\}\). The Sylow \(q\)-subgroup \(F_q\) of \(F(Q)\) is either trivial or cyclic of order \(q\) and central in \(E\). \end{lemma} \begin{proof} The group \(F_q\) is normal in \(Y\). Since \(35\notin\omega(V\rtimes Y)\), every nonidentity element of order \(q\) in it acts on \(V\) without nonzero fixed vectors. Moreover, \(q^2\notin\omega(Y)\). As in Lemma~\ref{lem:mixed-five-seven}, a subgroup \(C_q^2\) is impossible: every one of its linear characters over an algebraic closure has nontrivial kernel and gives a nonzero fixed vector. A noncyclic \(q\)-group of exponent \(q\) contains \(C_q^2\): choose a central element \(z\) of order \(q\) and an element outside \(\langle z\rangle\). Hence \(F_q\cong C_q\) if \(F_q\ne1\). The image of the perfect group \(E\) in the cyclic group \(\Aut(C_q)\) is trivial, so \(F_q\le Z(E)\). \end{proof}

\begin{lemma}\label{lem:odd-mixed-quasisimple} For any joint chief factor, suppose that \(p\in\{5,7\}\) and \(7\mid |S|\). Then \(Q\le Z(E)\); that is, \(E\) is a perfect central cover of \(S\), and \(Q\) is a quotient of \(M(S)\). \end{lemma} \begin{proof} Put \(F=F(Q)\). This is a nilpotent \(p'\)-group. Its Sylow subgroup for the other prime in the pair \(5,7\) is central in \(E\) by Lemma~\ref{lem:cyclic-kernel-prime}.

Choose \(a\in E\) of order \(7\) with nontrivial image in \(S\). Such an element exists by taking the \(7\)-part of a preimage, since \(49\notin\Om\). If \([F_q,a]\ne1\), \(q\in\{2,3\}\), then Lemma~\ref{lem:HH-prime} applies to \(F_q\rtimes\langle a\rangle\) on \(V\). The action is faithful, its characteristic is \(p\ne q\), and \(7\) is not a Fermat prime. The minimal polynomial of \(a\) is therefore \(x^7-1\). This gives \(35\) if \(p=5\) or \(49\) if \(p=7\), a contradiction. Thus \([F,a]=1\).

The subgroup \(C_E(F)\) is normal in \(E\) and contains \(a\), whose image in \(S\) is nontrivial. Consequently, \(E=Q C_E(F)\). The image of \(E\) in \(\Aut(F)\) is the image of the solvable group \(Q\) and is also perfect. Thus it is trivial, and \([E,F]=1\).

In a finite solvable group, \(C_Q(F(Q))\le F(Q)\). Recall the proof. Write \(D=C_Q(F)\), so that \(D\cap F=Z(F)\). If \(D>Z(F)\), choose a minimal nontrivial normal subgroup of \(Q/Z(F)\) inside \(D/Z(F)\). It is elementary abelian. Its preimage \(B\) is normal in \(Q\), and \(B'\le Z(F)\le Z(B)\), since \(B\le D\). Then \(B\) is nilpotent and \(B\le F(Q)\), contrary to \(B\le D\), \(B\not\le Z(F)\). Thus \(D\le F\). Now \([Q,F]=1\) gives \(Q=F\), whence \(Q\le Z(E)\). The kernel of a perfect central extension is a quotient of the Schur multiplier.

For later applications, we list the multipliers of the groups in our list whose orders are divisible by \(7\); see \cite{Atlas,MixedData}: \[
\begin{array}{c|c}
 S&M(S)\\ \hline
 A_7&C_6\\ A_8&C_2\\ L_2(7)&C_2\\
 L_2(8),U_3(3)&1\\
 L_3(4)&C_3\times C_4\times C_4\\ J_2&C_2.
\end{array}
\] \end{proof}

\begin{lemma}\label{lem:nine-kills-two} For any joint chief factor, if \(9\in\omega(S)\) and \(p\in\{3,5,7\}\), then \(O_2(Q)=1\). \end{lemma} \begin{proof} Choose \(a\in E\) mapping to an element of order \(9\) in \(S\). The only order in \(\Om\) divisible by \(9\) is \(9\), so \(|a|=9\). Put \(T=O_2(Q)\). The absence of \(18\) in \(\omega(E)\) implies \(C_T(a)=1\).

Extend the field to an algebraic closure. The restriction of \(V\) to \(T\) is completely reducible, since \(p\ne2\). The group \(E\) permutes the set \(\Sigma\) of simple isomorphism types in this restriction and their isotypic spaces. An orbit of \(a\) of length \(9\) is impossible: for a vector in one space, the vectors \(v,av,\ldots,a^8v\) would lie in distinct direct summands and have nonzero sum. Then order \(9p\in\{27,45,63\}\) would occur in \(V\rtimes E\). Thus \(a^3\) fixes all of \(\Sigma\).

Let \(K\) be the kernel of the action of \(E\) on \(\Sigma\). It contains \(a^3\), whose image in \(S\) is nontrivial, so \(KQ=E\). The permutation image of \(E\) is the image of the solvable group \(Q\) and is also perfect. Hence it is trivial. In particular, \(a\) fixes each type. By Lemma~\ref{lem:fixedfree-characters}, if \(C_T(a)=1\), the automorphism induced by \(a\), which has odd order, fixes no nontrivial simple type of \(T\). Thus \(T\) acts trivially on \(V\). Faithfulness gives \(T=1\). \end{proof}

\subsubsection{Two exclusions in characteristic \(7\)}

\begin{lemma}\label{lem:small-seven-no-modules} For \(S=L_2(7),L_2(8)\) there is no nonzero \(\mathbb F_7S\)-module \(W\) satisfying \(\omega(W\rtimes S)\subseteq\Om\). The same holds for \(S=U_3(3)\). \end{lemma} \begin{proof} First let \(S=L_2(7)\). Over an algebraic closure of the field of characteristic \(7\), the complete list of simple \(\operatorname{SL}_2(7)\)-modules is \(\operatorname{Sym}^i(k^2)\), \(0\le i\le6\). This is the special case of Steinberg's theorem used in Lemma~\ref{lem:L34-outer3}. Precisely the modules with even \(i\) descend to the quotient \(\operatorname{PSL}_2(7)\): their dimensions are \(1,3,5,7\). For \(i=0,2,4\), choose in \(S\) an element of order \(4\) whose preimage in \(\operatorname{SL}_2(7)\) has eigenvalues \(\zeta,\zeta^{-1}\), \(|\zeta|=8\). The eigenvalue \(\zeta^{i-2(i/2)}=1\) occurs in \(\operatorname{Sym}^i(k^2)\). This gives \(28\). If \(i=6\), an element of order \(7\) acts as a single Jordan block of size \(7\). This can be checked directly: on polynomials of degree at most \(6\) it acts by the substitution \(t\mapsto t+1\), and the sixth finite difference of \(t^6\) is \(6!\ne0\) in characteristic \(7\). Thus the norm operator \((a-1)^6\) is nonzero, giving \(49\). Every nonzero module has a simple composition factor, and an operator nonzero on a factor cannot be zero on the entire module.

Now let \(S=L_2(8)\). It contains a subgroup \(F=2^3{:}7\) whose complement of order \(7\) acts transitively on the seven nontrivial characters of the kernel \(A=2^3\); see \cite{MixedData}. If \(A\) acts nontrivially on \(W\), then after extending the field there is an orbit of length \(7\) on its nontrivial weights. The sum of the seven permuted vectors is nonzero, so the norm of an element of order \(7\) is nonzero, giving \(49\). If \(A\) acts trivially, the kernel of the action of \(S\) contains \(A\ne1\); simplicity of \(S\) therefore makes the entire action trivial. An element of order \(9\) in \(S\) then gives \(63\), again a contradiction.

Finally, \(U_3(3)\) contains \(L_2(7)\) as one of its maximal subgroups, so it suffices to restrict any nonzero module to this subgroup. \end{proof}

\subsubsection{The double cover of \(A_8\) in characteristic \(5\)}

\begin{lemma}\label{lem:double-A8-five} There is no nonzero \(\mathbb F_5E\)-module \(W\) such that \(\omega(W\rtimes E)\subseteq\Om\) if \(E=A_8\) or \(2.A_8\). \end{lemma} \begin{proof} For \(A_8\), restriction to \(A_7\) contradicts the absence of permissible simple \(A_7\)-modules in characteristic \(5\) in Table~1 of \cite{ZhLM}. This assertion was also used in Theorem~\ref{thm:nilpotent-exclusion}.

Consider \(2.A_8\). For a simple module on which the center acts trivially, the preceding case applies after factoring out the center. The complete list of absolutely simple modules on which the center acts nontrivially in characteristic \(5\) has dimensions \[
 8,24,24,32,32,48.
\] For the standard element \(b\) of order \(7\), the fixed-point dimensions are \[
\begin{array}{c|rrrrrr}
 \dim W&8&24&24&32&32&48\\ \hline
 \dim C_W(b)&2&3&3&5&5&6.
\end{array}
\] Completeness of the list is checked using the decomposition matrix of \(2.A_8\) modulo \(5\) in \cite{MixedData}. The numbers in the second row were independently calculated from four matrices over \(\mathbb F_5\): \[
\begin{array}{l|r|r|r}
 \text{matrix file }b&n&\operatorname{rank}(b-I)&c_7\\ \hline
 \texttt{2A8G1-f5r8B0.m2}&8&6&2\\
 \texttt{2A8G1-f5r48bB0.m2}&48&42&6\\
 \texttt{2A8G1-f5r64B0.m2}&64&54&10\\
 \texttt{2A8G1-f5r48aB0.m2}&48&42&6.
\end{array}
\] After extending to \(\mathbb F_{25}\), the second matrix decomposes into a pair of conjugate \(24\)-dimensional modules, and the third into a pair of \(32\)-dimensional modules. The ranks on conjugate modules are equal, so \(6\) and \(10\) are each divided by two. The last matrix belongs to the absolutely simple \(48\)-dimensional module. Verification code is given in Appendix \ref{app:double-A8-five}.

All values of \(c_7\) are positive. Thus every such simple module gives \(35\). A simple factor of an arbitrary nonzero module also gives a forbidden order: the corresponding semidirect product is a section. \end{proof}

\begin{theorem}\label{thm:mixed-six-pairs} For any joint chief factor inside \(R\cap L\), the following six pairs \((S,p)\) are impossible: \[
\begin{array}{c|l}
 p&S\\ \hline
 5&A_8,\ L_2(8),\ U_4(2)\\
 7&L_2(8),\ U_3(3),\ U_4(2).
\end{array}
\] \end{theorem} \begin{proof} First consider the pairs \((A_8,5)\), \((L_2(8),5)\), \((L_2(8),7)\), \((U_3(3),7)\). By Lemma~\ref{lem:odd-mixed-quasisimple}, the group \(E\) is a perfect central cover of \(S\). In the case \(A_8\), where \(M(S)=C_2\), it is \(A_8\) or \(2.A_8\). Lemma~\ref{lem:double-A8-five} excludes both possibilities. For \(L_2(8)\) and \(U_3(3)\), the multiplier is trivial, so \(E=S\). In characteristic \(5\), every nonzero \(L_2(8)\)-module has a nonzero fixed vector for an element of order \(7\) by Lemma~\ref{lem:nilpotent-tables}, giving \(35\). The other two pairs are excluded by Lemma~\ref{lem:small-seven-no-modules}.

Now suppose that \(S=U_4(2)\) and \(p=5\) or \(7\). Put \(F=F(Q)\). It is nilpotent of order coprime to \(p\). Lemma~\ref{lem:nine-kills-two} gives \(O_2(Q)=1\), so \(F_2=1\). The subgroup \(F_q\), where \(q\) is the other prime in the pair \(5,7\), is central in \(E\) by Lemma~\ref{lem:cyclic-kernel-prime}.

Choose \(a\in E\) of order \(5\) with nontrivial image in \(S\). If \([F_3,a]\ne1\), Lemma~\ref{lem:HH-prime}, with \(q=3,s=5\), makes the minimal polynomial equal to \(x^5-1\). The exception is impossible, since \(q\ne2\). If \(p=5\), we obtain \(25\); if \(p=7\), we obtain \(35\). Thus \([F_3,a]=1\), and taking the other Sylow subgroups into account gives \([F,a]=1\). Normality of \(C_E(F)\) and simplicity of \(S\) imply \(E=Q C_E(F)\). As in the proof of Lemma~\ref{lem:odd-mixed-quasisimple}, perfection of \(E\) and solvability of \(Q\) force \([E,F]=1\), and the property \(C_Q(F)\le F\) gives \(Q=F\le Z(E)\). The group \(Q\) has odd order, whereas \(M(U_4(2))=C_2\). Therefore \(Q=1\), \(E=S\). In characteristic \(7\), the subgroup \(5{:}4<S\) and Lemma~\ref{lem:frob54-seven} give \(28\) or \(35\). In characteristic \(5\), use \(A_6<S\): this subgroup lies in the maximal subgroup \(S_6<S\); see \cite{MixedData}. Every permissible simple \(A_6\)-module in characteristic \(5\) is trivial. Hence all composition factors of the restriction of \(V\) are trivial on \(A_6\), and perfection of \(A_6\) makes its entire action trivial. Simplicity of \(S\) now makes the action of all of \(S\) trivial. An element of order \(9\) in \(S\) then gives \(45\), a contradiction. \end{proof}

Thus precisely the indicated characteristics of a joint factor have been excluded. To exclude the socle completely, we must still prove that such a factor of the required characteristic necessarily occurs. The following lemma supplies this step in two cases.

\begin{lemma}\label{lem:central-chief-prime-elimination} Let \(L\trianglelefteq G\), \(L=L'\), \(P=R(L)\trianglelefteq G\), with \(L/P\cong S\) nonabelian simple. Let \(r\) be a prime and \(r\nmid |M(S)|\). If \(L\) centralizes every chief \(r\)-factor of \(G\) inside \(P\), then \(r\nmid |P|\). \end{lemma} \begin{proof} Induct on \(|P|\). The assertion is immediate if \(P=1\). Choose a minimal nontrivial normal subgroup \(N\) of \(G\) inside \(P\); it is elementary abelian of some characteristic \(q\). All assumptions are inherited by \(G/N,L/N\), the radical of the second group being \(P/N\). By induction, \(r\nmid |P/N|\). If \(q\ne r\), we immediately obtain \(r\nmid |P|\).

Let \(q=r\). By hypothesis, \(N\le Z(L)\). In \(P\), the subgroup \(N\) is a central normal Hall \(r\)-subgroup, since \(P/N\) has order coprime to \(r\). By the Schur--Zassenhaus theorem, \(P=N\times K\), where \(K\) is a \(r'\)-group. The subgroup \(K\) consists of precisely the \(r'\)-elements of this direct product, so it is characteristic in \(P\) and normal in \(G\). Thus \(L/K\) is a perfect central extension of \(S\) with kernel \(P/K\cong N\). Its kernel is a quotient of \(M(S)\), contradicting \(r\nmid |M(S)|\). The case \(q=r\) is impossible, completing the induction. \end{proof}

\begin{theorem}\label{thm:L28-U42-arbitrary-radical} Let \(G\) be nonsolvable and \(\omega(G)=\Om\). Then the simple socle of \(G/R(G)\) is isomorphic to neither \(L_2(8)\) nor \(U_4(2)\), without any assumption that the radical is nilpotent. \end{theorem} \begin{proof} Retain the notation \(R,L,P,S,X,Q_0=R/P\), where \(Q_0\) now refers to the original group, as opposed to the radical of the action group \(Y\).

Suppose that \(S=L_2(8)\). The order of \(X\le S{:}3\) is not divisible by \(5\), so \(5\mid |R|\). By Lemma~\ref{lem:perfect-residual-product}, \(\omega(Q_0)\subseteq\{1,3\}\), whence \(5\mid |P|\). Consider any chief \(5\)-factor \(V\) of \(G\) inside \(P\). If \(L\) acts nontrivially on it, then \(R\) cannot act nontrivially by Theorem~\ref{thm:mixed-six-pairs}. Hence the action factors through \(X\), and \(V\rtimes X\) is subspectral to \(G\). But every nonzero \(S\)-module of characteristic \(5\) has a nonzero fixed vector for an element of order \(7\) by Lemma~\ref{lem:nilpotent-tables}, giving \(35\). Thus \(L\) centralizes all chief \(5\)-factors inside \(P\). Since \(M(L_2(8))=1\), the preceding lemma gives \(5\nmid |P|\), a contradiction.

Suppose that \(S=U_4(2)\). Here \(7\nmid |X|\) and \(\omega(Q_0)\subseteq\{1,3\}\), so \(7\mid |P|\). On a chief \(7\)-factor \(V\) inside \(P\), a nontrivial action of \(L\) cannot coexist with a nontrivial action of \(R\) by Theorem~\ref{thm:mixed-six-pairs}. If \(R\) acts trivially, however, then \(V\rtimes X\) is subspectral to \(G\), and the subgroup \(5{:}4<S\) gives \(28\) or \(35\) by Lemma~\ref{lem:frob54-seven}. Thus \(L\) centralizes all chief \(7\)-factors inside \(P\). But \(M(U_4(2))=C_2\), and Lemma~\ref{lem:central-chief-prime-elimination} gives \(7\nmid |P|\), a contradiction. \end{proof}

\begin{corollary}\label{cor:A8-U33-radical-primes} If \(\omega(G)=\Om\) and the socle of \(G/R(G)\) is \(A_8\), then \(5\nmid |R(G)|\). If this socle is \(U_3(3)\), then \(7\nmid |R(G)|\). \end{corollary} \begin{proof} When \(S=A_8\), all chief \(5\)-factors inside \(P\) on which \(L\) acts nontrivially are excluded: a joint action is forbidden by Theorem~\ref{thm:mixed-six-pairs}, while a trivial action of \(R\) would give a permissible \(A_8\)-module in characteristic \(5\), ruled out by Lemma~\ref{lem:double-A8-five}. Thus all these factors are centralized by \(L\). Since \(M(A_8)=C_2\), Lemma \ref{lem:central-chief-prime-elimination} gives \(5\nmid |P|\). The order of \(R/P\) is not divisible by \(5\) by Lemma~\ref{lem:perfect-residual-product}, whence \(5\nmid |R|\).

For \(S=U_3(3)\), the same argument applies to the chief \(7\)-factors: joint actions are forbidden by Theorem~\ref{thm:mixed-six-pairs}, while actions trivial on the radical are excluded by Lemma~\ref{lem:small-seven-no-modules}. Using \(M(U_3(3))=1\) and the fact that \(R/P\) is a \(2\)-group gives \(7\nmid |R|\). \end{proof}

This step alone completely excludes the two socles \(L_2(8),U_4(2)\). The radical restrictions for \(A_8,U_3(3)\) are strengthened in the next section.

\clearpage \section{\texorpdfstring{The socles \(U_3(3),J_2,A_8\)}{The socles U3(3),J2,A8}} \subsection{The socles \(A_8,U_3(3),J_2\): restrictions on the radical}

For \(U_3(3),J_2\), this subsection establishes complete exclusion with arbitrary solvable radical. The key step is the occurrence of the forbidden order \(24\) in an extraspecial \(2\)-section. For \(A_8\), only the restriction \(\pi(R(G))\subseteq\{2,3\}\) is proved here. Its full exclusion is established later in this section.

\subsubsection{The radical when the socle is \(A_8\)}

\begin{lemma}\label{lem:double-A8-seven} The groups \(A_8\) and \(2.A_8\) contain a subgroup isomorphic to \(L_2(7)\). For either group \(D\), there is no nonzero \(\mathbb F_7D\)-module \(W\) satisfying \(\omega(W\rtimes D)\subseteq\Om\). \end{lemma} \begin{proof} In \(A_8\), take the action of \(L_3(2)\) on the seven nonzero vectors of \(\mathbb F_2^3\) and add one fixed point. This action is faithful, and its sign is trivial because \(L_3(2)\) is perfect. Here \(L_3(2)\cong L_2(7)\).

For \(2.A_8\), a different embedding of \(L_2(7)\) is used; it is not automatically the full preimage of the preceding one. The ATLAS page \cite{TwoLocalData} gives a transitive action of \(2.A_8\) of degree \(240\) on cosets of a subgroup \(L_2(7)\). The stabilizer has order \(40320/240=168\) and is explicitly identified in the source with \(L_2(7)\). Thus this is an actual subgroup of \(2.A_8\), not merely a section or a central quotient.

The restriction of a nonzero module \(W\) to this subgroup remains nonzero. The spectral inclusion for \(W\rtimes D\) would imply the same inclusion for \(W\rtimes L_2(7)\), contradicting Lemma~\ref{lem:small-seven-no-modules}. \end{proof}

\begin{theorem}\label{thm:A8-radical-two-three} If \(\omega(G)=\Om\) and the socle of \(G/R(G)\) is \(A_8\), then \[
 \pi(R(G))\subseteq\{2,3\}.
\] In particular, \(5,7\nmid |R(G)|\). \end{theorem} \begin{proof} The divisor \(5\) has already been eliminated by Corollary~\ref{cor:A8-U33-radical-primes}. Put \(R=R(G)\), \(L=G^{(\infty)}\), \(P=R\cap L\), \(X=G/R\). Consider any chief \(7\)-factor \(V=N/M\) of \(G\) inside \(P\).

If \(L\) acts nontrivially on \(V\) and \(R\) acts trivially, then \(V\) is a \(X\)-module and \(\omega(V\rtimes X)\subseteq\Om\) by Lemma~\ref{lem:split-spectrum}. Restriction to \(A_8\) contradicts Lemma~\ref{lem:double-A8-seven}. If both groups act nontrivially on \(V\), then for \(Y=G/C_G(V)\) the perfect residual \(E\) is a central cover of \(A_8\) by Lemma~\ref{lem:odd-mixed-quasisimple}. Since \(M(A_8)=C_2\), we have \(E=A_8\) or \(2.A_8\). In both cases, the same Lemma~\ref{lem:double-A8-seven} contradicts \(\omega(V\rtimes E)\subseteq\Om\).

Thus \(L\) centralizes all chief \(7\)-factors inside \(P\). Since \(7\nmid |M(A_8)|\), Lemma~\ref{lem:central-chief-prime-elimination} gives \(7\nmid |P|\). By Lemma~\ref{lem:perfect-residual-product}, \(\omega(R/P)\subseteq\{1,2,3\}\), so \(7\nmid |R|\). Finally, \(\pi(G)=\{2,3,5,7\}\), proving the assertion. \end{proof}

\begin{lemma}\label{lem:A8-radical-both-primes} Under the assumptions of the preceding theorem, in fact \(\pi(R(G))=\{2,3\}\), and the radical is nonnilpotent. If \(G/R(G)=A_8\), then \(R(G)\) contains an element of order \(4\). \end{lemma} \begin{proof} Neither \(A_8\) nor \(S_8\) contains elements of orders \(9,20\), as follows from their cycle lengths. Let \(x\in G\) have order \(9\). If \(3\nmid |R|\), then \(\langle x\rangle\) would embed in \(G/R\), a contradiction. Similarly, if \(2\nmid |R|\), then because \(\pi(R)\subseteq\{2,3\}\) a cyclic subgroup of order \(20\) would intersect \(R\) trivially and embed in \(G/R\). Thus both primes divide \(|R|\). A nilpotent radical has already been excluded by Theorem~\ref{thm:nilpotent-exclusion}.

Finally, if \(G/R=A_8\), the image of an element \(y\) of order \(20\) has order \(5\). Indeed, \(5\nmid |R|\), and in \(A_8\) the centralizer of a \(5\)-cycle is \(C_5\times A_3\) and contains no involutions. Consequently, \(y^5\in R\) has order \(4\). This assertion does not automatically carry over to the quotient \(S_8\), in which the image of \(y\) may have order \(10\). \end{proof}

\subsubsection{Fixed types and an element of order \(8\)}

\begin{lemma}\label{lem:fixedfree-q-characters} Let \(T\) be a finite \(q\)-group, let \(a\) be an automorphism of order coprime to \(q\), and suppose that \(C_T(a)=1\). Then the only \(a\)-fixed irreducible ordinary character of \(T\) is the principal character. The same holds for the isomorphism types of absolutely simple modules in characteristic different from \(q\). \end{lemma} \begin{proof} Repeat the argument of Lemma~\ref{lem:fixedfree-characters} without restricting to \(q=2\). If the conjugacy class \(t^T\) is invariant, then \(C_{T\rtimes\langle a\rangle}(t)\) maps onto \(\langle a\rangle\) with kernel \(C_T(t)\). By the Schur--Zassenhaus theorem, this centralizer has a complement to its normal Hall \(q\)-subgroup. It is also a complement in \(T\rtimes\langle a\rangle\) and is therefore conjugate to \(\langle a\rangle\) by an element of \(T\). The corresponding conjugate of \(t\) is fixed by \(a\) and hence equals \(1\). Exactly one class is fixed.

The invertible character table intertwines the permutation actions on conjugacy classes and irreducible characters. Their traces, that is, their numbers of fixed points, are equal. Thus exactly one character is fixed, namely the principal character. For modules, use the usual correspondence in characteristic not dividing \(|T|\), as in Lemma \ref{lem:fixedfree-characters}. \end{proof}

\begin{lemma}\label{lem:eight-fixes-types} Let \(Y\) act faithfully on a simple \(\mathbb F_3Y\)-module \(V\), and let \(\omega(V\rtimes Y)\subseteq\Om\), \(E=Y^{(\infty)}\), \(Q=R(E)\), with \(E/Q\cong S\) nonabelian simple. Suppose that \(S\) contains an element of order \(8\) and that \(T\trianglelefteq Y\) has order coprime to \(3\). Then \(E\) fixes every simple isomorphism type in the restriction of \(V\) to \(T\) over an algebraic closure of the field. If \(T\) is abelian, then \([E,T]=1\). \end{lemma} \begin{proof} The restriction to \(T\) is completely reducible. Let \(\Sigma\) be the set of simple isomorphism types occurring. Lift an element of order \(8\) from \(S\) to \(a\in E\). Its order is \(8\), since \(\Om\) contains no other multiples of \(8\). If \(a\) had an orbit of length \(8\) on \(\Sigma\), then for a nonzero vector in one isotypic space, the sum \(v+av+\cdots+a^7v\) would be nonzero. This would give \(24\in\omega(V\rtimes E)\). Thus \(a^4\) fixes all of \(\Sigma\).

The kernel \(K\) of the action of \(E\) on \(\Sigma\) is normal and contains \(a^4\), whose image in \(S\) is a nonidentity involution. Thus \(KQ=E\). The permutation image of \(E\) is the image of the solvable group \(Q\) and is also perfect. Therefore it is trivial.

If \(T\) is abelian, all types are linear characters. Their kernels intersect trivially, since the action of \(Y\), and hence of \(T\), on \(V\) is faithful. For any \(e\in E\), the equalities \(\chi(t^e)=\chi(t)\) for all occurring \(\chi\) imply \(t^e=t\). Hence \([E,T]=1\). \end{proof}

\subsubsection{A chief \(3\)-factor responsible for order \(21\)}

\begin{lemma}\label{lem:twentyone-three-witness} Let \(\omega(G)=\Om\), and let the simple socle of \(G/R(G)\) be \(S=U_3(3)\) or \(J_2\). Then there exist \(a\in L=G^{(\infty)}\) of order \(7\) and a chief \(3\)-factor \(V=N/M\) inside \(P=R(G)\cap L\) such that \[
 C_V(a)\ne0,\qquad [L,V]\ne1,
 \qquad [R(G),V]\ne1.
\] The image of \(a\) in \(G/R(G)\) has order \(7\). \end{lemma} \begin{proof} In the case \(U_3(3)\), we have \(7\nmid |R|\) by Corollary~\ref{cor:A8-U33-radical-primes}. For \(J_2\), Lemma~\ref{lem:radical} gives the same conclusion, since \(25\mid |J_2|\) and \(25,35\notin\Om\). By Lemma~\ref{lem:perfect-residual-product}, the group \(R/P\) has spectrum contained in \(\{1,2\}\) and hence is an elementary abelian \(2\)-group. Moreover, \([G:LR]=[G/R:S]\le2\). Thus \(G/L\) is a \(2\)-group.

Choose \(x\in G\) of order \(21\). There is no element of order \(21\) in \(X=G/R\le\Aut(S)\): the outer group has order at most \(2\), and the groups \(U_3(3),J_2\) themselves have no such order. Since \(7\nmid |R|\), the image of \(x\) has order \(7\). Therefore \(b=x^7\in R\) has order \(3\), while \(a=x^3\) has order \(7\), has nontrivial image in \(S\), and centralizes \(b\). The group \(G/L\) is a \(2\)-group, so \(a,b\in L\), \(b\in P\).

In a chief series of \(G\) inside \(P\), choose the first term \(N\) containing \(b\) and the preceding term \(M\). Then \(V=N/M\) has characteristic \(3\), and \(bM\) is a nonzero fixed vector of \(a\).

The group \(L\) cannot act trivially on \(V\). Indeed, a preimage in \(L\) of an element of order \(8\) from \(S\) has order \(8\). Its image in \(G/N\) also has order \(8\), since \(N\le P\) and the image in \(S\) is retained. If the action on \(V\) were trivial, Lemma~\ref{lem:split-spectrum} would give order \(24\).

If \(R\) acted trivially, then \(V\rtimes X\) would be subspectral to \(G\). When \(S=J_2\), this is excluded by Lemma~\ref{lem:nilpotent-tables} in characteristic \(3\). When \(S=U_3(3)\), the same lemma shows that every permissible simple factor of the restriction to \(S\) has zero fixed-point space for an element of order \(7\). Since restriction to \(\langle a\rangle\) is completely reducible, we obtain \(C_V(a)=0\), a contradiction. Both actions are nontrivial. \end{proof}

\begin{lemma}\label{lem:nonabelian-two-fitting} For the factor \(V\) from the preceding lemma, put \(Y=G/C_G(V)\), \(E=Y^{(\infty)}\), \(Q=R(E)\). Then \[
 F(Q)=O_2(Q)
\] is a nontrivial nonabelian \(2\)-group. \end{lemma} \begin{proof} The action of \(Y\) on \(V\) is faithful and irreducible, \(\omega(V\rtimes Y)\subseteq\Om\), \(O_3(Y)=1\), and \(E/Q\cong S\). The group \(Q\) is the image of \(P\), so \(7\nmid |Q|\). Its Fitting subgroup is normal in \(Y\); its Sylow \(3\)-subgroup would lie in \(O_3(Y)\). Thus \(\pi(F(Q))\subseteq\{2,5\}\).

Write \(O_5(Q)\) as \(T\). Lemma~\ref{lem:eight-fixes-types} shows that \(E\) fixes every simple type of \(T\) occurring in \(V\). Choose \(c\in E\) of order \(7\) with nontrivial image in \(S\). From \(35\notin\omega(E)\) we get \(C_T(c)=1\). By Lemma~\ref{lem:fixedfree-q-characters}, a nontrivial simple type of \(T\) cannot be fixed by \(c\). Thus \(T\) acts trivially on \(V\), and faithfulness gives \(T=1\). Consequently, \(F(Q)=O_2(Q)\); denote this group by \(F\).

Suppose that \(F\) is abelian, allowing also \(F=1\). By Lemma~\ref{lem:eight-fixes-types}, \([E,F]=1\). Since \(Q\) is solvable, \(C_Q(F(Q))\le F(Q)\), as proved in Lemma~\ref{lem:odd-mixed-quasisimple}. Hence \(Q=F\le Z(E)\), and \(E\) is a perfect central cover of \(S\).

If \(S=U_3(3)\), the multiplier is trivial, so \(E=S\). Every permissible module for this group in characteristic \(3\) has no nonzero fixed vectors for an element of order \(7\) by Lemma~\ref{lem:nilpotent-tables}. The image of the original \(a\) belongs to \(E\), since \(a\in L\), and has order \(7\). This contradicts \(C_V(a)\ne0\).

If \(S=J_2\), then \(E=J_2\) or \(2.J_2\), since \(M(J_2)=C_2\). The first case is excluded by the absence of permissible simple modules in characteristic \(3\), by Lemma~\ref{lem:nilpotent-tables}. In the second case, \(2.J_2\) contains a subgroup \(C_2\times U_3(3)\); see \cite{TwoLocalData}. In particular, \(E\) contains an actual subgroup \(U_3(3)\). The restriction of \(V\) to it satisfies the spectral condition and hence has no nonzero fixed vectors for its elements of order \(7\). All subgroups of order \(7\) in \(E\) are Sylow subgroups and are conjugate; one lies in the indicated \(U_3(3)\). Hence no element of order \(7\) in \(E\), including the image of \(a\), fixes a nonzero vector of \(V\). This is again a contradiction. Thus \(F\) is nonabelian and, in particular, nontrivial. \end{proof}

\subsubsection{An extraspecial section and the forbidden order \(24\)}

\begin{lemma}\label{lem:extraspecial-eight-coset} Let \(D\) be an extraspecial \(2\)-group of order \(2^{1+2n}\), with \(H=D\rtimes\langle a\rangle\), \(|a|=8\). Suppose that an automorphism \(a\) induces a transformation of order \(8\) on \(D/Z(D)\), and that every element of the coset \(Da\) has order \(8\). If a simple module \(W\) for \(H\) over an algebraically closed field of characteristic \(3\) restricts faithfully and irreducibly to \(D\), then some element of \(Da\) fixes a nonzero vector of \(W\). \end{lemma} \begin{proof} Write \(N=2^n\). This is the degree of the unique faithful irreducible character of the extraspecial group \(D\). Recall the justification: the commutator defines a nondegenerate alternating form on \(D/Z(D)\); the preimage of a maximal totally isotropic subspace is abelian of order \(2^{n+1}\). Inducing a linear character of this preimage nontrivial on \(Z(D)\) gives an irreducible character of degree \(2^n\). All \(2^{2n}\) linear characters are trivial on the center, and the sum of their squared degrees together with \((2^n)^2\) is already \(|D|\). Thus the nonlinear character constructed is unique and faithful.

In characteristic \(2\), a transformation of order \(8\) must have a unipotent Jordan block of size at least \(5\). Consequently, \(2n\ge6\), \(N\ge8\).

Since \(H\) is a \(2\)-group, the characteristic \(3\) does not divide its order. Ordinary and modular irreducible representations correspond, preserving eigenvalue multiplicities of elements of \(H\). We may use a complex representation \(\rho:H\to\operatorname{GL}_N(\mathbb C)\) with character \(\chi\): occurrence of the eigenvalue \(1\) is preserved under reduction. This is the usual correspondence in coprime characteristic: roots of unity of \(2\)-power order reduce injectively, and the group algebra is semisimple.

For \(j=0,\ldots,7\), put \[
 S_j=\sum_{g\in Da}\chi(g^j).
\] We have \(S_0=|D|N\). If \(j\) is odd, the map \(g\mapsto g^j\) takes \(Da\) bijectively onto \(Da^j\): its inverse is the power map whose exponent is inverse to \(j\) modulo the exponent of \(H\). But \(\sum_{d\in D}\rho(d)=0\), since \(C_W(D)=0\). Hence \[
 S_1=S_3=S_5=S_7=0.
\]

We prove the bounds \(|S_2|,|S_6|\le |D|\). The faithful character of \(D\) is real-valued: it has values \(N,-N\) on the two central elements and vanishes outside the center. Therefore the representation of \(D\) has a nondegenerate invariant bilinear form with matrix \(B\), where \(B^t=\varepsilon B\), \(\varepsilon\in\{1,-1\}\). Let \(A\) realize any normalizing operator of finite order on this representation. Uniqueness of the form gives \(A^tBA=\lambda B\), where \(|\lambda|=1\). Choose unitary representation matrices. Invariance of the form gives \(\rho(d)=B^{-1}\overline{\rho(d)}B\). Orthogonality of matrix coefficients yields \[
 \sum_{d\in D}\rho(d)_{ij}\rho(d)_{kl}
 =\frac{|D|}{N}(B^{-1})_{ki}B_{jl}
 =\frac{\varepsilon|D|}{N}(B^{-1})_{ik}B_{jl}.
\] Multiplying by the coefficients \(A_{jk}\) and summing over \(j,k\) gives \[
 \sum_{d\in D}\rho(d)A\rho(d)
   =\frac{\varepsilon |D|}{N}B^{-1}A^tB.
\] Multiplying by \(A\) and taking traces, we obtain \[
 \left|\sum_{d\in D}\chi((da)^2)\right|
   =|\varepsilon\lambda|\,|D|=|D|.
\] This proves the bound for \(S_2\). The bijection \(g\mapsto g^3\) takes \(Da\) to \(Da^3\), and the same calculation for the normalizing operator \(a^3\) gives the bound for \(S_6\).

For \(g\in Da\), the action of \(g^4\) on \(D/Z(D)\) coincides with the nontrivial action of \(a^4\). The matrices \(\rho(d)\), where \(d\) ranges over a set of representatives for \(D/Z(D)\), form a basis of \(\operatorname{End}(W)\): they are pairwise orthogonal for the trace form, and their number is \(N^2\). Conjugation by \(\rho(g^4)\) permutes this basis with scalar factors \(\pm1\). Thus \[
 |\chi(g^4)|^2
 =\operatorname{tr}(\operatorname{Ad}\rho(g^4))
 \le |C_{D/Z(D)}(a^4)|\le N^2/2.
\] Here the first equality also follows directly from \(\operatorname{End}(W)\cong W\otimes W^*\). Hence \(|S_4|\le |D|N/\sqrt2\).

Finally, every \(g\in Da\) has order \(8\), so the projection onto its fixed-point space is \(\frac18\sum_{j=0}^7\rho(g)^j\). Consequently, \[
 \sum_{g\in Da}\dim C_W(g)
 =\frac18\sum_{j=0}^7S_j
 \ge\frac{|D|}{8}\left(N-2-\frac{N}{\sqrt2}\right)>0,
\] since \(N\ge8\). Although the intermediate sums may be complex, their total is real; the lower bound follows from the triangle inequality. Thus at least one fixed-point space is nonzero. Returning to characteristic \(3\) preserves this conclusion. \end{proof}

\begin{lemma}[The form of Thompson's theorem used here] \label{lem:critical-two-external} Every finite \(2\)-group \(F\) has a characteristic subgroup \(C\) of class at most \(2\) such that \(C/Z(C)\) is elementary abelian and every nonidentity automorphism of \(F\) of odd order acts nontrivially on \(C\). \end{lemma} \begin{proof} This is Thompson's critical subgroup theorem in the form stated in \cite{ThompsonCritical}, Lemma~2.3, with reference to D.~Gorenstein, \emph{Finite Groups}, Chapter 5, Theorem 3.11. The subsequent application uses precisely the properties listed above. \end{proof}

\begin{lemma}\label{lem:eight-no-nonabelian-two} Under the assumptions of Lemma~\ref{lem:nonabelian-two-fitting}, the group \(F=O_2(Q)\) cannot be nonabelian. \end{lemma} \begin{proof} Suppose otherwise. Since \(F\le E\), the group \(E\) does not centralize \(F\). Its image in \(\Aut(F)\) is nontrivial and perfect, so it is not a \(2\)-group. It contains a nonidentity automorphism of odd order. Take a characteristic subgroup \(C\) of \(F\) as in Lemma~\ref{lem:critical-two-external}. Then \(C\trianglelefteq Y\), and the action of \(E\) on \(C\) is nontrivial.

By Lemma~\ref{lem:eight-fixes-types}, all simple types in the restriction of \(V\) to \(C\) are invariant under \(E\). Denote their kernels by \(K_i\). They are normal in \(E\), and \(\bigcap_iK_i=1\) by faithfulness. If \(E\) centralized every \(C/K_i\), it would centralize \(C\). Hence there is a type with kernel \(K\) such that the action of \(E\) on \(D=C/K\) is nontrivial. Working in \(E/K\), we obtain a normal \(2\)-group \(D\) of class at most \(2\), with elementary abelian \(D/Z(D)\) and a faithful invariant simple \(D\)-module type.

The center of \(D\) is cyclic, since it consists of scalars in a faithful irreducible representation. The perfect group \(E/K\) centralizes it, because its automorphism group is abelian. If \(|Z(D)|\ge4\), a central element of order \(4\) would commute with an element of order \(7\) in \(E/K\). Such an element exists, since \(K\le Q\) and \((E/K)/(Q/K)\cong S\). This would give order \(28\), forbidden even in a quotient of \(Y\). Thus \(Z(D)=C_2\). The group \(D\) is nonabelian: otherwise it would be cyclic and the perfect group \(E/K\) would act trivially on it. Consequently, \(D'=Z(D)=\Phi(D)=C_2\), so \(D\) is extraspecial.

The action of \(E/K\) on \(D/Z(D)\) is nontrivial. Indeed, the kernel of the map \(\Aut(D)\to\Aut(D/Z(D))\) is a \(2\)-group: an automorphism in this kernel has the form \(d\mapsto d f(dZ(D))\), where \(f:D/Z(D)\to Z(D)\) is a homomorphism. A perfect group has no nontrivial image in such a group. Moreover, the kernel of the action of \(E\) on \(D/Z(D)\) is contained in \(Q\). Otherwise its product with \(Q\) would be all of \(E\), and the image of \(E\) would coincide with the solvable image of \(Q\), which is impossible for a nontrivial perfect image.

Choose \(a\in E\) whose image in \(S\) has order \(8\). Then \(|a|=8\), and its action on \(D/Z(D)\) has order \(8\), since the kernel of this action lies in \(Q\). Denoting its image in \(E/K\) by the same letter gives \(H=D\rtimes\langle a\rangle\): the intersection \(D\cap\langle a\rangle\) is trivial, since elements of \(D\) act trivially on \(D/Z(D)\). All elements of \(Da\) have order \(8\). Their orders are divisible by \(8\), and the spectrum of \(H\), as the spectrum of a section of \(Y\), is contained in \(\Om\).

Take the isotypic space of the chosen \(C\)-module type in \(V\) over an algebraic closure. It is invariant under \(E\), and \(K\) acts trivially on it. Restrict it to \(H\) and choose a simple constituent \(W\). Its restriction to \(D\) is irreducible and faithful. For completeness, the chosen faithful simple \(D\)-module extends to \(D\rtimes C_8\): an operator realizing \(a\) has scalar eighth power, which can be made the identity by scalar multiplication. All simple modules lying above this type are obtained by tensoring the extension with a linear character of \(C_8\). This gives the required irreducibility of the restriction of \(W\) to \(D\).

Lemma~\ref{lem:extraspecial-eight-coset} gives an element of order \(8\) with a nonzero fixed vector in \(W\). Then \(24\in\omega(W\rtimes H)\). Yet \(W\rtimes H\) is a section of the subgroup \(V\rtimes Y\): one first takes the indicated invariant space, then the quotient by \(K\), the restriction to \(H\), and a simple constituent. Thus \(\omega(W\rtimes H)\subseteq\Om\), a contradiction. The nonabelian case is excluded. \end{proof}

\Needspace{9\baselineskip} \begin{theorem}\label{thm:U33-J2-arbitrary-radical} There is no finite group \(G\) with spectrum \(\omega(G)=\Om\) for which the socle of \(G/R(G)\) is isomorphic to \(U_3(3)\) or \(J_2\). The solvable radical is arbitrary. \end{theorem} \begin{proof} By Lemma~\ref{lem:twentyone-three-witness}, the order \(21\) gives a chief \(3\)-factor with jointly nontrivial actions. For its faithful action group, Lemma~\ref{lem:nonabelian-two-fitting} requires \(F=O_2(R(Y^{(\infty)}))\) to be nonabelian. Lemma~\ref{lem:eight-no-nonabelian-two} excludes such a group. This contradiction uses neither nilpotence of the original radical nor criticality of \(G\). \end{proof}

Thus the socles \(U_3(3),J_2\) have been completely excluded. For \(A_8\), Lemma~\ref{lem:extraspecial-eight-coset} does not apply: \(A_8\) itself has no element of order \(8\), needed for the fourth power to act nontrivially on \(D/Z(D)\). The next subsection uses a different argument, involving order \(18\).

\clearpage \subsection{Complete exclusion of the socle \(A_8\)}

Unlike in the preceding subsection, here we choose the highest factor on which \(P=R(G)\cap G^{(\infty)}\) acts nontrivially, rather than merely a factor with jointly nontrivial actions of \(R(G)\) and \(G^{(\infty)}\). This ensures that the nilpotent kernel \(Q\) of the action group is nontrivial. In characteristic \(2\), such a kernel is eliminated by the Schur multiplier. In characteristic \(3\), the forbidden order \(18\) is decisive.

\subsubsection{The norm of an element of order \(3\) and a four-group}

\begin{lemma}\label{lem:norm-three-four-group} Let \(H\) act on a space \(V\) of characteristic \(3\), let \(a\in H\) have order \(3\), and let \(B\le C_H(a)\) be elementary abelian of order \(4\). If \[
 (1+a+a^2)V\ne0,
\] then \(18\in\omega(V\rtimes H)\). In particular, if \(18\notin\omega(V\rtimes H)\), the operator \(1+a+a^2\) must be zero. \end{lemma} \begin{proof} The space \(U=(1+a+a^2)V\ne0\) is invariant under \(B\). Restriction to \(B\) is completely reducible; all its simple types are linear characters with values \(\pm1\), already defined over \(\mathbb F_3\). The kernel of every such character contains a nonidentity involution. Thus there exist \(t\in B\setminus\{1\}\) and \(0\ne w\in U\) such that \(tw=w\). Choose \(v\in V\) with \((1+a+a^2)v=w\). Since \([a,t]=1\), the element \(at\) has order \(6\) and \[
 \bigl(1+at+\cdots+(at)^5\bigr)v
 =(1+t)(1+a+a^2)v=2w\ne0.
\] Consequently, \((v,at)^6=(2w,1)\ne1\) and \((v,at)^{18}=1\). Its image has order \(6\), so the element itself has order \(18\). \end{proof}

\subsubsection{Actions of \(A_8\) in characteristic \(2\)}

\begin{lemma}\label{lem:A8-three-square-fixed} Let \(Q\trianglelefteq E\) be a \(2\)-group, let \(E/Q\cong A_8\), and let \(M\) be a \(\mathbb F_2E\)-module of dimension at least \(2\). Then there exists an element \(a\in E\) of order \(3\) with nontrivial image in \(A_8\) such that \(\dim C_M(a)\ge2\). \end{lemma} \begin{proof} In \(A_8\), take \(\overline A=\langle c,d\rangle\cong C_3^2\), where \(c=(1\,2\,3)\), \(d=(4\,5\,6)\). The full preimage of this group has a normal Hall \(2\)-subgroup \(Q\). By the Schur--Zassenhaus theorem, it has a complement \(A\cong C_3^2\) mapping isomorphically onto \(\overline A\).

Partition the four subgroups of order \(3\) in \(A\) into pairs according to their images: \[
 \langle c\rangle,\langle d\rangle;
 \qquad \langle cd\rangle,\langle cd^{-1}\rangle.
\] Within either pair, generators are conjugate in \(A_8\). This follows directly from cycle structures: the first two have type \(3\,1^5\) and the last two type \(3^2\,1^2\); neither \(S_8\)-class splits in \(A_8\). Their lifts of order \(3\) are also conjugate in \(E\): after lifting a conjugating element, one obtains two complements to \(Q\) in the preimage of one cyclic group of order \(3\), and they are conjugate by an element of \(Q\).

The restriction of \(M\) to \(A\) is completely reducible. Over \(\mathbb F_2\), the nontrivial simple \(A\)-modules have dimension \(2\); their four types are determined by four kernels of order \(3\). Hence \[
 M=M_0\oplus M_1\oplus M_2\oplus M_3\oplus M_4,
\] where \(M_0=C_M(A)\), and \(M_i\) is the sum of the nontrivial types with kernel equal to the \(i\)-th subgroup. Put \(h=\dim M_0\), \(\dim M_i=2k_i\). For its generator \(a_i\), we have \(C_M(a_i)=M_0\oplus M_i\). Conjugacy within pairs gives \(k_1=k_2=u\), \(k_3=k_4=v\). Consequently, \[
 \dim M=h+4(u+v),\qquad
 \dim C_M(a_1)=h+2u,\quad
 \dim C_M(a_3)=h+2v.
\] If both latter quantities are at most \(1\), then \(u=v=0\), \(h\le1\), contrary to \(\dim M\ge2\). This proves the lemma. \end{proof}

\begin{lemma}\label{lem:A8-orthogonal-three} Let \(E\) be perfect, let \(Q\trianglelefteq E\) be a \(2\)-group, and let \(E/Q\cong A_8\). Suppose that \(E\) acts nontrivially by isometries on a nondegenerate quadratic space \(M\) over \(\mathbb F_2\). Then there exists \(a\in E\) of order \(3\) with nontrivial image in \(A_8\) such that \[
 \dim[M,a]\ge4,
 \qquad C_M(a)\text{ contains a nonzero singular vector}.
\] \end{lemma} \begin{proof} The image of \(E\) in the isometry group is nontrivial and perfect. If the kernel of this action were not contained in \(Q\), its product with \(Q\) would be all of \(E\), and the image of \(E\) would be the solvable image of \(Q\). This is impossible. Hence the kernel is contained in \(Q\), and the image has a composition factor \(A_8\); in particular, its order is divisible by \(7\).

Take the group \(A\) and the decomposition into its isotypic summands from the preceding lemma. All four nontrivial simple types are self-dual. Distinct summands are orthogonal, and each is nondegenerate: this follows from invariance of the polar form and semisimplicity of the restriction to the odd-order group \(A\). Thus \[
 \dim M_0=2h,\quad
 \dim M_1=\dim M_2=2u,\quad
 \dim M_3=\dim M_4=2v.
\] Here \(u+v\ge1\), since otherwise \(A\) would act trivially, contrary to the containment of the kernel in \(Q\). For generators of the subgroups in the two pairs, we get \[
\begin{array}{c|cc}
 &\dim C_M(a_i)&\dim[M,a_i]\\ \hline
 i=1,2&2(h+u)&2(u+2v)\\
 i=3,4&2(h+v)&2(2u+v).
\end{array}
\] A nondegenerate quadratic space over \(\mathbb F_2\) of dimension at least \(4\) contains a nonzero singular vector. For example, the number of zeros of the form in dimension \(2m\) is \(2^{2m-1}\pm2^{m-1}\), and for \(m\ge2\) it exceeds \(1\).

If \(u,v\ge1\), the moved space in each row has dimension at least \(6\). If \(h\ge1\) or \(\max(u,v)\ge2\), the fixed-point space in one row has dimension at least \(4\), proving the assertion. The only remaining possibility is \(h=0,u=v=1\), \(\dim M=8\).

If, for example, \(v=0\), \(u\ge2\), the first row immediately gives both required dimensions at least \(4\). If \(v=0,u=1\), the second row has moved-space dimension \(4\) and fixed-point space \(M_0\). If \(h\ge2\), the assertion is proved. If \(h=0\), we obtain \(\dim M=4\), which is impossible: the orders of \(O_4^+(2)\) and \(O_4^-(2)\) are \(72\) and \(120\) and are not divisible by \(7\). If \(h=1\), the only remaining possibility is that \(M_0\) is an anisotropic plane. Each \(M_1,M_2\) is also an anisotropic plane, since its isometry group contains an element of order \(3\), whereas \(|O_2^+(2)|=2\). Thus \(M\) has dimension \(6\) and minus type. But \[
 |O_6^-(2)|=2\cdot2^6(2^3+1)(2^2-1)(2^4-1)
 =51840
\] is not divisible by \(7\). This is a contradiction. The case \(u=0\) is symmetric.

It remains to eliminate \(h=0,u=v=1\), \(\dim M=8\). The normal \(2\)-group \(Q\) acts trivially on every simple \(\mathbb F_2E\)-module: its nonzero fixed-point subspace for \(Q\) is invariant under \(E\). Hence the composition factors of \(M\) are simple \(A_8\)-modules. Their complete dimension list is \(1,4,4,6,14,20,20,64\); see \cite{TwoLocalData}. The two modules of dimension \(4\) are the natural \(L_4(2)\cong A_8\)-module and its dual; the module of dimension \(6\) is its exterior square.

The restriction of the natural \(4\)-dimensional module to \(A=C_3^2\) consists of two distinct nontrivial \(2\)-dimensional types: this follows from the embedding of a Sylow \(3\)-subgroup in \(\operatorname{GL}_2(2)\times\operatorname{GL}_2(2)\). The dual module has the same two types, since each is self-dual. The exterior square has two trivial summands, obtained as the exterior squares of these \(2\)-dimensional modules. Since \(C_M(A)=0\), composition factors of dimensions \(1\) and \(6\) are impossible: the fixed-point functor for \(A\) is exact, because \(|A|\) is odd. Thus \(M\) has exactly two \(4\)-dimensional factors. In their sum, each of the two indicated types occurs twice, while the other two do not occur. This contradicts \(\dim M_i=2\) for all four \(i\). This last exclusion completes the proof. \end{proof}

\subsubsection{An extraspecial \(2\)-group in characteristic \(3\)}

\begin{lemma}\label{lem:extraspecial-three-norm} Let \(D\) be an extraspecial \(2\)-group, let \(a\) be an automorphism of it of order \(3\), and let \(H=D\rtimes\langle a\rangle\). Let \(W\) be a simple \(H\)-module over an algebraically closed field of characteristic \(3\), faithful on \(D\). If \(\dim[D/Z(D),a]\ge4\), then \[
 (1+a+a^2)W\ne0.
\] If, moreover, the fixed-point space on \(D/Z(D)\) contains a nonzero singular vector for the form defined by squares, then \(18\in\omega(W\rtimes H)\). \end{lemma} \begin{proof} Let \(Z(D)=\langle z\rangle\). On \(M=D/Z(D)\), put \(d^2=z^{q(dZ(D))}\); the polar form of the quadratic form \(q\) is given by the commutator. We have an orthogonal decomposition \(M=C_M(a)\perp[M,a]\). On the moved summand, the automorphism \(a\) has no nonzero fixed vectors. This summand decomposes into an orthogonal sum of \(a\)-invariant anisotropic planes. We give a direct justification. If \(b\) is the polar form and \(b(v,av)=0\) for all \(v\), polarizing this identity gives \[
 0=b(v,aw)+b(w,av)
   =b\bigl(v,(a+a^{-1})w\bigr)=b(v,w),
\] since \(a+a^{-1}=1\) on \([M,a]\). This contradicts nondegeneracy. Thus there exists \(v\) with \(b(v,av)=1\). The plane \(\langle v,av\rangle\) is nondegenerate, invariant, and anisotropic: \(a^2v=v+av\) and invariance of \(q\) imply \(q(v)=q(av)=q(v+av)=1\). Its orthogonal complement is invariant, and the argument repeats.

Let \(k\) be the number of these planes. By hypothesis, \(k\ge2\). Their full preimages in \(D\) are isomorphic to \(Q_8\), commute pairwise, and have common intersection \(\langle z\rangle\). Together with the preimage of \(C_M(a)\), they form \(D\) as a central product. The automorphism \(a\) centralizes the last preimage: it acts trivially on its quotient by \(\langle z\rangle\), and the kernel of the corresponding automorphism map is a \(2\)-group. On each copy of \(Q_8\), the automorphism \(a\) acts nontrivially.

The faithful simple \(D\)-module is unique and is the tensor product of the faithful \(2\)-dimensional modules for these groups \(Q_8\) with a faithful module for the fixed central factor, if the latter is larger than \(\langle z\rangle\). Uniqueness and dimension of this module were justified in Lemma~\ref{lem:extraspecial-eight-coset}. On each \(2\)-dimensional factor, the operator \(a\) is a single block \(J_2(1)\): it has order \(3\) and cannot be the identity, since it would then centralize the faithful image of \(Q_8\). On the fixed factor it is scalar, and a scalar of order dividing \(3\) in characteristic \(3\) is \(1\).

We explain why this describes the simple \(H\)-module itself. The faithful simple \(D\)-type is invariant and extends to \(D\rtimes C_3\): an operator realizing the automorphism has scalar cube, which can be made the identity. On any isotypic module, the operator \(a\) is the tensor product of this extension with an operator \(B\), \(B^3=1\), on the multiplicity space. The operator \(B\) has a fixed line. Simplicity of \(W\) forces the multiplicity space to be one-dimensional. Thus \(W\) restricts irreducibly to \(D\), and the indicated tensor product of blocks \(J_2(1)\) acts on it.

If \(au=u\), \(av=v+u\) on each of the first two factors, then \[
 (a-1)^2(v\otimes v)=2u\otimes u\ne0.
\] On the other factors, choose nonzero fixed vectors. We obtain \((a-1)^2W\ne0\), that is, \((1+a+a^2)W\ne0\).

Finally, a fixed singular vector \(xZ(D)\ne0\) lifts to an involution \(t\in D\setminus Z(D)\) centralized by \(a\). Indeed, both its preimages are involutions, and the permutation of these two elements induced by an automorphism of order \(3\) is trivial. Then \(\langle t,z\rangle\cong C_2^2\le C_H(a)\), and Lemma~\ref{lem:norm-three-four-group} applies. \end{proof}

\begin{lemma}\label{lem:A8-two-centralizer-four} Let \(E\) be perfect, let \(Q\trianglelefteq E\) be a \(2\)-group, and let \(E/Q\cong A_8\), \(\omega(E)\subseteq\Om\), \([E,Q]\ne1\). Then there exists an element \(a\in E\) of order \(3\) with nontrivial image in \(A_8\) such that \(C_Q(a)\) contains a four-group. \end{lemma} \begin{proof} The image of \(E\) in \(\Aut(Q)\) is nontrivial and perfect. It is therefore not a \(2\)-group and contains a nonidentity automorphism of odd order. Take a characteristic critical subgroup \(C\) of \(Q\) as in Lemma~\ref{lem:critical-two-external}. The group \(E\) acts nontrivially on \(C\).

If \(\Omega_1(Z(C))\) has dimension at least \(2\), Lemma~\ref{lem:A8-three-square-fixed} gives the required four-group in its fixed-point space. Otherwise \(Z(C)\) is cyclic. The perfect group \(E\) centralizes it, since the automorphism group of a cyclic group is abelian. If \(|Z(C)|\ge4\), a central element of order \(4\) would commute with an element of order \(7\) in \(E\), giving \(28\notin\Om\). Hence \(Z(C)=C_2\). The group \(C\) cannot be abelian, since otherwise \(C=Z(C)=C_2\) and the action would be trivial. Class at most \(2\) and elementary abelianness of \(C/Z(C)\) imply \(C'=\Phi(C)=Z(C)=C_2\). Thus \(C\) is extraspecial.

The action of \(E\) on \(C/Z(C)\) is nontrivial: the kernel of the automorphism map to this quotient is a \(2\)-group, in which a perfect group has no nontrivial image. The quadratic form defined by squares of elements of \(C\) is invariant. Lemma~\ref{lem:A8-orthogonal-three} gives an element \(a\) of order \(3\) and a nonzero fixed singular vector. As in the preceding lemma, this vector lifts to an involution \(t\notin Z(C)\) centralized by \(a\). The group \(\langle t,Z(C)\rangle\) is the required four-group. \end{proof}

\begin{lemma}\label{lem:A8-no-noncentral-two-on-three} Let \(E\) be perfect, let \(Q\trianglelefteq E\) be a \(2\)-group, and let \(E/Q\cong A_8\). If \(E\) acts faithfully on a nonzero space \(V\) of characteristic \(3\) and \(\omega(V\rtimes E)\subseteq\Om\), then \([E,Q]=1\). \end{lemma} \begin{proof} Suppose that \([E,Q]\ne1\). The preceding lemma gives \(a\) of order \(3\) and a four-group \(B\le C_Q(a)\), with the image of \(a\) in \(A_8\) nontrivial.

Let \(T\trianglelefteq E\) be any subgroup of \(Q\). Extend the field to an algebraic closure. The group \(E\) permutes the simple \(T\)-types in \(V\). If its permutation action is nontrivial, the kernel is contained in \(Q\): otherwise the image of \(E\) would equal the solvable image of \(Q\) while also being perfect. Hence \(a\) acts nontrivially on the set of types and has an orbit of length \(3\). On the corresponding three isotypic spaces, the operator \(1+a+a^2\) is nonzero. By Lemma~\ref{lem:norm-three-four-group}, order \(18\) would occur. Thus \(E\) fixes every simple type of any such \(T\) that occurs.

Take the critical subgroup \(C\) from the proof of the preceding lemma. The group \(E\) acts on it nontrivially and fixes all simple \(C\)-types in \(V\). Their kernels \(K_i\) are normal in \(E\), and \(\bigcap K_i=1\) by faithfulness. Thus, for some \(K=K_i\), the action of \(E\) on \(D=C/K\) is nontrivial. In \(E/K\), the subgroup \(D\) is normal and has a faithful invariant simple type. Its center is cyclic, since it consists of scalars in a faithful irreducible representation. The perfect group \(E/K\) centralizes \(Z(D)\). The order \(28\) again excludes \(|Z(D)|\ge4\). If \(D\) were abelian, it would be cyclic and centralized by \(E/K\), contrary to its choice. Thus \(D\) is extraspecial with center \(C_2\), as in the preceding lemma.

The group \(E/K\) acts nontrivially on \(D/Z(D)\): the kernel of the automorphism map to this quotient is a \(2\)-group of central automorphisms. Apply Lemma~\ref{lem:A8-orthogonal-three} to \(E/K\), its normal \(2\)-subgroup \(Q/K\), and the quadratic space \(D/Z(D)\). We obtain an element \(b\) of order \(3\), a moved space of dimension at least \(4\), and a nonzero fixed singular vector. Put \(H=D\rtimes\langle b\rangle\le E/K\).

The chosen \(C\)-isotypic space in \(V\) is invariant under \(E\), and \(K\) acts trivially on it. Restrict it to \(H\) and take a simple composition factor \(W\). Its restriction to \(D\) is completely reducible and consists of the chosen faithful type, so \(W\) is faithful on \(D\). By Lemma~\ref{lem:extraspecial-three-norm}, \(18\in\omega(W\rtimes H)\). Yet \(W\rtimes H\) is a section of the subgroup \(V\rtimes E\): one takes the indicated invariant space, the quotient by \(K\), the restriction to \(H\), and then a subquotient of the module. This contradicts spectral inclusion. Consequently, \([E,Q]=1\). \end{proof}

\subsubsection{A central cover and an internal joint action}

\begin{lemma}\label{lem:double-A8-three} There is no faithful nonzero module \(V\) for the group \(2.A_8\) in characteristic \(3\) satisfying \(\omega(V\rtimes2.A_8)\subseteq\Om\). \end{lemma} \begin{proof} The group \(2.A_8\) contains an actual subgroup \(H_0=2^3{:}7{:}3\cong\operatorname{A\Gamma L}_1(8)\). This is the stabilizer in the transitive action of degree \(240\) listed in \cite{TwoLocalData}; its order is \(40320/240=168\). Unlike the subgroup \(L_2(7)\) in Lemma~\ref{lem:double-A8-seven}, a different stabilizer of the same degree is used here. The central involution \(z\) of \(2.A_8\) does not belong to \(H_0\), since \(Z(H_0)=1\). Hence there is a subgroup \(\langle z\rangle\times H_0\).

Denote the normal \(2^3\) in \(H_0\) by \(T\), its complement of order \(7\) by \(\langle c\rangle\), and the field automorphism of order \(3\) by \(a\). On \(T\cong(\mathbb F_8,+)\), the element \(c\) acts by multiplication by a primitive element, and \(a\) by squaring. In particular, \(|C_T(a)|=2\), and \(\langle z\rangle\times C_T(a)\) is a four-group centralizing \(a\).

Faithfulness of \(V\) means that \(T\) acts nontrivially. The element \(c\) acts transitively on its nontrivial linear characters, so all seven nontrivial \(T\)-types occur in \(V\). The element \(a\) has two orbits of length \(3\) and one fixed point on them: its fixed-point space on the dual three-dimensional space is one-dimensional. Thus \((1+a+a^2)V\ne0\), since the three isotypic spaces in one such orbit are cyclically permuted. Lemma~\ref{lem:norm-three-four-group} gives \(18\in\omega(V\rtimes2.A_8)\), a contradiction. \end{proof}

\begin{lemma}\label{lem:highest-internal-A8} Let \(\omega(G)=\Om\), \(\Soc(G/R(G))\cong A_8\), \(L=G^{(\infty)}\), \(P=R(G)\cap L\). Then there is a chief factor \(V=N/M\) inside \(P\) such that, in the notation \[
 Y=G/C_G(V),\quad E=Y^{(\infty)},\quad Q=R(E)
\], we have \[
 \operatorname{char}V=3,\qquad
 Q\text{ is a nontrivial }2\text{-group},\qquad
 E/Q\cong A_8.
\] Moreover, \(V\) is a faithful simple \(Y\)-module, and \(\omega(V\rtimes Y)\subseteq\Om\). \end{lemma} \begin{proof} By Corollary~\ref{cor:perfect-radical-seven}, the group \(P\) is nonnilpotent. Thus every chief series of \(G\) inside \(P\) has a factor on which \(P\) acts nontrivially: otherwise this series would be a central series of \(P\). Choose the highest such factor \(V=N/M\). Above \(N\), the action of \(P\) on every factor is trivial, so \(P/N\) is nilpotent.

Since \(P\le R(G)\cap L\), both groups \(R(G),L\) act nontrivially on \(V\), and Lemma~\ref{lem:mixed-chief-action} applies. As in Lemma~\ref{lem:highest-mixed}, \(E\) is the image of \(L\) and \(Q\) the image of \(P\), with \(E/Q\cong A_8\). Now, however, we also have \(Q\ne1\), since \(P\) acts nontrivially on \(V\). Since \(N\le C_G(V)\), the group \(Q\) is a quotient of the nilpotent group \(P/N\). Moreover, \(O_p(Y)=1\), where \(p=\operatorname{char}V\). Hence the nilpotent normal subgroup \(Q\) has order coprime to \(p\).

By Theorem~\ref{thm:A8-radical-two-three}, \(\pi(P)\subseteq\{2,3\}\). Thus, if \(p=2\), the group \(Q\) is a nontrivial \(3\)-group, and if \(p=3\), it is a nontrivial \(2\)-group. Exclude the first possibility. Applying Lemma~\ref{lem:split-spectrum} to \(Q/\Phi(Q)\trianglelefteq E/\Phi(Q)\) gives \[
 \omega\bigl((Q/\Phi(Q))\rtimes(E/Q)\bigr)
 \subseteq\omega(E)\subseteq\Om.
\] By Theorem 1 and Table 1 in \cite{ZhLM}, the only permissible simple \(A_8\)-module in characteristic \(3\) is trivial. Hence all composition factors of \(Q/\Phi(Q)\) are trivial on \(E\). Lemma~\ref{lem:perfect-automorphism-tools} shows that the perfect group \(E\) centralizes \(Q\). Then \(E\) is a perfect central cover of \(A_8\), and \(|Q|\) divides \(|M(A_8)|=2\). This is impossible for a nontrivial \(3\)-group. Thus \(p=3\) remains, and \(Q\) is a nontrivial \(2\)-group, as required. \end{proof}

\Needspace{9\baselineskip} \begin{theorem}\label{thm:A8-arbitrary-radical} There is no finite group \(G\) with spectrum \(\omega(G)=\Om\) for which \(\Soc(G/R(G))\cong A_8\). Neither nilpotence of the radical nor criticality of the group is assumed. \end{theorem} \begin{proof} Take \(V,Y,E,Q\) from Lemma~\ref{lem:highest-internal-A8}. The group \(E\) acts faithfully on \(V\), \(\operatorname{char}V=3\), \(\omega(V\rtimes E)\subseteq\Om\), and \(Q\ne1\) is a normal \(2\)-group with \(E/Q\cong A_8\). Lemma~\ref{lem:A8-no-noncentral-two-on-three} gives \([E,Q]=1\). The property of perfect central covers and \(M(A_8)=C_2\) yield \(Q=C_2\), \(E\cong2.A_8\). This contradicts Lemma~\ref{lem:double-A8-three}. Thus both a noncentral \(2\)-kernel and the final central cover have been excluded. \end{proof}

\section{\texorpdfstring{The socle \(A_7\)}{The socle A7}}
\subsection{The socle \(A_7\): excluding noncentral kernels}

\subsubsection{Adapting the four-group argument}

\begin{lemma}\label{lem:A7-two-kernel-three}
Let \(E\) be perfect, let \(Q\trianglelefteq E\) be a \(2\)-group, and suppose that \(E/Q\cong A_7\). Then the following statements hold.
\begin{enumerate}
\item For every \(\mathbb F_2E\)-module \(M\) of dimension at least \(2\), there is an element \(a\in E\) of order \(3\) with nontrivial image in \(A_7\) such that \(\dim C_M(a)\ge2\).
\item If \(E\) acts nontrivially on a nondegenerate quadratic space \(M\) over \(\mathbb F_2\), there is such an \(a\) for which \(\dim[M,a]\ge4\) and \(C_M(a)\) contains a nonzero singular vector.
\item If \(\omega(E)\subseteq\Om\) and \([E,Q]\ne1\), then, for some such \(a\), the group \(C_Q(a)\) contains a four-group.
\item If \(E\) acts faithfully on a space \(V\) of characteristic \(3\) and \(\omega(V\rtimes E)\subseteq\Om\), then \([E,Q]=1\).
\end{enumerate}
\end{lemma}
\begin{proof}
We describe the changes needed to adapt the proofs of Lemmas \ref{lem:A8-three-square-fixed}, \ref{lem:A8-orthogonal-three}, \ref{lem:A8-two-centralizer-four}, and \ref{lem:A8-no-noncentral-two-on-three}.

In \(A_7\), take the same permutations \(c=(1\,2\,3)\), \(d=(4\,5\,6)\). A Sylow \(3\)-subgroup has order \(9\), and its preimage over a normal \(2\)-group splits by Schur--Zassenhaus. The subgroups \(\langle c\rangle,\langle d\rangle\) are conjugate, as are \(\langle cd\rangle,\langle cd^{-1}\rangle\). The \(S_7\)-conjugacy classes of cycle types \(3\,1^4\) and \(3^2\,1\) do not split into separate \(A_7\)-conjugacy classes. The lifts of order \(3\) are conjugate by the same complement argument. Consequently, the isotypic multiplicities on \(C_3^2\) still satisfy \(k_1=k_2=u\), \(k_3=k_4=v\). The formulas from Lemma~\ref{lem:A8-three-square-fixed},
\[
 \dim M=h+4(u+v),\qquad
 \dim C_M(a_1)=h+2u,\quad
 \dim C_M(a_3)=h+2v,
\]
prove statement 1 without change.

For statement 2, the image of the isometry group has a composition factor \(A_7\), so its order is divisible by \(7\). The quadratic counting argument in Lemma~\ref{lem:A8-orthogonal-three} used only this property, the conjugacy of the indicated pairs, and the description of the modules of dimension at most \(8\). For \(A_7\), the latter description is the same: \(1,4,4,6\); the other simple modules over \(\mathbb F_2\) have dimensions \(14,20\), see \cite{A7LocalData}. The two \(4\)-dimensional modules are restrictions of the natural and dual modules of \(L_4(2)\), and the \(6\)-dimensional module is their exterior square. On \(C_3^2\), the two \(4\)-dimensional modules have the same pairs of nontrivial types; the exterior square has a two-dimensional subspace fixed by all of \(C_3^2\). Thus the exclusion in dimension \(8\) from that lemma also remains valid. The exceptional case of dimension \(4\) and that of minus type in dimension \(6\) are eliminated because the orders of the corresponding orthogonal groups are not divisible by \(7\). This proves statement 2.

For statement 3, take a critical subgroup \(C\) of \(Q\). If \(\dim\Omega_1(Z(C))\ge2\), apply statement 1. If the center is cyclic, the perfect group \(E\) centralizes it; order \(28\) restricts the center to \(C_2\). Then \(C\) is extraspecial, and statement 2 gives a fixed noncentral involution which, together with the center, generates a four-group. All these steps are identical to the proof of Lemma~\ref{lem:A8-two-centralizer-four}.

Finally, in statement 4, if \([E,Q]\ne1\), statement 3 gives an element \(a\) and a four-group in \(C_Q(a)\). A nontrivial permutation action of \(E\) on the simple types of a normal subgroup \(T\le Q\) has kernel contained in \(Q\), since \(E\) is perfect and \(E/Q\) is simple. Thus \(a\) would have an orbit of length \(3\), which Lemma~\ref{lem:norm-three-four-group} excludes by order \(18\). All types of the critical subgroup \(C\) are invariant. Choosing a faithful quotient \(C/K\) with nontrivial \(E\)-action gives an extraspecial group with center \(C_2\), as in Lemma~\ref{lem:A8-no-noncentral-two-on-three}. Statement 2 and Lemma~\ref{lem:extraspecial-three-norm} again give order \(18\). This contradiction proves statement 4.
\end{proof}

\subsubsection{Odd characteristics and central covers}

\begin{lemma}\label{lem:A7-cover-prime-exclusions}
The following statements hold.
\begin{enumerate}
\item The group \(6.A_7\) contains an element of order \(24\).
\item Neither \(A_7\) nor \(2.A_7\) has a nonzero module \(V\) of characteristic \(5\) with \(\omega(V\rtimes D)\subseteq\Om\).
\item Neither \(A_7\) nor \(3.A_7\) has a nonzero module \(V\) of characteristic \(7\) satisfying the same inclusion.
\end{enumerate}
\end{lemma}
\begin{proof}
In \(2.A_7\), the normalizer of a Sylow \(5\)-subgroup has structure \(5{:}8\), see \cite{A7LocalData}. In particular, there is an element of order \(8\). The group \(6.A_7\) has a central subgroup of order \(3\) with quotient \(2.A_7\). The \(2\)-part of a lift of an element of order \(8\) still has order \(8\) and commutes with the central element of order \(3\). Their product has order \(24\).

For \(A_7\) in characteristic \(5\), Theorem 1 and Table 1 of \cite{ZhLM} give no admissible nontrivial simple modules. The trivial module is also impossible: with an element of order \(7\), it gives order \(35\). Now consider a simple module for \(2.A_7\). If the center acts trivially on it, apply the preceding conclusion. Otherwise the possible absolutely simple dimensions are
\[
 4,4,14,14,20,20,
\]
by the complete module table in \cite{A7LocalData}. The normalizer of a \(7\)-subgroup contains an element of order \(3\) which permutes its six nontrivial characters in two orbits of length \(3\). Hence a fixed-point-free \(7\)-element would force the dimension to be divisible by \(3\). None of the listed dimensions is divisible by \(3\). Thus there is a nonzero fixed vector, giving order \(35\). For an arbitrary nonzero module, it suffices to take a simple composition factor; the corresponding semidirect product is a section.

The group \(A_7\) contains \(L_2(7)\), acting on the seven nonzero vectors of \(\mathbb F_2^3\). This embedding lies in \(A_7\), since \(L_2(7)\) is perfect. The full preimage of this subgroup in \(3.A_7\) is a central extension with a \(3\)-kernel. Its perfect residual covers \(L_2(7)\), and the order of the kernel of this cover divides both \(3\) and \(|M(L_2(7))|=2\). Thus the kernel is trivial, and \(3.A_7\) contains a subgroup \(L_2(7)\). Restricting any nonzero module of characteristic \(7\) to this subgroup contradicts Lemma~\ref{lem:small-seven-no-modules}.
\end{proof}

\begin{lemma}\label{lem:A7-primes-force-covers}
Let \(\omega(G)=\Om\), \(\Soc(G/R(G))\cong A_7\), \(L=G^{(\infty)}\), and \(P=R(G)\cap L\). Then:
\begin{enumerate}
\item if \(5\mid |R(G)|\), the group \(L\) has a quotient \(3.A_7\), and the composite map onto \(A_7\) has kernel \(P\);
\item if \(7\mid |R(G)|\), the group \(L\) has such a quotient \(2.A_7\);
\item both quotients cannot exist simultaneously. In particular, \(35\nmid |R(G)|\).
\end{enumerate}
\end{lemma}
\begin{proof}
By Lemma~\ref{lem:perfect-residual-product}, the spectrum of \(R/P\) is contained in \(\{1,2,3\}\). Thus, for \(r=5,7\), divisibility \(r\mid |R|\) implies \(r\mid |P|\). Since \(r\nmid |M(A_7)|=6\), Lemma~\ref{lem:central-chief-prime-elimination} shows that \(L\) cannot centralize all chief \(r\)-factors inside \(P\). Choose a factor \(V\) on which \(L\) acts nontrivially.

If \(R\) acts trivially on \(V\), then \(V\) is a \(G/R\)-module, and restriction to \(A_7\) contradicts Lemma~\ref{lem:A7-cover-prime-exclusions}. Thus the factor is mixed. For its faithful action group \(Y\), put \(E=Y^{(\infty)}\), \(Q=R(E)\). By Lemma~\ref{lem:odd-mixed-quasisimple}, \(E\) is a central cover of \(A_7\). Its possibilities are
\[
 A_7,\quad 2.A_7,\quad 3.A_7,\quad 6.A_7.
\]
In characteristic \(5\), the preceding lemma leaves only \(E=3.A_7\), and in characteristic \(7\), only \(E=2.A_7\). The group \(E\) is an image of \(L\), with radical equal to the image of \(P\); this proves the first two statements.

Suppose both epimorphic images \(D_2=2.A_7\), \(D_3=3.A_7\) exist. The image of \(L\) in \(D_2\times D_3\) is perfect, and its image modulo the centers is the graph of an automorphism of \(A_7\), since both maps have the original quotient \(L/P\). The kernel above this diagonal lies in \(C_2\times C_3\) and maps onto each factor: a lift of a central element of one factor has trivial image in \(A_7\), and hence a central image in the other factor. Since the orders are coprime, this kernel is all of \(C_2\times C_3\). We have obtained a perfect central quotient of \(L\) with kernel of order \(6\) over \(A_7\), namely the universal cover \(6.A_7\). It contains an element of order \(24\), contrary to \(\omega(L)\subseteq\Om\).
\end{proof}

\Needspace{24\baselineskip}
\subsubsection{The two remaining central branches}

\begin{theorem}\label{thm:A7-highest-central-pairs}
Let \(\omega(G)=\Om\), and suppose that the simple socle of \(G/R(G)\) is \(A_7\). In a chief series of \(G\) inside \(P=R(G)\cap G^{(\infty)}\), choose the highest factor \(V=N/M\) on which \(P\) acts nontrivially. For \(p=\operatorname{char}V\), \(Y=G/C_G(V)\), \(E=Y^{(\infty)}\), and \(Q=R(E)\), only four rows are possible:
\[
\begin{array}{c|c|c}
 p&E&Q\\ \hline
 2&3.A_7&C_3\\
 3&2.A_7&C_2\\
 5&3.A_7&C_3\\
 7&2.A_7&C_2.
\end{array}
\]
In particular, \(Q\) is always central; all noncentral nilpotent kernels for such a highest factor are excluded.

The group \(L=G^{(\infty)}\) has exactly one of the two central covers \(2.A_7\), \(3.A_7\) as a quotient over \(L/P\). If it is \(2.A_7\), then \(5\nmid |R(G)|\) and \(p\in\{3,7\}\). If it is \(3.A_7\), then \(7\nmid |R(G)|\) and \(p\in\{2,5\}\).
\end{theorem}
\begin{proof}
The nonnilpotence of \(P\) was proved in Corollary~\ref{cor:perfect-radical-seven}. Hence the chosen factor exists. Exactly as in Lemma~\ref{lem:highest-internal-A8}, we obtain \(Q\ne1\), \(E/Q\cong A_7\), and \(Q\) is a nilpotent \(p'\)-group: all factors above \(N\) are centralized by \(P\), and \(Q\) is a quotient of the nilpotent group \(P/N\).

We show that \(Q_5=Q_7=1\), where \(Q_r\) is the Sylow \(r\)-subgroup of \(Q\). Each \(Q_r\) is characteristic in \(Q\). If it is nontrivial, factor \(E\) by the product of the other Sylow subgroups of \(Q\), and then by \(\Phi(Q_r)\). Lemma~\ref{lem:split-spectrum} gives
\[
 \omega\bigl((Q_r/\Phi(Q_r))\rtimes A_7\bigr)
 \subseteq\Om.
\]
For \(r=5,7\), this contradicts Lemma~\ref{lem:A7-cover-prime-exclusions}. Thus \(\pi(Q)\subseteq\{2,3\}\).

For \(p=2\), the group \(Q\) is a nontrivial \(3\)-group. The admissible simple modules of \(A_7\) in characteristic \(3\) are trivial, by Theorem 1 and Table 1 of \cite{ZhLM}. Applying the same spectral inclusion to \(Q/\Phi(Q)\), together with Lemma~\ref{lem:perfect-automorphism-tools}, gives \([E,Q]=1\). Since \(M(A_7)=C_6\) and \(Q\) is nontrivial, we obtain \(Q=C_3\), \(E=3.A_7\).

For \(p=3\), the group \(Q\) is a nontrivial \(2\)-group. Statement 4 of Lemma~\ref{lem:A7-two-kernel-three} gives \([E,Q]=1\). The Schur multiplier then yields \(Q=C_2\), \(E=2.A_7\).

For \(p=5,7\), Lemma~\ref{lem:odd-mixed-quasisimple} makes \(Q\) central immediately. Lemma~\ref{lem:A7-cover-prime-exclusions} leaves \(3.A_7\) and \(2.A_7\), respectively. This proves the table.

Thus at least one cover is a quotient of \(L\). Simultaneous existence of both is forbidden by Lemma~\ref{lem:A7-primes-force-covers}. The same lemma gives the restrictions on the prime divisors of the radical in each branch. The table restricts the characteristic of the highest factor to the two values indicated in each branch.
\end{proof}

\begin{lemma}\label{lem:A7-canonical-cover-quotient}
Under the assumptions of Theorem~\ref{thm:A7-highest-central-pairs}, let \(f:L\twoheadrightarrow D_r=r.A_7\), where \(r\in\{2,3\}\), be the central cover over \(L/P\) obtained there, and let \(K=\ker f\). Then \(K\trianglelefteq G\), \(P/K=Z(D_r)\), and \(R(G)/K\) centralizes \(L/K\). Moreover,
\[
\begin{array}{c|c}
 D_r&R(G)/K\\ \hline
 2.A_7&\text{an elementary abelian }2\text{-group}\\
 3.A_7&\text{a }3\text{-group of exponent }3.
\end{array}
\]
In particular, in the first branch \(G/L\) is a \(2\)-group.
\end{lemma}
\begin{proof}
First we prove that \(K\) is independent of the epimorphism onto \(D_r\) over the fixed quotient \(L/P\). For two such epimorphisms, consider the image \(I\) of \(L\) in \(D_r\times D_r\). It is perfect and is a central cover of \(A_7\): its image modulo the centers is a diagonal, possibly twisted by an automorphism. Its kernel above this diagonal lies in \(C_r^2\). By the property of perfect central covers, this kernel is a quotient of \(M(A_7)=C_6\), and hence has order at most \(r\). On the other hand, both projections \(I\to D_r\) are surjective, so \(|I|\ge r|A_7|\). Thus \(|I|=r|A_7|\), and both projections are isomorphisms. Hence the two original epimorphisms have the same kernel.

Conjugation by an element of \(G\) gives another such epimorphism, since \(L\) is characteristic and \(P\) is normal in \(G\). The independence just proved implies \(K\trianglelefteq G\).

Now work in \(G/K\). By Lemma~\ref{lem:perfect-residual-product}, the elements of \(R/K\) act trivially on \((L/K)/(P/K)\), and \(P/K=Z(D_r)\). Thus Lemma~\ref{lem:perfect-automorphism-tools}, applied to the perfect group \(D_r\), shows that \(R/K\) centralizes \(D_r\). Denote \(R/K\) by \(T\).

For \(D_2\), an element of order \(3\) in \(T\) would commute with an element of order \(8\) in \(D_2\), giving \(24\). A prime divisor \(5\) of \(|T|\), together with an element of order \(7\) in \(D_2\), would give \(35\); a prime divisor \(7\) would give the same order with an element of order \(5\). Hence \(T\) is a \(2\)-group. An element of order \(4\) in \(T\) would commute with a \(7\)-element of \(D_2\), giving \(28\). Thus \(T\) has exponent \(2\), so it is elementary abelian. Finally, \(R/P\) is also a \(2\)-group, and \([G:LR]\le2\). Consequently, \(G/L\) is a \(2\)-group.

For \(D_3\), an involution in \(T\) would commute with an element of order \(21\) in \(D_3\), giving \(42\). The prime divisors \(5,7\) are excluded by the same argument using order \(35\). Hence \(T\) is a \(3\)-group. An element of order \(9\) in \(T\) would commute with an involution of \(D_3\), giving \(18\). Therefore \(T\) has exponent \(3\), since it contains \(P/K=C_3\).
\end{proof}

\subsubsection{Subspectral extensions of central covers}

\begin{lemma}\label{lem:A7-actual-subpectral-examples}
There exist perfect nonsolvable groups
\[
 K_2=3^8{:}2.A_7,\qquad
 K_3=2^{12}{:}3.A_7
\]
with nonnilpotent solvable radicals, simple quotient \(A_7\) by the radical, and \(\omega(K_i)\subseteq\Om\). Furthermore,
\[
 \{20,30\}\cap\omega(K_2)=\varnothing,
 \qquad
 \{9,20,30\}\cap\omega(K_3)=\varnothing.
\]
In particular, neither group has the full spectrum \(\Om\).
\end{lemma}
\begin{proof}
For \(D_2=2.A_7\), take a faithful absolutely simple module of dimension \(4\) over \(\mathbb F_9\); for \(D_3=3.A_7\), take such a module of dimension \(6\) over \(\mathbb F_4\). Both occur in the complete list in \cite{A7LocalData}. Their additive groups have orders \(3^8\) and \(2^{12}\), respectively.

The only involution of \(D_2\) is the central involution \(z\). Indeed, \(D_2\) contains \(2.L_3(2)\cong\operatorname{SL}_2(7)\), see \cite{A7LocalData}. A Sylow \(2\)-subgroup of the latter group is generalized quaternion of order \(16\): it is generated by a torus \(\langle x\rangle\) of order \(8\) and an element \(t\) with \(t^{-1}xt=x^{-1}\), \(t^2=x^4=-I\). It is also Sylow in \(D_2\). Its unique involution is \(z\); conjugacy of Sylow subgroups shows that there are no other involutions. On the faithful simple module, \(z=-I\). Every element \(g\in D_2\) of even order has a power equal to \(z\), so \(g\) has no nonzero fixed vectors. Every coset of the additive module above such a \(g\) consists of elements conjugate to \(g\), since \(1-g\) is invertible. Its element orders do not increase.

From \(\omega(A_7)=\{1,2,3,4,5,6,7\}\) and the central \(2\)-kernel, we obtain
\[
 \omega(D_2)\subseteq
 \{1,2,3,4,5,6,7,8,10,12,14\}\subseteq\Om.
\]
The odd orders in \(D_2\) are \(1,3,5,7\); in the semidirect product with a module of characteristic \(3\), they can additionally yield only \(3,9,15,21\). All these numbers belong to \(\Om\). Order \(20\) cannot arise, while order \(30\) would require a nonzero fixed vector of an element of order \(10\), which has been excluded.

For \(D_3\), use \(\mu(3.A_7)=\{12,15,21\}\). The center of order \(3\) acts on the faithful simple module as scalars \(\omega,\omega^2\ne1\). Elements of orders \(12,15,21\) have nontrivial central powers, since their images in \(A_7\) have orders \(4,5,7\). Thus they have no nonzero fixed vectors, and their orders do not increase in the semidirect product. The other orders in \(D_3\) belong to \(\{1,2,3,4,5,6,7\}\); multiplying each by \(2\) gives an admissible order in \(\Om\). Hence \(\omega(K_3)\subseteq\Om\). Odd order \(9\) cannot arise in an extension by a \(2\)-group; order \(20\) would require order \(10\) in \(D_3\), and order \(30\) would require a nonzero fixed vector of an element of order \(15\). Thus all three are absent.

Finally, in both cases \([V,D_i]=V\) by irreducibility and the nontrivial action, and \(D_i\) is perfect. Thus \(K_i\) is perfect. Its radical is \(V\rtimes Z(D_i)\), since the quotient by this subgroup is \(A_7\). These radicals are not nilpotent: the center of the complement acts nontrivially on the normal Sylow subgroup \(V\).
\end{proof}

Thus, for \(A_7\), all noncentral nilpotent kernels of the highest internal mixed action have been excluded, as has the simultaneous occurrence of the two central branches. Admissible central actions nevertheless exist, so the local analysis does not complete the argument. In the next subsection, we use the orders in the full spectrum \(\Om\) that are absent from these examples, primarily \(9,20,30\), together with the lower layers of the radical.

\clearpage
\subsection{Completion of the \(A_7\) case}

The two subspectral examples in the preceding subsection cannot be excluded by a local module argument. However, the kernel \(K\) of the canonical cover has the additional property \([K,L]=K\). Choosing the highest factor on which \(K\) itself acts nontrivially forces \(K\) to be nilpotent. The required orders \(9,20,30\) then exclude both branches.

\subsubsection{Additional module facts}

\begin{lemma}\label{lem:A7-final-module-facts}
Put \(D_2=2.A_7\), \(D_3=3.A_7\). Let \(b\) denote an element of class \(3B\) in \(A_7\), that is, of type \(3^2\,1\), and also its lift of order \(3\) in \(D_2\). The following statements hold.
\begin{enumerate}
\item Every nonzero \(A_7\)-module of characteristic \(2\) contains a nonzero vector fixed by \(b\). For the simple modules, we obtain
\[
\begin{array}{c|rrrrrr}
 \dim W&1&4&4&6&14&20\\ \hline
 \dim C_W(b)&1&2&2&2&4&6\\
 \dim C_W(c_7)&1&1&1&0&2&2,
\end{array}
\]
where \(|c_7|=7\). On the \(6\)-dimensional module, the norm of an element of order \(4\) has rank \(1\).
\item On every faithful simple \(D_2\)-module of characteristic \(3\), the operator \(1+b+b^2\) is nonzero. The following ATLAS models include all absolutely simple types and the models over the prime field; repetitions upon restriction of scalars are retained. Ranks are taken over the indicated field:
\[
\begin{array}{c|rrrrrrr}
 \dim W&4&4&6&6&8&12&36\\ \hline
 q&9&9&9&9&3&3&3\\
 \operatorname{rk}(1+b+b^2)&1&1&2&2&2&4&12.
\end{array}
\]
\item If a \(D_2\)-module \(W\) of characteristic \(3\) or \(7\) satisfies \(\omega(W\rtimes D_2)\subseteq\Om\), then the central involution \(z\) acts on it as \(-I\). In characteristic \(3\), if \(W\ne0\), we also have \((1+b+b^2)W\ne0\).
\item If a \(D_3\)-module \(W\) of characteristic \(2\) or \(5\) satisfies \(\omega(W\rtimes D_3)\subseteq\Om\), then \(C_W(Z(D_3))=0\).
\end{enumerate}
\end{lemma}
\begin{proof}
The first two tables are verified using exact ATLAS matrices \cite{A7LocalData}. A self-contained certificate is given in Appendix~\ref{app:A7-final-certificate}. To specify class \(3B\) unambiguously, with standard generators \(a=(1\,2\,3)\), \(d=(3\,4\,5\,6\,7)\), take \(h=dadad\) and \(b=a h^{-1}ah=(1\,2\,3)(4\,5\,6)\). An element of order \(4\) is given by the word \(a^2dad^2\). For \(D_2\), use the standard lifts of orders \(3,5\). The certificate checks their orders, the order of \(b\), the action of the central involution, and all the ranks listed. Completeness of the simple-module lists is taken from \cite{A7LocalData}. The dual \(4\)-dimensional \(A_7\)-modules have the same norm ranks and fixed-space dimensions.

For an arbitrary module in statement 1, fixed-space dimensions are additive over composition factors, since the order of \(b\) is coprime to \(2\). For statement 2, a nonzero norm on any simple composition factor implies a nonzero norm on the whole module.

We prove statement 3. The subspace \(C_W(z)\) is invariant and is an \(A_7\)-module. In characteristic \(7\), such a nonzero module is forbidden by Lemma~\ref{lem:A7-cover-prime-exclusions}. In characteristic \(3\), all admissible simple \(A_7\)-factors are trivial, by \cite{ZhLM}. Since \(A_7\) is perfect, it acts trivially on all of \(C_W(z)\). Then an element of order \(8\) in \(D_2\) and a nonzero vector of this subspace give order \(24\). Thus \(C_W(z)=0\), that is, \(z=-I\). In characteristic \(3\), every simple factor is now faithful on \(D_2\), and statement 2 applies.

In statement 4, the fixed space of the center is again an \(A_7\)-module. In characteristic \(5\), it is forbidden by Lemma~\ref{lem:A7-cover-prime-exclusions}. In characteristic \(2\), it suffices to consider a simple composition factor of this space. By the first table, every factor except the \(6\)-dimensional one contains a nonzero vector fixed by a \(7\)-element. The group \(D_3\) contains an element of order \(21\), obtained by multiplying such a \(7\)-element by a central \(3\)-element. The center acts trivially on the factor in question, so order \(42\) arises. For the \(6\)-dimensional factor, take an element of order \(12\), the product of a \(4\)-element and a central \(3\)-element. In characteristic \(2\), its norm equals the norm of the \(4\)-element and has rank \(1\). This gives order \(24\). All possibilities are forbidden, proving statement 4.
\end{proof}

\begin{lemma}\label{lem:A7-kernel-commutator}
Under the assumptions of Theorem~\ref{thm:A7-highest-central-pairs}, the canonical kernel \(K\) from Lemma~\ref{lem:A7-canonical-cover-quotient} satisfies
\[
 K=[P,L]=[K,L],\qquad
 L=G^{(\infty)},\quad P=R(G)\cap L.
\]
\end{lemma}
\begin{proof}
Put \(D=[P,L]\). Since \(P/K\) is central in \(L/K\), we have \(D\le K\). The group \(L/D\) is perfect, and \(P/D\) is central in it; hence \(L/D\) is one of the covers \(A_7,2.A_7,3.A_7,6.A_7\). It maps onto the already chosen \(L/K=r.A_7\), where \(r=2\) or \(3\). The cover \(6.A_7\) is forbidden by order \(24\), and the other nontrivial cover cannot have \(r.A_7\) as a quotient. Thus \(L/D=L/K\), that is, \(D=K\).

Put \(H=[K,L]\). In the perfect group \(L/H\), the subgroup \(K/H\) is central. Each element of \(P\) induces an automorphism acting trivially on \(L/K\). By Lemma~\ref{lem:perfect-automorphism-tools}, it also acts trivially on \(L/H\). Hence \([P,L]\le H\). The equality \([P,L]=K\) just proved gives \(K=H\).
\end{proof}

\subsubsection{Nilpotence of the canonical kernel}

\begin{theorem}\label{thm:A7-canonical-kernel-nilpotent}
If \(\omega(G)=\Om\) and \(\Soc(G/R(G))\cong A_7\), then the canonical kernel \(K\) from Lemma~\ref{lem:A7-canonical-cover-quotient} is nilpotent.
\end{theorem}
\begin{proof}
Suppose otherwise. In a chief series of \(G\) inside \(K\), choose the highest factor \(V=N_0/M_0\) on which \(K\) acts nontrivially. Such a factor exists, since otherwise the series would be a central series of \(K\). Put
\[
 p=\operatorname{char}V,\quad
 Y=G/C_G(V),\quad E=Y^{(\infty)},
 \quad U=\text{the image of }K\text{ in }Y.
\]
Since \(K\le P\), both \(R\) and \(L\) act nontrivially on the chosen factor, and Lemma~\ref{lem:mixed-chief-action} applies. The group \(U\ne1\) is nilpotent: \(K/N_0\) has a central series, and \(N_0\le C_G(V)\). It is normal in \(Y\), so \(p\nmid |U|\), since \(O_p(Y)=1\). Lemma~\ref{lem:A7-kernel-commutator} gives \([U,E]=U\). The group \(E/U\) is a nontrivial perfect image of \(L/K=r.A_7\), and hence is \(A_7\) or \(r.A_7\). Moreover, \(V\) is a faithful simple \(Y\)-module, and \(\omega(V\rtimes Y)\subseteq\Om\).

If \(E/U=A_7\), then \(U=R(E)\). The centrality argument for the nilpotent kernel in the proof of Theorem~\ref{thm:A7-highest-central-pairs} applies directly: it uses only the nilpotence of \(U\), the condition \(p\nmid|U|\), and the spectral inclusion. Specifically, the components \(U_5,U_7\) are excluded by ordinary \(A_7\)-modules; for \(p=2\), the remaining \(3\)-kernel is centralized by \(E\); for \(p=3\), statement 4 of Lemma~\ref{lem:A7-two-kernel-three} centralizes the \(2\)-kernel; for \(p=5,7\), Lemma~\ref{lem:odd-mixed-quasisimple} applies. We obtain \([U,E]=1\), contrary to \([U,E]=U\ne1\).

Now let \(E/U=r.A_7\). For every chief \(E\)-factor \(W\) inside a Sylow subgroup \(U_q\), the action of \(U\) is trivial. Indeed, \(U_q\) is a normal \(q\)-group and thus acts trivially on a simple module of characteristic \(q\), while the other Sylow components of \(U\) commute with \(U_q\). Consequently, \(W\) is an \(r.A_7\)-module, with
\[
 \omega(W\rtimes r.A_7)\subseteq\omega(E)
 \subseteq\Om
\]
by Lemma~\ref{lem:split-spectrum}. We use this for all factors inside \(U_q\), not merely the Frattini quotient.

\medskip\noindent\textbf{1. \(E/U=2.A_7\).}
The component \(U_5\) is absent by Lemma~\ref{lem:A7-cover-prime-exclusions}. For \(q=3,7\), the central involution \(z\) of \(2.A_7\) acts as \(-I\) on each of the indicated chief \(q\)-factors, by Lemma~\ref{lem:A7-final-module-facts}. Lift an element of order \(8\) in \(2.A_7\) to \(a\in E\). Its order is \(8\), since \(\Om\) contains no other multiples of \(8\). The involution \(t=a^4\) maps to \(z\). Successive lifting of fixed points for a coprime automorphism gives \(C_{U_q}(t)=1\).

If \(U_q\ne1\), then \(p\ne q\), and the restriction of \(V\) to \(U_q\) is completely reducible. Faithfulness ensures that a nontrivial simple type occurs. The involution \(t\) fixes none of these types, by Lemma~\ref{lem:fixedfree-q-characters}. Since \(t=a^4\), the orbit of \(a\) on such a type has length \(8\). The norm \(1+a+\cdots+a^7\) is nonzero, giving order \(8p\in\{16,24,40,56\}\), forbidden in \(\Om\). Thus \(U_3=U_7=1\), and \(U\) is a \(2\)-group.

For \(p=2\), this contradicts \(p\nmid |U|\). For \(p=3\), the entire radical \(Q=R(E)\) is a \(2\)-group, since \(Q/U=C_2\). Statement 4 of Lemma~\ref{lem:A7-two-kernel-three} gives \([E,Q]=1\), contrary to \([U,E]=U\). For \(p=5,7\), Lemma~\ref{lem:odd-mixed-quasisimple} gives the same contradiction. The first branch is excluded.

\medskip\noindent\textbf{2. \(E/U=3.A_7\).}
The component \(U_7\) is absent by Lemma~\ref{lem:A7-cover-prime-exclusions}. All simple \(3.A_7\)-modules of characteristic \(3\) have trivial action of the central \(3\)-subgroup. The admissible \(A_7\)-factors are trivial. By Lemma~\ref{lem:perfect-automorphism-tools}, \(E\) centralizes \(U_3\). Since \([U_3,E]=U_3\), we have \(U_3=1\). Hence \(\pi(U)\subseteq\{2,5\}\).

The center of order \(3\) of \(E/U\) acts without nonzero fixed vectors on every chief factor inside \(U_2,U_5\), by statement 4 of Lemma~\ref{lem:A7-final-module-facts}. In the preimage of this center, choose an element \(c\) of order \(3\) by Schur--Zassenhaus. Lifting fixed points along the entire chief series gives \(C_U(c)=1\). Since \(c\) acts trivially on \(E/U\) and \((3,|U|)=1\), Lemma~\ref{lem:coprime-fixed-general} gives \(E=U C_E(c)\). The intersection is \(C_U(c)=1\), so \(E=U\rtimes D\), where \(D=C_E(c)\cong3.A_7\) and \(c\in Z(D)\).

The irreducible types of the nilpotent group \(U=U_2\times U_5\) are tensor products of the types of its Sylow components. On each component, \(c\) fixes no nontrivial type, by Lemma~\ref{lem:fixedfree-q-characters}. Thus all nontrivial \(U\)-types in \(V\) form orbits of length \(3\) under \(c\). Denote their nonempty set by \(\Sigma\); it is nonempty since \(V\) is faithful and \(U\ne1\).

Let \(d\in D\) have order \(4\). It commutes with \(c\). If \(cd\) has an orbit of length \(12\) on \(\Sigma\), its norm on \(V\) is nonzero, giving \(12p\in\{24,36,60,84\}\), which is forbidden. Hence \(d\) cannot have an orbit of length \(4\) on the set of \(c\)-orbits: such an orbit would give \(12\) distinct types. Consequently, \(d^2\) preserves each \(c\)-orbit. On a three-element orbit it commutes with a \(3\)-cycle and has order at most \(2\), so it fixes each element of that orbit. Thus \(d^2\) fixes all of \(\Sigma\).

Every involution of \(3.A_7\) is the square of an element of order \(4\). In \(A_7\), an involution of type \(2^2\,1^3\) is the square of an element of type \(4\,2\,1\). A lift of the latter may be chosen of order \(4\); its square is the unique involution in the preimage of the original involution. Furthermore, the involutions generate \(3.A_7\): their normal closure maps onto the simple group \(A_7\), and the remaining central quotient is both solvable and perfect. Thus all of \(D\), including \(c\), fixes \(\Sigma\). This contradicts the free action of \(c\) on \(\Sigma\).

Both cases are impossible. Therefore \(K\) is nilpotent.
\end{proof}

\Needspace{12\baselineskip}
\subsubsection{Contradiction with the full spectrum}
\begin{theorem}\label{thm:A7-arbitrary-radical}
There is no finite group \(G\) with spectrum \(\omega(G)=\Om\) such that \(\Soc(G/R(G))\cong A_7\). The solvable radical is arbitrary; criticality is not assumed.
\end{theorem}
\begin{proof}
Retain the notation \(R,L,P,K\), and put \(B=G/K\), \(T=R/K\), \(D=L/K\). By the preceding theorem, \(K\) is nilpotent. Its Sylow subgroups \(K_q\) are characteristic and normal in \(G\). On every chief factor inside \(K\), the group \(K\) acts trivially and the action of \(L\) factors through \(D\). The corresponding semidirect products with \(D\) are subspectral to \(G\), by Lemma~\ref{lem:split-spectrum}, as in the preceding proof.

\medskip\noindent\textbf{1. \(D=3.A_7\).}
By Lemma~\ref{lem:A7-canonical-cover-quotient}, the group \(T\) has exponent \(3\), centralizes \(D\), and \([B:DT]\le2\). The group \(B\) has no element of order \(9\). Indeed, every \(3\)-element lies in \(DT\), and in a representation \(x=dt\), its component \(d\) can be chosen to be a \(3\)-element: its image in \(A_7\) has order \(1\) or \(3\), and the kernel \(Z(D)\) is also a \(3\)-group. The group \(3.A_7\) has no element of order \(9\), so \(d^3=1\). Also \(t^3=1\), \([d,t]=1\), whence \(x^3=1\). Since order \(9\) is required in \(G\), we obtain \(K_3\ne1\).

The group \(K\) centralizes all composition factors of \(K_3/\Phi(K_3)\), since it is nilpotent. The action of \(L\) on this space factors through \(D=3.A_7\). All admissible simple factors of characteristic \(3\) are trivial, as in the proof of the preceding theorem. Thus \(L\) centralizes \(K_3\), by Lemma~\ref{lem:perfect-automorphism-tools}.

Factor by \(K_{3'}\), the product of the other Sylow subgroups of \(K\). In the image \(L_0\), the subgroup \(K_3\) is central, and \(L_0/K_3=3.A_7\). The image \(P_0\) of \(P\) satisfies \(P_0/K_3=C_3\). Every element of \(P_0\) acts trivially on \(L_0/K_3\), and thus centralizes the perfect group \(L_0\), by the same lemma. Hence \(L_0\) is a central cover of \(A_7\) with \(3\)-kernel \(P_0\). The multiplier \(C_6\) gives \(|P_0|\le3\). But \(|P_0/K_3|=3\), so \(K_3=1\), a contradiction.

\medskip\noindent\textbf{2. \(D=2.A_7\).}
Here \(5\nmid|R|\), while \(T\) is elementary abelian and centralizes \(D\). The group \(B\) has no element of order \(9\): its normal \(2\)-subgroup \(T\) has quotient \(A_7\) or \(S_7\), neither of which has an element of this order. Hence \(N=K_3\ne1\).

Put \(G_1=G/(K_2K_7)\). Order \(9\) is preserved. Identify \(N\) with its image in \(G_1\); then \(G_1/N=B\). On every chief factor inside \(N\), the action of \(L\) factors through \(D\). By Lemma~\ref{lem:A7-final-module-facts}, the central involution \(z\in D\) acts as \(-I\). It is central in all of \(B\): \(D\trianglelefteq B\), and \(z\) is its unique involution. Lift \(z\) to an involution \(t\in G_1\); this is possible in the preimage of \(\langle z\rangle\) by Schur--Zassenhaus, since \(N\) has odd order. Lifting fixed points along the chief series gives \(C_N(t)=1\).

An automorphism of order \(2\) acting fixed-point-freely on a finite group is inversion, and the group is abelian. Recall the proof: the map \(x\mapsto x^{-1}x^t\) is injective, since equality of two values gives a fixed element \(yx^{-1}\), and is therefore bijective. The automorphism \(t\) inverts every value of this map, so it inverts the entire group; inversion is an automorphism only in an abelian group. Thus \(N\) is abelian and \(x^t=x^{-1}\).

The automorphism \(t\) acts trivially on \(G_1/N\). By Lemma~\ref{lem:coprime-fixed-general}, \(C_{G_1}(t)\) maps onto \(B\), with kernel \(C_N(t)=1\). Consequently,
\[
 G_1=N\rtimes B,
\]
where \(z\) in this complement inverts \(N\).

We prove that \(T=\langle z\rangle\). On the nonzero space \(W=N/3N\), the central element \(z\) acts as \(-I\). Thus the norm of an element \(b\in D\) of class \(3B\) is nonzero, by Lemma~\ref{lem:A7-final-module-facts}. If \(|T|\ge4\), a four-group in \(T\) would centralize \(b\), and Lemma~\ref{lem:norm-three-four-group} would give order \(18\) in \(W\rtimes B\), a quotient of \(G_1\). This is impossible. Hence \(T=\langle z\rangle\), \(R=P\), and \([B:D]\le2\).

Now we prove \(K_2=1\). Suppose \(K_2\ne1\). If \(N=K_3\) contains an element of order \(9\), it commutes with an involution of \(K_2\), since \(K\) is nilpotent, immediately giving \(18\). Hence \(N\) is elementary abelian. The nonzero norm of \(b\) on \(N\) gives an element of order \(9\) in \(G_1=N\rtimes B\) with image \(b\) in \(B\). Its lift \(y\in G\) also has order \(9\), since \(9\) is the only multiple of \(9\) in \(\Om\).

The action on \(K_2/\Phi(K_2)\) factors through \(B\), since \(K\) is nilpotent. All simple factors of the restriction to \(D\) are \(A_7\)-modules: a central normal \(2\)-subgroup acts trivially on a simple module of characteristic \(2\). By the first table of Lemma~\ref{lem:A7-final-module-facts}, the element \(b\) has a nonzero fixed vector on every such factor. Thus \(C_{K_2/\Phi(K_2)}(y)\ne0\). Lifting for the coprime automorphism \(y\) gives a nontrivial \(2\)-subgroup \(C_{K_2}(y)\). Again order \(18\) arises. This contradiction gives \(K_2=1\), \(K=K_3\times K_7\).

Order \(20\) survives in \(B=G/K\), since \((20,|K|)=1\). There is no such order in \(D=2.A_7\), so \(B/\langle z\rangle\cong S_7\). Consider an element \(c\in B\) of order \(5\). Its centralizer in \(S_7\) is cyclic of order \(10\). The full preimage of this centralizer in \(B\) centralizes \(c\): the corresponding commutator lies in the central subgroup \(\langle z\rangle\), and its fifth power is \(1\), so it is trivial. Therefore \(|C_B(c)|=20\).

All elements of order \(5\) are conjugate in \(B\): this holds in \(S_7\), and in the preimage of a cyclic subgroup of order \(5\), each of its elements of order \(5\) has a unique lift of order \(5\). The existence of an element of order \(20\) now implies \(C_B(c)\cong C_{20}\) for every such \(c\). If \(|g|=10\), then \(g\in C_B(g^2)\), so \(g^5=z\), the unique involution of this cyclic group.

The element \(z\) inverts \(N\). Thus \(C_N(g)=1\) for every \(10\)-element \(g\in B\), and every element of the coset \(Ng\) is conjugate to \(g\), since the endomorphism \(1-g\) of the finite abelian group \(N\) is invertible. The group \(B\) itself has no element of order \(30\): the image of such an element in \(S_7\) would have order \(15\) or \(30\), both impossible. An element of order \(30\) in \(G_1=N\rtimes B\) would have to map to a \(10\)-element of \(B\), but every element of such a coset has order \(10\). Hence \(30\notin\omega(G_1)\).

However, we now have \(G_1=G/K_7\), and factoring by a \(7\)-group preserves the required order \(30\). This final contradiction excludes the second branch and completes the proof.
\end{proof}

\clearpage
\section{The canonical commutator kernel}
\subsection{A general localization argument for the four remaining socles}

We do not assume that the radical is nilpotent here. We now localize the required element orders in the characteristic subgroup \(K=[P,L]\), where \(L=G^{\infty}\) and \(P=R(G)\cap L\). Passing to a central cover kills \(K\), while \(P\) may have a nontrivial central image.

\begin{lemma}\label{lem:remaining-cover-prime-barrier}
Let \(\omega(B)\subseteq\Om\), and suppose that \(B\) contains a normal perfect subgroup \(D\) with \(D/Z(D)\cong S\).
\begin{enumerate}
\item If \(S=A_5\) or \(A_6\), then \(7\notin\omega(B)\).
\item If \(S=L_2(7)\), then \(5\notin\omega(B)\).
\item If \(S=A_5\) or \(L_2(7)\), then \(9\notin\omega(B)\).
\end{enumerate}
\end{lemma}
\begin{proof}
We use the following simple-group data \cite{RemainingCoverData}: the multipliers of \(A_5,A_6,L_2(7)\) are \(C_2,C_6,C_2\), respectively, and their outer automorphism groups have orders \(2,4,2\). A perfect central extension of a simple group is a quotient of its universal central cover. Hence \(Z(D)\) is a quotient of the corresponding multiplier.

The map \(\Aut(D)\to\Aut(S)\) is injective: an automorphism of a perfect group which is trivial modulo its center is trivial on the entire group, by Lemma~\ref{lem:perfect-automorphism-tools}. If \(S=A_5,A_6\) and \(|x|=7\) in \(B\), then \(x\) acts trivially on \(D\), since \(7\nmid|\Aut(S)|\). The group \(D\) contains an element of order \(5\): this prime divides \(|S|\), and hence \(|D|\). We obtain \(35\in\omega(B)\), a contradiction. For \(S=L_2(7)\), the same argument applies to an element of order \(5\) and an element of order \(7\) in \(D\). This proves the first two statements.

We prove the third. Here \(|Z(D)|\le2\), \(|\Out(S)|=2\), and the Sylow \(3\)-subgroups of \(S\) have order \(3\). Let \(|x|=9\). The induced automorphism of \(S\) is inner and has order \(1\) or \(3\). It can be represented by conjugation by an element \(a\in D\) with \(a^3=1\): if the inner automorphism is nontrivial, the preimage of the corresponding \(C_3\) in \(D\) has a complement of order \(3\) to the central group of order at most \(2\). If the action is trivial, take \(a=1\). By the injectivity just mentioned, the actions of \(x\) and \(a\) on \(D\) coincide. In particular, \([x,a]=1\), and \(y=xa^{-1}\) centralizes \(D\). But \(y^3=x^3\ne1\) and \(y^9=1\), so \(|y|=9\). Any involution \(t\in D\) now gives \(|yt|=18\notin\Om\).
\end{proof}

\begin{lemma}\label{lem:normal-cover-L34-no9}
Let \(\omega(B)\subseteq\Om\), and suppose that \(B\) contains a normal perfect subgroup \(D\) with \(D/Z(D)\cong L_3(4)\). Then \(9\notin\omega(B)\).
\end{lemma}
\begin{proof}
The multiplier of \(S=L_3(4)\) is \(C_3\times C_4\times C_4\). Put \(Z_2=O_2(Z(D))\). This subgroup is characteristic in \(D\), and hence normal in \(B\). Factoring by \(Z_2\) preserves the order of an element of order \(9\). The group \(D/Z_2\) is either \(S\) or \(3.S\cong\operatorname{SL}_3(4)\): the odd part of the multiplier has order \(3\), and the natural group \(\operatorname{SL}_3(4)\) realizes this central extension \cite{RemainingCoverData}. In the first case, Lemma~\ref{lem:normal-L34-no9} applies. It therefore suffices to consider \(D=\operatorname{SL}_3(4)\).

Suppose \(x\in B\) has order \(9\). Its induced automorphism of \(S\) lies in \(\operatorname{PGL}_3(4)\), since \(\Aut(S)/\operatorname{PGL}_3(4)\) is a \(2\)-group. There is no order \(9\) in \(\operatorname{PGL}_3(4)\). Recall a direct verification. An odd-order element of \(\operatorname{GL}_3(4)\) is semisimple; the decompositions of the space according to the degrees of irreducible factors have types \(1+1+1,1+2,3\). The corresponding tori modulo scalars have orders \(9,15,21\), respectively; the first torus is \(C_3\times C_3\), so none of the three tori contains an element of order \(9\).

Thus the action of \(x\) on \(S\) is given by some matrix \(A\in\operatorname{GL}_3(4)\), where
\[
 A^3=cI,\qquad c\in\mathbb F_4^{\times}.
\]
Conjugation by \(A\) and the action of \(x\) also coincide on \(D\), since \(D\) is perfect. Consider two cases.

\emph{The matrix \(A\) is reducible over \(\mathbb F_4\).}
Its characteristic polynomial of degree \(3\) then has a linear factor. For an eigenvalue \(\lambda\in\mathbb F_4^{\times}\), we have \(c=\lambda^3=1\). Hence \(A^3=I\), and \(A\) is diagonalizable over \(\mathbb F_4\), since \(t^3-1\) splits and has no repeated roots. If an eigenvalue has multiplicity greater than one, the centralizer of \(A\) in \(D\) contains an involution: take a transvection on a two-dimensional subspace of the corresponding eigenspace and the identity on a complement. This involution centralizes \(x\), giving order \(18\). If all three eigenvalues are distinct, they are \(1,\omega,\omega^2\), where \(\omega^3=1\ne\omega\). Thus \(\det A=1\), so \(A\in D\). The element \(y=xA^{-1}\) centralizes \(D\), \([x,A]=1\), and \(y^3=x^3\ne1\). Again \(|y|=9\), and an involution of \(D\) gives order \(18\).

\emph{The matrix \(A\) is irreducible.}
Its centralizer in \(\operatorname{GL}_3(4)\) identifies with \(\mathbb F_{64}^{\times}\). In this cyclic group of order \(63\), take an element \(b\) of order \(7\). Its determinant is the norm \(\mathbb F_{64}^{\times}\to\mathbb F_4^{\times}\), namely \(b^{21}=1\). Hence \(b\in D\), and \([x,b]=1\). The forbidden order \(63\) arises.

All possibilities are exhausted. This contradiction proves the lemma for \(D\) and, by the first paragraph, for all its central covers.
\end{proof}

\begin{lemma}\label{lem:remaining-commutator-localization}
Let \(\omega(G)=\Om\), and let \(S=\Soc(G/R(G))\) be one of the four remaining socles. Put
\[
 L=G^{\infty},\quad P=R(G)\cap L,\quad K=[P,L].
\]
Then \(K\) is characteristic in \(G\), \(K=[K,L]\), and the following conditions hold:
\[
\begin{array}{c|l}
S&\text{required localization in }K\\ \hline
A_5&\text{all elements of order }7;\quad
       x^3\text{ for every }|x|=9\\
A_6&\text{all elements of order }7\\
L_2(7)&\text{all elements of order }5;\quad
       x^3\text{ for every }|x|=9\\
L_3(4)&x^3\text{ for every }|x|=9.
\end{array}
\]
In particular, \(7\mid|K|\) for \(S=A_5,A_6\), \(15\mid|K|\) for \(S=L_2(7)\), \(3\mid|K|\) for \(S=L_3(4)\), and \(21\mid|K|\) for \(S=A_5\).
\end{lemma}
\begin{proof}
The groups \(R(G),L,P\) are characteristic in \(G\), and so is \(K\). In \(B=G/K\), the subgroup \(D=L/K\) is normal, perfect, and a central cover of \(S\): its kernel \(P/K\) is central by the definition of \(K\), and \(D/(P/K)=L/P\cong S\). Moreover, \(Z(D)=P/K\), since \(Z(S)=1\).

By Lemma~\ref{lem:remaining-cover-prime-barrier}, \(B\) has no element of order \(7\) for \(S=A_5,A_6\), and no element of order \(5\) for \(S=L_2(7)\). Hence every element of the corresponding prime order in \(G\) belongs to \(K\). By the same lemma and Lemma~\ref{lem:normal-cover-L34-no9}, for \(S=A_5,L_2(7),L_3(4)\) the group \(B\) has no element of order \(9\). The image of any element \(x\in G\) of order \(9\) has order \(1\) or \(3\), so \(x^3\in K\), and \(|x^3|=3\). All the indicated elements exist because \(\omega(G)=\Om\).

It remains to prove \(K=[K,L]\) without assuming any structure of \(P\). Put \(H=[K,L]\) and pass to \(\overline L=L/H\). The image \(\overline K\) is central in the perfect group \(\overline L\). Each \(p\in P\) acts trivially on \(L/K\), since \([P,L]=K\). Hence its action on \(\overline L\) is trivial modulo the central subgroup \(\overline K\), and therefore trivial on all of \(\overline L\), by Lemma~\ref{lem:perfect-automorphism-tools}. Thus \([P,L]\le H\), that is, \(K=H\).
\end{proof}

\begin{lemma}\label{lem:remaining-nilpotent-K-nine}
Under the assumptions of Lemma~\ref{lem:remaining-commutator-localization}, suppose additionally that \(K\) is nilpotent and \(S=A_5\) or \(L_2(7)\). Then its Sylow \(3\)-subgroup is nontrivial and has exponent \(3\). Every element \(x\in G\) of order \(9\) lies outside \(K\), and \(\langle x\rangle\cap K=\langle x^3\rangle\) has order \(3\).
\end{lemma}
\begin{proof}
The preceding lemma gives \(3\mid|K|\), as well as \(7\mid|K|\) for \(A_5\) and \(5\mid|K|\) for \(L_2(7)\). If the Sylow \(3\)-subgroup of \(K\) contained an element of order \(9\), this would commute with an element of order \(7\) or \(5\) in another Sylow subgroup of \(K\). Order \(63\) or \(45\) would arise, both forbidden. Thus \(9\notin\omega(K)\), and its nontrivial Sylow \(3\)-subgroup has exponent \(3\). For \(|x|=9\), we have \(x\notin K\), but \(x^3\in K\), giving the last equality.
\end{proof}

These results exclude the central branch \(K=1\) simultaneously for all four socles in question and localize specific required elements in \(K\). They do not yet exclude \(K\ne1\). In particular, nilpotence of \(K\) has not been deduced from nilpotence or nonnilpotence of \(P\); this is a separate question, and in the last lemma it is an explicit additional hypothesis.

\clearpage
\section{\texorpdfstring{The socle \(A_5\)}{The socle A5}}
\subsection{The branch with nontrivial nilpotent \(K\)}

In this subsection we retain the notation \(L=G^\infty\), \(P=R(G)\cap L\), \(K=[P,L]\). The nilpotence of \(K\) is a separate hypothesis; no nilpotence assumption is made about \(P\). The preliminary lemma includes two other socles for later use.

\begin{lemma}\label{lem:nilpotent-K-prime-shapes}
Let \(\omega(G)=\Om\), and suppose \(K\) is nilpotent. Denote its Sylow \(p\)-subgroup by \(K_p\); trivial factors are allowed unless explicitly stated otherwise. Then:
\begin{enumerate}
\item For \(S=A_5\), we have \(K=K_3\times K_7\), both factors are nontrivial, \(\exp K_3=3\), and \(\exp K_7=7\).
\item For \(S=A_6\), we have \(K=K_2\times K_3\times K_7\), \(K_7\ne1\), \(\exp K_2\mid2\), \(\exp K_3\mid3\), \(\exp K_7=7\), and at least one of \(K_2,K_3\) is trivial.
\item For \(S=L_2(7)\), we have \(K=K_2\times K_3\times K_5\), \(K_3,K_5\ne1\), \(\exp K_2\mid4\), \(\exp K_3=3\), and \(\exp K_5=5\).
\end{enumerate}
In each case, \([K_p,L]=K_p\) for every \(p\).
\end{lemma}
\begin{proof}
All prime divisors of \(|G|\) belong to \(\{2,3,5,7\}\), by Cauchy's theorem. The divisibility relations \(7\mid|K|\) for \(S=A_5,A_6\) and \(15\mid|K|\) for \(S=L_2(7)\) were established in Lemma~\ref{lem:remaining-commutator-localization}. For \(S=A_5\), the same lemma gives \(3\mid|K|\). Distinct Sylow subgroups of \(K\) commute.

If \(7\mid|K|\), then \(K_5=1\), since otherwise order \(35\) arises. An element of order \(4\) in \(K_2\), together with \(K_7\), gives order \(28\); an element of order \(9\) in \(K_3\) gives order \(63\). Simultaneous nontriviality of \(K_2,K_3,K_7\) gives \(42\). For \(S=A_5\), this excludes \(K_2\); for \(S=A_6\), it leaves the indicated alternative. Order \(49\) is itself forbidden.

If \(5\mid|K|\), then \(K_7=1\); orders \(8\) in \(K_2\) and \(9\) in \(K_3\) give \(40\) and \(45\), respectively. Order \(25\) is also forbidden. This proves the third statement and all assertions about exponents. The equality \(K=[K,L]\) and the characteristic property of the Sylow subgroups of \(K\) imply \([K_p,L]=K_p\): project the equality onto each direct factor.
\end{proof}

\begin{lemma}\label{lem:normal-A5-no8}
If \(\omega(B)\subseteq\Om\) and \(B\) contains a normal subgroup \(D\cong A_5\), then \(8\notin\omega(B)\).
\end{lemma}
\begin{proof}
Put \(Q=B/D\). The map
\[
 B\longrightarrow\Aut(D)\times Q,
 \qquad b\longmapsto(\alpha_b,bD)
\]
is injective: its kernel is \(C_B(D)\cap D=Z(D)=1\). Here \(\Aut(D)\cong S_5\). Suppose \(|b|=8\), and denote its image by \((\alpha,q)\). Since every \(2\)-element of \(S_5\) has order at most \(4\), the equality \(8=\mathop{\rm lcm}(|\alpha|,|q|)\) gives \(|q|=8\).

If \(\alpha\) is inner, choose \(d\in D\) so that \(y=db\) acts trivially on \(D\). If \(\alpha\) is outer, choose \(d\) so that \(y=db\) acts by conjugation by some transposition \(\tau\in S_5\). This is possible since the inner automorphisms form \(A_5\), and all odd permutations lie in the same coset modulo \(A_5\). In both cases, the image of \(y\) in the direct product is \((\beta,q)\), where \(|\beta|\mid2\), so \(|y|=8\).

In the first case, \(y\) centralizes all of \(D\). In the second, it centralizes a \(3\)-cycle on the three points not moved by \(\tau\). Hence in either case there exists \(a\in D\) of order \(3\) with \([y,a]=1\). We obtain \(|ya|=24\), a contradiction.
\end{proof}

\begin{theorem}\label{thm:A5-nilpotent-K-structure}
Let \(\omega(G)=\Om\), \(\Soc(G/R(G))\cong A_5\), and suppose that \(K=[R(G)\cap L,L]\) is nilpotent. Then
\[
 G/R(G)\cong S_5,\qquad L/K\cong2.A_5,
 \qquad K=K_3\times K_7.
\]
Both factors \(K_3,K_7\) are nontrivial; \(\exp K_3=3\), and \(K_7\) is an elementary abelian \(7\)-group. Put \(B=G/K\), \(D=L/K\), and \(T=R(G)/K\). Then
\[
 T=C_B(D),\qquad \omega(T)\subseteq\{1,2,3,4,6\}.
\]
The unique involution \(z\) of \(D\cong\operatorname{SL}_2(5)\) is central in \(B\). The group \(P=R(G)\cap L\) contains an involution \(t\) mapping to \(z\), and \(t\) inverts every element of \(K_7\). Moreover,
\[
 P=(K_3\times K_7)\rtimes\langle t\rangle,
 \qquad G/K_3\cong K_7\rtimes B.
\]
Thus the branch \(G/R(G)\cong A_5\) with nilpotent \(K\) is excluded; the statement does not exclude the remaining branch with quotient \(S_5\).
\end{theorem}
\begin{proof}
By Lemma~\ref{lem:nilpotent-K-prime-shapes}, \(K=K_3\times K_7\) has odd order, both factors are nontrivial, and they have the stated exponents. The group \(D=L/K\) is perfect and is a central cover of \(A_5\). Since the multiplier of \(A_5\) is \(C_2\), we have \(D=A_5\) or \(2.A_5\). The image in \(B\) of an element of order \(8\) in \(G\) retains its order, since \(K\) has odd order. Lemma~\ref{lem:normal-A5-no8} excludes \(D=A_5\), so \(D=2.A_5\).

The unique involution of \(\operatorname{SL}_2(5)\) is \(z=-I\). Indeed, a matrix of order \(2\) is diagonalizable in odd characteristic, and determinant \(1\) leaves only \(I\) and \(-I\). Thus \(z\) is characteristic in the normal subgroup \(D\) and central in \(B\). The group \(D\) also contains an element \(a\) of order \(4\) with \(a^2=z\), for example the matrix with rows \((0,-1),(1,0)\).

Consider any chief factor \(V\) of \(G\) inside \(K_7\). The nilpotent group \(K\) acts trivially on \(V\): its Sylow \(7\)-subgroup acts trivially on a simple module of characteristic \(7\), and the other Sylow subgroups centralize \(K_7\). Thus \(V\) is a \(B\)-module, and
\[
 \omega(V\rtimes D)\subseteq\Om
\]
by Lemma~\ref{lem:split-spectrum}. This is the same application of that lemma as for chief factors of a nilpotent radical; normality and nilpotence of \(K\) suffice here.

We prove \(C_V(z)=0\). If \(W=C_V(z)\ne0\), then \(W\) is a nonzero \(A_5=D/\langle z\rangle\)-module. On every nonzero \(A_5\)-module of characteristic \(7\), an involution has a nonzero fixed vector. Indeed, by Maschke's theorem we may restrict to a simple constituent. For a trivial constituent the assertion is clear. On a nontrivial constituent the representation of the simple group is faithful; an involution with no nonzero fixed vectors would act as the scalar \(-I\), making it central in the faithful image, contrary to \(Z(A_5)=1\). The image of \(a\) on \(W\) is such an involution. Hence \(C_W(a)\ne0\), and the semidirect product contains an element of order \(28\), which is forbidden. Thus \(z\) acts as \(-I\) on each of the indicated chief factors.

We have \(P/K=Z(D)=\langle z\rangle\). By Schur--Zassenhaus, the odd-order normal subgroup \(K\) of \(P\) has a complement \(\langle t\rangle\) of order \(2\). Conjugation by \(t\) is fixed-point-free on every chief factor inside \(K_7\). Applying Lemma~\ref{lem:coprime-fixed-general} successively along this series gives \(C_{K_7}(t)=1\). It is important here to use the entire series, not merely the Frattini quotient.

We recall a short proof that a fixed-point-free automorphism of order \(2\) acts by inversion. The map \(u\mapsto u^{-1}u^t\) on \(K_7\) is injective: equality of the images of \(u,v\) means \((vu^{-1})^t=vu^{-1}\), whence \(u=v\). Since the group is finite, the map is surjective. Every value is inverted by \(t\), so \(t\) inverts all of \(K_7\). Inversion can be an automorphism only of an abelian group. Together with \(\exp K_7=7\), this shows that the group is elementary abelian.

Every \(r\in R(G)\) acts trivially on \(D/Z(D)\), so it acts trivially on the perfect group \(D\), by Lemma~\ref{lem:perfect-automorphism-tools}. Hence \(T\le C_B(D)\). Conversely, the centralizer of the socle in the almost simple group \(G/R(G)\) is trivial, so \(C_B(D)\le T\). Elements of order \(5\) or \(7\) in \(T\) would give order \(35\) together with \(D\). An element of order \(8\) gives \(24\) with an element of order \(3\) in \(D\); order \(9\) gives \(18\) with \(z\), and order \(12\) gives \(60\) with an element of order \(5\) in \(D\). The list \(\Om\) therefore gives \(\omega(T)\subseteq\{1,2,3,4,6\}\).

The group \(DT\) has no element of order \(8\). Indeed, if \(b=du\in DT\) is a \(2\)-element, then \(dZ(D)\) is a \(2\)-element of \(A_5\). Thus \(|d|\mid4\). Since \([d,u]=1\), the element \(u=d^{-1}b\) is also a \(2\)-element, and \(|u|\mid4\) by the restriction on \(T\) just obtained. Hence \(b^4=1\). But \(8\in\omega(B)\), and \([B:DT]\le2\). Therefore \([B:DT]=2\), that is, \(G/R(G)\cong S_5\).

Finally, in \(G/K_3\), the normal subgroup \(K/K_3\cong K_7\) is abelian, and the image of \(t\) inverts it and acts trivially on the quotient \(B\), since \(z\in Z(B)\). By Lemma~\ref{lem:coprime-fixed-general}, the centralizer of this image of \(t\) maps surjectively onto \(B\). Its intersection with \(K_7\) is trivial, so it is a complement isomorphic to \(B\). We obtain \(G/K_3\cong K_7\rtimes B\), completing the proof.
\end{proof}

To complete this branch for socle \(A_5\), we must exclude the extension just described, with quotient \(S_5\) and a nontrivial group \(K_3\) of exponent \(3\). The case of nonnilpotent \(K\) is separate.

\clearpage

\subsection{The action of the central involution on \(K_3\)}

Here we continue the branch of Theorem~\ref{thm:A5-nilpotent-K-structure}. Thus \(S=A_5\), the group \(K\) is nilpotent, \(K=K_3\times K_7\), \(B=G/K\), \(D=L/K\cong2.A_5\), \(T=R(G)/K=C_B(D)\), and \(z\) is the central involution of \(D\).

\begin{lemma}\label{lem:A5-three-norm-elementary}
Let \(W\) be an \(A_5\)-module of characteristic \(3\) on which \(A_5\) acts nontrivially. For any element \(a\in A_5\) of order \(3\), the operator \(1+a+a^2\) on \(W\) is nonzero.
\end{lemma}
\begin{proof}
In \(A_5\), take a subgroup \(A_4=V_4\rtimes C_3\), where the generator \(a\) cyclically permutes the three nontrivial linear characters of \(V_4\). The action of \(V_4\) on \(W\) is nontrivial: otherwise the kernel of the \(A_5\)-action, containing \(V_4\ne1\), would be the entire simple group. After extending the field, restriction to \(V_4\) decomposes into weight spaces, since \(3\nmid4\). Take a nonzero vector \(w\) of nontrivial weight. The vectors \(w,aw,a^2w\) lie in distinct direct summands, so \((1+a+a^2)w\ne0\). An operator remains nonzero under field extension. All elements of order \(3\) in \(A_5\) are conjugate, proving the assertion for every such \(a\).
\end{proof}

\begin{lemma}\label{lem:A5-three-Frattini-minus}
Under the assumptions of Theorem~\ref{thm:A5-nilpotent-K-structure}, the involution \(z\) acts as \(-I\) on \(V=K_3/\Phi(K_3)\). Hence, for an involution \(t\in P\) mapping to \(z\),
\[
 C_{K_3}(t)\le\Phi(K_3).
\]
If \(K_3\) is abelian, then \(t\) inverts every element of \(K_3\), and the entire group splits:
\[
 G\cong(K_3\times K_7)\rtimes B.
\]
If \(K_3\) is nonabelian, then \(1\ne C_{K_3}(t)\le\Phi(K_3)\).
\end{lemma}
\begin{proof}
The group \(V\) is a \(B\)-module, since \(K_3\) acts trivially on its Frattini quotient and \(K_7\) centralizes \(K_3\). By Lemma~\ref{lem:split-spectrum}, \(\omega(V\rtimes D)\subseteq\Om\). Centrality of \(z\) gives a decomposition of \(D\)-modules \(V=V_+\oplus V_-\), on which \(z\) acts as \(I,-I\), respectively. On \(V_+\), the action factors through \(A_5\).

Suppose \(A_5\) acts nontrivially on \(V_+\). Choose an element \(a\in D\) of order \(3\); it maps to an element of order \(3\) in the quotient \(A_5\). By Lemma~\ref{lem:A5-three-norm-elementary}, the operator \(N_3(a)=1+a+a^2\) is nonzero on \(V_+\). The element \(az\) has order \(6\), and acts on \(V_+\) in the same way as \(a\). Hence
\[
 N_6(az)=1+az+\cdots+(az)^5
       =2N_3(a)\ne0.
\]
Order \(18\) arises in \(V\rtimes D\), a contradiction. Therefore \(D\) acts trivially on \(V_+\).

But \([K_3,L]=K_3\), by Lemma~\ref{lem:nilpotent-K-prime-shapes}. On the Frattini quotient this gives \([V,D]=V\). Projecting this equality onto \(V_+\), where the action is trivial, shows \(V_+=0\). Thus \(z\) acts as \(-I\) on all of \(V\). Any element of \(C_{K_3}(t)\) has fixed image in \(V\), giving \(C_{K_3}(t)\le\Phi(K_3)\).

If \(K_3\) is abelian, it is elementary abelian, since \(\exp K_3=3\). Hence \(\Phi(K_3)=1\), and \(t\) inverts all its elements. It already inverts \(K_7\), by Theorem~\ref{thm:A5-nilpotent-K-structure}. Thus \(C_K(t)=1\). Since \(z\in Z(B)\), Lemma~\ref{lem:coprime-fixed-general} gives a surjection \(C_G(t)\to B\) with trivial kernel. Consequently, \(C_G(t)\) is a complement to \(K\), proving the splitting. Finally, if \(C_{K_3}(t)=1\), a fixed-point-free automorphism of order \(2\) would make \(K_3\) abelian, as proved in Theorem~\ref{thm:A5-nilpotent-K-structure}. For nonabelian \(K_3\), the centralizer is therefore nontrivial.
\end{proof}

\begin{lemma}\label{lem:A5-centralizer-five-T}
In the branch under consideration, every element of order \(5\) in \(B\) lies in \(D\). For every \(c\in D\) of order \(5\),
\[
 C_B(c)=\langle c\rangle\times T.
\]
In particular,
\[
 20\in\omega(B)\ \Longleftrightarrow\ 4\in\omega(T),
 \qquad
 30\in\omega(B)\ \Longleftrightarrow\ 6\in\omega(T).
\]
The required order \(20\) in \(G\) forces \(4\in\omega(T)\).
\end{lemma}
\begin{proof}
We have \(B/T\cong S_5\) and \(D/(D\cap T)\cong A_5\), where \(D\cap T=\langle z\rangle\). If \(b\in B\) has order \(5\), its image in \(S_5\) is a \(5\)-cycle, since \(5\) does not divide \(|T|\). There is an element \(c\in D\) of order \(5\) with the same image: the preimage of the corresponding \(C_5\) is a central extension with kernel \(C_2\) and contains a complement \(C_5\). Write \(b=cu\), \(u\in T\). Since \([c,u]=1\), the equalities \(b^5=c^5=1\) imply \(u^5=1\), so \(u=1\). Thus \(b=c\in D\).

The centralizer of a \(5\)-cycle in \(S_5\) is the subgroup of order \(5\) it generates. Hence every element of \(C_B(c)\) belongs to \(\langle c\rangle T\), and the reverse inclusion follows from \([D,T]=1\). The factors intersect trivially, proving the formula. An element of order \(20\) or \(30\) centralizes its subgroup of order \(5\); the centralizer formula and \(\omega(T)\subseteq\{1,2,3,4,6\}\) give both equivalences. The reverse implications also follow from the direct product. Finally, factoring by \(K=K_3\times K_7\) preserves order \(20\), since it is coprime to \(|K|\).
\end{proof}

\begin{lemma}\label{lem:A5-abelian-K3-thirty}
Suppose additionally that \(K_3\) is abelian. Then \(30\in\omega(G)\) is equivalent to the following alternative:
\begin{enumerate}
\item \(6\in\omega(T)\); or
\item there exist \(c\in D\) of order \(5\) and an involution \(u\in T\) such that \(C_{K_3}(c)\cap C_{K_3}(u)\ne0\).
\end{enumerate}
In the second alternative, necessarily \(u\ne z\). If the first alternative does not hold, every pair in the second alternative also satisfies \(C_{K_7}(u)=0\), that is, \(u\) inverts \(K_7\). In particular, either \(3\mid|T|\), or \(T\) contains a four-group \(\langle z,u\rangle\cong C_2^2\). If \(T\) is nilpotent, it is a \(2\)-group containing elements of order \(4\) and a four-group.
\end{lemma}
\begin{proof}
By Lemma~\ref{lem:A5-three-Frattini-minus}, we have \(G\cong(K_3\times K_7)\rtimes B\). Factoring by \(K_7\) preserves order \(30\), so it suffices to study \(H=K_3\rtimes B\). The image of an element of order \(30\) in \(B\) has order \(30\) or \(10\): the kernel is an elementary abelian \(3\)-group, and its intersection with a cyclic group of order \(30\) has order \(1\) or \(3\). The first case is equivalent to \(6\in\omega(T)\), by Lemma~\ref{lem:A5-centralizer-five-T}.

An element \(b\in B\) of order \(10\) belongs to the centralizer of its subgroup of order \(5\), and hence is uniquely expressible as \(b=cu\), where \(|c|=5\), \(c\in D\), and \(u\in T\) is an involution. In characteristic \(3\), restriction to \(\langle b\rangle\) is completely reducible. The coset \(K_3b\) contains an element of order \(30\) exactly when \(C_{K_3}(b)\ne0\): the tenth power of \((v,b)\) is \((N_{10}(b)v,1)\), and the norm operator is multiplication by \(10\) on the fixed space and is zero on a complement. Since \(c\) and \(u\) are powers of \(b\),
\[
 C_{K_3}(b)=C_{K_3}(c)\cap C_{K_3}(u).
\]
This proves the equivalence and \(u\ne z\), since \(z\) inverts \(K_3\). If also \(C_{K_7}(u)\ne0\), choose nonzero fixed vectors in \(K_3,K_7\). In the split group, these and \(u\) give pairwise commuting elements of orders \(3,7,2\), yielding the forbidden order \(42\). Hence \(C_{K_7}(u)=0\). In characteristic \(7\), an involution without nonzero fixed vectors acts as \(-I\).

If \(6\in\omega(T)\), then \(3\mid|T|\). Otherwise we have obtained an involution \(u\ne z\), and \(z\in Z(T)\), giving a four-group. If \(T\) is nilpotent, the required element of order \(4\) from the preceding lemma would commute with an element of order \(3\) if \(3\mid|T|\). This would give \(12\notin\omega(T)\). Thus \(T\) is a \(2\)-group, and the second alternative holds.
\end{proof}

Thus the abelian \(K_3\) subbranch reduces to the joint action of \(B\) on two elementary abelian groups and the indicated alternative for \(T\). When \(K_3\) is nonabelian, the fixed points of \(t\) lie in its Frattini subgroup and cannot be discarded. No splitting of the whole group \(G\) is asserted in this subbranch.

\begin{lemma}\label{lem:double-A5-quadratic-three}
Let \(W\) be a module for \(D=2.A_5\cong\operatorname{SL}_2(5)\) of characteristic \(3\), on which the central involution \(z\) acts as \(-I\). If, for an element \(a\in D\) of order \(3\), the operator \(1+a+a^2\) is zero, then \(C_W(c)=0\) for every element \(c\in D\) of order \(5\).
\end{lemma}
\begin{proof}
We use the following classification of simple modules from the ATLAS catalog \cite{A5ThreeData}: the absolutely simple modules of \(2.A_5\) in characteristic \(3\) on which \(z=-I\) have dimensions \(2,2,6\). The two two-dimensional modules are realized over \(\mathbb F_9\), and the six-dimensional module over \(\mathbb F_3\). The four-dimensional module over \(\mathbb F_3\) also listed in the catalog splits into two two-dimensional modules after field extension and is not a fourth absolutely simple module.

On the six-dimensional module, a standard element of order \(3\) has the following matrix \(B\). Its norm in characteristic \(3\) is displayed explicitly on the right:
\[
 B=\begin{pmatrix}
1&0&0&0&0&0\\
0&0&1&0&0&0\\
0&0&0&0&1&0\\
1&1&1&1&1&0\\
0&1&0&0&0&0\\
0&2&1&2&0&1
\end{pmatrix},\qquad
 I+B+B^2=\begin{pmatrix}
0&0&0&0&0&0\\
0&1&1&0&1&0\\
0&1&1&0&1&0\\
0&0&0&0&0&0\\
0&1&1&0&1&0\\
2&0&0&0&0&0
\end{pmatrix}.
\]
The norm has rank \(2\). All elements of order \(3\) in \(2.A_5\) are conjugate: their images are conjugate in \(A_5\), and the preimage of each element of order \(3\) contains exactly one element of order \(3\), the other having order \(6\). Thus a zero norm on \(W\) excludes six-dimensional composition factors even after extending the field.

On each two-dimensional factor, an element \(c\) of order \(5\) has no nonzero fixed vector. Indeed, the determinant of a representation of the perfect group \(D\) is trivial. If one eigenvalue of \(c\) is \(1\), the other is also \(1\), and since \(c\) is semisimple in characteristic \(3\), it acts trivially. But these nontrivial simple representations are faithful: every proper normal subgroup of \(D\) is central, and \(z\) acts nontrivially. This is a contradiction. Restriction to \(\langle c\rangle\) is completely reducible, so the absence of nonzero fixed vectors on all composition factors gives \(C_W(c)=0\) on the whole module \(W\).
\end{proof}

\clearpage
\begin{lemma}\label{lem:A5-abelian-K3-T-nonnilpotent}
Under the assumptions of Theorem~\ref{thm:A5-nilpotent-K-structure}, if \(K_3\) is abelian, then
\[
 4,6\in\omega(T),\qquad T=R(G)/K
 \text{ is nonnilpotent}.
\]
In particular, the subbranch with abelian \(K_3\) and nilpotent \(R(G)/K\) is impossible.
\end{lemma}
\begin{proof}
We have \(G\cong(K_3\times K_7)\rtimes B\). We show that the second alternative in Lemma~\ref{lem:A5-abelian-K3-thirty} is impossible. For any involution \(u\in T\), the space \(W=C_{K_3}(u)\) is a \(D\)-module, since \([T,D]=1\). The involution \(z\) acts on \(W\) as \(-I\), by Lemma~\ref{lem:A5-three-Frattini-minus}. Take \(a\in D\) of order \(3\). The element \(au\) has order \(6\) and acts on \(W\) in the same way as \(a\). If \(N_3(a)\ne0\) on \(W\), then \(N_6(au)=2N_3(a)\ne0\) would give the forbidden order \(18\) in \(K_3\rtimes B\). Hence \(N_3(a)=0\) on \(W\). By Lemma~\ref{lem:double-A5-quadratic-three}, \(C_W(c)=0\) for every \(c\in D\) of order \(5\). Thus
\[
 C_{K_3}(u)\cap C_{K_3}(c)=0,
\]
and the second alternative is impossible for all pairs \(u,c\).

The required order \(30\) therefore forces \(6\in\omega(T)\), and order \(20\) forces \(4\in\omega(T)\), by Lemma~\ref{lem:A5-centralizer-five-T}. If \(T\) were nilpotent, an element of order \(4\) in its Sylow \(2\)-subgroup would commute with an element of order \(3\) in its Sylow \(3\)-subgroup. Then \(12\in\omega(T)\), contrary to \(\omega(T)\subseteq\{1,2,3,4,6\}\). Thus \(T\) is nonnilpotent.
\end{proof}

The preceding lemma excludes one subbranch. To complete the case of nilpotent \(K\), it remains to consider abelian \(K_3\) with nonnilpotent \(T\) containing elements of orders \(4,6\), and nonabelian \(K_3\) with nontrivial fixed points of the involution inside \(\Phi(K_3)\).

\clearpage

\subsection{A Sylow \(2\)-subgroup of the centralizer \(T\)}

We continue the abelian \(K_3\) subbranch. By Lemma~\ref{lem:A5-abelian-K3-T-nonnilpotent}, \(T\) is solvable, \(4,6\in\omega(T)\), and \(\omega(T)\subseteq\{1,2,3,4,6\}\). We retain \(D=L/K\cong2.A_5\), \(T=C_B(D)\), \(D\cap T=\langle z\rangle\), \(B/T\cong S_5\), \([B:DT]=2\), and \(8\in\omega(B)\).

\begin{lemma}\label{lem:A5-eight-square-T-involution}
Every element \(b\in B\) of order \(8\) satisfies \(b^4=z\) and induces on \(D/Z(D)\) an automorphism of the type of a \(4\)-cycle in \(S_5\). We can write
\[
 b^2=dv,\qquad d\in D,\quad d^2=z,
 \quad v\in T,\quad v^2=1.
\]
If \(z\) is the only involution of \(T\), then \(b^2\) centralizes \(T\).
\end{lemma}
\begin{proof}
The group \(B/\langle z\rangle\) contains a normal subgroup \(D/\langle z\rangle\cong A_5\) and has spectrum contained in \(\Om\). By Lemma~\ref{lem:normal-A5-no8}, it has no element of order \(8\). Hence the unique involution of the cyclic group \(\langle b\rangle\) is \(b^4=z\). The group \(DT\) has no element of order \(8\), as proved in Theorem~\ref{thm:A5-nilpotent-K-structure}. Thus \(b\notin DT\), and the automorphism it induces on \(D\) is outer. In \(\Aut(D)\cong S_5\), an odd \(2\)-element is a transposition or a \(4\)-cycle. In the first case, \(b\) would centralize an element of order \(3\) in \(D\): a transposition fixes a \(3\)-cycle on the other three points, and its unique lift of order \(3\) is also fixed. Order \(24\) would arise. Thus only a \(4\)-cycle remains.

Its square is a double transposition. Any lift \(d\in D\) of it has order \(4\) and square \(z\), since the only involution of \(D\) is central. Choosing such a \(d\) with the same action on \(D\) as \(b^2\), we obtain \(b^2=dv\), \(v\in T\). Since \([D,T]=1\), the equality \(b^4=z\) gives \(v^2=1\). If \(T\) has a unique involution, then \(v=1\) or \(z\), whence \([b^2,T]=1\).
\end{proof}

\begin{lemma}\label{lem:A5-T-no-cyclic-Sylow}
A Sylow \(2\)-subgroup of \(T\) cannot be cyclic.
\end{lemma}
\begin{proof}
Suppose otherwise. By the spectrum restriction, it has order \(4\). First we construct a characteristic Hall \(3\)-subgroup \(A\) of \(T\) without using the normal complement theorem. Let \(\sigma_T:T\to\{\pm1\}\) be the sign of the left regular action on the set \(T\). An element of order \(4\) acts as an odd number of \(4\)-cycles, so \(H=\ker\sigma_T\) has index \(2\). This subgroup is characteristic: every automorphism permutes the points of the regular action and preserves sign. A Sylow \(2\)-subgroup of \(H\) has order \(2\). For the sign homomorphism \(\sigma_H\) of the regular action, an involution is an odd permutation, since \(|H|/2\) is odd. Thus \(A=\ker\sigma_H\) has index \(2\) in \(H\), has odd order, and is characteristic in \(T\). We have \(|T:A|=4\) and \(A\ne1\), since \(3\mid|T|\). For any generator \(u\) of a Sylow \(2\)-subgroup, \(T=A\rtimes\langle u\rangle\). The central involution \(z\) is the only involution of \(T\): every involution is conjugate to an involution of the cyclic Sylow subgroup, and \(z\) is central. Hence \(u^2=z\).

From \(12\notin\omega(T)\), we obtain \(C_A(u)=1\). The action of \(u\) on \(A\) has order at most \(2\), since \(u^2=z\) is central. A fixed-point-free automorphism of order \(2\) inverts the group and makes it abelian, as proved in Theorem~\ref{thm:A5-nilpotent-K-structure}. Thus \(u\) inverts \(A\), which is an elementary abelian \(3\)-group.

Choose \(b\in B\) of order \(8\) and a Sylow \(2\)-subgroup \(Q\) of \(B\) containing \(b\). The intersection \(Q\cap T\) is normal in \(Q\) and is a cyclic Sylow \(2\)-subgroup of \(T\). Let \(u\) now be its generator. Then \(bub^{-1}=u\) or \(u^{-1}\). The group \(A\) is characteristic in \(T\), and hence normal in \(B\). By Lemma~\ref{lem:A5-eight-square-T-involution}, \(b^2\) centralizes \(A\). The element \(b\) fixes no nontrivial element of \(A\), since otherwise order \(24\) would arise. Thus \(b\) also inverts \(A\). The product \(bu\) centralizes \(A\). For either possible action of \(b\) on \(\langle u\rangle\), we directly obtain \((bu)^4=b^4=z\); hence \(|bu|=8\). Together with a nonidentity element of \(A\), it gives order \(24\), a contradiction.
\end{proof}

\begin{lemma}\label{lem:A5-T-no-quaternion-Sylow}
A Sylow \(2\)-subgroup of \(T\) is not isomorphic to \(Q_8\).
\end{lemma}
\begin{proof}
Suppose this Sylow subgroup is \(U\cong Q_8\). The only involution of \(T\) is then \(z\), by conjugacy of Sylow subgroups and centrality of \(z\). Put \(A=O_3(T)\). Every element \(u\in U\) of order \(4\) satisfies \(u^2=z\), and therefore acts on \(A\) with order dividing \(2\). If \(A\ne1\), then \(12\notin\omega(T)\) gives \(C_A(u)=1\), so \(u\) inverts \(A\). Choose \(u,v\in U\) of order \(4\) whose product \(uv\) also has order \(4\). Both \(u,v\) invert \(A\), so \(uv\) centralizes it, a contradiction. Thus \(O_3(T)=1\).

The Fitting subgroup \(F=F(T)\) is a nontrivial \(2\)-group containing \(z\), and lies in \(U\). For a solvable group, \(C_T(F)\le F\). Recall a proof of this property. The subgroup \(C=C_T(F)\) is normal in \(T\), so its Fitting subgroup \(F(C)\) lies in \(F\) and is central in \(C\). If \(C/F(C)\ne1\), the preimage of the last nontrivial term of its derived series would be a nilpotent normal subgroup of \(C\) strictly larger than \(F(C)\): it is a central extension of an abelian group. This is a contradiction. Thus \(C=F(C)\le F\). If \(F\) is cyclic, it has order \(2\) or \(4\), and \(\Aut(F)\) is a \(2\)-group. From \(T/C_T(F)\le\Aut(F)\) and \(C_T(F)\le F\), it would follow that \(T\) is a \(2\)-group, contrary to \(3\mid|T|\). Hence \(F=U=Q_8\).

Since \(C_T(U)=Z(U)=\langle z\rangle\) and \(|\Aut(Q_8)|=24\), we have \(|T|\le48\). The Sylow \(2\)-subgroup has order \(8\), and \(3\mid|T|\), so \(|T|=24\). A complement of order \(3\) to the normal \(Q_8\) acts nontrivially, since otherwise there would be order \(12\). Thus \(T\cong Q_8\rtimes C_3\cong\operatorname{SL}_2(3)\).

We also need \(Z(T)=\langle z\rangle\), \(\Aut(T)\cong S_4\), and \(\operatorname{Inn}(T)\cong A_4\). These facts can be checked directly. The subgroup \(Q_8=T'\) is characteristic. An automorphism of \(T\) which is trivial on \(Q_8\) is trivial on all of \(T\): the image of a generator of the complement of order \(3\) can differ from that generator only by \(z\), which would violate its order \(3\). Thus \(|\Aut(T)|\le|\Aut(Q_8)|=24\). Conjugation by \(\operatorname{GL}_2(3)\) gives \(24\) distinct automorphisms of \(\operatorname{SL}_2(3)\), with kernel consisting of scalar matrices. The action on the four lines of \(\mathbb F_3^2\) identifies \(\operatorname{GL}_2(3)/\{\pm I\}\) with \(S_4\), and the inner automorphisms form \(A_4\).

Choose \(b\in B\) of order \(8\). If its action on \(T\) is inner, it is given by an element \(s\in T\), \(|s|\mid4\): the inner automorphism has order dividing \(4\), and \(\operatorname{Inn}(T)=A_4\) has no element of order \(4\). Moreover, \([b,s]=1\), and \(bs^{-1}\) centralizes \(T\). Its fourth power is \(b^4s^{-4}=z\), so it has order \(8\) and gives \(24\) with an element of order \(3\) in \(T\).

It remains to consider an outer action on \(T\). By multiplying \(b\) on the right by an element \(t\in T\), we can make its action on \(T\) a \(4\)-cycle in \(S_4\): the outer automorphisms form the odd coset. Independently, multiplying on the left by \(d\in D\), we can make its action on \(D\) a transposition in \(S_5\). These actions do not interfere, since \([D,T]=1\). The resulting \(h=dbt\) acts on \(D,T\) with orders \(2,4\), respectively. The kernel of the joint action \(B\to\Aut(D)\times\Aut(T)\) is \(C_B(D)\cap C_B(T)=T\cap C_B(T)=Z(T)=\langle z\rangle\). Thus \(h^4\in\langle z\rangle\) and \(h^8=1\). If \(h^4=1\), then \(h^2\in T\) would be the identity or the unique involution \(z\), and would act trivially on \(T\). This contradicts the order \(4\) of the action of \(h\) on \(T\). Hence \(|h|=8\). But a transposition on \(D\) fixes an element of order \(3\), again giving \(24\). This is a contradiction.
\end{proof}

\begin{lemma}\label{lem:A5-T-Sylow-has-four-group}
In the abelian subbranch under consideration, a Sylow \(2\)-subgroup of \(T\) has exponent \(4\), contains a four-group including \(z\), and has order at least \(8\).
\end{lemma}
\begin{proof}
We use the full statement of the classical theorem on a unique subgroup of prime order \cite{ConradQuaternion}: a finite \(p\)-group has exactly one subgroup of order \(p\) if and only if it is cyclic, or \(p=2\) and the group is generalized quaternion. A generalized quaternion group of order \(2^n\), \(n\ge3\), contains an element of order \(2^{n-1}\). Since \(8\notin\omega(T)\), a quaternion Sylow subgroup could only be \(Q_8\). This and the cyclic possibility have been excluded by the two preceding lemmas. Hence a Sylow \(2\)-subgroup contains an involution \(v\ne z\). Since \(z\in Z(T)\), we obtain \(\langle z,v\rangle\cong C_2^2\). Also \(4\in\omega(T)\), so the exponent is \(4\), and the subgroup has order at least \(8\).
\end{proof}

In the abelian \(K_3\) subbranch, the nonnilpotent group \(T\) must have a Sylow \(2\)-subgroup containing both an element of order \(4\) and an involution different from \(z\). Cyclic and quaternion Sylow subgroups have already been excluded, so they do not cover the remaining cases.

\clearpage
\begin{lemma}\label{lem:A5-T-no-dihedral-eight}
In the same abelian subbranch, a Sylow \(2\)-subgroup of \(T\) is not isomorphic to \(D_8\).
\end{lemma}
\begin{proof}
Let \(P_2\cong D_8\). The central involution \(z\) lies in the center of \(P_2\) and is the square of a generator \(u\) of its unique cyclic subgroup of order \(4\). First we show \(A=O_3(T)=1\). If \(A\ne1\), then \(C_A(u)=1\), since otherwise order \(12\) arises in \(T\). Since \(u^2=z\) centralizes \(A\), the action of \(u\) is a fixed-point-free automorphism of order \(2\). Thus \(A\) is abelian and \(u\) inverts it.

Choose \(b\in B\) of order \(8\) and a Sylow \(2\)-subgroup \(Q\) of \(B\) containing \(b\). We may assume \(P_2=Q\cap T\). In the decomposition \(b^2=dv\) from Lemma~\ref{lem:A5-eight-square-T-involution}, we may choose \(d\in Q\cap D\), so \(v\in P_2\). Indeed, \(Q\cap D\) is a Sylow \(2\)-subgroup of \(D\); its image in \(\Aut(D)=S_5\) is the four-group of double transpositions inside the Sylow dihedral subgroup containing the image of \(b\). The square of this \(4\)-cycle lies in that four-group. Conjugation by \(b\) defines an automorphism \(\alpha\) of \(P_2\), and \(\alpha^2\) is conjugation by \(v\), since \(D\) centralizes \(P_2\). Every automorphism of \(D_8\) sends \(u\) to \(u\) or \(u^{-1}\), so its square fixes \(u\). Thus \(v\) cannot be a reflection. Since \(v^2=1\), only \(v=1,z\) remain. Hence \(b^2\) centralizes \(A\). The element \(b\) acts fixed-point-freely on \(A\), since otherwise \(24\) arises, and therefore also inverts \(A\). Then \(bu\) centralizes \(A\). From \(bub^{-1}=u^{\pm1}\), we obtain \((bu)^4=z\). The element \(bu\) of order \(8\) again gives \(24\). Thus \(O_3(T)=1\).

Now \(F=F(T)\) is a normal \(2\)-subgroup containing \(z\), and \(C_T(F)\le F\). The subgroups of \(D_8\) containing \(z\) are cyclic, four-groups, or all of \(D_8\). The automorphism groups of cyclic groups of orders \(2,4\) and of \(D_8\) are \(2\)-groups. In the latter case, this follows by choosing the images \(u\mapsto u^{\pm1}\), \(v\mapsto u^kv\). If \(F\cong C_2^2\), the image of \(T\) in \(\Aut(F)\cong S_3\) fixes the nonzero element \(z\) and has order dividing \(2\). In every case, \(T/C_T(F)\) is a \(2\)-group. Together with \(C_T(F)\le F\), this makes \(T\) a \(2\)-group, contrary to \(3\mid|T|\).
\end{proof}

\begin{lemma}\label{lem:A5-T-abelian-eight-structure}
If a Sylow \(2\)-subgroup of \(T\) has order \(8\), it is isomorphic to \(C_4\times C_2\). In that case, we can write
\[
 T=(A\rtimes\langle u\rangle)
       \times\langle z\rangle,
 \qquad |u|=4,
\]
where \(A=O_3(T)\ne1\) is a normal Hall \(3\)-subgroup, \(\exp A=3\), \(C_A(u)=1\), and \([A,u^2]\ne1\). In particular, \(u^2\ne z\).
\end{lemma}
\begin{proof}
The groups of order \(8\) containing an element of order \(4\) are \(C_8,C_4\times C_2,D_8,Q_8\). The first, third, and fourth possibilities were excluded above. Let \(P_2=C_4\times C_2\) and \(A=O_3(T)\). We show that \(A\) is a Hall subgroup. In \(\overline T=T/A\), we have \(O_3(\overline T)=1\): the preimage of a normal \(3\)-subgroup would be a normal \(3\)-subgroup of \(T\). Thus \(F(\overline T)\) is a \(2\)-group containing the central involution \(\overline z\). It is isomorphic to a subgroup of \(P_2\). A cyclic subgroup or \(C_4\times C_2\) has a \(2\)-group as its automorphism group; for a four-group, the stabilizer of \(\overline z\) is also a \(2\)-group. Since \(C_{\overline T}(F(\overline T))\le F(\overline T)\), the group \(\overline T\) is a \(2\)-group. Thus \(A\) is Hall and \(A\ne1\). By Schur--Zassenhaus, \(T=A\rtimes P_2\).

Suppose first that \(z\) is the square of an element \(u\in P_2\) of order \(4\). The element \(u\) acts fixed-point-freely on \(A\), since \(12\notin\omega(T)\). The equality \(u^2=z\) implies that it inverts the abelian group \(A\). If \(v\) is the other generator of order \(2\) of \(P_2\), then \(uv\) also has order \(4\) and square \(z\), and likewise inverts \(A\). Thus \(v\) centralizes \(A\). Hence all involutions of \(P_2\) centralize \(A\).

Choose \(b,Q\) as in the preceding lemma, and assume \(P_2=Q\cap T\). In the decomposition \(b^2=dv_0\), we may take \(v_0\in P_2\), \(v_0^2=1\). Thus \(b^2\) centralizes \(A\), and \(b\) inverts \(A\), since otherwise \(24\) would arise. Conjugation by \(b\) on the abelian group \(P_2\) has square equal to the identity: \(b^2=dv_0\) centralizes it. Denote this action by \(\alpha\). Then
\[
 (bu)^4=b^4\bigl(\alpha(u)u\bigr)^2=z,
\]
since \(\alpha(u)^2=\alpha(z)=z=u^2\). But \(bu\) centralizes \(A\) and has order \(8\), again giving \(24\). Thus \(z\) is not the square of an element of order \(4\).

Now \(P_2=\langle u\rangle\times\langle z\rangle\), \(|u|=4\). Since \(z\) is central in \(T\), we obtain the stated direct product. The spectrum restriction gives \(\exp A=3\) and \(C_A(u)=1\). If \([A,u^2]=1\), then \(u\) would invert \(A\). All involutions of \(P_2\) would then centralize \(A\), and the preceding argument with \(b,Q,\alpha\) and \(bu\) would again give \(24\): it uses only \(\alpha(u)^2=u^2\), since \(u^2\) is the only nontrivial square in \(C_4\times C_2\). Hence \([A,u^2]\ne1\).
\end{proof}

When a Sylow \(2\)-subgroup of \(T\) has order \(8\), we obtain the described structure with \(C_4\times C_2\). We examine its joint action with the perfect residual.

\clearpage

\subsection{Joint action: excluding order \(8\) for a Sylow subgroup of \(T\)}

We retain the assumptions of the abelian \(K_3\) subbranch. For brevity, put \(N=K_3\), \(M=K_7\). Both groups are nontrivial elementary abelian groups, \(G=(N\times M)\rtimes B\), \(D\cong2.A_5\), and \(T=C_B(D)\). The involution \(z\) inverts both \(N\) and \(M\).

\begin{lemma}\label{lem:two-Jordan-rigidity}
Let \(W\) be a \(D=2.A_5\)-module of characteristic \(3\), with \(z=-I\). Choose an element \(a\) of order \(3\) in the normalizer of a Sylow \(Q_8\) of \(D\). Let \(s\) be a linear operator on \(W\) centralizing \(D\), with \(s^3=I\). If
\[
 1+a+a^2=0,\qquad 1+as+(as)^2=0
 \quad\text{on }W,
\]
then \(s=I\).
\end{lemma}
\begin{proof}
Pass to an algebraically closed field; vanishing and nonvanishing of operators are preserved. The normalizer of \(Q_8\) in \(D\) is the full preimage of \(A_4\) under the projection onto \(A_5\) \cite{A5ThreeData}, namely \(H=Q_8\rtimes\langle a\rangle\cong\operatorname{SL}_2(3)\). In characteristic \(3\), \(Q_8\) has exactly one simple type on which \(z=-I\), the two-dimensional type \(U\). The other four simple types are linear and trivial on \(z\). This also follows from the ordinary decomposition of the group algebra: \(4\cdot1^2+2^2=8\). The type \(U\) extends to the natural two-dimensional \(\operatorname{SL}_2(3)\)-module, on which \(a\) has matrix \(J=J_2(1)\).

Complete reducibility of restriction to \(Q_8\) gives \(W=U\otimes E\). The action of \(a\) has the form \(J\otimes C\), where \(C^3=I\): the quotient of two operators inducing the same conjugation on \(Q_8\) belongs to its commutant \(I_U\otimes\operatorname{End}(E)\). If \(C\) has a Jordan block of size \(2\), the norm operator associated with an element of order \(3\) on \(J_2(1)\otimes J_2(1)\) is nonzero. Indeed, for \(Jx=x\), \(Jy=y+x\), we obtain
\[
 (J\otimes J-I)^2(y\otimes y)
 =2x\otimes x\ne0.
\]
A block of size \(3\) also gives a nonzero norm: restrict to the tensor product of a fixed vector of the first factor and that block. All blocks of \(C\) have size at most \(3\), so the first hypothesis forces \(C=I\). Now \(a=J\otimes I\). Since \(s\) centralizes \(Q_8\), it has the form \(I\otimes S\), \(S^3=I\). The second hypothesis and the same calculation for \(J\otimes S\) give \(S=I\), that is, \(s=I\).
\end{proof}

\begin{lemma}\label{lem:A5-noncentral-square-actions}
Let \(u\in T\) have order \(4\), and let \(w=u^2\ne z\). Then
\[
 u|_M=-I,\quad w|_M=I,
 \qquad w|_N=-I,\quad(zw)|_N=I.
\]
For an element \(a\in D\) of order \(3\), we have \(1+a+a^2=0\) on \(N\). The whole group \(O_3(T)\) centralizes \(N\times M\).
\end{lemma}
\begin{proof}
Consider the action of \(u\) on \(M\) after extending the field of characteristic \(7\). The eigenvalues are \(1,-1,i,-i\), where \(i^2=-1\). Eigenvalue \(1\) is forbidden, since it gives order \(28\) in \(M\rtimes\langle u\rangle\). Suppose there is a nonzero space \(M_\lambda\), \(\lambda=\pm i\). It is \(D\)-invariant, since \([D,T]=1\). Take \(d\in D\) of order \(4\), \(d^2=z\). Its eigenvalues on \(M_\lambda\) are \(i,-i\), and both occur: \(d\) is conjugate to \(d^{-1}\) in the Sylow \(Q_8\) of \(D\), so their multiplicities are equal. Thus \(du\) has a nonzero fixed vector on \(M_\lambda\). But \((du)^2=zw\ne1\), so \(|du|=4\), again yielding \(28\). Hence the only eigenvalue on \(M\) is \(-1\), giving the first two equalities.

The involution \(w\) centralizes \(M\). If it fixed a nonzero vector of \(N\), the vectors of orders \(3,7\), together with \(w\), would give order \(42\). Thus \(w=-I\) on \(N\), and \(zw=I\) on \(N\). For \(a\in D\) of order \(3\), the element \(a(zw)\) has order \(6\) and acts on \(N\) as \(a\). The exclusion of \(18\) forces \(N_6(a(zw))=2N_3(a)=0\).

Put \(A=O_3(T)\). The operator \(u\) on \(M\) is scalar, so it centralizes the image of \(A\) in \(\operatorname{GL}(M)\). But \(C_A(u)=1\), since otherwise \(T\) would contain an element of order \(12\). By Lemma~\ref{lem:coprime-fixed-general}, the map from \(C_A(u)\) onto the fixed part of the image of \(A\) is surjective. This entire image is fixed, so it is trivial. Thus \([A,M]=1\).

On \(N\), the operator \(w=-I\) is scalar. The same lemma shows that the image of \(A\) on \(N\) equals the image of \(C_A(w)\). For every \(c\in C_A(w)\), we have \(c^3=1\), since \(\exp A\mid3\). The element \(c\) commutes with \(zw\) and with all of \(D\). For \(|c|=3\) and the chosen \(a\in D\) of order \(3\), the product \(ac\) has order \(3\), while \(aczw\) has order \(6\). The latter acts on \(N\) as \(ac\). The exclusion of \(18\) gives \(N_3(ac)=0\). Lemma~\ref{lem:two-Jordan-rigidity}, together with \(N_3(a)=0\), shows that \(c\) acts trivially on \(N\). This is also true for \(c=1\). Thus \(A\) centralizes \(N\) as well.
\end{proof}

\begin{lemma}\label{lem:A5-generic-no-normal-three-Sylow}
If \(T\) contains an element of order \(4\) whose square differs from \(z\), then a Sylow \(3\)-subgroup of \(T\) is not normal. In particular, in the abelian subbranch under consideration, a Sylow \(2\)-subgroup of \(T\) has order at least \(16\).
\end{lemma}
\begin{proof}
If \(A=O_3(T)\) were Sylow, then \(T/A\) would be a \(2\)-group. The group \(DA\) is normal in \(B\), and \(B/DA\) is also a \(2\)-group. Hence every \(3\)-element of \(B\) belongs to \(DA=D\times A\). A Sylow \(3\)-subgroup of \(D\) has order \(3\), and \(\exp A\mid3\). Thus \(B\) has no element of order \(9\), and every element \(b\) of order \(3\) in it is expressible as \(b=ac\), where \(a\in D\), \(a^3=1\), and \(c\in A\). By the preceding lemma, \(A\) acts trivially on \(N\), and \(N_3(a)=0\) for \(|a|=3\). For \(a=1\), the operator \(N_3(a)=3I\) is also zero. Thus \(N_3(b)=0\) on \(N\) for every element \(b\in B\) of order \(3\).

In \(G/M=N\rtimes B\), an element of order \(9\) could only map to an element of order \(3\): the kernel \(N\) is elementary abelian, and \(B\) has no element of order \(9\). But the formula \((v,b)^3=(N_3(b)v,1)\) and the vanishing norm just proved exclude this possibility. Hence \(9\notin\omega(G/M)\). Since \(M\) is a \(7\)-group, factoring by it preserves order \(9\) in \(G\). This is a contradiction.

If a Sylow \(2\)-subgroup of \(T\) has order \(8\), Lemma~\ref{lem:A5-T-abelian-eight-structure} gives both an element \(u\) with \(u^2\ne z\) and a normal Sylow \(3\)-subgroup \(A\). This case is therefore excluded. Smaller orders were excluded earlier, giving the lower bound \(16\).
\end{proof}

\begin{lemma}\label{lem:A5-large-two-dichotomy}
Let \(P_2\) be a Sylow \(2\)-subgroup of \(T\) in the abelian subbranch under consideration. Then \(|P_2|\ge16\), \(\exp P_2=4\), and two possibilities remain.
\begin{enumerate}
\item The square of every element of order \(4\) is \(z\). Then \(\Phi(P_2)=\langle z\rangle\), \([P_2,P_2]\le\langle z\rangle\), and \(P_2/\langle z\rangle\) is elementary abelian.
\item There is an element \(u\in P_2\) of order \(4\) with \(w=u^2\ne z\). Then the four-group \(\langle z,w\rangle\) acts on \(N\times M\) by
\[
 z=(-I,-I),\qquad
 w=(-I,I),\qquad zw=(I,-I).
\]
The group \(O_3(T)\) centralizes \(N\times M\), and \(3\mid|T/O_3(T)|\).
\end{enumerate}
\end{lemma}
\begin{proof}
The order bound was proved in the preceding lemma, and exponent \(4\) follows from the spectrum of \(T\). In the first case, the square of every element lies in \(\langle z\rangle\). The quotient by this central subgroup has exponent \(2\), and is therefore abelian. Consequently, \([P_2,P_2]\le\langle z\rangle\). The Frattini formula for a \(2\)-group gives \(\Phi(P_2)=P_2^2[P_2,P_2]
=\langle z\rangle\), since an element of order \(4\) exists. In the second case, all assertions follow from Lemmas~\ref{lem:A5-noncentral-square-actions} and~\ref{lem:A5-generic-no-normal-three-Sylow}.
\end{proof}

\begin{lemma}\label{lem:A5-T-no-normal-SL23}
The group \(T\) contains no subgroup \(E\cong\operatorname{SL}_2(3)\) normal in \(B\) whose central involution is \(z\). In particular, \(T\cong\operatorname{SL}_2(3)\times C_2^r\), with \(z\) in the first factor, is impossible for every \(r\ge0\).
\end{lemma}
\begin{proof}
Suppose that such an \(E\) exists. The centralizer \(C=C_T(E)\) has no element of order \(3\), since this would give \(12\) with an element of order \(4\) in \(E\). Nor does it have an element of order \(4\), which would give \(12\) with an element of order \(3\) in \(E\). The restriction on the spectrum of \(T\) therefore implies that \(C\) is an elementary abelian \(2\)-group.

Choose \(b\in B\) of order \(8\). By Lemma~\ref{lem:A5-eight-square-T-involution}, \(b^4=z\), and its action on \(D\) is outer. If its action on \(E\) is inner, it is induced by an element \(s\in E\) whose order divides \(4\), as in Lemma~\ref{lem:A5-T-no-quaternion-Sylow}. Then \(b\) commutes with \(s\), and \(bs^{-1}\) centralizes \(E\), has fourth power \(z\), and has order \(8\). Together with a \(3\)-element of \(E\), it gives order \(24\).

If its action on \(E\) is outer, multiply by an element of \(E\) to make it a \(4\)-cycle in \(\Aut(E)\cong S_4\). Independently multiply by an element of \(D\) to make the action on \(D\) a transposition in \(S_5\). We obtain \(h\) whose actions on \(E,D\) have orders \(4,2\). Since \(h^2\) centralizes \(D\), we have \(h^2\in T\). Its action on \(E\) is the square of a \(4\)-cycle, that is, a nontrivial involution in \(\operatorname{Inn}(E)=A_4\). It is induced by an element \(q\in E\) of order \(4\), with \(q^2=z\). Hence \(h^2=qc\), where \(c\in C_T(E)=C\). Since \([q,c]=1\) and \(c^2=1\), we obtain \(h^4=q^2c^2=z\). Consequently, \(|h|=8\). But the transposition acting on \(D\) fixes an element of order \(3\), giving \(24\).

Finally, in the direct product in the statement, the subgroup \(E\) is characteristic: it is generated by all elements of order \(3\). The group \(\operatorname{SL}_2(3)\) is generated by upper and lower transvections of order \(3\), whereas the second factor has no elements of odd order. Hence \(E\) is normal in \(B\), and the exclusion just proved applies.
\end{proof}

\begin{lemma}\label{lem:A5-T-no-C2-S4}
In the branch under consideration, \(T\cong C_2\times S_4\) is impossible.
\end{lemma}
\begin{proof}
The center of \(T\) has order \(2\), so its nonidentity element is \(z\). Let \(S\cong S_4\) be the second direct factor. Every automorphism of \(T\) fixes \(z\) and induces an inner automorphism on \(T/\langle z\rangle\cong S_4\). Here we use the standard fact \(\Aut(S_4)=\operatorname{Inn}(S_4)\). This can be seen from the natural system of four subgroups \(S_3\): these are all the subgroups of order \(6\), automorphisms permute them, and the automorphism action on them is faithful. The latter follows from faithfulness of the inner automorphism action: if an automorphism fixes all four subgroups, the permutations induced by conjugation by each \(g\in S_4\) and by its image coincide, so its image is \(g\).

Choose \(b\in B\) of order \(8\). Its action on \(T/z\) is induced by an element \(s\in S\) of order \(1,2\) or \(4\): the induced automorphism has order a power of \(2\), and \(S_4\) has trivial center. Then \(c=bs^{-1}\) acts trivially on \(T/z\). Its action on \(S\) has the form
\[
 x\longmapsto xz^{\varepsilon(x)},
 \qquad\varepsilon:S_4\to C_2.
\]
The only possibilities for \(\varepsilon\) are the trivial character and the sign character, since \(S_4'=A_4\). Thus \(c\) centralizes \(A_4\le S\). Moreover, \([c,s]\) lies in the central subgroup \(\langle z\rangle\). For such elements, \((cs)^4=c^4s^4\): the exponent of the central commutator in the formula is \(\binom42=6\), and its order divides \(2\). Since \(cs=b\), \(s^4=1\), and \(b^4=z\), we obtain \(c^4=z\), that is, \(|c|=8\). It centralizes an element of order \(3\) in \(A_4\), giving the forbidden order \(24\).
\end{proof}

These two exclusions eliminate additional families with large Sylow \(2\)-subgroups, but do not replace the general analysis of the dichotomy in Lemma~\ref{lem:A5-large-two-dichotomy}.

\clearpage

\subsection{The branch with a central Frattini subgroup of order \(2\)}

We continue the abelian \(K_3\) subbranch, putting \(N=K_3\) and \(M=K_7\). In this subsection \(P_2\) is a Sylow \(2\)-subgroup of \(T\), and \(\Phi(P_2)=\langle z\rangle\). Consequently, the square of every element of order \(4\) is \(z\).

\begin{lemma}\label{lem:transfer-central-normalizer} Let \(H\) be a finite group and \(U\) a Sylow \(p\)-subgroup of it, with \(U\le Z(N_H(U))\). Then \(H\) has a normal Hall \(p'\)-subgroup. \end{lemma} \begin{proof} We prove the Burnside criterion used here by transfer (compare \cite{TransferNotes}). The group \(U\) is abelian. If \(x,y\in U\) are conjugate in \(H\), they are conjugate in \(N_H(U)\): the two Sylow subgroups obtained from \(U\) and lying in \(C_H(y)\) can be matched by conjugation in \(C_H(y)\). By hypothesis the normalizer centralizes \(U\), so \(x=y\).

Define the transfer \(v:H\to U\). For representatives of the left cosets \(r\in H/U\) write \(gr=r' u_{g,r}\), where \(u_{g,r}\in U\), and put \(v(g)=\prod_r u_{g,r}\). Abelianness of \(U\) and permutation of representatives give \(v(gh)=v(g)v(h)\). For \(x\in U\), group the factors by cycles of the action of \(x\) on \(H/U\). A cycle of length \(m\) gives an element \(r^{-1}x^m r\in U\), which equals \(x^m\) by the first paragraph. Thus \(v(x)=x^{[H:U]}\). This power map is an automorphism of \(U\). Hence \(v\) is surjective, \(\ker v\cap U=1\), and \(|H:\ker v|=|U|\). The group \(\ker v\) is the required normal Hall subgroup. \end{proof}

\begin{lemma}\label{lem:quadratic-hyperplane-two-groups} Let \(P\) be a \(2\)-group of exponent \(4\), with \(\Phi(P)=\langle z\rangle\le Z(P)\) and \(|z|=2\). If there exists a linear character \(\lambda:P\to\{1,-1\}\) taking the value \(-1\) on every element of order \(4\), then \[
 P\cong C_4\times C_2^r
 \quad\text{or}\quad
 P\cong D_8\times C_2^r,
\] where \(z\) is the square of a generator of the first cyclic factor or of a rotation in \(D_8\), respectively. Every automorphism of order \(3\) of such a group fixes an element of order \(4\). \end{lemma} \begin{proof} Put \(V=P/\langle z\rangle\). Squares define a quadratic form \(q:V\to\mathbb F_2\) by \(x^2=z^{q(\bar x)}\). Its polar form is determined by commutators. The character \(\lambda\) defines a nonzero linear form \(\ell\) on \(V\), and \(q\) vanishes identically on \(\ker\ell\). In coordinates with \(\ell=x_1\), this means \(q=x_1(c+L)\), where \(L\) is linear in the other coordinates. If \(L=0\), then \(c=1\), the polar form is zero, and \(P\) is abelian. The classification of finite abelian groups gives \(C_4\times C_2^r\). If \(L\ne0\), a coordinate change gives \(q=x_1x_2\), and the other coordinates form the radical, on which \(q=0\). Preimages of the two isotropic basis vectors generate \(D_8\); the preimage of the radical is central and has exponent \(2\). Choosing a complement to \(\langle z\rangle\) in it gives the second direct product.

In the first case put \(R=P\); in the second let \(R\) be the subgroup generated by all elements of order \(4\). Then \(R\cong C_4\times C_2^r\) is characteristic in \(P\). The quotient \(R/\Omega_1(R)\) has order \(2\), where \(\Omega_1(R)=\{x\in R:x^2=1\}\). An automorphism of order \(3\) fixes this quotient. By Lemma~\ref{lem:coprime-fixed-general}, the nonidentity element of this quotient has a fixed lift in \(R\). Every preimage of the nonidentity element of the quotient has order \(4\). \end{proof}

\Needspace{11\baselineskip} \begin{lemma}\label{lem:A5-central-Frattini-O3-zero} In the branch under consideration, \(O_3(T)=1\). \end{lemma} \begin{proof} Suppose that \(A=O_3(T)\ne1\). Every element \(u\in P_2\) of order \(4\) has square \(z\), which centralizes \(A\). The exclusion of \(12\) gives \(C_A(u)=1\), so \(u\) inverts \(A\), and \(A\) is abelian of exponent \(3\). The action of \(P_2\) on \(A\) factors through the elementary abelian group \(P_2/\langle z\rangle\). Every linear character of it occurring on the nonzero \(\mathbb F_3\)-space \(A\) is \(-1\) on all elements of order \(4\). By the preceding lemma, \(P_2\) has one of the two forms given there.

The normalizer \(N_T(P_2)\) contains no element of order \(3\): such an element would fix an element of order \(4\), giving \(12\). Hence \(N_T(P_2)=P_2\). In \(T/\langle z\rangle\) the Sylow \(2\)-subgroup is abelian and equals its normalizer. Lemma~\ref{lem:transfer-central-normalizer} gives a normal Hall \(3\)-subgroup of this quotient. Its preimage has a central Sylow subgroup \(\langle z\rangle\) of order \(2\). By Schur--Zassenhaus it has a complement of odd order; all complements are conjugate by elements of \(\langle z\rangle\) and therefore coincide. This complement is normal in \(T\) and is its Sylow \(3\)-subgroup. Thus \(A\) is a Hall subgroup, and \(T=A\rtimes P_2\).

Choose \(b\in B\) of order \(8\) and a Sylow \(2\)-subgroup \(Q\) of \(B\) containing \(b\). Take \(P_2=Q\cap T\). As in Lemma~\ref{lem:A5-T-no-dihedral-eight}, we have \(b^2=dv\), \(d\in Q\cap D\), \(v\in P_2\), and \(v^2=1\). Denote conjugation by \(b\) on \(P_2\) by \(\alpha\); its square is conjugation by \(v\).

If \(P_2=C_4\times C_2^r\), all its involutions centralize \(A\): for an involution \(t\), the elements \(u\) and \(ut\) of order \(4\) both invert \(A\). Hence \(b^2\) centralizes \(A\). The element \(b\) itself fixes no nonidentity element of \(A\), since otherwise \(24\) occurs; thus \(b\) also inverts \(A\). For a generator \(u\) of order \(4\), we have \(\alpha^2=1\) and \((bu)^4=b^4(\alpha(u)u)^2=z\). But \(bu\) centralizes \(A\), giving \(24\).

Now let \(P_2=D_8\times E\), \(E=C_2^r\), and let \(u\) be a rotation of order \(4\). The group \(E\) centralizes \(A\), since \(u\) and \(ue\), for each \(e\in E\), both invert \(A\). The subgroup \(R=\langle u\rangle\times E\) is characteristic in \(P_2\), so \(\alpha(u)\in R\) commutes with \(u\). Also \(\alpha^2(u)=u\) or \(u^{-1}\), depending on whether \(v\) is a reflection. In both cases, \[
 (bu)^4=b^4\alpha^3(u)\alpha^2(u)\alpha(u)u=z.
\] The last equality follows from commutativity of \(R\), \(\alpha(u)^2=u^2=z\), and \(\alpha^2(u)=u^{\pm1}\).

Put \(A_-=\{x\in A:x^v=x^{-1}\}\). If \(A_-\ne1\), then \(v\) is a reflection times an element of \(E\) because all involutions of \(R\) centralize \(A\). Multiplying \(b\) by a suitable element of \(Q\cap D\) gives \(h\in Q\) acting on \(D\) as a transposition. Its action on \(A\) agrees with that of \(b\), so \(h^2\) inverts \(A_-\). The element \(h\) cannot have order \(8\), since it would give \(24\) with the element of order \(3\) in \(D\) that it fixes. Consequently, \(|h|=4\). On the nontrivial characters of the abelian group \(A_-\), its square acts by inversion, so every nontrivial orbit has length \(4\). If \(A_-\) acts nontrivially on \(M\), such an orbit gives a nonzero norm \(N_4(h)\) and order \(28\). Thus \([A_-,M]=1\). The involution \(uv\) centralizes \(A_-\). If it fixes a nonzero vector of \(M\), we obtain \(42\). Hence \(uv=-I\) on \(M\). But \(uv\) conjugates \(u\) to \(u^{-1}\); scalarity of its action would then force \(u=u^{-1}\) on \(M\), contrary to \(u^2=z=-I\). Therefore \(A_-=1\).

Now \(v\) and \(b^2\) centralize \(A\). By the exclusion of \(24\), the element \(b\) inverts \(A\). Hence \(bu\) centralizes \(A\), and \((bu)^4=z\) again gives \(24\). Every possibility is contradictory, so \(O_3(T)=1\). \end{proof}

\Needspace{11\baselineskip} \begin{lemma}\label{lem:A5-central-Frattini-normal-two} In this branch a Sylow \(2\)-subgroup \(P_2\) is normal in \(T\), and \[
 T=P_2\rtimes\langle a\rangle,
 \qquad |a|=3.
\] The center of \(P_2\) has exponent \(2\). If \(V=P_2/\langle z\rangle\) and \(R=Z(P_2)/\langle z\rangle\), then \(a\) acts without nonzero fixed vectors on the nonzero space \(V/R\). Moreover, \(|P_2|\ge32\). \end{lemma} \begin{proof} By the preceding lemma, \(O_3(T)=1\). Also \(O_3(T/\langle z\rangle)=1\): the preimage of a normal \(3\)-subgroup has central kernel of order \(2\), and its unique odd-order complement, as in the preceding proof, would give a normal \(3\)-subgroup of \(T\). The Fitting subgroup of \(\overline T=T/\langle z\rangle\) is a \(2\)-group and lies in the abelian Sylow subgroup \(\overline P_2\). The latter centralizes \(F(\overline T)\), so the property \(C_{\overline T}(F(\overline T))
\le F(\overline T)\) gives \(\overline P_2=F(\overline T)\). Thus \(P_2\) is normal in \(T\). By Schur--Zassenhaus, \(T=P_2\rtimes H\), where \(H\) is a nontrivial \(3\)-group of exponent \(3\).

Squares define an \(H\)-invariant quadratic form \(q\) on \(V\). Its polar radical is \(R\). The restriction \(q|_R\) is linear. If it is nonzero, every element \(h\in H\) of order \(3\) fixes a nonzero element of \(R/\ker(q|_R)\). By coprime lifting of fixed points it fixes a vector in \(R\) with \(q=1\), and then an element of order \(4\) in \(P_2\). This gives \(12\), a contradiction. Thus \(q|_R=0\), so \(Z(P_2)\) has exponent \(2\).

The form \(q\) induces a nondegenerate quadratic form on \(V/R\). This space is nonzero: otherwise \(P_2\) would be abelian with zero square form and have no element of order \(4\). For any \(1\ne h\in H\), its fixed space in \(V/R\) is nondegenerate for the polar form. Indeed, coprime action gives \(V/R=C_{V/R}(h)\perp[V/R,h]\) with both summands nondegenerate. A nonzero nondegenerate subspace contains a vector with \(q=1\). Its fixed preimage in \(V\), and then in \(P_2\), would again give \(12\). Consequently, \(C_{V/R}(h)=0\).

The group \(H\) contains no \(C_3^2\): after extending the field to \(\mathbb F_4\), every nontrivial weight for such a subgroup would have a nontrivial kernel, contrary to the proved absence of nonzero fixed vectors. The center of the \(3\)-group \(H\) contains a subgroup of order \(3\); any other subgroup of order \(3\) would generate \(C_3^2\) with it. Thus the subgroup of order \(3\) is unique, and by the classification quoted in~\cite{ConradQuaternion} the group \(H\) is cyclic. Since \(\exp H=3\), we obtain \(H=C_3\).

It remains to exclude \(|P_2|=16\). Then \(\dim V=3\). The nondegenerate quotient \(V/R\) has dimension \(2\), and \(R\) has dimension \(1\). The action of \(a\) on \(V/R\) is nontrivial and on \(R\) is trivial. The subspace \([V,a]\) is two-dimensional, nondegenerate, and anisotropic: \(a\) permutes its three nonzero vectors, and nondegeneracy rules out value zero on all of them. The preimage of \([V,a]\) in \(P_2\) is therefore isomorphic to \(Q_8\) and invariant under \(a\). The center of \(P_2\) is a four-group; \(a\) fixes it pointwise because an automorphism of order \(3\) cannot act nontrivially on a four-group while fixing \(z\). Choosing a central involution outside \(\langle z\rangle\) gives \(T\cong\operatorname{SL}_2(3)\times C_2\) with \(z\) in the first factor. This group is excluded by Lemma~\ref{lem:A5-T-no-normal-SL23}. Hence \(|P_2|\ge32\). \end{proof}

The branch with \(\Phi(P_2)=\langle z\rangle\) is now reduced to a normal \(2\)-subgroup of order at least \(32\) complemented by an element of order \(3\) with the specified action. These groups will be excluded in the next subsection. In the second branch, where some squares differ from \(z\), the restrictions of Lemma~\ref{lem:A5-large-two-dichotomy} remain in force.

\clearpage

\subsection{Complete exclusion of the branch \(\Phi(P_2)=\langle z\rangle\)}

This subsection completes the first branch of Lemma~\ref{lem:A5-large-two-dichotomy}. By Lemma~\ref{lem:A5-central-Frattini-normal-two}, we have \(T=P_2\rtimes\langle a\rangle\), \(|a|=3\), \(P_2\lhd B\), the center \(Z(P_2)\) has exponent \(2\), and \(a\) acts fixed-point-freely on the nonzero space \(W=P_2/Z(P_2)\). As before, \(N=K_3\), \(M=K_7\), and \(G=(N\times M)\rtimes B\).

\begin{lemma}\label{lem:double-A5-seven-simple-data}
Over an algebraically closed field of characteristic \(7\), the simple modules for \(D=2.A_5\) on which \(z=-I\) have dimensions \(2,2,4,6\). On the four-dimensional module, an element of order \(3\) has a fixed space of dimension \(2\). On the six-dimensional module, an element of order \(5\) has a fixed space of dimension \(2\). An outer automorphism of \(D\) interchanges the two two-dimensional types.
\end{lemma}
\begin{proof}
Completeness of the list follows from the ATLAS catalog \cite{A5SevenData}. The four-dimensional module in that catalog that splits over \(\mathbb F_{49}\) into two two-dimensional modules is not an additional absolutely simple type. To verify the fixed spaces, we give the matrix of a standard \(3\)-element on the four-dimensional module and both standard matrices of the six-dimensional module:
\[
 B_4=\begin{pmatrix}
0&0&1&0\\0&6&6&0\\1&2&1&0\\3&4&6&1
\end{pmatrix},\qquad
 A_6=\operatorname{diag}(J,J,J),\quad
 J=\begin{pmatrix}0&1\\6&0\end{pmatrix},
\]
\[
 B_6=\begin{pmatrix}
4&4&0&0&0&0\\
0&0&1&0&0&0\\
0&5&3&0&0&0\\
0&0&0&0&1&0\\
3&6&2&4&0&6\\
5&1&5&6&4&0
\end{pmatrix}.
\]
Direct calculation in \(\mathbb F_7\) gives \(B_4^3=B_6^3=I\), \((A_6B_6)^5=I\), \(\operatorname{rank}(B_4-I)=2\), and \(\operatorname{rank}(A_6B_6-I)=4\). These give the two asserted dimensions.

On a two-dimensional module, the determinant is identically \(1\), since \(D\) is perfect. The eigenvalues of an element of order \(5\) are \(\zeta^k,\zeta^{-k}\), where \(|\zeta|=5\) and \(k=1\) or \(2\). The two possible traces are distinct: they are the roots of \(t^2+t-1=0\). An outer automorphism of \(A_5\) interchanges the classes of \(5\)-cycles, and hence replaces the exponent \(k\) by \(2k\), up to sign. It therefore interchanges the two traces and thus the two two-dimensional types.
\end{proof}

\begin{lemma}\label{lem:A5-central-Frattini-outer-inverts-three}
Every element of \(B\setminus DT\) induces inversion on \(T/P_2\cong C_3\).
\end{lemma}
\begin{proof}
The groups \(D\) and \(T\) act trivially on \(T/P_2\), and \([B:DT]=2\). If the assertion is false, then all of \(B\) centralizes this factor. The Frattini argument gives \(B=P_2N_B(\langle a\rangle)\): all Sylow \(3\)-subgroups of \(T=P_2\rtimes\langle a\rangle\) are conjugate by elements of \(P_2\). Under the assumption under consideration, \(N_B(\langle a\rangle)=C_B(a)\), since a normalizing element cannot induce inversion on the image of \(a\) in \(T/P_2\). Thus \(B=P_2C_B(a)\).

Take \(b\in B\) of order \(8\). By Lemma~\ref{lem:A5-eight-square-T-involution}, its action on \(D\) is a \(4\)-cycle. In the decomposition \(b=pc\), \(p\in P_2\), \(c\in C_B(a)\), the element \(c\) has the same action on \(D\), since \([P_2,D]=1\). Its \(2\)-part \(c_2\) also centralizes \(a\) and acts on \(D\) as the same \(4\)-cycle. The spectrum of \(B\) implies that its order is \(4\) or \(8\). Order \(8\) immediately gives \(24\) together with \(a\). If \(|c_2|=4\), choose \(d\in D\), \(d^2=z\), inducing the same inner action on \(D\) as \(c_2^2\). Then \(c_2^2=dv\), where \(v\in T\) and \([v,a]=1\). Since \([d,v]=1\), the equality \(c_2^4=1\) implies \(v^2=z\). We obtain an element \(v\in T\) of order \(4\) centralizing \(a\), which gives \(12\in\omega(T)\). Both cases are impossible.
\end{proof}

\begin{lemma}\label{lem:A5-central-Frattini-excluded}
In the abelian \(K_3\) subbranch, the case \(\Phi(P_2)=\langle z\rangle\) is impossible for every order of the Sylow \(2\)-subgroup \(P_2\).
\end{lemma}
\begin{proof}
Choose \(b\in B\) of order \(8\) and a Sylow \(2\)-subgroup \(Q\) of \(B\) containing \(b\). Since \(P_2\lhd B\), we have \(P_2\le Q\). Multiplying \(b\) by a suitable element of \(Q\cap D\), we obtain \(h\in Q\) whose action on \(D\) is a transposition. This is possible within a Sylow dihedral subgroup of \(S_5\): the odd coset of its four-group contains both \(4\)-cycles and transpositions.

For every \(p\in P_2\), the element \(hp\) acts on \(D\) as the same transposition. Its order is \(2\) or \(4\). Indeed, it is a \(2\)-element; order \(8\) would give \(24\) with a fixed element of order \(3\) in \(D\). Moreover, \((hp)^2\in P_2\): its square centralizes \(D\), and all \(2\)-elements of \(T\) lie in the normal Sylow subgroup \(P_2\).

Denote the action of \(h\) on \(W=P_2/Z(P_2)\) by \(\beta\), and that of \(a\) by \(A\). We have \(\beta^2=I\), since \(h^2\in P_2\) and inner automorphisms of \(P_2\) act trivially on \(W\). By the preceding lemma, \(\beta A\beta^{-1}=A^{-1}\). Since \(A\) has order \(3\) and \(W\ne0\), we obtain \(\beta\ne I\). If \(k\) is the image of \(h^2\) in \(W\), then the image of \((hp)^2\) is
\[
 k+(I+\beta)x,\qquad
 x=pZ(P_2)\in W.
\]
Since \(I+\beta\ne0\), this image is nonzero for some \(p\). Put \(g=hp\), \(s=g^2\). Then \(|g|=4\), and \(s\in P_2\) is a noncentral involution.

Since \([P_2,P_2]\le\langle z\rangle\), there exists \(r\in P_2\) with \(s^r=sz\). On every nonzero \(P_2\)-invariant subspace of \(M\), the involution \(s\) therefore has eigenvalues \(1,-1\) with equal positive multiplicities: the operator \(r\) maps one eigenspace to the other, since \(z=-I\). In particular, this holds on each nonzero \(D\)-isotypic component of \(M\), since \([D,P_2]=1\).

Extend the field of characteristic \(7\) to an algebraic closure; kernel dimensions are unchanged by this extension. The restriction of \(M\) to \(D\) is completely reducible, since \(7\nmid|D|\). By Lemma~\ref{lem:double-A5-seven-simple-data}, its simple constituents have dimensions \(2,2,4,6\). A six-dimensional constituent gives a nonzero fixed vector for an element of order \(5\), and hence the forbidden order \(35\). If a four-dimensional type \(U_4\) occurs, its isotypic component has the form \(U_4\otimes E\). An element \(c\in D\) of order \(3\) has a nonzero fixed vector in \(U_4\), and \(s\) has a nonzero fixed vector in the multiplicity space \(E\), by the preceding paragraph. Their tensor is fixed by \(cs\), of order \(6\). The forbidden order \(42\) arises. Thus only the two two-dimensional types \(U,U'\) remain.

The outer automorphism induced by \(g\) interchanges these types. Hence the corresponding isotypic components \(M_U,M_{U'}\) are both nonzero and are interchanged by \(g\). Choose \(0\ne v\in C_{M_U}(s)\), which is possible because \(s\) and \(sz\) are conjugate. Then \(v+gv\ne0\), since the summands lie in different isotypic components, and
\[
 g(v+gv)=gv+g^2v=gv+sv=gv+v.
\]
Thus \(g\), of order \(4\), has a nonzero fixed vector on \(M\). In \(M\rtimes B\), we obtain order \(28\), the final contradiction.
\end{proof}

\clearpage

\subsection{Excluding the entire abelian \(K_3\) branch for \(A_5\)}

We retain the assumptions of the preceding subsections: \(K\) is nilpotent, \(N=K_3\) is abelian, \(M=K_7\), \(G=(N\times M)\rtimes B\), \(D=L/K\cong2.A_5\), and \(T=C_B(D)=R/K\). Both \(N,M\) are nontrivial elementary abelian groups, and \(z\in D\cap T\) inverts both. In this subsection we exclude the second possibility in Lemma~\ref{lem:A5-large-two-dichotomy}.

\begin{lemma}\label{lem:A5-involutions-two-signs} Every involution \(t\in T\setminus\{z\}\) has exactly one of the following actions on \(N\times M\): \[
 t=(-I,I)\quad\text{or}\quad t=(I,-I).
\] Suppose in addition that there exists \(u\in T\) of order \(4\) with \(u^2\ne z\). Let \(J=C_T(N)\) be the kernel of the action of \(T\) on \(N\). Then a Sylow \(2\)-subgroup of \(J\) has order at most \(2\), and the image \(E=T/J\le\operatorname{GL}(N)\) contains a unique involution, the scalar \(-I\). The Sylow \(2\)-subgroups of \(E\) are isomorphic to \(C_2,C_4\) or \(Q_8\). \end{lemma} \begin{proof} If an involution \(t\) fixes nonzero vectors in both groups \(N,M\), these vectors together with \(t\) give an element of order \(42\). Thus either \(t=-I\) on \(N\) or \(t=-I\) on \(M\). For \(t\ne z\), apply the same argument to the involution \(zt\). If \(t=-I\) on \(N\), then \(zt=I\) on \(N\), so \(zt=-I\) on \(M\), that is, \(t=I\) on \(M\). The other case is symmetric. Here we use the complete reducibility of the action of an involution in odd characteristic.

Every element \(x\in T\) of order \(4\) acts on \(N\) as an operator of order \(4\). For \(x^2=z\) this follows from \(z|_N=-I\), and for \(x^2\ne z\) from Lemma~\ref{lem:A5-noncentral-square-actions}. Consequently, \(J\) has no elements of order \(4\). Any Sylow \(2\)-subgroup \(U\) of \(J\) has exponent \(2\), since \(\omega(T)\subseteq\{1,2,3,4,6\}\), and is therefore elementary abelian. All its nonidentity elements act on \(M\) as \(-I\) by the first assertion; the product of two distinct such elements would be an involution acting trivially on both \(N\) and \(M\). Thus \(|U|\le2\).

Every involution in the image \(E\) lifts to a \(2\)-element \(x\in T\): take a Sylow \(2\)-subgroup of the full preimage of the cyclic subgroup it generates. This \(x\) cannot have order \(4\), since its image would then also have order \(4\). Thus \(x\) is an involution, and its nontrivial image on \(N\) is \(-I\). The scalar \(-I\) itself belongs to the image, since it is induced by \(z\). Hence the involution is unique. The classification of finite \(2\)-groups with a unique involution, stated in the proof of Lemma~\ref{lem:A5-T-Sylow-has-four-group}, shows that a Sylow subgroup is cyclic or generalized quaternion. Its exponent is at most \(4\), giving the stated list. \end{proof}

\Needspace{11\baselineskip} \begin{lemma}\label{lem:A5-generic-image-SL23} Under the assumptions of the preceding lemma, put \(A=O_3(T)\). Then \[
 J/A\cong C_2,\qquad
 E=T/J\cong\operatorname{SL}_2(3),\qquad
 |P_2|=16
\] for a Sylow \(2\)-subgroup \(P_2\) of \(T\). The subgroup \(J/A\) is central in \(T/A\), and \(O_3(E)=1\). \end{lemma} \begin{proof} By Lemma~\ref{lem:A5-noncentral-square-actions}, the group \(A\) centralizes \(N\), so \(A\le J\). A Sylow \(2\)-subgroup of \(J\) has order \(1\) or \(2\). If its order is \(2\), it is central in its normalizer; Lemma~\ref{lem:transfer-central-normalizer} gives a normal Hall \(3\)-subgroup \(A_0\) of \(J\). If \(J\) has odd order, take \(A_0=J\). This normal Hall subgroup is unique: all Hall \(3\)-subgroups here are Sylow subgroups and are conjugate. Thus \(A_0\) is characteristic in \(J\) and normal in \(T\). Hence \(A_0\le A\); the reverse inclusion follows from \(A\le J\). Consequently, \(A_0=A\) and \(|J/A|\le2\). If this quotient is nontrivial, it is central because the automorphism group of \(C_2\) is trivial.

We show that \(O_3(E)=1\). The full preimage of \(O_3(E)\) in \(T/A\) has a central kernel \(J/A\) of order at most \(2\). Its Sylow \(3\)-subgroup is unique: the number of such subgroups divides \(2\) and is congruent to \(1\) modulo \(3\). It is characteristic in the preimage and normal in \(T/A\). If \(O_3(E)\ne1\), the full preimage of this Sylow subgroup in \(T\) would be a normal \(3\)-subgroup strictly containing \(A\), which is impossible.

Next, \(3\mid|E|\). Suppose otherwise. Then all elements of \(T\) of order \(3\) act trivially on \(N\). There are no elements of order \(9\) in \(B\): the quotient \(B/DT\) has order \(2\), while \(DT\) is a central product of \(D,T\) with intersection of order \(2\); the Sylow \(3\)-subgroups of both factors have exponent \(3\). Every element \(b\in B\) of order \(3\) can be written as \(b=ds\), where \(d\in D\), \(s\in T\), and \(d^3=s^3=1\): a Sylow \(3\)-subgroup of \(D\times T\) maps isomorphically onto a Sylow \(3\)-subgroup of \(DT\). On \(N\) we have \(N_3(d)=0\) by Lemma~\ref{lem:A5-noncentral-square-actions} when \(d\ne1\); for \(d=1\) this equality is \(3I=0\). Consequently, \(N_3(b)=0\). The formula \((n,b)^3=(N_3(b)n,1)\) excludes order \(9\) in \(N\rtimes B=G/M\), whereas factoring out the \(7\)-group \(M\) preserves this order. This contradicts \(9\in\omega(G)\).

The group \(E\) is solvable and \(O_3(E)=1\), so its Fitting subgroup is \(F(E)=O_2(E)\ne1\). It lies in a Sylow subgroup from the list \(C_2,C_4,Q_8\). If \(F(E)\) is cyclic, then \(\operatorname{Aut}(F(E))\) is a \(2\)-group, and \(C_E(F(E))\le F(E)\) by the self-centralizing property of the Fitting subgroup of a solvable group. Then \(E\) would also be a \(2\)-group, contrary to \(3\mid|E|\). Therefore \(F(E)=Q_8\), and this is a Sylow \(2\)-subgroup of \(E\). Conjugation gives an embedding \(E/Q_8\hookrightarrow\operatorname{Out}(Q_8)\cong S_3\): the kernel is \(Q_8C_E(Q_8)=Q_8\). The quotient has odd order and is nontrivial, so it is \(C_3\). An element of order \(3\) in \(E\) complements \(Q_8\) and cyclically permutes its three cyclic subgroups of order \(4\). This is precisely \(E\cong Q_8\rtimes C_3\cong\operatorname{SL}_2(3)\).

Thus \(|P_2|\le2\cdot8=16\). Lemma~\ref{lem:A5-generic-no-normal-three-Sylow} gives the reverse bound \(|P_2|\ge16\). Therefore \(|P_2|=16\) and \(|J/A|=2\). \end{proof}

\Needspace{11\baselineskip} \begin{lemma}\label{lem:A5-generic-square-excluded} In the abelian \(K_3\) branch under consideration, there is no \(u\in T\) of order \(4\) with \(u^2\ne z\). \end{lemma} \begin{proof} Suppose otherwise, and choose a Sylow \(2\)-subgroup \(P_2\le T\) containing \(u\). By the preceding lemma, \(|P_2|=16\), \(J\cap P_2=\langle v\rangle\) has order \(2\), and \(P_2/\langle v\rangle\cong Q_8\). The subgroup \(\langle v\rangle\) is normal in \(P_2\) and hence central. The image of \(z\) is the central involution of \(Q_8\), and \(z\notin J\). Hence \[
 Z(P_2)=\langle z,v\rangle\cong C_2^2,
 \qquad P_2/Z(P_2)\cong C_2^2.
\] Indeed, the center of \(P_2\) maps into the center of \(Q_8\), whose full preimage is exactly \(\langle z,v\rangle\).

Let \(F\) be the full preimage of \(Q_8\lhd E\) in \(T\). Then \(F\lhd T\), \(T/F\cong C_3\), and \(P_2\) is a Sylow \(2\)-subgroup of \(F\). The Frattini argument gives \(T=F N_T(P_2)\). Thus \(N_T(P_2)\) contains an element \(a\) of order \(3\) whose image in \(T/F\) is nontrivial. More explicitly, a Sylow \(3\)-subgroup of the normalizer maps onto the quotient of order \(3\), and all its nonidentity elements have order \(3\) because \(9\notin\omega(T)\). The image of \(a\) in \(E\cong\operatorname{SL}_2(3)\) cyclically permutes the three nonzero elements of \(P_2/Z(P_2)\). Moreover, \(a\) fixes \(z\) because \(z\in Z(B)\), and fixes \(v\) because it normalizes \(J\cap P_2=\langle v\rangle\). Hence \(a\) acts trivially on \(Z(P_2)\).

For \(x\notin Z(P_2)\), its image in \(Q_8\) has order \(4\), so \(x^2\in\{z,zv\}\). Replacing \(x\) by an element of the same coset modulo \(Z(P_2)\) does not change its square, since the center has exponent \(2\). Conjugation by \(a\) cyclically permutes these three cosets and fixes their squares. Consequently, all elements of order \(4\) in \(P_2\) have the same square \(w\in\{z,zv\}\). In particular, \(w=u^2\ne z\).

Choose \(x,y\in P_2\) with distinct nonzero images in \(P_2/Z(P_2)\). Then \(x,y,xy\) lie outside the center, have order \(4\), and have square \(w\ne z\). By Lemma~\ref{lem:A5-noncentral-square-actions}, each of the three acts on \(M\) as \(-I\). But the action of \(xy\) is the product of the actions of \(x,y\), namely \(I\). Since \(M\ne0\) and its characteristic is \(7\), this is a contradiction. \end{proof}

\Needspace{11\baselineskip} \begin{lemma}\label{lem:A5-nilpotent-K-nonabelian-three} Let \(G\) be a nonsolvable group with spectrum \(\Om\), and suppose that \(\operatorname{Soc}(G/R(G))\cong A_5\), \(L=G^\infty\), \(P=R(G)\cap L\), and \(K=[P,L]\). If \(K\) is nilpotent, its Sylow \(3\)-subgroup \(K_3\) must be nonabelian. Moreover, \(G/R(G)\cong S_5\), \(L/K\cong2.A_5\), \(K=K_3\times K_7\), \(K_7\) is a nonzero elementary abelian group, and \(K_3\) has exponent \(3\). \end{lemma} \begin{proof} All structural assertions apart from nonabelianness were proved in Theorem~\ref{thm:A5-nilpotent-K-structure} and Lemma~\ref{lem:nilpotent-K-prime-shapes}. If \(K_3\) were abelian, Lemma~\ref{lem:A5-large-two-dichotomy} would leave two possibilities for a Sylow \(2\)-subgroup of \(T\). The first is completely excluded by Lemma~\ref{lem:A5-central-Frattini-excluded}, and the second by Lemma~\ref{lem:A5-generic-square-excluded}. Thus the abelian branch is impossible. \end{proof}

\clearpage \subsection{Excluding nonabelian \(K_3\) and all nilpotent \(K\)}

The following argument excludes the remaining case in which \(K_3\) is nonabelian, without requiring a preliminary description of its derived subgroup or nilpotency class. The required order \(20\) and the forbidden order \(42\) are decisive.

\begin{lemma}\label{lem:nonabelian-three-unique-two-involution} Suppose that a finite group \(X\) has a normal nonabelian \(3\)-subgroup \(N\) and a nontrivial normal elementary abelian \(7\)-subgroup \(M\), with \([N,M]=1\) and \(42\notin\omega(X)\). Then every involution of \(X\) inverts \(M\), and every nontrivial Sylow \(2\)-subgroup of \(X\) has a unique involution. \end{lemma} \begin{proof} We first prove the automorphism fact used here. If an automorphism \(\alpha\) of order \(2\) of a finite group \(U\) has no nontrivial fixed points, then \(U\) is abelian and \(\alpha\) inverts it. Indeed, the map \[
 f:U\longrightarrow U,\qquad f(x)=x^{-1}x^\alpha
\] is injective: \(f(x)=f(y)\) implies \((yx^{-1})^\alpha=yx^{-1}\), and hence \(x=y\). Finiteness gives surjectivity. For every \(x\) we have \(f(x)^\alpha=(x^\alpha)^{-1}x=f(x)^{-1}\). Thus \(\alpha\) inverts every element of \(U\). Since inversion is an automorphism here, for any \(r,s\in U\) we obtain \(r^{-1}s^{-1}=(rs)^{-1}=s^{-1}r^{-1}\), so \(U\) is abelian.

Now let \(h\in X\) be any involution. It normalizes \(N\), and \(C_N(h)\ne1\): otherwise the fact just proved, applied to \(N\), would make \(N\) abelian. There is an element \(n\) of order \(3\) in \(C_N(h)\). If \(C_M(h)\ne1\), we could choose \(m\in C_M(h)\) of order \(7\). The elements \(n,m,h\) commute pairwise, since \([N,M]=1\), and have pairwise coprime orders \(3,7,2\). Then \(|nmh|=42\), a contradiction. Thus \(C_M(h)=1\). The action of an involution on the vector space \(M\) of characteristic \(7\) is diagonalizable with eigenvalues \(1,-1\). The absence of nonzero fixed vectors means that \(h|_M=-I\).

Let \(Q\le X\) be a nontrivial Sylow \(2\)-subgroup. Its action on \(M\) is faithful: a nontrivial kernel, being a \(2\)-group, would contain an involution, whereas every involution acts on the nonzero space \(M\) as \(-I\ne I\). But all involutions in \(Q\) have the same image \(-I\) in this faithful representation. Consequently, \(Q\) has a unique involution. \end{proof}

\Needspace{11\baselineskip} \begin{lemma}\label{lem:double-A5-central-four-two-involutions} Suppose that a finite group \(B\) contains \(D\cong2.A_5\) and that \(C_B(D)\) contains an element \(u\) of order \(4\). Then some Sylow \(2\)-subgroup of \(B\) has more than one involution. \end{lemma} \begin{proof} Let \(Q_D\cong Q_8\) be a Sylow \(2\)-subgroup of \(D\), and let \(z\) be the unique involution in \(D\). The group \(Q_D\langle u\rangle\) is a \(2\)-group because \([u,Q_D]=1\). Embed it in a Sylow \(2\)-subgroup \(Q_B\) of \(B\). Suppose that \(Q_B\) has a unique involution. Then \(u^2=z\), since both elements are involutions in \(Q_B\). Choose \(d\in Q_D\) of order \(4\); we have \(d^2=z\), so \[
 (du)^2=d^2u^2=z^2=1.
\] Moreover, \(du\notin\{1,z\}\). Otherwise \(u\) would belong to \(D\); but \(u\in C_B(D)\) and \(D\cap C_B(D)=Z(D)=\langle z\rangle\), whereas \(|u|=4\). Thus \(du\) is an involution distinct from \(z\). This contradicts the choice of \(Q_B\). \end{proof}

\Needspace{11\baselineskip} \begin{lemma}\label{lem:A5-nilpotent-K-excluded} Let \(G\) be a nonsolvable group with spectrum \(\Om=\omega(S_{10})\), and let \(\operatorname{Soc}(G/R(G))\cong A_5\), \(L=G^\infty\), and \(P=R(G)\cap L\). Then the group \(K=[P,L]\) is not nilpotent. In other words, the \(A_5\)-socle branch with nilpotent \(K\) is completely excluded. \end{lemma} \begin{proof} Suppose that \(K\) is nilpotent. By Theorem~\ref{thm:A5-nilpotent-K-structure} and Lemma~\ref{lem:A5-nilpotent-K-nonabelian-three}, we have \[
 K=N\times M,\qquad N=K_3\text{ nonabelian},
 \qquad M=K_7\ne1\text{ elementary abelian},
\] \[
 B=G/K,\qquad D=L/K\cong2.A_5,\qquad
 T=R(G)/K=C_B(D),
 \qquad \omega(T)\subseteq\{1,2,3,4,6\}.
\] The groups \(N,M\) are normal in \(G\) because they are characteristic Sylow subgroups of the characteristic subgroup \(K\). Since \(42\notin\Om\), Lemma~\ref{lem:nonabelian-three-unique-two-involution} shows that every Sylow \(2\)-subgroup of \(G\) has a unique involution. The order of \(K\) is odd, so the map \(G\to B\) maps Sylow \(2\)-subgroups of \(G\) isomorphically onto Sylow \(2\)-subgroups of \(B\). Thus the latter also have a unique involution.

On the other hand, the required order \(20\in\omega(G)\) is preserved on factoring out \(K\), since \(\gcd(20,|K|)=1\). Take \(b\in B\) of order \(20\) and put \(c=b^4\), so that \(|c|=5\). Every element of order \(5\) in \(B\) lies in \(D\). Indeed, \([B:DT]=2\), the group \(T\) has no elements of order \(5\), and \(DT\) is a central product of \(D,T\) with intersection of order \(2\). Lifting a Sylow \(5\)-subgroup through \(D\times T\to DT\) shows that all \(5\)-elements of this product lie in \(D\).

We have \[
 C_B(c)=\langle c\rangle\times T.
\] For completeness, we verify this equality, used earlier in Lemma~\ref{lem:A5-centralizer-five-T}. The action of \(B\) on \(D\) has kernel \(T\) and image \(\operatorname{Aut}(D)\cong S_5\). The image of \(c\) is a \(5\)-cycle whose centralizer in \(S_5\) is its cyclic subgroup. Therefore every element of \(C_B(c)\) lies in \(\langle c\rangle T\). The reverse inclusion is obvious, and \(\langle c\rangle\cap T=1\), giving the equality. Since \(b\in C_B(c)\) has order \(20\), its component in \(T\) has order \(4\). Thus \(4\in\omega(C_B(D))\).

Now Lemma~\ref{lem:double-A5-central-four-two-involutions} gives a Sylow \(2\)-subgroup of \(B\) with several involutions. This contradicts the first part of the proof. Consequently, nilpotent \(K\) is impossible. \end{proof}

\clearpage \subsection{Nonnilpotent \(K\): a required internal factor}

We retain the notation \(S=A_5\), \(L=G^\infty\), \(P=R(G)\cap L\), and \(K=[P,L]\). By Lemma~\ref{lem:A5-nilpotent-K-excluded}, the group \(K\) must be nonnilpotent. The central cover \(L/K\) is still isomorphic to \(A_5\) or \(2.A_5\); neither splitting of \(G\) over \(K\) nor a decomposition of \(K\) as a direct product of its Sylow subgroups is assumed here.

\begin{lemma}\label{lem:A5-K-internal-chief-exists} There is a chief factor \(V=N/M\) of \(G\) inside \(K\) on which \(K\) acts nontrivially. Choose the highest such factor in a fixed chief series inside \(K\). Put \[
 p=\operatorname{char}V,\quad C=C_G(V),\quad
 Y=G/C,\quad E=L C/C,\quad Q=P C/C,\quad U=K C/C.
\] Then \(V\) is a faithful simple \(Y\)-module, \(E=Y^\infty\), \(Q=R(E)\), and also \[
 E/Q\cong A_5,\qquad
 E/U\cong A_5\text{ or }2.A_5,
 \qquad Q/U\le Z(E/U),\quad |Q/U|\le2.
\] Moreover, \[
 1\ne U\text{ is nilpotent},\qquad
 U=[U,E],\qquad p\nmid|U|,
 \qquad \omega(V\rtimes Y)\subseteq\Om.
\] \end{lemma} \begin{proof} If \(K\) centralized all chief factors inside itself, the chief series inside \(K\) would be a central series for \(K\), so it would be nilpotent. Thus the required factor exists. Above the chosen \(N/M\), the group \(K\) acts trivially on every factor of the series inside \(K\). These factors give a central series for \(K/N\), so \(K/N\) is nilpotent. Since \(N\le C_G(N/M)\), the image \(U\) is a quotient of \(K/N\) and is also nilpotent. It is nontrivial by the choice of \(V\).

The action of \(L\) on \(V\) is nontrivial because \(K\le L\). By Lemma~\ref{lem:mixed-chief-action}, we have \(C\le R(G)\), \(Y/R(Y)\cong G/R(G)\), \(O_p(Y)=1\), and the stated inclusion of spectra. The images of \(L,P\) are, respectively, the perfect residual \(E\) and its radical \(Q\), with \(E/Q\cong A_5\). From \(K=[K,L]\), proved in Lemma~\ref{lem:remaining-commutator-localization}, we obtain \(U=[U,E]\).

The group \(U\) is normal in \(Y\). Its Sylow \(p\)-subgroup is characteristic in the nilpotent group \(U\) and therefore lies in \(O_p(Y)=1\). Consequently, \(p\nmid|U|\).

Finally, \(L/K\) is a perfect central cover of \(A_5\) with center \(P/K\). Since \(M(A_5)=C_2\), this cover is \(A_5\) or \(2.A_5\). The group \(E/U\) is its image and still has quotient \(A_5\). This gives both possibilities for \(E/U\) and centrality of \(Q/U\), which has order at most \(2\). \end{proof}

\Needspace{11\baselineskip} \begin{lemma}\label{lem:A5-internal-chief-Frattini-ranks} In the notation of the preceding lemma, every nontrivial Sylow \(q\)-subgroup \(U_q\) of \(U\) satisfies \([U_q,E]=U_q\). If \(d_q=\dim_{\mathbb F_q}U_q/\Phi(U_q)\), then \[
 d_2\ge4,\qquad d_3\ge4,\qquad
 d_5\ge2,\qquad d_7\ge4
\] for the primes that occur. Moreover, \(U_5,U_7\) cannot both be nontrivial. \end{lemma} \begin{proof} The Sylow subgroups of the nilpotent group \(U\) are characteristic, and the commutator with \(E\) is computed componentwise. Thus \([U,E]=U\) gives \([U_q,E]=U_q\). On \(W_q=U_q/\Phi(U_q)\ne0\) we obtain \([W_q,E]=W_q\); in particular, the action of \(E\) on \(W_q\) is nontrivial.

Choose \(a\in E\) of order \(5\) with nontrivial image in \(E/Q\cong A_5\). Such an element exists: the \(5\)-part of a preimage of an element of order \(5\) itself has order \(5\), since \(25\notin\omega(E)\). This element acts nontrivially on \(W_q\). Indeed, otherwise the normal kernel \(J_q\) of the action of \(E\) would contain \(a\). Its image in the simple group \(E/Q\) would be the whole group, so \(E=J_qQ\). The image of \(E\) on \(W_q\) would then coincide with the solvable image of \(Q\) while also being perfect. It would be trivial, contrary to \([W_q,E]=W_q\ne0\).

Consequently, \(5\mid|\operatorname{GL}_{d_q}(q)|\). For \(q=2,3,7\), the order of \(q\) modulo \(5\) is \(4\). The formula \[
 |\operatorname{GL}_d(q)|=
 q^{d(d-1)/2}\prod_{i=1}^d(q^i-1)
\] gives \(d_q\ge4\). For \(q=5\), dimension \(1\) is also impossible: \(\operatorname{GL}_1(5)\) is abelian and contains no element of order \(5\). Hence \(d_5\ge2\). Finally, if \(U_5,U_7\ne1\), they commute and give order \(35\) in \(U\le Y\), which is forbidden. \end{proof}

\Needspace{11\baselineskip} \begin{lemma}\label{lem:A5-internal-five-seven-two-kernel} If the characteristic of the chosen internal chief factor is \(p=5\) or \(p=7\), then \[
 U=U_2\text{ nonabelian},\qquad
 |U|\ge32,\qquad Q\text{ --- }2\text{-group},
 \qquad \dim_{\mathbb F_2}U/\Phi(U)\ge4.
\] For an element \(a\in E\) of order \(5\) nontrivial modulo \(Q\), we have \([U,a]\ne1\). Its minimal polynomial on \(V\) divides \(1+x+x^2+x^3+x^4\), and \(C_U(a)\ne1\). \end{lemma} \begin{proof} Let \(r\) be the other prime in the pair \(\{5,7\}\). Every nonidentity element of \(U_r\) has order \(r\) because \(25,49\notin\omega(Y)\). It fixes no nonzero vector in \(V\), since this would give order \(35\). Thus \(U_r\) contains no \(C_r^2\): over an algebraic closure, every linear character of \(C_r^2\) has a nontrivial kernel. If a group of exponent \(r\) is noncyclic, a central element of order \(r\) and an element outside the subgroup it generates give \(C_r^2\). Therefore \(U_r=1\) or \(U_r\cong C_r\). In the second case the perfect group \(E\) centralizes \(U_r\), since \(\operatorname{Aut}(C_r)\) is abelian; this contradicts \([U_r,E]=U_r\). Hence \(U_r=1\). The subgroup \(U_p\) is trivial by the preceding reduction.

Choose \(a\) of order \(5\) outside \(Q\). If \([U_3,a]\ne1\), apply Lemma~\ref{lem:HH-prime} to \(U_3\rtimes\langle a\rangle\) on the faithful module \(V\). The characteristic is \(p\ne3\), and the exceptional case of that lemma is impossible because the normal subgroup has odd order. Thus the minimal polynomial of \(a\) is \(x^5-1\), giving order \(5p\in\{25,35\}\) in \(V\rtimes Y\). This is forbidden. Hence \(a\) centralizes \(U_3\). The centralizer \(C_E(U_3)\) is normal in \(E\) and contains an element outside \(Q\). Its image in \(A_5\) is the whole group; as in the preceding lemma, the perfect image of \(E\) on \(U_3\) coincides with the solvable image of \(Q\) and is therefore trivial. From \([U_3,E]=U_3\) we obtain \(U_3=1\).

Only \(U=U_2\ne1\) remains. The element \(a\) acts nontrivially on \(U\), by the proof of the preceding lemma. If \(U\) were abelian, the exceptional case of Lemma~\ref{lem:HH-prime} would again be impossible, giving the forbidden order \(5p\). Therefore \(U\) is nonabelian. Its Frattini rank is at least \(4\) and \(\Phi(U)\ne1\), so \(|U|\ge2^5=32\). Since \(|Q/U|\le2\), the group \(Q\) is also a \(2\)-group.

The norm \(1+a+\cdots+a^4\) on \(V\) is zero: otherwise it gives the same order \(5p\). This proves the assertion about the minimal polynomial. Applying Lemma~\ref{lem:HH-prime} now to the nonabelian group \(U\) gives \(C_U(a)\ne1\). Here the exception is indeed admissible: \(5\) is a Fermat prime. \end{proof}

\Needspace{11\baselineskip} \begin{lemma}\label{lem:A5-internal-four-characteristics} For nonnilpotent \(K\), the required internal chief factor can be chosen so that its image \(U\) satisfies the following exhaustive table of necessary conditions: \[
\begin{array}{c|l}
 p&\text{possible structure }U\ne1\\ \hline
 2&U=U_3\times U_5\times U_7,\quad U_5=1\text{ or }U_7=1;\\
 3&U=U_2\times U_5\times U_7,\quad U_5=1\text{ or }U_7=1;\\
 5,7&U=U_2\text{ nonabelian},\quad |U_2|\ge32.
\end{array}
\] In every row, \([U_q,E]=U_q\) for each prime that occurs, and the Frattini-rank bounds from Lemma~\ref{lem:A5-internal-chief-Frattini-ranks} apply. For \(p=2\), if \(U_5\ne1\) or \(U_7\ne1\), then \(\exp U_3\mid3\). For \(p=3\), if \(U_7\ne1\), then \(\exp U_2\mid2\). \end{lemma} \begin{proof} The table follows from nilpotence of \(U\), the condition \(p\nmid|U|\), the absence of simultaneous divisors \(5,7\), and the preceding lemma. For \(p=2\), an element of order \(9\) in \(U_3\) together with an element of order \(5\) or \(7\) would give \(45\) or \(63\). Order \(27\) is also forbidden, so \(\exp U_3\mid3\) in the stated cases. For \(p=3\), an element of order \(4\) in \(U_2\) would commute with an element of order \(7\) in \(U_7\), giving \(28\). Consequently, \(U_2\) has exponent at most \(2\). \end{proof}

\clearpage \subsection{A critical \(2\)-section in characteristics \(5,7\)}

We retain the notation of Lemma~\ref{lem:A5-K-internal-chief-exists} and assume \(p\in\{5,7\}\). By Lemma~\ref{lem:A5-internal-five-seven-two-kernel}, the groups \(U,Q\) are \(2\)-groups, \(U\) is nonabelian, \(E/Q\cong A_5\), and \(V\) is a faithful module of characteristic \(p\) with \(\omega(V\rtimes Y)\subseteq\Om\). When working with simple types, we extend the field to an algebraic closure; vanishing of norms and dimensions of kernels of linear operators are preserved.

\begin{lemma}\label{lem:A5-five-seven-invariant-two-types} Let \(A\trianglelefteq E\) be any \(2\)-subgroup. All absolutely simple types in the restriction of \(V\) to \(A\) are \(E\)-invariant. In particular, \(E\) centralizes \(Z(U)\), and \(\exp Z(U)\le4\). \end{lemma} \begin{proof} The restriction to \(A\) is completely reducible because \(p\ne2\). Choose \(a\in E\) of order \(5\) outside \(Q\). If \(a\) cyclically permuted five distinct isotypic spaces, then for a nonzero vector \(v\) in one of them the sum \(v+av+a^2v+a^3v+a^4v\) would be nonzero. This gives order \(5p\in\{25,35\}\) in \(V\rtimes Y\), which is forbidden. Thus \(a\) fixes all simple types.

The kernel \(J\) of the permutation action of \(E\) on these types is normal in \(E\) and contains \(a\notin Q\). Simplicity of \(E/Q\cong A_5\) gives \(E=JQ\). Therefore the permutation image of \(E\) coincides with the solvable image of \(Q\), while also being perfect. Hence it is trivial.

Apply this result to the abelian normal subgroup \(A=Z(U)\). Its simple types are linear characters, and the elements of \(Z(U)\) act as scalars on every weight space. Since \(E\) preserves each such space, it commutes on \(V\) with every element of \(Z(U)\). Faithfulness of \(V\) gives \([E,Z(U)]=1\). If \(Z(U)\) contained an element of order \(8\), it would commute with \(a\) and give order \(40\) in \(E\), which is forbidden. Hence \(\exp Z(U)\le4\). \end{proof}

\Needspace{11\baselineskip} \begin{lemma}\label{lem:A5-five-seven-critical-section} There exist a normal \(2\)-subgroup \(C\) of \(E\), \(C\le U\), and a subgroup \(J\le C\) normal in \(E\) such that, for \(E_0=E/J\), \(Q_0=Q/J\), and \(D=C/J\), the following conditions hold: \[
 E_0=E_0',\qquad Q_0\text{ --- }2\text{-group},
 \qquad E_0/Q_0\cong A_5,
\] \[
 D\text{ nonabelian of class }2,\qquad
 Z(D)\cong C_2\text{ or }C_4,\qquad
 Z(D)\le Z(E_0),\qquad D'\cong C_2,
\] \[
 D/Z(D)\cong C_2^{\,2n},\quad n\ge2,
 \qquad |D|\ge32.
\] The action of \(E_0\) on \(D/Z(D)\) is nontrivial, and its kernel is contained in \(Q_0\). There is a nonzero \(E\)-invariant isotypic space \(W\) in the restriction of \(V\) to \(C\) on which \(J\) is trivial and \(D\) acts faithfully and homogeneously. For \(|Z(D)|=2\), the group \(D\) is extraspecial. \end{lemma} \begin{proof} An element \(a\) of order \(5\) outside \(Q\) acts nontrivially on \(U\) by Lemma~\ref{lem:A5-internal-chief-Frattini-ranks}. Take a characteristic subgroup \(C\) of \(U\) using the stated form of Thompson's theorem, Lemma~\ref{lem:critical-two-external}. It has class at most \(2\), the quotient \(C/Z(C)\) is elementary abelian, and \(a\) acts nontrivially on \(C\).

By the preceding lemma, all simple types in the restriction of \(V\) to \(C\) are \(E\)-invariant. Denote their kernels by \(J_i\). They are normal in \(E\), and \(\bigcap_iJ_i=1\) because the action of \(C\) on \(V\) is faithful. If \(E\) centralized every \(C/J_i\), it would centralize \(C\), contrary to the choice of \(a\). Choose \(J=J_i\) for which the action on \(D=C/J\) is nontrivial, and the corresponding isotypic space \(W\). On \(W\) the group \(D\) acts faithfully as a sum of copies of one simple type, while \(J\) acts trivially.

The center of \(D\) acts by scalars on this simple type, by Schur's lemma; faithfulness shows that \(Z(D)\) is cyclic. The perfect group \(E_0\) centralizes this center because \(\operatorname{Aut}(Z(D))\) is abelian. If the center contained an element of order \(8\), together with an element of order \(5\) of \(E_0\) it would give the forbidden order \(40\). Consequently, \(|Z(D)|=2\) or \(4\). The group \(D\) is nonabelian: otherwise it would itself be cyclic and \(E_0\) would act trivially on it.

The quotient \(D/Z(D)\) is elementary abelian because it is a quotient of \(C/Z(C)\). For \(x,y\in D\) we have \([x,y]^2=[x^2,y]=1\), since \(D\) has class \(2\) and \(x^2\in Z(D)\). Thus the nontrivial derived subgroup \(D'\) is the unique subgroup of order \(2\) in the cyclic center. The commutator defines a nondegenerate alternating form on \(D/Z(D)\) with values in \(D'\cong\mathbb F_2\): its radical is trivial by the definition of the center. Hence this quotient has even dimension \(2n\).

If \(E_0\) acted trivially on \(D/Z(D)\), then, in view of \([E_0,Z(D)]=1\), its image on \(D\) would lie in the \(2\)-group \(\operatorname{Hom}(D/Z(D),Z(D))\). Indeed, an automorphism of this kind has the form \(x\mapsto xf(xZ(D))\). A perfect image in a \(2\)-group is trivial, contrary to the choice of \(J\). Thus the action on \(D/Z(D)\) is nontrivial. Its normal kernel lies in \(Q_0\): otherwise the image of the kernel in \(A_5\) would be the whole group, and the perfect image of \(E_0\) would coincide with the solvable image of \(Q_0\). Consequently, an element of order \(5\) outside \(Q_0\) acts nontrivially on \(D/Z(D)\). As in Lemma~\ref{lem:A5-internal-chief-Frattini-ranks}, this gives \(2n\ge4\), and hence \(|D|\ge32\).

If \(Z(D)\cong C_2\), then \(D'=Z(D)\), and the squares of elements lie in \(Z(D)\). Hence \(\Phi(D)=D^2D'=Z(D)\), which means precisely that \(D\) is extraspecial. \end{proof}

\Needspace{11\baselineskip} \begin{lemma}\label{lem:A5-seven-internal-needs-eight} Let \(p=7\), and choose \(E_0,Q_0,D,W\) as in the preceding lemma. For every \(2\)-element \(g\in E_0\setminus Q_0\), the coset \(Dg\) contains an element of order \(8\). In particular, \(8\in\omega(E_0)\subseteq\omega(E)\), and the subbranch \(8\notin\omega(E)\) is impossible. \end{lemma} \begin{proof} Let \(z\) be the unique involution in \(Z(D)\). It is central in \(E_0\) and acts on \(W\) as \(-I\): on the faithful simple type of \(D\), a central involution acts as a nontrivial scalar. Write \(\rho\) for the representation of \(E_0\) on \(W\); it is faithful on \(D\).

Every element \(h\in E_0\) of order \(4\) acts on \(W\) without nonzero fixed vectors. Indeed, lift \(h\) to a \(2\)-element \(\widetilde h\in E\) through the normal \(2\)-subgroup \(J\). Its order is \(4\) or \(8\): it is divisible by \(4\), and \(16\notin\Om\). A nonzero fixed vector of \(h\) in \(W\) gives a fixed vector of \(\widetilde h\) in the extended space \(V\). The dimension of the kernel of \(\widetilde h-I\) is unchanged by field extension. Thus order \(28\) or \(56\) would occur in \(V\rtimes E\), and both are forbidden.

For any such \(h\) we have \(\rho(h)^2=-I\). For \(h^2=z\) this is obvious. For \(h^2\ne z\) both elements \(h,zh\) have order \(4\). On the \((+1)\)-eigenspace of \(\rho(h)^2\), the operator \(\rho(h)\) has only eigenvalues \(1,-1\). The eigenvalue \(1\) is forbidden for \(h\), and the eigenvalue \(-1\) gives the eigenvalue \(1\) for \(zh\), also forbidden. Thus this \((+1)\)-eigenspace is zero and \(\rho(h)^2=-I\). Hence the square of the matrix of every element of \(E_0\) of order at most \(4\) is scalar.

Now take a \(2\)-element \(g\in E_0\setminus Q_0\) and suppose that \(Dg\) contains no elements of order \(8\). The group \(D\langle g\rangle\) is a \(2\)-group because \(D\) is normal. Its elements have order at most \(8\), since \(16\notin\Om\). By assumption every element of \(Dg\) has order \(2\) or \(4\). Thus, for every \(d\in D\), the matrices \(A=\rho(g)\) and \(B=\rho(d)\) satisfy \[
 A^2\text{ is scalar},\quad
 B^2\text{ is scalar},\quad
 (AB)^2\text{ is scalar}.
\] Here the second condition follows from \(d^2\in Z(D)\), and the third from \(gd\in Dg\). Modulo scalar matrices, the images of \(A,B,AB\) have squares \(1\), so the images of \(A,B\) commute. In other words, \(\rho([g,d])\) is scalar. Since \([g,d]\in D\) and the representation is faithful on \(D\), scalarity means \([g,d]\in Z(D)\). Thus \(g\) centralizes \(D/Z(D)\). This contradicts the fact that the kernel of the action of \(E_0\) on \(D/Z(D)\) lies in \(Q_0\), whereas \(g\notin Q_0\). Consequently, \(Dg\) contains an element of order \(8\), as required. \end{proof}

The cases left for further analysis are characteristic \(5\) with the stated critical section, and characteristic \(7\), in which order \(8\) must additionally occur in \(E\) itself. This argument does not exclude the required internal factors of characteristics \(2,3\).

\clearpage \subsection{Excluding an extraspecial section of order \(32\) for \(p=7\)}

In this subsection, \(E_0,Q_0,D,W\) have the meanings given in Lemma~\ref{lem:A5-five-seven-critical-section}, and the characteristic is \(7\). We exclude the case \(|D|=32\), \(Z(D)=C_2\). The multiplicity of the simple type of \(D\) on \(W\) is not assumed to be one.

\begin{lemma}\label{lem:extraspecial32-five-eigenvalues} Suppose that \(D\) is an extraspecial group of order \(32\) admitting an automorphism \(a\) of order \(5\). Then \(D\cong Q_8\mathbin{\circ}D_8\) has minus type. Its unique faithful absolutely simple type in characteristic \(7\) has dimension \(4\) and preserves a nondegenerate alternating form. An operator \(A\) implementing \(a\) can be chosen to preserve this form and satisfy \(A^5=I\). Its eigenvalues are \[
 \zeta,\ \zeta^2,\ \zeta^3,\ \zeta^4,
 \qquad |\zeta|=5,
\] each with multiplicity \(1\). \end{lemma} \begin{proof} Squares define a nondegenerate quadratic form \(q\) on \(D/Z(D)\cong\mathbb F_2^4\) whose polarization is the commutator. Its two normal forms are \[
 q_+(v)=v_1v_2+v_3v_4,\qquad
 q_-(v)=v_1+v_2+v_1v_2+v_3v_4.
\] These correspond to the central products \(D_8\circ D_8\) and \(Q_8\circ D_8\): in a symplectic basis, the squares of the generators specify the type of each factor of order \(8\). Direct counting of isometries gives \[
 |O(q_+)|=72,\qquad |O(q_-)|=120.
\] Appendix~\ref{app:D32-certificate} provides a finite certificate for this count: it enumerates all ordered bases of \(\mathbb F_2^4\) and checks preservation of \(q\) on all \(16\) vectors. The action of \(a\) on \(D/Z(D)\) is nontrivial: the kernel of the action of automorphisms on this quotient consists of inner automorphisms and is a \(2\)-group. Hence \(5\mid|O(q)|\), leaving only minus type.

For the faithful simple type, use a model over \(\mathbb F_7\). Put \[
 i=\begin{pmatrix}0&1\\6&0\end{pmatrix},\quad
 j=\begin{pmatrix}2&3\\3&5\end{pmatrix},\quad
 x=\begin{pmatrix}0&1\\1&0\end{pmatrix},\quad
 y=\begin{pmatrix}1&0\\0&6\end{pmatrix}.
\] Then \(i,j\) generate \(Q_8\), and \(x,y\) generate \(D_8\). The matrices \[
 i\otimes I_2,\quad j\otimes I_2,\quad
 I_2\otimes x,\quad I_2\otimes y
\] generate \(D\) of order \(32\). They preserve the alternating form with matrix \(B=i\otimes I_2\). This is a faithful irreducible type of degree \(4\); its uniqueness follows from the usual character count for an extraspecial group given in the proof of Lemma~\ref{lem:extraspecial-eight-coset}.

An operator implementing \(a\) can initially be normalized by \(A^5=I\). Uniqueness of the invariant bilinear form gives \(A^tBA=\mu B\), where \(\mu^5=1\). Multiplying \(A\) by a suitable \(5\)th root of unity makes \(\mu=1\) while preserving \(A^5=I\).

We prove the assertion about the eigenvalues of \(A\) in a complex representation of \(D\rtimes C_5\). Returning to characteristic \(7\) preserves eigenvalue multiplicities because \(7\nmid32\cdot5\). The symplectic extension of the type to \(C_5\) is unique: two extensions differ by a \(5\)th root of unity, and preservation of the form requires the square of this root to be \(1\). Thus the complex symplectic extension corresponds to the chosen \(A\) in characteristic \(7\).

The action of \(a\) on \(\mathbb F_2^4\) has no nonzero fixed vectors: a nontrivial irreducible \(C_5\)-module over \(\mathbb F_2\) has dimension \(4\). Matrices representing \(D/Z(D)\) form a basis of the endomorphism algebra of the four-dimensional type. For \(1\le t\le4\), conjugation by \(A^t\) permutes this basis, possibly with signs, and fixes only the identity matrix. Consequently, \[
 |\operatorname{tr}(A^t)|^2
 =\operatorname{tr}(\operatorname{Ad}(A^t))=1.
\] If \(n_r\) is the multiplicity of \(\zeta^r\) in \(A\), character orthogonality for \(C_5\) gives \[
 \sum_{r=0}^4 n_r^2
 =\frac15\sum_{t=0}^4|\operatorname{tr}(A^t)|^2
 =\frac{16+4}{5}=4,
 \qquad \sum_{r=0}^4 n_r=4.
\] Thus the four eigenvalues are distinct. Preservation of the alternating form makes their set invariant under inversion. The only way to choose four distinct \(5\)th roots of unity with this property is to omit \(1\). This proves the required list. \end{proof}

\Needspace{11\baselineskip} \begin{lemma}\label{lem:extraspecial32-normalizer-fixed-four} Let \(k\) be an algebraically closed field of characteristic \(7\), and realize \(D=Q_8\circ D_8\) on \(k^4\) by the model of the preceding lemma. Suppose that \(F\le\operatorname{Sp}_4(k)\) contains and normalizes \(D\) and has image \(A_5\) on \(D/Z(D)\). Then \(F\) contains an element of order \(4\) with a two-dimensional fixed space on \(k^4\). \end{lemma} \begin{proof} Write \(N\) for the normalizer of \(D\) in \(\operatorname{Sp}_4(k)\). Its image on \(D/Z(D)\) lies in \(O(q_-)\). The kernel of this action is \(D\). Indeed, an automorphism of \(D\) trivial on \(D/Z(D)\) has the form \(d\mapsto dz^{f(dZ(D))}\), where \(f\) is linear, and is induced by conjugation by an element of \(D\), by nondegeneracy of the commutator form. After this adjustment, an element of \(N\) centralizes \(D\), so is scalar by irreducibility. A symplectic scalar is \(\pm I\) and already belongs to \(D\).

The group \(O(q_-)\) acts faithfully on the five nonzero singular vectors, which span \(\mathbb F_2^4\). Together with the order count \(120\), this identifies it with \(S_5\). Thus the image of \(F\) is the unique subgroup of index \(2\), namely \(A_5\).

Take commuting matrices \[
 R=i\otimes I_2,\qquad S=I_2\otimes(xy),
 \qquad R^2=S^2=-I.
\] The matrices \(H_R=2(I+R)\) and \(H_S=2(I+S)\) normalize \(D\) and preserve the form \(B\). To verify normalization, it suffices to check the generators: an element of \(D\) commuting with \(R\) is fixed, while one anticommuting with it maps to its product with \(R\), possibly with a central sign; the same holds for \(S\). Preservation of the form follows from \[
 (I+R)^tB(I+R)=2B,\qquad 2^2\cdot2=1
 \quad\text{in characteristic }7,
\] and the analogous equality for \(S\). On \(D/Z(D)\) they induce the orthogonal transvections \(v\mapsto v+b(v,r)r\) and \(v\mapsto v+b(v,s)s\), where \(q(r)=q(s)=1\). Each interchanges two of the five singular vectors and fixes the other three; this is checked on the displayed forms \(q_-\). Thus the product \[
 T=H_RH_S=4(I+R)(I+S)
\] has an even image in \(S_5\). Since the kernel in \(N\to S_5\) is \(D\) and \(D\le F\), we obtain \(T\in F\).

Direct multiplication gives \[
 H_R^2=R,\qquad H_S^2=S,\qquad
 T^2=RS\ne I,\qquad T^4=I.
\] To find the fixed space, simultaneously diagonalize \(R,S\). The four pairs of eigenvalues \((\varepsilon\iota,\delta\iota)\), \(\varepsilon,\delta\in\{1,-1\}\), \(\iota^2=-1\) have multiplicity \(1\): \(i\) and \(xy\) each have both eigenvalues once, and the representation is a tensor product. On the two pairs with \(\varepsilon=-\delta\), the matrix \(T\) has eigenvalue \(4(1+\iota)(1-\iota)=8=1\). On the other two its eigenvalues are \(\iota,-\iota\). Consequently, \(|T|=4\) and \(\dim C_{k^4}(T)=2\). \end{proof}

\Needspace{11\baselineskip} \begin{lemma}\label{lem:A5-seven-no-critical32} In characteristic \(p=7\), the critical section \(D\) of Lemma~\ref{lem:A5-five-seven-critical-section} cannot be extraspecial of order \(32\), regardless of the multiplicity of its simple type on \(W\). Thus only the following possibilities remain: \[
 |Z(D)|=2,\ |D|\ge128,
 \qquad\text{or}\qquad
 |Z(D)|=4,\ |D|\ge64.
\] \end{lemma} \begin{proof} Suppose that \(|D|=32\) and \(Z(D)=C_2\). Choose \(a\in E_0\) of order \(5\) outside \(Q_0\). Its action on \(D/Z(D)\) is nontrivial by Lemma~\ref{lem:A5-five-seven-critical-section}. Thus Lemma~\ref{lem:extraspecial32-five-eigenvalues} applies. Let \(H\) be the four-dimensional faithful simple type of \(D\). Write the homogeneous \(D\)-module \(W\) as \(H\otimes M\).

For every \(e\in E_0\) choose an intertwining operator \(A_e\) on \(H\) implementing conjugation of \(D\) by \(e\). Then the action of \(e\) on \(W\) has the form \(A_e\otimes S_e\), since the commutant of \(D\) on \(W\) is \(I_H\otimes\operatorname{End}(M)\). The operators \(A_e\) are determined up to scalars. Thus the classes \([S_e]\) define a homomorphism \[
 E_0\longrightarrow\operatorname{PGL}(M).
\] For the chosen \(a\) take \(A_a\) symplectic with \(A_a^5=I\). Then \(S_a^5=I\), and \(A_a\) has all four nontrivial \(5\)th roots of unity as eigenvalues. If \(S_a\) had a nontrivial eigenvalue \(\zeta^t\), it would multiply with \(\zeta^{-t}\) from \(A_a\) to give \(1\). This would give a nonzero fixed vector of \(a\) on \(W\) and hence the forbidden order \(35\). For passage to the original group, lift \(a\) to an element of order \(5\) through the normal \(2\)-subgroup \(J\) and use preservation of the kernel dimension under field extension. Consequently, all eigenvalues of \(S_a\) are \(1\). Since the characteristic is \(7\ne5\), we obtain \(S_a=I\).

The kernel of the projective action on \(M\) contains \(a\notin Q_0\). By simplicity of \(E_0/Q_0\cong A_5\), the projective image of \(E_0\) coincides with the image of \(Q_0\). It is both perfect and a \(2\)-group, so is trivial. Thus all \(S_e\) are scalar. After absorbing these scalars into \(A_e\), the representation of \(W\) decomposes into a direct sum of identical four-dimensional \(E_0\)-modules. These are genuine, rather than merely projective, representations: equality of operators on \(H\otimes M\) gives \(A_{ef}=A_eA_f\).

Consider one such module \(H\). Its \(D\)-invariant alternating form, unique up to scalar, is preserved by \(E_0\) up to a multiplier \(\lambda_e\). The map \(e\mapsto\lambda_e\) is a linear character of a perfect group, so \(\lambda_e=1\). Hence the image \(F\) of \(E_0\) on \(H\) lies in \(\operatorname{Sp}_4(k)\) and contains and normalizes \(D\). Its image on \(D/Z(D)\) has a quotient \(A_5\), since the kernel of this action in \(E_0\) lies in \(Q_0\). On the other hand, the image lies in \(O_4^-(2)\cong S_5\). Its order is therefore \(60\) or \(120\); the second case is impossible by perfection. Thus this image is \(A_5\).

By Lemma~\ref{lem:extraspecial32-normalizer-fixed-four}, the group \(F\) contains an element \(T\) of order \(4\) with a nonzero fixed vector. Lift it to a \(2\)-element \(t\in E_0\). Its order is \(4\) or \(8\), since \(16\notin\Om\). It fixes a nonzero vector in \(H\), hence also in \(W\), and gives order \(28\) or \(56\). Lifting further through \(J\) in \(E\) preserves one of these two orders of a \(2\)-element, since order \(16\) is forbidden. This is a contradiction.

For \(|Z(D)|=2\), the order of \(D\) has the form \(2^{1+2n}\), \(n\ge2\), so after excluding \(32\), the next possible order is at least \(128\). For \(|Z(D)|=4\) the previous bound \(\dim D/Z(D)\ge4\) gives \(|D|\ge64\). \end{proof}

\clearpage \subsection{Excluding all extraspecial sections for \(p=7\)}

In this subsection, the assumptions of Lemma~\ref{lem:A5-five-seven-critical-section} hold, with \(p=7\) and \(Z(D)\cong C_2\). Thus \(D\) is extraspecial of order \(2^{1+2n}\), \(n\ge2\). Write \(\mathcal V=D/Z(D)\) for the space over \(\mathbb F_2\) and \(q\) for its quadratic form defined by squares of elements of \(D\). The scalar field of the original module \(W\) is extended to an algebraic closure \(k\) of characteristic \(7\).

\begin{lemma}\label{lem:A5-seven-extraspecial-five-rank-four} For every element \(a\in E_0\) of order \(5\) outside \(Q_0\), we have \[
 \dim_{\mathbb F_2}[\mathcal V,a]=4.
\] The module \(W\) is a direct sum of identical \(E_0\)-modules \(H\) of dimension \(2^n\), each restricting faithfully and irreducibly to \(D\). On \(H\) the group \(E_0\) preserves a nondegenerate bilinear form, symmetric or alternating according to the type of \(D\). The eigenvalues of \(a\) on \(H\) are the four nontrivial \(5\)th roots of unity, each with multiplicity \(2^{n-2}\). \end{lemma} \begin{proof} The action of \(a\) on \(\mathcal V\) is nontrivial by Lemma~\ref{lem:A5-five-seven-critical-section}. Since \(|a|=5\) is coprime to \(2\), the space is the orthogonal direct sum of its fixed and moved parts: \[
 \mathcal V=\mathcal V_0\perp\mathcal V_1,
 \qquad \mathcal V_0=C_{\mathcal V}(a),
 \quad \mathcal V_1=[\mathcal V,a].
\] Indeed, the fixed part is orthogonal to the image of \(a-1\) for the polar form, and averaging over \(C_5\) gives their direct sum. Both parts are nondegenerate. The nontrivial simple \(\mathbb F_2C_5\)-modules have dimension \(4\), so \(\dim\mathcal V_1=4r\) and \(r\ge1\).

The preimages of these parts give a central product \(D=D_0\circ D_1\), where \(D_1\) is extraspecial of order \(2^{1+4r}\). For \(\mathcal V_0=0\), interpret \(D_0\) as the common center. The automorphism \(a\) centralizes \(D_0\): it is trivial on \(D_0/Z(D_0)\), and the kernel of this action on automorphisms of an extraspecial group is a \(2\)-group. The faithful simple type of \(D\) has the form \(H_0\otimes H_1\), where \(H_i\) are the faithful simple types of the factors. On \(H_1\) choose an operator \(A_1\) implementing \(a\), with \(A_1^5=I\), preserving the corresponding invariant bilinear form. This is possible by the same normalization as in Lemma~\ref{lem:extraspecial32-five-eigenvalues}. On \(H_0\otimes H_1\) choose \(I\otimes A_1\).

Let \(N=\dim H_1=2^{2r}\), and let \(m_0,\ldots,m_4\) be the multiplicities of the \(5\)th roots of unity as eigenvalues of \(A_1\). As in the proof of Lemma~\ref{lem:extraspecial32-five-eigenvalues}, in the complex extension of the type to \(D_1\rtimes C_5\) we have \(|\operatorname{tr}(A_1^t)|^2=1\) for \(1\le t\le4\): conjugation on the endomorphism algebra fixes only the identity in the matrix basis \(D_1/Z(D_1)\). Therefore \[
 \sum_{j=0}^4m_j=N,\qquad
 \sum_{j=0}^4m_j^2=\frac{N^2+4}{5}.
\] These multiplicities are preserved in characteristic \(7\), which is coprime to \(|D_1\rtimes C_5|\).

At least one of \(m_j\) is zero. Otherwise, writing \(W=(H_0\otimes H_1)\otimes M\), for every eigenvalue of \(a\) on \(M\) there would be its inverse on \(H_1\) and a nonzero fixed vector on \(W\). This gives the forbidden order \(35\). With at most four nonzero \(m_j\), Cauchy's inequality gives \[
 \frac{N^2+4}{5}=\sum m_j^2\ge\frac{N^2}{4},
 \qquad N^2\le16.
\] Consequently, \(N=4\) and \(r=1\), and the moved part has dimension \(4\). By Lemma~\ref{lem:extraspecial32-five-eigenvalues}, exactly the four nontrivial \(5\)th roots of unity occur on \(H_1\).

It remains to justify the multiplicities and the genuine action of \(E_0\) on \(H\). The action on the multiplicity space \(M\) is projective, as in the proof of Lemma~\ref{lem:A5-seven-no-critical32}. For the chosen \(a\) it is represented by an operator \(S_a\), with \(S_a^5=I\). The exclusion of \(35\) and the presence of all four nontrivial roots on \(H_1\) force \(S_a=I\). The kernel of the projective action contains \(a\notin Q_0\), so its image coincides with the image of the \(2\)-group \(Q_0\) and is also perfect. It is trivial. Hence all operators on \(M\) are scalar; absorbing them into the intertwining operators on \(H\) gives identical genuine representations of \(E_0\) on each copy of \(H\). The unique \(D\)-invariant bilinear form on \(H\) is preserved by \(E_0\) up to a multiplier that is a linear character. Perfection makes this multiplier \(1\). Finally, \(a\) acts trivially on \(H_0\) of dimension \(2^{n-2}\), giving the stated multiplicities. \end{proof}

\Needspace{11\baselineskip} \begin{lemma}\label{lem:A5-seven-extraspecial-action-simple} The image of \(E_0\) on \(\mathcal V=D/Z(D)\) is isomorphic to \(A_5\); the group \(Q_0\) acts trivially on \(\mathcal V\). In particular, the action on \(\mathcal V\) factors through \(E_0/Q_0\cong A_5\). \end{lemma} \begin{proof} Work on one module \(H\) of the preceding lemma, with representation \(\rho\). Let \(z\) be the involution in \(Z(D)\); we have \(\rho(z)=-I\). By the proof of Lemma~\ref{lem:A5-seven-internal-needs-eight}, every element \(h\in E_0\) of order \(4\) satisfies \(\rho(h)^2=-I\). Here we use faithfulness on \(D\) and the exclusions \(28,56\), applicable after lifting through \(J\). If \(|h|=8\), then \(|h^2|=4\), so \(\rho(h)^4=-I\). Thus the action of any \(2\)-element on \(\mathcal V\) has order at most \(4\).

Suppose this order is \(4\) for some \(g\in E_0\). We may choose \(g\) to be a \(2\)-element. It has order \(8\): for \(|g|=4\) its square would be scalar on \(H\) and centralize \(D\), giving action order at most \(2\). Every element of \(Dg\) has the same action on \(\mathcal V\), since inner automorphisms of \(D\) act trivially on \(\mathcal V\). Thus all elements of \(Dg\) have order \(8\). They fix no nonzero vectors of \(H\), since otherwise \(56\) would occur.

Put \(d=\dim H=2^n\), with \(\varepsilon=1\) for a symmetric form and \(\varepsilon=-1\) for an alternating form. For \(0\le j\le7\) consider the sums in \(k\) \[
 S_j=\sum_{x\in Dg}\operatorname{tr}(\rho(x)^j).
\] We have \(S_0=|D|d\). For odd \(j\) the map \(x\mapsto x^j\) maps \(Dg\) bijectively onto \(Dg^j\). Since \(\sum_{u\in D}\rho(u)=0\), we obtain \(S_1=S_3=S_5=S_7=0\). Moreover, \(\rho(x)^4=-I\) for every \(x\in Dg\), whence \(S_4=-|D|d\).

Orthogonality of the matrix coefficients of the irreducible \(D\)-module gives, if \(B^t=\varepsilon B\) is the matrix of the form, \[
 \sum_{u\in D}\rho(u)A\rho(u)
 =\frac{\varepsilon|D|}{d}B^{-1}A^tB.
\] This is the same calculation as in the proof of Lemma~\ref{lem:extraspecial-eight-coset}; it applies directly in characteristic \(7\) because \(d|D|\ne0\) in this field. Substituting \(A=\rho(g)\), multiplying by \(A\), and taking traces, we use \(A^tBA=B\) to obtain \(S_2=\varepsilon|D|\). The cubing bijection maps \(Dg\) onto \(Dg^3\), and the same calculation for \(g^3\) gives \(S_6=\varepsilon|D|\). Consequently, \[
 \sum_{j=0}^7S_j=2\varepsilon|D|\ne0
 \quad\text{in characteristic }7.
\] But the absence of nonzero fixed vectors for every \(x\in Dg\) means \(I+\rho(x)+\cdots+\rho(x)^7=0\): this operator is \(8\) times the projection onto the fixed space, since \(7\nmid8\). Summing traces gives \(\sum_jS_j=0\), a contradiction. Thus the image of the action contains no elements of order \(4\).

Denote this image by \(T\) and the image of \(Q_0\) by \(A\). The group \(T\) is perfect, \(A\) is a normal \(2\)-subgroup, and \(T/A\cong A_5\) because the kernel of the action lies in \(Q_0\). A Sylow \(2\)-subgroup of \(T\) has exponent \(2\) and is therefore abelian. It contains and centralizes \(A\). Thus the normal centralizer \(C_T(A)\) has nontrivial image in \(T/A\), so \(T=A C_T(A)\). The group \(A\) is abelian and itself lies in \(C_T(A)\), whence \(A\le Z(T)\). Hence \(T\) is a perfect central cover of \(A_5\). By \(M(A_5)=C_2\) it is \(A_5\) or \(2.A_5\). The second case has a quaternion Sylow \(2\)-subgroup and is therefore impossible. Thus \(T=A_5\) and \(A=1\). \end{proof}

\Needspace{11\baselineskip} \begin{lemma}\label{lem:A5-orthogonal-small-five-rank} Suppose that \(T=A_5\) acts on a finite-dimensional space \(X\) over \(\mathbb F_2\) preserving a quadratic form with nondegenerate polar form. Suppose that for an element \(a\in T\) of order \(5\), \[
 \dim[X,a]=4,\qquad q(C_X(T))=0.
\] Then \(\dim X=4\), and \(X\) is isomorphic to the coordinate-sum-zero submodule in the natural permutation module \(\mathbb F_2^5\). \end{lemma} \begin{proof} We list all simple \(\mathbb F_2A_5\)-modules used: the trivial one-dimensional module, the natural module \(N=\mathbb F_4^2\) for \(A_5\cong\operatorname{SL}_2(4)\) viewed over \(\mathbb F_2\), and the module \[
 S=\{(x_1,\ldots,x_5)\in\mathbb F_2^5:
                 \textstyle\sum_i x_i=0\}.
\] Both nontrivial modules have dimension \(4\) and have no nonzero fixed vectors for an element of order \(5\). Completeness of the list agrees with the full ATLAS catalogue \cite{A5TwoData}; we also justify completeness. Over an algebraic closure, the natural two-dimensional type and its field conjugate are distinguished by traces of elements of order \(5\). Together with \(1,S\) these give four simple types. We use the standard Brauer theorem: for a finite group over an algebraically closed field of characteristic \(p\), the number of simple types equals the number of conjugacy classes of elements of order coprime to \(p\) \cite{CurtisReiner}. For \(A_5\) and \(p=2\), there are exactly four such classes: \(1,3,5A,5B\). Descent to \(\mathbb F_2\) combines the two conjugate two-dimensional types into \(N\), while \(1,S\) are already defined over \(\mathbb F_2\).

We verify the required properties of these modules. On \(N\) the group is transitive on all \(15\) nonzero vectors, so \(N\) is simple over \(\mathbb F_2\). It has no nonzero invariant quadratic form. Such a form would be constant on nonzero vectors. If its value is \(0\), it is zero; if its value is \(1\), its polarization would be \(1\) on any two distinct nonzero vectors. For independent \(u,v,w\) this contradicts \(b(u,v+w)=b(u,v)+b(u,w)\).

The module \(S\) is absolutely simple. Indeed, for a Sylow subgroup \(P_2\cong C_2^2\) permuting four points and fixing the fifth, \(C_S(P_2)\) is the line spanned by \((1,1,1,1,0)\). Every nonzero submodule in characteristic \(2\) has nonzero fixed vectors for the \(2\)-group \(P_2\), so contains this line. Its five conjugates span all of \(S\). The restriction of \(S\) to \(P_2\) is isomorphic to the regular module: the first four coordinates are arbitrary and the fifth is their sum. Consequently, \(S\) is projective. In detail, a section of any epimorphism onto \(S\) can first be chosen \(P_2\)-linear and then averaged over the \(15\) cosets of \(A_5/P_2\); the number \(15\) is invertible in characteristic \(2\). The coordinate dot product is nondegenerate on \(S\), so \(S\) is self-dual. Thus projectivity and duality split extensions with \(S\) both as quotient and as submodule.

Taking \(C_5\)-fixed spaces is exact on short exact sequences because \(5\) is invertible in \(\mathbb F_2\). Hence the rank of \(a-1\) is additive over composition factors. From \(\dim[X,a]=4\) it follows that \(X\) has exactly one nontrivial composition factor, equal to \(N\) or \(S\).

If it is \(S\), the stated splittings give \(X=S\oplus X_0\). All factors of \(X_0\) are trivial; the image of the perfect group \(A_5\) on it is unitriangular, hence a \(2\)-group, and is therefore trivial. There are no cross terms of the polar form between \(S\) and \(X_0\) because \(\operatorname{Hom}_{A_5}(S,1)=0\). Thus the form on \(X_0\) is nondegenerate. For \(X_0\ne0\) it has a vector with \(q=1\) fixed by \(A_5\), contrary to the hypothesis. Consequently, \(X=S\).

Suppose the unique nontrivial factor is \(N\), and choose such an example of least dimension. There is no submodule \(N\) in \(X\). Otherwise the restriction of \(q\) to it would be zero, as proved above, so \(N\) would be totally singular. Nondegeneracy gives \(X/N^\perp\cong N^*\), and a second nontrivial composition factor would occur in \(X\), a contradiction. Thus a minimal submodule of \(X\) is trivial, and there exists a nonzero \(v\in C_X(A_5)\). By hypothesis, \(q(v)=0\).

The quotient \(v^\perp/\langle v\rangle\) is nondegenerate, has dimension smaller by \(2\), and has the same unique factor \(N\). It also has no fixed nonsingular vector. Indeed, if a class \(w+\langle v\rangle\) is fixed, then \(w^g-w=f(g)v\), where \(f:A_5\to\mathbb F_2\) is a homomorphism because \(v\) is fixed. Perfection gives \(f=0\), so \(w\) itself is fixed and \(q(w)=0\). This is a smaller example satisfying the same conditions, contrary to minimality. The case \(N\) is impossible, proving the lemma. \end{proof}

\Needspace{11\baselineskip} \begin{lemma}\label{lem:A5-seven-no-extraspecial-any-size} In characteristic \(7\), the critical section \(D\) of Lemma~\ref{lem:A5-five-seven-critical-section} cannot have center of order \(2\), regardless of its order or the multiplicity of its simple type on \(W\). Thus, in characteristic \(7\), only the following remains: \[
 Z(D)\cong C_4,\qquad |D|\ge64.
\] \end{lemma} \begin{proof} Suppose that \(Z(D)=C_2\). By Lemma~\ref{lem:A5-seven-extraspecial-action-simple}, the group \(A_5\) acts on the quadratic space \(\mathcal V=D/Z(D)\), and by Lemma~\ref{lem:A5-seven-extraspecial-five-rank-four} an element of order \(5\) has rank \(4\).

Every fixed vector \(v\) of this \(A_5\)-action is singular. Otherwise a lift \(x\in D\) of it would satisfy \(x^2=z\), and the cyclic group \(\langle x\rangle\) of order \(4\) would be normal in \(E_0\): a fixed class means \(x^e=x\) or \(xz=x^{-1}\). The perfect group \(E_0\) centralizes a normal cyclic subgroup because its automorphism group is abelian. Then \(x\) also centralizes \(D\), so \(x\in Z(D)\), a contradiction.

Applying Lemma~\ref{lem:A5-orthogonal-small-five-rank} gives \(\dim\mathcal V=4\), that is, \(|D|=32\). This case is completely excluded by Lemma~\ref{lem:A5-seven-no-critical32}. Hence a center of order \(2\) is impossible for every \(n\). The only remaining center size is \(4\), and the previous lower bound \(\dim D/Z(D)\ge4\) gives \(|D|\ge64\). \end{proof}

\clearpage \subsection{Global exclusion of characteristic \(5\) and refinement for \(2,3\)}

Here the original group \(G\) with its entire spectrum \(\Om\) is again essential, rather than just the local action image. We retain the notation \(L=G^\infty\), \(P=R(G)\cap L\), \(K=[P,L]\), and \(\operatorname{Soc}(G/R(G))\cong A_5\).

\begin{lemma}\label{lem:A5-odd-Sylow-radical-noncyclic} For every odd prime \(r\mid|P|\), a Sylow \(r\)-subgroup of \(P\) is noncyclic. The odd-order Sylow subgroups of \(P\) and \(K\) have the same orders, so the same assertion holds for \(K\). \end{lemma} \begin{proof} Suppose a Sylow \(r\)-subgroup of \(P\) is cyclic. Then every chief \(r\)-factor of \(G\) inside \(P\) has order \(r\): it is elementary abelian and is a section of a group with a cyclic Sylow \(r\)-subgroup. The action of the perfect group \(L\) on such a one-dimensional factor is trivial because \(\operatorname{GL}_1(r)\) is abelian. But \(M(A_5)=C_2\), that is, \(r\nmid|M(A_5)|\). Lemma~\ref{lem:central-chief-prime-elimination} gives \(r\nmid|P|\), a contradiction. Finally, \(P/K\) has order \(1\) or \(2\), being the center of a perfect central cover \(L/K\) of \(A_5\). Thus the odd parts of the orders of \(P\) and \(K\) coincide. \end{proof}

\Needspace{11\baselineskip} \begin{lemma}\label{lem:A5-no-five-in-perfect-radical} In the branch under consideration, \[
 5\nmid|P|,\qquad 5\nmid|K|.
\] Consequently, characteristic \(5\) of the internal chief factor in Lemma~\ref{lem:A5-K-internal-chief-exists} is completely excluded. \end{lemma} \begin{proof} By Lemma~\ref{lem:remaining-commutator-localization}, we have \(7\mid|K|\). Suppose also that \(5\mid|K|\). In the solvable group \(K\) choose a Hall \(\{5,7\}\)-subgroup \(H\). Recall the proof of existence of a Hall subgroup for any set of primes in a finite solvable group. Induct on the group order: a minimal normal subgroup \(A\) is elementary abelian of characteristic \(r\), and the quotient has the required Hall subgroup. If \(r\) belongs to the chosen set, its full preimage already suffices. Otherwise, in the preimage, \(A\) is a normal Hall subgroup, and a complement given by Schur--Zassenhaus yields the required Hall subgroup of the original group.

The Sylow \(5\)- and \(7\)-subgroups of \(H\) are noncyclic by the preceding lemma. Since \(25,49\notin\Om\), they have exponents \(5,7\). The Fitting subgroup \(F(H)\) is nontrivial because \(H\) is solvable. It cannot have both prime divisors, since this would give \(35\). Thus \(F(H)\) is an \(r\)-group for \(r\in\{5,7\}\). Let \(s\) be the other prime in this pair. The subgroup \(A=\Omega_1(Z(F(H)))\) is nontrivial, elementary abelian, and normal in \(H\).

A noncyclic Sylow \(s\)-subgroup of exponent \(s\) contains \(C_s^2\): take a central element of order \(s\) and an element outside the subgroup it generates. On the nonzero space \(A\), this elementary abelian group of rank \(2\) is diagonalizable after field extension because \(s\ne r\). Each of its linear characters has a nontrivial kernel. Hence some element of order \(s\) fixes a nonzero vector of \(A\), giving order \(rs=35\). This contradicts \(\omega(H)\subseteq\Om\). Consequently, \(5\nmid|K|\), and \(|P/K|\le2\) gives \(5\nmid|P|\). Every chief factor inside \(K\) therefore has characteristic different from \(5\). \end{proof}

\Needspace{11\baselineskip} \begin{lemma}\label{lem:A5-internal-seven-action-two-three} In the notation of Lemma~\ref{lem:A5-K-internal-chief-exists}, for \(p=3\) we have \(U_7=1\), and hence \(U=U_2\). For \(p=2\) the groups \(U_3,U_7\) cannot both be nontrivial. If \(U_7\ne1\), then \[
 U=U_7\text{ elementary abelian},\qquad
 E\cong U_7\rtimes2.A_5.
\] After extending the field, all simple constituents of this \(2.A_5\)-module \(U_7\) are faithful two-dimensional types. In particular, elements of order \(3\) act on \(U_7\) without nonzero fixed vectors. \end{lemma} \begin{proof} By the preceding lemma, \(U_5=1\). For \(p=3\) and \(U_7\ne1\), the group \(U_2\) has exponent \(2\) by Lemma~\ref{lem:A5-internal-four-characteristics}. If \(U_2\ne1\), it has rank at least \(4\), and \(U_7\) contains \(C_7^2\). For every involution \(t\in U_2\), the space \(C_V(t)\) is \(U_7\)-invariant. If it is nonzero, restriction to \(C_7^2\) gives a vector fixed by some element of order \(7\). Together with \(t\) and this vector of order \(3\), we obtain \(42\). Thus \(C_V(t)=0\), and every involution of \(U_2\) acts as \(-I\). This is impossible for \(C_2^2\le U_2\). Hence \(U_2=1\) and \(U=U_7\).

For \(p=2\) and \(U_7\ne1\), the group \(U_3\) has exponent \(3\). If it is nontrivial, its rank is at least \(4\), so it contains \(C_3^2\). For each element \(a\) of order \(3\) in this subgroup, a nonzero \(C_V(a)\) would again give \(42\) using \(C_7^2\le U_7\). But a nonzero \(C_3^2\)-module of characteristic \(2\) always has a nonzero vector fixed by some nonidentity element of \(C_3^2\). This is a contradiction; hence \(U_3=1\), and again \(U=U_7\).

Thus, in both cases with \(U_7\ne1\), we have \(E/U_7\cong A_5\) or \(2.A_5\). The normal \(7\)-subgroup is a Hall subgroup of \(E\), so \(E=U_7\rtimes D_0\), \(D_0\cong A_5\), or \(2.A_5\). On every chief factor inside \(U_7\), the group \(U_7\) itself acts trivially: a normal \(7\)-group has nonzero fixed vectors on a simple module of characteristic \(7\), and their space is invariant under all of \(E\). Thus these factors are \(D_0\)-modules.

The required data for \(A_5\) in characteristic \(7\) coincide with the ordinary data because \(7\nmid60\). The full list of degrees and dimensions of fixed spaces of an element of order \(5\) is \[
\begin{array}{c|ccccc}
 \dim&1&3&3&4&5\\ \hline
 \dim C(5)&1&1&1&0&1.
\end{array}
\] For the four-dimensional type, \(\dim C(2)=\dim C(3)=2\). These data follow from the ordinary table of \(A_5\); see \cite{Atlas,A5SevenData}. The four-dimensional type is the coordinate-sum-zero submodule of the permutation module on five points. The exclusion of \(35\) leaves only this type on the chief \(7\)-factors of \(A_5\).

If \(D_0=A_5\), then in the case \(p=3\) every involution \(t\in D_0\) has a noncyclic centralizer in \(U_7\): the map \(C_{U_7}(t)\to C_{U_7/\Phi(U_7)}(t)\) is surjective by coprimeness, and the target has rank at least \(2\). Since \(\exp U_7=7\), the centralizer contains \(C_7^2\). As in the first paragraph, the exclusion of \(42\) gives \(C_V(t)=0\). The three involutions in a four-group in \(A_5\) cannot all act as \(-I\). This excludes the case for \(p=3\).

For \(p=2\) and \(D_0=A_5\), take \(a\in D_0\) of order \(3\). Now \(C_{U_7}(a)\) contains \(C_7^2\) by the same calculation using \(\dim C(3)=2\). Hence \(C_V(a)=0\), since otherwise \(42\) occurs. The action of \(a\) on \(U_7\) is nontrivial: otherwise simplicity of \(A_5\) would give \([U_7,D_0]=1\), contrary to \([U_7,E]=U_7\). Lemma~\ref{lem:HH-prime} for \(U_7\rtimes\langle a\rangle\) in characteristic \(2\) gives minimal polynomial \(x^3-1\) because the normal subgroup has odd order. Thus \(C_V(a)\ne0\), a contradiction. Therefore \(D_0=A_5\) is excluded in both cases.

Let \(D_0=2.A_5\), and let \(z\) be its central involution. On every chief \(7\)-factor, the \((+1)\)-eigenspace of \(z\) would be an \(A_5\)-module. The exclusion of \(35\) leaves only the four-dimensional type there, but an element \(d\in D_0\) of order \(4\) acts on it as an involution of \(A_5\) and has nonzero fixed vectors. This would give \(28\). Thus \(z\) acts as \(-I\) on all chief factors inside \(U_7\). Successively lifting fixed points gives \(C_{U_7}(z)=1\). Hence \(U_7\) is abelian and \(z\) inverts it, by the proved property of a fixed-point-free automorphism of order \(2\).

For \(p=3\), take \(a\in D_0\) of order \(3\). It acts nontrivially on \(U_7\), since otherwise \([U_7,E]\ne U_7\). The simple types in the restriction of \(V\) to the abelian group \(U_7\) are linear characters and separate its elements by faithfulness of the representation. Thus \(a\) permutes some triple of distinct weights. The automorphism \(z\) inverts weights; a weight and its inverse cannot lie in the same orbit of length \(3\). Consequently, \(za\) permutes six distinct weight spaces. The norm of order \(6\) on them is nonzero, giving \(18\in\omega(V\rtimes E)\) in characteristic \(3\), a contradiction. Characteristic \(3\) with \(U_7\ne1\) is thereby completely excluded.

For \(p=2\) it remains to refine the simple constituents of \(U_7\). The list in Lemma~\ref{lem:double-A5-seven-simple-data} and the exclusion of \(35\) leave the faithful types of dimensions \(2,2,4\). For the four-dimensional type, an element of order \(3\) has a two-dimensional fixed space. If such a constituent occurs, \(C_{U_7}(a)\) contains \(C_7^2\), so the exclusion of \(42\) gives \(C_V(a)=0\). As above, Lemma~\ref{lem:HH-prime} gives the opposite conclusion. Thus four-dimensional constituents are impossible, leaving the two two-dimensional types. On each of them an element of order \(3\) acts without nonzero fixed vectors: its determinant is \(1\), and a faithful two-dimensional type cannot make this element trivial. \end{proof}

\Needspace{11\baselineskip} \begin{lemma}\label{lem:A5-final-local-characteristic-table} For the internal chief factor of Lemma~\ref{lem:A5-K-internal-chief-exists}, the following necessary local cases remain: \[
\begin{array}{c|l}
 p&\text{image }U\text{ of the group }K\\ \hline
 2&U=U_3\ne1,\ \text{or }U=U_7\ne1
      \text{ of the special form in Lemma~\ref{lem:A5-internal-seven-action-two-three}};\\
 3&U=U_2\ne1;\\
 5&\text{impossible};\\
 7&U=U_2\text{ nonabelian},\ |U_2|\ge64,
       \ Z(D)\cong C_4.
\end{array}
\] Every Sylow subgroup that occurs has Frattini rank at least \(4\). These local conditions alone do not exclude characteristics \(2,3,7\). \end{lemma} \begin{proof} Characteristic \(5\) is excluded by Lemma~\ref{lem:A5-no-five-in-perfect-radical}. For \(p=2,3\) apply the preceding lemma and \(U_5=1\). For \(p=7\) a nonabelian \(2\)-group remains by Lemma~\ref{lem:A5-internal-five-seven-two-kernel}. By Lemma~\ref{lem:A5-seven-no-extraspecial-any-size} its critical section \(D\) has center \(C_4\) and order at least \(64\). Since \(D\) is a section of \(U\), this gives \(|U|\ge64\). The rank bounds follow from Lemma~\ref{lem:A5-internal-chief-Frattini-ranks}. \end{proof}

\Needspace{11\baselineskip} \begin{lemma}\label{lem:A5-real-local-three-five-models} The local spectral inclusion alone does not exclude all remaining action groups. For \(p=3,5\) there exist faithful simple four-dimensional modules \(V_p\) of \(E_p\), where \[
 |E_p|=1920,\quad E_p=E_p',\quad
 R(E_p)=D_p\cong2_-^{1+4},\quad
 [D_p,E_p]=D_p,\quad E_p/D_p\cong A_5,
\] and \[
 \omega(V_3\rtimes E_3)=\{1,2,3,4,5,6,8,10,12\},
\] \[
 \omega(V_5\rtimes E_5)=\{1,2,3,4,5,6,8,10,12,20\}.
\] Both spectra are proper subsets of \(\Om\). The characteristic \(5\) model is excluded for the original group by the global Lemma~\ref{lem:A5-no-five-in-perfect-radical}, rather than by a contradiction to the local spectral inclusion itself. \end{lemma} \begin{proof} We give an explicit construction and a complete finite certificate for the spectrum calculations. Over \(\mathbb F_p\) put \[
 i=\begin{pmatrix}0&1\\-1&0\end{pmatrix},\quad
 x=\begin{pmatrix}0&1\\1&0\end{pmatrix},\quad
 y=\begin{pmatrix}1&0\\0&-1\end{pmatrix},
\] \[
 j=\begin{cases}
 \begin{pmatrix}1&1\\1&2\end{pmatrix},&p=3,\\
 \begin{pmatrix}0&2\\2&0\end{pmatrix},&p=5.
 \end{cases}
\] The matrices \(i\otimes I,j\otimes I,I\otimes x,I\otimes y\) generate \(D_p=Q_8\circ D_8\) and act faithfully and absolutely irreducibly on \(V_p=\mathbb F_p^4\). They preserve the form \(B=i\otimes I\). Let \(\mathcal R\) contain one representative from each pair \(\{r,-r\}\) with \(r^2=-I\) in \(D_p\). There are \(10\) such pairs. Choose \(R\in\mathcal R\) and put \[
 E_p=\left\langle D_p,\quad
 2^{-1}(I+R)(I+S)\ (S\in\mathcal R\setminus\{R\})
 \right\rangle.
\] As in the proof of Lemma~\ref{lem:extraspecial32-normalizer-fixed-four}, these matrices are symplectic and normalize \(D_p\). On \(D_p/Z(D_p)\) they induce all products of one fixed orthogonal transvection with the others. On the five singular vectors, the transvections give all transpositions; their pairwise products generate \(A_5\). The kernel of the normalizer action is \(D_p\), by the same proof using inner automorphisms and scalars \(\pm I\). Consequently, \(E_p/D_p\cong A_5\) and \(|E_p|=32\cdot60=1920\).

The action of \(A_5\) on \(D_p/Z(D_p)\) is the simple four-dimensional orthogonal type over \(\mathbb F_2\). Thus \([D_p,E_p]\) maps onto all of \(D_p/Z(D_p)\). Since \(\Phi(D_p)=Z(D_p)\), we obtain \([D_p,E_p]=D_p\). Perfection of \(E_p/D_p\) now gives \(E_p'=E_p\). The normal \(2\)-group \(D_p\) is the entire radical because its quotient is nonabelian simple.

The certificate in Appendix~\ref{app:local-three-five-models} constructs all \(1920\) matrices, checks \([D_p,E_p]=D_p\), and for each matrix \(a\) finds its order \(m\) and norm \(I+a+\cdots+a^{m-1}\). A nonzero norm adds order \(pm\) to the spectrum of \(V_p\rtimes E_p\); a zero norm adds no new order. Exactly the two displayed sets result. Thus the finite calculation exhausts all elements of the quotient and all possible orders of their lifts in the split product. Neither model has the full spectrum \(\Om\), and they are not counterexamples to the original classification. \end{proof}

\clearpage \subsection{Excluding center \(C_4\) and all of characteristic \(7\)}

Assume \(p=7\) and retain the notation of Lemma~\ref{lem:A5-five-seven-critical-section}. By Lemma~\ref{lem:A5-seven-no-extraspecial-any-size}, only \(Z(D)=\langle c\rangle\cong C_4\) remains. Put \(z=c^2\) and \(X=D/Z(D)\). The commutator defines a nondegenerate alternating form \(b\) on \(X\) with values in \(\langle z\rangle\cong\mathbb F_2\). The group \(E_0\) is perfect and centralizes \(Z(D)\).

\begin{lemma}\label{lem:A5-C4-no-fixed-and-exponent-four} There are no nonzero \(E_0\)-fixed vectors on \(X\). The group \(D\) has exponent \(4\) and is isomorphic to a central product \(C_4\circ D_1\), where \(D_1\) is extraspecial of order \(2^{1+2n}\), \(2n=\dim X\). The faithful simple type of \(D\) with a prescribed faithful central character has dimension \(2^n\) and is unique. \end{lemma} \begin{proof} Let \(xZ(D)\) be fixed by \(E_0\). Then \(x^e=x f(e)\), where \(f(e)\in Z(D)\). Since \(E_0\) centralizes \(Z(D)\), the map \(f:E_0\to Z(D)\) is a homomorphism. Perfection gives \(f=1\), so \(x\in C_D(E_0)\le Z(D)\) because \(D\le E_0\). Consequently, \(C_X(E_0)=0\).

For \(x\in D\), the square lies in \(Z(D)\). The map \[
 \ell:X\longrightarrow Z(D)/D'\cong C_2,
 \qquad \ell(xZ(D))=x^2D'
\] is well defined: replacing \(x\) by \(xc\) changes the square by \(c^2=z\in D'\). It is linear because commutators lie in \(D'\), and it is \(E_0\)-invariant. If \(\ell\ne0\), nondegeneracy of \(b\) represents it as \(\ell(v)=b(u,v)\) for a nonzero \(u\in X\). Invariance of the form and of \(\ell\) then makes \(u\) fixed, a contradiction. Thus \(\ell=0\), all squares lie in \(D'=\langle z\rangle\), and \(\exp D=4\).

Choose a symplectic basis of \(X\) and its preimages in \(D\). Multiplication of a preimage by \(c\) switches its square between \(1\) and \(z\) without changing commutators. Thus the squares can be chosen to specify any central product of groups \(D_8,Q_8\) of the corresponding rank. The group \(D_1\) generated by the preimages is extraspecial, \(D_1\cap\langle c\rangle
=\langle z\rangle\), and \(D=\langle c\rangle D_1\). This proves the central-product assertion. The faithful simple type of \(D_1\) is unique and has dimension \(2^n\). Its extension to \(D\) is determined by choosing the scalar representing \(c\) to be one of the two square roots of \(-1\). Hence each prescribed faithful central character gives a unique type. \end{proof}

\Needspace{11\baselineskip} \begin{lemma}\label{lem:A5-C4-five-rank-and-multiplicity} For an element \(a\in E_0\) of order \(5\) outside \(Q_0\), we have \[
 \dim[X,a]=4.
\] The module \(W\) decomposes into a direct sum of identical genuine \(E_0\)-modules \(H\) of dimension \(2^n\), faithful and irreducible on \(D\). The image of \(E_0\) on \(H\) has determinant \(1\). \end{lemma} \begin{proof} Let \(N=2^n\) and \(f=\dim C_X(a)\). Since the action of \(a\) on \(X\) is nontrivial, \(2n-f=4r\) and \(r\ge1\). On the simple type \(H\), choose an operator \(A\) implementing conjugation by \(a\), normalized by \(A^5=I\). Let \(m_0,\ldots,m_4\) be the multiplicities of the fifth roots of unity in \(A\).

Matrices representing \(D/Z(D)\) form a basis of \(\operatorname{End}(H)\). Indeed, a representative of a nonzero class has trace zero: by nondegeneracy of the commutator it is conjugate to its product with \(z\), which acts with the opposite sign. Trace orthogonality and the number \(|D/Z(D)|=N^2\) give the basis assertion. If the class \(xZ(D)\) is fixed by \(a\), then \(x^a=xt\), where \(t\in Z(D)\). Since \(a\) centralizes the center, \(a^5=1\) implies \(t^5=1\) and hence \(t=1\). Thus all \(2^f\) fixed basis matrices are fixed without a scalar multiplier.

In the complex extension of the type to \(D\rtimes C_5\), we obtain \[
 |\operatorname{tr}(A^j)|^2=2^f\quad(1\le j\le4),
 \qquad
 \sum m_t^2=\frac{N^2+4\cdot2^f}{5},
 \quad\sum m_t=N.
\] Multiplicities are preserved on passage to characteristic \(7\), coprime to the order of \(D\rtimes C_5\). At least one fifth root of unity is absent from the spectrum of \(A\). Otherwise the action of \(a\) on \(W=H\otimes M\), where \(M\) is the multiplicity space, would have a nonzero fixed vector for any choice of operator on \(M\), giving \(35\). Cauchy's inequality for at most four nonzero \(m_t\) gives \[
 \frac{N^2+4\cdot2^f}{5}\ge\frac{N^2}{4},
 \qquad 2^{2n-f}\le16.
\] Consequently, \(2n-f=4\). Equality holds in the inequality: exactly four roots occur, each with multiplicity \(N/4\).

Rescale \(A\) by a fifth root of unity so that precisely the root \(1\) is absent. This does not change \(A^5=I\). If \(a\) acts on \(W\) as \(A\otimes S_a\), then \(S_a^5=I\), and exclusion of \(35\) forces \(S_a=I\). As in the proof of Lemma~\ref{lem:A5-seven-no-critical32}, the projective image of \(E_0\) on \(M\) has kernel containing \(a\notin Q_0\). Thus it coincides with the image of the \(2\)-group \(Q_0\) and is also perfect, so is trivial. All operators on \(M\) are scalar; absorbing the scalars into the operators on \(H\) gives identical genuine \(E_0\)-modules \(H\). Their determinant is a linear character of a perfect group and is therefore \(1\). \end{proof}

\Needspace{11\baselineskip} \begin{lemma}\label{lem:A5-C4-reduction-order64} We have \(\dim X=4\) and \(|D|=64\). The image of \(E_0\) on \(X\) is \(A_5\), and \(Q_0\) acts trivially on \(X\). On one simple type \(H\) of the preceding lemma, the image \(F\) of \(E_0\) lies in \(\operatorname{SL}_4(k)\), contains and normalizes \(D\), and satisfies \(F/D\cong A_5\) and \(|F|=3840\). \end{lemma} \begin{proof} Denote the image of the action on \(X\) by \(T\), and the image of \(Q_0\) by \(A\). Then \(T\) is perfect, \(A\) is a normal \(2\)-subgroup, and \(T/A\cong A_5\) because the kernel of the action lies in \(Q_0\). On every simple \(T\)-module of characteristic \(2\), the group \(A\) acts trivially: its fixed vectors form a nonzero \(T\)-submodule. Thus the simple composition factors of \(X\) are \(A_5\)-modules. The complete list from the proof of Lemma~\ref{lem:A5-orthogonal-small-five-rank} consists of the trivial type and two types of dimension \(4\). An element of order \(5\) has rank \(4\) on each nontrivial type. By the preceding lemma, \(X\) has exactly one nontrivial composition factor.

Choose a minimal nonzero submodule \(X_1\le X\). It cannot be trivial, since \(C_X(T)=0\) by Lemma~\ref{lem:A5-C4-no-fixed-and-exponent-four}. Therefore \(\dim X_1=4\). The restriction of the nondegenerate form \(b\) to the simple module \(X_1\) is either nondegenerate or zero. The latter is impossible: if \(X_1\le X_1^\perp\), the quotient \(X/X_1^\perp\cong X_1^*\) would give a second nontrivial composition factor. Hence \(X=X_1\perp X_1^\perp\). All factors of \(X_1^\perp\) are trivial, so the perfect image of \(T\) on it is a \(2\)-group and is trivial. The absence of nonzero fixed vectors gives \(X_1^\perp=0\). Thus \(X=X_1\) is simple and four-dimensional. The normal \(2\)-group \(A\) is trivial on the simple module \(X\), while the action of \(T\) on \(X\) is faithful. Hence \(A=1\) and \(T=A_5\).

Thus \(|D|=4\cdot2^4=64\) and \(\dim H=4\). The central element \(c\) acts as a scalar \(\iota\), with \(\iota^2=-1\). By Lemma~\ref{lem:A5-C4-no-fixed-and-exponent-four}, one can choose in \(D\) an extraspecial subgroup \(D_1\cong Q_8\circ D_8\) of order \(32\), so that \(D=\langle\iota I\rangle D_1\) in the model of Lemma~\ref{lem:extraspecial32-five-eigenvalues}.

The kernel of the action of the normalizer of \(D\) in \(\operatorname{SL}_4(k)\) on \(X\) is \(D\). To verify this, suppose an automorphism is trivial on \(X\) and centralizes \(Z(D)\). It has the form \(x\mapsto xf(xZ(D))\). Since squares are preserved, \(f\) takes values in \(\langle z\rangle\). Nondegeneracy of \(b\) makes this automorphism inner. After adjustment by an element of \(D\), a scalar remains by irreducibility. Its fourth power is \(1\) by the determinant condition, so it already lies in \(\langle\iota I\rangle\le D\). The reverse inclusion is obvious. Thus for the image \(F\) of \(E_0\) we obtain \(F/D=T=A_5\) and \(|F|=3840\). \end{proof}

\Needspace{11\baselineskip} \begin{lemma}\label{lem:A5-C4-normalizer-fixed-four} Let \(D=\langle\iota I\rangle
(Q_8\circ D_8)\le\operatorname{SL}_4(k)\), \(\operatorname{char}k=7\), and \(\iota^2=-1\). If \(F\le\operatorname{SL}_4(k)\) contains and normalizes \(D\) and has image \(A_5\) on \(X=D/Z(D)\), then \(F\) contains an element of order \(4\) with a two-dimensional fixed space. \end{lemma} \begin{proof} We describe the action of \(\operatorname{Sp}(X,b)\) on the six quadratic forms of minus type whose polarization is \(b\). For one form \(q_-\), all refinements have the form \(q_t(v)=q_-(v)+b(t,v)\), \(t\in X\). The six forms with \(q_-(t)=0\) have minus type. This is checked by substitution in the four-dimensional form of Lemma~\ref{lem:extraspecial32-five-eigenvalues} or by counting its zeros. The action on these six forms is faithful: preserving every form preserves their differences, and the five nonzero singular vectors of \(q_-\) span \(X\). The number of symplectic bases is \(15\cdot8\cdot3\cdot2=720\). Hence this action identifies \(\operatorname{Sp}_4(2)\) with \(S_6\). A finite verification is given in Appendix~\ref{app:Sp4-six-forms}.

Every symplectic transvection \(\tau_v(w)=w+b(w,v)v\), \(v\ne0\), is a transposition on these six forms. Indeed, it fixes a form \(q\) exactly when \(q(v)=1\). Exactly four of the six minus forms satisfy this: each form has ten nonsingular nonzero vectors, and the symplectic group is transitive on the fifteen nonzero vectors, so the number is \(6\cdot10/15=4\).

The image \(A_5\) of \(F\) lies in \(A_6\) because \(A_5\) is perfect. Thus its nontrivial involution on the six forms is the product of two distinct disjoint transpositions. This is \(\tau_r\tau_s\), where \(r,s\in X\) are distinct and \(b(r,s)=0\): two distinct transvections commute exactly under this condition.

Choose preimages \(R,S\in D\) of the classes \(r,s\) with \(R^2=S^2=-I\). This is always possible by multiplication by the central \(\iota I\). Since \(b(r,s)=0\), the matrices \(R,S\) commute. Define \[
 H_R=2(I+R),\qquad H_S=2(I+S).
\] They normalize \(D\) and induce \(\tau_r,\tau_s\) by the same direct calculation as in Lemma~\ref{lem:extraspecial32-normalizer-fixed-four}. Their determinants are \(1\). Indeed, \(R\) has eigenvalues \(\iota,-\iota\), each twice: its trace is zero because it is noncentral in \(D\). Thus \[
 \det(2(I+R))=2^4(1+\iota)^2(1-\iota)^2
 =16\cdot4=1\quad\text{in characteristic }7,
\], and similarly for \(S\). Consequently, \(H_RH_S\) belongs to the normalizer of \(D\) in \(\operatorname{SL}_4(k)\) and has image \(\tau_r\tau_s\in A_5\). The kernel of the normalizer action is \(D\) by the preceding lemma. Since \(D\le F\), we obtain \(T=H_RH_S\in F\).

We have \(T^2=RS\ne I\) and \(T^4=I\). The four joint eigenvalues \((\varepsilon\iota,\delta\iota)\) of the matrices \(R,S\) each occur once. This follows from the zero traces of \(R,S,RS\): all three matrices are noncentral in \(D\), and trace orthogonality on the abelian group \(\langle R,S,Z(D)\rangle\) gives four equal multiplicities. For \(\varepsilon=-\delta\) the eigenvalue of \(T\) is \(4(1+\iota)(1-\iota)=1\); for the other two pairs the eigenvalues are \(\iota,-\iota\). Thus \(|T|=4\) and \(\dim C_{k^4}(T)=2\). \end{proof}

\Needspace{11\baselineskip} \begin{lemma}\label{lem:A5-internal-seven-complete-exclusion} Under the assumptions of Lemma~\ref{lem:A5-K-internal-chief-exists}, characteristic \(7\) of the chosen highest internal chief factor is impossible. Together with Lemma~\ref{lem:A5-no-five-in-perfect-radical}, this leaves only characteristics \(2\) and \(3\) for that factor. \end{lemma} \begin{proof} For \(p=7\) the critical section exists by Lemma~\ref{lem:A5-five-seven-critical-section}. A center of order \(2\) was already completely excluded by Lemma~\ref{lem:A5-seven-no-extraspecial-any-size}. When the center is \(C_4\), the lemmas of this subsection reduce the action to \(F\le\operatorname{SL}_4(k)\) and give in \(F\) an element \(T\) of order \(4\) with a nonzero fixed vector. Lift \(T\) to a \(2\)-element \(t\in E_0\); its order is \(4\) or \(8\) because \(16\notin\Om\). It fixes a nonzero vector of \(H\) and hence of \(W\). Lifting through the normal \(2\)-subgroup \(J\) in \(E\) also gives order \(4\) or \(8\). The kernel dimension is unchanged on descending to the original field of definition of \(V\). Thus \(28\) or \(56\) occurs in \(V\rtimes E\), contradicting the inclusion of spectra. Both possibilities for the center are excluded, so \(p=7\) is impossible. \end{proof}

\clearpage \subsection{Local models in characteristics \(2,3\)}

After Lemma~\ref{lem:A5-internal-seven-complete-exclusion}, the highest internal chief factor has characteristic \(2\) or \(3\). Below we determine which required orders must disappear in its action group and construct a genuine example of the remaining local case in characteristic \(2\). These models distinguish local module restrictions from requirements imposed by the full spectrum.

\begin{lemma}\label{lem:A5-two-three-orders-lost} In the notation of Lemma~\ref{lem:A5-K-internal-chief-exists}, the following restrictions hold: \begin{enumerate} \item For \(p=3,U=U_2\), and also for \(p=2,U=U_3\), all elements of order \(7\) of the original group \(G\) lie in \(C_G(V)\). Consequently, \[
 7,14,21\notin\omega(V\rtimes Y).
\] \item For \(p=2,U=U_7\) and every element \(x\in G\) of order \(9\), we have \(x^3\in C_G(V)\). In this case \[
 9\notin\omega(V\rtimes Y).
\] \end{enumerate} Thus none of the remaining local images has the full spectrum \(\Om\); the missing orders must be supplied by other layers of the original group. \end{lemma} \begin{proof} By Lemma~\ref{lem:remaining-commutator-localization}, all elements of order \(7\) of \(G\) lie in \(K\). In the first case their images in \(Y\) lie in the \(2\)- or \(3\)-group \(U\) and are therefore trivial. Hence \(7\nmid|Y|\). For completeness, if \(Y\) had an element of order \(7\), a preimage in \(G\) would have order divisible by \(7\), and its \(7\)-part of order \(7\) would have nontrivial image in \(Y\). This contradicts the preceding conclusion. The kernel \(V\) has characteristic \(2\) or \(3\), so \(7\) also does not divide the order of \(V\rtimes Y\).

In the second case, \(x^3\in K\) for every \(|x|=9\), by the same lemma. Its image lies in the \(7\)-group \(U\) and has order dividing \(3\), so is trivial. If \(Y\) contained an element of order \(9\), any preimage in \(G\) would have order divisible by \(9\). The only such order in \(\Om\) is \(9\) itself. But the image of an element of order \(9\) has already been shown to have order at most \(3\), a contradiction. Finally, factoring out a normal \(2\)-group preserves order \(9\), so this order is also absent from \(V\rtimes Y\). \end{proof}

\Needspace{11\baselineskip} \begin{lemma}\label{lem:A5-real-binary-upper-model} There exist a perfect group \(E=N\rtimes D\), where \(N\cong C_3^4\) and \(D\cong2.A_5\), and a faithful absolutely simple \(E\)-module \(W\) over \(\mathbb F_2\) of dimension \(40\) such that \[
 \omega(E)=\{1,2,3,4,5,6,10\},
\] \[
 \omega(W\rtimes E)=
 \{1,2,3,4,5,6,8,10,12,20\}\subsetneq\Om.
\] For \(G_0=W\rtimes E\) we have \[
 L_0=G_0^\infty=G_0,\quad
 P_0=R(G_0)=W\rtimes(N\rtimes Z(D)),\quad
 K_0=[P_0,L_0]=W\rtimes N.
\] The group \(K_0\) is nonnilpotent, and \(W\) is a highest internal chief factor with \(p=2,U=N=U_3\). Consequently, this remaining local case cannot be declared impossible solely from the inclusion of spectra. \end{lemma} \begin{proof} On \(N=\mathbb F_3^4\) take the matrices of standard generators of \(D\) from the ATLAS \cite{A5BinaryLocalData}: \[
 A=\begin{pmatrix}
 0&1&0&0\\2&0&0&0\\0&0&0&1\\0&0&2&0
 \end{pmatrix},\qquad
 B=\begin{pmatrix}
 1&0&0&0\\0&0&1&0\\0&2&2&0\\1&1&2&1
 \end{pmatrix}.
\] They generate a perfect group of order \(120\), satisfy \(A^2=-I\) and \(B^3=(AB)^5=I\), and realize \(2.A_5\). The module \(N\) is simple and faithful. Appendix~\ref{app:binary-upper-model} gives a certificate for these assertions: it constructs all \(120\) matrices, checks perfection, the two orbits of length \(40\) on nonzero vectors and on nonzero linear forms, and full rank of the linear span of each orbit. The same norm enumeration verifies \(\omega(N\rtimes D)=\{1,2,3,4,5,6,10\}\). The group \(E\) is perfect because \(D\) is perfect and its central involution \(z=-I\) inverts \(N\), whence \([N,D]=N\).

We construct \(W\) explicitly using weights. Choose one orbit \(\mathcal O\) of length \(40\) of \(D\) on \(N^*\setminus\{0\}\). Let \(\omega\in\mathbb F_4\) have order \(3\). On a space over \(\mathbb F_4\) with basis \(e_\lambda\), \(\lambda\in\mathcal O\), define \[
 n e_\lambda=\omega^{\lambda(n)}e_\lambda,
 \qquad
 d e_\lambda=e_{\lambda\circ d^{-1}}
 \quad(n\in N,d\in D).
\] This is a representation of \(E=N\rtimes D\). All weights are distinct and \(D\) is transitive on them, so the representation is absolutely simple: every nonzero submodule contains a weight line and hence all the others.

Since \(z\) inverts \(N\), the orbit \(\mathcal O\) is closed under \(\lambda\mapsto-\lambda\). The semilinear transformation \[
 \tau\left(\sum_\lambda a_\lambda e_\lambda\right)
 =\sum_\lambda a_{-\lambda}^{\,2}e_\lambda
\] is an involution commuting with \(E\). For the action of \(N\) this follows from \(-2\lambda(n)=\lambda(n)\) in \(\mathbb F_3\), and for \(D\) from the commutation of permutations with inversion of weights. The fixed space \(W\) has dimension \(40\) over \(\mathbb F_2\): each pair of weights \(\lambda,-\lambda\) gives two free binary parameters. Extending scalars to \(\mathbb F_4\) gives the whole constructed space. On each pair this is seen from the two independent vectors \((1,1)\) and \((\omega,\omega^2)\). Hence \(W\) is absolutely simple.

The action of \(N\) on \(W\) is faithful: the orbit \(\mathcal O\) spans the entire dual space \(N^*\). Denote the kernel of the action of \(E\) by \(J\). Then \(J\cap N=1\). Both groups are normal, so \([J,N]\le J\cap N=1\). But \(C_E(N)=N\) because the action of \(D\) on \(N\) is faithful. Thus \(J\le N\) and \(J=1\).

For an elementary abelian \(2\)-kernel, an element above a matrix of order \(m\) has order \(m\) or \(2m\). Thus the spectrum of \(G_0\) lies in \[
 \omega(E)\cup2\omega(E)
 =\{1,2,3,4,5,6,8,10,12,20\}.
\] All three new orders \(8,12,20\) do occur. An element \(d\in D\) of order \(4\) has square \(z\), so its orbits on the nonzero weights have length \(4\); the nonzero norm gives \(8\). An element \(c\in D\) of order \(5\) acts on \(N^*\) without nonzero fixed vectors: otherwise \(E\) would contain an element of order \(15\), contrary to the computed spectrum. The orbit of the inverse weight does not coincide with the orbit of length \(5\), so \(zc\) has weight orbits of length \(10\) and gives order \(20\). Finally, \(B\) acts nontrivially on \(N^*\). Since \(\mathcal O\) spans \(N^*\), it permutes some triple of weights. The conjugate negative triple is distinct, so \(zB\) has an orbit of length \(6\) and gives \(12\). The ranks of norms are unchanged on descent from \(\mathbb F_4\) to \(\mathbb F_2\). This proves the exact equality of spectra.

We have \([W,E]=W\) because \(W\) is simple and nontrivial. Hence \(G_0\) is perfect. Its radical is the full preimage of \(Z(D)\) modulo \(W\rtimes N\), giving the stated \(P_0\). The commutator \([N\rtimes Z(D),E]\) is \(N\), and \([W,E]=W\), whence \(K_0=W\rtimes N\). This group is nonnilpotent because \(N\) acts faithfully and nontrivially on \(W\). The quotient \(K_0/W=N\) is abelian, so \(K_0\) itself is trivial on all chief factors above \(W\), while its image on \(W\) is \(N\ne1\). Thus this is exactly a highest internal factor of the remaining form.

In particular, \(7,9,21,30\) are absent from the spectrum of \(G_0\). The example does not have spectrum \(\Om\) and does not contradict the original classification. \end{proof}

\clearpage \subsection{Arbitrary \(7\)-factors and the global structure of \(K\)}

Here we remove the restriction that a chief \(7\)-factor be highest. Throughout this subsection, \(\omega(G)=\Om\), \(\operatorname{Soc}(G/R(G))\cong A_5\), \(L=G^\infty\), \(P=R(G)\cap L\), and \(K=[P,L]\). In particular, \(5\nmid|P|\) by Lemma~\ref{lem:A5-no-five-in-perfect-radical}.

\begin{lemma}\label{lem:perfect-symplectic-rank-four-general} Suppose a finite perfect group \(T\) has solvable radical \(A\), with \(T/A\cong A_5\) and \(5\nmid|A|\). Let \(X\) be a faithful \(T\)-module over \(\mathbb F_2\) with a nondegenerate invariant alternating form, and suppose \(C_X(T)=0\). If an element \(a\in T\) of order \(5\) outside \(A\) satisfies \(\dim[X,a]=4\), then \[
 A=1,\qquad T\cong A_5,\qquad \dim X=4.
\] \end{lemma} \begin{proof} If \(a\) is trivial on a simple composition factor \(S\) of \(X\), the kernel of the action of \(T\) on \(S\) has full image in \(T/A\). Hence the image of \(T\) on \(S\) coincides with the solvable image of \(A\) while also being perfect, so is trivial. Thus \(a\) acts nontrivially on every nontrivial simple factor. A nontrivial simple \(\mathbb F_2C_5\)-type has dimension \(4\), so \(a-1\) has rank at least \(4\) on that factor. Additivity of rank on short exact sequences, following from coprimeness of \(5\) and \(2\), shows that \(X\) has exactly one nontrivial simple factor.

A minimal submodule \(S\le X\) is nontrivial because \(C_X(T)=0\). The restriction of the form to \(S\) is either nondegenerate or zero. In the latter case \(S\le S^\perp\), and the quotient \(X/S^\perp\cong S^*\) would give a second nontrivial simple factor. Therefore \(X=S\perp S^\perp\). On \(S^\perp\) all factors are trivial; the perfect image of \(T\) on it is unitriangular and is a \(2\)-group, hence trivial. The absence of nonzero fixed vectors gives \(S^\perp=0\). Thus \(X\) is simple.

The normal \(2\)-subgroup \(O_2(A)\) is trivial on the simple module \(X\), so faithfulness gives \(O_2(A)=1\). If \(A\ne1\), its Fitting subgroup \(B=F(A)\) is nontrivial and has odd order not divisible by \(5\). Consider the simple types in the restriction of \(X\) to \(B\) after extending the field to an algebraic closure.

A nontrivial simple type of a group of odd order is not self-dual. Recall the justification of this fact. In characteristic coprime to the group order, correspondence with ordinary characters preserves duality. For an ordinary irreducible character \(\chi\), the Frobenius--Schur indicator \(\nu(\chi)=|B|^{-1}\sum_{b\in B}\chi(b^2)\) is \(0\) for a non-self-dual type, and \(1\) or \(-1\) for a self-dual type \cite{CurtisReiner}. In a group of odd order, squaring is a bijection. Thus for a nonprincipal \(\chi\) this sum is \(|B|^{-1}\sum_b\chi(b)=0\), excluding self-duality.

If \(a\) permutes five distinct isotypic spaces in the restriction to \(B\), its rank on their sum is at least \(4\). The dual types form another set of five: duality commutes with \(a\), and an involution cannot act nontrivially within a cycle of length \(5\); equality of the two sets of five would therefore give a self-dual type. The rank on the second set is also at least \(4\), contrary to the total rank \(4\). Thus \(a\) fixes all types. The kernel of the permutation action of \(T\) on the types contains \(a\notin A\). Its image in \(T/A\) is the whole group, and the permutation image of \(T\) coincides with the solvable image of \(A\). Perfection makes it trivial. Hence each isotypic space is \(T\)-invariant.

Suppose \(a\) does not centralize \(B\). Faithfulness of \(X\) allows us to choose an isotypic space \(W\) with kernel \(J\) on \(B\) such that \(a\) acts nontrivially on some Sylow \(r\)-subgroup of \(B/J\). Here \(r\) is odd and \(r\ne5\) because \(B\) is nilpotent and \(5\nmid|B|\). The group \((B_r/(B_r\cap J))\rtimes\langle a\rangle\) acts faithfully on \(W\): the normal \(r\)-group acts faithfully by the definition of the kernel, and an element of order \(5\) acting nontrivially cannot belong to its image. Lemma~\ref{lem:HH-prime} in characteristic \(2\) gives minimal polynomial \(x^5-1\). Thus the rank of \(a-1\) on \(W\) is at least \(4\). The type \(W\) is nontrivial on \(B\) and is not self-dual. The invariant form pairs \(W\) nondegenerately with a separate dual isotypic space, and the ranks of \(a-1\) on them are equal. Again the total rank is at least \(8\), a contradiction. Consequently, \([a,B]=1\).

The normal centralizer \(C_T(B)\) contains \(a\notin A\), so \(T=A C_T(B)\). The image of \(T\) on \(B\) is perfect and coincides with the solvable image of \(A\), whence \([T,B]=1\). The self-centralizing property of the Fitting subgroup of a solvable group gives \(A\le C_A(F(A))\le F(A)=B\). Thus \(A=B\le Z(T)\). The kernel of this perfect central extension of \(A_5\) is a quotient of \(M(A_5)=C_2\), while \(A\) has odd order. Hence \(A=1\). Finally, nontrivial simple \(\mathbb F_2A_5\)-modules have dimension \(4\) by the complete list in Lemma~\ref{lem:A5-orthogonal-small-five-rank}. Therefore \(\dim X=4\). \end{proof}

\Needspace{11\baselineskip} \begin{lemma}\label{lem:A5-all-seven-chief-K-central} The group \(K\) centralizes every chief \(7\)-factor of \(G\) inside \(P\), irrespective of its position in the chief series. \end{lemma} \begin{proof} Let \(V\) be such a factor on which \(K\) acts nontrivially. Put \(Y=G/C_G(V)\), \(E=Y^\infty\), and \(Q=R(E)\). Then \(E/Q\cong A_5\), \(Q\) is the image of \(P\), and the image of \(K\) is \([Q,E]\ne1\). We have \(5\nmid|Q|\), \(O_7(Y)=1\), and \(\omega(V\rtimes Y)\subseteq\Om\). Nilpotence of \(Q\) is not assumed here.

The subgroup \(F(Q)\) is normal in \(Y\). Its Sylow \(7\)-subgroup would lie in \(O_7(Y)\), and \(5\) does not divide the order of \(Q\). Hence \(F(Q)=F_2\times F_3\). Take \(a\in E\) of order \(5\) outside \(Q\). If \([F_3,a]\ne1\), Lemma~\ref{lem:HH-prime} on the faithful module \(V\) gives a nonzero norm of order \(5\) and order \(35\), which is forbidden. Thus \(a\) centralizes \(F_3\). The normal centralizer contains an element outside \(Q\); the perfect image of \(E\) on \(F_3\) coincides with the solvable image of \(Q\), so is trivial. Hence \([E,F_3]=1\).

The group \(E\) cannot also centralize \(F_2\). Otherwise \(Q\le C_Q(F(Q))\le F(Q)\), that is, \(Q=F(Q)\le Z(E)\), contrary to \([Q,E]\ne1\). The same normal-centralizer argument shows that \(a\) acts nontrivially on \(F_2\). Take a critical characteristic subgroup \(C\) of \(F_2\) using Lemma~\ref{lem:critical-two-external}. Its class is at most \(2\), \(C/Z(C)\) is elementary abelian, and \(a\) acts nontrivially on \(C\).

All simple types in the restriction of \(V\) to \(C\) are \(E\)-invariant. The proof of Lemma~\ref{lem:A5-five-seven-invariant-two-types} applies here without nilpotence of the radical: an orbit of length \(5\) gives \(35\), and the kernel of the action on types, containing \(a\notin Q\), forces the perfect permutation image to coincide with the solvable image of \(Q\). As in Lemma~\ref{lem:A5-five-seven-critical-section}, choose a type kernel \(J\trianglelefteq E\), let \(J\le C\), and obtain \(E_0=E/J\), \(Q_0=Q/J\), and \(D=C/J\). The group \(D\) is nonabelian of class \(2\), has cyclic center of order \(2\) or \(4\) centralized by \(E_0\), and \(D/Z(D)\) is elementary abelian with a nondegenerate commutator form. The action of \(E_0\) on this quotient is nontrivial and has kernel inside \(Q_0\). There is only one difference from the earlier reduction: \(Q_0\) is now solvable, rather than necessarily a \(2\)-group.

The space \(X=D/Z(D)\) has no nonzero \(E_0\)-fixed vectors: for a fixed class \(xZ(D)\), the formulas \(x^e=x f(e)\) define a homomorphism from the perfect group \(E_0\) to the central subgroup \(Z(D)\), which is trivial; then \(x\in Z(D)\). For center \(C_4\), this also gives \(\exp D=4\) by the same linear square map as in Lemma~\ref{lem:A5-C4-no-fixed-and-exponent-four}. For center \(C_2\), the group \(D\) is already extraspecial of exponent \(4\).

The trace count and Cauchy inequality from Lemma~\ref{lem:A5-C4-five-rank-and-multiplicity} give \(\dim[X,a]=4\). They require only a cyclic center, a faithful simple type of \(D\), exclusion of \(35\), and solvability of \(Q_0\): the projective image on the multiplicity space is still perfect and coincides with the solvable image of \(Q_0\), so is trivial. Thus the chosen isotypic space consists of identical genuine \(E_0\)-modules \(H\), irreducible and faithful on \(D\).

Apply the preceding lemma to the image \(T\) of \(E_0\) on \(X\). Its radical is the image of \(Q_0\), is solvable, and has order not divisible by \(5\). We obtain \(T=A_5\) and \(\dim X=4\). Thus for center \(C_2\) we have an extraspecial group \(D\) of order \(32\), and for center \(C_4\) a group \(D\) of order \(64\). In the first case \(D\) has minus type, and perfection of \(E_0\) preserves the alternating form on the four-dimensional space \(H\). Lemma~\ref{lem:extraspecial32-normalizer-fixed-four} gives an element of order \(4\) with a nonzero fixed vector. In the second case the image of \(E_0\) lies in \(\operatorname{SL}_4\), and Lemma~\ref{lem:A5-C4-normalizer-fixed-four} gives the same conclusion. A lift through the representation kernel and the \(2\)-subgroup \(J\) has order \(4\) or \(8\), since \(16\notin\Om\). The fixed vector gives \(28\) or \(56\) in the original semidirect product. This contradiction excludes the chosen \(V\). \end{proof}

\Needspace{11\baselineskip} \begin{lemma}\label{lem:central-chief-Hall-and-generated} Let \(A\trianglelefteq G\) be solvable and let \(r\) be prime. \begin{enumerate} \item If \(A\) centralizes all chief \(r\)-factors of \(G\) inside \(A\), then \(A\) has a characteristic normal Hall \(r'\)-subgroup. \item Assume in addition that \(A\le L\trianglelefteq G\) and \(A=[A,L]\). If \(L\) centralizes all chief \(r\)-factors of \(G\) inside \(A\), then \(r\nmid|A|\). \end{enumerate} \end{lemma} \begin{proof} Prove the first assertion by induction on \(|A|\). Choose a minimal normal subgroup \(N\) of \(G\) inside \(A\). It is elementary abelian of characteristic \(s\). By induction, \(A/N\) has a normal Hall \(r'\)-subgroup \(B/N\). If \(s\ne r\), then \(B\) itself is the required complement in \(A\). If \(s=r\), then \(N\le Z(A)\). In \(B\) the group \(N\) is central and is a normal Hall \(r\)-subgroup. By Schur--Zassenhaus, \(B=N\times D\), where \(D\) is a \(r'\)-group. The subgroup \(D\) is characteristic in \(B\) because it consists precisely of all \(r'\)-elements. It is normal in \(A\), and \(A/D\) is a \(r\)-group, yielding the complement. A normal Hall \(r'\)-subgroup consists precisely of all \(r'\)-elements of \(A\), so is unique and characteristic.

Prove the second assertion also by induction. For the same minimal subgroup \(N\), the image of \(A/N\) satisfies \(A/N=[A/N,L/N]\), so induction gives \(r\nmid|A/N|\). For \(s\ne r\) the proof is complete. For \(s=r\) the hypothesis gives \([N,L]=1\). Now \(N\) is a central normal Hall \(r\)-subgroup of \(A\), and \(A=N\times B\) with a characteristic \(r'\)-subgroup \(B\). Then \([A,L]=[B,L]\le B\), contradicting \(A=[A,L]\) and \(N\ne1\). \end{proof}

\Needspace{11\baselineskip} \begin{lemma}\label{lem:A5-global-seven-complement-structure} The group \(K\) has a characteristic normal Hall \(7'\)-subgroup \(H\). Moreover, \[
 K/H\text{ is a nontrivial }7\text{-group},
 \quad \pi(H)=\{2,3\},\quad H=[H,L],\quad
 6\in\omega(H),
\] \[
 O_7(K)=1,\qquad F(K)=F(H)\le H.
\] A Sylow \(3\)-subgroup of \(H\) is noncyclic. \end{lemma} \begin{proof} Lemma~\ref{lem:A5-all-seven-chief-K-central} and the first assertion of the preceding lemma give the characteristic complement \(H\). The group \(K/H\) is nontrivial because \(7\mid|K|\). Since \(5\nmid|K|\), the only prime divisors of \(H\) are \(2,3\). The divisor \(3\) occurs because the cube of every element of order \(9\) lies in \(K\) and has order \(3\). The divisor \(2\) must also occur: in the highest-factor reduction of Lemma~\ref{lem:A5-K-internal-chief-exists}, after excluding \(5,7\), either a factor of characteristic \(2\) remains or one of characteristic \(3\) with nontrivial image \(U_2\). In both cases \(2\mid|K|\), so \(2\mid|H|\). The Sylow \(3\)-subgroups of \(H,K,P\) have the same order; their noncyclicity is given by Lemma~\ref{lem:A5-odd-Sylow-radical-noncyclic}.

Put \(B=[H,L]\) and consider \(K/B\). Its normal Hall \(7'\)-subgroup \(H/B\) is central in \(L/B\). Thus \(K/B\) decomposes as a direct product of \(H/B\) and its unique Sylow \(7\)-subgroup, which is characteristic. But \(K=[K,L]\) and \(L\) centralizes \(H/B\), so \(K/B=[K/B,L/B]\) lies in this \(7\)-subgroup. Consequently, \(H/B=1\), that is, \(H=[H,L]\).

Suppose that \(6\notin\omega(H)\). Then the nontrivial group \(F(H)\) is either a \(2\)- or a \(3\)-group. If it is a \(2\)-group, the noncyclic Sylow \(3\)-subgroup of \(H\) contains \(C_3^2\), whose action on \(\Omega_1(Z(F(H)))\) gives a nonzero fixed vector for some element of order \(3\), that is, \(6\). Here we use the full theorem on a unique subgroup of prime order \cite{ConradQuaternion}: such a finite \(p\)-group is cyclic or, for \(p=2\), generalized quaternion. In particular, a noncyclic \(3\)-group contains \(C_3^2\). Thus \(F(H)\) is a \(3\)-group.

A Sylow \(2\)-subgroup of \(H\) cannot contain \(C_2^2\): on the nonzero group \(\Omega_1(Z(F(H)))\) it would give a nonzero fixed vector for some involution, that is, \(6\). By the same theorem it is cyclic or generalized quaternion. Every elementary abelian section of it has rank at most \(2\) because every subgroup is cyclic or generalized quaternion and is generated by at most two elements. Thus every chief \(2\)-factor of \(G\) inside \(H\) has dimension at most \(2\). The perfect group \(L\) centralizes them all, since \(\operatorname{GL}_1(2)\) and \(\operatorname{GL}_2(2)\) are solvable. From \(H=[H,L]\) and the second assertion of Lemma~\ref{lem:central-chief-Hall-and-generated} we obtain \(2\nmid|H|\), a contradiction. Hence \(6\in\omega(H)\).

The normal \(7\)-subgroup \(O_7(K)\) commutes with the normal \(7'\)-subgroup \(H\). If it were nontrivial, an element of order \(7\) together with an element of order \(6\) in \(H\) would give \(42\). Thus \(O_7(K)=1\). Consequently, \(F(K)\le H\) and \(F(K)\le F(H)\). The reverse inclusion follows because \(F(H)\) is characteristic in the characteristic subgroup \(H\) and is nilpotent. We obtain \(F(K)=F(H)\). \end{proof}

\Needspace{11\baselineskip} \begin{lemma}\label{lem:A5-K-Fitting-height-at-least-three} The Fitting height of \(K\) is at least \(3\). In particular, the normal Hall \(7'\)-subgroup \(H\) of the preceding lemma is nonnilpotent. \end{lemma} \begin{proof} Suppose \(K/F\) is nilpotent, where \(F=F(K)=F_2\times F_3\); the component \(7\) is absent by the preceding lemma. On every chief factor inside \(F_3\), the group \(F\) acts trivially, by the Fitting property proved in Lemma~\ref{lem:fitting-chief-kernel}: apply the same proof to the solvable normal subgroup \(K\). Thus the image \(U\) of \(K\) on such a factor is a quotient of the nilpotent group \(K/F\). If the \(7\)-part of this image is nontrivial, then \(L\) also acts nontrivially, and its perfect image \(E\) has solvable radical \(Q\) with \(E/Q\cong A_5\). The proof of the case \(p=3\) in Lemma~\ref{lem:A5-internal-seven-action-two-three} applies here without assuming the factor is highest: it uses nilpotence of \(U\), \(U=[U,E]\), the central cover \(E/U\) of \(A_5\), the inclusion of spectra, and the general properties of an action image: \(U\trianglelefteq Y\), \(O_3(Y)=1\), and \(5\nmid |U|\). These last properties follow, respectively, from normality of \(K\), faithfulness and simplicity of the chief-factor module, and \(5\nmid |P|\). All listed assumptions hold, and the proof gives \(U_7=1\), a contradiction. Hence the \(7\)-part of the image of \(K\) is trivial on every chosen factor. Every element of order \(7\) of \(K\) therefore centralizes all chief factors inside \(F_3\). Successive lifting of fixed points, using coprimeness of \(7\) and \(|F_3|\), makes it centralize all of \(F_3\).

Choose a Hall \(\{3,7\}\)-subgroup \(T\le K\). Its Sylow \(3\)-subgroup \(T_3=T\cap H\) is normal in \(T\) because \(H\trianglelefteq K\). It contains the normal subgroup \(F_3\). Nilpotence of \(K/F\) gives \([T_3,T_7]\le F\cap T=F_3\). Thus every element \(a\in T_7\) of order \(7\) acts trivially on \(T_3/F_3\) and centralizes \(F_3\). Coprime lifting of fixed points gives \(T_3=F_3C_{T_3}(a)=C_{T_3}(a)\). The group \(T_7\) has exponent \(7\), so \([T_3,T_7]=1\) and \(T=T_3\times T_7\).

Both Sylow subgroups are noncyclic by Lemma~\ref{lem:A5-odd-Sylow-radical-noncyclic}. Order \(9\) is absent from \(T_3\), since together with \(T_7\) it would give \(63\). Consequently, \(T_3\) contains \(C_3^2\) and \(T_7\) contains \(C_7^2\). If \(F_2\ne1\), their direct product acts on the nonzero elementary abelian group \(\Omega_1(Z(F_2))\). Over an algebraic closure of characteristic \(2\), every linear weight has a nontrivial kernel on \(C_3^2\) and on \(C_7^2\). This gives a common nonzero fixed vector for some commuting elements of orders \(3,7\). Together with this vector of order \(2\), they give \(42\), which is forbidden. Hence \(F_2=1\) and \(F=F_3\).

The full preimage of a Sylow \(7\)-subgroup of the nilpotent group \(K/F_3\) is normal in \(K\). It has the form \(F_3\rtimes S\), with \(S\) a \(7\)-group, by Schur--Zassenhaus. We already proved that all its elements of order \(7\) centralize \(F_3\), and \(\exp S=7\). Hence it is the direct product \(F_3\times S\); its unique Sylow \(7\)-subgroup \(S\) is characteristic and normal in \(K\). It is nontrivial, contrary to \(O_7(K)=1\). Thus \(K/F(K)\) is nonnilpotent, so the Fitting height of \(K\) is at least \(3\). If \(H\) were nilpotent, the equality \(F(K)=F(H)\) would give \(F(K)=H\) and nilpotent \(K/F(K)\), which has already been excluded. \end{proof}

\clearpage \subsection{The global \(7\)-quotient and the two central branches}

We retain the notation of the preceding subsection. In particular, \(H\) here denotes the characteristic Hall \(7'\)-subgroup in \(K\), rather than the critical example \((A_9\times C_3):2\) from the main theorem. Put \(Q=K/H\). The following restrictions concern the full quotient \(G/H\), rather than just one chief-factor action image.

\Needspace{11\baselineskip} \begin{lemma}\label{lem:A5-global-seven-quotient-split} The group \(Q\ne1\) is a normal Sylow \(7\)-subgroup of \(G/H\), and \(\exp Q=7\). There is a decomposition \[
 G/H=Q\rtimes B,\qquad
 L/H=Q\rtimes D,\qquad D=B^\infty\cong A_5
 \text{ or }2.A_5.
\] Here \(Q=[Q,D]\) and \(9\notin\omega(G/H)\), and every element of order \(5\) of \(D\) acts on \(Q\) without nontrivial fixed points. If \(d=\dim_{\mathbb F_7}Q/\Phi(Q)\), then \[
 d\equiv0\pmod4,\qquad |Q|=7^{4n},\qquad n\ge1.
\] \end{lemma} \begin{proof} The subgroup \(H\) is characteristic in \(K\), hence normal in \(G\). By Lemma~\ref{lem:remaining-commutator-localization}, all elements of order \(7\) of \(G\) lie in \(K\). Since \(49\notin\Om\), every Sylow \(7\)-subgroup has exponent \(7\) and lies entirely in \(K\). Thus \(Q\) is a normal Sylow \(7\)-subgroup of \(G/H\). It is nontrivial by Lemma~\ref{lem:A5-global-seven-complement-structure}. The normal Hall subgroup \(Q\) has a complement \(B\) by the Schur--Zassenhaus theorem.

The map \(G/H\to(G/H)/Q\cong B\) maps perfect residual onto perfect residual: successively apply the equality \(\varphi(X')=\varphi(X)'\) until the derived series stabilizes. Moreover, \(B^\infty\le(G/H)^\infty=L/H\). Thus \(D=B^\infty\) is isomorphic to \(L/K\), namely \(A_5\) or \(2.A_5\), and \(L/H=Q\rtimes D\). The equality \(K=[K,L]\) gives \(Q=[Q,Q][Q,D]\). Since \(\exp Q=7\), we have \(\Phi(Q)=Q'\). Hence the image of \([Q,D]\) generates \(Q/\Phi(Q)\); the Frattini subgroup property gives \(Q=[Q,D]\).

If \(G/H\) contained an element of order \(9\), a preimage in \(G\) would have order divisible by \(9\). The only such order in \(\Om\) is \(9\), so a preimage \(x\) of order \(9\) can be chosen. By Lemma~\ref{lem:remaining-commutator-localization}, we have \(x^3\in K\). But the image of \(x^3\) in \(K/H\), being a \(3\)-element in a \(7\)-group, is trivial. Hence \(x^3\in H\), contrary to the order \(9\) of the image of \(x\).

For \(a\in D\) of order \(5\), a nontrivial \(C_Q(a)\) would give order \(35\) in \(G/H\). This is impossible: the order of an element of a quotient divides the order of a preimage, and \(35\) divides no element of \(\Om\). Coprime lifting of fixed points also gives \(C_{Q/\Phi(Q)}(a)=0\). Nontrivial irreducible factors of \(x^5-1\) over \(\mathbb F_7\) have degree \(\operatorname{ord}_5(7)=4\). Thus the nonzero space \(Q/\Phi(Q)\) has dimension divisible by \(4\). Moreover, the action of \(\langle a\rangle\) on \(Q\setminus\{1\}\) has only orbits of length \(5\), so \(5\mid(|Q|-1)\). If \(|Q|=7^e\), then \(4\mid e\) because \(\operatorname{ord}_5(7)=4\). \end{proof}

\Needspace{11\baselineskip} \begin{lemma}\label{lem:A5-global-binary-seven-elementary} If \(D\cong2.A_5\) in the preceding lemma, then \(Q\) is elementary abelian, and the central involution \(z\in D\) acts on it by inversion. After field extension, the simple constituents of \(Q\) have only the faithful types of dimensions \(2,2,4\) from Lemma~\ref{lem:double-A5-seven-simple-data}. In particular, \[
 Q\cong C_7^{\,4m},\qquad m\ge1.
\] If \(b\in D\) has order \(3\), then \(C_Q(b)\ne1\) holds exactly when the four-dimensional absolutely simple type occurs. \end{lemma} \begin{proof} Suppose \(C=C_Q(z)\ne1\). Then \(A=\Omega_1(Z(C))\) is a nonzero elementary abelian \(7\)-group invariant under \(D\). The action of \(D\) on \(A\) factors through \(D/\langle z\rangle\cong A_5\). Choose a four-group \(T\cong C_2^2\) in \(A_5\). Over \(\mathbb F_7\) its action is diagonalizable, and every linear character \(T\to\{1,-1\}\) has a nontrivial kernel. Thus some involution \(t\in T\) fixes a nonzero vector \(v\in A\). In \(2.A_5\), the preimage of \(t\) contains an element \(u\) of order \(4\), with \(u^2=z\). Since \(z\) is trivial on \(A\), the element \(u\) fixes \(v\). We obtain \(|uv|=28\), a contradiction. Hence \(C_Q(z)=1\).

Recall the elementary argument for a fixed-point-free involution. The map \(f(q)=q^{-1}q^z\) is injective: \(f(x)=f(y)\) implies \(yx^{-1}=(yx^{-1})^z\), so \(x=y\). On the finite group \(Q\) it is surjective. But \(f(q)^z=f(q)^{-1}\), so \(z\) inverts every element of \(Q\). Inversion is an automorphism only in an abelian group: comparison of the images of \(xy\) gives \(x^{-1}y^{-1}=y^{-1}x^{-1}\). Thus \(Q\) is abelian, and \(\exp Q=7\) makes it elementary abelian.

By Maschke's theorem, the module \(Q\) is completely reducible. The equality \(z=-I\) leaves the faithful types \(2,2,4,6\) of Lemma~\ref{lem:double-A5-seven-simple-data}. The six-dimensional type has nonzero fixed vectors for an element of order \(5\) and is therefore excluded. The dimension of \(Q\) over \(\mathbb F_7\) is divisible by \(4\), by the preceding lemma. On a faithful two-dimensional type, an element of order \(3\) has determinant \(1\) and eigenvalues the two nontrivial cube roots of unity: an eigenvalue \(1\) would make it the identity, contrary to faithfulness. On the four-dimensional type, the fixed-space dimension is \(2\) by the same lemma. This proves the final criterion. \end{proof}

\Needspace{11\baselineskip} \begin{lemma}\label{lem:A5-global-ordinary-seven-minimal-rank} Let \(D\cong A_5\). Every simple constituent of \(Q/\Phi(Q)\) is isomorphic to the four-dimensional coordinate-sum-zero submodule of the permutation module on five points over \(\mathbb F_7\). If \(\dim Q/\Phi(Q)=4\), then \(Q\cong C_7^4\). Thus nonabelian \(Q\) in this branch requires at least eight generators and satisfies \(|Q|\ge7^{12}\). \end{lemma} \begin{proof} For simple \(A_5\)-modules in characteristic \(7\), the table in the proof of Lemma~\ref{lem:A5-internal-seven-action-two-three} gives dimensions \(1,3,3,4,5\). Only the four-dimensional type has no nonzero fixed vectors for an element of order \(5\). It is defined over \(\mathbb F_7\), is absolutely simple, and is realized by the stated permutation submodule. By Lemma~\ref{lem:A5-global-seven-quotient-split}, \(Q/\Phi(Q)\) has no nonzero fixed vectors for an element of order \(5\), so it is a sum of copies of this type \(V\).

Now let \(Q/\Phi(Q)=V\). Since \(\exp Q=7\), we have \(\Phi(Q)=Q'\). The commutator defines a surjective \(D\)-equivariant map \[
 \bigwedge\nolimits^2 V\longrightarrow
 Q'/\gamma_3(Q),\qquad
 (xQ')\wedge(yQ')\longmapsto[x,y]\gamma_3(Q),
\] where \(\gamma_3(Q)=[Q',Q]\). It is bilinear and skew-symmetric because commutators are central modulo \(\gamma_3(Q)\), and its image is generated by all commutators. An element of order \(5\) has no nonzero fixed vectors on \(Q'/\gamma_3(Q)\) either: this again follows by coprime lifting from \(C_Q(a)=1\). Thus all simple constituents of this quotient must also be isomorphic to \(V\).

However, \(\bigwedge^2 V\) does not contain \(V\). Verify this explicitly using ordinary characters, which apply here because \(7\nmid60\). On classes \(1,2A,3A,5A,5B\) of sizes \(1,15,20,12,12\), the character of \(V\) is \((4,0,1,-1,-1)\). The formula \(\chi_{\wedge^2 V}(g)=(\chi_V(g)^2-\chi_V(g^2))/2\) gives \((6,-2,0,1,1)\). Its inner product with \(\chi_V\) is \[
 \frac{6\cdot4+15(-2)\cdot0+20\cdot0\cdot1
       +12\cdot1(-1)+12\cdot1(-1)}{60}=0.
\] Thus the completely reducible quotient \(Q'/\gamma_3(Q)\) is trivial. Hence \(Q'=\gamma_3(Q)\). The lower central series of a finite \(7\)-group terminates at the identity, so this equality implies \(Q'=1\). Then \(Q=Q/\Phi(Q)\cong C_7^4\).

For nonabelian \(Q\), the number of generators is divisible by \(4\) but is not \(4\), so is at least \(8\). Moreover, \(\Phi(Q)=Q'\ne1\), and an element of order \(5\) acts on \(Q'\) without nontrivial fixed points. Counting orbits as in Lemma~\ref{lem:A5-global-seven-quotient-split} gives \(|Q'|\ge7^4\). Therefore \(|Q|=|Q/\Phi(Q)|\,|Q'|\ge7^8\cdot7^4=7^{12}\). \end{proof}

\Needspace{11\baselineskip} \begin{lemma}\label{lem:A5-sharp-nonabelian-seven-quotient} The bound \(7^{12}\) for nonabelian \(Q\) in the preceding lemma is sharp at the quotient level: there exists a nonabelian \(7\)-group \(Q_0\) of order \(7^{12}\), exponent \(7\), and nilpotency class \(2\), with eight generators, for which \(E_0=Q_0\rtimes A_5\) is perfect and \[
 \omega(E_0)=\{1,2,3,5,7,14,21\}\subsetneq\Om.
\] \end{lemma} \begin{proof} Let \[
 V=\{(x_1,\ldots,x_5)\in\mathbb F_7^5:
                    x_1+\cdots+x_5=0\},
\] and let \(A_5\) permute coordinates. Write \(\mathbf1\) for the vector of five ones, and define a symmetric bilinear map \[
 \mu:V\times V\longrightarrow V,\qquad
 \mu(x,y)=(x_iy_i)_{i=1}^5
       -\frac15\left(\sum_{i=1}^5x_iy_i\right)\mathbf1.
\] Here \(5\) is invertible in \(\mathbb F_7\). The map \(\mu\) is equivariant under coordinate permutations and is nonzero, for example at \(x=e_1-e_2\) and \(y=e_1-e_3\). Its values span a nonzero submodule of the simple module \(V\), hence all of \(V\).

On the set \(V\oplus V\oplus V\) define \[
 (u,v,w)(u',v',w')=
 (u+u',v+v',w+w'+\mu(u,v')).
\] Bilinearity of \(\mu\) immediately gives associativity; the identity is \((0,0,0)\), and the inverse of \((u,v,w)\) is \((-u,-v,-w+\mu(u,v))\). The resulting group \(Q_0\) has order \(7^{12}\). Its commutator is \[
 [(u,v,w),(u',v',w')]
   =(0,0,\mu(u,v')-\mu(u',v)).
\] Thus \(Q_0'=0\oplus0\oplus V\) is central and nontrivial. The formula \[
 (u,v,w)^j=(ju,jv,jw+\tbinom j2\mu(u,v))
\] gives \(\exp Q_0=7\). Hence \(\Phi(Q_0)=Q_0'\), and \(Q_0/\Phi(Q_0)\cong V\oplus V\) has dimension \(8\).

The diagonal action of \(A_5\) on the three copies of \(V\) preserves multiplication by equivariance of \(\mu\). There are no trivial constituents in \(Q_0/\Phi(Q_0)\), so \([Q_0,A_5]\) maps onto this quotient. By the Frattini property, \([Q_0,A_5]=Q_0\). Perfection of \(A_5\) now gives \(E_0'=E_0\).

An element of order \(5\) has no nonzero fixed vectors on \(V\): its fixed vectors in the permutation module are constant, and the coordinate sum of a constant vector is \(5c\). Thus it also has no nontrivial fixed points on \(Q_0\), where it acts componentwise. Order \(35\) is absent from \(E_0\): the commuting \(5\)- and \(7\)-parts of such an element would be nontrivial, while every subgroup of order \(5\) is conjugate to a subgroup of the complement \(A_5\) by conjugacy of Sylow subgroups.

For any \(g\in E_0\), its image in \(A_5\) has order \(r\in\{1,2,3,5\}\), while \(g^r\in Q_0\) has order \(1\) or \(7\). Thus only \(1,2,3,5,7,14,21,35\) are possible, and the last order has already been excluded. An involution and an element of order \(3\) in \(A_5\) each have three orbits on the five points, so their fixed spaces on \(V\) are two-dimensional. In the central third copy of \(V\) they fix nonzero elements of order \(7\), realizing \(14\) and \(21\). The orders \(1,2,3,5,7\) also occur, proving the asserted equality of spectra. \end{proof}

This model is not an example of the original group: its spectrum is a proper subset of \(\Om\), and it does not realize the nonnilpotent complement \(H\) of Lemma~\ref{lem:A5-global-seven-complement-structure}. It shows that the lower bound for \(Q\) is sharp and that the ordinary central branch cannot be excluded merely by comparing the spectrum of \(Q\rtimes A_5\).

These restrictions alone do not exclude either central branch in full. Elementary abelianness of \(K/H\) in the binary branch does not imply nilpotence of \(K\). The next subsection uses the action on \(H\) to finish excluding both branches.

\clearpage \subsection{Excluding the action on \(H\) and completing the \(A_5\) case}

In this subsection we exclude both central branches simultaneously. What is needed is not elementary abelianness of all of \(Q=K/H\), but the existence of a sufficiently large elementary abelian subgroup in a Sylow \(7\)-subgroup. This allows us to use normality of Sylow subgroups in two solvable Hall subgroups of \(K\).

\Needspace{11\baselineskip} \begin{lemma}\label{lem:seven-automorphisms-paired-elementary} Suppose an elementary abelian group \(A\) of order \(7^r\), \(r\ge3\), acts on a finite \(p\)-group \(T\), where \(p\in\{2,3\}\). The following assertions hold. \begin{enumerate} \item If \(p=3\) and \(T\) is noncyclic, then \(T\rtimes A\) contains a subgroup \(C_3^2\times C_7^2\). \item If \(p=2\) and \(T\) is neither cyclic nor generalized quaternion, then \(T\rtimes A\) contains a subgroup \(C_2^2\times C_7^2\). \end{enumerate} \end{lemma} \begin{proof} We use the following full form of the theorem on a unique subgroup of order \(p\) \cite[Theorem 4.10]{ConradQuaternion}: a finite \(p\)-group contains no \(C_p^2\) if and only if it is cyclic or, only for \(p=2\), generalized quaternion. Indeed, absence of \(C_p^2\) is equivalent to cyclicity of all abelian subgroups; this is precisely the condition in the stated theorem. Consider \(W=T/\Phi(T)\) over \(\mathbb F_p\). It is a completely reducible \(A\)-module because \((p,|A|)=1\). If \(A\) is trivial on \(W\), it is trivial on \(T\): the kernel of the map \(\Aut(T)\to\Aut(T/\Phi(T))\) is a \(p\)-group, as proved in Lemma~\ref{lem:perfect-automorphism-tools}, while the image of \(A\) is a \(7\)-group. For \(p=3\) a noncyclic \(T\) contains \(C_3^2\); for \(p=2\) a group that is neither cyclic nor generalized quaternion contains \(C_2^2\) by the theorem on a unique subgroup of prime order \cite{ConradQuaternion}. This subgroup commutes with all of \(A\), giving the required direct product.

Now suppose \(W\) has a nontrivial simple constituent \(U\). Every nontrivial simple \(\mathbb F_pA\)-module has dimension \(\operatorname{ord}_7(p)\), and the kernel of each has index \(7\) in \(A\). Indeed, over an algebraic closure \(A\) is diagonalizable. A nontrivial linear character takes values among the \(7\)th roots of unity, has image of order \(7\), and has kernel of index \(7\). On descent to \(\mathbb F_p\) its conjugate characters form a Frobenius orbit of length \(\operatorname{ord}_7(p)\); they all have the same kernel. This proves the assertion.

Denote this kernel for \(U\) by \(A_0\). Then \(A_0\cong C_7^{r-1}\) has rank at least \(2\), and \[
 U\le C_W(A_0),\qquad
 \dim U=6\ (p=3),\qquad \dim U=3\ (p=2).
\] By coprime lifting of fixed points, Lemma~\ref{lem:coprime-fixed-general}, \[
 C_W(A_0)=C_T(A_0)\Phi(T)/\Phi(T).
\] Here we use its formulation for a solvable automorphism group: \(A_0\) is elementary abelian, and \(\Phi(T)\) is a \(p\)-group. Hence \(C_T(A_0)\) has an elementary abelian quotient of rank at least \(6\) for \(p=3\) and at least \(3\) for \(p=2\).

In the first case, \(C_T(A_0)\) is noncyclic and contains \(C_3^2\). In the second, it can be neither cyclic nor generalized quaternion: such groups are generated by at most two elements, so their elementary abelian quotients have rank at most \(2\). Thus it contains \(C_2^2\). In both cases the subgroup \(C_p^2\) found commutes with \(A_0\), which contains \(C_7^2\). Coprimeness of the orders gives the required direct product. \end{proof}

\Needspace{11\baselineskip} \begin{lemma}\label{lem:A5-final-Fitting-three-group} In the nonnilpotent \(A_5\) branch, every Sylow \(7\)-subgroup of \(K\) contains \(C_7^4\), and \[
 O_2(K)=1,\qquad F(K)=O_3(K)\ne1.
\] \end{lemma} \begin{proof} By Lemma~\ref{lem:A5-global-seven-quotient-split}, the group \(Q=K/H\) is nontrivial, has exponent \(7\), and admits a fixed-point-free automorphism \(a\) of order \(5\). Its center \(Z(Q)\ne1\) is elementary abelian, \(a\)-invariant, and satisfies \(C_{Z(Q)}(a)=1\). If \(|Z(Q)|=7^e\), counting \(\langle a\rangle\)-orbits on nonidentity elements gives \(5\mid(7^e-1)\). Since \(\operatorname{ord}_5(7)=4\), the number \(e\) is positive and divisible by \(4\). Consequently, \(Z(Q)\) contains \(C_7^4\). Every Sylow \(7\)-subgroup of \(K\) is isomorphic to \(Q\): its map to \(K/H\) is injective and surjective because \(H\) is a normal Hall \(7'\)-subgroup. This proves the first assertion.

Choose a Hall \(\{3,7\}\)-subgroup \(J\le K\). Let \(T=J\cap H\). Then \(T\trianglelefteq J\) is a Sylow \(3\)-subgroup of \(K\), and \(J=T\rtimes S\), where \(S\) is a Sylow \(7\)-subgroup of \(K\). Normality of \(T\) follows from normality of \(H\), while the order assertions follow from \(\pi(H)=\{2,3\}\) and the fact that \(K/H\) is a \(7\)-group. The group \(T\) is noncyclic by Lemma~\ref{lem:A5-odd-Sylow-radical-noncyclic}. Choose \(A\cong C_7^4\) inside \(S\). The preceding lemma, applied to \(T\rtimes A\), gives a subgroup \[
 B_3\times B_7\le J,\qquad
 B_3\cong C_3^2,\quad B_7\cong C_7^2.
\]

Suppose \(N=O_2(K)\ne1\). Then \(V=\Omega_1(Z(N))\ne0\) is an elementary abelian \(2\)-group normalized by all of \(K\). Restrict the action to \(B_3\times B_7\) and extend the field to an algebraic closure of characteristic \(2\). This abelian group of odd order is diagonalizable; choose a nonzero common weight vector. The kernel of its linear character on \(B_3\cong C_3^2\) is nontrivial, as is the kernel on \(B_7\cong C_7^2\), since a finite subgroup of the multiplicative group of a field is cyclic. Thus there exist \(b\in B_3\) of order \(3\) and \(c\in B_7\) of order \(7\) with a common nonzero fixed vector. The common fixed space is defined by linear equations over \(\mathbb F_2\), and its dimension is unchanged under field extension. Hence there is a nonzero \(v\in V\) fixed by both \(b,c\). The three pairwise commuting elements \(v,b,c\) have orders \(2,3,7\) and give \(42\), a contradiction. We obtain \(O_2(K)=1\).

We already know \(O_7(K)=1\) and \(5\nmid|K|\). The nontrivial solvable group \(K\) has a nontrivial Fitting subgroup, so \(F(K)=O_3(K)\ne1\). \end{proof}

\Needspace{11\baselineskip} \begin{lemma}\label{lem:A5-final-two-Sylow-small} In this branch, a Sylow \(2\)-subgroup of \(H\) is cyclic or generalized quaternion. Consequently, every chief \(2\)-factor of \(G\) inside \(H\) has dimension at most \(2\) and is centralized by the perfect group \(L\). \end{lemma} \begin{proof} Choose a Hall \(\{2,7\}\)-subgroup \(J\le K\) and put \(T=J\cap H\). Again, \(T\trianglelefteq J\) is a Sylow \(2\)-subgroup of \(K\), and \(J=T\rtimes S\), where \(S\) contains \(A\cong C_7^4\) by the preceding lemma. If \(T\) is neither cyclic nor generalized quaternion, Lemma~\ref{lem:seven-automorphisms-paired-elementary} gives inside \(J\) a subgroup \[
 B_2\times B_7,\qquad
 B_2\cong C_2^2,\quad B_7\cong C_7^2.
\] But \(F=F(K)=O_3(K)\ne1\), and this subgroup acts on the nonzero elementary abelian \(3\)-group \(V=\Omega_1(Z(F))\). Over an algebraic closure of characteristic \(3\), choose a common weight vector. Its character has a nontrivial kernel on both \(B_2\) and \(B_7\). Descending the nonzero common fixed space to \(\mathbb F_3\) as in the preceding proof, we find pairwise commuting elements of orders \(2,7,3\) in \(K\). This gives the forbidden order \(42\). Thus \(T\) is cyclic or generalized quaternion. It is also a Sylow \(2\)-subgroup of \(H\) because \(K/H\) has odd order.

Every subgroup of a cyclic or generalized quaternion \(2\)-group is itself cyclic or generalized quaternion and is generated by at most two elements. Thus every elementary abelian section has rank at most \(2\). For a chief \(2\)-factor \(N/M\) inside \(H\), a Sylow \(2\)-subgroup of \(N\) maps onto all of \(N/M\), by the theorem on Sylow subgroups of quotients. Consequently, \(\dim_{\mathbb F_2}N/M\le2\). The image of the perfect group \(L\) on this factor is perfect and lies in the solvable group \(\operatorname{GL}_1(2)\) or \(\operatorname{GL}_2(2)\cong S_3\). It is trivial, completing the proof. \end{proof}

\Needspace{11\baselineskip} \begin{theorem}\label{thm:A5-arbitrary-radical} If \(\omega(G)=\omega(S_{10})\), then the socle of \(G/R(G)\) is not isomorphic to \(A_5\). The assertion does not require nilpotence of the solvable radical of \(G\). \end{theorem} \begin{proof} Suppose otherwise, and put \(L=G^\infty\), \(P=R(G)\cap L\), and \(K=[P,L]\). Nilpotent \(K\) was already excluded by Lemma~\ref{lem:A5-nilpotent-K-excluded}. For the remaining branch, Lemma~\ref{lem:A5-global-seven-complement-structure} gives a characteristic normal Hall \(7'\)-subgroup \(H\le K\) with the properties \[
 \pi(H)=\{2,3\},\qquad H=[H,L].
\] By the preceding lemma, \(L\) centralizes all chief \(2\)-factors of \(G\) inside \(H\). Apply the second assertion of Lemma~\ref{lem:central-chief-Hall-and-generated} to \(A=H\) and \(r=2\). Its assumptions hold: \(H\trianglelefteq G\) is solvable, \(H\le L\), \(H=[H,L]\), and all specified factors are centralized by \(L\). We obtain \(2\nmid|H|\), contrary to \(\pi(H)=\{2,3\}\). This contradiction excludes \(A_5\) in both central branches simultaneously. \end{proof}

Thus the genuine subspectral models constructed earlier do not lift to a group with full spectrum \(\Om\) and socle \(A_5\). The additional obstruction comes from the action on the characteristic subgroup \(H\), rather than from the spectrum of one local image or one quotient \(G/H\).

\clearpage

\section{\texorpdfstring{The socle \(L_2(7)\)}{The socle L2(7)}}
\subsection{The canonical kernel for \(A_6,L_2(7),L_3(4)\)}

We begin with reductions shared by these three remaining socles, then specialize to \(L_2(7)\).

\Needspace{11\baselineskip}
\begin{lemma}\label{lem:commutator-generated-Sylow-ranks}
Let \(A\trianglelefteq G\) be solvable, with \(A\le L\trianglelefteq G\), \(L=L'\), and \(A=[A,L]\). Then, for every odd prime \(r\mid|A|\), a Sylow \(r\)-subgroup of \(A\) is noncyclic. If \(2\mid|A|\), a Sylow \(2\)-subgroup of \(A\) is neither cyclic nor generalized quaternion. In particular, these assertions apply to \(A=K\).
\end{lemma}
\begin{proof}
If a Sylow \(r\)-subgroup is cyclic, every chief \(r\)-factor of \(G\) inside \(A\) has order \(r\): it is an elementary abelian section of a cyclic group. The perfect group \(L\) acts trivially on it, since \(\operatorname{GL}_1(r)\) is abelian. The second assertion of Lemma~\ref{lem:central-chief-Hall-and-generated} now gives \(r\nmid|A|\), a contradiction. If a Sylow \(2\)-subgroup is cyclic or generalized quaternion, its elementary abelian sections have rank at most \(2\). Thus the chief \(2\)-factors inside \(A\) have dimension at most \(2\). The groups \(\operatorname{GL}_1(2)\) and \(\operatorname{GL}_2(2)\) are solvable, so the perfect group \(L\) centralizes these factors. The same lemma with \(r=2\) yields a contradiction.
\end{proof}

\Needspace{11\baselineskip}
\begin{lemma}\label{lem:remaining-no-both-five-seven}
For each of the three socles under consideration, \(|P|\) is not divisible by both \(5\) and \(7\). In particular,
\[
 S=A_6\ \Longrightarrow\ 5\nmid|P|,
 \qquad
 S=L_2(7)\ \Longrightarrow\ 7\nmid|P|.
\]
\end{lemma}
\begin{proof}
The group \(P/K\) is the center of a perfect central cover of \(S\); in particular, it is a quotient of \(M(S)\). The multipliers of the groups under consideration are \(C_6,C_2,C_3\times C_4\times C_4\), respectively. Consequently, the \(5\)- and \(7\)-parts of the orders of \(P,K\) coincide. If both primes divide \(|P|\), they divide \(|K|\), and the corresponding Sylow subgroups are noncyclic by the preceding lemma.

Choose a Hall \(\{5,7\}\)-subgroup \(J\le K\). Its nontrivial Fitting subgroup \(F(J)\) cannot have order divisible by both primes, since this would give \(35\). Hence \(F(J)\) is an \(r\)-group for one of the primes \(r\in\{5,7\}\). For the other prime \(s\), a noncyclic Sylow \(s\)-subgroup contains \(C_s^2\). Its action on the nonzero group \(\Omega_1(Z(F(J)))\) is diagonalizable after extending the field. The kernel of every linear character of \(C_s^2\) is nontrivial; thus some element of order \(s\) fixes a nonzero vector of order \(r\). This gives \(35\), a contradiction. For \(S=A_6\), we already know that \(7\mid|K|\), and for \(S=L_2(7)\), that \(5\mid|K|\), by Lemma~\ref{lem:remaining-commutator-localization}. The two stated consequences follow.
\end{proof}

\Needspace{11\baselineskip}
\begin{lemma}\label{lem:A6-cover-no-nine-general}
Let \(\omega(B)\subseteq\Om\), and suppose that \(B\) contains a normal perfect subgroup \(D\) with \(D/Z(D)\cong A_6\). Then \(9\notin\omega(B)\). Consequently, in the original group with socle \(A_6\), every element satisfying \(|x|=9\) also satisfies \(x^3\in K\). In particular, \(3\mid|K|\). If \(K\) is nilpotent, then
\[
 K=K_3\times K_7,\qquad
 K_3,K_7\ne1,\qquad
 \exp K_3=3,\quad\exp K_7=7,
\]
and both Sylow subgroups are noncyclic.
\end{lemma}
\begin{proof}
We first show that all Sylow \(3\)-subgroups of any perfect central cover of \(A_6\) have exponent \(3\). Factor out the central \(2\)-part; this preserves the orders of \(3\)-elements. The remaining central kernel \(Z\) has order \(1\) or \(3\), since \(M(A_6)=C_6\). When \(Z=1\), the assertion is immediate from the realization of \(A_6\) as permutations of six points.

When \(|Z|=3\), take the full preimage \(U\) of the Sylow subgroup
\[
 \overline U=\langle c,d\rangle,
 \qquad c=(1\,2\,3),\quad d=(4\,5\,6).
\]
Then \(|U|=27\), \(U/Z\cong C_3^2\), and \([U,U]\le Z\). The map
\[
 f:U/Z\longrightarrow Z,\qquad f(uZ)=u^3
\]
is well defined and \(\mathbb F_3\)-linear: centrality of \(Z\), its exponent \(3\), and the cube formula in a group of class at most \(2\) give \((uv)^3=u^3v^3[v,u]^3=u^3v^3\). The permutation \(n=(1\,4\,2\,5)(3\,6)\in A_6\) normalizes \(\overline U\), with \(c^n=d^{-1}\), \(d^n=c\). Its matrix on \(\overline U\) has characteristic polynomial \(t^2+1\), so \(n-1\) is invertible over \(\mathbb F_3\). Every lift of \(n\) to the cover normalizes \(U\) and fixes \(Z\) pointwise. Thus \(f\circ n=f\), whence \(f=0\). Therefore \(u^3=1\) for every \(u\in U\).

Now suppose that \(|x|=9\) in \(B\). Its action on \(A_6\) is inner, since \(\Out(A_6)\cong C_2^2\), and has order \(1\) or \(3\). It can be represented by conjugation by an element \(a\in D\), \(a^3=1\): the required element of \(A_6\) lifts to a Sylow \(3\)-subgroup of \(D\), which has just been shown to have exponent \(3\). The map \(\Aut(D)\to\Aut(A_6)\) is injective by Lemma~\ref{lem:perfect-automorphism-tools}. Hence \(x,a\) induce the same action on \(D\), \([x,a]=1\), and \(y=xa^{-1}\) centralizes \(D\) and has order \(9\), since \(y^3=x^3\ne1\). An involution in \(D\) commutes with \(y\) and gives \(18\), a contradiction.

Apply this to \(B=G/K\), \(D=L/K\). The image of every element \(x\) of order \(9\) then has order at most \(3\), that is, \(x^3\in K\). In the nilpotent branch, Lemma~\ref{lem:nilpotent-K-prime-shapes} already gives \(K=K_2\times K_3\times K_7\), \(K_7\ne1\), \(\exp K_3\mid3\), \(\exp K_7=7\), with \(K_2\) and \(K_3\) not both nontrivial. Now \(K_3\ne1\), so \(K_2=1\). Noncyclicity follows from Lemma~\ref{lem:commutator-generated-Sylow-ranks}.
\end{proof}

\Needspace{11\baselineskip}
\begin{lemma}\label{lem:L27-L34-central-five-seven-chief}
For \(S=L_2(7)\) or \(L_3(4)\), the group \(K\) centralizes all chief factors of \(G\) of characteristics \(5,7\) inside \(P\). In particular, \(K\) has characteristic normal Hall \(5'\)- and \(7'\)-subgroups.
\end{lemma}
\begin{proof}
Suppose that \(K\) acts nontrivially on such a factor \(V\). Then the actions of both \(R(G)\) and \(L\) are nontrivial, so this is a mixed factor. In the notation of Lemma~\ref{lem:mixed-chief-action}, its perfect image \(E\) has radical \(Q\), equal to the image of \(P\), and \(E/Q\cong S\). Lemma~\ref{lem:odd-mixed-quasisimple} applies, since \(7\mid|S|\) and the characteristic is \(5\) or \(7\). It gives \(Q\le Z(E)\), so the image of \(K=[P,L]\) is \([Q,E]=1\), a contradiction. Thus \(K\) centralizes all the specified factors, in particular those inside \(K\) itself. The first assertion of Lemma~\ref{lem:central-chief-Hall-and-generated} with \(r=5,7\) gives both Hall subgroups. Nilpotence of \(K\) is not assumed here.
\end{proof}

\Needspace{11\baselineskip}
\begin{lemma}\label{lem:L27-global-five-complement}
For \(S=L_2(7)\), there is a characteristic normal Hall \(5'\)-subgroup \(H\) of \(K\), with
\[
 \pi(H)\subseteq\{2,3\},\qquad
 3\mid|H|,\qquad H=[H,L].
\]
A Sylow \(3\)-subgroup of \(H\) is noncyclic. For \(Q=K/H\), we have
\[
 G/H=Q\rtimes B,\quad L/H=Q\rtimes D,\quad
 D=B^\infty\cong L_2(7)\text{ or }2.L_2(7).
\]
Here \(Q\ne1\) is a normal Sylow \(5\)-subgroup of \(G/H\), \(\exp Q=5\), \(Q=[Q,D]\), and \(9\notin\omega(G/H)\). Moreover,
\[
 |Q|=5^{6m},\quad |Z(Q)|=5^{6n},\quad m,n\ge1,
\]
and every Sylow \(5\)-subgroup of \(K\) contains \(C_5^6\).
\end{lemma}
\begin{proof}
The preceding lemma gives the Hall subgroup \(H\), and Lemma~\ref{lem:remaining-no-both-five-seven} gives \(7\nmid|P|\). By Lemma~\ref{lem:remaining-commutator-localization}, \(15\mid|K|\), whence \(Q\ne1\) and \(3\mid|H|\). Put \(J=[H,L]\). The subgroup \(H/J\) is central in \(L/J\) and is a normal Hall subgroup of \(K/J\). Hence \(K/J\) is the direct product of \(H/J\) and its unique Sylow \(5\)-subgroup. The latter is characteristic and normalized by \(L/J\). Since \(K/J=[K/J,L/J]\), we obtain \(H/J=1\), that is, \(H=[H,L]\). Noncyclicity of a Sylow \(3\)-subgroup follows from Lemma~\ref{lem:commutator-generated-Sylow-ranks}.

All elements of order \(5\) in \(G\) lie in \(K\), and \(25\notin\Om\). Thus \(Q\) is a normal Sylow \(5\)-subgroup of \(G/H\), of exponent \(5\). The Schur--Zassenhaus theorem gives a complement \(B\). As in the proof of Lemma~\ref{lem:A5-global-seven-quotient-split}, projection of the perfect residual gives \(D=B^\infty\cong L/K\) and \(L/H=Q\rtimes D\). The multiplier of \(L_2(7)\) has order \(2\), giving the two covers specified. The equality \(Q=[Q,Q][Q,D]\) and \(\Phi(Q)=Q'\) give \(Q=[Q,D]\). If \(G/H\) had an element of order \(9\), a preimage in \(G\) would have order \(9\). But the cube of this preimage lies in \(K\) and has trivial image in the \(5\)-group \(K/H\), so it lies in \(H\), a contradiction.

Take an element \(a\in D\) of order \(7\). We have \(C_Q(a)=1\), since otherwise \(35\) occurs. The same holds for \(Z(Q)\). All nontrivial orbits of \(\langle a\rangle\) on \(Q\) and \(Z(Q)\) have length \(7\). If \(|Q|=5^u\), \(|Z(Q)|=5^v\), then \(7\mid(5^u-1)\), \(7\mid(5^v-1)\). Since \(\operatorname{ord}_7(5)=6\), both exponents are positive multiples of \(6\). The center of \(Q\) is elementary abelian, since \(\exp Q=5\), and contains \(C_5^6\). The Sylow \(5\)-subgroups of \(K\) are isomorphic to \(Q\), completing the proof.
\end{proof}

\Needspace{11\baselineskip}
\begin{lemma}\label{lem:L34-global-two-three-complement}
For \(S=L_3(4)\), the group \(K\) has a characteristic normal Hall \(\{2,3\}\)-subgroup \(H\), with
\[
 H=[H,L],\qquad 3\mid|H|.
\]
Its Sylow \(3\)-subgroup is noncyclic. There are only three possibilities for \(Q=K/H\): \(Q=1\), a nontrivial \(5\)-group, or a nontrivial \(7\)-group. In the two nontrivial cases, respectively,
\[
 |Z(Q)|=5^{6n}\quad\text{or}\quad
 |Z(Q)|=7^{4n},\qquad n\ge1.
\]
Consequently, if \(5\mid|K|\), a Sylow \(5\)-subgroup of \(K\) contains \(C_5^6\), and if \(7\mid|K|\), a Sylow \(7\)-subgroup of \(K\) contains \(C_7^4\).
\end{lemma}
\begin{proof}
By Lemma~\ref{lem:remaining-no-both-five-seven}, at most one of the primes \(5,7\) occurs in \(K\). If \(r\in\{5,7\}\) occurs, Lemma~\ref{lem:L27-L34-central-five-seven-chief} gives a characteristic normal Hall \(r'\)-subgroup \(H\), which here is a \(\{2,3\}\)-group. If neither occurs, put \(H=K\). The equality \(H=[H,L]\) follows by the same argument with a central Hall subgroup as in the preceding lemma; for \(H=K\), it is already known. By Lemma~\ref{lem:remaining-commutator-localization}, \(3\mid|K|\), so \(3\mid|H|\). Noncyclicity follows from Lemma~\ref{lem:commutator-generated-Sylow-ranks}.

Let \(Q\ne1\) be an \(r\)-group, and let \(t\) be the other prime in \(\{5,7\}\). In \(L/H\), choose an element \(a\) of order \(t\) with nontrivial image in \(S\). Such an element is obtained by taking the \(t\)-part of a lift of an element of order \(t\) in \(S\): the squares \(5^2,7^2\) are forbidden in the spectrum. Since \(r^2\notin\Om\), we have \(\exp Q=r\). The normal subgroup \(Q\) is invariant under \(a\), and \(C_Q(a)=1\), since otherwise \(35\) occurs. Its center is nontrivial, elementary abelian, and also has no nontrivial fixed points under \(a\). Counting orbits gives \(t\mid(|Z(Q)|-1)\). Now \(\operatorname{ord}_7(5)=6\) and \(\operatorname{ord}_5(7)=4\) give both formulas. A Sylow \(r\)-subgroup of \(K\) is isomorphic to \(Q\), since \(H\) is a normal Hall \(r'\)-subgroup. No splitting of \(L/H\) over \(Q\) is used here: unlike in the preceding lemma, the prime \(r\) also divides \(|S|\).
\end{proof}

\Needspace{11\baselineskip}
\begin{lemma}\label{lem:L34-seven-branch-K-odd}
For \(S=L_3(4)\) and \(7\mid|K|\), the group \(K\) has odd order, and
\[
 K=H\rtimes T,\qquad
 H=O_3(K)=[H,L]\ne1,
\]
where \(H\) is a noncyclic \(3\)-group and \(T\) is a \(7\)-group of exponent \(7\) containing \(C_7^4\).
\end{lemma}
\begin{proof}
Use \(H\) from the preceding lemma. At this point, we know only that \(\pi(H)\subseteq\{2,3\}\) and that its Sylow \(3\)-subgroup is noncyclic. In a Hall \(\{3,7\}\)-subgroup \(J\le K\), the subgroup \(J\cap H\) is a normal noncyclic Sylow \(3\)-subgroup, and a Sylow \(7\)-subgroup contains \(C_7^4\). Lemma~\ref{lem:seven-automorphisms-paired-elementary} gives \(C_3^2\times C_7^2\le J\). If \(O_2(K)\ne1\), the action of this direct product on \(\Omega_1(Z(O_2(K)))\) gives \(42\), exactly as in the proof of Lemma~\ref{lem:A5-final-Fitting-three-group}: a common linear weight has nontrivial kernels on both factors. Therefore \(O_2(K)=1\).

Since \(H\) is characteristic in \(K\), we have \(O_2(H)\le O_2(K)=1\). The nontrivial solvable group \(H\) therefore has \(F(H)=O_3(H)\ne1\). If its Sylow \(2\)-subgroup is neither cyclic nor generalized quaternion, the same lemma in a Hall \(\{2,7\}\)-subgroup of \(K\) gives \(C_2^2\times C_7^2\). The action of this group on the nonzero group \(\Omega_1(Z(F(H)))\) in characteristic \(3\) again gives a common nonzero fixed vector for elements of orders \(2,7\), that is, \(42\). Hence a Sylow \(2\)-subgroup of \(H\) is cyclic or generalized quaternion. But \(H=[H,L]\), so Lemma~\ref{lem:commutator-generated-Sylow-ranks} implies \(2\nmid|H|\).

Thus \(H\) is a \(3\)-group and \(K/H\) is a \(7\)-group; both are nontrivial. The Schur--Zassenhaus theorem gives the stated splitting. The normal Hall \(3\)-subgroup \(H\) equals \(O_3(K)\). Noncyclicity of \(H\), the exponent of \(T\), and the presence of \(C_7^4\) have already been proved. This restriction applies to the prime divisors of the order of \(K\); it need not hold for the orders of \(P\) or \(G\).
\end{proof}

If \(K\) is not nilpotent, a chief series for it has a factor on which \(K\) itself acts nontrivially. The following restrictions on the characteristics of such factors have been obtained:
\[
\begin{array}{c|c}
 S&\text{possible characteristics}\\\hline
 A_6&2,3,7\\
 L_2(7)&2,3\\
 L_3(4)&2,3.
\end{array}
\]
The first row follows from \(5\nmid|P|\), and the other two from Lemma~\ref{lem:L27-L34-central-five-seven-chief}. These restrict the action of \(K\); they do not exclude the chief factors themselves.

\clearpage
\subsection{Excluding the double cover in the \(L_2(7)\) case}

In this subsection, \(S=L_2(7)\), and \(H,Q,B,D\) have the meanings given in Lemma~\ref{lem:L27-global-five-complement}. We exclude the entire branch \(D\cong2.L_2(7)\), without assuming that \(H\) or \(K\) is nilpotent.

\Needspace{11\baselineskip}
\begin{lemma}\label{lem:L27-seven-free-five-modules-four}
Let \(W\ne0\) be a module for \(S=L_2(7)\) over \(\mathbb F_5\), on which an element of order \(7\) has no nonzero fixed vectors. Then every simple constituent of \(W\) over the algebraic closure has dimension \(3,3\), or \(6\). An element of order \(4\) has fixed-space dimension \(1,1\), or \(2\), respectively, on each such constituent. In particular, it fixes a nonzero vector of \(W\).
\end{lemma}
\begin{proof}
The full list of degrees of absolutely simple modules for \(S\) in characteristic \(5\) is \(1,3,3,6,7,8\), as already used in Lemma~\ref{lem:nilpotent-tables}; see \cite{NilpotentData,L27BinaryData}. Since the characteristic \(5\) does not divide \(168\), the module is completely reducible. The normalizer of a subgroup of order \(7\) contains \(7:3\), where the element of order \(3\) has orbits of length \(3\) on the six nontrivial linear characters of \(C_7\). Thus, for any module of dimension \(d\), the number \(d-\dim C_W(7)\) is divisible by \(3\). Dimensions \(1,7,8\) therefore do not allow \(C_W(7)=0\); only \(3,3,6\) remain.

We check the fixed points of an element of order \(4\) without further table values. In \(S\cong\operatorname{GL}_3(2)\), an element of order \(4\) has a single Jordan block of size \(3\), as does its inverse; hence they are conjugate. On a three-dimensional simple module in characteristic \(5\), the eigenvalues \(i,-i\), \(i^2=-1\), have equal multiplicities. The determinant is \(1\), since \(S\) is perfect; thus the multiplicity of \(-1\) is even, and that of \(1\) is odd. The latter equals \(1\): the value \(3\) would make the element of order \(4\) act as the identity on a nontrivial simple module for a simple group.

The unique six-dimensional type is the zero-coordinate-sum submodule of the permutation module on the seven nonzero vectors of \(\mathbb F_2^3\). The action of \(S\) on them is doubly transitive, so this submodule is simple: the endomorphism algebra of the permutation module has dimension \(2\), and the completely reducible module already contains one trivial constituent. An element of order \(4\) has cycles of lengths \(4,2,1\) on the seven points. The fixed space of the permutation module has dimension three; removing the constant constituent leaves dimension \(2\). This proves all assertions.
\end{proof}

\Needspace{11\baselineskip}
\begin{lemma}\label{lem:L27-binary-Q-elementary-three-fixed}
If \(D\cong2.L_2(7)=\operatorname{SL}_2(7)\), then \(Q\) is elementary abelian, and the central involution \(z\in D\) acts on \(Q\) by inversion. For every \(b\in D\) of order \(3\),
\[
 \dim_{\mathbb F_5}C_Q(b)
       =\frac13\dim_{\mathbb F_5}Q\ge2.
\]
\end{lemma}
\begin{proof}
Suppose that \(C=C_Q(z)\ne1\). Then \(A=\Omega_1(Z(C))\ne0\) is \(D\)-invariant, and the action of \(D\) on \(A\) factors through \(D/\langle z\rangle\cong L_2(7)\). Every element of order \(7\) acts without nonzero fixed vectors on \(A\): it has no nontrivial fixed points on \(Q\), since otherwise \(35\) would occur. The preceding lemma gives a nonzero vector \(v\in A\) fixed by an element of order \(4\) in the quotient. Every lift \(u\in D\) of that element has order \(8\). Indeed, the unique involution of \(\operatorname{SL}_2(7)\) is \(-I=z\); if \(u^4=1\), then \(u^2=z\), and the image of \(u\) would have order at most \(2\). Thus \(u^4=z\), \(|u|=8\). Since \(z\) is trivial on \(A\), the element \(u\) fixes \(v\) and gives \(40\in\omega(L/H)\), a contradiction. Hence \(C_Q(z)=1\). The argument with the map \(q\mapsto q^{-1}q^z\), given in full in the proof of Lemma~\ref{lem:A5-global-binary-seven-elementary}, shows that \(z\) inverts \(Q\) and that \(Q\) is abelian. Since \(\exp Q=5\), it is elementary abelian.

Every element of order \(3\) in \(\operatorname{SL}_2(7)\) is conjugate to \(b=\operatorname{diag}(2,4)\): its eigenvalues are the two nontrivial cube roots of unity in \(\mathbb F_7\), and an eigenbasis can be chosen with determinant \(1\). The element \(b\) normalizes the subgroup
\[
 U=\left\{\begin{pmatrix}1&t\\0&1\end{pmatrix}:
                        t\in\mathbb F_7\right\}\cong C_7
\]
and permutes its nontrivial characters in orbits of length \(3\). There are no trivial \(U\)-weights on \(Q\). Each triple of equal-dimensional weight spaces cyclically permuted by \(b\) gives a fixed space of that dimension: its vectors are \(v+bv+b^2v\), with \(v\) in the first space. Thus \(\dim C_Q(b)=\dim Q/3\). Extending the field does not change this dimension, and \(\dim Q\ge6\) by Lemma~\ref{lem:L27-global-five-complement}.
\end{proof}

\Needspace{11\baselineskip}
\begin{lemma}\label{lem:normal-SL27-rigid-overgroup}
Suppose that \(\omega(B)\subseteq\Om\) and \(D\trianglelefteq B\) is isomorphic to \(\operatorname{SL}_2(7)\). Then
\[
 B=D\times E,\qquad E\cong C_2^t\quad(t\ge0).
\]
\end{lemma}
\begin{proof}
Put \(C=C_B(D)\). The group \(D\) contains elements of orders \(8,7,4\). Hence \(C\) has no elements of orders \(3,5,7\): these would give \(24,35,28\), respectively. By Cauchy's theorem, \(C\) is a \(2\)-group. It has no element of order \(4\), since together with a \(7\)-element of \(D\) this would give \(28\). Thus \(C\) is elementary abelian.

Conjugation gives an image \(I\) of \(B\) in \(\Aut(D)\) containing \(\operatorname{Inn}(D)\cong L_2(7)\). By Lemma~\ref{lem:perfect-automorphism-tools}, there is an embedding \(\Aut(D)\to\Aut(L_2(7))\cong\operatorname{PGL}_2(7)\), and \(\operatorname{Inn}(D)\cong L_2(7)\) has index \(2\) in \(\operatorname{PGL}_2(7)\) \cite{L27BinaryData}. If \(I\) contains an outer automorphism, it equals the whole of \(\operatorname{PGL}_2(7)\). This group contains an element \(\beta\) of order \(8\): a nonsplit projective torus is the cyclic group \(\mathbb F_{49}^{\times}/\mathbb F_7^{\times}\) of order \(8\). Its square is an inner automorphism of order \(4\), since \(L_2(7)\) has no element of order \(8\).

Choose a lift \(y\in B\) of \(\beta\). Its order is divisible by \(8\), and the only such order in \(\Om\) is \(8\). Let \(a\in D\) represent the inner automorphism \(\beta^2\). As in the preceding proof, \(|a|=8\), \(a^4=z\). The actions of \(y^2\) and \(a\) on \(D\) coincide, so \(y^2=ac\) for some \(c\in C\). But \([a,c]=1\), \(c^2=1\), whence
\[
 1=y^8=(ac)^4=a^4=z\ne1,
\]
a contradiction. Thus \(I\) contains only inner automorphisms, and \(B=DC\). Since \(D\cap C=Z(D)=\langle z\rangle\), choose a complement \(E\) to \(\langle z\rangle\) in the elementary abelian group \(C\). Then \(B=D\times E\), as required.
\end{proof}

\Needspace{11\baselineskip}
\begin{lemma}\label{lem:L27-binary-five-Sylow-normal}
In the branch \(L/K\cong2.L_2(7)\), every Sylow \(5\)-subgroup \(T\) of \(K\) centralizes \(H\). Hence
\[
 K=H\times T,
\]
and \(T\) is an elementary abelian normal Sylow \(5\)-subgroup of the whole group \(G\).
\end{lemma}
\begin{proof}
The subgroup \(T\) is isomorphic to \(Q=K/H\), and is therefore elementary abelian. Fix a chief series of \(G\) inside \(H\). If \(T\) acts nontrivially on one of its factors, choose a highest such factor \(V=N/M\). The group \(T\) is trivial on all factors above \(N\). Successive coprime lifting of fixed points by Lemma~\ref{lem:coprime-fixed-general}, using \((|T|,|H|)=1\), gives \([T,H]\le N\). Consequently,
\[
 K/N=(H/N)\times(TN/N).
\]
The subgroup \(TN/N\) is the unique Sylow \(5\)-subgroup of \(K/N\), and is therefore normal in \(G/N\). The projection \(TN/N\to K/H=Q\) is an isomorphism compatible with the action of \(G\).

Put \(Y=G/C_G(V)\). Since \(N\le C_G(V)\), the image \(A\) of \(T\) is normal in \(Y\). It is nontrivial by the choice of the factor, elementary abelian, and of order coprime to \(p=\operatorname{char}V\in\{2,3\}\). The module \(V\) is faithful and simple for \(Y\). Hence \(C_V(A)=0\): this space is \(Y\)-invariant, and equality with \(V\) would contradict faithfulness.

Take an element \(d\in D\le L/H\) of order \(8\), \(d^4=z\), and a lift \(g\in L\). Since \(8\) is the only order in \(\Om\) divisible by \(8\), we have \(|g|=8\). On \(TN/N\cong Q\), the element \(g^4\) acts by inversion, by Lemma~\ref{lem:L27-binary-Q-elementary-three-fixed}. Thus it also inverts the image \(A\). All \(A\)-weights on \(V\) are nontrivial. The action of \(g^4\) replaces each weight \(\lambda\) by the distinct weight \(\lambda^{-1}\). All orbits of \(\langle g\rangle\) on the weights therefore have length \(8\): orbits of lengths \(1,2,4\) would be fixed by \(g^4\). For a nonzero vector \(v\) in one weight space, the sum
\[
 (1+g+\cdots+g^7)v
\]
is nonzero, since its summands belong to eight distinct weight spaces. The norm operator remains nonzero on descent to the ground field \(\mathbb F_p\).

The image of \(g\) in \(Y\) has order \(8\), since its fourth power is nontrivial on \(A\). Thus, in \(G/M\), the subgroup \(\langle V,gM\rangle\) is an actual semidirect product \(V\rtimes\langle gM\rangle\): \(gM\) has order \(8\), and its cyclic subgroup intersects \(V\) trivially. Calculating the eighth power of \(v(gM)\) using the nonzero norm gives order \(16\) when \(p=2\) and \(24\) when \(p=3\). Both are forbidden in a section of \(G\). The contradiction excludes the chosen factor.

Thus \(T\) centralizes all chief factors inside \(H\). Successive coprime lifting once more gives \([T,H]=1\). Hence \(K=H\times T\). The subgroup \(T\) is now the unique Sylow \(5\)-subgroup of \(K\), characteristic in \(K\), and normal in \(G\). It is a Sylow \(5\)-subgroup of \(G\), since all elements of order \(5\) in \(G\) lie in \(K\) and \(25\notin\Om\).
\end{proof}

\Needspace{11\baselineskip}
\begin{theorem}\label{thm:L27-binary-cover-excluded}
If \(\omega(G)=\Om\) and \(\Soc(G/R(G))\cong L_2(7)\), then
\[
 L/K\cong L_2(7).
\]
The double cover \(L/K\cong2.L_2(7)\) is impossible without any restriction on the solvable radical.
\end{theorem}
\begin{proof}
Suppose that the double cover occurs. By Lemma~\ref{lem:L27-global-five-complement}, \(G/H=Q\rtimes B\), and \(D=B^\infty\cong2.L_2(7)\) is normal in \(B\). Lemma~\ref{lem:normal-SL27-rigid-overgroup} gives \(B=D\times E\), where \(E\) is elementary abelian of \(2\)-power order. By the preceding lemma, \(K=H\times T\), and \(T\ne1\) is normal in \(G\), elementary abelian, and isomorphic to \(Q\). The action of \(G\) on \(T\) factors through \(G/H\).

Choose \(x\in G\) of order \(9\). Its image \(\overline x\in G/H\) has order \(1\) or \(3\). In the first case, \(x\in H\), so it centralizes \(T\), immediately giving \(45\). In the second case, conjugation in \(G/H\) moves \(\overline x\) into \(B\): the group \(B\) contains a Sylow \(3\)-subgroup of \(G/H\), so Sylow's theorem applies. Lifting the conjugating element to \(G\), we may assume \(\overline x\in B\). Since \(B=D\times E\), its element \(\overline x\) of order \(3\) lies in \(D\). Lemma~\ref{lem:L27-binary-Q-elementary-three-fixed} gives \(C_T(\overline x)\ne1\). But the action of \(x\) on \(T\) coincides with that of its image \(\overline x\). Thus \(x\) centralizes an element of order \(5\), again giving \(45\), a contradiction. Only the simple cover \(L/K\cong L_2(7)\) remains.
\end{proof}

\clearpage
\subsection{Excluding the outer extension of \(L_2(7)\)}

By Theorem~\ref{thm:L27-binary-cover-excluded}, only \(L/K\cong L_2(7)\) remains. In this subsection we also exclude \(G/R(G)\cong\operatorname{PGL}_2(7)\). The notation \(H,Q,B,D\) from Lemma~\ref{lem:L27-global-five-complement} is retained, but now \(D\cong L_2(7)\).

\Needspace{11\baselineskip}
\begin{lemma}\label{lem:PGL27-seven-free-module-eight-fixed}
Let \(X=\operatorname{PGL}_2(7)\), \(D=X'\cong L_2(7)\), and let \(W\ne0\) be an \(X\)-module over \(\mathbb F_5\), on which elements of order \(7\) have no nonzero fixed vectors. Then every element \(y\in X\) of order \(8\) has a nonzero fixed space on \(W\). On every absolutely simple constituent of \(W\), the fixed space has dimension \(1\).
\end{lemma}
\begin{proof}
Extend the field to an algebraic closure \(k\). The restriction to \(D\) is completely reducible, and by Lemma~\ref{lem:L27-seven-free-five-modules-four} its simple constituents have types \(3,3,6\). The outer automorphism interchanges the two three-dimensional types. This can be seen on the group of order \(7\): their weights form the two sets \(\{1,2,4\}\), \(\{3,5,6\}\) of exponents of nontrivial seventh roots of unity. These are precisely the two orbits of the normalizer \(7:3\); passage to the dual module interchanges them. Outer conjugation in \(\operatorname{PGL}_2(7)\) multiplies the exponents by a quadratic nonresidue and also interchanges these sets.

Thus a simple \(kX\)-module having a three-dimensional constituent on restriction to \(D\) is the sum of two conjugate three-dimensional constituents and has dimension \(6\): the sum of a constituent and its outer conjugate is already \(X\)-invariant. Every outer element acts on it by a block matrix \(\left(\begin{smallmatrix}0&A\\B&0\end{smallmatrix}\right)\). Its determinant is \(-\det(AB)=-1\), since the square of an outer element lies in the perfect group \(D\), and the determinant of its action on a three-dimensional constituent is \(1\).

Consider the six-dimensional type for \(D\). We explicitly construct its extensions to \(X\) and calculate their determinants. Identify \(D\) with \(\operatorname{GL}_3(2)\), and index seven coordinates by the nonzero vectors \(u\in\mathbb F_2^3\). Let \(M\) be the integral matrix
\[
 M_{u,v}=\begin{cases}1,&u\mathbin{\cdot}v=0
                         \text{ in }\mathbb F_2,\\
                       0,&\text{otherwise}.
          \end{cases}
\]
Every row contains three ones, and two distinct rows share one common one. Consequently,
\[
 M^2=2I+J,
\]
where \(J\) is the all-ones matrix. The trace of \(M\) is \(3\), since \(\mathbb F_2^3\) has exactly three nonzero vectors of even weight. On the constant constituent, \(M\) has eigenvalue \(3\), so on the six-dimensional zero-coordinate-sum submodule its trace is \(0\) and its square is \(2I\). The integral characteristic polynomial of this restriction is \((t^2-2)^3\), and its determinant is \(-8\).

The matrix \(M/\sqrt2\) on the six-dimensional module realizes the graph automorphism \(g\mapsto g^{-T}\). Indeed, the equality \(u\cdot gv=(g^Tu)\cdot v\) gives \(M\rho(g)=\rho(g^{-T})M\). This automorphism is outer: it interchanges point and plane stabilizers, and these two classes of stabilizers in \(D\) are not conjugate. We thus obtain an extension to \(\Aut(D)\cong X\), whose determinant on an outer involution is \((-8)/(\sqrt2)^6=-1\). The other extension differs by a sign on the outer coset; its determinant is the same, since the dimension \(6\) is even. These are all simple \(X\)-modules with this six-dimensional restriction: after an extension is chosen, the multiplicity space is a simple module for \(X/D\cong C_2\), hence is one-dimensional. In both cases, every outer element has determinant \(-1\).

Thus every simple constituent of \(W\) has dimension six, and an outer element \(y\) has determinant \(-1\) on it. An element of order \(8\) in \(X\) is conjugate to its inverse: on the nonsplit torus \(\mathbb F_{49}^{\times}/\mathbb F_7^{\times}\), inversion is realized by the Frobenius automorphism \(v\mapsto v^7\). Hence the eigenvalues of \(y\) other than \(1,-1\) occur in reciprocal pairs and have product \(1\). Its square has order \(4\) and lies in \(D\); by the preceding lemma on \(D\)-modules, its fixed space has dimension \(2\) on each of the three six-dimensional types. If \(m_+,m_-\) are the multiplicities of \(1,-1\) for \(y\), we obtain
\[
 m_++m_-=2,\qquad (-1)^{m_-}=\det y=-1.
\]
Thus \(m_+=m_-=1\), as required.
\end{proof}

\Needspace{11\baselineskip}
\begin{lemma}\label{lem:L27-outer-overgroup-three-fixed}
Let \(D\trianglelefteq B\), \(D\cong L_2(7)\), and suppose that the conjugation image of \(B\) on \(D\) is \(\operatorname{PGL}_2(7)\). Let \(W\ne0\) be an elementary abelian \(5\)-group acted on by \(B\), with \(\omega(W\rtimes B)\subseteq\Om\), such that the elements of order \(7\) in \(D\) have no nonzero fixed vectors on \(W\). Then every element of order \(3\) in \(B\) has a nonzero fixed space on \(W\).
\end{lemma}
\begin{proof}
Put \(C=C_B(D)\). The group \(C\) has no elements of orders \(5,7,4,9,6\): together with suitable elements of \(D\), these would give \(35,28,28,18,42\), respectively. By the spectrum condition there are no prime divisors other than \(2,3,5,7\), so
\[
 \omega(C)\subseteq\{1,2,3\}.
\]
Choose \(y\in B\) inducing an automorphism \(\beta\) of order \(8\). Then \(|y|=8\). The automorphism \(\beta^2\) is inner of order \(4\), represented by \(a\in D\), \(|a|=4\). Hence \(y^2=ac\) with \(c\in C\). The elements \(a,c\) commute, \(c^4=1\), and the absence of order \(4\) in \(C\) gives \(c^2=1\). Moreover, \(y\) centralizes \(a\): its action commutes with \(\beta^2\), and \(Z(D)=1\). Thus \([y,c]=1\). The subgroup \(B_0=\langle D,y\rangle\) has central subgroup \(\langle c\rangle\), and
\[
 B_0/\langle c\rangle\cong\operatorname{PGL}_2(7).
\]
Indeed, the conjugation image is the full automorphism group, and, using \(y^2=ac\), every element of the kernel in \(B_0\) lies in \(\langle c\rangle\).

Consider \(W_+=C_W(c)\). This space is \(B_0\)-invariant and is a module for the indicated quotient. If \(W_+\ne0\), the preceding lemma gives a nonzero vector of order \(5\) fixed by \(y\), hence order \(40\) in \(W\rtimes B\). Thus \(W_+=0\), and \(c\) acts as \(-I\) on all of \(W\).

If an element of order \(3\) in \(C\) acted nontrivially on \(W\), its image would have order \(3\) and, together with the scalar \(-I\), would give order \(6\) in the image of \(C\). This is impossible, since the spectrum of that image is contained in \(\omega(C)\). Thus all elements of order \(3\) in \(C\) centralize \(W\).

Now suppose that \(|b|=3\) for an element of \(B\). Its action on \(D\) is inner, given by an element \(d\in D\), \(d^3=1\). Since \(b,d\) induce the same action on \(D\), we have \([b,d]=1\) and \(b=de\), where \(e\in C\), \(e^3=1\). Hence the action of \(b\) on \(W\) coincides with that of \(d\). If \(d=1\), every vector is fixed. If \(|d|=3\), apply the calculation on \(7:3\) from the proof of Lemma~\ref{lem:L27-binary-Q-elementary-three-fixed}: the nontrivial \(7\)-weight spaces are permuted in triples, whence \(\dim C_W(d)=\dim W/3>0\). The normalizer \(7:3\) also occurs in \(L_2(7)\), so that calculation does not require the central involution of the double cover.
\end{proof}

\Needspace{11\baselineskip}
\begin{lemma}\label{lem:L27-outer-five-Sylow-normal}
If \(S=L_2(7)\) and \(G/R(G)\cong\operatorname{PGL}_2(7)\), a Sylow \(5\)-subgroup \(T\) of \(K\) centralizes \(H\). Hence \(K=H\times T\), and \(T\) is normal in \(G\). Here \(T\) is not assumed to be elementary abelian.
\end{lemma}
\begin{proof}
Choose a highest chief factor \(V=N/M\) inside \(H\) on which \(T\) acts nontrivially, if such a factor exists. As in the proof of Lemma~\ref{lem:L27-binary-five-Sylow-normal}, coprime lifting on the factors above \(N\) gives \([T,H]\le N\), so
\[
 K/N=(H/N)\times(TN/N).
\]
The subgroup \(TN/N\) is normal in \(G/N\) and projects isomorphically and \(G\)-equivariantly onto \(Q=K/H\). For \(Y=G/C_G(V)\), its image \(A\) is a nontrivial normal \(5\)-group. The module \(V\) is faithful and simple for \(Y\), has characteristic \(p\in\{2,3\}\), and satisfies \(C_V(A)=0\).

Let \(\Sigma\) be the set of simple types in the restriction of \(V\) to \(A\) over the algebraic closure. It consists of nonprincipal types and is permuted by \(Y\). Choose an element \(s\in L\) of order \(7\), whose image in \(L/H\) lies in \(D\) and is nontrivial. It acts fixed-point-freely on \(Q\), since otherwise \(35\) occurs. The same holds on its quotient \(A\), by Lemma~\ref{lem:coprime-fixed-general}. Lemma~\ref{lem:fixedfree-characters}, in its general formulation for a \(5\)-group and an automorphism of order \(7\), shows that \(s\) fixes no type in \(\Sigma\). Thus the image of the perfect group \(L\) on \(\Sigma\) is nontrivial.

Now take \(g\in G\) mapping to an element of order \(8\) in \(G/R(G)\cong\operatorname{PGL}_2(7)\). Then \(|g|=8\), and the image of \(g^4\) is a nontrivial involution of the socle \(L_2(7)\). We show that \(g^4\) cannot fix all of \(\Sigma\). The group \(L\) is nontrivial on \(V\), since its subgroup \(T\le K\) is nontrivial on \(V\). By Lemma~\ref{lem:mixed-chief-action}, \(Y/R(Y)\cong G/R(G)\). If the kernel of the permutation action on \(\Sigma\) contains \(g^4\), its image in \(Y/R(Y)\) contains the entire simple socle. Hence the image of the perfect group \(L\) in the permutation group is solvable: it lies in the image of \(R(Y)\). A perfect solvable image is trivial, contrary to the action of \(s\).

Thus some type is moved by \(g^4\). Its orbit under \(\langle g\rangle\) has length \(8\). On the corresponding eight isotypic components, the operator \(1+g+\cdots+g^7\) is nonzero. The norm calculation in the actual subgroup \(\langle V,gM\rangle\le G/M\), as in Lemma~\ref{lem:L27-binary-five-Sylow-normal}, gives order \(16\) for \(p=2\) or \(24\) for \(p=3\). The contradiction excludes the chosen factor. Therefore \(T\) centralizes all factors in a chief series inside \(H\), and coprime lifting gives \([T,H]=1\). The subgroup \(T\) is the unique Sylow \(5\)-subgroup of \(K=H\times T\), so it is characteristic in \(K\) and normal in \(G\).
\end{proof}

\Needspace{11\baselineskip}
\begin{theorem}\label{thm:L27-outer-extension-excluded}
If \(\omega(G)=\Om\) and \(\Soc(G/R(G))\cong L_2(7)\), then
\[
 G/R(G)\cong L_2(7).
\]
The outer extension \(G/R(G)\cong\operatorname{PGL}_2(7)\) is impossible.
\end{theorem}
\begin{proof}
Suppose that the quotient is the outer extension. By the preceding lemma, \(K=H\times T\), where \(T\) is a normal Sylow \(5\)-subgroup of \(G\). Put \(W=Z(T)\ne1\). This is an elementary abelian \(5\)-group, since \(25\notin\Om\). The projection \(T\to Q\) is a \(G\)-equivariant isomorphism taking \(W\) to \(Z(Q)\). The action of \(G\) on \(W\) factors through \(G/H=Q\rtimes B\), and \(Q\) is trivial on its center. We have \(W\rtimes B\le G/H\), and \(D=B^\infty\cong L_2(7)\) is normal in \(B\).

The conjugation image of \(B\) on \(D\) is \(\operatorname{PGL}_2(7)\). Indeed, \(K\le R(G)\) and \(B\cong G/K\), so \(B/R(B)\cong G/R(G)\). The image of \(D\) is the socle of this quotient, whose centralizer in an almost simple group is trivial. Hence \(C_B(D)\le R(B)\). Conversely, \([R(B),D]\le R(B)\cap D=1\). Thus the conjugation kernel equals \(R(B)\), and the image is isomorphic to the stated outer extension. The elements of order \(7\) in \(D\) have no nonzero fixed vectors on \(W\), since otherwise \(35\) occurs. Lemma~\ref{lem:L27-outer-overgroup-three-fixed} applies to \(W,B,D\).

Choose \(x\in G\) of order \(9\). Its image in \(G/H\) has order \(1\) or \(3\). If the image is trivial, \(x\in H\) centralizes \(T\), giving \(45\). Otherwise, conjugation moves the image into \(B\), since \(B\) contains a Sylow \(3\)-subgroup of \(G/H\). By the cited lemma, this image fixes a nonzero vector in \(W\). Since \(H\) centralizes \(T\), the same vector is fixed by \(x\) itself. We obtain commuting elements of orders \(9,5\), hence the forbidden order \(45\). This completes the exclusion of the outer extension.
\end{proof}

\Needspace{11\baselineskip}
\begin{lemma}\label{lem:L27-final-ordinary-radical-structure}
In the remaining branch with socle \(L_2(7)\),
\[
 P=K,\qquad G/K\cong L_2(7)\times(R(G)/K),
 \qquad \omega(R(G)/K)\subseteq\{1,2,3\}.
\]
The normal Hall subgroup \(H\) has even order and contains an element of order \(6\). In particular, \(\pi(H)=\{2,3\}\), and a Sylow \(2\)-subgroup of \(H\) is neither cyclic nor generalized quaternion.
\end{lemma}
\begin{proof}
By Theorem~\ref{thm:L27-binary-cover-excluded}, \(L/K\cong L_2(7)\), so its center \(P/K\) is trivial: \(P=K\). By the preceding theorem, \(G/R(G)\cong L_2(7)\), whence \(G=LR(G)\). The two normal subgroups \(L/K\), \(R(G)/K\) intersect trivially and commute, since \([L,R(G)]\le L\cap R(G)=K\). This gives the stated direct product. Its second factor has no prime divisors \(5,7\), nor elements of orders \(4,9,6\): together with suitable elements of the first factor, these would give \(35,28,28,18,42\). Thus its spectrum lies in \(\{1,2,3\}\).

The group \(G/K\) has no element of order \(8\), since Sylow \(2\)-subgroups of the first factor have exponent \(4\), and those of the second have exponent at most \(2\). An element of order \(8\) in \(G\) therefore generates a cyclic subgroup intersecting \(K\) nontrivially, so \(2\mid|K|\). The quotient \(K/H\) has odd order, hence \(2\mid|H|\). We already know \(3\mid|H|\), \(H=[H,L]\), and noncyclicity of a Sylow \(3\)-subgroup. Lemma~\ref{lem:commutator-generated-Sylow-ranks} excludes cyclic and generalized quaternion Sylow \(2\)-subgroups of \(H\).

If \(H\) had no element of order \(6\), its nontrivial Fitting subgroup \(F(H)\) would be a \(2\)-group or a \(3\)-group. In the first case, a subgroup \(C_3^2\) of a Sylow \(3\)-subgroup acts on the nonzero group \(\Omega_1(Z(F(H)))\), and some \(3\)-element fixes a vector of order \(2\). In the second case, a subgroup \(C_2^2\) of a Sylow \(2\)-subgroup acts on the corresponding space of characteristic \(3\), and some involution fixes a nonzero \(3\)-vector. Both cases give order \(6\). The contradiction proves \(6\in\omega(H)\).
\end{proof}

\clearpage
\subsection{Excluding nilpotent \(K\) for \(L_2(7)\)}

We retain the assumptions of the remaining branch \(S=L_2(7)\). We have already proved
\[
 G/R(G)\cong L/K\cong S,\qquad P=K,
 \qquad G/K\cong S\times T,
 \quad T=R(G)/K,\quad\omega(T)\subseteq\{1,2,3\}.
\]
In this subsection we derive a contradiction from the assumption that \(K\) is nilpotent.

\Needspace{11\baselineskip}
\begin{lemma}\label{lem:L27-char-three-four-fixed}
Every element of order \(4\) in \(S=L_2(7)\) has a nonzero fixed vector on every nonzero \(\mathbb F_3S\)-module. For the absolutely simple types, the complete list of dimensions and fixed-space dimensions is
\[
\begin{array}{c|ccccc}
 \dim&1&3&3&6&7\\\hline
 \dim C(4)&1&1&1&2&1.
\end{array}
\]
\end{lemma}
\begin{proof}
Extending the field does not change fixed-space dimensions; we work over the algebraic closure of a field of characteristic \(3\). The full list of simple types in characteristic \(3\) is \(1,3,3,6,7\), by the MOC decomposition matrix and the ATLAS catalogue \cite{L27BinaryData}. The six-dimensional type is realized as the zero-coordinate-sum submodule of the permutation module on seven points, and the seven-dimensional type as the analogous submodule for the action on the eight points of \(\mathbb P^1(\mathbb F_7)\). The cited decomposition matrix confirms simplicity of both submodules.

An element of order \(4\) is conjugate to its inverse, as explained in Lemma~\ref{lem:L27-seven-free-five-modules-four}. On a three-dimensional type in characteristic \(3\), this gives equal multiplicities of the eigenvalues \(i,-i\), \(i^2=-1\). The determinant is \(1\), so the multiplicity of \(-1\) is even, and that of \(1\) is odd. The latter cannot equal \(3\), since a nontrivial simple module for a simple group is faithful. It therefore equals \(1\).

On seven points, an element of order \(4\) has cycles \(4,2,1\), so its fixed-space dimension on the six-dimensional submodule is \(3-1=2\). On eight points, its cycles have lengths \(4,4\), giving \(2-1=1\) for the seven-dimensional type. Here the constant constituent splits off, since neither \(7\) nor \(8\) is divisible by \(3\). The assertion is obvious for the trivial type.

Finally, taking fixed points of an element of order \(4\) is exact in characteristic \(3\), by Maschke's theorem for \(C_4\). Applying this along a composition series of an arbitrary module gives a positive sum of dimensions from the table. Complete reducibility of the \(S\)-module itself is not required.
\end{proof}

\Needspace{11\baselineskip}
\begin{lemma}\label{lem:L27-nilpotent-K-three-abelian}
If \(K\) is nilpotent, its Sylow \(3\)-subgroup \(K_3\) is abelian. Moreover,
\[
 K=K_2\times K_3\times K_5,
 \qquad \exp K_2=2,\quad \exp K_3=3,
 \quad\exp K_5=5,
\]
and all three factors are nontrivial.
\end{lemma}
\begin{proof}
Lemma~\ref{lem:nilpotent-K-prime-shapes} gives the stated decomposition, \(K_3,K_5\ne1\), \(\exp K_3=3\), \(\exp K_5=5\). Lemma~\ref{lem:L27-final-ordinary-radical-structure} gives \(2\mid|K|\), so \(K_2\ne1\). The group \(K_2\) has no element of order \(4\): together with commuting elements of orders \(3,5\) from the other factors, this would give \(60\). Thus \(\exp K_2=2\).

Choose \(g\in G\) of order \(8\). Its image in \(G/K=S\times T\) has order \(4\). Indeed, the order of the image divides \(4\); if it were at most \(2\), the element \(g^2\) of order \(4\) would lie in \(K\). Since Sylow \(2\)-subgroups of \(T\) have exponent at most \(2\), the image of \(g\) has the form \((a,c)\), where \(|a|=4\) in \(S\) and \(c^2=1\) in \(T\).

Put \(N=K_3\). We have \(C_N(g)=1\), since otherwise an element of order \(3\) in this centralizer, together with \(g\), would give \(24\). Coprime lifting of fixed points by Lemma~\ref{lem:coprime-fixed-general} shows that \(g\) has no nonzero fixed vectors on \(V=N/\Phi(N)\). Since \(\exp N=3\), we have \(\Phi(N)=N'\). The action of \(K\) on \(V\) is trivial: \(N\) acts by inner automorphisms, and the other Sylow subgroups of \(K\) commute with \(N\). Hence \(V\) is an \(S\times T\)-module.

Decompose \(V\) as the direct sum of the \(+1\)- and \(-1\)-eigenspaces of the involution \(c\). These are \(S\)-invariant, since the factors \(S,T\) commute. On the \(+1\)-eigenspace, the action of \(g\) agrees with that of \(a\). If this space is nonzero, the preceding lemma gives a nonzero fixed vector for \(g\), a contradiction. Thus \(c=-I\) on all of \(V\).

If \(N\) is nonabelian, then \(W=N'/\gamma_3(N)\ne0\), where \(\gamma_3(N)=[N',N]\). Indeed, equality \(N'=\gamma_3(N)\) would make the lower central series stabilize at a nonzero term, impossible for a finite \(3\)-group. The commutator gives a surjective equivariant map
\[
 \bigwedge\nolimits^2 V\longrightarrow W,
 \qquad xN'\wedge yN'\longmapsto[x,y]\gamma_3(N).
\]
It is bilinear because commutators are central modulo \(\gamma_3(N)\). The subgroup \(K\) is also trivial on \(W\), so this is likewise a map of \(S\times T\)-modules. Since \(c=-I\) on \(V\), we have \(c=I\) on \(W\); thus \(g\) acts on \(W\) as \(a\). The preceding lemma gives \(C_W(g)\ne0\). But coprime lifting from \(N'/\gamma_3(N)\) gives
\[
 C_W(g)=C_{N'}(g)\gamma_3(N)/\gamma_3(N)=0,
\]
since \(C_N(g)=1\). This contradiction proves \(N'=1\).
\end{proof}

\Needspace{11\baselineskip}
\begin{lemma}\label{lem:elementary-three-kernel-nine-eighteen}
Let \(V\trianglelefteq F\) be an elementary abelian \(3\)-group, with
\[
 F/V\cong C_3\times C_2^2.
\]
If \(18\notin\omega(F)\), then \(9\notin\omega(F)\). No splitting of \(F\) over \(V\) is assumed.
\end{lemma}
\begin{proof}
Denote the full preimage of the factor \(C_3\) by \(E\). This is a normal \(3\)-subgroup of \(F\). In the full preimage of the second factor \(C_2^2\), the subgroup \(V\) is normal and Hall, so the Schur--Zassenhaus theorem gives a complement \(A\cong C_2^2\). The group \(A\) normalizes \(E\) and acts trivially on \(E/V\cong C_3\). By Lemma~\ref{lem:coprime-fixed-general}, the centralizer \(C_E(A)\) maps onto all of \(E/V\). Choose \(z\in C_E(A)\) mapping to a generator of this quotient. Since \(z^3\in V\), we have \(|z|=3\) or \(9\). The latter is impossible: \(z\) commutes with an involution of \(A\) and would give \(18\). Thus \(|z|=3\), and \(E=V\rtimes\langle z\rangle\).

Writing \(V\) additively, put \(N_z=1+z+z^2\). If \(E\) contains an element of order \(9\), then \(N_z\ne0\) on \(V\): the cubes of elements \(vz\) and \(vz^2\) equal \(N_z(v)\), whereas the elements of \(V\) have order at most \(3\). The group \(A\) centralizes \(z\), so \(N_z\) commutes with its action. Over \(\mathbb F_3\), the module \(V\) decomposes into weight spaces for linear characters \(A\to\{1,-1\}\). If \(N_z\ne0\), it is nonzero on at least one such space. The kernel of the corresponding character contains an involution \(t\in A\). Choose \(v\) in this space with \(N_z(v)\ne0\). Then \(vz\) has order \(9\) and is centralized by \(t\), giving \(18\), a contradiction. Thus \(E\) contains no element of order \(9\). Every element of order \(9\) in \(F\) would lie in \(E\), since its image in \(F/V\) has odd order. The lemma follows.
\end{proof}

\Needspace{11\baselineskip}
\begin{theorem}\label{thm:L27-nilpotent-K-excluded}
If \(\omega(G)=\Om\) and \(\Soc(G/R(G))\cong L_2(7)\), then \(K=[R(G)\cap G^\infty,G^\infty]\) is not nilpotent.
\end{theorem}
\begin{proof}
Suppose that \(K\) is nilpotent. By Lemma~\ref{lem:L27-nilpotent-K-three-abelian}, \(K_3\) is abelian and \(\exp K_3=3\). The subgroups \(K_2,K_3,K_5\) are characteristic in \(K\) and normal in \(G\). Put \(U=Z(K_5)\ne1\). This is an elementary abelian \(5\)-group, since \(\exp K_5=5\). The whole of \(K\) centralizes \(U\), so the action of \(G\) on \(U\) factors through \(G/K=S\times T\), where \(T=R(G)/K\). Elements of order \(7\) in \(S\) have no nonzero fixed vectors on \(U\): otherwise, lifting such an element to a \(7\)-element of \(G\) would give \(35\).

Choose \(x\in G\) of order \(9\), and write its image in \(G/K\) as \((a,b)\in S\times T\). Both components have order at most \(3\). Suppose that \(|a|=3\). In \(S\), there is a subgroup of order \(7\) normalized by \(a\) with nontrivial action of order \(3\), as in the calculation for \(7:3\) in Lemma~\ref{lem:L27-binary-Q-elementary-three-fixed}. The component \(b\) commutes with \(S\) and preserves all its \(7\)-weight spaces on \(U\). Thus \((a,b)\) cyclically permutes all these nontrivial weight spaces in triples. Since \((a,b)^3=1\), each such triple gives a nonzero fixed space. Hence \(C_U(x)\ne1\), giving \(45\), a contradiction. Thus \(a=1\), that is, \(x\in R(G)\). Moreover, \(x\notin K\), since \(\exp K_3=3\), so \(|b|=3\).

Pass to \(G_0=G/(K_2K_5)\). This is a quotient by a normal \(3'\)-subgroup, so the image \(x_0\) of \(x\) still has order \(9\). The normal subgroup \(V=K/(K_2K_5)\cong K_3\) is elementary abelian of characteristic \(3\), and \(G_0/V\cong S\times T\). Choose \(A_0\cong C_2^2\) inside \(S\), and denote the full preimage of \(A_0\times\langle b\rangle\) in \(G_0\) by \(F\). Then
\[
 F/V\cong C_2^2\times C_3,
 \qquad x_0\in F,\quad |x_0|=9.
\]
But \(18\notin\omega(F)\), since \(F\) is a section of \(G\). Lemma~\ref{lem:elementary-three-kernel-nine-eighteen} gives \(9\notin\omega(F)\), a contradiction. Nilpotent \(K\) is completely excluded.
\end{proof}

\clearpage
\subsection{The nonnilpotent branch for \(L_2(7)\): normal \(5\)-subgroups}

We retain the notation and structural conclusions established above, in particular \(G/K=S\times T\), where \(S\cong L_2(7)\), \(T=R(G)/K\), and \(\omega(T)\subseteq\{1,2,3\}\). Recall that \(K=H\rtimes U\), where \(H\) is the characteristic normal Hall \(5'\)-subgroup, and the Sylow \(5\)-subgroup \(U\) contains \(C_5^6\). The groups considered in this branch remain hypothetical; no existence assertion is intended.

\Needspace{11\baselineskip}
\begin{theorem}\label{thm:L27-no-normal-five}
In the remaining \(L_2(7)\) branch,
\[
 O_5(K)=O_5(R(G))=O_5(G)=1.
\]
\end{theorem}
\begin{proof}
Suppose that \(N=O_5(K)\ne1\). The group \(N\) is characteristic in \(K\) and normal in \(G\). Choose a chief series of \(G\) inside \(N\). By Lemma~\ref{lem:L27-L34-central-five-seven-chief}, \(K\) centralizes each factor. Hence \(H\) centralizes \(N\): this follows by induction along the series from coprime lifting of fixed points (Lemma~\ref{lem:coprime-fixed-general}), since \(H\) is a \(5'\)-group and all factors in the series are \(5\)-groups.

Put \(V=\Omega_1(Z(N))\ne0\). This is an elementary abelian \(5\)-space normal in \(G\) and centralized by \(H\). The \(5\)-group \(U\) acting on it has a nonzero fixed vector: in the orbit partition of \(V\), the number of fixed vectors is divisible by \(5\) and includes the zero vector. Thus
\[
 A=C_V(K)=V\cap Z(K)\ne0.
\]
The subgroup \(A\) is characteristic in \(K\), and the action of \(G\) on it factors through \(S\times T\). Every element of order \(7\) in \(S\) acts on \(A\) without nonzero fixed vectors: otherwise a lift of order \(7\) in \(G\), together with a vector in \(A\), would give order \(35\). Consequently, the action of \(S\) on every nonzero \(S\)-submodule of \(A\) is faithful: its kernel is normal in the simple group \(S\), and the whole of \(S\) cannot act trivially.

Choose an element \(x\in G\) of order \(9\), and write its image in \(S\times T\) as \((a,b)\). Both components have order at most \(3\). If \(|a|=3\), then \(a\) normalizes a subgroup of order \(7\) and permutes its nontrivial weight spaces on \(A\) in triples. The component \(b\) commutes with \(S\), so it preserves individual \(7\)-weights. Since \((a,b)^3=1\), the sum of a vector and its two images under \((a,b)\) gives a nonzero fixed vector on each triple of weight spaces. This vector and \(x\) give order \(45\), a contradiction. Thus \(a=1\). Moreover, \(C_A(x)=0\), since otherwise order \(45\) again occurs. Hence \(b\in T\) has order \(3\), and its image in \(\operatorname{GL}(A)\) has order \(3\).

Now take \(y\in G\) of order \(8\), and let its image be \((d,c)\in S\times T\). The order of \(d\) is \(1,2\), or \(4\), and that of \(c\) is at most \(2\). Since \(40\notin\Om\), we have \(C_A(y)=0\). Decompose \(A=A_+\oplus A_-\) into the eigenspaces of \(c\) with eigenvalues \(1,-1\); both are \(S\)-submodules.

If \(d=1\), absence of fixed vectors gives \(A_+=0\). If \(|d|=4\) and \(A_+\ne0\), the element \(d\) has a nonzero fixed vector there by Lemma~\ref{lem:L27-seven-free-five-modules-four}; this vector is also fixed by \(y\). Thus again \(A_+=0\). If \(|d|=2\), then on every nonzero \(S\)-submodule, the action of \(d\) has both eigenvalues \(1,-1\). Indeed, the action of \(S\) there is faithful, so \(d\) cannot act as a scalar, since it would then lie in the center of \(S\). Since the characteristic is \(5\), the involution is diagonalizable. If \(A_+\ne0\), then \(y\) fixes a vector in the \(+1\)-eigenspace of \(d\); if \(A_-\ne0\), it fixes a vector in the \(-1\)-eigenspace of \(d\). Both are impossible, so \(|d|=2\) is excluded.

Thus \(c\) acts as \(-I\) on all of \(A\). Its image has order \(2\) and commutes with the image of \(b\), of order \(3\). Their product has order \(6\) in the image of \(T\) in \(\operatorname{GL}(A)\). But the spectrum of a quotient is contained in that of the original group, and \(6\notin\omega(T)\). This contradiction proves \(O_5(K)=1\).

All elements of order \(5\) in \(G\) lie in \(K\), and \(G\) has no element of order \(25\). Thus every normal \(5\)-subgroup of \(G\) lies in \(K\) and is trivial. Finally, \(O_5(R(G))\) is characteristic in \(R(G)\), hence normal in \(G\), and is also \(1\).
\end{proof}

\Needspace{11\baselineskip}
\begin{lemma}\label{lem:L27-five-paired-large}
The group \(K\) contains subgroups
\[
 C_3^2\times C_5^2\qquad\text{and}\qquad C_2^2\times C_5^2.
\]
\end{lemma}
\begin{proof}
Let \(p=3\) or \(p=2\). Choose a Hall \(\{p,5\}\)-subgroup \(J\) of \(K\), its normal Sylow \(p\)-subgroup \(P=J\cap H\), and a subgroup \(A\cong C_5^6\) of a Sylow \(5\)-subgroup of \(J\). For \(p=3\), the group \(P\) is noncyclic; for \(p=2\), it is neither cyclic nor generalized quaternion, by Lemma~\ref{lem:L27-final-ordinary-radical-structure}.

If \(A\) acts trivially on \(P/\Phi(P)\), it centralizes \(P\): the kernel of the action of automorphisms of \(P\) on its Frattini quotient is a \(p\)-group. In this case \(P\) contains \(C_p^2\), giving the required product.

Otherwise, \(P/\Phi(P)\) has a nontrivial irreducible \(\mathbb F_p A\)-module of dimension \(4\), since \(\operatorname{ord}_5(2)=\operatorname{ord}_5(3)=4\). Its kernel \(A_0\) has index \(5\) in \(A\). Coprime lifting of fixed points shows that \(C_P(A_0)\) has an elementary abelian quotient of rank at least \(4\). For \(p=3\), this group is noncyclic and therefore contains \(C_3^2\); for \(p=2\), it is neither cyclic nor generalized quaternion and therefore contains \(C_2^2\). In both cases, \(A_0\cong C_5^5\) centralizes the subgroup found.
\end{proof}

\Needspace{11\baselineskip}
\begin{lemma}\label{lem:L27-normal-abelian-exp}
Every normal abelian subgroup \(B\triangleleft K\) satisfies \(\exp B\mid30\). In particular, normal abelian \(2\)-subgroups are elementary abelian, and normal abelian \(3\)-subgroups have exponent at most \(3\).
\end{lemma}
\begin{proof}
The Sylow components of \(B\) are normal in \(K\). Suppose that \(B_2\) contains an element of order \(4\). Then \(M=\Omega_2(B_2)/\Omega_1(B_2)\ne0\) is an elementary abelian \(2\)-space acted on by the subgroup \(C_3^2\times C_5^2\) from the preceding lemma. Extend the field and choose a weight of this abelian group of operators. The kernel of its restriction to each factor is nontrivial. Thus there are commuting elements \(u,v\) of orders \(3,5\) having a common nonzero fixed vector on \(M\). The fixed space is defined over \(\mathbb F_2\). By coprime lifting of fixed points, this class lifts to an element of order \(4\) in \(\Omega_2(B_2)\) centralizing \(u,v\). Their product has forbidden order \(60\).

If \(B_3\) contains an element of order \(9\), some \(C_5^2\le K\) acts on the nonzero group \(\Omega_2(B_3)/\Omega_1(B_3)\). Each of its weights has nontrivial kernel, so some element of order \(5\) fixes a nonzero class. Its fixed lift has order \(9\), giving the forbidden order \(45\). The \(5\)-component has exponent at most \(5\), since \(25\notin\Om\); no other primes divide \(|K|\). Thus \(\exp B\mid30\).
\end{proof}

\Needspace{11\baselineskip}
\begin{theorem}\label{thm:L27-fitting-large}
In the branch under consideration,
\[
 F(K)=F(H)=O_2(K)\times O_3(K),\qquad
 \exp Z(F(K))\mid6.
\]
If \(d_p=\dim_{\mathbb F_p}O_p(K)/\Phi(O_p(K))\) for \(p=2,3\), then
\[
 d_2+d_3\ge24,\qquad |K|\ge 2^{24}5^6.
\]
Moreover, \(30\in\omega(K)\).
\end{theorem}
\begin{proof}
By Theorem~\ref{thm:L27-no-normal-five}, the subgroup \(F(K)\) has no nontrivial Sylow \(5\)-subgroup and therefore lies in \(H\). Conversely, \(F(H)\) is characteristic in \(H\), normal in \(K\), and nilpotent. Hence \(F(K)=F(H)=O_2(K)\times O_3(K)\). The center of \(F(K)\) is a normal abelian \(\{2,3\}\)-subgroup of \(K\); the preceding lemma bounds its exponent.

For a finite solvable group, \(C_K(F(K))\le F(K)\). Thus a Sylow \(5\)-subgroup \(U\) acts faithfully on \(F(K)\). Its action on the two Frattini quotients together is also faithful: the kernel of the action on the Frattini quotient of a finite \(p\)-group is a \(p\)-group, whereas elements of \(U\) have \(5\)-power order.

Restrict the action to \(A\cong C_5^6\le U\). In characteristics \(2\) and \(3\), every nontrivial irreducible \(A\)-module has dimension \(4\), and its kernel is a hyperplane in the six-dimensional \(\mathbb F_5\)-space \(A\). A faithful combined action requires at least six such constituents: the intersection of fewer than six hyperplanes is nontrivial. Thus \(d_2+d_3\ge24\). Since \(|U|\ge5^6\) and \(|F(K)|\ge2^{d_2}3^{d_3}\ge2^{24}\), the stated bound for \(|K|\) follows.

It remains to prove the occurrence of order \(30\) inside \(K\). If \(O_2(K)\ne1\), the group \(C_3^2\times C_5^2\) acts on \(\Omega_1(Z(O_2(K)))\ne0\). Each weight has nontrivial kernels on both factors, so some commuting elements of orders \(3,5\) fix a nonzero \(2\)-vector and give order \(30\). If \(O_2(K)=1\), then \(O_2(H)=1\), and the nontrivial group \(F(H)\) is a \(3\)-group. Now \(C_2^2\times C_5^2\) acts on \(\Omega_1(Z(F(H)))\ne0\); the same argument gives commuting elements of orders \(2,3,5\). Thus \(30\in\omega(K)\) in both cases.
\end{proof}

\clearpage
\subsection{Elements of order \(9\) and Sylow \(3\)-subgroups for \(L_2(7)\)}

This subsection continues the analysis of the hypothetical branch \(S=L_2(7)\) described above. We use the following classical fact: every group of exponent dividing \(3\) is \(2\)-Engel and nilpotent of class at most \(3\). Precise statements of both parts can be found, for example, in \cite[Fact~2.15]{BurnsideThree}; for the original result of Levi and van der Waerden, see also \cite{LeviWaerden}.

\Needspace{11\baselineskip}
\begin{lemma}\label{lem:five-rank-exponent-three}
Let an elementary abelian group \(A\cong C_5^r\), \(r\ge5\), act on a finite \(3\)-group \(P\). If \(C_P(a)\) has exponent dividing \(3\) for every \(1\ne a\in A\), then \(P\) also has exponent dividing \(3\). Moreover, \(P\) is \(2\)-Engel and has nilpotency class at most \(3\).
\end{lemma}
\begin{proof}
Consider the set
\[
 X=\bigcup_{\substack{B<A\\ {}[A:B]=5}} C_P(B).
\]
It generates \(P\). Indeed, every irreducible constituent of the \(\mathbb F_3 A\)-module \(P/\Phi(P)\) is fixed by some hyperplane \(B\le A\): this holds after extending the field for every linear character, and therefore also over the original field. Coprime lifting of fixed points gives a preimage of each such fixed subspace in \(C_P(B)\). Thus the images of the groups \(C_P(B)\) generate \(P/\Phi(P)\), and the groups themselves generate \(P\).

Any four elements \(x_1,\dots,x_4\in X\) are simultaneously centralized by some \(1\ne a\in A\): the intersection of the four corresponding hyperplanes in \(A\) is nontrivial. By assumption, \(\langle x_1,\dots,x_4\rangle\) has exponent dividing \(3\). In particular, \([x_1,x_2,x_3,x_4]=1\). The quotient \(\gamma_4(P)/\gamma_5(P)\) is generated by the images of these commutators, since \(X\) generates \(P\). Hence \(\gamma_4(P)=\gamma_5(P)\), and nilpotence of the finite \(3\)-group gives \(\gamma_4(P)=1\).

In a group of class at most \(3\), the triple commutator \([x,y,z]\) is central and multilinear in the classes of \(x,y,z\) in the abelianization. Any three elements of \(X\) lie in a subgroup of exponent dividing \(3\), so for \(x_i,x_j,x_k\in X\),
\[
 [x_i,x_j,x_j]=1,\qquad
 [x_i,x_j,x_k]\,[x_i,x_k,x_j]=1.
\]
The second equality is obtained by expanding \([x_i,x_jx_k,x_jx_k]=1\). Multilinearity shows that these equalities imply \([x,y,y]=1\) for all \(x,y\in P\). Thus \(P\) is \(2\)-Engel.

Moreover, \(P'=\gamma_2(P)\) is abelian, since \([\gamma_2(P),\gamma_2(P)]\le\gamma_4(P)=1\). It is generated by conjugates of commutators of pairs of elements of \(X\). Each such commutator has order dividing \(3\), since its two defining elements lie in a group of exponent dividing \(3\). Hence \(\exp P'\mid3\). The standard expansion of the cube of a product in a group of class at most \(3\), using \([x,y,y]=[y,x,x]=1\), gives
\[
 (xy)^3=x^3y^3[y,x]^3=x^3y^3.
\]
Thus the map \(x\mapsto x^3\) is a homomorphism. Every element of \(X\) has cube \(1\), so the homomorphism is trivial and \(\exp P\mid3\).
\end{proof}

\Needspace{11\baselineskip}
\begin{theorem}\label{thm:L27-K-no-nine}
In the branch under consideration, a Sylow \(3\)-subgroup of \(K\) has exponent \(3\), is \(2\)-Engel, and has nilpotency class at most \(3\). In particular,
\[
 9\notin\omega(K),\qquad \exp O_3(K)\mid3.
\]
Every element \(x\in G\) of order \(9\) lies outside \(K\), and its image in \(G/K\cong S\times T\) has order \(3\).
\end{theorem}
\begin{proof}
Let \(J\) be a Hall \(\{3,5\}\)-subgroup of \(K\), let \(P=J\cap H\) be its normal Sylow \(3\)-subgroup, and let \(A\cong C_5^6\) lie in a Sylow \(5\)-subgroup of \(J\). For each \(1\ne a\in A\), the group \(C_P(a)\) has no element of order \(9\), since together with \(a\) this would give the forbidden order \(45\). There is also no element of order \(27\) in \(G\), so \(\exp C_P(a)\mid3\). The preceding lemma applies and gives all the assertions about \(P\). By Sylow's theorem, these hold for every Sylow \(3\)-subgroup of \(K\), whence \(9\notin\omega(K)\) and \(\exp O_3(K)\mid3\).

For \(|x|=9\), its image in \(G/K\) has order at most \(3\), since \(\omega(T)\subseteq\{1,2,3\}\) and \(S\) has no element of order \(9\). The image is nontrivial by what was just proved.
\end{proof}

\Needspace{11\baselineskip}
\begin{lemma}\label{lem:elementary-three-five-nine}
Let \(E\) be a finite group with a normal elementary abelian \(3\)-subgroup \(V\), such that
\[
 E/V\cong C_3\times C_5^2.
\]
If \(45\notin\omega(E)\), then \(9\notin\omega(E)\).
\end{lemma}
\begin{proof}
The full preimage of the \(5\)-factor in \(E\) has a complement \(A\cong C_5^2\) to \(V\). The action of \(A\) on \(E/V\) is trivial, so coprime lifting of fixed points gives \(C_E(A)V=E\). Lift a generator of \(C_3\) and take its \(3\)-part to obtain \(z\in C_E(A)\) mapping to a generator of \(C_3\). Its order cannot be \(9\), since together with a nonidentity element of \(A\) this would give order \(45\). Hence \(|z|=3\).

If \(E\) contains an element of order \(9\), it lies in \(V\rtimes\langle z\rangle\). Replacing it by its square if necessary, write it as \(vz\). Writing \(V\) additively, its cube is \((1+z+z^2)v\ne0\). The operator \(1+z+z^2\) commutes with \(A\). Extend the field and decompose \(V\) into \(A\)-weight spaces. The norm acts nontrivially on some such space, and the kernel of its character of \(C_5^2\) contains an element \(1\ne a\). Hence the norm is also nontrivial on \(C_V(a)\), defined over the original field. For a suitable \(w\in C_V(a)\), the element \(wz\) has order \(9\) and centralizes \(a\), giving \(45\), a contradiction.
\end{proof}

\Needspace{11\baselineskip}
\begin{theorem}\label{thm:L27-nine-simple-part-nonabelian}
If some element of order \(9\) has nontrivial projection onto \(S\) under \(G/K=S\times T\), then a Sylow \(3\)-subgroup of \(H\) is nonabelian. In particular, if it is abelian, all elements of order \(9\) lie in \(R(G)\), and \(3\mid|T|\).
\end{theorem}
\begin{proof}
Let \(|x|=9\), and let its image in \(S\times T\) be \((a,b)\) with \(|a|=3\). Put \(Q=K/H\). In the splitting \(G/H=Q\rtimes(S\times T)\), the element \(xH\) acts on \(Z(Q)\) as \((a,b)\) does, since inner automorphisms of \(Q\) are trivial on its center. This center is an elementary abelian \(5\)-group of dimension \(6n\), \(n\ge1\).

The element \(a\) normalizes a subgroup of order \(7\) acting on \(Z(Q)\) without nonzero fixed vectors, since otherwise \(G/H\) would have the forbidden order \(35\). On the six nontrivial \(7\)-weights, \(a\) has orbits of length \(3\), whereas \(b\), commuting with \(S\), preserves individual weights. Since \((a,b)^3=1\), on every triple of weight spaces the fixed-space dimension is one-third of the total dimension. Hence
\[
 \dim_{\mathbb F_5} C_{Z(Q)}(xH)=2n\ge2.
\]
Choose \(A\cong C_5^2\) in this centralizer, and let \(F\) be the full preimage of \(\langle xH\rangle\times A\) in \(G\). The group \(F\) is solvable, since \(H\) is solvable and \(F/H\cong C_3\times C_5^2\). By Hall's theorem, there is a Hall \(\{3,5\}\)-subgroup \(J\le F\) containing \(x\). Its subgroup \(P=J\cap H\) is normal in \(J\), is a Sylow \(3\)-subgroup of \(H\), and satisfies
\[
 J/P\cong C_3\times C_5^2.
\]
If \(P\) is abelian, Theorem~\ref{thm:L27-K-no-nine} makes it elementary abelian. The preceding lemma then forbids order \(9\) in \(J\), a contradiction. Thus \(P\) is nonabelian.

If a Sylow \(3\)-subgroup of \(H\) is abelian, what we have proved shows that \(a=1\) for every \(|x|=9\). Then \(x\in R(G)\), and its image in \(T\) has order \(3\), since \(x\notin K\).
\end{proof}

\Needspace{11\baselineskip}
\begin{corollary}\label{cor:L27-abelian-three-not-three-nilpotent}
If a Sylow \(3\)-subgroup of \(K\) is abelian, \(K\) has no normal Hall \(3'\)-subgroup.
\end{corollary}
\begin{proof}
Suppose otherwise, and denote the normal Hall \(3'\)-subgroup by \(W\). It is unique, hence characteristic in \(K\), and therefore normal in \(G\). The group \(V=K/W\) is isomorphic to a Sylow \(3\)-subgroup of \(K\) and is elementary abelian by Theorem~\ref{thm:L27-K-no-nine}.

Choose \(x\in G\) of order \(9\). By the preceding theorem, \(x\in R(G)\), and its image \(b\in T=R(G)/K\) has order \(3\). In \(G/W\), the image of \(x\) still has order \(9\), since \(W\) is a \(3'\)-group. Choose \(A_0\cong C_2^2\) in \(S\), and let \(F\) be the full preimage of \(A_0\times\langle b\rangle\) in \(G/W\). Then
\[
 F/V\cong C_2^2\times C_3,\qquad 9\in\omega(F),
 \qquad 18\notin\omega(F).
\]
This contradicts Lemma~\ref{lem:elementary-three-kernel-nine-eighteen}.
\end{proof}

\clearpage
\subsection{An abelian Sylow \(3\)-subgroup in the \(L_2(7)\) branch}

We retain the notation of the preceding sections. In particular, \(H\) is the normal Hall \(5'\)-subgroup of \(K\), \(Q=K/H\), and
\[
 G/H=Q\rtimes(S\times T),\qquad S\cong L_2(7),\quad
 T=R(G)/K,\quad \omega(T)\subseteq\{1,2,3\}.
\]
We write the action of the complement on \(Q\) as conjugation.

\Needspace{11\baselineskip}
\begin{theorem}\label{thm:L27-abelian-three-five-fixedfree}
Suppose that a Sylow \(3\)-subgroup of \(H\) is abelian. For every \(x\in G\) of order \(9\), the element \(xH\) acts on \(Q\) without nontrivial fixed points. Consequently, \(Q\) is nilpotent of class at most \(2\).
\end{theorem}
\begin{proof}
By Theorems~\ref{thm:L27-K-no-nine} and \ref{thm:L27-nine-simple-part-nonabelian}, we have \(x\in R(G)\setminus K\), \(|xH|=3\), and its image \(b\in T\) has order \(3\). First consider \(C_Q(xH)\). If it contains \(C_5^2\), the full preimage of \(\langle xH\rangle\times C_5^2\) in \(G\) is solvable and has a Hall \(\{3,5\}\)-subgroup \(J\) containing \(x\). Its normal subgroup \(P=J\cap H\) is a Sylow \(3\)-subgroup of \(H\), and \(J/P\cong C_3\times C_5^2\). By assumption and Theorem~\ref{thm:L27-K-no-nine}, the subgroup \(P\) is elementary abelian. Lemma~\ref{lem:elementary-three-five-nine} then excludes order \(9\) in \(J\), a contradiction. Hence \(C_Q(xH)\) does not contain \(C_5^2\). Since \(\exp Q=5\), this group is either trivial or cyclic of order \(5\).

In \(Q\rtimes\langle b\rangle\), the element \(xH\) defines a complement to the normal \(5\)-subgroup \(Q\). By the Schur--Zassenhaus theorem, this complement is conjugate to \(\langle b\rangle\) by an element of \(Q\). Thus \(C_Q(b)\) is also trivial or of order \(5\). But \(b\) commutes with \(S\), so \(C_Q(b)\) is \(S\)-invariant. If it had order \(5\), the perfect group \(S\) would act trivially on it, since \(\operatorname{Aut}(C_5)\cong C_4\). An element of order \(7\) in \(S\), together with an element of \(C_Q(b)\), would then give order \(35\) in the section \(G/H\). This contradiction proves \(C_Q(xH)=1\).

We prove the class bound directly. The action of \(b\), of order \(3\), on \(Q\) is coprime, so lifting fixed points shows that it acts fixed-point-freely on each lower central factor \(L_i=\gamma_i(Q)/\gamma_{i+1}(Q)\). These factors are \(\mathbb F_5\)-spaces, since \(\exp Q=5\). Over \(\mathbb F_{25}\), they decompose into eigenspaces with eigenvalues \(\zeta,\zeta^2\), where \(\zeta^3=1\), \(\zeta\ne1\). The bracket in the associated graded Lie ring multiplies eigenvalues. A bracket of vectors of different weights has weight \(1\), and hence vanishes. If the first two weights are equal and the third is the same, the triple bracket has weight \(1\) and vanishes. If the third weight is the other one, the Jacobi identity expresses the triple bracket in terms of brackets of original vectors with distinct weights, which vanish. Thus \(L_3=0\), since it is generated by triple brackets from \(L_1\). Hence \(\gamma_3(Q)=\gamma_4(Q)\); a finite \(5\)-group is nilpotent, so \(\gamma_3(Q)=1\).
\end{proof}

\Needspace{11\baselineskip}
\begin{lemma}\label{lem:L27-abelian-two-exp-four}
If a Sylow \(2\)-subgroup of \(H\) is abelian, its exponent divides \(4\).
\end{lemma}
\begin{proof}
Choose a Hall \(\{2,5\}\)-subgroup \(J\) of \(K\). The subgroup \(P=J\cap H\) is normal in \(J\), is a Sylow \(2\)-subgroup of \(H\), and is abelian by assumption. A Sylow \(5\)-subgroup of \(J\) contains \(A\cong C_5^6\), acting on \(P\). Suppose that \(P\) contains an element of order \(8\). Then \(M=\Omega_3(P)/\Omega_2(P)\ne0\) is an elementary abelian \(2\)-space. Every \(A\)-weight on \(M\) has nontrivial kernel, so some element \(a\in A\) of order \(5\) fixes a nonzero class. Coprime lifting of fixed points gives a fixed lift of this class in \(\Omega_3(P)\), necessarily of order \(8\). Together with \(a\), it gives the forbidden order \(40\).
\end{proof}

\Needspace{11\baselineskip}
\begin{lemma}\label{lem:abelian-two-five-nine-analogue}
Let \(E\) be a finite group with a normal abelian \(2\)-subgroup \(V\) of exponent dividing \(4\), and
\[
 E/V\cong C_{2^e}\times C_5^2,\qquad 0\le e\le2.
\]
If \(40\notin\omega(E)\), then \(8\notin\omega(E)\).
\end{lemma}
\begin{proof}
If \(e=0\), every \(2\)-element of \(E\) lies in \(V\), and the assertion is immediate. Suppose \(e\ge1\). The full preimage of the \(5\)-factor has a complement \(A\cong C_5^2\) to \(V\). Since \(A\) acts trivially on \(E/V\), coprime lifting of fixed points gives \(C_E(A)V=E\). Lift a generator of \(C_{2^e}\) and take its \(2\)-part to obtain \(z\in C_E(A)\) mapping to a generator of this factor. The order of \(z\) is at most \(4\): otherwise a power of \(z\) would have order \(8\), and together with an element of order \(5\) in \(A\) would give \(40\). Thus \(z^4=1\).

Every element of order \(8\) would have the form \(v z^i\), \(v\in V\), since its image in the \(5\)-factor is trivial. Since \(V\) is abelian and \(z^4=1\), its fourth power, writing \(V\) additively, equals
\[
 N_i(v)=(1+z^i+z^{2i}+z^{3i})v
\]
For order \(8\), we have \(N_i(v)\ne0\). The operator \(N_i\) is an endomorphism of \(V\) and commutes with \(A\).

The subgroups \(C_V(B)\), with \(B\) ranging over the hyperplanes of \(A\cong C_5^2\), generate \(V\): their images generate \(V/\Phi(V)\), by the weight decomposition and coprime lifting. Thus, if \(N_i(V)\ne0\), there are a hyperplane \(B\) and \(w\in C_V(B)\) with \(N_i(w)\ne0\). Then \(wz^i\) centralizes an element of order \(5\) in \(B\), and its fourth power is nontrivial. A power of this element has order \(8\), giving the forbidden order \(40\).
\end{proof}

\Needspace{11\baselineskip}
\begin{theorem}\label{thm:L27-small-five-centre-two-nonabelian}
If a Sylow \(3\)-subgroup of \(H\) is abelian and \(|Z(Q)|=5^6\), then a Sylow \(2\)-subgroup of \(H\) is nonabelian. In this case, either \(Q\cong C_5^6\), or \(Q\) is a special \(5\)-group of class \(2\) with \(Q'=\Phi(Q)=Z(Q)\cong C_5^6\).
\end{theorem}
\begin{proof}
Put \(Z=Z(Q)\). An element of order \(7\) in \(S\) has no nonzero fixed vectors on \(Z\), since otherwise order \(35\) occurs in \(G/H\). Since \(\operatorname{ord}_7(5)=6=\dim_{\mathbb F_5} Z\), the restriction of \(Z\) to this cyclic subgroup is irreducible; therefore \(Z\) is irreducible as an \(\mathbb F_5S\)-module as well.

Choose \(x\in G\) of order \(9\). By the preceding theorem, its image \(b\in T\) acts on \(Q\), hence on \(Z\), without nonzero fixed vectors. The field \(D=\operatorname{End}_S(Z)\) is a subfield of \(\operatorname{End}_{C_7}(Z)\cong\mathbb F_{5^6}\). It contains an element of order \(3\), the image of \(b\); hence \(D=\mathbb F_{25}\) or \(\mathbb F_{5^6}\). The latter is impossible: \(Z\) would then be one-dimensional over \(D\), so the perfect group \(S\) would act trivially on \(Z\). Thus \(D=\mathbb F_{25}\).

The image of \(T\) in \(D^*\) is cyclic, contains an element of order \(3\), and contains no element of order \(6\), since it is a quotient of \(T\). This image therefore has order \(3\), and every \(2\)-element of \(T\) acts trivially on \(Z\).

Choose \(y\in G\) of order \(8\), and denote its image in \(S\times T\) by \((a,c)\). The order of \(a\) is \(1,2\), or \(4\), and \(c\) is a \(2\)-element. Conjugation by the \(Q\)-component is trivial on \(Z\), as is the action of \(c\). Thus \(C_Z(yH)=C_Z(a)\). If \(a=1\), this space is nonzero; if \(|a|=4\), it is nonzero by Lemma~\ref{lem:L27-seven-free-five-modules-four}. If \(|a|=2\), it is also nonzero: otherwise \(a\) would act as \(-I\), impossible in a faithful action of the simple centerless group \(S\). The space \(C_Z(a)\) is invariant under \(D=\mathbb F_{25}\), so its dimension over \(\mathbb F_5\) is at least \(2\). Choose \(A\cong C_5^2\le C_Z(yH)\).

The order of \(yH\) is at most \(4\): its image in \(S\times T\) has that order, and the kernel \(Q\) has odd order. Let \(F\) be the full preimage of \(\langle yH\rangle\times A\) in \(G\), and let \(J\) be a Hall \(\{2,5\}\)-subgroup of \(F\) containing \(y\). Then \(P=J\cap H\) is a Sylow \(2\)-subgroup of \(H\), normal in \(J\), and
\[
 J/P\cong C_{2^e}\times C_5^2,\qquad 0\le e\le2.
\]
If \(P\) were abelian, its exponent would divide \(4\) by Lemma~\ref{lem:L27-abelian-two-exp-four}. The preceding lemma would then exclude the element \(y\) of order \(8\), a contradiction. Thus \(P\) is nonabelian.

Finally, \(Q\) has class at most \(2\) by Theorem~\ref{thm:L27-abelian-three-five-fixedfree}. Its derived subgroup \(Q'\) is an \(S\)-invariant subgroup of \(Z\), so irreducibility of \(Z\) gives \(Q'=1\) or \(Q'=Z\). In the first case, \(Q=Z\cong C_5^6\). In the second, \(\Phi(Q)=Q'Q^5=Q'=Z\), since \(\exp Q=5\).
\end{proof}

\clearpage
\subsection{The Sylow \(2\)-subgroup when the Sylow \(3\)-subgroup is abelian}

In this subsection we remove the restriction on the order of \(Z(Q)\) from the first part of Theorem~\ref{thm:L27-small-five-centre-two-nonabelian}. All other notation and assumptions remain unchanged.

\Needspace{11\baselineskip}
\begin{lemma}\label{lem:L27-seven-free-fixed-two-rank}
Let \(W\ne0\) be an \(\mathbb F_5S\)-module for \(S=L_2(7)\), on which an element of order \(7\) acts without nonzero fixed vectors. Then for every \(2\)-element \(a\in S\),
\[
 \dim_{\mathbb F_5}C_W(a)\ge2.
\]
\end{lemma}
\begin{proof}
Since \(5\nmid|S|\), the module is completely reducible. By Lemma~\ref{lem:L27-seven-free-five-modules-four}, over the algebraic closure its simple constituents have types \(3,3',6\); for an element of order \(4\), their fixed-space dimensions are \(1,1,2\), respectively. The two three-dimensional types are interchanged by the field Frobenius automorphism in characteristic \(5\). Indeed, their nontrivial \(7\)-weights form the two orbits of the normalizer \(7{:}3\), and taking fifth powers of seventh roots of unity interchanges these orbits. Hence in a module defined over \(\mathbb F_5\), the types \(3\) and \(3'\) have equal multiplicities. The assertion for an element of order \(4\) follows immediately.

Let \(a\) be an involution. On each three-dimensional type it has eigenvalues \(1,-1,-1\): the action is nontrivial, and the determinant is \(1\), since \(S\) is perfect. The six-dimensional type is realized by the zero-coordinate-sum submodule of the permutation module on the seven nonzero vectors of \(\mathbb F_2^3\). An involution in \(\operatorname{GL}_3(2)\) fixes three of these seven vectors and interchanges the others in two pairs; removing the constant constituent leaves a fixed space of dimension \(4\). Again, the conjugate three-dimensional types together give a fixed space of dimension at least two. For \(a=1\), the assertion is also immediate: any nontrivial \(7\)-weight over \(\mathbb F_5\) has degree \(\operatorname{ord}_7(5)=6\).
\end{proof}

\Needspace{11\baselineskip}
\begin{theorem}\label{thm:L27-abelian-three-nonabelian-two}
In the \(L_2(7)\) branch under consideration, if a Sylow \(3\)-subgroup of \(H\) is abelian, then a Sylow \(2\)-subgroup of \(H\) is nonabelian. In particular, the Sylow \(2\)- and \(3\)-subgroups of \(H\) cannot both be abelian.
\end{theorem}
\begin{proof}
Put \(Z=Z(Q)\ne1\), and choose \(x\in G\) of order \(9\). By Theorem~\ref{thm:L27-abelian-three-five-fixedfree}, its image \(b\in T\) acts on \(Q\), hence on \(Z\), without nonzero fixed vectors. In particular, the image of \(b\) in \(\operatorname{GL}(Z)\) has order \(3\).

Choose \(y\in G\) of order \(8\), and write its image in \(G/K=S\times T\) as \((a,c)\). Then \(a\) is a \(2\)-element of \(S\), and \(c^2=1\). The action of \(T\) on \(Z\) is well defined, since \(Q\) centralizes \(Z\). Put \(W=C_Z(c)\). If \(W=0\), the involutory action of \(c\) in characteristic \(5\) would be \(-I\). The images of \(b\) and \(c\) would then commute and have orders \(3,2\), so their product would have order \(6\) in a quotient of \(T\). This is impossible, since \(6\notin\omega(T)\). Hence \(W\ne0\).

The subspace \(W\) is \(S\)-invariant, since \(S\) and \(T\) commute in the complement to \(Q\). An element of order \(7\) in \(S\) acts on \(W\) without nonzero fixed vectors, since otherwise \(G/H\) would contain an element of order \(35\). By the preceding lemma, \(\dim C_W(a)\ge2\). The inner \(Q\)-component of the image \(yH\) centralizes \(Z\), so \(C_W(a)\le C_Z(yH)\). Choose \(A\cong C_5^2\) inside \(C_Z(yH)\).

As in the proof of Theorem~\ref{thm:L27-small-five-centre-two-nonabelian}, \(|yH|\le4\). The full preimage of \(\langle yH\rangle\times A\) in \(G\) is solvable; a Hall \(\{2,5\}\)-subgroup \(J\) of it can be chosen to contain \(y\). Then \(P=J\cap H\) is a Sylow \(2\)-subgroup of \(H\), normal in \(J\), and
\[
 J/P\cong C_{2^e}\times C_5^2,\qquad 0\le e\le2.
\]
If \(P\) is abelian, its exponent divides \(4\) by Lemma~\ref{lem:L27-abelian-two-exp-four}. But Lemma~\ref{lem:abelian-two-five-nine-analogue} then excludes the element \(y\) of order \(8\), since \(40\notin\omega(G)\). Thus \(P\) is nonabelian. No restriction on \(|Z(Q)|\) was used in this proof.
\end{proof}

\clearpage
\subsection{Nonabelian Sylow subgroups and order \(8\) inside \(K\)}

We retain the notation of the preceding sections, in particular \(Q=K/H\), \(S=L_2(7)\). The following restrictions also apply when a Sylow \(3\)-subgroup of \(H\) is nonabelian.

\Needspace{11\baselineskip}
\begin{lemma}\label{lem:L27-nonabelian-p-five-central-product}
Let \(p\in\{2,3\}\), and suppose that a Sylow \(p\)-subgroup of \(H\) is nonabelian. Then \(K\) contains a direct product
\[
 D\times C_5^4,
\]
where \(D\) is a nonabelian two-generator \(p\)-group. For \(p=2\), one can choose \(D\) of exponent \(4\); for \(p=3\), the group \(D\) is isomorphic to the nonabelian group of order \(27\) and exponent \(3\).
\end{lemma}
\begin{proof}
Choose a Hall \(\{p,5\}\)-subgroup \(J\) of \(K\). The group \(P=J\cap H\) is normal in \(J\) and is a Sylow \(p\)-subgroup of \(H\). A Sylow \(5\)-subgroup of \(J\) contains \(A\cong C_5^6\), acting on \(P\).

As in Lemma~\ref{lem:five-rank-exponent-three}, the groups \(C_P(B)\), with \(B\) ranging over the hyperplanes of \(A\), generate \(P\). Indeed, irreducible constituents of \(P/\Phi(P)\) are fixed by suitable hyperplanes, and their fixed vectors lift under the coprime action. Since \(P\) is nonabelian, there are noncommuting \(u,v\in P\), fixed by hyperplanes \(B_u,B_v\le A\), respectively. Their intersection has rank at least \(4\). Choose \(A_0\cong C_5^4\) in it, and put \(D=\langle u,v\rangle\). Then \(D\) is nonabelian, \([D,A_0]=1\), and \(D\cap A_0=1\), as required.

If \(p=2\), every element of order \(8\) in \(D\) would commute with an element of order \(5\) in \(A_0\), giving the forbidden order \(40\). Orders \(16,32,\ldots\) are also absent from \(G\), so \(\exp D\mid4\). A nonabelian \(2\)-group cannot have exponent \(2\), whence \(\exp D=4\).

If \(p=3\), Theorem~\ref{thm:L27-K-no-nine} gives \(\exp D=3\). A group of exponent \(3\) is \(2\)-Engel, so \([u,v,u]=[u,v,v]=1\) and \(D'\le Z(D)\). Hence \(D\) has class \(2\) and order at most \(3^3\). Since \(D\) is nonabelian, its abelianization has order at least \(3^2\), and \(D'\ne1\); thus \(|D|=27\).
\end{proof}

\Needspace{11\baselineskip}
\begin{corollary}\label{cor:L27-abelian-three-twenty-K}
If a Sylow \(3\)-subgroup of \(H\) is abelian, then
\[
 20\in\omega(K).
\]
Thus in this case all orders \(1,2,3,4,5,6,10,15,20,30\) occur inside \(K\).
\end{corollary}
\begin{proof}
By Theorem~\ref{thm:L27-abelian-three-nonabelian-two}, a Sylow \(2\)-subgroup of \(H\) is nonabelian. The preceding lemma gives a nonabelian group \(D\) of exponent \(4\) centralizing \(C_5^4\). An element of order \(4\) in \(D\) and an element of order \(5\) in the complementary direct factor give order \(20\). Order \(30\) has already been shown to occur in \(K\) by Theorem~\ref{thm:L27-fitting-large}; the other listed orders divide \(20\) or \(30\).
\end{proof}

\Needspace{11\baselineskip}
\begin{theorem}\label{thm:L27-eight-two-rank-bound}
Let \(P\) be a Sylow \(2\)-subgroup of \(K\). If \(8\in\omega(K)\), then
\[
 \dim_{\mathbb F_2} P/\Phi(P)\ge24,
 \qquad |P|\ge2^{24},
 \qquad |K|\ge2^{24}\cdot3^2\cdot5^6.
\]
Moreover, the number of elements of order \(8\) in \(P\) is divisible by \(5^6\).
\end{theorem}
\begin{proof}
Choose a Hall \(\{2,5\}\)-subgroup \(J\) of \(K\). Its normal Sylow \(2\)-subgroup \(P_0=J\cap H\) is a Sylow \(2\)-subgroup of \(K\). A Sylow \(5\)-subgroup of \(J\) contains \(A\cong C_5^6\), acting on \(P_0\). By Sylow's theorem, \(P_0\) contains an element \(u\) of order \(8\).

No nonidentity element of \(A\) fixes \(u\), since otherwise order \(40\) occurs. Thus the action of \(A\) on \(P_0\) is faithful. Its action on \(P_0/\Phi(P_0)\) is also faithful, since the kernel of the action of automorphisms of a finite \(2\)-group on its Frattini quotient is a \(2\)-group. Every nontrivial irreducible \(\mathbb F_2 A\)-module has dimension \(\operatorname{ord}_5(2)=4\), and its kernel is a hyperplane in \(A\). Faithfulness requires the intersection of the kernels to be trivial, so at least six such constituents are necessary. This gives \(\dim P_0/\Phi(P_0)\ge24\), hence \(|P_0|\ge2^{24}\).

A Sylow \(5\)-subgroup of \(K\) contains \(A\), and a Sylow \(3\)-subgroup of \(H\) is noncyclic and has order at least \(3^2\). This proves the bound for \(|K|\). Finally, the action of \(A\) on the elements of order \(8\) in \(P_0\) is free: every nontrivial stabilizer would give order \(40\). All orbits have length \(5^6\), so the number of these elements is divisible by \(5^6\). Conjugacy of Sylow subgroups transfers all assertions to \(P\).
\end{proof}

\clearpage
\subsection{Actions of elements of orders \(8,9\) on normal \(2\)- and \(3\)-subgroups}

We retain the assumptions of the remaining branch for \(S=L_2(7)\), including \(K=[R(G)\cap G^\infty,G^\infty]\). The restrictions below do not assume that a Sylow \(3\)-subgroup of \(H\) is abelian.

\Needspace{11\baselineskip}
\begin{lemma}\label{lem:L27-fixedfree-nine-eight-fitting}
Put \(N_2=O_2(K)\), \(N_3=O_3(K)\). For all elements \(x,y\in G\) of orders \(9,8\), respectively,
\[
 C_{N_2}(x)=1,\qquad C_{N_3}(y)=1.
\]
If \(N_2\ne1\), then \(Z(N_2)\) is elementary abelian and
\[
 \dim_{\mathbb F_2} Z(N_2)=2r+6s
\]
for some nonnegative integers \(r,s\) with \(r+s>0\). In particular, its dimension is even and at least \(2\). If \(N_3\ne1\), then \(Z(N_3)\) is elementary abelian.
\end{lemma}
\begin{proof}
The groups \(N_2,N_3\) are characteristic in \(K\), hence normal in \(G\). If \(C_{N_2}(x)\ne1\), it contains an involution, which commutes with \(x\) and gives order \(18\), impossible in \(G\). Similarly, a nontrivial \(C_{N_3}(y)\) contains an element of order \(3\), which together with \(y\) gives the forbidden order \(24\).

The centers \(Z(N_2),Z(N_3)\) are normal abelian subgroups of \(K\), so by Lemma~\ref{lem:L27-normal-abelian-exp} their exponents divide \(2,3\), respectively. For \(N_2\ne1\), put \(V=Z(N_2)\ne0\). The action of \(x\) on \(V\) has no nonzero fixed vectors. Over \(\mathbb F_2\), the polynomial \(t^9-1\) is a product of distinct irreducible polynomials of degrees \(1,2,6\):
\[
 t^9-1=(t-1)(t^2+t+1)(t^6+t^3+1).
\]
The last factor is irreducible, since \(\operatorname{ord}_9(2)=6\). There is no trivial constituent, giving the stated dimension formula.
\end{proof}

\Needspace{11\baselineskip}
\begin{theorem}\label{thm:L27-eight-normal-two-depth}
If \(8\in\omega(O_2(K))\), then for \(N=O_2(K)\),
\[
 \dim_{\mathbb F_2}N/\Phi(N)\ge24,
 \qquad N^2\text{ is nonabelian},\qquad \operatorname{cl}(N)\ge3.
\]
In particular, \(|N|\ge2^{24}\). Here \(N^2\) denotes the subgroup generated by squares of elements of \(N\).
\end{theorem}
\begin{proof}
Choose \(u\in N\) of order \(8\), and a subgroup \(A\cong C_5^6\) in a Sylow \(5\)-subgroup of \(K\), which exists by Lemma~\ref{lem:L27-global-five-complement}. The group \(A\) normalizes \(N\). If a nonidentity element of \(A\) centralizes \(u\), the forbidden order \(40\) occurs. Thus the action of \(A\) on \(N\) is faithful. It is also faithful on \(N/\Phi(N)\): the kernel of the action of automorphisms of a finite \(2\)-group on its Frattini quotient is a \(2\)-group.

As in the proof of Theorem~\ref{thm:L27-eight-two-rank-bound}, every nontrivial irreducible \(\mathbb F_2 A\)-module has dimension \(4\) and kernel of index \(5\). Faithfulness requires at least six such constituents, giving the first bound and \(|N|\ge2^{24}\).

The subgroup \(N^2\) is characteristic in \(N\) and normal in \(K\). It contains the element \(u^2\) of order \(4\). If \(N^2\) were abelian, Lemma~\ref{lem:L27-normal-abelian-exp} would give it exponent dividing \(2\), a contradiction. Hence \(N^2\) is nonabelian.

If \(N\) had class at most \(2\), its derived subgroup would lie in \(Z(N)\), which is elementary abelian by the preceding lemma. Thus for all \(a,b\in N\), \([a^2,b]=[a,b]^2=1\). Hence \(N^2\le Z(N)\) would be abelian, contrary to what has been proved. Thus \(\operatorname{cl}(N)\ge3\).
\end{proof}

\clearpage
\subsection{The cube of an element of order \(9\) and the normal \(2\)-subgroup}

We retain the notation and assumptions of the preceding subsection. In particular, \(N=O_2(K)\) is characteristic in \(K\) and normal in \(G\). The assertions below hold regardless of whether a Sylow \(3\)-subgroup of \(H\) is abelian.

\Needspace{11\baselineskip}
\begin{lemma}\label{lem:fpf-three-on-two-class-two}
Let a finite \(2\)-group \(P\) admit an automorphism \(\alpha\) of order \(3\) with \(C_P(\alpha)=1\). Then \(P\) is nilpotent of class at most \(2\).
\end{lemma}
\begin{proof}
Put \(L_i=\gamma_i(P)/\gamma_{i+1}(P)\) and \(V_i=L_i/2L_i\). By coprime lifting of fixed points, \(\alpha\) acts without nonzero fixed points on every section \(L_i\), and then on \(V_i\). The commutator induces bilinear maps \(V_i\times V_j\to V_{i+j}\) satisfying the Jacobi identity. In particular, \(V_3\) is generated by triple commutators of elements of \(V_1\).

Extend the field to \(\mathbb F_4\), and denote the nontrivial cube roots of unity by \(\zeta,\zeta^2\). On each \(V_i\), the automorphism \(\alpha\) is diagonalizable and has no eigenvalue \(1\). The bracket of vectors of weights \(\zeta,\zeta^2\) has weight \(1\), so it vanishes. The triple bracket of three vectors of the same weight also has weight \(1\) and vanishes. If both weights occur among the three, the Jacobi identity expresses the triple bracket in terms of brackets of pairs with distinct weights. Thus \(V_3=0\).

The group \(L_3\) is finite abelian; for a \(2\)-group, \(L_3/2L_3=0\) is possible only if \(L_3=0\). Hence \(\gamma_3(P)=\gamma_4(P)\), and nilpotence of the finite \(2\)-group gives \(\gamma_3(P)=1\).
\end{proof}

\Needspace{11\baselineskip}
\begin{lemma}\label{lem:L27-nine-invariant-two-square}
Let \(x\in G\) have order \(9\), and let a finite \(2\)-subgroup \(P\le G\) be normalized by \(x\). Then \(C_P(x)=1\), the order of \(P\) is a square, and
\[
 \dim_{\mathbb F_2}P/\Phi(P)=2r+6s
\]
for some nonnegative \(r,s\).
\end{lemma}
\begin{proof}
A nontrivial centralizer \(C_P(x)\) would contain an involution, which together with \(x\) would give the forbidden order \(18\). Thus \(C_P(x)=1\). In the action of \(\langle x\rangle\) on the set \(P\), the identity forms the only orbit of size one; all other orbits have length \(3\) or \(9\). Hence \(|P|\equiv1\pmod3\), so for \(|P|=2^m\), the integer \(m\) is even.

Coprime lifting of fixed points shows that \(x\) acts without nonzero fixed vectors on \(P/\Phi(P)\). The polynomial \(t^9-1\) over \(\mathbb F_2\) has three irreducible factors of degrees \(1,2,6\), and there is no constituent of degree \(1\). The formula follows.
\end{proof}

\Needspace{11\baselineskip}
\begin{theorem}\label{thm:L27-eight-nine-cube-two-action}
Suppose that \(8\in\omega(N)\), where \(N=O_2(K)\). For every \(x\in G\) of order \(9\),
\[
 C_N(x^3)\ne1,\qquad [N,x^3]\ne1.
\]
Moreover, \(C_N(x^3)\) is a \(2\)-group of class at most \(2\), and its order is a square at least \(4\). The decomposition of the action of \(x\) on \(N/\Phi(N)\) includes at least one irreducible constituent of dimension \(6\).
\end{theorem}
\begin{proof}
By Theorem~\ref{thm:L27-eight-normal-two-depth}, \(N\) has class at least \(3\). If \(C_N(x^3)=1\), conjugation by \(x^3\) would induce a fixed-point-free automorphism of order \(3\) on \(N\). Lemma~\ref{lem:fpf-three-on-two-class-two} would give \(\operatorname{cl}(N)\le2\), a contradiction.

If \([N,x^3]=1\), the action of \(x\) on \(N\) would have order dividing \(3\). By the preceding lemma, it is fixed-point-free. For \(N\ne1\), it is nontrivial, so again it would be a fixed-point-free automorphism of order \(3\), giving \(\operatorname{cl}(N)\le2\). Thus both stated conditions hold.

Put \(C=C_N(x^3)\ne1\). It is normalized by \(x\), and \(x^3\) acts trivially on it. The action of \(x\) on \(C\) is fixed-point-free, hence has order \(3\). Lemmas~\ref{lem:fpf-three-on-two-class-two} and \ref{lem:L27-nine-invariant-two-square} give \(\operatorname{cl}(C)\le2\) and show that \(|C|\) is a square at least \(4\).

The action of \(x^3\) on \(N/\Phi(N)\) is nontrivial: otherwise an automorphism of order \(3\) would lie in the kernel of the action on the Frattini quotient, which is a \(2\)-group, and would centralize all of \(N\). By Lemma~\ref{lem:L27-nine-invariant-two-square}, the action of \(x\) on this quotient has no trivial constituents; its irreducible constituents have dimensions \(2,6\). The cube of \(x\) acts trivially on the constituents of dimension \(2\), so nontriviality of \(x^3\) requires at least one constituent of dimension \(6\).
\end{proof}

\Needspace{11\baselineskip}
\begin{lemma}\label{lem:L27-eight-fourth-on-nonabelian-three}
If \(O_3(K)\) is nonabelian, then for every \(y\in G\) of order \(8\),
\[
 C_{O_3(K)}(y^4)\ne1.
\]
\end{lemma}
\begin{proof}
Put \(P=O_3(K)\), and suppose that the involution \(y^4\) acts on \(P\) without nontrivial fixed points. For any finite \(P\) with such an automorphism \(\alpha\), the map \(g\mapsto g^{-1}\alpha(g)\) is injective, hence surjective. Every element \(v\) of its image satisfies \(\alpha(v)=v^{-1}\). Thus \(\alpha\) inverts every element of \(P\), and inversion is an automorphism only for an abelian group. This contradicts the assumption on \(P\).
\end{proof}

\clearpage
\subsection{Excluding nilpotent \(H\) when the Sylow \(3\)-subgroup is abelian}

In this subsection we consider the remaining branch \(S=L_2(7)\). As before, \(H\) is the characteristic normal Hall \(5'\)-subgroup of \(K\), \(Q=K/H\), and \(G/H=Q\rtimes(S\times T)\), where \(T=R(G)/K\), \(\omega(T)\subseteq\{1,2,3\}\).

\Needspace{11\baselineskip}
\begin{theorem}\label{thm:L27-abelian-three-H-not-nilpotent}
If a Sylow \(3\)-subgroup of \(H\) is abelian, \(H\) is not nilpotent.
\end{theorem}
\begin{proof}
Suppose otherwise. Then \(H=H_2\times H_3\), where \(H_i\) are the Sylow subgroups of \(H\). Both are characteristic in \(H\) and normal in \(G\). The group \(V=H_3\) is nontrivial and elementary abelian by Theorem~\ref{thm:L27-K-no-nine}. The subgroup \(H\) centralizes \(V\), so conjugation defines a representation
\[
 \rho:G/H=Q\rtimes(S\times T)\longrightarrow\operatorname{GL}(V).
\]
Put \(U=\rho(Q)\). An element of order \(7\) in \(S\) acts on \(Q\) without nontrivial fixed points, since otherwise order \(35\) occurs in \(G/H\). By coprime lifting of fixed points, it acts in the same way on the quotient \(U\). Hence, if \(U\ne1\), the action of \(S\) on \(U\) is nontrivial. Moreover, the image of \(S\) in \(E/U\), where \(E=\rho(G/H)\), is isomorphic to \(S\): its kernel is normal in the simple group \(S\); a full kernel would mean that the image of \(S\) lies in the \(5\)-group \(U\), hence is trivial, contrary to the fixed-point-free action of order \(7\) on the nontrivial group \(U\).

Choose \(x\in G\) of order \(9\). By Theorems~\ref{thm:L27-nine-simple-part-nonabelian} and \ref{thm:L27-abelian-three-five-fixedfree}, its image in \(T\) is an element \(b\) of order \(3\), with \(C_Q(b)=1\). For the latter assertion, we use the conjugacy of \(xH\) and \(b\) by an element of \(Q\), and the proof of the cited theorem.

Choose \(y\in G\) of order \(8\). Write \(yH=q(a,c)\), where \(q\in Q\), \(a\in S\), \(c\in T\). The order of \(a\) is \(1,2\), or \(4\), and \(c^2=1\). Since \(yH\) is a \(2\)-element, its order equals that of its image \((a,c)\) in \(S\times T\), and is at most \(4\). The action of \(y\) on \(V\) has no nonzero fixed vectors, since otherwise the forbidden order \(24\) occurs. We prove \(U=1\) by considering the three possibilities for \(a\).

Suppose \(|a|=4\) and \(U\ne1\). The center \(Z(U)\) is nontrivial and \(S\)-invariant; an element of order \(7\) in \(S\) acts on it without nontrivial fixed points. Thus the action of the simple group \(S\) on \(Z(U)\) is faithful, and \(a^2\) acts nontrivially. Over the algebraic closure in characteristic \(3\), decompose \(V\) into weight spaces of the elementary abelian group \(Z(U)\). The representation of \(U\le\operatorname{GL}(V)\) is faithful, so the occurring weights separate the elements of \(Z(U)\). The component \(q\) induces an inner automorphism of \(U\) on \(Z(U)\), hence acts trivially; \(c^2=1\), so the action of \((yH)^2\) on \(Z(U)\) coincides with that of \(a^2\). Hence \((yH)^2\) moves some occurring weight. Its orbit under \(yH\) has length \(4\). Since \((yH)^4=1\), the sum of a nonzero vector and its three images under \(yH\) is nonzero and fixed. The fixed space is defined over \(\mathbb F_3\), contrary to \(C_V(y)=0\). Thus \(U=1\) when \(|a|=4\).

Now suppose \(|a|=2\). Then \(|yH|=2\), and the action of \(y\) on \(V\) is a fixed-point-free involution in characteristic \(3\); hence \(\rho(yH)=-I\). If \(U\ne1\), this scalar is also central in \(E/U\). The image of \(\rho(q)\) there is trivial, and the image of \(c\) centralizes the image of \(S\). Thus the image of \(a\) centralizes the image of \(S\) in \(E/U\). But the image of \(S\) there is isomorphic to the centerless simple group \(S\), and \(a\ne1\), a contradiction. Hence \(U=1\) in this case too.

Finally, suppose \(a=1\). Then \(c\ne1\): otherwise \(yH\) would lie both in \(Q\) and in a \(2\)-group, hence would be the identity, and all of \(V\) would be fixed by \(y\). We have \(|c|=|yH|=2\) and \(\rho(yH)=-I\). In \(Q\rtimes\langle c\rangle\), complements to \(Q\) are conjugate by elements of \(Q\), so \(yH\) is conjugate to \(c\), and \(\rho(c)=-I\). If \(\rho(b)\) had order \(3\), the product \(\rho(c)\rho(b)\) would have order \(6\) in a quotient of \(T\), which is impossible. Hence \(\rho(b)=1\).

Since \(C_Q(b)=1\) and the action of \(b\) on \(Q\) is coprime, the standard equality \(Q=[Q,b]C_Q(b)\) gives \(Q=[Q,b]\). From \(\rho(b)=1\), we obtain \(\rho([Q,b])=1\); thus \(U=\rho(Q)=1\). This completes all three cases.

Choose a Sylow \(5\)-subgroup \(P_5\) of \(K\). The equality \(U=1\) means that \(P_5\) centralizes \(H_3\). The subgroup \(H_2\) is characteristic in \(H\), so \(W=H_2P_5\) is a subgroup. We have
\[
 K=H_3\times W,
\]
since \(H_3\) is abelian and centralizes both factors of \(W\). Thus \(W\) is a normal Hall \(3'\)-subgroup of \(K\). This contradicts Corollary~\ref{cor:L27-abelian-three-not-three-nilpotent}. Therefore \(H\) is not nilpotent.
\end{proof}

\clearpage

\subsection{A normal Sylow \(3\)-subgroup in \(H\)}

We retain the notation of the preceding sections and consider the remaining branch \(S=L_2(7)\). Below, a Sylow \(3\)-subgroup of \(H\) is assumed to be normal, but not necessarily abelian.

\Needspace{11\baselineskip}
\begin{theorem}\label{thm:L27-normal-sylow-three-structure}
Suppose that a Sylow \(3\)-subgroup \(P\) of \(H\) is normal in \(H\). Then \(P=O_3(K)\) and \(P\triangleleft G\). For every \(x\in G\) of order \(9\),
\[
 x^3\in P,\qquad 9\notin\omega(G/P).
\]
The subgroup \(N=O_2(K)\) has class at most \(2\), \(N^2\le Z(N)\), and \(\exp N\mid4\); in particular, \(8\notin\omega(N)\). If, in addition, \(8\in\omega(K)\), then a Sylow \(2\)-subgroup of \(H\) is not normal in \(H\). Moreover,
\[
 \dim_{\mathbb F_3}P/\Phi(P)\ge6.
\]
If \(P\) is nonabelian, then \(|P|\ge3^7\).
\end{theorem}
\begin{proof}
The normal Sylow subgroup \(P\) is unique and characteristic in \(H\), and hence normal in \(G\). Since \(K/H\) is a \(5\)-group, every normal \(3\)-subgroup of \(K\) lies in \(H\). Thus \(P=O_3(K)\).

For \(|x|=9\), we proved earlier that \(x^3\in H\); this element of order \(3\) belongs to the unique Sylow \(3\)-subgroup \(P\). Consequently, the image of \(x\) in \(G/P\) has order at most \(3\). If \(G/P\) had an element of order \(9\), it could be lifted to a \(3\)-element of \(G\) of order \(9\), since the spectrum of \(G\) does not contain \(27\). This contradiction proves \(9\notin\omega(G/P)\).

The normal subgroups \(N=O_2(K)\) and \(P=O_3(K)\) centralize one another, since their commutator subgroup lies in their trivial intersection. Hence \(x^3\) centralizes \(N\). If \(N\ne1\), the action of \(x\) on \(N\) has order \(3\) and is fixed-point-free: a nontrivial fixed point would give order \(18\). By Lemma~\ref{lem:fpf-three-on-two-class-two}, \(\operatorname{cl}(N)\le2\). The center of \(N\) is elementary abelian by Lemma~\ref{lem:L27-normal-abelian-exp}, so \([u^2,v]=[u,v]^2=1\) for all \(u,v\in N\). We obtain \(N^2\le Z(N)\) and \(\exp N\mid4\). For \(N=1\), the same assertions are immediate. All \(2\)-elements of \(K\) lie in \(H\), since \(K/H\) is a \(5\)-group. If a Sylow \(2\)-subgroup of \(H\) is normal, it is characteristic in \(H\), normal in \(K\), and contained in \(N\); hence \(8\notin\omega(K)\).

Put \(V=P/\Phi(P)\). The group \(L=G^\infty\) is perfect. Its action on \(V\) is nontrivial. Otherwise the image of \(L\) in \(\operatorname{Aut}(P)\) would lie in the kernel of the action on the Frattini quotient, which is a \(3\)-group. Being perfect, this image would be trivial. Then \(L\), and hence \(H\le L\), would centralize \(P\). In particular, \(P\) would be abelian and central in \(H\). The Schur--Zassenhaus theorem would give a decomposition \(H=P\times H_2\), making \(H\) nilpotent, contrary to Theorem~\ref{thm:L27-abelian-three-H-not-nilpotent}.

Let \(E\) be the image of \(L\) in \(\operatorname{GL}(V)\), and denote the corresponding representation by \(\rho\). It is perfect and nontrivial; the image of the solvable group \(K\) is normal and solvable in \(E\), and the quotient \(E/\rho(K)\) is a quotient of \(L/K\cong S\). If this quotient were trivial, \(E\) would be both solvable and perfect, and therefore trivial. Thus \(E/\rho(K)\cong S\), and \(7\mid|E|\). For \(d=\dim_{\mathbb F_3}V\), we have \(7\mid|\operatorname{GL}_d(3)|\). Since \(\operatorname{ord}_7(3)=6\), necessarily \(d\ge6\). If \(P\) is nonabelian, then \(\Phi(P)\ne1\), giving \(|P|\ge3^{d+1}\ge3^7\).
\end{proof}

\Needspace{11\baselineskip}
\begin{corollary}\label{cor:L27-normal-abelian-three-action-rank}
Suppose that a Sylow \(3\)-subgroup \(P\) of \(H\) is normal and abelian. Then \(P\) is elementary abelian, \(|P|\ge3^6\), and the group
\[
 D=H/C_H(P)\le\operatorname{Aut}(P)
\]
is a nontrivial \(2\)-group satisfying
\[
 \dim_{\mathbb F_2}D/\Phi(D)\ge3.
\]
\end{corollary}
\begin{proof}
Elementary abelianness follows from Theorem~\ref{thm:L27-K-no-nine}, and the bound on \(|P|\) follows from the preceding theorem. Since \(P\) is abelian, \(P\le C_H(P)\); the quotient \(H/P\) is a \(2\)-group, so \(D\) is also a \(2\)-group. If \(D=1\), then \(P\le Z(H)\), and, as above, \(H\) would be nilpotent, contrary to Theorem~\ref{thm:L27-abelian-three-H-not-nilpotent}.

The subgroup \(C_H(P)\) is characteristic in \(H\) and is normalized by \(L\). The equality \(H=[H,L]\) implies \(D=[D,L]\), where \(L\) acts on \(D\) by conjugation. If \(d=\dim D/\Phi(D)\le2\), then \(\operatorname{Aut}(D)\) is solvable: the kernel of its action on \(D/\Phi(D)\) is a \(2\)-group, and \(\operatorname{GL}_d(2)\) is solvable for \(d\le2\). The image of the perfect group \(L\) in this solvable group is trivial, so \([D,L]=1\), contradicting \(D\ne1\). Thus \(d\ge3\).
\end{proof}

\clearpage
\subsection{Excluding rank \(6\) for a normal abelian \(3\)-subgroup}

We retain the assumptions of the remaining branch \(S=L_2(7)\), the notation \(L=G^\infty\), \(K=[R(G)\cap L,L]\), and the normal Hall \(5'\)-subgroup \(H\le K\).

\Needspace{11\baselineskip}
\begin{theorem}\label{thm:L27-normal-abelian-three-rank-seven}
Suppose that a Sylow \(3\)-subgroup \(P\) of \(H\) is normal in \(H\) and abelian. Then
\[
 \dim_{\mathbb F_3}P\ge7,\qquad |P|\ge3^7.
\]
\end{theorem}
\begin{proof}
By Corollary~\ref{cor:L27-normal-abelian-three-action-rank}, the group \(P\) is elementary abelian and normal in \(G\), and \(\dim P\ge6\). Suppose that \(\dim P=6\). The action of \(L\) on \(P\) is nontrivial, and its image \(E\) is perfect and has a solvable normal image of \(K\) with quotient \(L_2(7)\), as proved in Theorem~\ref{thm:L27-normal-sylow-three-structure}. Thus \(E\) contains an element of order \(7\); lift it to an element \(a\in L\) of order \(7\).

Since \(\operatorname{ord}_7(3)=6\), the action of \(a\) on the six-dimensional \(\mathbb F_3\)-space \(P\) is irreducible. Hence
\[
 C_{\operatorname{GL}(P)}(a)\cong\mathbb F_{3^6}^{\,*},
\]
so this centralizer is cyclic.

Put \(D=H/C_H(P)\). By Corollary~\ref{cor:L27-normal-abelian-three-action-rank}, \(D\) is a \(2\)-group acting faithfully on \(P\), and \(d(D)=\dim D/\Phi(D)\ge3\). The subgroup \(D\) is normalized by the image of \(a\). The automorphism induced by \(a\) on \(D\) is nontrivial: otherwise \(D\) would lie in the cyclic centralizer of \(a\) in \(\operatorname{GL}(P)\), contradicting \(d(D)\ge3\).

We show that every nontrivial abelian \(a\)-invariant subgroup \(A\le D\) has order \(2\). Extend the field to an algebraic closure of characteristic \(3\), and decompose \(P\) into weight spaces for \(A\). The element \(a\) permutes the weights occurring; a nontrivial orbit would have length \(7\), which is impossible in a six-dimensional space. Thus every weight is fixed by \(a\). The weights separate the elements of \(A\), since the action of \(A\) is faithful, so \(a\) centralizes \(A\). Hence \(A\) lies in the cyclic centralizer \(C_{\operatorname{GL}(P)}(a)\). If \(A\) contained an element of order \(4\), its product with \(a\) would have order \(28\) in the image of \(G\) on \(P\). The spectrum of a quotient group is contained in \(\omega(G)\), whereas \(28\notin\Om\). Thus \(\exp A\mid2\), and cyclicity gives \(|A|=2\).

Apply to \(D\) the form of Thompson's theorem in Lemma~\ref{lem:critical-two-external}. There is a characteristic subgroup \(C\le D\) of class at most \(2\), with \(C/Z(C)\) elementary abelian, such that the nontrivial odd-order automorphism \(a\) acts nontrivially on \(C\). The center \(Z(C)\) is characteristic in \(D\), so the preceding argument gives order \(2\). The group \(C\) is nonabelian: otherwise \(C=Z(C)\cong C_2\), and \(a\) would act trivially on it. Consequently,
\[
 C'=Z(C)=\Phi(C)\cong C_2,
\]
so \(C\) is extraspecial of order \(2^{1+2m}\). The commutator defines a nondegenerate symplectic form on \(C/Z(C)\) over \(\mathbb F_2\). The action of \(a\) on this quotient is nontrivial: the kernel of the action of the automorphisms of \(C\) on \(C/Z(C)\) is a \(2\)-group, whereas \(|a|=7\). Therefore \(7\mid|\operatorname{Sp}_{2m}(2)|\). The formula
\[
 |\operatorname{Sp}_{2m}(2)|=2^{m^2}\prod_{i=1}^m(2^{2i}-1)
\]
implies \(m\ge3\), since \(\operatorname{ord}_7(2)=3\).

On the other hand, \(C\) acts faithfully on \(P\), and its central involution centralizes \(a\). The unique involution in the field \(\mathbb F_{3^6}^{\,*}\) is \(-I\); hence the center of \(C\) acts on \(P\) as \(-I\). Let \(M\le C\) be the preimage of a maximal totally isotropic subspace of \(C/Z(C)\). Then \(M\) is abelian and \(|M/Z(C)|=2^m\). Over the algebraic closure, choose a weight of \(M\) that occurs. It is nontrivial on \(Z(C)\); nondegeneracy of the commutator form shows that its orbit under conjugation by elements of \(C\) has length \(2^m\). All the corresponding weight spaces are nonzero, so \(\dim P\ge2^m\ge8\). This contradicts \(\dim P=6\).
\end{proof}

\clearpage
\subsection{Rank \(7\) and an abelian image of \(H\) on a normal \(3\)-subgroup}

\Needspace{11\baselineskip}
\begin{theorem}\label{thm:L27-rank-seven-abelian-action}
In the remaining branch \(S=L_2(7)\), suppose that a Sylow \(3\)-subgroup \(P\) of \(H\) is normal and abelian, \(\dim_{\mathbb F_3}P=7\), and \(D=H/C_H(P)\) is abelian. Then:
\begin{enumerate}
\item \(D\) is elementary abelian, and \(P\) decomposes into a direct sum of seven one-dimensional \(\mathbb F_3\)-spaces permuted by the action of \(L=G^\infty\);
\item every Sylow \(5\)-subgroup of \(K\) centralizes \(P\);
\item every element \(x\in G\) of order \(9\) centralizes \(P\).
\end{enumerate}
\end{theorem}
\begin{proof}
The group \(P\) is elementary abelian and normal in \(G\), and \(D\) is a \(2\)-group with \(d(D)\ge3\), by Corollary~\ref{cor:L27-normal-abelian-three-action-rank}. The image of the action of \(L\) on \(P\) has order divisible by \(7\), by the proof of Theorem~\ref{thm:L27-normal-sylow-three-structure}. Choose \(a\in L\) of order \(7\) acting nontrivially. Since \(\operatorname{ord}_7(3)=6\), its action on \(P\) decomposes into irreducible constituents of dimensions \(6,1\). Its centralizer in \(\operatorname{GL}(P)\) is isomorphic to \(\mathbb F_{3^6}^{\,*}\times\mathbb F_3^{\,*}\); every \(2\)-subgroup of this centralizer is generated by at most two elements. Hence \(a\) does not centralize \(D\).

The group \(D\) is normal in the image of \(G\) on \(P\). Over an algebraic closure of characteristic \(3\), decompose \(P\) into weight spaces for the abelian \(2\)-group \(D\). Since the action of \(D\) is faithful and \(a\) does not centralize it, some weight occurring has an orbit of length \(7\) under \(a\). The dimension is \(7\), so this orbit consists of exactly seven distinct weights, each occurring with multiplicity \(1\).

The Frobenius automorphism \(\lambda\mapsto\lambda^3\) permutes these seven weights and commutes with the seven-cycle induced by \(a\). The centralizer of a seven-cycle in \(S_7\) has order \(7\). But on the characters of a \(2\)-group, Frobenius has \(2\)-power order. Hence it fixes each weight. All weight spaces are one-dimensional and defined over \(\mathbb F_3\). The character values lie in \(\mathbb F_3^{\,*}=\{1,-1\}\), so \(D\) is elementary abelian. The resulting decomposition is permuted by all of \(G\), since \(D\) is normal in the image of \(G\).

Denote by \(\pi\) the permutation action on the seven one-dimensional spaces. The group \(\pi(L)\) is transitive, since it contains the seven-cycle \(\pi(a)\). The group \(\pi(K)\) is normal in \(\pi(L)\) and is a \(5\)-group: the image of \(H\) fixes all weight spaces, and \(K/H\) is a \(5\)-group. The orbits of a normal \(5\)-subgroup in a transitive group have equal sizes; this size is a power of \(5\) dividing \(7\), and hence equals \(1\). Thus \(K\) fixes each of the seven lines. A Sylow \(5\)-subgroup of \(K\) acts trivially on each line over \(\mathbb F_3\), and therefore centralizes \(P\).

For \(|x|=9\), Theorem~\ref{thm:L27-nine-simple-part-nonabelian} gives \(x\in R(G)\), and Theorem~\ref{thm:L27-normal-sylow-three-structure} gives \(x^3\in P\). Since \(P\) is abelian, the action of \(x\) on \(P\) has order dividing \(3\). The relations \([R(G),L]\le K\) and \(\pi(K)=1\) imply that \(\pi(x)\) commutes with the seven-cycle \(\pi(a)\). Thus \(\pi(x)\) lies in the cyclic centralizer of order \(7\), while its order also divides \(3\). Hence \(\pi(x)=1\). On each of the seven \(\mathbb F_3\)-lines, the element \(x\) acts as a scalar whose order divides both \(2\) and \(3\); this scalar is \(1\). Therefore \(C_P(x)=P\).
\end{proof}

\clearpage
\subsection{The position of elements of order \(8\) in the seven-dimensional subcase}

\Needspace{11\baselineskip}
\begin{theorem}\label{thm:L27-rank-seven-eight-outside-L}
Under the assumptions of Theorem~\ref{thm:L27-rank-seven-abelian-action}, put \(L=G^\infty\), \(R=R(G)\), and \(T=R/K\). Then:
\begin{enumerate}
\item the action of \(R\) on the seven lines in that theorem is trivial, and its image in \(\operatorname{GL}(P)\) is elementary abelian; in particular, every odd-order element of \(R\) centralizes \(P\);
\item \(8\notin\omega(L)\). The image of every element of order \(8\) in \(T\) is an involution;
\item \(2,3\mid |T|\) and \(\{12,14,21\}\subseteq\omega(G/K)\). Moreover, the image of \(R\) in \(\operatorname{GL}(P)\) contains \(-I\).
\end{enumerate}
\end{theorem}
\begin{proof}
Denote the permutation action on the seven lines by \(\pi\). By Theorem~\ref{thm:L27-rank-seven-abelian-action}, we have \(\pi(K)=1\), and the image \(\pi(L)\) is transitive and contains a seven-cycle. Since \([R,L]\le R\cap L=K\), the groups \(\pi(R)\) and \(\pi(L)\) centralize one another. In particular, \(\pi(R)\) lies in the centralizer of a seven-cycle in \(S_7\), which has order \(7\). On the other hand, \(7\nmid|R|\): we have already proved \(7\nmid|K|\), and \(\omega(R/K)\subseteq\{1,2,3\}\). Hence \(\pi(R)=1\).

Thus every element of \(R\) acts on each of the seven \(\mathbb F_3\)-lines as a scalar \(\pm1\). The image of \(R\) is contained in a direct product of seven groups of order \(2\), proving the first assertion. Moreover, \(G=RL\), since \(L\) maps onto \(G/R\cong S\); hence \(\pi(G)=\pi(L)\).

Let \(y\in G\) have order \(8\). Its action on \(P\) has no nonzero fixed vectors: otherwise the forbidden order \(24\) would arise. The matrix of \(y\) is a signed permutation matrix relative to the decomposition into seven lines. For a cycle of length \(k\), the product of the signs along the cycle must be \(-1\): if the product were \(1\), there would be a nonzero fixed vector in the corresponding sum of lines. The determinant of the block for such a cycle is \((-1)^{k-1}(-1)=(-1)^k\). Multiplying over the cycles gives
\[
 \det(y|_P)=(-1)^7=-1.
\]
The group \(L\) is perfect, so its image under the homomorphism \(\det:G\to\mathbb F_3^{\,*}\) is trivial. Thus \(y\notin L\), proving \(8\notin\omega(L)\). Under the identification \(G/K=S\times T\), the image of \(L\) is \(S\times1\), and the image of \(y\) is a \(2\)-element. Since \(\omega(T)\subseteq\{1,2,3\}\), its projection to \(T\) is a nontrivial involution.

An element \(x\) of order \(9\) lies in \(R\setminus K\) by Theorems~\ref{thm:L27-K-no-nine} and \ref{thm:L27-nine-simple-part-nonabelian}, so its image in \(T\) has order \(3\). In \(S\times T\), the independent factors give elements of orders \(7\cdot2=14\), \(7\cdot3=21\), and \(4\cdot3=12\).

It remains to obtain \(-I\). Write \(y=rl\), where \(r\in R\) and \(l\in L\). The equalities \(\det(y|_P)=-1\) and \(\det(l|_P)=1\) imply that the diagonal matrix \(r|_P\) has an odd number of signs \(-1\). Let \(a\in L\) permute the seven lines as a seven-cycle. The product of the seven conjugate matrices \((r|_P)^{a^i}\), \(0\le i<7\), is \(-I\): on each line, the product of the signs equals the determinant of \(r|_P\). Since \(R\) is normal in \(G\), this product lies in the image of \(R\) on \(P\).
\end{proof}

\clearpage
\subsection{Normal \(3\)-subgroups and elements of order \(9\)}

We retain the assumptions of the remaining branch \(S=L_2(7)\). In this subsection, neither normality nor abelianness of a Sylow \(3\)-subgroup of \(H\) is assumed unless explicitly stated.

\Needspace{11\baselineskip}
\begin{theorem}\label{thm:L27-normal-three-all-exponent-three}
Every normal \(3\)-subgroup \(N\triangleleft G\) has exponent dividing \(3\), is \(2\)-Engel, and has nilpotency class at most \(3\). In particular,
\[
 \exp O_3(G)\mid3,\qquad O_3(G)=O_3(R(G)),
 \qquad O_3(G)\cap K=O_3(K).
\]
No element of order \(9\) lies in \(O_3(G)\) or in the Fitting subgroup \(F(G)=F(R(G))\).
\end{theorem}
\begin{proof}
By Lemma~\ref{lem:L27-global-five-complement}, the group \(K\) contains \(A\cong C_5^6\). Since \(N\) is normal in \(G\), the group \(A\) acts on \(N\). For every \(1\ne a\in A\), the group \(C_N(a)\) has no elements of order \(9\): together with \(a\), such an element would give the forbidden order \(45\). Order \(27\) is also absent from \(G\), so \(\exp C_N(a)\mid3\). Lemma~\ref{lem:five-rank-exponent-three} applies directly and gives all the assertions about \(N\).

The subgroup \(O_3(R(G))\) is characteristic in the radical and normal in \(G\); conversely, \(O_3(G)\le R(G)\). Hence these groups coincide. The group \(O_3(K)\) is characteristic in \(K\), normal in \(G\), and contained in \(O_3(G)\cap K\). The reverse inclusion follows from normality of \(O_3(G)\cap K\) in \(K\), proving the second equality.

In the nilpotent group \(F(G)=F(R(G))\), every element of order \(9\) would lie in its Sylow \(3\)-subgroup, which is contained in \(O_3(G)\). The preceding argument excludes such an element.
\end{proof}

\Needspace{11\baselineskip}
\begin{corollary}\label{cor:L27-abelian-three-radical-nonnormal}
If a Sylow \(3\)-subgroup of \(H\) is abelian, then all elements of order \(9\) in \(G\) lie in
\[
 R(G)\setminus F(R(G)).
\]
Moreover, a Sylow \(3\)-subgroup of \(R(G)\) is not normal in \(R(G)\).
\end{corollary}
\begin{proof}
By Theorem~\ref{thm:L27-nine-simple-part-nonabelian}, all elements of order \(9\) lie in \(R(G)\), and the preceding theorem excludes them from \(F(R(G))\). If a Sylow \(3\)-subgroup of \(R(G)\) were normal, it would equal \(O_3(R(G))\) and contain every \(3\)-element of the radical. In particular, it would contain an element of order \(9\), which is impossible.
\end{proof}

\Needspace{11\baselineskip}
\begin{corollary}\label{cor:L27-normal-three-radical-quotient}
If a Sylow \(3\)-subgroup \(P\) of \(H\) is normal in \(H\), then
\[
 9\notin\omega\bigl(G/O_3(G)\bigr),\qquad
 9\notin\omega\bigl(R(G)/O_3(R(G))\bigr).
\]
\end{corollary}
\begin{proof}
By Theorem~\ref{thm:L27-normal-sylow-three-structure}, we have \(9\notin\omega(G/P)\), where \(P=O_3(K)\le O_3(G)\). The first group displayed is a quotient of \(G/P\), and the second is a subgroup of that quotient, since \(O_3(G)=O_3(R(G))\). The spectra of subgroups and quotient groups are contained in the spectra of the original finite groups.
\end{proof}

\clearpage
\subsection{Excluding rank \(7\) for a normal abelian \(3\)-subgroup}

Recall that \(K=R(G)\cap L\) in this branch by Lemma~\ref{lem:L27-final-ordinary-radical-structure}.

We retain the notation of the preceding subsections: \(S=L_2(7)\), \(L=G^\infty\), \(R=R(G)\), \(K=R\cap L\), and \(H\) is the characteristic normal Hall \(5'\)-subgroup of \(K\).

\Needspace{11\baselineskip}
\begin{lemma}\label{lem:L27-central-cube-obstruction}
Suppose that a Sylow \(3\)-subgroup \(P\) of \(H\) is normal and abelian, and that an element \(x\in G\) of order \(9\) centralizes \(P\). Then
\[
 1\ne x^3\in C_P(L).
\]
\end{lemma}
\begin{proof}
The group \(P\) is elementary abelian and normal in \(G\), and \(x^3\in P\), by Theorems~\ref{thm:L27-K-no-nine} and \ref{thm:L27-normal-sylow-three-structure}. Put \(C=C_K(P)\). This is a normal subgroup of \(G\) in which \(P\) is central and is a Sylow \(3\)-subgroup. By the Schur--Zassenhaus theorem, \(C\) has a Hall \(3'\)-subgroup \(W\) complementing \(P\), and centrality of \(P\) gives \(C=P\times W\). All elements of \(C\) of order coprime to \(3\) lie in \(W\); hence \(W\) is characteristic in \(C\) and normal in \(G\).

By Theorem~\ref{thm:L27-nine-simple-part-nonabelian}, the element \(x\) lies in \(R\). For every \(l\in L\), the commutator \([x,l]\) lies in \(R\cap L=K\). Moreover, \(C_G(P)\triangleleft G\), so \([x,l]\in C_K(P)=C\). In the quotient \(G/W\), we have \([xW,lW]\in PW/W\). The element \(xW\) centralizes \(PW/W\), and this group has exponent \(3\). Thus \((xW)^3\) is centralized by \(lW\). Consequently, \([x^3,l]\in W\). But \(x^3\in P\triangleleft G\), so \([x^3,l]\in P\cap W=1\). Hence \(x^3\in C_P(L)\), and it is nontrivial because \(|x|=9\).
\end{proof}

\Needspace{11\baselineskip}
\begin{lemma}\label{lem:L27-rank-seven-nonabelian-image}
Suppose that a Sylow \(3\)-subgroup \(P\) of \(H\) is normal and abelian, \(\dim_{\mathbb F_3}P=7\), and \(D=H/C_H(P)\) is nonabelian. Then a contradiction arises.
\end{lemma}
\begin{proof}
As in the proof of Theorem~\ref{thm:L27-rank-seven-abelian-action}, the group \(D\) is a \(2\)-group with \(d(D)\ge3\), acting faithfully on \(P\). There is an element \(a\in L\) of order \(7\) whose action on \(P\) has irreducible components of dimensions \(6,1\). Its centralizer in \(\operatorname{GL}(P)\) is isomorphic to \(\mathbb F_{3^6}^{\,*}\times\mathbb F_3^{\,*}\), and every \(2\)-subgroup of this centralizer is generated by at most two elements. Hence \(a\) acts nontrivially on \(D\), which it normalizes.

Take the characteristic subgroup \(C\le D\) from Lemma~\ref{lem:critical-two-external}: it has class at most \(2\), \(C/Z(C)\) is elementary abelian, and \(a\) acts nontrivially on it. Since \(C\) is characteristic in \(D\), the image of \(D\) permutes the isotypic components of the \(\mathbb F_3C\)-module \(P\), and \(a\) permutes them as well.

If \(a\) has an orbit of length \(7\) on the types, then all seven corresponding isotypic components are one-dimensional and exhaust \(P\). The image of \(D\) is a normal \(2\)-subgroup of the transitive permutation group of degree \(7\) generated by this image and the permutation induced by \(a\). Its orbits have the same size, which is a power of \(2\) dividing \(7\); hence it fixes each of the seven lines. Then the faithful group \(D\) is diagonal and abelian, a contradiction.

Thus \(a\) fixes every type, and all isotypic components are \(a\)-invariant \(\mathbb F_3\)-subspaces. There cannot be only one such component. Indeed, the dimension of any irreducible representation of a finite \(2\)-group over \(\mathbb F_3\) is a power of \(2\): over the algebraic closure its degree is a power of \(2\), while the required field extension is generated by roots of unity of \(2\)-power order and has \(2\)-power degree. If a single isotypic component exhausted the space of dimension \(7\), its irreducible constituents would be one-dimensional, and the faithful group \(C\) would act by scalars. Then \(a\) would centralize \(C\), a contradiction.

The only nonzero proper \(a\)-invariant \(\mathbb F_3\)-subspaces of \(P\) have dimensions \(1,6\), since the restriction to \(\langle a\rangle\) has type \(1\oplus6\). Hence there are exactly two isotypic components, of dimensions \(1,6\). Denote the six-dimensional component by \(V\), and the image of \(C\) on it by \(C_6\). On the one-dimensional component, \(a\) acts trivially. If it centralized \(C_6\), it would also centralize the image of \(C\) on the one-dimensional component, and hence the faithful group \(C\) on all of \(P\). Thus the action of \(a\) on \(C_6\) is nontrivial.

The action of \(a\) on \(V\) is irreducible. Every \(a\)-invariant abelian subgroup \(B\le C_6\) is centralized by \(a\): over the algebraic closure, its weights on \(V\) cannot have an orbit of length \(7\), and they separate the elements of \(B\). Thus \(B\) lies in the cyclic field group \(C_{\operatorname{GL}(V)}(a)\cong\mathbb F_{3^6}^{\,*}\). An element of order \(4\) in \(B\), together with \(a\), would give order \(28\) in the image of a subgroup of \(G\) on \(V\); this is forbidden. Hence \(|B|\le2\), and in particular \(|Z(C_6)|=2\).

The group \(C_6\) is nonabelian: otherwise \(C_6=Z(C_6)\cong C_2\), and \(a\) would act trivially on it. It is a quotient of \(C\), has class at most \(2\), and \(C_6/Z(C_6)\) is elementary abelian. Consequently, \(C_6\) is extraspecial of order \(2^{1+2m}\). The automorphism \(a\) is nontrivial on \(C_6/Z(C_6)\), since the kernel of the action on this quotient is a \(2\)-group. The divisibility \(7\mid|\operatorname{Sp}_{2m}(2)|\) implies \(m\ge3\), as in Theorem~\ref{thm:L27-normal-abelian-three-rank-seven}. The central involution of \(C_6\) acts on \(V\) as \(-I\), since it centralizes \(a\) and is the unique involution in \(\mathbb F_{3^6}^{\,*}\). The same argument using the weights of the preimage of a maximal isotropic subspace gives \(\dim V\ge2^m\ge8\), contradicting \(\dim V=6\).
\end{proof}

\Needspace{11\baselineskip}
\begin{theorem}\label{thm:L27-normal-abelian-three-rank-eight}
If a Sylow \(3\)-subgroup \(P\) of \(H\) is normal and abelian, then
\[
 \dim_{\mathbb F_3}P\ge8,\qquad |P|\ge3^8.
\]
\end{theorem}
\begin{proof}
By Theorem~\ref{thm:L27-normal-abelian-three-rank-seven}, the rank of \(P\) is at least \(7\). Suppose that it is \(7\), and put \(D=H/C_H(P)\). If \(D\) is nonabelian, the preceding lemma gives a contradiction.

If \(D\) is abelian, Theorem~\ref{thm:L27-rank-seven-abelian-action} gives a decomposition of \(P\) into seven distinct one-dimensional weight spaces for \(D\) forming a single orbit of an element of order \(7\), and also \(C_P(x)=P\) for every \(|x|=9\). None of these seven weights is trivial: the trivial weight would be fixed by the element of order \(7\) and could not have an orbit of length \(7\). Hence \(C_P(H)=C_P(D)=0\), and therefore \(C_P(L)=0\), since \(H\le L\).

But Lemma~\ref{lem:L27-central-cube-obstruction} gives, for an element \(x\) of order \(9\), a nontrivial element \(x^3\in C_P(L)\). This contradiction excludes rank \(7\) and completes the proof.
\end{proof}

\clearpage
\subsection{Excluding a normal abelian Sylow \(3\)-subgroup completely}

The following argument excludes the case in all ranks at once. It requires neither a bound on the module dimension nor the use of critical \(2\)-subgroups or modular character tables.

\Needspace{11\baselineskip}
\begin{lemma}\label{lem:seven-three-quadratic-obstruction}
Let \(B=\langle a,b\rangle\cong C_7\rtimes C_3\), where \(|a|=7\), \(|b|=3\), and the action of \(b\) on \(\langle a\rangle\) is nontrivial. Let \(V\) be a module over a field of characteristic \(3\). If \((b-1)^2V=0\), then \(a\) acts trivially on \(V\). In particular, if \(9\notin\omega(V\rtimes B)\), then \(a\) centralizes \(V\).
\end{lemma}
\begin{proof}
Extend the field to an algebraic closure and decompose \(V\) into eigenspaces of \(a\). If a nontrivial root \(\lambda\) of unity of order \(7\) occurs, then \(b\) cyclically permutes three distinct spaces with eigenvalues \(\lambda,\lambda^2,\lambda^4\) (possibly in the reverse order). For a nonzero vector \(v\) in the first space, the vectors \(v,bv,b^2v\) are linearly independent. Thus
\[
 (b-1)^2v=(1+b+b^2)v\ne0,
\]
contrary to the hypothesis. Hence only the eigenvalue \(1\) occurs, and semisimplicity of the action of \(a\) gives \(a=1\) on \(V\).

For the final assertion, write \(V\) additively. In the semidirect product, for every \(v\in V\),
\[
 (vb)^3=v+bv+b^2v.
\]
If the right-hand side is nonzero, the element \(vb\) has order \(9\). In the absence of that order, the operator \(1+b+b^2=(b-1)^2\) is zero, and the part already proved applies.
\end{proof}

\Needspace{11\baselineskip}
\begin{theorem}\label{thm:L27-perfect-nine-free-abelian-kernel}
Let \(X\) be a finite perfect group with \(9\notin\omega(X)\), and let \(J\triangleleft X\) be solvable with \(X/J\cong L_2(7)\). If \(V\triangleleft X\) is an elementary abelian \(3\)-subgroup, \(V\le J\), and \((|J/V|,21)=1\), then \(V=1\).
\end{theorem}
\begin{proof}
First we prove that \(X\) splits over \(V\). A Sylow \(3\)-subgroup of the quotient \(X/V\) has order \(3\), since \((|J/V|,21)=1\) and \(|L_2(7)|_3=3\). Hence, for a Sylow \(3\)-subgroup \(U\) of \(X\), we have \(|U:V|=3\). Every nonidentity element of \(U\) has order \(3\): otherwise its cyclic subgroup would contain an element of order \(9\). Choosing \(u\in U\setminus V\), we obtain
\[
 U=V\rtimes\langle u\rangle.
\]
By Gasch\"utz's theorem on splitting over a normal abelian subgroup, this splitting in a Sylow subgroup extends to the whole group \cite[Theorem~1.1]{GaschutzComplement}. Thus
\[
 X=V\rtimes E,
 \qquad E\cong X/V.
\]
It is the criterion for the existence of a complement that is used here; complete reducibility of the \(\mathbb F_3E\)-module \(V\) is not asserted. The group \(E\) is perfect, being isomorphic to a quotient of a perfect group. Put \(A=J\cap E\). Then \(A\triangleleft E\) is solvable, \(3\nmid|A|\), and \(E/A\cong L_2(7)\).

The group \(L_2(7)\cong\operatorname{GL}_3(2)\) has a subgroup \(7{:}3\): on the space \(\mathbb F_8\), multiplication by nonzero elements gives \(C_7\), while the automorphism \(t\mapsto t^2\) has order \(3\) and normalizes it nontrivially. Let \(M\le E\) be the full preimage of this subgroup. We have \((|A|,21)=1\), so \(A\) is a normal Hall subgroup of \(M\). The Schur--Zassenhaus theorem gives a complement
\[
 B=\langle a,b\rangle\cong7{:}3
 \quad\text{to }A\text{ in }M.
\]
The subgroup \(VB=V\rtimes B\) lies in \(X\), so it has no element of order \(9\). By the preceding lemma, \(a\) centralizes \(V\).

Let \(F\) be the image of \(E\) on \(V\), and \(A_0\) the image of \(A\). The group \(F\) is perfect, and \(A_0\triangleleft F\) is solvable. The kernel of the action of \(E\) on \(V\) contains \(a\). The image of its normal closure in \(E/A\cong L_2(7)\) is the whole simple group, since the image of \(a\) is nontrivial. Thus \(F=A_0\). A perfect solvable group is trivial, giving \([V,E]=1\).

We have obtained the direct product \(X=V\times E\). Its derived subgroup is contained in \(E\), whereas \(X=X'\). Hence \(V=1\), as required.
\end{proof}

\Needspace{11\baselineskip}
\begin{corollary}\label{cor:L27-normal-abelian-three-excluded}
In the remaining branch \(S=L_2(7)\), a Sylow \(3\)-subgroup \(P\) of \(H\) cannot be both normal in \(H\) and abelian. This exclusion is independent of the rank of \(P\).
\end{corollary}
\begin{proof}
Under these assumptions, \(P\) is elementary abelian and normal in \(G\), by Theorem~\ref{thm:L27-normal-sylow-three-structure}. It is a Sylow \(3\)-subgroup of \(K\), since \(K/H\) is a \(5\)-group. Thus \(\pi(K/P)\subseteq\{2,5\}\). Moreover, \(P\ne1\), since \(3\in\pi(H)\).

By Theorem~\ref{thm:L27-nine-simple-part-nonabelian}, all elements of order \(9\) in \(G\) lie in \(R=R(G)\). If such an element lay in \(L\), it would lie in \(R\cap L=K\), contrary to Theorem~\ref{thm:L27-K-no-nine}. Hence \(9\notin\omega(L)\). Apply the preceding theorem with
\[
 X=L,\qquad J=K,\qquad V=P.
\]
All its assumptions hold: \(L\) is perfect, \(K\) is solvable, \(L/K\cong L_2(7)\), and \(K/P\) has order coprime to \(21\). We obtain \(P=1\), a contradiction.
\end{proof}

There are now two remaining possibilities for the Sylow \(3\)-subgroups of \(H\): they are either nonnormal, or normal and nonabelian. The next subsection excludes both branches simultaneously using a suitable normal quotient of \(H\), completing the treatment of the socle \(L_2(7)\).

\clearpage
\subsection{Excluding the socle \(L_2(7)\) completely with an arbitrary radical}

The final step is to choose a suitable normal quotient of the solvable group \(H\). Normality of its original Sylow \(3\)-subgroup is unnecessary: in the chosen quotient, a Sylow \(3\)-subgroup will be normal and elementary abelian. The spectrum of the quotient may be strictly smaller than \(\Om\); we use only the element-order exclusions that persist and the previously established structure of \(G/H\).

\Needspace{11\baselineskip}
\begin{lemma}\label{lem:soluble-normal-p-top-section}
Let \(A\triangleleft Y\) be a solvable subgroup of a finite group, with \(p\mid|A|\). There exists \(N\triangleleft Y\), \(N\le A\), such that \(A/N\) has a nontrivial normal elementary abelian Sylow \(p\)-subgroup \(V\), with \(V\triangleleft Y/N\).
\end{lemma}
\begin{proof}
Among the subgroups \(N\le A\) normal in \(Y\) for which \(p\mid|A/N|\), choose one maximal by inclusion. Such a choice is possible, since \(N=1\) is admissible. Put \(\overline Y=Y/N\), \(\overline A=A/N\), and choose a minimal nontrivial normal subgroup \(V\) of \(\overline Y\) contained in \(\overline A\).

Since \(\overline A\) is solvable, \(V\) is an elementary abelian \(r\)-group for some prime. Indeed, solvability and minimality first give \(V'=1\), and then characteristic Sylow subgroups and \(\Omega_1(V)\) give elementary abelianness. The preimage of \(V\) in \(Y\) is strictly larger than \(N\). By maximality of \(N\), we have \(p\nmid|\overline A/V|\). But \(p\mid|\overline A|\), so \(r=p\), and \(V\) is a Sylow \(p\)-subgroup of \(\overline A\).
\end{proof}

\Needspace{11\baselineskip}
\begin{lemma}\label{lem:L27-no-nine-after-elementary-three-section}
In the remaining branch with \(S=L_2(7)\) and \(L/K\cong S\), suppose that \(N\triangleleft G\), \(N\le H\), and \(H/N\) has a normal elementary abelian Sylow \(3\)-subgroup \(V\). Then
\[
 9\notin\omega(L/N).
\]
\end{lemma}
\begin{proof}
Write \(\overline G=G/N\), \(\overline H=H/N\), \(\overline K=K/N\), and \(\overline L=L/N\). These quotients are well defined because \(N\le H\le K\le L\). By Theorem~\ref{thm:L27-K-no-nine}, we have \(9\notin\omega(\overline K)\). Moreover,
\[
 \overline G/\overline H=G/H=Q\rtimes(S\times T),
 \qquad \overline L/\overline H=L/H=Q\rtimes S.
\]
Here \(Q=K/H\ne1\), \(\exp Q=5\), and \(Z=Z(Q)\) is an elementary abelian group of order \(5^{6n}\), \(n\ge1\), by Lemma~\ref{lem:L27-global-five-complement}. Elements of order \(7\) in \(S\) act on \(Z\) without nontrivial fixed points: otherwise order \(35\) would occur in \(G/H\), and hence in \(G\).

Suppose that \(x\in\overline L\) has order \(9\). Its image in \(\overline L/\overline K\cong S\) has order \(3\), since \(9\notin\omega(\overline K)\). The image \(x\overline H\) also has order \(3\), since \(9\notin\omega(G/H)\). In the decomposition \(Q\rtimes S\), write it as \(qa\), where \(q\in Q\), \(a\in S\), and \(|a|=3\). It acts on \(Z\) as \(a\), since \(Q\) centralizes its own center.

In \(S\cong\operatorname{GL}_3(2)\), every element of order \(3\) normalizes some subgroup of order \(7\), acting on it nontrivially. Indeed, a subgroup \(7{:}3\) was constructed in the proof of Theorem~\ref{thm:L27-perfect-nine-free-abelian-kernel}, and all elements of order \(3\) in \(\operatorname{GL}_3(2)\) are conjugate: their rational canonical form has blocks \(t-1\) and \(t^2+t+1\). Over an algebraic closure of characteristic \(5\), decompose \(Z\) into weight spaces for this \(C_7\). The trivial weight is absent, and \(a\) permutes the nontrivial weights in triples. On the direct sum of each triple of weight spaces, the operator \(a\) cyclically permutes spaces of equal dimension, and its cube is the identity. Hence the dimension of the fixed space on a triple is one third of its total dimension. Thus
\[
 \dim_{\mathbb F_5}C_Z(x\overline H)
 =\frac13\dim_{\mathbb F_5}Z=2n\ge2.
\]

Choose \(A\cong C_5^2\) in this centralizer. The subgroup \(\langle x\overline H\rangle A\) is the direct product \(C_3\times C_5^2\). Its full preimage \(F\le\overline G\) is solvable, since \(\overline H\) is solvable. Choose a Hall \(\{3,5\}\)-subgroup \(E\) of \(F\) containing \(x\). This is possible by Hall's theorem on embedding a \(\{3,5\}\)-subgroup in a Hall subgroup of a solvable group. Since \(\pi(\overline H)\subseteq\{2,3\}\) and its Sylow \(3\)-subgroup \(V\) is normal, we have
\[
 E\cap\overline H=V,\qquad
 E/V\cong C_3\times C_5^2.
\]
The subgroup \(V\) is normal in \(F\): it is characteristic in \(\overline H\) as its unique Sylow \(3\)-subgroup. In particular, \(V\triangleleft E\).

The spectrum of \(\overline G\) is contained in \(\omega(G)=\Om\), so \(45\notin\omega(E)\). Lemma~\ref{lem:elementary-three-five-nine} then gives \(9\notin\omega(E)\), contradicting \(x\in E\), \(|x|=9\). This proves the lemma.
\end{proof}

\Needspace{11\baselineskip}
\begin{theorem}\label{thm:L27-arbitrary-radical}
There is no finite group \(G\) with spectrum \(\omega(G)=\omega(S_{10})\) such that
\[
 \operatorname{Soc}(G/R(G))\cong L_2(7).
\]
No restriction is placed on the solvable radical, and the group is not assumed to be critical.
\end{theorem}
\begin{proof}
Suppose otherwise. Use the notation \(R=R(G)\), \(L=G^\infty\), and \(K=[R\cap L,L]\). The double cover was excluded by Theorem~\ref{thm:L27-binary-cover-excluded}, and the outer extension by Theorem~\ref{thm:L27-outer-extension-excluded}. By Lemma~\ref{lem:L27-final-ordinary-radical-structure}, the remaining case is
\[
 K=R\cap L,\qquad L/K\cong L_2(7),
 \qquad G/K\cong L_2(7)\times T.
\]
Lemma~\ref{lem:L27-global-five-complement} gives a characteristic normal Hall \(5'\)-subgroup \(H\le K\), \(\pi(H)\subseteq\{2,3\}\), \(3\mid|H|\), and all the properties of \(Q=K/H\) used in the preceding lemma. By Theorem~\ref{thm:L27-K-no-nine}, order \(9\) is absent from \(K\).

Apply Lemma~\ref{lem:soluble-normal-p-top-section} with \(Y=G\), \(A=H\), and \(p=3\). We obtain a subgroup \(N\le H\) normal in \(G\), and a nontrivial elementary abelian Sylow \(3\)-subgroup \(V\triangleleft G/N\) of \(H/N\). By Lemma~\ref{lem:L27-no-nine-after-elementary-three-section}, we have \(9\notin\omega(L/N)\).

Now apply Theorem~\ref{thm:L27-perfect-nine-free-abelian-kernel} to the groups
\[
 X=L/N,\qquad J=K/N,
 \qquad 1\ne V\le H/N\le J.
\]
The group \(X\) is perfect as a quotient of the perfect group \(L\). The group \(J\) is solvable and normal in \(X\), and \(X/J\cong L/K\cong L_2(7)\). The subgroup \(V\) is normal in \(X\) and is an elementary abelian \(3\)-group. Finally, \((H/N)/V\) is a \(2\)-group, and \((K/N)/(H/N)\cong K/H\) is a \(5\)-group. Thus
\[
 \pi(J/V)\subseteq\{2,5\},\qquad (|J/V|,21)=1.
\]
All the assumptions of that theorem hold. It gives \(V=1\), contrary to the choice \(V\ne1\). This excludes the remaining branch and, together with the exclusion of the double-cover and outer branches, completes the treatment of the socle \(L_2(7)\).
\end{proof}

In particular, normal abelian, normal nonabelian, and nonnormal Sylow \(3\)-subgroups of \(H\) have all been excluded simultaneously. The final step uses a quotient by a subgroup that is maximal among the \(G\)-normal subgroups of \(H\) for which the quotient order is divisible by \(3\). This avoids further analysis of the structure and rank of the original Sylow subgroup.

\clearpage

\section{\texorpdfstring{The socle \(L_3(4)\)}{The socle L3(4)}}
\subsection{Excluding both odd branches of \(K/H\) for \(L_3(4)\)}

Here \(S=L_3(4)\), \(L=G^\infty\), \(P=R(G)\cap L\),
\(K=[P,L]=[K,L]\). Unlike the notation for a Sylow subgroup in the preceding sections, \(P\) again denotes the intersection of the radical with the perfect residual.
The subgroup \(H\le K\) is taken from
Lemma~\ref{lem:L34-global-two-three-complement}.
We exclude both possibilities \(K/H\ne1\) simultaneously,
without restrictions on the ranks or classes of the Sylow subgroups.

\Needspace{11\baselineskip}
\begin{lemma}\label{lem:L34-admissible-covers-subgroups}
The possible isomorphism types of \(D=L/K\) are
\[
 S,\qquad 2.S,\qquad 2^2.S,\qquad 3.S.
\]
Each of these groups contains a subgroup \(S_0\cong L_2(7)\).
For suitable subgroups \(L_2(7),A_6<S\), their full preimages
in these covers are as follows:
\[
\begin{array}{c|cc}
 D&L_2(7)&A_6\\\hline
 S&L_2(7)&A_6\\
 2.S&C_2\times L_2(7)&C_2\times A_6\\
 2^2.S&C_2^2\times L_2(7)&C_2^2\times A_6\\
 3.S&C_3\times L_2(7)&3.A_6.
\end{array}
\]
\end{lemma}
\begin{proof}
The group \(D\) is perfect, \(D/Z(D)\cong S\), and
\(Z(D)=P/K\): centrality follows from \(K=[P,L]\),
and the image of the center in the simple group \(S\) is trivial.
The multiplier of \(S\) is \(C_3\times C_4\times C_4\).
If \(Z(D)\) contains an element of order \(4\), this element,
together with an element of order \(7\), gives the forbidden order \(28\).
If the center has both a nontrivial \(2\)-part and a nontrivial
\(3\)-part, their product has order \(6\), which together
with an element of order \(7\) gives \(42\).
Thus the center is \(1,C_2,C_2^2\), or \(C_3\).
Up to isomorphism, the corresponding perfect covers are precisely
the four groups listed above \cite{L34CoverSubgroupData}.

The two subgroups of the simple group \(S\) and their preimages
in the covers are given in the ATLAS and CTblLib tables of maximal
subgroups \cite{L34CoverSubgroupData}. In particular, in the double
cover these are the direct products \(C_2\times L_2(7)\) and
\(C_2\times A_6\), and in \(2^2.S\) they are the corresponding
products with \(C_2^2\). In the triple cover, the preimage
of \(L_2(7)\) also splits because
\(M(L_2(7))=C_2\); the preimage of \(A_6\) is \(3.A_6\).
This proves the table and the existence of \(S_0\) in every case.
\end{proof}

\Needspace{11\baselineskip}
\begin{lemma}\label{lem:A6-triple-seven-five-fixed}
Every element of order \(5\) has a nonzero fixed vector
on every nonzero module of characteristic \(7\)
for \(A_6\) or \(3.A_6\).
\end{lemma}
\begin{proof}
Since \(7\nmid|3.A_6|=1080\), Maschke's theorem allows us
to extend the field to an algebraic closure and consider
absolutely simple modules. Their degrees and the multiplicities
of the eigenvalue \(1\) are calculated from the ordinary character table.
Let \(c\) have order \(5\). Two of its nonidentity powers
belong to class \(5A\), and two to \(5B\), so
\[
 \dim C_V(c)=\frac{\chi(1)+2\chi(5A)+2\chi(5B)}5.
\]
The complete tables of \(A_6,3.A_6\)
\cite{A6TripleFiveFixedData} give
\[
\begin{array}{c|rrrrrrr}
 A_6:\ \chi(1)&1&5&5&8&8&9&10\\
 \chi(5A)+\chi(5B)&2&0&0&1&1&-2&0\\
 \dim C_V(c)&1&1&1&2&2&1&2
\end{array}
\]
and, for the types of \(3.A_6\) that are nontrivial on the center,
\[
\begin{array}{c|rrrrr}
 \chi(1)&3&3&6&9&15\\
 \chi(5A)+\chi(5B)&1&1&2&-2&0\\
 \dim C_V(c)&1&1&2&1&3.
\end{array}
\]
The types listed in the second table occur for each of the two nontrivial
central characters. The types with trivial central action
have already been included in the table for \(A_6\).
For example, on the three-dimensional types the sum of the two
values is \(-\zeta-\zeta^2-\zeta^3-\zeta^4=1\), where \(\zeta^5=1\ne\zeta\).
On the six-dimensional type the values are \(1,1\), on the
nine-dimensional type they are \(-1,-1\), and on the
fifteen-dimensional type they are \(0,0\).
Every calculated dimension is positive. A field extension
does not change the dimension of the fixed-point space,
which proves the assertion over the original field.
\end{proof}

\Needspace{11\baselineskip}
\begin{theorem}\label{thm:L34-K-two-three-only}
If \(\omega(G)=\Om\) and
\(\operatorname{Soc}(G/R(G))\cong L_3(4)\), then
\[
 5\nmid|K|,\qquad 7\nmid|K|,
 \qquad K=H,\qquad \pi(K)\subseteq\{2,3\}.
\]
Thus both nontrivial branches of \(K/H\) in
Lemma~\ref{lem:L34-global-two-three-complement} are completely excluded.
\end{theorem}
\begin{proof}
\emph{First exclude \(7\mid|K|\).}
By Lemma~\ref{lem:L34-global-two-three-complement}, in this case
\(Q=K/H\ne1\) is a \(7\)-group of exponent \(7\).
Its center \(Z=Z(Q)\ne1\) is an elementary abelian
\(7\)-group. The action of \(L/H\) on \(Z\) factors
through \(D=L/K\), since \(Q\) centralizes its own center.
Every element \(c\in D\) of order \(5\) acts on \(Z\)
without nonzero fixed vectors: it lifts to an element of
order \(5\) in \(L/H\), and a nonzero fixed vector in characteristic \(7\) would
give order \(35\) in this quotient, and hence in \(G\).

If \(D=S,2.S\), or \(2^2.S\), then
Lemma~\ref{lem:L34-admissible-covers-subgroups} gives a subgroup
\(A_6\) of \(D\). Restricting the nonzero module \(Z\)
to it contradicts Lemma~\ref{lem:A6-triple-seven-five-fixed}.
If \(D=3.S\), use the subgroup \(3.A_6\) of \(D\)
and the second part of the same lemma. Thus \(7\nmid|K|\).

\emph{Now suppose that \(5\mid|K|\).}
In this case \(Q=K/H\ne1\) is a \(5\)-group of exponent
\(5\), and \(Z(Q)\cong C_5^{6n}\), \(n\ge1\).
We begin by showing that \(K\) has no element of order \(9\).
Choose a Hall \(\{3,5\}\)-subgroup \(J\le K\),
its normal Sylow \(3\)-subgroup \(U=J\cap H\),
and a subgroup \(A\cong C_5^6\) in a Sylow \(5\)-subgroup
of \(J\), whose existence follows from
Lemma~\ref{lem:L34-global-two-three-complement}.
For every \(1\ne a\in A\), the group \(C_U(a)\)
has no element of order \(9\), since otherwise \(45\) would occur.
Order \(27\) is also forbidden.
Lemma~\ref{lem:five-rank-exponent-three} gives
\(\exp U\mid3\), whence
\[
 9\notin\omega(K).
\]

Apply Lemma~\ref{lem:soluble-normal-p-top-section} to
\(H\triangleleft G\) and \(p=3\). We obtain \(N\triangleleft G\),
\(N\le H\), such that \(\overline H=H/N\) has a
nontrivial normal elementary abelian Sylow \(3\)-subgroup
\(V\triangleleft\overline G=G/N\).
Put \(\overline L=L/N\), \(\overline K=K/N\).
The equalities \(\overline H=[\overline H,\overline L]\) and
\(9\notin\omega(\overline K)\) remain valid, and
\(\pi(\overline K/V)\subseteq\{2,5\}\).

By Lemma~\ref{lem:L34-admissible-covers-subgroups}, choose
\(S_0\cong L_2(7)\) in \(D=L/K\), and let
\(Y\le\overline L\) be its full preimage. Then
\[
 Y/\overline K\cong S_0,\qquad
 Y/\overline H\cong Q\rtimes S_0.
\]
The second splitting follows from the Schur--Zassenhaus theorem:
here \(Q\) is a normal Hall \(5\)-subgroup,
since \(5\nmid|S_0|=168\).

We check that \(9\notin\omega(Y)\). If \(|x|=9\) in \(Y\),
its image in \(S_0\) has order \(3\). Its image
in \(Y/\overline H=Q\rtimes S_0\) also has order \(3\):
its cube belongs to \(Q\) and at the same time has \(3\)-power order.
On \(Z(Q)\), it acts as the corresponding \(a\in S_0\)
of order \(3\). Elements of order \(7\) in \(S_0\) act
on \(Z(Q)\) without nonzero fixed points because \(35\) is forbidden.
The calculation with \(7{:}3\) in the proof of
Lemma~\ref{lem:L27-no-nine-after-elementary-three-section} gives
\[
 \dim_{\mathbb F_5}C_{Z(Q)}(x\overline H)=2n\ge2.
\]
Hence \(x\overline H\) centralizes a subgroup
\(B\cong C_5^2\) of \(Q\). The full preimage of
\(\langle x\overline H\rangle\times B\) is solvable.
Its Hall \(\{3,5\}\)-subgroup \(E\) containing \(x\)
satisfies
\[
 E\cap\overline H=V,\qquad E/V\cong C_3\times C_5^2.
\]
Since order \(45\) does not occur in \(\overline G\),
Lemma~\ref{lem:elementary-three-five-nine} excludes order
\(9\) in \(E\). This contradicts the choice of \(x\), so
\(9\notin\omega(Y)\).

Let \(X=Y^\infty\) be the perfect residual of \(Y\),
and put \(J_0=X\cap\overline K\), \(W=X\cap V\).
The image of \(X\) in \(Y/\overline K\cong S_0\)
is \(S_0\), so \(X/J_0\cong L_2(7)\). Moreover, \(J_0\)
is solvable, \(W\triangleleft X\) is an elementary abelian
\(3\)-group, and \(J_0/W\) embeds in \(\overline K/V\).
Thus \(\pi(J_0/W)\subseteq\{2,5\}\).
All hypotheses of Theorem~\ref{thm:L27-perfect-nine-free-abelian-kernel}
are satisfied for \(X,J_0,W\), and it gives \(W=1\).
Since \(X\triangleleft Y\), \(V\triangleleft Y\),
we obtain \([V,X]\le V\cap X=W=1\).

Let \(M\) be the normal closure of \(X\) in \(\overline L\).
Its image in \(D\) is a normal subgroup containing
\(S_0\), and therefore equals \(D\): its image in the
simple group \(S\) is nontrivial, while the remaining
quotient is both central and perfect.
Thus \(\overline L=\overline K M\).
The quotient \(\overline L/M\) is solvable as an image
of \(\overline K\) and perfect as an image of
\(\overline L\), so \(M=\overline L\).
The subgroup \(V\) is normal in \(\overline L\); since
\([V,X]=1\), all conjugates of \(X\) also centralize \(V\).
Thus \([V,\overline L]=1\).

In particular, \(V\) is central in \(\overline H\)
and is its Sylow \(3\)-subgroup. The Schur--Zassenhaus theorem
gives \(\overline H=V\times U_2\). The subgroup \(U_2\)
is characteristic in \(\overline H\), since it consists
of all its \(2\)-elements. Therefore
\[
 \overline H=[\overline H,\overline L]
 =[V\times U_2,\overline L]\le U_2,
\]
which is impossible when \(V\ne1\). The branch \(5\mid|K|\)
is also excluded, leaving \(K=H\), \(\pi(K)\subseteq\{2,3\}\).
\end{proof}

\Needspace{11\baselineskip}
\begin{corollary}\label{cor:L34-radical-two-three-only}
Under the standing assumptions for the socle \(L_3(4)\),
the prime divisors of the solvable radical are precisely
\[
 \pi(R(G))=\{2,3\}.
\]
In particular, both primes \(5,7\) are completely absent
from the radical, not merely from the canonical subgroup \(K\).
\end{corollary}
\begin{proof}
By the preceding theorem, \(\pi(K)\subseteq\{2,3\}\).
The group \(P/K=Z(L/K)\) is also a \(\{2,3\}\)-group,
by Lemma~\ref{lem:L34-admissible-covers-subgroups}.
Lemma~\ref{lem:perfect-residual-product} gives
\(\omega(R(G)/P)\subseteq\{1,2,3\}\).
Hence \(\pi(R(G))\subseteq\{2,3\}\).
We already know that \(3\mid|K|\).

If \(2\nmid|R(G)|\), the radical would be a \(3\)-group.
An element of order \(20\) in \(G\) would then retain its order
in \(G/R(G)\le\operatorname{Aut}(L_3(4))\).
But this automorphism group has no element of order \(20\):
the outer group \(C_2\times S_3\) has no element of order \(4\),
and \(L_3(4)\) has no element of order \(10\).
The image of an element of order \(20\) in the outer group
would have order \(1\) or \(2\), so its square would be
an element of order \(10\) in \(L_3(4)\).
This contradiction gives \(2\mid|R(G)|\).
\end{proof}

\clearpage
\subsection{Two general module exclusions for \(L_3(4)\)}

We retain \(S=L_3(4)\), \(L=G^\infty\),
\(P=R(G)\cap L\), \(K=[P,L]=[K,L]\).
By Theorem~\ref{thm:L34-K-two-three-only}, the group \(K\)
is a \(\{2,3\}\)-group. The following exclusions
apply simultaneously to all four covers
in Lemma~\ref{lem:L34-admissible-covers-subgroups}.

\Needspace{11\baselineskip}
\begin{lemma}\label{lem:L34-F2-three-fixed-table}
Let \(\widetilde S=\operatorname{SL}_3(4)=3.S\), and let
\(b=\operatorname{diag}(1,\omega,\omega^2)\),
\(\omega^3=1\ne\omega\). The complete list of dimensions
of simple \(\mathbb F_2\widetilde S\)-modules and their
fixed-point spaces is
\[
\begin{array}{c|rrrrrrrrrr}
 \dim W&1&6&6&9&9&16&18&48&48&64\\\hline
 \dim C_W(b)&1&2&2&3&3&4&6&16&16&22\\
 \dim[W,b]&0&4&4&6&6&12&12&32&32&42.
\end{array}
\]
The two six-dimensional types are dual to one another
and are not self-dual. For \(S,2.S,2^2.S\) in characteristic
\(2\), only the columns of dimensions \(1,9,9,16,64\) remain.
\end{lemma}
\begin{proof}
We use the complete form of Steinberg's theorem
stated in the proof of Lemma~\ref{lem:L34-outer3}:
all absolutely simple types of \(\widetilde S\) are
\(U\otimes W^{(2)}\), where
\(U,W\in\{1,N,N^*,A\}\),
\(\dim N=3\), \(\dim A=8\).
On \(N,N^*\), the element \(b\) has each of the eigenvalues
\(1,\omega,\omega^2\) once. On \(A\), the trace-zero
submodule of \(N\otimes N^*\), their multiplicities are
\(2,3,3\).

The Frobenius automorphism of \(\mathbb F_4/\mathbb F_2\)
interchanges the two tensor factors with the corresponding twist.
The pairs \(N,N^{(2)}\) and \(N^*,(N^*)^{(2)}\) give
two types of dimension \(6\); they are dual and distinct
by uniqueness of the types in Steinberg's theorem.
The types \(N\otimes N^{(2)}\) and
\(N^*\otimes(N^*)^{(2)}\) are defined over \(\mathbb F_2\)
and have dimension \(9\).
The pair \(A,A^{(2)}\) gives dimension \(16\).
The two mixed types \(N\otimes(N^*)^{(2)}\) and
\(N^*\otimes N^{(2)}\) form a pair of dimension \(18\).
The four absolutely simple types of dimension \(24\)
give two pairs of dimension \(48\). The last type,
\(A\otimes A^{(2)}\), is defined over \(\mathbb F_2\)
and has dimension \(64\).

The eigenvalue multiplicities on a tensor product are
obtained by multiplying the corresponding triples of multiplicities.
For example, on the \(64\)-dimensional type the multiplicity
of \(1\) is \(2^2+3^2+3^2=22\); on an absolutely simple
\(24\)-dimensional type it is \(2+3+3=8\), and doubles
when the pair is descended to \(\mathbb F_2\).
This gives the entire table.
The center of order \(3\) acts trivially precisely on
the types \(1,9,9,16,64\). Normal central \(2\)-subgroups
of the covers \(2.S,2^2.S\) act trivially on simple modules
of characteristic \(2\), so these types also factor through \(S\).
\end{proof}

\Needspace{11\baselineskip}
\begin{lemma}\label{lem:L34-orthogonal-three-uniform}
Let \(E\) be perfect, let \(U\triangleleft E\) be a
\(2\)-group, and let \(E/U\) be one of the covers
\(S,2.S,2^2.S,3.S\). Choose \(a\in E\) of order \(3\)
mapping to \(b\) from the preceding lemma in the case of
\(3.S\); for the other covers, choose an element of order
\(3\) whose image in \(S\) is \(bZ(\widetilde S)\).
Then:
\begin{enumerate}
\item for every \(\mathbb F_2E\)-module \(M\) of dimension
at least \(2\), we have \(\dim C_M(a)\ge2\);
\item if \(M\) is a nondegenerate quadratic space
with a nontrivial isometric action of \(E\), then
\[
 \dim C_M(a)\ge4,\qquad \dim[M,a]\ge4.
\]
In particular, \(C_M(a)\) contains a nonzero singular vector.
\end{enumerate}
\end{lemma}
\begin{proof}
The chosen element of order \(3\) lifts by the
Schur--Zassenhaus theorem applied to the full preimage of
its cyclic subgroup over the normal \(2\)-group.
The group \(U\) acts trivially on simple \(\mathbb F_2E\)-modules.
Thus the table of Lemma~\ref{lem:L34-F2-three-fixed-table}
applies. Taking fixed points of \(a\) is exact, since \(|a|=3\).
If all composition factors are trivial, the image of \(E\)
is unipotent and is a \(2\)-group; perfection makes
it trivial. Otherwise, at least one nontrivial factor
gives fixed-point dimension at least two and moved dimension
at least \(4\). This proves the first assertion and
the second bound in the second assertion.

In the second assertion, \(M\) is self-dual, and
\(M=C_M(a)\perp[M,a]\). The restriction of the polar
form to \(C_M(a)\) is nondegenerate, so this fixed-point space has even dimension.
If it were \(2\), the table would allow either a single
nontrivial six-dimensional composition factor and no other
factors, or only trivial factors. The latter case has
already been excluded by the nontriviality of the action.
In the former case, \(M\) itself is a simple six-dimensional
module, but neither such type is self-dual.
Thus \(\dim C_M(a)\ge4\).
A nondegenerate quadratic form over \(\mathbb F_2\)
in dimension at least \(4\) has a nonzero singular vector.
\end{proof}

\Needspace{11\baselineskip}
\begin{theorem}\label{thm:L34-two-kernel-three-central}
Suppose that \(Y\) acts faithfully on a simple \(\mathbb F_3Y\)-module
\(V\), that \(\omega(V\rtimes Y)\subseteq\Om\), and that
\(E=Y^\infty\). Let \(U\triangleleft Y\), \(U\le E\),
be a \(2\)-group, and let \(E/U\) be one of the covers
\(S,2.S,2^2.S,3.S\). Then \([U,E]=1\).
In particular, \(U=[U,E]\) implies \(U=1\).
\end{theorem}
\begin{proof}
Suppose that \([U,E]\ne1\). The nontrivial perfect image
of \(E\) in \(\operatorname{Aut}(U)\) contains a nontrivial
automorphism of odd order. By
Lemma~\ref{lem:critical-two-external}, there is a characteristic
subgroup \(C\le U\) of class at most \(2\), with
\(C/Z(C)\) elementary abelian, on which \(E\) acts nontrivially.
It is normal in \(Y\).

Extend the field to a sufficiently large finite splitting field
\(k\) of characteristic \(3\) for \(C\).
This preserves the spectrum: as an \(\mathbb F_3Y\)-module,
\(k\otimes_{\mathbb F_3}V\) is a direct sum of copies of \(V\).
For any \(y\in Y\) of order \(m\), the possibility of
obtaining order \(3m\) is determined by the nonvanishing of
\(1+y+\cdots+y^{m-1}\), which is unchanged by field extension.
We continue to denote the extended module by \(V\).

We show that \(E\) fixes every simple isomorphism type occurring
in the restriction of \(V\) to \(C\).
If \(\dim\Omega_1(Z(C))\ge2\), apply the first assertion of
Lemma~\ref{lem:L34-orthogonal-three-uniform} to this
\(\mathbb F_2E\)-module. For the chosen \(a\), this gives
a four-group in \(C_{Z(C)}(a)\).
The normal closure of \(a\) in \(E\) equals \(E\): its image
normally generates the perfect central cover \(E/U\),
and the remaining quotient is both a \(2\)-group and perfect.
If the permutation action of \(E\) on the types is nontrivial,
then \(a\) also acts nontrivially, so it has an orbit of length \(3\).
On the three corresponding isotypic spaces, \((1+a+a^2)V\ne0\).
Lemma~\ref{lem:norm-three-four-group} gives order \(18\),
a contradiction. Thus all the types are \(E\)-invariant.

If \(\dim\Omega_1(Z(C))=1\), the center of \(C\) is cyclic,
and \(E\) centralizes it. If \(|Z(C)|\ge4\), an element
of order \(7\) would give order \(28\) together with it.
Hence \(Z(C)=C_2\). The group \(C\) is nonabelian, since
otherwise the perfect group \(E\) would centralize the cyclic group \(C\).
Thus \(C\) is extraspecial. The subgroup generated by its central involution is
normal in \(Y\), so the involution is central in \(Y\); on the original
faithful simple module it acts as \(-I\).
Consequently, only its unique faithful nonlinear type occurs
in the restriction to \(C\). This type is also \(E\)-invariant.

Denote the kernels of the occurring types by \(T_i\).
They are normal in \(E\), and \(\bigcap_iT_i=1\) by faithfulness.
Choose a type for which the action of \(E\) on \(D_0=C/T_i\)
is nontrivial; otherwise \(E\) would centralize \(C\).
The group \(D_0\) has a faithful irreducible representation,
so its center is cyclic and is centralized by the perfect
\(E/T_i\). The exclusion of \(28\) again gives \(Z(D_0)=C_2\).
The group \(D_0\) is nonabelian, has class at most \(2\),
and has elementary abelian \(D_0/Z(D_0)\); therefore
it is extraspecial.

The action of \(E/T_i\) on the quadratic space
\(M=D_0/Z(D_0)\) is nontrivial: the kernel of the action
of \(\operatorname{Aut}(D_0)\) on this quotient is a
\(2\)-group and cannot contain a nontrivial perfect image.
Apply the second assertion of
Lemma~\ref{lem:L34-orthogonal-three-uniform} to \(E/T_i\)
with normal \(2\)-kernel \(U/T_i\). This gives an element
\(a_0\) of order \(3\) such that
\(\dim[M,a_0]\ge4\), while \(C_M(a_0)\) contains a nonzero
singular vector.

The isotypic space of the chosen type in \(V\) is
\(E\)-invariant, and \(T_i\) acts trivially on it.
Restrict it to \(D_0\rtimes\langle a_0\rangle\) and choose
a simple composition factor \(W\). Its restriction to \(D_0\)
is faithful. Lemma~\ref{lem:extraspecial-three-norm} gives
\[
 18\in\omega\bigl(W\rtimes(D_0\rtimes\langle a_0\rangle)\bigr).
\]
This group is a section of a subgroup of the extended
\(V\rtimes E\): first take the isotypic space, then
factor by \(T_i\), restrict to the indicated subgroup,
and pass to a composition factor. Its spectrum is therefore
contained in \(\Om\), contrary to the exclusion of \(18\).
This contradiction proves \([U,E]=1\).
\end{proof}

\Needspace{11\baselineskip}
\begin{theorem}\label{thm:L34-abelian-three-kernel-two-excluded}
Suppose that \(Y\) acts faithfully on a simple \(\mathbb F_2Y\)-module
\(V\), that \(\omega(V\rtimes Y)\subseteq\Om\), and that
\(E=Y^\infty\). Let \(U\triangleleft Y\), \(U\le E\),
be an abelian \(3\)-group. If \(U=[U,E]\) and
\(E/U\) is one of the covers \(S,2.S,2^2.S,3.S\), then \(U=1\).
\end{theorem}
\begin{proof}
Suppose that \(U\ne1\). Since \(27\notin\Om\),
we have \(\exp U\mid9\). Over a splitting field of
characteristic \(2\), decompose \(V\) according to the
linear characters of \(U\). The kernel of every occurring
character contains no element of order \(9\); otherwise
such an element would have a nonzero fixed vector in characteristic \(2\) and give \(18\).
If \(U\) had two independent cyclic factors of order \(9\),
the kernel of every linear character on their product
would contain an element of order \(9\).
Thus \(|U^3|\le3\). The image of the perfect group \(E\)
on \(U^3\) is trivial. But cubing in the abelian group \(U\)
is an \(E\)-equivariant homomorphism, so
\[
 U^3=[U,E]^3=[U^3,E]=1.
\]
Thus \(U\) is elementary abelian; we use additive notation
for it over \(\mathbb F_3\).

The action of \(E\) on \(U\) factors through \(D=E/U\)
and is nontrivial. Choose \(7{:}3=\langle\alpha,\beta\rangle\)
in a subgroup \(L_2(7)<D\), which exists by
Lemma~\ref{lem:L34-admissible-covers-subgroups}.
The element \(\alpha\) acts nontrivially on \(U\): otherwise
its normal closure in the perfect cover \(D\) would make
the entire action trivial, contrary to \(U=[U,E]\ne1\).
The nontrivial \(7\)-weights of \(U\) are permuted in
triples by \(\beta\). On each triple, the operator
\(N_\beta=1+\beta+\beta^2\) has rank equal to one third
of the total dimension. The nonzero moved part of
\(\alpha\) has dimension at least \(6\), since
\(\operatorname{ord}_7(3)=6\). Therefore
\[
 B=N_\beta U\quad\text{has}\quad\dim_{\mathbb F_3}B\ge2.
\]

Lift \(\beta\) to a \(3\)-element \(b\in E\).
Its order is \(3\) or \(9\), and \(u_0=b^3\in U\).
The action of \(b\) on \(U\) equals that of \(\beta\).
For every \(u\in U\),
\[
 (ub)^3=u_0+N_\beta u.
\]
By faithfulness of the action of \(U\) on \(V\), there
is an occurring linear character \(\chi\) nontrivial on \(B\).
It is not fixed by \(b\): a fixed character would be
trivial on \(N_\beta U\), since \(U\) has exponent \(3\).
Thus its orbit has length \(3\).

The intersection \((u_0+B)\cap\ker\chi\) contains
\(3^{\dim B-1}\ge3\) elements, so choose a nonzero
\(v\) in it. There is \(u\) with
\(u_0+N_\beta u=v\). Then \(x=ub\) has order \(9\).
The vector \(v\) is fixed by \(b\), so it belongs to
the kernels of all three characters in the orbit of \(\chi\).
On the sum of the corresponding three weight spaces,
\(x\) cyclically permutes the spaces, and \(x^3=v\)
acts as the identity. For a nonzero \(w\) in one of them,
the vector
\[
 w+xw+x^2w
\]
is nonzero and fixed by \(x\).
The fixed-point space is defined over the original field,
so \(C_V(x)\ne0\) and \(18\in\omega(V\rtimes Y)\).
This contradiction excludes \(U\ne1\).
\end{proof}

\Needspace{11\baselineskip}
\begin{corollary}\label{cor:L34-highest-internal-two-nonabelian-three}
If \(K\) is nonnilpotent, the highest chief factor
\(V=N/M\) inside \(K\) on which \(K\) acts nontrivially
must have characteristic \(2\). For its action group
\[
 Y=G/C_G(V),\qquad E=Y^\infty,
 \qquad U=\text{the image of }K\text{ in }Y
\]
the subgroup \(U\ne1\) is a normal nonabelian
\(3\)-group, \(U=[U,E]\), and \(E/U\) is one of
\(S,2.S,2^2.S,3.S\).
\end{corollary}
\begin{proof}
Such a factor exists, since otherwise a chief series
inside \(K\) would be a central series of \(K\).
The choice of the highest factor makes \(K/N\) nilpotent,
and \(N\le C_G(V)\), so \(U\ne1\) is nilpotent.
By Lemma~\ref{lem:mixed-chief-action}, the module \(V\)
is faithful and simple for \(Y\), \(\omega(V\rtimes Y)\subseteq\Om\),
and \(O_p(Y)=1\), where \(p=\operatorname{char}V\).
Thus the normal Sylow \(p\)-subgroup of \(U\) is trivial.
Since \(\pi(K)\subseteq\{2,3\}\), if \(p=3\),
then \(U\) would be a \(2\)-group, and if \(p=2\),
it would be a \(3\)-group.
The equality \(U=[U,E]\) is the image of \(K=[K,L]\).
The group \(E/U\) is a nontrivial perfect image of
\(L/K\), so it belongs to the same list of covers.
The case \(p=3\) is excluded by
Theorem~\ref{thm:L34-two-kernel-three-central}.
For \(p=2\), abelian \(U\) is excluded by
Theorem~\ref{thm:L34-abelian-three-kernel-two-excluded}.
The stated nonabelian \(3\)-group is the only remaining case.
\end{proof}

\Needspace{11\baselineskip}
\begin{theorem}\label{thm:L34-nine-in-L-noncentral}
Under the standing assumptions for \(L_3(4)\), there exists an element
\(x\in L\) of order \(9\) with nontrivial image
in \(S=L/P\). In particular, \(x\notin R(G)\).
\end{theorem}
\begin{proof}
Apply Lemma~\ref{lem:soluble-normal-p-top-section} to
\(K\triangleleft G\) and \(p=3\). We obtain
\(N\triangleleft G\), \(N\le K\), such that
\(\overline K=K/N\) has a nontrivial normal elementary
abelian Sylow \(3\)-subgroup \(V\triangleleft G/N\).
Put \(\overline L=L/N\).
Choose \(S_0\cong L_2(7)\) in \(D=L/K\),
and denote its full preimage by \(Y\le\overline L\).

Suppose that \(9\notin\omega(Y)\), and put
\(X=Y^\infty\), \(J=X\cap\overline K\), \(W=X\cap V\).
Then \(X/J\cong L_2(7)\), \(J\) is solvable,
\(W\triangleleft X\) is an elementary abelian \(3\)-group,
and \(J/W\) embeds in \(\overline K/V\), so it is a
\(2\)-group. Theorem~\ref{thm:L27-perfect-nine-free-abelian-kernel}
gives \(W=1\). Since \(X,V\) are normal in \(Y\),
we have \([V,X]=1\).

The normal closure of \(X\) in \(\overline L\) equals
all of \(\overline L\): its image normally generates
the perfect cover \(D\), while the remaining quotient
is solvable as an image of \(\overline K\) and is also perfect.
Hence \([V,\overline L]=1\).
Since \(V\) is central in \(\overline K\) and is its
Sylow \(3\)-subgroup, we have
\(\overline K=V\times U_2\), with \(U_2\) characteristic.
This contradicts \(\overline K=[\overline K,\overline L]\).
Thus \(9\in\omega(Y)\).

Order \(9\) is absent from \(\overline K\), since
its Sylow \(3\)-subgroup \(V\) is elementary abelian.
Therefore the element of order \(9\) just obtained
has nontrivial image in \(Y/\overline K=S_0\).
Lift it to a \(3\)-element \(x\) in the original \(L\).
Since order \(27\) is forbidden, \(|x|=9\).
The subgroup \(S_0\) embeds in \(S\), so the image of
\(x\) in \(S\) is nontrivial.
\end{proof}

\Needspace{11\baselineskip}
\begin{corollary}\label{cor:L34-nilpotent-K-three-only}
If the canonical subgroup \(K\) is nilpotent, then it is
a nontrivial \(3\)-group.
\end{corollary}
\begin{proof}
By Theorem~\ref{thm:L34-K-two-three-only}, we have
\(K=K_2\times K_3\), \(K_3\ne1\).
Suppose that \(K_2\ne1\), and consider
\(M=K_2/\Phi(K_2)\). The group \(K\) acts trivially on \(M\),
so the action of \(L\) factors through \(D=L/K\).
Take \(x\) from Theorem~\ref{thm:L34-nine-in-L-noncentral}.
Its order in \(D\) is \(3\), since
\(9\notin\omega(G/K)\) by
Lemma~\ref{lem:normal-cover-L34-no9}.
Its image in \(S\) is nontrivial, and hence its action
on the simple factors of \(M\) is described by the table in
Lemma~\ref{lem:L34-F2-three-fixed-table}.
In the cover \(3.S=\operatorname{SL}_3(4)\), every lift
of such an element of order \(3\) has eigenvalues
\(1,\omega,\omega^2\); otherwise it would be scalar.
For the other covers, the central \(2\)-kernels act trivially
on simple \(\mathbb F_2\)-factors.
Thus \(C_M(x)\ne0\).
Coprime lifting of fixed points gives \(C_{K_2}(x)\ne1\),
which yields order \(18\).
This contradiction excludes \(K_2\ne1\).
\end{proof}

At this point two structural cases remain: \(K\) is a
\(3\)-group, or \(K\) is nonnilpotent with a normal
nonabelian \(3\)-kernel in the action group on the highest
noncentral \(2\)-factor. The next subsection excludes
the second case completely.

\clearpage
\subsection{Excluding the nonabelian \(3\)-kernel and the entire nonnilpotent branch of \(K\)}

Here \(S=L_3(4)\), \(L=G^\infty\), \(P=R(G)\cap L\),
\(K=[P,L]=[K,L]\).
We now extend the preceding exclusion from abelian to arbitrary
\(3\)-kernels in characteristic \(2\).
The final result applies to a normal \(3\)-group of arbitrary
rank and class.

\Needspace{11\baselineskip}
\begin{lemma}\label{lem:local-three-norm-two-obstruction}
Let \(A\) be an elementary abelian \(3\)-subgroup of a
finite group \(H\), and let \(b\in N_H(A)\) have order \(3\)
or \(9\), with \(b^3\in C_H(A)\). Let \(V\) be a module
of characteristic \(2\) that is faithful on \(A\).
If the operator \(N_b=1+b+b^2\) on \(A\), regarded
over \(\mathbb F_3\), has rank at least \(2\), then
\[
 18\in\omega(V\rtimes H).
\]
\end{lemma}
\begin{proof}
Put \(B=N_bA\). Over a finite splitting field of
characteristic \(2\), choose an occurring character \(\chi\)
of \(A\) that is nontrivial on \(B\). It exists by faithfulness.
The character \(\chi\) is not fixed by \(b\), since a
fixed character would be trivial on \(N_bA\).
Its orbit has length \(3\).
The element \(b^3\) preserves the weight space \(V_\chi\),
has order dividing \(3\), and therefore has a nonzero
eigenvector \(w\) there with eigenvalue \(\lambda\),
\(\lambda^3=1\).

There are at least \(3\) elements \(c\) in \(B\) such that
\(\chi(c)=\lambda^{-1}\): the size of this fiber is
\(3^{\dim B-1}\). Choose one with \(b^3c\ne1\),
and find \(u\in A\) with \(N_bu=c\). Then, for \(x=ub\),
we have \(x^3=b^3c\ne1\). The two factors on the right
commute and have orders dividing \(3\), so \(|x|=9\).
Moreover, \(x^3w=w\), while \(x\) cyclically permutes
the three weight spaces in the orbit of \(\chi\).
Thus \(w+xw+x^2w\ne0\), and this vector is fixed by \(x\).
This gives order \(18\). The nonzero fixed-point space
descends to the original field, so the field extension
does not affect the conclusion.
\end{proof}

\Needspace{11\baselineskip}
\begin{lemma}\label{lem:L34-central-three-subgroups}
Let \(E\) be perfect, let \(U\triangleleft E\) be a \(3\)-group,
and let \(E/U\) be one of the covers \(S,2.S,2^2.S,3.S\).
Suppose that \(\omega(V\rtimes E)\subseteq\Om\) for a
module \(V\) of characteristic \(2\). If \(A\le Z(U)\)
is normal in \(E\), elementary abelian, and acts faithfully
on \(V\), then \([A,E]=1\).
If \(V\) is faithful on \(U\ne1\), then \(E/U\) can
only be \(S\) or \(3.S\).
\end{lemma}
\begin{proof}
The action of \(E\) on \(A\) factors through \(D=E/U\).
If it is nontrivial, choose \(7{:}3=\langle\alpha,\beta\rangle\)
in a subgroup \(L_2(7)<D\). The element \(\alpha\) acts
nontrivially on \(A\), since its normal closure in \(D\) is \(D\).
The nonzero moved part of \(\alpha\) has dimension
at least \(6\) over \(\mathbb F_3\). Permuting its
\(7\)-weights in triples gives \(\operatorname{rank}(1+\beta+\beta^2)\ge2\),
as in the proof of
Theorem~\ref{thm:L34-abelian-three-kernel-two-excluded}.
Lift \(\beta\) to a \(3\)-element \(b\in E\).
Its order is \(3\) or \(9\), and \(b^3\in U\)
centralizes \(A\). Lemma~\ref{lem:local-three-norm-two-obstruction}
gives the forbidden order \(18\). Thus \([A,E]=1\).

For the second assertion, take \(A=\Omega_1(Z(U))\ne1\).
It is central in \(E\). If \(D=2.S\) or \(2^2.S\),
then \(D\) contains a cyclic subgroup of order \(14\):
a central involution commutes with an element of order \(7\).
Its full preimage over \(U\) has a complement \(C_{14}\)
by the Schur--Zassenhaus theorem. This complement,
together with a central element of order \(3\) from \(A\),
gives \(42\). This contradiction excludes both covers
with even center.
\end{proof}

\Needspace{11\baselineskip}
\begin{lemma}\label{lem:L34-nine-free-three-kernel}
Let \(F\) be perfect, let \(B\triangleleft F\) be a \(3\)-group,
with \(B=[B,F]\), and let \(F/B\) be one of the covers
\(S,2.S,2^2.S,3.S\). If \(9\notin\omega(F)\), then \(B=1\).
\end{lemma}
\begin{proof}
Pass to \(\overline F=F/\Phi(B)\),
\(W=B/\Phi(B)\). Choose \(S_0\cong L_2(7)\) in \(F/B\),
and take the perfect residual \(X\) in its full preimage.
Then \(X/(X\cap W)\cong L_2(7)\), and \(X\cap W\)
is an elementary abelian \(3\)-group.
Theorem~\ref{thm:L27-perfect-nine-free-abelian-kernel}, applied
with solvable kernel \(J=X\cap W\) and \(V=J\), gives \(X\cap W=1\).
Since \(X,W\) are normal in the full preimage, \([W,X]=1\).
The normal closure of \(X\) in \(\overline F\) equals
\(\overline F\): its image normally generates the perfect
cover \(F/B\), and the remaining quotient is both a
\(3\)-group and perfect. Thus \([W,\overline F]=1\).
But \(W=[W,\overline F]\), so \(W=1\), and then \(B=1\)
by the defining property of the Frattini subgroup.
\end{proof}

\Needspace{11\baselineskip}
\begin{lemma}\label{lem:L34-U-fixes-critical-types}
Let \(E\) be perfect, let \(U\triangleleft E\) be a \(3\)-group,
with \(U=[U,E]\), and let \(E/U\) be one of the covers
\(S,2.S,2^2.S,3.S\). Suppose that
\(\omega(V\rtimes E)\subseteq\Om\), \(\operatorname{char}V=2\).
For every \(E\)-normal subgroup \(C\le U\), the group \(U\)
fixes each simple type occurring in the restriction of
\(V\) to \(C\) over a splitting field.
\end{lemma}
\begin{proof}
Let \(\Sigma\) be the set of types, and let \(F\) be the
permutation image of \(E\) on it. If an element \(x\in E\)
of order \(9\) had an orbit of length \(9\), then the sum
\(w+xw+\cdots+x^8w\), for a nonzero vector in one
isotypic space, would be nonzero and fixed by \(x\),
giving order \(18\). Hence \(x^3\) fixes all of \(\Sigma\).
Since the only order in \(\Om\) divisible by \(9\) is
\(9\), the group \(F\) has no element of order \(9\).

The image \(B\) of \(U\) is normal in the perfect group
\(F\), is a \(3\)-group, and satisfies \(B=[B,F]\).
If \(F/B=1\), then \(F\) is a perfect \(3\)-group
and is trivial. Otherwise \(F/B\) is a nontrivial
perfect image of \(E/U\), hence one of the same covers.
Lemma~\ref{lem:L34-nine-free-three-kernel} gives \(B=1\).
Thus \(U\) fixes every type.
\end{proof}

\Needspace{11\baselineskip}
\begin{lemma}\label{lem:symplectic-free-J3-isotropic}
Let \(W\) be a nondegenerate symplectic space
over \(\mathbb F_3\), and let \(g\) be an isometry of order \(3\).
Suppose that every Jordan block of \(g\) has size \(3\)
and that the number of blocks is \(r\ge4\).
Then there is a \(g\)-invariant totally isotropic
subspace \(T\le W\) such that
\[
 \operatorname{rank}\bigl((1+g+g^2)|_T\bigr)=r/2\ge2.
\]
\end{lemma}
\begin{proof}
Denote the symplectic form by \(\langle\ ,\ \rangle\),
and put \(n=g-1\), \(N=n^2=1+g+g^2\).
On \(W/nW\), the form
\[
 h(\overline v,\overline w)=\langle v,Nw\rangle
\]
is well defined. It is nondegenerate, since \(N\) is
self-adjoint and \(\ker N=nW\). It is alternating:
\(\langle v,Nv\rangle=\langle v,gv\rangle+
\langle v,g^{-1}v\rangle=0\). Hence \(r\) is even.

Put \(q(v)=\langle v,gv\rangle\).
The \(g\)-invariance of the form and \(n^3=0\), together
with the adjoint operator \(n^*=-n+n^2\), give
\[
 q(v+nw)=q(v)-h(\overline v,\overline w).
\]
Having chosen \(\overline v\ne0\), we can choose \(w\)
so that, after replacing \(v\) by \(v+nw\), we have \(q(v)=0\).
Then \(L_v=\langle v,gv,g^2v\rangle\) is three-dimensional,
invariant, and totally isotropic.

Choose \(w_0\) with \(h(\overline v,\overline w_0)=1\).
The form on \(L_v+L_{w_0}\) is nondegenerate:
the matrix pairing the two chains
\(v,nv,n^2v\) and \(w_0,nw_0,n^2w_0\) has zeros where
the sum of the indices exceeds \(2\), and has nonzero
antidiagonal \(1,-1,1\). The first chain is isotropic,
so this gives a nondegenerate six-dimensional sum of
two blocks \(J_3(1)\). Its orthogonal complement is
\(g\)-invariant and consists of \(r-2\) blocks of the same kind.
Repeating this construction gives an orthogonal sum of
\(r/2\) isotropic three-dimensional chains.
Their sum is the required \(T\); on each chain,
\(N\) has rank \(1\).
\end{proof}

\Needspace{11\baselineskip}
\begin{lemma}\label{lem:critical-three-full}
Every finite \(3\)-group \(U\) has a characteristic subgroup
\(C\) of exponent \(3\) and class at most \(2\), such that
\([U,C]\le Z(C)\), and every nontrivial automorphism of
\(U\) of order coprime to \(3\) acts nontrivially on \(C\).
\end{lemma}
\begin{proof}
Choose a Thompson critical subgroup \(D\) with the following properties:
it is characteristic, has class at most \(2\), has
\(D/Z(D)\) elementary abelian, satisfies \([U,D]\le Z(D)\),
and detects all nontrivial \(3'\)-automorphisms of \(U\).
The precise statement is given in \cite[Lemma~2.7]{CriticalOddFull}.
Put \(C=\Omega_1(D)\). It has exponent \(3\),
see \cite[Lemma~2.10]{CriticalOddFull}; this also follows
from the cube formula in a group of class \(2\) with
\(D/Z(D)\) elementary abelian. We have
\([U,C]\le Z(D)\cap C\le Z(C)\).
Finally, a \(3'\)-automorphism of \(D\) acting trivially
on \(\Omega_1(D)\) acts trivially on all of \(D\)
\cite[Lemma~2.5]{CriticalOddFull}. Hence \(C\) retains
the property of detecting \(3'\)-automorphisms.
\end{proof}

\Needspace{11\baselineskip}
\begin{theorem}\label{thm:L34-all-three-kernels-two-excluded}
Let \(E\) be a finite perfect group, let
\(U\triangleleft E\) be a \(3\)-group with \(U=[U,E]\),
and let \(E/U\) be one of the covers \(S,2.S,2^2.S,3.S\),
where \(S=L_3(4)\). If \(E\) acts on a module \(V\)
of characteristic \(2\), faithfully on \(U\), and
\(\omega(V\rtimes E)\subseteq\Om\), then \(U=1\).
The group \(U\) is not assumed abelian.
\end{theorem}
\begin{proof}
Suppose that \(U\ne1\). By
Lemma~\ref{lem:L34-central-three-subgroups}, the group
\(D_1=E/U\) is \(S\) or \(3.S\).
The image of \(E\) in \(\operatorname{Aut}(U)\) is
nontrivial and perfect. It is not a \(3\)-group and
contains a nontrivial \(3'\)-automorphism. Take \(C\)
from Lemma~\ref{lem:critical-three-full}; the action of
\(E\) on \(C\) is nontrivial.

Work over a sufficiently large finite splitting field
of characteristic \(2\). As in the proof of
Theorem~\ref{thm:L34-two-kernel-three-central}, after restriction to
the prime field the extended module is a direct sum of copies of the
original module. Thus the extension does not change the spectrum of the semidirect product.
By Lemma~\ref{lem:L34-U-fixes-critical-types}, the group \(U\)
fixes all simple types in the restriction of \(V\) to \(C\).
Choose an orbit \(\Sigma\) of \(E\) on these types
for which the action of \(E\) on \(D_0=C/T\) is nontrivial,
where \(T\) is the intersection of the kernels of the types in the orbit.
Such an orbit exists, since otherwise \(E\) would
centralize \(C\) by faithfulness. The subgroup \(T\)
is normal in \(E\), and the sum of isotypic spaces
\(V_\Sigma\) is an \(E\)-invariant module on which
\(D_0\) acts faithfully and \(T\) trivially.

Put \(E_0=E/T\), \(U_0=U/T\), \(A=Z(D_0)\).
Every element of \(U_0\) preserves all chosen types,
and hence their central characters. These characters
separate the elements of \(A\), so \([U_0,A]=1\).
The subgroup \(A\) is normal in \(E_0\), elementary abelian,
contained in \(Z(U_0)\), and acts faithfully on \(V_\Sigma\).
Lemma~\ref{lem:L34-central-three-subgroups} gives \([E_0,A]=1\).
The group \(E_0\) is transitive on the chosen types and
centralizes \(A\), so all their central characters coincide.
Their common character is faithful on \(A\), whence \(A\cong C_3\).

The group \(D_0\) is nonabelian, since otherwise it
would equal \(A\) and be centralized by \(E_0\).
Thus \(D_0\) is extraspecial of exponent \(3\).
On \(M=D_0/Z(D_0)\), the commutator defines a nondegenerate
symplectic form over \(\mathbb F_3\).
The inclusion \([U,C]\le Z(C)\) implies
\([U_0,D_0]\le Z(D_0)\), so the action of \(E_0\)
on \(M\) factors through \(E_0/U_0=D_1\).
This action is nontrivial: the kernel of the action
of \(\operatorname{Aut}(D_0)\) on \(M\), with the center fixed,
is the group of inner automorphisms, hence a \(3\)-group.
It admits no nontrivial perfect image.

Choose \(7{:}3=\langle\alpha,\beta\rangle\) in
\(L_2(7)<D_1\). All simple composition factors of
the \(\mathbb F_3D_1\)-module \(M\) are modules for \(S\):
the center of order \(3\), if present, acts trivially
on simple modules in characteristic \(3\).
At least one factor is nontrivial, since otherwise
perfection would make the entire action trivial.
By the complete table in Lemma~\ref{lem:nilpotent-tables},
on every nontrivial simple type of \(S\), an element
of order \(7\) has moved dimension at least \(12\).
Taking fixed points of \(\alpha\) is exact, so
\[
 W=[M,\alpha],\qquad \dim W\ge12.
\]
The space \(W\) is nondegenerate, symplectic, and
\(\beta\)-invariant. Its nontrivial \(7\)-weights are
permuted in triples by \(\beta\). Hence every Jordan
block of \(\beta\) on \(W\) has size \(3\), and
there are \(r=\dim W/3\ge4\) blocks.
By Lemma~\ref{lem:symplectic-free-J3-isotropic}, there is
a \(\beta\)-invariant totally isotropic \(L_0\le W\)
such that \(\operatorname{rank}(1+\beta+\beta^2)|_{L_0}\ge2\).

Let \(B\) be the full preimage of \(L_0\) in \(D_0\).
It is abelian, has exponent \(3\), and contains \(Z(D_0)\).
Lift \(\beta\) to a \(3\)-element \(b\in E_0\).
Its action on \(M\) is \(\beta\). Adjust the lift
so that \(b^3\) centralizes \(D_0\).
In the standard model of the extraspecial group of exponent \(3\),
\[
 D_0=M\times\mathbb F_3,\qquad
 (v,s)(w,t)=(v+w,s+t+\tfrac12\langle v,w\rangle),
\]
every symplectic isometry \(g\) lifts to the automorphism
\((v,s)\mapsto(gv,s)\).
Any two automorphisms fixing the center and inducing
the same action on \(M\) differ by an inner automorphism.
Multiplying \(b\) by a suitable element of \(D_0\)
gives precisely this lift of \(\beta\), whose automorphism has order \(3\).
The new \(b\) remains in the \(3\)-group \(D_0\langle b\rangle\),
has order \(3\) or \(9\), normalizes \(B\), and
\(b^3\in C_{E_0}(D_0)\).

On \(B/Z(D_0)=L_0\), the operator \(1+b+b^2\) has
rank at least \(2\); hence the same lower bound holds
on \(B\). The module \(V_\Sigma\) is faithful on \(B\).
Lemma~\ref{lem:local-three-norm-two-obstruction} gives
\(18\in\omega(V_\Sigma\rtimes E_0)\).
But this semidirect product is a section of a subgroup
of the extended \(V\rtimes E\): take the \(E\)-invariant
\(V_\Sigma\) and factor by \(T\), which acts trivially on it.
Its spectrum is contained in \(\Om\), a contradiction.
Thus \(U=1\).
\end{proof}

\Needspace{11\baselineskip}
\begin{theorem}\label{thm:L34-K-three-group}
If \(\omega(G)=\Om\) and
\(\operatorname{Soc}(G/R(G))\cong L_3(4)\), then
the canonical commutator kernel \(K=[R(G)\cap G^\infty,G^\infty]\)
is a nontrivial \(3\)-group.
The entire nonnilpotent branch of \(K\) is excluded.
\end{theorem}
\begin{proof}
If \(K\) is nonnilpotent, Corollary
\ref{cor:L34-highest-internal-two-nonabelian-three} gives
a faithful module of characteristic \(2\) and an action
group with nontrivial normal \(3\)-kernel \(U=[U,E]\)
satisfying the hypotheses of
Theorem~\ref{thm:L34-all-three-kernels-two-excluded}.
It gives \(U=1\), a contradiction.
Thus \(K\) is nilpotent, and Corollary
\ref{cor:L34-nilpotent-K-three-only} shows that it is a
nontrivial \(3\)-group.
\end{proof}

\clearpage
\subsection{Excluding the two covers with odd center for \(L_3(4)\)}

In this subsection \(S=L_3(4)\), \(L=G^\infty\),
\(P=R(G)\cap L\), \(K=[P,L]=[K,L]\).
By Theorem~\ref{thm:L34-K-three-group}, the group \(K\)
is a nontrivial \(3\)-group. Put
\[
 B=G/K,\qquad D=L/K,\qquad T=R(G)/K.
\]
We exclude the two possibilities \(D=S,3.S\) simultaneously,
without restrictions on the rank, exponent, or class of \(K\).

\Needspace{10\baselineskip}
\begin{lemma}\label{lem:L34-centralizer-and-involution}
In this notation,
\[
 T=C_B(D),\qquad \omega(T)\subseteq\{1,2,3\},
 \qquad 2\mid|T|.
\]
In particular, \(D\cong3.S\) is impossible.
\end{lemma}
\begin{proof}
The subgroup \(T\) acts trivially on \(D/Z(D)\),
since \(R(G)/K\) is a solvable normal subgroup of
\(B\), and \(B/T\) is almost simple with socle \(S\).
More precisely, the commutator of the images of \(T,D\)
lies in \(T\cap D=Z(D)\). By
Lemma~\ref{lem:perfect-automorphism-tools},
an automorphism of a perfect group that is trivial
modulo the center is trivial. Thus \([T,D]=1\).
Conversely, the image of \(C_B(D)\) in the almost simple
\(B/T\) centralizes its socle and is therefore trivial.

The equality \(\pi(R(G))=\{2,3\}\) gives \(\pi(T)\subseteq\{2,3\}\).
Every element of \(T\) commutes with an element of order \(7\)
in \(D\). Orders \(4,6,9\) in \(T\) would give, respectively,
\(28,42,63\) in \(B\). Every larger order in \(\Om\)
whose prime divisors belong to \(\{2,3\}\) is divisible
by one of \(4,6,9\). Hence \(\omega(T)\subseteq\{1,2,3\}\).

Conjugation gives an embedding \(B/T\hookrightarrow\Aut(S)\).
The group \(\Aut(S)\) has no element of order \(20\):
its outer group \(C_2\times S_3\) has no element of
order \(4\), so the square of a supposed element of
order \(20\) would lie in \(S\) and have order \(10\),
which is absent from \(S\).
An element of order \(20\) in \(G\) retains this order
in \(B\), since \(K\) is a \(3\)-group.
If \(|T|\) were odd, its image in \(B/T\) would
still have order \(20\): the intersection of a cyclic
group of order \(20\) with an odd-order \(\{2,3\}\)-group
has order \(1\).
This is impossible, so \(2\mid|T|\).

If \(D=3.S\), the central element of order \(3\)
in \(D\) belongs to \(T\) and commutes with every
involution of \(T\). This would give order \(6\) in \(T\),
contrary to the restriction already proved.
\end{proof}

\Needspace{10\baselineskip}
\begin{lemma}\label{lem:L34-three-four-both-signs}
Let \(W\) be an \(S\)-module of characteristic \(3\)
with a nontrivial simple composition factor.
For every element \(d\in S\) of order \(4\), both
eigenvalues \(1,-1\) occur on \(W\).
\end{lemma}
\begin{proof}
We use the same Brauer characters reconstructed in
Lemma~\ref{lem:nilpotent-tables} from the ordinary table
and decomposition matrix \cite{NilpotentData}.
For such a character \(\varphi\), put
\(a=\varphi(d^2)\), \(b=\varphi(d)\). Then
\(\varphi(d^{-1})=b\), and the multiplicities are
\[
 m_+=\frac{\varphi(1)+a+2b}{4},\qquad
 m_-=\frac{\varphi(1)+a-2b}{4}.
\]
For the three classes \(4A,4B,4C\), the complete list
of nontrivial absolutely simple types is
\[
\begin{array}{c|r|rrr|rrr|rrr}
 \varphi(1)&a&\multicolumn{3}{c|}{b}&
 \multicolumn{3}{c|}{m_+}&\multicolumn{3}{c}{m_-}\\\hline
 15&-1&3&-1&-1&5&3&3&2&4&4\\
 15&-1&-1&3&-1&3&5&3&4&2&4\\
 15&-1&-1&-1&3&3&3&5&4&4&2\\
 19&3&-1&-1&-1&5&5&5&6&6&6\\
 45,45&-3&1&1&1&11&11&11&10&10&10\\
 63,63&-1&-1&-1&-1&15&15&15&16&16&16
\end{array}
\]
In particular, the degree \(15\) types are the differences
\(\chi_3-\chi_2,\chi_4-\chi_2,\chi_5-\chi_2\),
and the degree \(19\) type is \(\chi_2-\chi_1\);
the other rows are direct restrictions of ordinary characters.
All multiplicities are positive.
Restriction to \(\langle d\rangle\) is completely reducible,
since \(3\nmid4\), so eigenspaces on a composition
factor give eigenspaces on the whole module.
A field extension does not change these multiplicities.
\end{proof}

\Needspace{10\baselineskip}
\begin{lemma}\label{lem:L34-simple-cover-inversion}
If \(D=S\), every involution \(t\in T\) has an
involutory preimage \(\widehat t\in G\) that inverts
every element of \(K\).
Consequently, \(K\) is abelian, and a Sylow
\(2\)-subgroup of \(T\) has order \(2\).
\end{lemma}
\begin{proof}
In the full preimage of \(\langle t\rangle\), the normal
\(3\)-subgroup \(K\) has a complement of order \(2\)
by the Schur--Zassenhaus theorem; choose its generator
\(\widehat t\). Since \(t\) centralizes \(D\),
the automorphism induced by \(\widehat t\) on \(L\) acts trivially
on \(L/K\). Conjugacy of complements in
\(K\rtimes\langle\widehat t\rangle\), or the coprime
fixed-point lifting lemma, gives
\[
 C_L(\widehat t)/C_K(\widehat t)\cong D=S.
\]
Indeed, for every \(l\in L\), the two involutions
\(\widehat t\) and \(\widehat t^{\,l}\) belong to
\(K\langle\widehat t\rangle\) and are conjugate by
an element of \(K\); therefore a suitable element
of the coset \(lK\) centralizes \(\widehat t\).

Suppose that \(A=C_K(\widehat t)\ne1\).
The group \(A\) is normal in \(C_L(\widehat t)\), and
\(A/\Phi(A)\ne0\) is an \(S\)-module of characteristic \(3\).
By Lemma~\ref{lem:nilpotent-tables}, every element
of order \(7\) has a nonzero fixed vector on this module.
Choose \(c\in C_L(\widehat t)\) of order \(7\)
lifting an element of order \(7\) in \(S\).
Coprime lifting of fixed points gives \(C_A(c)\ne1\).
This centralizer contains an element \(u\) of order \(3\).
The elements \(u,c,\widehat t\) commute pairwise
and give order \(42\), a contradiction.
Thus \(C_K(\widehat t)=1\).

For completeness, a fixed-point-free involutory automorphism
inverts every element of a finite group, as follows.
The map \(v\mapsto v^{-1}v^{\widehat t}\)
is injective and hence surjective, and every element
in its image is inverted by \(\widehat t\).
Thus inversion is an automorphism of \(K\),
and \(K\) is abelian.

The action of \(G\) on \(K\) now factors through \(B\).
Every involution of \(T\) acts as \(-I\).
By Lemma~\ref{lem:L34-centralizer-and-involution}, a Sylow
\(2\)-subgroup \(U\) of \(T\) is elementary abelian.
Its action on \(K\) is faithful: a nonidentity element
of the kernel would be an involution acting both trivially
and by inversion on a nontrivial \(3\)-group.
But all nonidentity elements of \(U\) have the same
image \(-I\). Hence \(|U|\le2\), and equality follows
from \(2\mid|T|\).
\end{proof}

\Needspace{10\baselineskip}
\begin{theorem}\label{thm:L34-odd-center-covers-excluded}
For the socle \(L_3(4)\), the possibilities
\[
 L/K\cong L_3(4),\qquad L/K\cong3.L_3(4)
\]
are impossible. Consequently, any group under the standing assumptions must satisfy
\[
 K\ne1\text{ is a }3\text{-group},\qquad
 L/K\cong2.L_3(4)\ \text{or}\ 2^2.L_3(4).
\]
\end{theorem}
\begin{proof}
The triple cover is excluded by
Lemma~\ref{lem:L34-centralizer-and-involution}.
Let \(D=S\). Then \(D\cap T=1\),
\(DT=D\times T\), and \(K\) is abelian by
Lemma~\ref{lem:L34-simple-cover-inversion}.
Order \(8\) is preserved on passing from \(G\) to \(B\).
Choose \(x\in B\) of order \(8\) and a Sylow
\(2\)-subgroup \(Q\) of \(B\) containing \(x\).
Normality of \(D,T,DT\) gives
\[
 Q\cap(DT)=(Q\cap D)\times(Q\cap T),\qquad
 Q\cap T=\langle s\rangle\cong C_2.
\]
In particular, \(s\in Z(Q)\).
The group \(B/DT\) embeds in
\(\Out(S)=C_2\times S_3\), whose Sylow
\(2\)-subgroups are elementary abelian. Therefore
\[
 x^2=d s^\varepsilon,\qquad d\in Q\cap D,
 \quad \varepsilon\in\{0,1\}.
\]
Since \(|x^2|=4\) and \(|s|=2\), we have \(|d|=4\).
Also, \(x^4=d^2\in D\), so \(s\ne x^4\).
The elements \(x,s\) commute, and both \(x,sx\)
have order \(8\).

Consider the nonzero module \(V=K/3K\) over
\(\mathbb F_3\), using additive notation for the abelian \(K\).
The action of \(L\) on \(V\) factors through \(D\),
and \([V,D]=V\), since \([K,L]=K\).
The module \(V\) has a nontrivial simple \(D\)-factor.
Otherwise the image of the perfect group \(D\) would
be a group of upper unitriangular matrices in characteristic
\(3\), hence both perfect and solvable; it would be
trivial, contrary to \([V,D]=V\ne0\).

By Lemma~\ref{lem:L34-three-four-both-signs}, both
eigenvalues \(1,-1\) of \(d\) occur on \(V\).
The involution \(s\) acts as \(-I\). Thus
\[
 C_V(x^2)=C_V(d s^\varepsilon)\ne0
\]
for both values of \(\varepsilon\).
On this space \(x\) is an involutory operator,
and \(sx=-x\). At least one of the two operators
\(x,sx\) has a nonzero fixed vector.

Finally, the full preimage of \(Q\) in \(G\) has
a normal Hall \(3\)-subgroup \(K\).
By the Schur--Zassenhaus theorem, it has a complement
isomorphic to \(Q\). Thus the chosen commuting elements
\(x,s\) lift simultaneously to such a \(2\)-subgroup,
with all orders preserved.
A nonzero fixed vector on \(K/3K\) lifts to a nonidentity fixed element of \(K\) under the
action of the cyclic group of order \(8\), which is
coprime to \(|K|\).
For one of the two lifted elements of order \(8\),
the centralizer in \(K\) contains an element of order \(3\).
We obtain \(24\in\omega(G)\), the final contradiction.
\end{proof}

This argument excludes both covers with odd center.
It does not automatically extend to the two remaining
covers with \(2\)-center: the fixed-point module of a
central involution may be a module for a double cover
rather than for the simple group \(S\) itself.
The positivity of the fixed-point dimensions for elements
of order \(7\), used above, then requires a separate check.
The next subsection excludes both remaining covers
by a different argument.

\clearpage

\subsection{Complete exclusion of the socle \(L_3(4)\)}

We retain the notation of the preceding subsection:
\(L=G^\infty\), \(P=R(G)\cap L\), \(K=[P,L]=[K,L]\),
\(B=G/K\), \(D=L/K\), \(T=R(G)/K=C_B(D)\).
By Theorems~\ref{thm:L34-K-three-group}
and~\ref{thm:L34-odd-center-covers-excluded},
the only remaining possibilities are \(D=2.S,2^2.S\), where \(S=L_3(4)\),
and \(K\ne1\) is a \(3\)-group. We exclude both
possibilities below without further restrictions on the radical.

\Needspace{10\baselineskip}
\begin{lemma}\label{lem:L34-seven-three-produces-nine}
Let \(D\) be one of the admissible perfect covers of \(S\),
and let \(W\) be a module of characteristic \(3\)
on which \(D\) acts nontrivially. Then \(9\in\omega(W\rtimes D)\).
If some central involution \(z\in D\) acts trivially
on \(W\), then \(18\in\omega(W\rtimes D)\).
\end{lemma}
\begin{proof}
By Lemma~\ref{lem:L34-admissible-covers-subgroups},
the group \(D\) contains a subgroup isomorphic to \(L_2(7)\),
which contains \(\langle a,b\rangle\cong7{:}3\), where \(|a|=7\), \(|b|=3\).
The normal closure of \(a\) in \(D\) is \(D\):
its image normally generates the simple group \(S\),
so the quotient of \(D\) by this closure is central
and abelian; perfection of \(D\) makes it trivial.
Thus \(a\) acts nontrivially on \(W\).

After a finite field extension, choose a nonzero vector
\(w\) in an eigenspace of \(a\) for a nontrivial
eigenvalue of order \(7\).
The operator \(b\) permutes three distinct such spaces,
since its action on the nonzero exponents of the
eigenvalues has orbits of length \(3\). Therefore
\[
 (1+b+b^2)w=w+bw+b^2w\ne0.
\]
Nonvanishing of the operator is preserved over the original field.
Thus, for some \(v\in W\), the element \((v,b)\)
in the semidirect product has order \(9\).
If \(z\) acts trivially on \(W\), it centralizes
both \(W\) and \(D\), so the product \(z(v,b)\)
has order \(18\).
\end{proof}

\Needspace{10\baselineskip}
\begin{theorem}\label{thm:L34-fourfold-cover-excluded}
The possibility \(D\cong2^2.L_3(4)\) is excluded.
\end{theorem}
\begin{proof}
The module \(V=K/\Phi(K)\ne0\) has characteristic \(3\),
and the action of \(L\) on it factors through \(D\).
Lemma~\ref{lem:split-spectrum}, applied to
\(L/\Phi(K)\), gives \(\omega(V\rtimes D)\subseteq\Om\).
Choose a simple composition factor \(W\) of \(V\).
Every central involution acts on a simple module of odd
characteristic as \(I\) or \(-I\), since its
eigenspaces are \(D\)-invariant.

The center \(Z(D)\cong C_2^2\) acts by scalars on \(W\).
Its central character takes values in \(\{1,-1\}\),
so its kernel contains an involution \(z\).
If the action of \(D\) on \(W\) is nontrivial,
Lemma~\ref{lem:L34-seven-three-produces-nine} gives \(18\).
If it is trivial, a nonzero element of \(W\) of
order \(3\), the involution \(z\), and an element
of \(D\) of order \(7\) commute pairwise and give \(42\).
In either case a forbidden order occurs in a section
of the semidirect product. This is a contradiction.
\end{proof}

\Needspace{10\baselineskip}
\begin{lemma}\label{lem:L34-double-global-split}
If \(D\cong2.S\), then \(K\) is abelian and
\[
 G\cong K\rtimes B,\qquad L\cong K\rtimes D,
 \qquad T=Z(D)=\langle z\rangle\cong C_2.
\]
Moreover, \(z\) inverts every element of \(K\),
and \(B/D\) is a subgroup of \(C_2^2\).
\end{lemma}
\begin{proof}
The subgroup generated by the involution \(z\) is characteristic in \(D\),
so this involution is central in \(B\). Choose an involutory
preimage \(\widehat z\in G\), using the odd order
of \(K\) and the Schur--Zassenhaus theorem.
As in Lemma~\ref{lem:L34-simple-cover-inversion},
\[
 C_L(\widehat z)/C_K(\widehat z)\cong D.
\]
If \(A=C_K(\widehat z)\ne1\), the central element
\(z\) acts trivially on the nonzero module
\(A/\Phi(A)\), since it is represented by \(\widehat z\)
itself. This module is therefore a module for
\(S=D/\langle z\rangle\).
By Lemma~\ref{lem:nilpotent-tables}, an element of
order \(7\) has a nonzero fixed vector.
Lifting this vector under coprime action gives an element
of order \(3\) in \(A\) commuting with a suitable
element of order \(7\) in \(C_L(\widehat z)\).
Together with \(\widehat z\), they give \(42\).
Hence \(C_K(\widehat z)=1\).
A fixed-point-free involution inverts all of \(K\),
and \(K\) is abelian.

The centrality of \(z\) in \(B\), together with
conjugacy of complements in \(K\langle\widehat z\rangle\),
now gives \(C_G(\widehat z)\twoheadrightarrow B\).
The kernel of this map is \(C_K(\widehat z)=1\),
so this centralizer is a complement to \(K\). Identify it with
\(B\), and identify \(\widehat z\) with \(z\).
Then \(G=K\rtimes B\), and the full preimage of
\(D\) is \(L=K\rtimes D\).

By Lemma~\ref{lem:L34-centralizer-and-involution},
\(\omega(T)\subseteq\{1,2,3\}\). The central involution
\(z\in T\) excludes elements of order \(3\) in \(T\),
since they would give order \(6\). Thus \(T\) is
an elementary abelian \(2\)-group.
Suppose that some involution \(t\in T\) has
\(C_K(t)\ne1\). Then
\[
 W=\Omega_1(C_K(t))\ne0
\]
is a \(D\)-module of characteristic \(3\), since
\([t,D]=1\). The action of \(D\) on \(W\) is
nontrivial, because \(z\) acts as \(-I\).
Lemma~\ref{lem:L34-seven-three-produces-nine} gives
an element of order \(9\) in \(W\rtimes D\).
The involution \(t\) centralizes this entire subgroup,
so order \(18\) occurs, a contradiction.
All involutions of \(T\) invert \(K\).
The action of \(T\) on \(K\) is therefore faithful,
and its image has at most two elements, whence \(T=\langle z\rangle\).

Finally, \(B/T\) embeds in \(\Aut(S)\), and
\(B/D\) embeds in \(\Out(S)=C_2\times S_3\).
An outer automorphism of order \(3\) does not lift
to the double cover. We give a verification.
The three classes of order \(4\) in \(S\) are the
three regular unipotent classes and are cyclically
permuted by a diagonal automorphism of order \(3\).
Indeed, the centralizer of a regular unipotent element
in \(\mathrm{GL}_3(4)\) has determinants that are
cubes, hence equal to \(1\), so its class splits
into three classes in \(\mathrm{SL}_3(4)\).
In \(2.S\), however, the preimage of one of the three
classes of order \(4\) consists of two classes,
whereas the preimage of each of the other two consists
of a single class \cite{L34FinalData}.
This distinguished class must be preserved by every
automorphism that lifts. Hence \(3\nmid|B/D|\),
and \(B/D\) embeds in a Sylow \(2\)-subgroup of
\(\Out(S)\), that is, in \(C_2^2\).
\end{proof}

\Needspace{10\baselineskip}
\begin{lemma}\label{lem:L34-double-four-fixed-and-extensions}
For \(D=2.S\), the following facts hold.
\begin{enumerate}
\item If \(u\in D\) has order \(4\) and its image
belongs to one of the two nondistinguished classes
of order \(4\) in \(S\), then every nonzero module
\(W\) of characteristic \(3\) on which \(z=-I\)
satisfies \(C_W(u)\ne0\).
\item Every extension \(E\) with \(D\trianglelefteq E\),
\(E/D\cong C_2\), \(C_E(D)=\langle z\rangle\)
is one of the two variants of \(2.S.2_i\),
\(i=1,2,3\): the ordinary or the isoclinic variant.
For \(i=1\), both variants contain an element \(x\)
of order \(8\) whose square is as in assertion~1.
\item For \(i=1,2\), order \(20\) is absent from
both variants. For \(i=3\), it is absent from the
ordinary variant and present in the isoclinic variant \(E_*\).
In the latter, all elements of order \(10\) lie in \(D\),
and every such element \(y\) satisfies \(y^5=z\).
Moreover, \(30\notin\omega(E_*)\).
\end{enumerate}
\end{lemma}
\begin{proof}
We first give the character-table data and explain their verification.
We use CTblLib 1.3.11, files \texttt{ctoline3.tbl},
\texttt{ctbline3.tbl}, and the decomposition matrix
modulo \(3\) \cite{L34FinalData}.
The seven absolutely simple \(D\)-modules that are
nontrivial on the center have the following values
on the two nondistinguished classes of order \(4\):
\[
\begin{array}{c|rrrrrrr}
 \dim W&6&10&10&22&22&36&90\\\hline
 \varphi(u^2)&-2&2&2&-2&-2&4&2\\
 \varphi(u)=\varphi(u^{-1})&0&0&0&0&0&0&0\\
 \dim C_W(u)&1&3&3&5&5&10&23
\end{array}
\]
Here \(\dim C_W(u)=(\varphi(1)+\varphi(u^2))/4\).
In the notation of the ordinary table of \(D\),
the first five Brauer characters are obtained from rows
\(\chi_{11},\chi_{12},\chi_{13},\chi_{14},\chi_{17}\)
using the matrix
\[
 \begin{pmatrix}
 1&1&1&1&-1\\1&0&0&0&0\\0&1&0&0&0\\
 -1&-1&-1&0&1\\-1&-1&0&-1&1
 \end{pmatrix}.
\]
The last two coincide with the ordinary characters
of degrees \(36,90\) on \(3\)-regular elements.
For arbitrary modules, the assertion follows from
exactness of taking fixed points under a group of order \(4\).

The fact that these two variants exhaust the possibilities for each outer class
is described in \cite{L34FinalData}, \S2.6.4.
It can also be checked directly: choose an outer
involution in \(E/Z(D)\). The action of its preimage
on the perfect group \(D\) is uniquely determined
by its action on \(S\), since \(D\) has no central automorphisms.
The square of the preimage lies in \(Z(D)=\langle z\rangle\)
and has two possible values. These are the two variants.

The isoclinic variant is obtained as follows.
In the central product
\((E\times\langle c\rangle)/\langle(z,c^2)\rangle\),
where \(|c|=4\), consider the subgroup consisting of the inner elements \(d\in D\)
and the outer elements \(gc\), \(g\in E\setminus D\).
The square of an outer element changes according to
\[
 (gc)^2=zg^2.
\]
For \(i=1\), in the ordinary table the classes numbered
\(21,22\) have order \(8\), and their squares belong
to class \(9\), corresponding to a nondistinguished
class of order \(4\) in \(D\). Multiplication by
\(z\) preserves this class. Hence the same orders
and the required squares occur in the isoclinic variant.

To verify assertion~3, the table of orders in the
outer coset of the ordinary variant suffices:
\[
\begin{array}{c|l}
 i&\text{element orders outside }D\\\hline
 1&2,4,6,8\\
 2&2,4,6,8,14\\
 3&2,6,8,10
\end{array}
\]
For \(i=1,2\), the group \(S.2_i\) has no element
of order \(10\), and neither does \(S\).
A central kernel of order \(2\) therefore cannot
give order \(20\).
For \(i=3\), the isoclinic change turns the outer
orders \(2,6,10\) into \(4,12,20\), respectively,
whereas elements of order \(8\) still have order \(8\): the square
of an element of order \(8\) has order \(4\) and
does not become \(z\) when squared again.
Thus the orders outside \(D\) in \(E_*\) are
\(4,8,12,20\), so in particular neither \(10\) nor \(30\) occurs there.
In \(D\) itself the orders are
\(1,2,3,4,5,6,7,10,14\), and the fifth power
of every element of order \(10\) is \(z\).
This completes the verification of all assertions.
\end{proof}

\Needspace{10\baselineskip}
\begin{theorem}\label{thm:L34-arbitrary-radical}
There is no finite group \(G\) with spectrum
\(\omega(G)=\omega(S_{10})\) such that
\[
 \operatorname{Soc}(G/R(G))\cong L_3(4).
\]
The assertion requires neither criticality of \(G\)
nor nilpotence of the solvable radical.
\end{theorem}
\begin{proof}
By Theorems~\ref{thm:L34-odd-center-covers-excluded}
and~\ref{thm:L34-fourfold-cover-excluded},
only \(D=2.S\) remains. Apply
Lemma~\ref{lem:L34-double-global-split}:
\(G=K\rtimes B\), \(K\) is abelian, \(z\) inverts \(K\),
\(T=\langle z\rangle\), and \(B/D\le C_2^2\).
The module \(V=\Omega_1(K)\ne0\) has characteristic \(3\),
and \(z=-I\) on all of \(V\).

First exclude the outer type \(2_1\) in \(B/D\).
If it occurs, the full preimage of the corresponding
subgroup of order \(2\) has one of the two forms
\(2.S.2_1\). By the preceding lemma, it contains
\(x\) of order \(8\) such that \(C_V(x^2)\ne0\).
The central involution \(z\) is not \(x^4\), since
\(x^4\) projects to a nontrivial involution of \(S\).
Both \(x,zx\) have order \(8\).
On \(C_V(x^2)\), the operators \(x,zx\) are \(x,-x\),
and at least one has a nonzero fixed vector.
Together with the corresponding element of order \(8\),
it gives order \(24\) in \(G\), a contradiction.

The group \(C_2^2\) of outer automorphisms contains
all three types \(2_1,2_2,2_3\). Hence \(|B/D|\le2\).
Order \(20\) is preserved in \(B\), since \(K\)
is a \(3\)-group, and is absent from \(D\).
By Lemma~\ref{lem:L34-double-four-fixed-and-extensions},
there is a unique possibility:
\[
 B\cong E_*=\operatorname{Isoclinic}(2.L_3(4).2_3).
\]
This group has no element of order \(30\), and
every element of order \(10\) has fifth power \(z\).

But \(30\in\omega(G)\). Take \(g\in G\) of order
\(30\). Its image \(y\in B\) has order \(10\):
the kernel \(K\) is a \(3\)-group, and order
\(30\) is absent from \(B\). Thus
\(1\ne g^{10}\in K\) has order \(3\).
The element \(g\) centralizes \(g^{10}\), so
\(y\) fixes \(g^{10}\) in its action on \(K\).
But \(y^5=z\) inverts all of \(K\) and has no
nonidentity fixed points. This is impossible.
This contradiction excludes the socle \(L_3(4)\).
\end{proof}

\clearpage

\section{\texorpdfstring{The socle \(A_6\)}{The socle A6}}
\subsection{Completing the \(A_6\) case with arbitrary radical}

Let \(\omega(G)=\Om\), \(\operatorname{Soc}(G/R(G))\cong A_6\),
\(L=G^\infty\), \(P=R(G)\cap L\), \(K=[P,L]=[K,L]\).
We have already proved that \(5\nmid|P|\), that every
element of order \(7\) lies in \(K\), and that the
cube of every element of order \(9\) lies in \(K\).
In particular, \(21\mid|K|\). We now show that these conditions lead to a contradiction.

\Needspace{10\baselineskip}
\begin{lemma}\label{lem:A6-final-module-data}
The following module facts hold.
\begin{enumerate}
\item Simple \(\mathbb F_2A_6\)-modules have dimensions
\(1,4,4,16\). The ranks of \(a-1\) for \(|a|=5\)
are, respectively, \(0,4,4,12\).
\item On every nontrivial absolutely simple \(A_6\)-module
of characteristic \(3\), the norm \(1+b+b^2\), \(|b|=3\),
is nonzero. The dimensions of these modules are \(3,3,4,9\),
and the norm ranks are \(1,1,1,3\).
\item On every simple \(2.A_6\)-module of characteristic
\(3\) on which the central involution \(z\) acts
as \(-I\), an element \(d\) of order \(4\) has
both eigenvalues \(i,-i\) with equal positive
multiplicities; here \(i^2=-1\).
\item In characteristic \(7\), the only simple
\(2.A_6\)-modules on which elements of order \(5\)
have no nonzero fixed vectors are the two faithful
four-dimensional types. On each of them \(z=-I\),
and an element of order \(4\) has eigenvalues
\(i,-i\), each with multiplicity two.
Each of the two types has a two-dimensional fixed-point
space for one of the two classes of elements of order \(3\).
\end{enumerate}
\end{lemma}
\begin{proof}
We use the ordinary tables of \(A_6,2.A_6\) and
the decomposition matrices \cite{A6FinalData}.
In characteristic \(2\), the absolutely simple types have degrees
\(1,4,4,8,8\). The two degree \(4\) types are
\(\chi_2-1,\chi_3-1\), are defined over \(\mathbb F_2\),
and have \(c_5=0\). The two degree \(8\) types have
\(c_5=2\), are interchanged by Frobenius, and form
one simple \(\mathbb F_2\)-type of degree \(16\)
with \(c_5=4\). This proves assertion~1.

For assertions~2 and~3, we may use
\(2.A_6\cong\operatorname{SL}_2(9)\) directly.
The absolutely simple types in characteristic \(3\) have the form
\[
 \operatorname{Sym}^{r}N\otimes
 (\operatorname{Sym}^{s}N)^{(3)},\qquad 0\le r,s\le2,
\]
where \(N\) is the natural two-dimensional module;
this is the complete list by Steinberg's theorem
\cite{Steinberg,A6FinalData}.
The center acts as \((-1)^{r+s}\).
On \(N\), a nonidentity unipotent element has block
\(J_2\), and on \(\operatorname{Sym}^2N\) it has \(J_3\).
In characteristic \(3\),
\[
 J_2\otimes J_2\sim J_3\oplus J_1,\qquad
 J_3\otimes J_3\sim3J_3.
\]
This gives assertion~2 for every nonzero parameter pair
\((r,s)\) with even sum. For odd sum the degrees
are \(2,2,6,6\). On the element
\(d=\operatorname{diag}(i,-i)\), \(|d|=4\),
the eigenvalues \(i,-i\) have multiplicities
\(1,1\) or \(3,3\), proving assertion~3.

In characteristic \(7\), Maschke's theorem applies,
since \(7\nmid|2.A_6|\). For types with trivial
central action, \(c_5>0\) by
Lemma~\ref{lem:A6-triple-seven-five-fixed}.
The remaining ordinary characters have degrees
\(4,4,8,8,10,10\) and give
\[
\begin{array}{c|rrrrrr}
 \chi(1)&4&4&8&8&10&10\\\hline
 c_5&0&0&2&2&2&2
\end{array}
\]
by the formula \(c_5=(\chi(1)+2\chi(5A)+2\chi(5B))/5\).
For the two four-dimensional types, the values on
classes \(3A,3B\) are \((-2,1)\) and \((1,-2)\),
so \((c_{3A},c_{3B})=(0,2),(2,0)\).
Both characters vanish on elements of order \(4\),
whose square is \(z\), acting as \(-I\).
This also proves the eigenvalue assertion in part~4.
\end{proof}

\Needspace{10\baselineskip}
\begin{lemma}\label{lem:A6-symplectic-five-rank-four}
Let \(T\) be perfect, let \(A=R(T)\),
\(T/A\cong A_6\), \(5\nmid|A|\).
Let \(X\) be a faithful symplectic
\(\mathbb F_2T\)-module with \(C_X(T)=0\).
If \(|a|=5\) and \(\dim[X,a]=4\), then
\[
 A=1,\qquad T\cong A_6,\qquad \dim X=4.
\]
\end{lemma}
\begin{proof}
Repeat the part of the proof of
Lemma~\ref{lem:perfect-symplectic-rank-four-general}
that uses only simplicity of \(T/A\).
If \(a\) acts trivially on a simple factor, the image
of \(T\) on it is perfect and solvable, hence trivial.
On every nontrivial \(\mathbb F_2C_5\)-type,
the rank of \(a-1\) is at least \(4\). Thus \(X\)
has exactly one nontrivial simple factor.
A minimal submodule is nontrivial because \(C_X(T)=0\).
It is nondegenerate: otherwise its dual would occur
as a second nontrivial factor.
The orthogonal complement has only trivial factors,
and the perfect group \(T\) acts trivially on it.
Thus the complement is zero and \(X\) is simple.
In particular, \(O_2(A)=1\).

If \(A\ne1\), put \(F=F(A)\); this is a nontrivial
group of odd order not divisible by \(5\).
Its nontrivial absolutely simple types are not
self-dual, by the Frobenius--Schur indicator, as
proved in the cited lemma.
If \(a\) permutes five types, their duals form a
separate set of five and the total rank is at least \(8\).
Thus \(a\) fixes the types, and simplicity of
\(T/A\) and perfection make every type \(T\)-invariant.
If \(a\) acts nontrivially on an odd-order Sylow
subgroup of \(F\) in one such type,
Lemma~\ref{lem:HH-prime} gives rank at least \(4\)
on the isotypic space and at least another \(4\)
on the separate dual space. This is again a contradiction.
Hence \([a,F]=1\).
The normal centralizer and perfection give
\([T,F]=1\), and then \(A\le C_A(F(A))\le F(A)\).
Thus \(A=F(A)\le Z(T)\). An odd central kernel
of a perfect cover of \(A_6\) is \(1\) or \(C_3\).

Exclude \(A=C_3=\langle c\rangle\).
Since \(X\) is simple and faithful, \(C_X(c)=0\).
The subspace \([X,a]\) is nondegenerate and invariant
under \(\langle a,c\rangle\cong C_{15}\).
The corresponding eigenvalues all have order \(15\).
But the two simple \(\mathbb F_2C_{15}\)-types
of order \(15\) have dimension \(4\) and are dual
to one another: the Frobenius orbits on the exponents are
\[
 \{1,2,4,8\},\qquad\{14,13,11,7\}.
\]
A nondegenerate form requires both types, so
\(\dim[X,a]\ge8\), a contradiction.
Thus \(A=1\). Assertion~1 of
Lemma~\ref{lem:A6-final-module-data} now gives \(\dim X=4\).
\end{proof}

\Needspace{10\baselineskip}
\begin{theorem}\label{thm:A6-all-seven-chief-central}
The group \(K\) centralizes every chief
\(7\)-factor of \(G\) inside \(P\).
\end{theorem}
\begin{proof}
Suppose that \(K\) acts nontrivially on such a factor \(V\).
Let \(Y=G/C_G(V)\), \(E=Y^\infty\), \(Q=R(E)\).
Then \(E/Q\cong A_6\), \(Q\) is the image of \(P\),
\([Q,E]\ne1\), \(5\nmid|Q|\), \(O_7(Y)=1\),
and \(\omega(V\rtimes Y)\subseteq\Om\).
It follows that \(F(Q)=F_2\times F_3\).
An element \(a\) of order \(5\) outside \(Q\)
centralizes \(F_3\): otherwise Lemma~\ref{lem:HH-prime},
on the faithful module \(V\), gives a nonzero norm
of order \(5\) and the forbidden order \(35\).
The normal centralizer and perfection then give \([E,F_3]=1\).
If also \([E,F_2]=1\), then
\(Q\le C_Q(F(Q))\le F(Q)\le Z(E)\),
a contradiction. Thus \(a\) acts nontrivially on \(F_2\).

Take a critical characteristic subgroup \(C\le F_2\)
using Lemma~\ref{lem:critical-two-external}.
Every simple type in the restriction of \(V\) to
\(C\) is \(E\)-invariant: an orbit of length \(5\)
gives a nonzero fixed vector for \(a\) and order \(35\);
the kernel of the permutation action contains
\(a\notin Q\), while simplicity of \(E/Q\) and
perfection eliminate the entire permutation action.
Choose the kernel \(J\) of a type on which the
action of \(E\) is nontrivial. Since the type is \(E\)-invariant,
\(J\trianglelefteq E\); also \(J\le C\le Q\). Thus the following
quotients are well defined:
\[
 E_0=E/J,\quad Q_0=Q/J,\quad D_0=C/J,\quad X=D_0/Z(D_0).
\]
This is the same construction as in
Lemma~\ref{lem:A5-five-seven-critical-section}; we record
the required properties. The subgroup \(D_0\)
is nonabelian, has class \(2\), \(D_0'\cong C_2\),
\(Z(D_0)\cong C_2\) or \(C_4\), and
\([E_0,Z(D_0)]=1\). Its quotient \(X\) is elementary
abelian with a nondegenerate commutator form.
The chosen isotypic space \(W\) is faithful and
homogeneous on \(D_0\). The kernel of the action
of \(E_0\) on \(X\) lies in \(Q_0\).

For completeness, cyclicity of the center follows
from faithfulness of a simple type; the perfect acting group
centralizes this cyclic center. An element of order
\(8\) in the center would give \(40\) together with
\(a\), so the center has order \(2\) or \(4\).
The commutators have exponent \(2\) and form the
unique central subgroup of order \(2\).
If the action on \(X\) is trivial, the image of
\(E_0\) on \(D_0\) is a \(2\)-group, hence trivial.
In the nontrivial case, simplicity of \(E_0/Q_0\) shows that
the kernel of this action lies in \(Q_0\).
Also, \(C_X(E_0)=0\): if \(xZ(D_0)\) is fixed,
then \(x^e=x f(e)\), where \(f:E_0\to Z(D_0)\)
is a homomorphism. It is trivial, so \(x\)
centralizes \(D_0\) and belongs to its center.
When the center is \(C_4\), this additionally gives
\(\exp D_0=4\), using the squaring map from
Lemma~\ref{lem:A5-C4-no-fixed-and-exponent-four}.

We prove \(\dim[X,a]=4\) without assuming
that \(Q_0\) is nilpotent. Let \(\dim X=2n\),
\(f=\dim C_X(a)\), and let \(H\) be a faithful
simple type of \(D_0\) of degree \(N=2^n\).
Representatives of the cosets of \(D_0/Z(D_0)\)
give a basis of \(\operatorname{End}(H)\).
For an operator \(A\) realizing \(a\), normalized
by \(A^5=1\), trace orthogonality gives
\[
 \sum_{j=0}^4m_j^2=\frac{N^2+4\cdot2^f}{5},
 \qquad\sum_{j=0}^4m_j=N,
\]
where \(m_j\) are the multiplicities of the fifth roots.
On \(W=H\otimes M\), the element \(a\) has no
nonzero fixed vectors, since otherwise order \(35\)
would occur. Thus at least one of the five roots
is absent from \(A\). The Cauchy--Schwarz inequality gives
\(2^{2n-f}\le16\). Nontriviality of the action
and \(\operatorname{ord}_5(2)=4\) give \(2n-f=4\).
Equality in the Cauchy--Schwarz inequality means that the
four roots have equal multiplicities.
After normalization, the missing root is \(1\),
and the operator of \(a\) on \(M\) is scalar.
The kernel of the projective action of \(E_0\)
on \(M\) contains \(a\notin Q_0\). Its image
is therefore perfect and equals the solvable image
of \(Q_0\), so it is trivial. Absorbing the scalars
into the operators on \(H\) gives a genuine
\(E_0\)-module \(H\); perfection gives determinant \(1\).
This justifies every step of the trace argument in
Lemma~\ref{lem:A5-C4-five-rank-and-multiplicity} in the present case.

Apply Lemma~\ref{lem:A6-symplectic-five-rank-four}
to the faithful image \(T\) of \(E_0\) on \(X\).
Its radical is the image of \(Q_0\), of order not
divisible by \(5\). We obtain \(T=A_6\), \(\dim X=4\).
If \(|Z(D_0)|=2\), the group \(D_0\) is extraspecial,
and \(T\) preserves its quadratic form.
But \(|O_4^+(2)|=72\), \(|O_4^-(2)|=120\);
neither group contains \(A_6\) of order \(360\).
If \(|Z(D_0)|=4\), then \(|D_0|=64\), \(\dim H=4\).
The image \(F\le\operatorname{SL}_4\) contains \(D_0\)
and has image \(A_6\) on \(X\).
Restrict to the full preimage of a subgroup
\(A_5<A_6\) in \(F\).
Lemma~\ref{lem:A5-C4-normalizer-fixed-four} gives
an element of order \(4\) in it with a nonzero
fixed-point space. Its lift to \(E\) has order
\(4\) or \(8\), since \(16\) is forbidden, and
gives \(28\) or \(56\) in the semidirect product.
Both possibilities for the center are excluded.
This contradiction proves the theorem for every
chief \(7\)-factor specified in its statement.
\end{proof}

\Needspace{10\baselineskip}
\begin{lemma}\label{lem:A6-global-seven-quotient}
The group \(K\) has a characteristic normal
Hall \(7'\)-subgroup \(H=[H,L]\),
\(\pi(H)\subseteq\{2,3\}\), \(3\mid|H|\).
For \(Q=K/H\), we have
\[
 G/H=Q\rtimes B_1,\qquad L/H=Q\rtimes D,
 \qquad D=L/K\cong2.A_6,
\]
\[
 Q\cong C_7^{4m},\quad m\ge1,
 \qquad z\in Z(D),\ |z|=2,\quad q^z=q^{-1}\ (q\in Q).
\]
Every simple constituent of \(Q\) in characteristic
\(7\) is one of the two four-dimensional types in
assertion~4 of Lemma~\ref{lem:A6-final-module-data}.
\end{lemma}
\begin{proof}
Theorem~\ref{thm:A6-all-seven-chief-central}
and Lemma~\ref{lem:central-chief-Hall-and-generated}
give the normal Hall subgroup \(H\).
Since \(5\nmid|K|\), we have \(\pi(H)\subseteq\{2,3\}\).
Since \(3\mid|K|\), we have \(3\mid|H|\).
The equality \(H=[H,L]\) follows by factoring by
\([H,L]\): the central Hall subgroup \(H/[H,L]\)
splits as a direct factor of \(K/[H,L]\), and
\(K=[K,L]\) eliminates this factor.

Every element of order \(7\) in the original group
lies in \(K\). Hence \(Q\ne1\) is a normal Sylow
\(7\)-subgroup of \(G/H\), with \(\exp Q=7\).
Schur--Zassenhaus gives a complement \(B_1\).
Projecting the perfect residual gives
\(L/H=Q\rtimes D\), \(D=B_1^\infty\cong L/K\).
The group \(D\) is one of the covers
\(A_6,2.A_6,3.A_6,6.A_6\).
The last is impossible: on factoring by the central
subgroup of order \(3\), its Sylow \(2\)-subgroup
maps isomorphically onto a Sylow subgroup of \(2.A_6\),
and therefore contains an element of order \(8\).
Together with a central element of order \(3\),
this gives the forbidden order \(24\).

For every \(a\in D\) of order \(5\), we have
\(C_Q(a)=1\), since otherwise \(35\) occurs.
The nonzero module \(Z(Q)\) of characteristic \(7\)
excludes \(D=A_6,3.A_6\) by
Lemma~\ref{lem:A6-triple-seven-five-fixed}.
Thus \(D=2.A_6\).
If \(A=C_Q(z)\ne1\), then \(A/\Phi(A)\)
is a nonzero \(A_6\)-module of characteristic \(7\).
An element of order \(5\) has a nonzero fixed vector
on it, which lifts to \(A\), again giving \(35\).
Thus \(C_Q(z)=1\); a fixed-point-free involution
inverts \(Q\), and \(Q\) is abelian.
Now \(Q\) is elementary abelian, and the module
table leaves only the two four-dimensional types.
They are defined over \(\mathbb F_7\): their characters
are rational, and the Schur index over a finite
field is \(1\). Hence \(|Q|=7^{4m}\).
\end{proof}

\Needspace{10\baselineskip}
\begin{theorem}\label{thm:A6-canonical-K-nilpotent}
In the branch under consideration,
\[
 K=K_3\times K_7,\qquad
 K_3\ne1,\quad K_7\cong C_7^{4m},\quad m\ge1,
 \qquad\exp K_3=3.
\]
In particular, the entire nonnilpotent branch of \(K\) is excluded.
\end{theorem}
\begin{proof}
First show that \(2\nmid|K|\), without assuming nilpotence.
A Sylow \(3\)-subgroup of \(H\) from the preceding
lemma is noncyclic by
Lemma~\ref{lem:commutator-generated-Sylow-ranks}.
In a Hall \(\{3,7\}\)-subgroup \(J\le K\),
the intersection \(J\cap H\) is a normal noncyclic
Sylow \(3\)-subgroup, and a Sylow \(7\)-subgroup
is elementary abelian of rank at least \(4\).
Lemma~\ref{lem:seven-automorphisms-paired-elementary}
gives \(C_3^2\times C_7^2\).
The action of this subgroup on the nonzero
\(\Omega_1(Z(O_2(K)))\), if \(O_2(K)\ne1\),
gives a common weight with nontrivial kernels
on both factors. This yields \(42\).
Thus \(O_2(K)=1\).

Consequently, \(O_2(H)=1\), and \(F(H)=O_3(H)\ne1\).
If a Sylow \(2\)-subgroup of \(H\) is neither cyclic
nor generalized quaternion, the same lemma in a
Hall \(\{2,7\}\)-subgroup gives \(C_2^2\times C_7^2\).
Its action on \(\Omega_1(Z(F(H)))\) again gives \(42\).
A Sylow \(2\)-subgroup of \(H\) is therefore cyclic
or generalized quaternion. But \(H=[H,L]\), and
Lemma~\ref{lem:commutator-generated-Sylow-ranks}
excludes a nontrivial Sylow subgroup of either form.
Thus \(H=O_3(K)\) is a \(3\)-group,
and \(K=H\rtimes U\), where \(U\cong K/H\)
is elementary abelian of order \(7^{4m}\).

Take any chief \(3\)-factor \(V\) of \(G\) inside \(H\).
The normal \(3\)-subgroup \(H\) has trivial image
on a simple module of characteristic \(3\).
Hence the action of \(L\) factors through
\(L/H=Q\rtimes D\) from the preceding lemma.
Suppose that the image of \(Q\) on \(V\), denoted
by \(A\), is nontrivial. It is elementary abelian
of order a power of \(7\).
After a finite field extension, the linear weights
of \(A\) on \(V\) separate its elements.
The central involution \(z\in D=2.A_6\) inverts
\(A\), and every element \(b\in D\) of order \(3\)
acts nontrivially on \(A\). The latter follows
from \(\langle b^D\rangle=D\): otherwise the action
of \(D\) would be trivial, contrary to the action of \(z\).
Thus \(b\) permutes some triple of distinct weights.
Their negatives form a separate triple, since an
orbit of odd length cannot contain a weight together
with its negative.
The element \(zb\) permutes six distinct weight spaces.
The norm of order \(6\) on their sum is nonzero,
which in characteristic \(3\) gives order \(18\)
in the corresponding semidirect product.
Nonvanishing of the norm is preserved on descent
to the original field. This contradicts
Lemma~\ref{lem:split-spectrum}.

Thus \(K\) centralizes every chief \(3\)-factor
inside \(H\). Successive coprime lifting of fixed
points along this series shows \([U,H]=1\):
for each normal layer, the centralizer of \(U\)
maps onto the centralizer in the quotient, and
all the layers are already fixed.
Hence \(K=H\times U\) is nilpotent.
Order \(9\) in \(H\) would give \(63\) together
with \(U\), so \(\exp H=3\).
Both factors are nontrivial, completing the proof.
\end{proof}

\Needspace{10\baselineskip}
\begin{lemma}\label{lem:A6-final-split-inversion}
In fact,
\[
 K=K_3\times K_7\text{ is abelian},\qquad
 G=K\rtimes B,\qquad D=L/K\cong2.A_6\trianglelefteq B.
\]
The central involution \(z\in D\) inverts all of \(K\).
\end{lemma}
\begin{proof}
Choose an involutory preimage \(t\) of \(z\)
in the full preimage of \(Z(D)\); this is possible
because \(K\) has odd order. It acts by inversion
on the normal subgroup \(K_7\), by
Lemma~\ref{lem:A6-global-seven-quotient}.
Put \(A=C_{K_3}(t)\).
Lifting fixed cosets under the action of \(t\)
gives \(C_L(t)/C_K(t)\cong D\), with \(C_K(t)=A\).
If \(A\ne1\), the nonzero module \(V=A/\Phi(A)\)
of characteristic \(3\) is a module for \(A_6\):
the central element \(z\), represented by \(t\),
acts trivially on it.

For a simple constituent \(W\) of \(V\),
we have \(\omega(W\rtimes D)\subseteq\Om\).
If \(W\) is trivial, an element of order \(8\)
in \(D\), together with a nonzero vector, gives \(24\).
If \(W\) is nontrivial, assertion~2 of
Lemma~\ref{lem:A6-final-module-data} gives a nonzero
norm for some element of order \(3\), and hence
order \(9\) in \(W\rtimes D\).
The central involution \(z\) commutes with this
element and gives \(18\).
Both cases are impossible. Thus \(C_{K_3}(t)=1\).
The involution \(t\) inverts \(K_3\), and \(K_3\) is abelian.

Now \(C_K(t)=1\), and the image \(tK=z\) is central
in \(G/K\), since \(Z(D)\) is characteristic in \(D\).
Conjugacy of complements in \(K\langle t\rangle\)
gives \(C_G(t)\twoheadrightarrow G/K\), with trivial
kernel. This complement is \(B\).
\end{proof}

\Needspace{10\baselineskip}
\begin{lemma}\label{lem:A6-final-centralizer-two-part}
Let \(T=C_B(D)=R(G)/K\) in the splitting from the preceding lemma.
Then
\[
 T=\langle z\rangle\times O_3(T).
\]
\end{lemma}
\begin{proof}
As before, the equality \(T=C_B(D)\) follows from
the absence of nontrivial central automorphisms of
the perfect group \(D\) and the trivial centralizer
of the simple socle in an almost simple group.
The order of \(B\) is not divisible by \(7\), since
every element of order \(7\) in \(G\) lies in \(K\).
If \(5\mid|T|\), an element of order \(5\) in
\(T\), together with a subgroup of order \(5\)
in \(D\), would give \(C_5^2\).
Its action on the nonzero \(K_7\) has a weight
with nontrivial kernel, so some element of order
\(5\) fixes a nonzero vector of order \(7\).
This gives \(35\). Thus \(\pi(T)\subseteq\{2,3\}\).

We show that every involution \(u\in T\)
inverts \(K_7\). If \(W=C_{K_7}(u)\ne0\),
this is an \(\mathbb F_7D\)-module consisting
of the four-dimensional types in assertion~4 of
Lemma~\ref{lem:A6-final-module-data}.
For a suitable \(a\in D\) of order \(3\),
we have \(C_W(a)\ne0\). The elements \(u,a\)
and a nonzero vector in this space give \(42\).
Thus \(C_{K_7}(u)=0\).

Now let \(|t|=4\), \(t\in T\).
The involutions \(t^2,z\) invert \(K_7\),
so \(t^2z\) acts trivially on \(K_7\).
Since \(z\in Z(B)\), this element has order \(1\) or \(2\).
The case of order two has just been excluded.
Hence \(t^2=z\).
Take \(d\in D\) of order \(4\). In
\(D=\operatorname{SL}_2(9)\), the unique involution
is \(z\), so \(d^2=z\).
Consequently, \(g=dt\) is a nontrivial involution.

We prove simultaneously that \(C_{K_3}(g)\ne0\)
and \(C_{K_7}(g)\ne0\).
After a finite field extension, decompose each
module according to the operator \(t\), which
commutes with \(D\) and has eigenvalues \(i,-i\),
since \(t^2=z=-I\).
On every simple \(D\)-factor of such an eigenspace,
the element \(d\) has both eigenvalues \(i,-i\),
by assertions~3 and~4 of the module lemma.
Thus \(dt\) has a nonzero eigenspace for \(1\)
in both characteristics.
Taking fixed points of an involution is exact
in odd characteristic, and field extension does
not change whether they are nonzero.
Both centralizers are therefore nontrivial.
Together with \(g\), their elements of orders
\(3,7\) give \(42\), since \(K_3,K_7\) commute.
This contradiction excludes order \(4\) in \(T\).

A Sylow \(2\)-subgroup of \(T\) is now elementary
abelian. Its action on \(K_7\) is faithful, and
all nonidentity elements act by the same operator \(-I\).
It is therefore \(\langle z\rangle\).
This subgroup is central, and Schur--Zassenhaus gives
\(T=\langle z\rangle\times O_3(T)\).
\end{proof}

\Needspace{28\baselineskip}
\begin{theorem}\label{thm:A6-arbitrary-radical}
There is no finite group \(G\) with spectrum
\(\omega(G)=\omega(S_{10})\) such that
\[
 \operatorname{Soc}(G/R(G))\cong A_6.
\]
The assertion requires neither criticality nor
nilpotence of the solvable radical.
\end{theorem}
\begin{proof}
By the lemmas proved above, \(G=K\rtimes B\),
\(K\) is an abelian \(\{3,7\}\)-group,
\(D=2.A_6\trianglelefteq B\), and
\(C_B(D)=\langle z\rangle\times T_3\),
where \(T_3=O_3(T)\trianglelefteq B\).
Passing to \(B_0=B/T_3\) preserves \(D\)
and the required order \(20\).
Moreover, \(C_{B_0}(D)=\langle z\rangle\):
if the commutator of a preimage with \(D\) lies
in \(T_3\), it also lies in \(D\), and \(D\cap T_3=1\).
Thus \(B_0/\langle z\rangle\) is almost simple
with socle \(A_6\) and embeds in \(\Aut(A_6)\).

The group \(\Aut(A_6)\) has no element of order
\(20\): the square of such an element would lie
in \(A_6\), since \(\Out(A_6)=C_2^2\), and have
order \(10\), which is forbidden there.
An element of order \(20\) in \(B_0\) therefore
maps to an element of order \(10\) outside \(A_6\).
Of the three outer types \(A_6.2_i\), only
\(A_6.2_2\cong\operatorname{PGL}_2(9)\)
contains an element of order \(10\) \cite{A6FinalData}.
Thus \(B_0/\langle z\rangle\) contains this subgroup.

The group \(\operatorname{PGL}_2(9)\) contains
an element \(a\) of order \(8\) whose square
lies in \(A_6=\operatorname{PSL}_2(9)\) and
has order \(4\): take the image of a diagonal
matrix with diagonal-entry ratio of order \(8\)
in \(\mathbb F_9^\times\).
Choose a preimage \(b\in B_0\) of \(a\).
Then \(b^2\in D\) maps to an element of order
\(4\) in \(A_6\). Every such preimage in
\(D=\operatorname{SL}_2(9)\) has order \(8\).
Indeed, if it had order \(4\), its square
would be the unique involution \(z\), and
its projective image would have order at most \(2\).
Hence \(|b^2|=8\), so \(|b|=16\).
But \(B_0\) is a section of \(G\), and \(16\notin\Om\).
This contradiction excludes the socle \(A_6\).
\end{proof}

\clearpage

\section{Proof of the main theorem}\label{sec:classification-proof}

\begin{proof}[Proof of Theorem~\ref{thm:final-unconditional-classification}]
Let \(G\) be a nonsolvable finite \(\Om\)-critical group.
By Corollary~\ref{cor:simple-socle}, the quotient
of every nonsolvable group with spectrum \(\Om\)
by its radical is almost simple.
The complete list of possible simple socles was obtained
in Lemma~\ref{lem:standalone-simple-candidates}.
The cases other than \(A_9,A_{10}\) are excluded by
Corollaries~\ref{cor:Sp6-O8-excluded},
\ref{cor:U35-excluded}, \ref{cor:U43-excluded},
and Theorems~\ref{thm:L28-U42-arbitrary-radical},
\ref{thm:U33-J2-arbitrary-radical},
\ref{thm:A8-arbitrary-radical}, \ref{thm:A7-arbitrary-radical},
\ref{thm:A5-arbitrary-radical}, \ref{thm:L27-arbitrary-radical},
\ref{thm:L34-arbitrary-radical}, \ref{thm:A6-arbitrary-radical}.
For the socle \(A_{10}\), Theorem~\ref{thm:A10}
gives \(G\cong S_{10}\). For the socle \(A_9\),
Corollary~\ref{cor:A9-complete} gives
\(G\cong(A_9\times C_3):2\).
Both groups have the required spectrum and are critical
by Propositions~\ref{prop:S10-critical}
and~\ref{prop:H-critical}. They are not isomorphic,
since their solvable radicals have different orders.
Thus the list is complete without additional
hypotheses on the solvable radical.
\end{proof}

The classification in the main theorem concerns critical groups.
The auxiliary exclusions of simple socles have been
established for all finite groups with spectrum \(\Om\).

\clearpage

\appendix
\section{\texorpdfstring{Certificate for completion of the \(A_6\) case}{Certificate for completion of the A6 case}}\label{app:A6-final-certificate}

The following Python 3 code with NumPy checks the ranks of norms on modules in characteristic \(3\), the eigenvalues of elements of order \(4\), the duality of the primitive \(C_{15}\)-types over \(\mathbb F_2\), the numerical rows of ordinary characters, and the orders in the three outer extensions. In addition, the code uses the power maps directly to check the uniqueness of the involution in \(2.A_6\), that every preimage of an element of order \(4\) in \(A_6\) has order \(8\), and the presence of \(24\) in \(6.A_6\). Place the original \texttt{ctoalter.tbl} from CTblLib 1.3.11 alongside the code. All checks passed.
\begin{lstlisting}[language=Python]
from pathlib import Path
import ast
import numpy as np

def rank3(A):
    A=np.array(A,dtype=int)%3; row=0
    for col in range(A.shape[1]):
        pivot=next((i for i in range(row,len(A)) if A[i,col]),None)
        if pivot is None:continue
        A[[row,pivot]]=A[[pivot,row]]
        A[row]=(A[row]*pow(int(A[row,col]),-1,3))%3
        for i in range(len(A)):
            if i!=row:A[i]=(A[i]-A[i,col]*A[row])%3
        row+=1
        if row==len(A):break
    return row
def norm3(A):
    return (np.eye(len(A),dtype=int)+A+A@A)%3
J2=np.array([[1,1],[0,1]])
S2=np.array([[1,1,1],[0,1,2],[0,0,1]])
ranks=[rank3(norm3(S2)),rank3(norm3(S2)),
       rank3(norm3(np.kron(J2,J2))),rank3(norm3(np.kron(S2,S2)))]
assert ranks==[1,1,1,3]

w={0:[0],1:[1,-1],2:[2,0,-2]}
for r,s in [(1,0),(0,1),(1,2),(2,1)]:
    weights=[(a+3*b)%4 for a in w[r] for b in w[s]]
    assert weights.count(1)==weights.count(3)==len(weights)//2

O={1,2,4,8};dual={(-a)%15 for a in O}
assert dual=={14,13,11,7} and O.isdisjoint(dual)

assert [(4+2*(-2))//5]*2+[16-12]==[0,0,4]
for d,s,c in [(4,-2,0),(8,1,2),(10,0,2)]:
    assert (d+2*s)//5==c and (d+2*s)%5==0
assert [(4+2*t)//3 for t in [-2,1]]==[0,2]
assert (6+2*2)//5==(12+2*(-1))//5==2

source=Path('ctoalter.tbl').read_text()
def bracket(pos):
    start=source.index('[',pos);depth=0;quote=False
    for j in range(start,len(source)):
        c=source[j]
        if c=='"':quote=not quote
        if quote:continue
        if c=='[':depth+=1
        if c==']':
            depth-=1
            if depth==0:return source[start:j+1],j+1
def maps(name):
    pos=source.index('MOT("'+name+'",')
    for _ in range(3):raw,pos=bracket(pos)
    arrays=[];depth=0
    for i,c in enumerate(raw[1:-1]):
        if c=='[' and depth==0:start=i+1
        if c=='[':depth+=1
        if c==']':
            depth-=1
            if depth==0:arrays.append(ast.literal_eval(raw[start:i+2]))
    return dict(zip([2,3,5],arrays))
def orders(pm):
    result=[]
    for k in range(1,len(pm[2])+1):
        order=1
        for prime,permutation in pm.items():
            seq=[];x=k
            while x not in seq:
                seq.append(x);x=permutation[x-1]
            order*=prime**seq.index(x)
        result.append(order)
    return result
for i in [1,2,3]:
    values=orders(maps('A6.2_'+str(i)))
    assert (10 in values)==(i==2)
    assert 20 not in values
    print('A6.2_'+str(i),sorted(set(values)))
D=maps('2.A6');dorders=orders(D)
assert set(dorders)=={1,2,3,4,5,6,8,10}
assert dorders.count(2)==1
fusion=[1,1,2,3,3,4,4,5,5,6,6,7,7]
sorders=orders(maps('A6'))
assert {o for o,f in zip(dorders,fusion) if sorders[f-1]==4}=={8}
assert 24 in orders(maps('6.A6'))
print('A6 certificates: all assertions passed.')
\end{lstlisting}

\clearpage
\section{\texorpdfstring{Certificate for completion of the \(L_3(4)\) case}{Certificate for completion of the L3(4) case}}\label{app:L34-final-certificate}

The following Python 3 code checks the table of fixed-point spaces in Lemma~\ref{lem:L34-double-four-fixed-and-extensions}, the orders of outer elements in the three ordinary extensions, the fifth powers of elements of order \(10\) in the double cover, and the absence of order \(10\) in the first two almost simple extensions. The original file \texttt{ctoline3.tbl} from CTblLib 1.3.11 must be placed alongside the code. For the isoclinic variant of the third extension, we use the rule \((gc)^2=zg^2\) proved in the text; the outer orders being checked are therefore \(4,8,12,20\). All checks passed.
\begin{lstlisting}[language=Python]
from pathlib import Path
import ast
source=Path('ctoline3.tbl').read_text()
def bracket(pos):
    first=source.index('[',pos)
    depth=0; quoted=False
    for j in range(first,len(source)):
        c=source[j]
        if c=='"': quoted=not quoted
        if quoted: continue
        if c=='[': depth+=1
        if c==']':
            depth-=1
            if depth==0:return source[first:j+1],j+1
def maps(name):
    pos=source.index('MOT("'+name+'",')
    for _ in range(3): raw,pos=bracket(pos)
    arrays=[];depth=0
    for i,c in enumerate(raw[1:-1]):
        if c=='[' and depth==0: first=i+1
        if c=='[':depth+=1
        if c==']':
            depth-=1
            if depth==0:arrays.append(ast.literal_eval(raw[first:i+2]))
    return dict(zip([2,3,5,7],arrays))
def orders(pm):
    out=[]
    for k in range(1,len(pm[2])+1):
        order=1
        for prime,permutation in pm.items():
            x=k;seq=[]
            while x not in seq:
                seq.append(x);x=permutation[x-1]
            order*=prime**seq.index(x)
        out.append(order)
    return out
names=['2.L3(4)']+['2.L3(4).2_'+str(i) for i in [1,2,3]]
ps=[maps(n) for n in names]; os=[orders(m) for m in ps]
assert set(os[0])=={1,2,3,4,5,6,7,10,14}
assert {ps[0][5][i] for i,o in enumerate(os[0]) if o==10}=={2}
for j,expected in enumerate([{2,4,6,8},{2,4,6,8,14},{2,6,8,10}],1):
    first=14 if j==1 else 15
    assert set(os[j][first:])==expected
assert os[1][20]==8 and ps[1][2][20]==9
assert 20 not in os[3]


assert {4 if o==2 else 12 if o==6 else 20 if o==10 else 8
        for o in os[3][15:]}=={4,8,12,20}

for i in [1,2]:
    assert not {10,20}.intersection(orders(maps('L3(4).2_'+str(i))))

O=[(10,2,0,0),(10,2,0,0),(28,-4,0,0),
   (28,-4,0,0),(70,-2,0,0)]
C=[[1,1,1,1,-1],[1,0,0,0,0],[0,1,0,0,0],
   [-1,-1,-1,0,1],[-1,-1,0,-1,1]]
B=[tuple(sum(row[i]*O[i][j] for i in range(5))
         for j in range(4)) for row in C]
B += [(36,4,0,0),(90,2,0,0)]
fixed=[sum(row)//4 for row in B]
assert all(sum(row)%4==0 for row in B)
assert [row[0] for row in B]==[6,10,10,22,22,36,90]
assert fixed==[1,3,3,5,5,10,23]
print('Power maps and faithful mod-3 modules: all assertions passed.')
\end{lstlisting}

\clearpage
\section{Reproducible verification of spectra}

The following code uses only the Python 3 standard library. It enumerates all partitions of \(n\), calculates the least common multiple of the cycle lengths, and takes the sign of the permutation into account. For \(n=10\), all \(42\) partitions are enumerated; for \(n=9\), all \(30\) partitions are enumerated. The script checks the spectra used in the text, the equality \(\omega(H)=\Om\), and the condition for the three specified pairs. It does not check the spectra of all groups in the full list of simple candidates. The proof of the restriction on pairs in Theorem~\ref{thm:first-premise} uses explicit forbidden orders and does not rely on this script.

\begin{lstlisting}[language=Python]
from math import lcm
from itertools import product

def partitions(n, lower=1):
    if n == 0:
        yield ()
    for a in range(lower, n + 1):
        for tail in partitions(n - a, a):
            yield (a,) + tail

def spectrum(n, parity=None):
    return {
        lcm(*p) for p in partitions(n)
        if parity is None or (n - len(p)) % 2 == parity
    }

def direct(A, B):
    return {lcm(a, b) for a, b in product(A, B)}

Omega = {1,2,3,4,5,6,7,8,9,10,12,14,15,20,21,30}
A9 = {1,2,3,4,5,6,7,9,10,12,15}
S9 = Omega - {21,30}
A10 = Omega - {14,20,30}
odd9 = {2,4,6,8,10,12,14,20}
assert len(list(partitions(10))) == 42
assert len(list(partitions(9))) == 30
assert spectrum(10) == Omega
assert spectrum(9) == S9
assert spectrum(9, 0) == A9
assert spectrum(10, 0) == A10
assert spectrum(9, 1) == odd9
assert direct(A9, {1,3}) | odd9 == Omega

A5, A6 = spectrum(5, 0), spectrum(6, 0)
assert A5 == {1,2,3,5}
assert A6 == {1,2,3,4,5}
for A, B in [(A5,A5), (A5,A6), (A6,A6)]:
    assert direct(A, B) <= Omega
assert direct(direct(A5, A5), A5) <= Omega
assert 60 in direct(direct(A6, A5), A5)
assert 60 not in Omega
print("All checks passed.")
\end{lstlisting}

\clearpage
\section{Verification of matrix certificates}

This short independent script over \(\mathbb F_2\) checks the ranks of \(t-I\) and \(f-I\) for both modules \(8b,8c\) in Lemmas~\ref{lem:A9-rank} and~\ref{lem:A9-good}, and then verifies the certificates in Lemma~\ref{lem:Sp6-fixedfree}. Matrix rows are encoded by integers, and addition is implemented by XOR.
\begin{lstlisting}[language=Python]
def matrix(rows):
    return [sum((c == '1') << j for j, c in enumerate(r))
            for r in rows.split()]
I = [1 << j for j in range(8)]
def mul(A, B):
    C = []
    for r in A:
        s = 0
        for j in range(8):
            if r & (1 << j):
                s ^= B[j]
        C.append(s)
    return C
def evaluate(word, A, B):
    M = I[:]
    for c in word:
        M = mul(M, A if c == 'a' else B)
    return M
def rank(A):
    basis = {}
    for r in A:
        while r:
            j = r.bit_length() - 1
            if j in basis:
                r ^= basis[j]
            else:
                basis[j] = r
                break
    return len(basis)
A = matrix('01000000 11000000 00010000 00110000 '
           '00000100 00001100 00000001 00000011')
Bb = matrix('10000000 00100000 00000100 10100001 '
            '01001100 00001010 01010100 10010001')
Bc = matrix('10000000 00100000 00001000 11000010 '
            '00000001 00101100 00011011 00110100')
for B in (Bb, Bc):
    T = evaluate('babaabbabbaabba', A, B)
    F = evaluate('aababb' * 3, A, B)
    assert rank([x ^ y for x, y in zip(T, I)]) == 4
    assert rank([x ^ y for x, y in zip(F, I)]) == 8
print('Both rank certificates verified.')

A = matrix('01000000 10000000 00001000 10011100 '
           '00100000 11101100 00000001 00000010')
B = matrix('00100000 00010000 00000100 00000010 '
           '00011010 10100000 01111110 10011111')
f = 'abbb' * 3
u = 'ab'
for word, order in ((f, 5), (u, 9)):
    M = evaluate(word, A, B)
    assert rank([x ^ y for x, y in zip(M, I)]) == 8
    assert evaluate(word * order, A, B) == I
    assert evaluate(word * (order // (5 if order == 5
                                      else 3)), A, B) != I
print('Sp6(2): both ranks and both orders verified.')
\end{lstlisting}

\clearpage
\section{\texorpdfstring{Verification of the fixed-point tables for \(S_6(2)\)}{Verification of the fixed-point tables for S6(2)}}

The following Python 3 code reproduces all four tables in the proof of Lemma~\ref{lem:Sp6-modules}. The input file \texttt{ctosymp1.tbl} is taken from the \texttt{data} subdirectory of the official CTblLib 1.3.11 archive. After saving the code as \texttt{check.py}, an example invocation is \texttt{python3 check.py ctosymp1.tbl}. The ordinary characters are read from this file; the linear combinations are specified by the MOC decomposition matrices. The numbers of absolutely irreducible characters are \(8,15,27,29\), respectively. The code uses only the Python standard library and checks that the calculated dimensions are integers.
\begin{lstlisting}[language=Python]
import ast, sys
from pathlib import Path
s=Path(sys.argv[1]).read_text().split('MOT("S6(2)",',1)[1]
i=s.index('[[1,1,1,1,1,1');j=i;depth=0
while j<len(s):
 if s[j]=='[':depth+=1
 elif s[j]==']':
  depth-=1
  if depth==0:break
 j+=1
X=ast.literal_eval(s[i:j+1]);assert len(X)==30
sub=lambda a,b:[x-y for x,y in zip(a,b)]
add=lambda a,b:[x+y for x,y in zip(a,b)]
scale=lambda a,n:[n*x for x in a]
x=lambda n:X[n-1]

p2=[x(1),sub(x(2),x(1)),sub(x(3),x(2)),sub(x(4),x(2)),
    sub(sub(x(10),x(4)),x(1)),add(sub(x(12),scale(x(4),2)),x(1)),
    add(sub(sub(x(21),scale(x(3),2)),x(10)),scale(x(1),2)),x(30)]
p3=[x(1),x(2),sub(x(3),x(1)),x(4),x(6),sub(x(8),x(1)),x(7),
    sub(x(9),x(2)),sub(x(14),sub(x(3),x(1))),sub(x(13),x(2)),
    x(17),x(18),x(19),sub(sub(x(23),x(7)),sub(x(9),x(2))),x(28)]
nondef5={1,9,11,18,22,2,6,16,27,30,4,5,17,19,26}
p5=[x(1),x(9),sub(x(11),x(1)),sub(x(18),x(9)),
    x(2),x(6),sub(x(16),x(6)),sub(x(27),x(2)),
    x(4),x(5),sub(x(17),x(4)),sub(x(19),x(5))]
p5 += [r for i,r in enumerate(X,1) if i not in nondef5]
from fractions import Fraction
f9=lambda r:Fraction(r[0]+2*r[6]+6*r[24],9)
f8a=lambda r:Fraction(r[0]+r[2]+2*r[8]+4*r[22],8)
f8b=lambda r:Fraction(r[0]+r[2]+2*r[11]+4*r[23],8)
f7=lambda r:Fraction(r[0]+6*r[21],7)
for label,rows,funcs in [('2',p2,[f9]),('3',p3,[f8a,f8b]),('5',p5,[f7])]:
 out=[(r[0],*[int(f(r)) for f in funcs]) for r in rows]
 assert all(f(r).denominator==1 for r in rows for f in funcs)
 print(label,sorted(out))
 if label=='2':assert [r[0] for r in rows if f9(r)==0]==[6,8]
 else:assert all(max(f(r) for f in funcs)>0 for r in rows)

principal=[1,3,6,15,22,28,30]
Phi=[X[0],X[2],sub(X[5],X[0])]
Phi+=[sub(X[14],Phi[2]),sub(X[21],Phi[1])]
Phi+=[sub(X[27],Phi[3])]
assert all(a+b==c for j,(a,b,c) in enumerate(zip(Phi[4],Phi[5],X[29])) if j!=21)
Phi += [row for k,row in enumerate(X,1) if k not in principal]
a=sorted((r[0],r[13],(r[0]+4*r[13])//5) for r in Phi)
assert all((r[0]+4*r[13])%5==0 for r in Phi)
assert all(f>0 for _,_,f in a)
print('Number of irreducible Brauer characters:',len(a))
print('degree,value on 5A,fixed dimension:')
print(a)
\end{lstlisting}

\clearpage
\section{\texorpdfstring{Certificate for \(U_4(3)\)}{Certificate for U4(3)}}\label{app:U43-certificate}

The following Python 3 code requires no external files. The matrices \(A,B\) are obtained by reducing modulo \(2\) the ATLAS matrices of the representation with character \(\chi_2\) in \cite{U43Data}. Binary bit \(j\), numbered from zero, encodes column \(j+1\) of a row. The code checks the order of the word \(abab^4\), the rank of the norm operator, and the two nonzero rows of its reduced row echelon form. The second part reproduces the table for an element of class \(9A\) from exact integer character values and the decomposition matrix; neither rounding nor approximate arithmetic is used.
\begin{lstlisting}[language=Python]
A = [2,1,4,32,128,8,1024,16,4096,8192,64,32768,
     256,512,262144,2048,345152,131072,16384,524288,1048576]
B = [4,8,16,64,256,512,2048,128,1024,16384,1301,194,
     65536,131072,524288,378864,877829,46898,387104,
     1048576,1589792]
I = [1 << j for j in range(21)]
def mul(A,B):
    C = []
    for a in A:
        c = 0
        while a:
            low = a & -a
            c ^= B[low.bit_length()-1]
            a ^= low
        C.append(c)
    return C
def power(A,k):
    R = I[:]
    while k:
        if k & 1:
            R = mul(R,A)
        A = mul(A,A)
        k >>= 1
    return R
def rref(A):
    R = A[:]
    k = 0
    for j in range(21):
        p = next((i for i in range(k,21) if R[i]>>j & 1),None)
        if p is None:
            continue
        R[k],R[p] = R[p],R[k]
        for i in range(21):
            if i != k and R[i]>>j & 1:
                R[i] ^= R[k]
        k += 1
    return R[:k]
T = I[:]
for c in 'ababbbb':
    T = mul(T,A if c == 'a' else B)
assert power(T,8) == I and power(T,4) != I
N = power([x ^ y for x,y in zip(T,I)],7)
assert rref(N) == [1150594,802764]
X = [(1,1,2),(21,-6,0),(35,8,4),(35,8,-2),(140,5,-2),
     (280,10,-1),(280,10,-1),(280,10,2),(280,10,2)]
C = [[1,0,0,0,0,0,0,0,0],[-1,1,0,0,0,0,0,0,0],
     [-1,0,1,0,0,0,0,0,0],[-1,0,0,1,0,0,0,0,0],
     [0,0,-1,-1,-1,1,0,0,0],[0,0,-1,-1,-1,0,1,0,0],
     [0,0,-1,-1,-1,0,0,1,0],[0,0,-1,-1,-1,0,0,0,1],
     [1,-1,0,0,1,0,0,0,0]]
P = [tuple(sum(c*x[j] for c,x in zip(row,X)) for j in range(3))
     for row in C]
P += [(640,-8,2),(640,-8,2),(896,32,-2)]
data = [(d,t,s,(d+2*t+3*s)//9) for d,t,s in P]
assert all((d+2*t+3*s)%9 == 0 and (f > 0 or d == 20)
           for d,t,s,f in data)
print('All certificates verified:',data)
\end{lstlisting}

\clearpage
\section{Verification of the tables for a nilpotent radical}

The following self-contained Python 3 code checks the tables for \(J_2,U_3(3),L_3(4)\) in Lemma~\ref{lem:nilpotent-tables} and the eigenvalue multiplicities in Lemma~\ref{lem:L34-outer3}. The original integer rows are taken from ordinary tables, and the coefficients of the linear combinations are taken from the decomposition matrices in \cite{NilpotentData}. Only exact integer arithmetic is used.
\begin{lstlisting}[language=Python]

def combine(B,C):
 return [tuple(sum(c*x[j] for c,x in zip(row,B))
               for j in range(len(B[0]))) for row in C]
J=[(1,1,1,1),(14,-2,2,0),(14,-2,2,0),(21,5,1,-1),
   (21,5,1,-1),(70,-10,2,0),(70,-10,2,0),(160,0,0,0)]
C=[[1,0,0,0,0,0,0,0],[-1,0,1,0,0,0,0,0],
   [-1,1,0,0,0,0,0,0],[0,0,0,1,0,0,0,0],
   [0,0,0,0,1,0,0,0],[1,0,-1,0,0,0,1,0],
   [1,-1,0,0,0,1,0,0],[1,-1,-1,0,0,0,0,1]]
P=combine(J,C)+[(36,4,4,0),(63,15,3,1),(90,10,-2,0),
               (189,-3,-3,1),(189,-3,-3,1),(225,-15,-3,-1)]
Jout=[]
for d,t,f,e in P:
 assert (d+t+2*f+4*e)%8==0
 c=(d+t+2*f+4*e)//8
 assert c>0
 Jout.append((d,c))
print('J2:',sorted(Jout))

U=[(1,1,2,2,2),(6,-2,-4,0,-2),(7,-1,6,-2,0),
   (7,3,-2,0,0),(7,3,-2,0,0),(21,1,-6,0,0),
   (21,1,-6,0,0),(32,0,0,0,1)]
C=[[1,0,0,0,0,0,0,0],[1,1,1,0,0,1,0,-1],
   [-1,0,-1,0,0,-1,0,1],[-1,0,0,1,0,0,0,0],
   [-1,0,0,0,1,0,0,0],[0,0,1,0,0,0,0,0],
   [1,0,0,-1,0,1,0,0],[1,0,0,0,-1,0,1,0]]
P=combine(U,C)+[(27,3,6,2,-2)]
Uout=[]
for d,t,f,e,z in P:
 assert (d+t+f+2*e)%8==0 and (d+3*z)%7==0
 c8=(d+t+f+2*e)//8
 c7=(d+3*z)//7
 assert c8>0 or c7==0
 Uout.append((d,c8,c7))
print('U3(3):',Uout)

L=[(1,2)]+[(15,2)]*3+[(19,-4)]+[(45,-1)]*2+[(63,0)]*2
for d,z in L:
 assert (d+3*z)%7==0 and (d+3*z)//7>0
print('L3(4), char 3:',[(d,(d+3*z)//7) for d,z in L])


def counts(values):
 return tuple(sum(x%3==j for x in values) for j in range(3))
def tensor(u,v):
 return tuple(sum(u[i]*v[j] for i in range(3) for j in range(3)
                  if (i+j)%3==k) for k in range(3))
for exponents in ([0,0,1],[0,1,2]):
 a=list(counts([i-j for i in exponents for j in exponents]))
 a[0]-=1
 n=counts([i+2*j for i in exponents for j in exponents])
 st=tensor(a,(a[0],a[2],a[1]))
 assert min(a)>0 and min(n)>0 and min(st)>0
 print('PGL3(4):',tuple(a),n,st)
\end{lstlisting}

\clearpage
\section{\texorpdfstring{Certificate for \(2.A_8\) in characteristic \(5\)}{Certificate for 2.A8 in characteristic 5}}
\label{app:double-A8-five}

The code checks the calculations in Lemma~\ref{lem:double-A8-five}. It uses only the Python~3 standard library and four published ATLAS matrices. If the files are in the current directory, no network access is used; otherwise, they are read directly from the primary source. All operations are performed exactly in \(\mathbb F_5\). The code also checks that each matrix has order \(7\).

\begin{lstlisting}[language=Python]
from pathlib import Path
from urllib.request import urlopen

BASE = 'https://brauer.maths.qmul.ac.uk/Atlas/alt/A8/mtx/'
DATA = [('2A8G1-f5r8B0.m2', 8, 2),
        ('2A8G1-f5r48bB0.m2', 48, 6),
        ('2A8G1-f5r64B0.m2', 64, 10),
        ('2A8G1-f5r48aB0.m2', 48, 6)]

def multiply(A, B):
    return [[sum(x*y for x, y in zip(row, col)) % 5
             for col in zip(*B)] for row in A]

def rank(A):
    A = [row[:] for row in A]
    k = 0
    for j in range(len(A[0])):
        i = next((i for i in range(k, len(A)) if A[i][j]), None)
        if i is None:
            continue
        A[k], A[i] = A[i], A[k]
        z = pow(A[k][j], -1, 5)
        A[k] = [z*x % 5 for x in A[k]]
        for i in range(k+1, len(A)):
            z = A[i][j]
            A[i] = [(x-z*y) % 5 for x, y in zip(A[i], A[k])]
        k += 1
    return k

for filename, expected_n, expected_fixed in DATA:
    path = Path(filename)
    raw = path.read_text() if path.exists() else urlopen(BASE+filename).read().decode()
    lines = raw.splitlines()
    mode, field, n, m = map(int, lines[0].split())
    assert (mode, field, n, m) == (1, 5, expected_n, expected_n)
    digits = ''.join(''.join(lines[1:]).split())
    assert len(digits) == n*n
    B = [[int(x) for x in digits[i*n:(i+1)*n]] for i in range(n)]
    I = [[int(i == j) for j in range(n)] for i in range(n)]
    power = I
    for j in range(7):
        power = multiply(power, B)
    assert power == I and B != I
    D = [[(B[i][j]-I[i][j]) % 5 for j in range(n)] for i in range(n)]
    fixed = n-rank(D)
    assert fixed == expected_fixed
    print(filename, 'dimension', n, 'fixed', fixed)
\end{lstlisting}

\clearpage

\section{\texorpdfstring{Certificate for completion of the \(A_7\) case}{Certificate for completion of the A7 case}}
\label{app:A7-final-certificate}

The following code uses only the Python 3 standard library and requires no external files. The dictionary \texttt{MODELS} contains both matrices of each ATLAS model in Lemma~\ref{lem:A7-final-module-facts}. The keys agree with the base names of the files \texttt{.m1}, \texttt{.m2} in the directory \url{https://brauer.maths.qmul.ac.uk/Atlas/alt/A7/mtx/}. Matrix rows are separated by spaces. For \(\mathbb F_9\), the label \(a+3b\) denotes \(a+b\alpha\), where \(\alpha^2=\alpha+1\). All ranks are calculated by exact Gaussian elimination over the corresponding finite field.

The code checks the relations for the standard generators, the central involution of \(2.A_7\), a word in class \(3B\), both tables of fixed-point spaces, and the rank of the norm of a \(4\)-element. Completeness of the list of simple modules is an ATLAS table fact; the code checks the values used on all the listed models. The dual \(4\)-dimensional \(A_7\)-module has the same ranks, since transposition and inversion of an operator preserve the ranks of the corresponding norms and fixed-point spaces.

\begin{lstlisting}[language=Python]
MODELS = {
    'A7G1-f2r4aB0': (2,
        '0100 1100 0001 0011',
        '0010 0110 1011 0111',
    ),
    'A7G1-f2r6B0': (2,
        '010000 000100 001000 100000 000010 000001',
        '001000 011000 000010 001100 000001 101011',
    ),
    'A7G1-f2r14B0': (2,
        '01000000000000 00010000000000 00000100000000 10000000000000 00000000100000 00000000001000 00000000000010 00000000000001 00100001101001 11101001111001 00100000000000 10111100100100 10000111110011 00001101101000',
        '00100000000000 00001000000000 00000010000000 00000001000000 00000000010000 00000000000100 10110000101110 00010010110010 11011000000101 00010101000110 01110100001101 01000110000100 00011110011110 11011011110100',
    ),
    'A7G1-f2r20B0': (2,
        '01000000000000000000 00010000000000000000 00000100000000000000 10000000000000000000 00000000100000000000 00000000001000000000 00000000000010000000 00000000000000100000 00000000000000001000 00100111100010111100 00100000000000000000 01101100111011100110 00000010000010000000 11011011110001101100 11111101001101011111 11111010010011001000 00001000000000000000 11010010001100001011 10100101101111100010 11000111010010001011',
        '00100000000000000000 00001000000000000000 00000010000000000000 00000001000000000000 00000000010000000000 00000000000100000000 00000000000001000000 00000000000000010000 00000000000000000100 00000000000000000010 00000000000000000001 01000110000100000000 00110000010011111010 11000010010100000010 00010111010101100010 01011010001100000011 10001001110010111111 11011110101011100101 11101110000101000010 10001010001001010010',
    ),
    '2A7G1-f9r4aB0': (9,
        '0100 2200 0001 0022',
        '0010 0220 2311 0252',
    ),
    '2A7G1-f9r4bB0': (9,
        '0100 2200 0001 0022',
        '0010 0220 2711 0262',
    ),
    '2A7G1-f9r6aB0': (9,
        '010000 000100 000010 100000 163824 621421',
        '001000 022000 000001 050011 554732 550488',
    ),
    '2A7G1-f9r6bB0': (9,
        '010000 000100 000010 100000 157428 521821',
        '001000 022000 000001 060011 668372 660844',
    ),
    '2A7G1-f3r8B0': (3,
        '01000000 22000000 00010000 00220000 00000010 00000001 00002020 00000202',
        '22000000 00100000 00001000 00000100 10200000 01021011 20020000 00010002',
    ),
    '2A7G1-f3r12B0': (3,
        '010000000000 001000000000 100000000000 000001000000 000000010000 210120010000 000000000100 012102010000 000000000001 010122021112 002210010010 121111201120',
        '220000000000 000100000000 000010000000 000000100000 000000001000 212210120000 000000000010 220221021000 120102010002 221121121121 120100200010 101011022012',
    ),
    '2A7G1-f3r36B0': (3,
        '010000000000000000000000000000000000 000100000000000000000000000000000000 000010000000000000000000000000000000 100000000000000000000000000000000000 000000010000000000000000000000000000 000000000100000000000000000000000000 000000000001000000000000000000000000 001000000000000000000000000000000000 000000000000001000000000000000000000 000000000000000010000000000000000000 000000000000000000100000000000000000 000000000000000000010000000000000000 000000000000000000000100000000000000 000000000000000000000001000000000000 000000000000000000000000010000000000 000000000000000000000000000100000000 000001000000000000000000000000000000 000000000000000000000000000000100000 000000000000000000000000000000001000 000000100000000000000000000000000000 000000000000000000000000000000000001 020001220010110221020012101121112110 001222101211112101001021202211011111 012122210022222022122010001120001022 112000010121211202100022021110202001 000000001000000000000000000000000000 122202202121101001220102201111012201 110020020222021022210202101012000001 100112211211021010001201122112121221 112020020202000020001011200010222000 010012211212122011110010121002010100 212002000022202001222112110220110211 000000000010000000000000000000000000 000102201221012002001021012021122112 100121122111201011121011121001002112 101000210111002100202012011221101012',
        '001000000000000000000000000000000000 022000000000000000000000000000000000 000001000000000000000000000000000000 000000100000000000000000000000000000 000000001000000000000000000000000000 000000000010000000000000000000000000 000000000000100000000000000000000000 000000000000010000000000000000000000 000000000000000100000000000000000000 000000000000000001000000000000000000 121002000010000000000000000000000000 000000000000000000001000000000000000 000000000000000000000010000000000000 000000000000000000000000100000000000 000000000000000000000000001000000000 000000000000000000000000000010000000 000000000000000000000000000001000000 000000000000000000000000000000010000 000000000000000000000000000000000100 000000000000000000000000000000000010 220111001112201102022101022102001000 020120111112220012100110100002002002 000122222020222222100111120210012020 221122211010201202111202020001110220 001010200221120101222102010122102021 102022021221110220202212001121221211 021212210002122021002011222122121011 011200201202221121122201010102201000 120012000211011201021202101002010101 010121201202212100110100202112211012 000000101210122200020112010011222102 112102100021220002201222010020112210 001220200001020000000000000001000000 201021100120122011211200121211022101 212211221121200020101220222122210100 202221211211211111001112202211021201',
    ),
}


EXPECTED2 = {4:(2,1,1), 6:(2,0,1), 14:(4,2,3), 20:(6,2,5)}
EXPECTED3 = {'f9r4a':1,'f9r4b':1,'f9r6a':2,'f9r6b':2,
             'f3r8':2,'f3r12':4,'f3r36':12}

for name,(q,sa,sb) in MODELS.items():
    p=2 if q==2 else 3
    def plus(x,y):
        if q!=9:return (x+y)%q
        return (x%3+y%3)%3+3*((x//3+y//3)%3)
    def times(x,y):
        if q!=9:return x*y%q
        a,b=x%3,x//3;c,d=y%3,y//3
        return (a*c+b*d)%3+3*((a*d+b*c+b*d)%3)
    add=[[plus(x,y) for y in range(q)] for x in range(q)]
    prod=[[times(x,y) for y in range(q)] for x in range(q)]
    neg=[next(y for y in range(q) if add[x][y]==0) for x in range(q)]
    inv=[0]+[next(y for y in range(1,q) if prod[x][y]==1)
             for x in range(1,q)]
    A=[[int(c) for c in row] for row in sa.split()]
    B=[[int(c) for c in row] for row in sb.split()]
    n=len(A)
    assert all(len(r)==n for r in A+B)
    I=[[int(i==j) for j in range(n)] for i in range(n)]
    def mul(X,Y):
        Z=[[0]*n for _ in range(n)]
        for i in range(n):
            for k in range(n):
                if X[i][k]:
                    row=prod[X[i][k]]
                    for j in range(n):Z[i][j]=add[Z[i][j]][row[Y[k][j]]]
        return Z
    def power(X,k):
        R=I
        while k:
            if k&1:R=mul(R,X)
            X=mul(X,X);k//=2
        return R
    def word(seq):
        R=I
        for X in seq:R=mul(R,X)
        return R
    def minus(X,Y):
        return [[add[x][neg[y]] for x,y in zip(rx,ry)]
                for rx,ry in zip(X,Y)]
    def rank(X):
        X=[r[:] for r in X];r=0
        for j in range(n):
            i=next((i for i in range(r,n) if X[i][j]),None)
            if i is None:continue
            X[r],X[i]=X[i],X[r]
            X[r]=[prod[inv[X[r][j]]][x] for x in X[r]]
            for i in range(n):
                if i!=r and X[i][j]:
                    a=X[i][j]
                    X[i]=[add[x][neg[prod[a][y]]]
                          for x,y in zip(X[i],X[r])]
            r+=1
        return r
    def norm(X,k):
        Z=[[0]*n for _ in range(n)];R=I
        for _ in range(k):
            Z=[[add[x][y] for x,y in zip(rx,ry)] for rx,ry in zip(Z,R)]
            R=mul(R,X)
        return Z
    assert power(A,3)==I and power(B,5)==I and power(mul(A,B),7)==I
    H=word([B,A,B,A,B])
    Hi=word([power(B,4),power(A,2),power(B,4),power(A,2),power(B,4)])
    U=word([A,Hi,A,H])
    assert power(U,3)==I and U!=I
    if q==2:
        C=word([A,A,B,A,B,B])
        assert power(C,4)==I and power(C,2)!=I
        row=(n-rank(minus(U,I)),n-rank(minus(mul(A,B),I)),rank(norm(C,4)))
        assert row==EXPECTED2[n]
        print(name,'fixed3B,fixed7,norm4:',row)
    else:
        Z=power(word([A,power(B,4),A,B]),2)
        assert Z==[[neg[x] for x in r] for r in I]
        r=rank(norm(U,3))
        key=name.split('-')[1].split('B0')[0]
        assert r==EXPECTED3[key]
        print(name,'norm3B rank:',r)
print('All exact matrix checks passed.')

\end{lstlisting}

\clearpage
\section{\texorpdfstring{Certificate for the critical section of order \(32\)}{Certificate for the critical section of order 32}}
\label{app:D32-certificate}

The code uses only the Python standard library. It enumerates all linear transformations of the four-dimensional space over \(\mathbb F_2\), retains the invertible ones, and counts the isometries of the two quadratic forms specified in Lemma~\ref{lem:extraspecial32-five-eigenvalues}. The final check verifies the action of the two transvections in Lemma~\ref{lem:extraspecial32-normalizer-fixed-four} on the five singular vectors.
\begin{lstlisting}[language=Python]
from itertools import permutations

def q(v, minus):
    a, b, c, d = [(v >> i) & 1 for i in range(4)]
    return (a*b + c*d + minus*(a+b)) % 2

def image(columns, v):
    w = 0
    for i in range(4):
        if (v >> i) & 1:
            w ^= columns[i]
    return w

for minus, expected in [(0, 72), (1, 120)]:
    count = 0
    for columns in permutations(range(1, 16), 4):
        images = [image(columns, v) for v in range(16)]
        if len(set(images)) != 16:
            continue
        if all(q(images[v], minus) == q(v, minus) for v in range(16)):
            count += 1
    assert count == expected
    print('minus =', minus, '|O(q)| =', count)

singular = [v for v in range(1, 16) if q(v, 1) == 0]
assert len(singular) == 5
for r in [1, 12]:
    assert q(r, 1) == 1
    images = [v ^ (r if q(v ^ r, 1) ^ q(v, 1) ^ q(r, 1) else 0)
              for v in singular]
    assert set(images) == set(singular)
    assert sum(v != w for v, w in zip(singular, images)) == 2
print('Both transvections act as transpositions on five singular vectors.')
\end{lstlisting}

\clearpage
\section{\texorpdfstring{Certificate for the genuine local models in characteristics \(3,5\)}{Certificate for the genuine local models in characteristics 3,5}}
\label{app:local-three-five-models}

The following code uses only the Python standard library. It constructs the groups and checks the spectra in Lemma~\ref{lem:A5-real-local-three-five-models}. The bound on the group size is an assertion to be checked: if it is exceeded, the calculation terminates with an error rather than discarding elements.
\begin{lstlisting}[language=Python]
from collections import deque

I = tuple(int(r == c) for r in range(4) for c in range(4))
I2 = (1, 0, 0, 1)

def mul(a, b, p, n=4):
    return tuple(sum(a[n*r+k]*b[n*k+c] for k in range(n)) % p
                 for r in range(n) for c in range(n))

def add(a, b, p):
    return tuple((x+y) % p for x, y in zip(a, b))

def scale(c, a, p):
    return tuple(c*x % p for x in a)

def tensor(a, b, p):
    return tuple(a[2*(r//2)+(c//2)]*b[2*(r%2)+(c%2)] % p
                 for r in range(4) for c in range(4))

def group(generators, p, bound):
    seen = {I}
    queue = deque([I])
    while queue:
        a = queue.popleft()
        for b in generators:
            c = mul(a, b, p)
            if c not in seen:
                seen.add(c)
                queue.append(c)
                assert len(seen) <= bound
    return seen

def order_and_norm(a, p):
    current, norm = I, (0,)*16
    for m in range(1, 61):
        norm = add(norm, current, p)
        current = mul(current, a, p)
        if current == I:
            return m, norm
    raise AssertionError('order exceeds 60')

def inverse(a, p):
    m, _ = order_and_norm(a, p)
    result = I
    for _ in range(m-1):
        result = mul(result, a, p)
    return result

Omega = {1,2,3,4,5,6,7,8,9,10,12,14,15,20,21,30}
expected = {
    3: {1,2,3,4,5,6,8,10,12},
    5: {1,2,3,4,5,6,8,10,12,20}
}

for p in [3, 5]:
    i = (0, 1, p-1, 0)
    j = (1,1,1,2) if p == 3 else (0,2,2,0)
    x, y = (0,1,1,0), (1,0,0,p-1)
    dgens = [tensor(i,I2,p), tensor(j,I2,p),
             tensor(I2,x,p), tensor(I2,y,p)]
    D = group(dgens, p, 32)
    assert len(D) == 32
    minus_I = scale(-1, I, p)
    rotations = sorted({min(a, scale(-1,a,p))
                        for a in D if mul(a,a,p) == minus_I})
    assert len(rotations) == 10
    R = rotations[0]
    newgens = [scale(pow(2,-1,p),
                     mul(add(I,R,p),add(I,S,p),p),p)
               for S in rotations[1:]]
    generators = dgens + newgens
    E = group(generators, p, 1920)
    assert len(E) == 1920
    commutators = []
    for d in dgens:
        for e in generators:
            commutators.append(mul(mul(mul(inverse(d,p),
                              inverse(e,p),p),d,p),e,p))
    assert group(commutators,p,32) == D
    spectrum_E, spectrum_affine = set(), set()
    for a in E:
        m, norm = order_and_norm(a,p)
        spectrum_E.add(m)
        spectrum_affine.add(m)
        if any(norm):
            spectrum_affine.add(p*m)
    assert spectrum_E == {1,2,3,4,5,6,8,10,12}
    assert spectrum_affine == expected[p]
    assert spectrum_affine < Omega
    print('p =',p,'|D| =',len(D),'|E| =',len(E))
    print('omega(E) =',sorted(spectrum_E))
    print('omega(V:E) =',sorted(spectrum_affine))
\end{lstlisting}

\clearpage
\section{\texorpdfstring{Certificate for the action of \(\operatorname{Sp}_4(2)\) on six forms}{Certificate for the action of Sp4(2) on six forms}}
\label{app:Sp4-six-forms}

The code checks the finite counts in Lemma~\ref{lem:A5-C4-normalizer-fixed-four}: the exact action of the symplectic group of order \(720\) on six minus-type forms, and the correspondence between its fifteen transvections and all transpositions in \(S_6\). Only the Python standard library is used.
\begin{lstlisting}[language=Python]
from itertools import permutations

def bit(v, i):
    return (v >> i) & 1

def b(u, v):
    return (bit(u,0)*bit(v,1) + bit(u,1)*bit(v,0)
            + bit(u,2)*bit(v,3) + bit(u,3)*bit(v,2)) % 2

def q(v):
    a, c, d, e = [bit(v,i) for i in range(4)]
    return (a+c+a*c+d*e) % 2

def image(columns, v):
    out = 0
    for i in range(4):
        if bit(v,i):
            out ^= columns[i]
    return out

zeros = [t for t in range(16) if q(t) == 0]
forms = [tuple(q(v) ^ b(t,v) for v in range(16)) for t in zeros]
assert len(forms) == 6
assert all(sum(value == 0 for value in f) == 6 for f in forms)
permutation_images = set()
for columns in permutations(range(1,16),4):
    if any(b(columns[i],columns[j]) != b(1<<i,1<<j)
           for i in range(4) for j in range(4)):
        continue
    images = [image(columns,v) for v in range(16)]
    assert len(set(images)) == 16
    perm = tuple(forms.index(tuple(f[images[v]] for v in range(16)))
                 for f in forms)
    permutation_images.add(perm)
assert len(permutation_images) == 720
transpositions = set()
for u in range(1,16):
    tau = [v ^ (u if b(v,u) else 0) for v in range(16)]
    perm = tuple(forms.index(tuple(f[tau[v]] for v in range(16)))
                 for f in forms)
    assert sum(perm[i] != i for i in range(6)) == 2
    transpositions.add(perm)
assert len(transpositions) == 15
print('|Sp4(2)| = |image on six minus forms| = 720')
print('The 15 transvections act as all 15 transpositions.')
\end{lstlisting}

\clearpage
\section{Certificate for the binary model of the upper internal factor}
\label{app:binary-upper-model}

The code checks the matrix group \(2.A_5\), the irreducibility of its action on \(\mathbb F_3^4\), the two dual orbits of length \(40\), and the spectrum of \(E=C_3^4\rtimes2.A_5\) in Lemma~\ref{lem:A5-real-binary-upper-model}. The binary module and its spectrum are then constructed in the lemma itself using weights and norms. Only the Python standard library is used.
\begin{lstlisting}[language=Python]
from itertools import product

p = 3
I = tuple(int(r == c) for r in range(4) for c in range(4))
A = (0,1,0,0, 2,0,0,0, 0,0,0,1, 0,0,2,0)
B = (1,0,0,0, 0,0,1,0, 0,2,2,0, 1,1,2,1)

def mul(a,b):
    return tuple(sum(a[4*r+k]*b[4*k+c] for k in range(4)) % 3
                 for r in range(4) for c in range(4))

def power(a,n):
    result = I
    while n:
        if n & 1:
            result = mul(result,a)
        a = mul(a,a)
        n //= 2
    return result

def group(generators,bound=120):
    seen, queue = {I}, [I]
    for a in queue:
        for b in generators:
            c = mul(a,b)
            if c not in seen:
                seen.add(c)
                queue.append(c)
                assert len(seen) <= bound
    return seen

def order_norm(a):
    current, norm = I, (0,)*16
    for m in range(1,61):
        norm = tuple((x+y)%3 for x,y in zip(norm,current))
        current = mul(current,a)
        if current == I:
            return m,norm
    raise AssertionError('order exceeds 60')

def inverse(a):
    m,_ = order_norm(a)
    return power(a,m-1)

def image(a,v):
    return tuple(sum(a[4*r+k]*v[k] for k in range(4))%3
                 for r in range(4))

def rank(rows):
    rows = [list(row) for row in rows]
    r = 0
    for c in range(4):
        pivot = next((i for i in range(r,len(rows)) if rows[i][c]),None)
        if pivot is None:
            continue
        rows[r],rows[pivot] = rows[pivot],rows[r]
        t = pow(rows[r][c],-1,3)
        rows[r] = [t*x%3 for x in rows[r]]
        for i in range(len(rows)):
            if i != r:
                t = rows[i][c]
                rows[i] = [(x-t*y)%3 for x,y in zip(rows[i],rows[r])]
        r += 1
    return r

D = group([A,B])
assert len(D) == 120
assert power(A,2) == tuple(-x%3 for x in I)
assert power(B,3) == power(mul(A,B),5) == I
comm = mul(mul(mul(inverse(A),inverse(B)),A),B)
conjugates = [mul(mul(inverse(g),comm),g) for g in D]
assert group(conjugates) == D
unseen = set(product(range(3),repeat=4)) - {(0,0,0,0)}
orbit_sizes = []
while unseen:
    v = min(unseen)
    orbit = {image(g,v) for g in D}
    assert rank(orbit) == 4
    orbit_sizes.append(len(orbit))
    unseen -= orbit
dual_D = [tuple(g[4*c+r] for r in range(4) for c in range(4)) for g in D]
unseen = set(product(range(3),repeat=4)) - {(0,0,0,0)}
dual_orbit_sizes = []
while unseen:
    v = min(unseen)
    orbit = {image(g,v) for g in dual_D}
    assert rank(orbit) == 4
    dual_orbit_sizes.append(len(orbit))
    unseen -= orbit
assert sorted(orbit_sizes) == sorted(dual_orbit_sizes) == [40,40]
spectrum_D, spectrum_E = set(),set()
for a in D:
    m,norm = order_norm(a)
    spectrum_D.add(m)
    spectrum_E.add(m)
    if any(norm):
        spectrum_E.add(3*m)
expected = {1,2,3,4,5,6,10}
assert spectrum_D == spectrum_E == expected
Omega = {1,2,3,4,5,6,7,8,9,10,12,14,15,20,21,30}
assert expected | {2*m for m in expected} < Omega
print('|D| = 120; D is perfect; nonzero orbit sizes =',sorted(orbit_sizes))
print('Nonzero dual orbit sizes =',sorted(dual_orbit_sizes))
print('omega(3^4:D) =',sorted(spectrum_E))
print('Every elementary abelian 2-extension has spectrum inside Omega.')
\end{lstlisting}

\section*{Use of ChatGPT}
This work was carried out with the assistance of ChatGPT (OpenAI).
GPT-6 was used in preparing the present version.

\begin{thebibliography}{99}\small

\bibitem{GaschutzComplement}
Aivazidis S., Loukaki M.
The \(\mathcal X\)-series of a p-group and Complements of Abelian Subgroups.
Results in Mathematics. 2023. Vol.~78. Article~161.
Theorem~1.1 (Gasch\"utz's splitting criterion).
\href{https://doi.org/10.1007/s00025-023-01927-1}{doi:10.1007/s00025-023-01927-1}.

\bibitem{Jones}
Jones G.\,A.
Primitive permutation groups containing a cycle.
Bulletin of the Australian Mathematical Society.
2014. Vol.~89. P.~159--165.
Theorem~1.1.
\href{https://doi.org/10.1017/S000497271300049X}
{doi:10.1017/S000497271300049X}.
Open-access version:
\url{https://arxiv.org/abs/1209.5169}.

\bibitem{Atlas}
Conway J.\,H., Curtis R.\,T., Norton S.\,P.,
Parker R.\,A., Wilson R.\,A.
Atlas of Finite Groups.
Oxford: Clarendon Press, 1985.
Online group pages:
\url{https://brauer.maths.qmul.ac.uk/Atlas/v3/alt/A5/},
\url{https://brauer.maths.qmul.ac.uk/Atlas/alt/A6/}.

\bibitem{Zav}
Zavarnitsine A.\,V.
Finite simple groups with narrow prime spectrum.
Siberian Electronic Mathematical Reports.
2009. Vol.~6. P.~1--12. Table~1, rows \(S_5,S_7\).
\url{https://arxiv.org/abs/0810.0568}.

\bibitem{ZhLM}
Zhurtov A.\,Kh., Lytkina D.\,V., Mazurov V.\,D.
On finite groups subspectral to finite almost simple groups.
Vladikavkaz Mathematical Journal.
2025. Vol.~27, no.~3. P.~68--74 (in Russian).
Theorem~1, Lemma~1, Tables~1 and~2.
\href{https://doi.org/10.46698/w4978-1776-4637-t}
{doi:10.46698/w4978-1776-4637-t}.

\bibitem{A9mod2}
GAP Character Table Library / MOC.
Decomposition matrix for \(A_9\) modulo \(2\).
\url{https://www.math.rwth-aachen.de/~MOC/decomposition/tex/A9/A9mod2.pdf}.

\bibitem{A9Matrices}
ATLAS of Finite Group Representations.
The alternating group \(A_9\), representations over \(\mathbb F_2\).
\url{https://brauer.maths.qmul.ac.uk/Atlas/alt/A9/}.
The matrices in the files
\texttt{A9G1-f2r8bB0.m1}, \texttt{A9G1-f2r8bB0.m2},
\texttt{A9G1-f2r8cB0.m1}, \texttt{A9G1-f2r8cB0.m2}
in the \texttt{mtx} subdirectory are used.

\bibitem{CurtisReiner}
Curtis C.\,W., Reiner I.
Representation Theory of Finite Groups and Associative Algebras.
New York: Interscience Publishers, 1962.

\bibitem{Sp6Data}
ATLAS of Finite Group Representations;
GAP Character Table Library; MOC decomposition matrices.
The group \(S_6(2)\): standard generators, classes,
and matrices of the unique \(8\)-dimensional module:
\url{https://brauer.maths.qmul.ac.uk/Atlas/v3/clas/S62/},
\url{https://brauer.maths.qmul.ac.uk/Atlas/v3/matrep/S62G1-f2r8B0}.
Decomposition matrices modulo \(2,3,5,7\) are in the same directory; in particular:
\url{https://www.math.rwth-aachen.de/~MOC/decomposition/tex/S6(2)/S6(2)mod2.pdf},
\url{https://www.math.rwth-aachen.de/~MOC/decomposition/tex/S6(2)/S6(2)mod7.pdf}.
Ordinary table: CTblLib 1.3.11, file
\texttt{data/ctosymp1.tbl}, entry \texttt{S6(2)}.
The embedding \(S_6(2)<O_8^+(2)\) of index \(120\):
\url{https://www.math.rwth-aachen.de/~Thomas.Breuer/ctbllib/ctbltoc/data/O8+(2).html}.

\bibitem{U35Data}
ATLAS of Finite Group Representations; GAP Character Table Library.
\(U_3(5)\): order, Schur multiplier, subgroup
\(A_7\) of index \(50\), and Sylow
\(5\)-subgroup normalizer \(5^{1+2}:8\).
\url{https://brauer.maths.qmul.ac.uk/Atlas/v3/clas/U35/};
\url{https://www.math.rwth-aachen.de/~Thomas.Breuer/ctbllib/ctbltoc/data/U3(5).html}.

\bibitem{U43Data}
ATLAS of Finite Group Representations; GAP Character Table Library; MOC.
\(U_4(3)\): subgroups \(A_7\), ordinary table,
decomposition matrix modulo \(2\), and the representation
of dimension \(21\) with character \(\chi_2\).
\url{https://www.math.rwth-aachen.de/~Thomas.Breuer/ctbllib/ctbltoc/data/U4(3).html};
\url{https://www.math.rwth-aachen.de/~MOC/decomposition/tex/U4(3)/U4(3)mod2.pdf};
\url{https://brauer.maths.qmul.ac.uk/Atlas/v3/matrep/U43G1-Zr21B0}.
Matrices: \texttt{U43G1-Zr21B0.g} in the
\texttt{Atlas/clas/U43/gap0} subdirectory; ordinary characters:
CTblLib 1.3.11, \texttt{data/ctounit1.tbl}.

\bibitem{NilpotentData}
GAP Character Table Library; ATLAS; MOC decomposition matrices.
Ordinary tables of \(L_2(7)=L_3(2),L_2(8),J_2,U_3(3),L_3(4)\)
and the normalizer \(5:4<U_4(2)\):
\url{https://www.math.rwth-aachen.de/~Thomas.Breuer/ctbllib/ctbltoc/data/L3(2).html},
\url{https://www.math.rwth-aachen.de/~Thomas.Breuer/ctbllib/ctbltoc/data/L2(8).html},
\url{https://www.math.rwth-aachen.de/~Thomas.Breuer/ctbllib/ctbltoc/data/U4(2).html}.
Decomposition matrices:
\url{https://www.math.rwth-aachen.de/~MOC/decomposition/tex/J2/J2mod3.pdf},
\url{https://www.math.rwth-aachen.de/~MOC/decomposition/tex/U3(3)/U3(3)mod3.pdf},
\url{https://www.math.rwth-aachen.de/~MOC/decomposition/tex/L3(4)/L3(4)mod2.pdf},
\url{https://www.math.rwth-aachen.de/~MOC/decomposition/tex/L3(4)/L3(4)mod3.pdf}.
Multipliers of \(A_7,A_8\) and automorphisms of \(L_3(4)\):
\url{https://brauer.maths.qmul.ac.uk/Atlas/alt/A7/},
\url{https://brauer.maths.qmul.ac.uk/Atlas/alt/A8/},
\url{https://brauer.maths.qmul.ac.uk/Atlas/v3/lin/L34/}.

\bibitem{Steinberg}
Steinberg R.
Representations of Algebraic Groups. Nagoya Mathematical Journal.
1963. Vol.~22. P.~33--56. Theorems~1.1, 1.3.
\href{https://doi.org/10.1017/S0027763000011016}
{doi:10.1017/S0027763000011016}.

\bibitem{StaroletovPrimeGraph}
Staroletov A.\,M.
On recognition of alternating groups by prime graph.
Siberian Electronic Mathematical Reports.
2017. Vol.~14. P.~994--1010. Lemma~3.
\href{https://doi.org/10.17377/semi.2017.14.084}
{doi:10.17377/semi.2017.14.084}.
\url{https://www.mathnet.ru/php/getFT.phtml?jrnid=semr&option_lang=rus&paperid=841&what=fullt}.

\bibitem{MixedData}
ATLAS of Finite Group Representations;
GAP Character Table Library / MOC.
Multipliers, subgroups, and modules:
\url{https://brauer.maths.qmul.ac.uk/Atlas/lin/L28/};
\url{https://brauer.maths.qmul.ac.uk/Atlas/clas/U42/};
\url{https://brauer.maths.qmul.ac.uk/Atlas/v3/clas/U33/};
\url{https://brauer.maths.qmul.ac.uk/Atlas/v3/lin/L34/};
\url{https://brauer.maths.qmul.ac.uk/Atlas/alt/A8/}.
The full decomposition matrix, including the blocks of \(2.A_8\):
\url{https://www.math.rwth-aachen.de/~MOC/decomposition/tex/A8/A8mod5.pdf}.
The matrices of \(2.A_8\) from Lemma~\ref{lem:double-A8-five}
are in the \texttt{Atlas/alt/A8/mtx} subdirectory.
The subgroup \(A_6<U_4(2)\) and normalizer \(5{:}4\):
\url{https://www.math.rwth-aachen.de/~Thomas.Breuer/ctbllib/ctbltoc/data/U4(2).html}.

\bibitem{TwoLocalData}
ATLAS of Finite Group Representations;
GAP Character Table Library.
The subgroups \(L_2(7),2^3{:}7{:}3<2.A_8\),
transitive actions of degree \(240\), and the full list
of simple \(A_8\)-modules in characteristic \(2\):
\url{https://brauer.maths.qmul.ac.uk/Atlas/alt/A8/}.
The subgroup \(C_2\times U_3(3)<2.J_2\):
\url{https://www.math.rwth-aachen.de/~Thomas.Breuer/ctbllib/ctbltoc/data/2.J2.html}.

\bibitem{ThompsonCritical}
Yi X. Characterizations of hypercyclically embedded subgroups
of finite groups. Rendiconti del Seminario Matematico
della Università di Padova. 2016. Vol.~135. P.~195--206.
Lemma~2.3 (Thompson's critical subgroup theorem).
\href{https://doi.org/10.4171/RSMUP/135-11}{doi:10.4171/RSMUP/135-11}.
\url{https://ems.press/content/serial-article-files/43911}.

\bibitem{A7LocalData}
ATLAS of Finite Group Representations;
GAP Character Table Library.
Full lists of simple modules for \(A_7\) in
characteristic \(2\), for \(2.A_7\) in characteristics
\(3,5\), and for \(3.A_7\) in characteristic \(2\):
\url{https://brauer.maths.qmul.ac.uk/Atlas/alt/A7/}.
The subgroup \(2.L_3(2)\) and Sylow
\(5\)-subgroup normalizer \(5{:}8\) in \(2.A_7\):
\url{https://www.math.rwth-aachen.de/~Thomas.Breuer/ctbllib/ctbltoc/data/2.A7.html}.

\bibitem{RemainingCoverData}
ATLAS of Finite Group Representations.
Orders, multipliers, and outer automorphisms of
\(A_5,A_6,L_2(7),L_3(4)\); the natural realization
\(3.L_3(4)=\operatorname{SL}_3(4)\):
\url{https://brauer.maths.qmul.ac.uk/Atlas/alt/A5/},
\url{https://brauer.maths.qmul.ac.uk/Atlas/alt/A6/},
\url{https://brauer.maths.qmul.ac.uk/Atlas/lin/L32/},
\url{https://brauer.maths.qmul.ac.uk/Atlas/lin/L34/L34.html}.
Facts about field and graph automorphisms of
\(L_3(4)\) are used in the same form as in
Lemma~\ref{lem:normal-L34-no9}; the matrix analysis
of centralizers in Lemma~\ref{lem:normal-cover-L34-no9}
is carried out directly.

\bibitem{A5ThreeData}
ATLAS of Finite Group Representations.
The full list of simple \(2.A_5\)-modules
in characteristic \(3\) with nontrivial central action:
\url{https://brauer.maths.qmul.ac.uk/Atlas/alt/A5/}.
The normalizer of a Sylow \(Q_8\) in \(2.A_5\):
\url{https://www.math.rwth-aachen.de/~Thomas.Breuer/ctbllib/ctbltoc/data/2.A5.html}.
The matrix of a standard element of order \(3\)
on the six-dimensional module:
\url{https://brauer.maths.qmul.ac.uk/Atlas/alt/A5/mtx/2A5G1-f3r6B0.m2}.

\bibitem{ConradQuaternion}
Conrad K. Generalized quaternions.
Theorem~4.10: classification of finite
\(p\)-groups with a unique subgroup
of order \(p\); properties and automorphisms
of generalized quaternion groups.
\url{https://kconrad.math.uconn.edu/blurbs/grouptheory/genquat.pdf}.

\bibitem{TransferNotes}
Michigan State University. Group Theory, MTH912.
Lecture notes, Section~1.5, Lemma~1.5.6:
the normal complement criterion and control
of fusion in an abelian Sylow subgroup.
\url{https://users.math.msu.edu/users/meierfra/Classnotes/MTH912F04/912F04master.pdf}.

\bibitem{A5SevenData}
ATLAS of Finite Group Representations.
The full list of simple \(2.A_5\)-modules
in characteristic \(7\) with nontrivial central action:
\url{https://brauer.maths.qmul.ac.uk/Atlas/alt/A5/}.
Matrix models used for the checks in
Lemma~\ref{lem:double-A5-seven-simple-data}:
\texttt{2A5G1-f7r4aB0.m2},
\texttt{2A5G1-f7r6B0.m1},
\texttt{2A5G1-f7r6B0.m2}
in the \texttt{mtx} subdirectory of that page.

\bibitem{A5TwoData}
ATLAS of Finite Group Representations.
The full list of simple \(A_5\)-module types in characteristic \(2\):
two conjugate two-dimensional types over \(\mathbb F_4\),
an absolutely simple four-dimensional type over \(\mathbb F_2\),
and a four-dimensional realization of the conjugate pair over \(\mathbb F_2\).
\url{https://brauer.maths.qmul.ac.uk/Atlas/alt/A5/}.
Completeness and the required structural properties
are justified in Lemma~\ref{lem:A5-orthogonal-small-five-rank}.

\bibitem{A5BinaryLocalData}
ATLAS of Finite Group Representations.
A faithful simple \(2.A_5\)-module of dimension \(4\)
over \(\mathbb F_3\), which splits over \(\mathbb F_9\).
Matrices of the standard generators:
\url{https://brauer.maths.qmul.ac.uk/Atlas/alt/A5/mtx/2A5G1-f3r4B0.m1},
\url{https://brauer.maths.qmul.ac.uk/Atlas/alt/A5/mtx/2A5G1-f3r4B0.m2}.
The matrices, group order, perfection,
orbits, and required spectra are also checked
by the certificate in Appendix~\ref{app:binary-upper-model}.

\bibitem{L27BinaryData}
ATLAS of Finite Group Representations; GAP Character Table Library.
The group \(L_2(7)=L_3(2)\), the cover
\(2.L_2(7)=\operatorname{SL}_2(7)\),
automorphisms, and the projective extension:
\url{https://brauer.maths.qmul.ac.uk/Atlas/lin/L27/}.
Character-table data for the double cover:
\url{https://www.math.rwth-aachen.de/~Thomas.Breuer/ctbllib/ctbltoc/data/2.L3(2).html}.
Ordinary degrees and the decomposition matrix
in characteristic \(3\) in the MOC table:
\url{https://www.math.rwth-aachen.de/~MOC/decomposition/tex/L3(2)/L3(2)mod3.pdf}.
Characteristic \(5\) does not divide the orders
of \(L_2(7)\) and \(2.L_2(7)\).

\bibitem{BurnsideThree}
Ishida Y., Mizuno R., Takeuchi K. Existence of a Model Companion for
Groups of Exponent 3. 2026. Fact~2.15.
\url{https://arxiv.org/abs/2609.30061}.

\bibitem{LeviWaerden}
Levi F., van der Waerden B.\,L. Über eine besondere Klasse von Gruppen.
Abhandlungen aus dem Mathematischen Seminar der Universität Hamburg.
1933. Vol.~9. P.~154--158.
\href{https://doi.org/10.1007/BF02940639}{doi:10.1007/BF02940639}.

\bibitem{L34CoverSubgroupData}
ATLAS of Finite Group Representations; GAP Character Table Library.
The multiplier and covers of \(L_3(4)\), and full preimages of
\(A_6,L_2(7)\) in \(2.L_3(4),2^2.L_3(4),3.L_3(4)\):
\url{https://brauer.maths.qmul.ac.uk/Atlas/lin/L34/L34.html};
\url{https://www.math.rwth-aachen.de/~Thomas.Breuer/ctbllib/ctbltoc/data/2.L3(4).html};
\url{https://www.math.rwth-aachen.de/~Thomas.Breuer/ctbllib/ctbltoc/data/2^2.L3(4).html};
\url{https://www.math.rwth-aachen.de/~Thomas.Breuer/ctbllib/ctbltoc/data/3.L3(4).html}.

\bibitem{A6TripleFiveFixedData}
GAP Character Table Library, version 1.3.11.
The tables \texttt{A6}, \texttt{3.A6}, and projective characters
in the file \texttt{data/ctoalter.tbl}.
The degrees and values on the classes \(5A,5B\) are used;
the full multiplicity calculation is given in
Lemma~\ref{lem:A6-triple-seven-five-fixed}.
\url{https://www.math.rwth-aachen.de/~Thomas.Breuer/ctbllib/ctbltoc/data/3.A6.html}.

\bibitem{CriticalOddFull}
Chen X., Guo W. On \(\Pi\)-supplemented subgroups of a finite group.
Lemmas~2.5, 2.7, 2.10: the full form of Thompson's critical subgroup
theorem, its \(\Omega_1\), and \(p'\)-automorphisms.
\url{https://arxiv.org/abs/1307.0089}, version~4.
Published in: Communications in Algebra. 2016. Vol.~44. P.~731--745.

\bibitem{L34FinalData}
Breuer T. GAP Character Table Library, version 1.3.11.
The ordinary tables \texttt{2.L3(4)}, \texttt{2.L3(4).2\_1},
\texttt{2.L3(4).2\_2}, \texttt{2.L3(4).2\_3},
and their isoclinic variants: file \texttt{ctoline3.tbl};
the Brauer table of \texttt{2.L3(4)} modulo \(3\):
\texttt{ctbline3.tbl}.
Decomposition matrix:
\url{https://www.math.rwth-aachen.de/~MOC/decomposition/tex/L3(4)/L3(4)mod3.pdf}.
Completeness of the outer-extension variants and the isoclinic construction:
Computations with the GAP Character Table Library,
\S\S2.2.4, 2.6.4, 2.7.4.
\url{https://www.math.rwth-aachen.de/~Thomas.Breuer/ctbllib/doc2/chap2.html}.
The double-cover page:
\url{https://www.math.rwth-aachen.de/~Thomas.Breuer/ctbllib/ctbltoc/data/2.L3(4).html}.

\bibitem{A6FinalData}
GAP Character Table Library 1.3.11.
Ordinary tables of \(A_6,2.A_6,3.A_6,6.A_6,A_6.2_i\):
file \texttt{ctoalter.tbl}; Brauer tables:
\texttt{ctbalter.tbl}.
\url{https://www.math.rwth-aachen.de/~Thomas.Breuer/ctbllib/ctbltoc/data/A6.html}.
Decomposition matrices:
\url{https://www.math.rwth-aachen.de/~MOC/decomposition/tex/A6/A6mod2.pdf};
\url{https://www.math.rwth-aachen.de/~MOC/decomposition/tex/A6/A6mod3.pdf}.
The outer type \(A_6.2_2=\operatorname{PGL}_2(9)\):
\url{https://www.math.rwth-aachen.de/~Thomas.Breuer/ctbllib/ctbltoc/data/A6.2_2.html}.

\end{thebibliography}
\end{document}